%% file: __the_exponential_world.tex
\documentclass[11pt,oneside]{book}
\usepackage{geometry}                		
\usepackage{graphicx}									
\usepackage{amssymb}\vspace{0.25cm}
\usepackage{mathtools}
\usepackage{setspace}
\usepackage{tikz-cd}
\usepackage[mathscr]{euscript}
\newcommand{\abs}[1]{\left\vert#1\right\vert}
\usepackage[normalem]{ulem}
\usepackage{manfnt}
\usepackage{pgfplots}
\usepackage{pifont}
\usepackage[safe]{tipa}
\usepackage{marvosym}
\usepackage{scalerel,stackengine}
\usepackage{titlesec}
\usepackage{makeidx}
\usepackage{centernot}
\usepackage{xspace}
\usepackage[noauto]{chappg}
\usepackage[OT2,T1]{fontenc}
\usepackage{hyperref}
\usepackage{amsmath}
\usepackage[amsmath,thmmarks]{ntheorem}
\usepackage{tabularx}
\usepackage{yfonts}
\newcolumntype{Y}{>{\raggedright\arraybackslash}X}

\newcommand\mA{%
$A$\xspace
}
\newcommand\mB{%
$B$\xspace
}
\newcommand\mC{%
$C$\xspace
}
\newcommand\mD{%
$D$\xspace
}
\newcommand\mE{%
$E$\xspace
}
\newcommand\mF{%
$F$\xspace
}
\newcommand\mG{%
$G$\xspace
}

\newcommand\mJ{%
$J$\xspace
}

\newcommand\mL{%
$L$\xspace
}
\newcommand\mM{%
$M$\xspace
}
\newcommand\mN{%
$N$\xspace
}

\newcommand\mP{%
$P$\xspace
}
\newcommand\mQ{%
$Q$\xspace
}
\newcommand\Qbar{%
\overline{\mathbb{Q}\raisebox{.29cm}{\hspace{0pt}}}\xspace
}
\newcommand\mR{%
$R$\xspace
}
\newcommand\mS{%
$S$\xspace
}
\newcommand\mT{%
$T$\xspace
}
\newcommand\mU{%
$U$\xspace
}
\newcommand\mV{%
$V$\xspace
}

\newcommand\mX{%
$X$\xspace
}
\newcommand\mY{%
$Y$\xspace
}

\newcommand\mfc{%
\mathfrak{c}\xspace
}

\newcommand\C{%
\mathbb{C}
}
\newcommand\E{%
\mathbb{E}
}
\newcommand\F{%
\mathbb{F}
}

\newcommand\K{%
\mathbb{K}
}
\newcommand\LL{%
\mathbb{L}
}

\newcommand\N{%
\mathbb{N}
}
\newcommand\PP{%
\mathbb{P}
}
\newcommand\Q{%
\mathbb{Q}
}
\newcommand\R{%
\mathbb{R}
}
\newcommand\T{%
\mathbb{T}
}

\newcommand\Z{%
\mathbb{Z}
}

\newcommand\ba{%
\textbf{a}\xspace
}

\newcommand\bb{%
\textbf{b}
}

\newcommand\be{%
\textbf{e}\xspace
}

\newcommand\bk{%
\textbf{k}\xspace
}
\newcommand\bL{%
\textbf{L}\xspace
}

\newcommand\bq{%
\textbf{q}\xspace
}

\newcommand\bss{%
\textbf{s}\xspace
}
\newcommand\bu{%
\textbf{u}\xspace
}
\newcommand\bv{%
\textbf{v}\xspace
}
\newcommand\bw{%
\textbf{w}\xspace
}
\newcommand\bx{%
\textbf{x}\xspace
}
\newcommand\by{%
\textbf{y}\xspace
}
\newcommand\bz{%
\textbf{z}\xspace
}
\newcommand\bzw{%
\textbf{zw}\xspace
}

\newcommand\elln{%
\ell\text{n}
}
\newcommand\ecl{%
\text{ec}\ell
}
\newcommand\lindim{%
\text{lindim}
}
\newcommand\Log{%
\text{Log}\hspace{0.05cm}
}

\newcommand\sB{%
\mathcal{B}
}

\newcommand\sD{%
\mathcal{D}
}
\newcommand\sE{%
\mathcal{E}
}

\newcommand\sH{%
\mathcal{H}
}
\newcommand\sM{%
\mathcal{M}
}

\newcommand\sP{%
\mathcal{P}
}

\newcommand\sV{%
\mathcal{V}
}
\newcommand\sW{%
\mathcal{W}
}

\newcommand\EXP{%
\text{EXP}\hspace{0.05cm}
}
\newcommand\EX{%
\text{EX}\hspace{0.05cm}
}

\newcommand\Ker{%
\text{Ker}\hspace{0.03cm}
}

\newcommand\spanxs[1]{%
\text{span}_{\hspace{0.03cm}{#1}}\hspace{0.05cm}
}

\newcommand\diag{%
\text{diag}\hspace{0.05cm}
}

\newcommand\Gal{%
\text{Gal}
}

\newcommand\glb{%
\text{glb}\hspace{0.05cm}
}
\newcommand\lub{%
\text{lub}\hspace{0.05cm}
}

\newcommand\PV{%
\text{PV}
}
\newcommand\rank{%
\text{rank}\hspace{0.05cm}
}
\newcommand\spt{%
\text{spt}\hspace{0.05cm}
}

\newcommand\id{%
\text{id}
}

\newcommand\nth{%
\text{th}
}
\newcommand\tD{%
\text{D}
}
\newcommand\td{%
\text{d}
}
\newcommand\tH{%
\text{H}
}

\newcommand\tL{%
\text{L}
}
\newcommand\tmul{%
\text{mul}
}
\newcommand\tN{%
\text{N}
}

\newcommand\tS{%
\text{S}
}

\newcommand\tU{%
\text{U}
}

\newcommand\trdeg{%
\text{trdeg}
}
\newcommand\trdegQ{%
\text{trdeg}_{\mathbb{Q}} \hspace{0.05cm}
}
\newcommand\trdegQbar{%
\text{trdeg}_{\overline{\mathbb{Q}}} \hspace{0.1cm}
}
\newcommand\Tee{%
\mathsf{T}
}

\newcommand\ra{%
\rightarrow
}
\newcommand\lra{%
\longrightarrow
}

\newcommand\ds{%
\displaystyle
}

\newcommand\fL{%
\frak{L}
}

\newcommand\Imx{%
\text{Im}\hspace{0.05cm}
}
\newcommand\Rex{%
\text{Re}\hspace{0.05cm}
}
\newcommand\Arg{%
\text{Arg}\hspace{0.05cm}
}
\newcommand\card{%
\text{card}\hspace{0.05cm}
}

\newcommand\modx{%
\text{mod}\hspace{0.05cm}
}

\newcommand\tr{%
\text{tr}\hspace{0.05cm}
}

\newcommand\un[1]{%
\underline{#1}\xspace
}

\newcommand\ov[1]{%
\overline{#1}
}

\newcommand\restr[2]{%
{#1}|{#2}
}

\newcommand\fraOneAlphaone{%
\Large\text{\big($\frac{1}{\alpha_1}$\big)}
}
\newcommand\LogfraOneAlphaone{%
 \Large\text{\large{Log} \Large{\big($\frac{1}{\alpha_1}$\big)}}
}

\newcommand\hsx{%
\hspace{0.05cm}
}
\newcommand\hsy{%
\hspace{0.03cm}
}
\newcommand\hsm{%
\hspace{0.3cm}
}

\newcommand\bzero{%
{\boldsymbol{0}}
}

\newcommand\bcdot{%
{\boldsymbol{\ \cdot \ }}
}
\newcommand\bbeta{%
{\boldsymbol{\beta}}
}
\newcommand\bomega{%
{\boldsymbol{\omega}}
}
\newcommand\blambda{%
{\boldsymbol{\lambda}}
}

\newcommand\strutx{\rule{0pt}{13pt}}
\newcommand\struty{\rule{0pt}{11pt}}
\newcommand\strutz{\rule{0pt}{8pt}}
\newcommand\strutv{\rule{0pt}{8pt}} 
\newcommand\strutw{\rule{0pt}{9pt}} 

\newcommand\upx{\raisebox{.42cm}{\hspace{0pt}}}
\newcommand\upy{\raisebox{.35cm}{\hspace{0pt}}}

\newif\ifcomments
\commentstrue
\newif\ifdmc
\dmcfalse
\newif\ifdmcx
\dmcxtrue
\newif\ifgw
\gwfalse

\usepackage{scalerel,stackengine}
\stackMath
\newcommand\reallywidehat[1]{%
\savestack{\tmpbox}{\stretchto{%
  \scaleto{%
    \scalerel*[\widthof{\ensuremath{#1}}]{\kern-.6pt\bigwedge\kern-.6pt}%
    {\rule[-\textheight/2]{1ex}{\textheight}}
  }{\textheight}%
}{0.5ex}}%
\stackon[1pt]{#1}{\tmpbox}%
}
\newcommand*{\wtilde}[2][0pt]{
  \raisebox{#1}{$\widetilde{\raisebox{-#1}{$#2$}}$}%
}%

\makeatletter
\DeclareRobustCommand\widecheck[1]{{\mathpalette\@widecheck{#1}}}
\def\@widecheck#1#2{%
    \setbox\z@\hbox{\m@th$#1#2$}%
    \setbox\tw@\hbox{\m@th$#1%
       \widehat{%
          \vrule\@width\z@\@height\ht\z@
          \vrule\@height\z@\@width\wd\z@}$}%
    \dp\tw@-\ht\z@
    \@tempdima\ht\z@ \advance\@tempdima2\ht\tw@ \divide\@tempdima\thr@@
    \setbox\tw@\hbox{%
       \raise\@tempdima\hbox{\scalebox{1}[-1]{\lower\@tempdima\box
\tw@}}}%
    {\ooalign{\box\tw@ \cr \box\z@}}}
\makeatother

\newtheoremstyle{xx}
  {4pt}
  {0pt}
  {\upshape}
  {\bfseries}
  {}
  { }
  {}
  
\makeatletter 
 \newtheoremstyle{myu}%
  {\upshape\item[ \indent\indent\bf\underline{\theorem@headerfont ##2:}]}%
\makeatother
\makeatletter 
 \newtheoremstyle{myn}%
  {\item[\hskip\labelsep \ \bf ##1 \theorem@headerfont ##2.]}%
\makeatother
\theoremstyle{myn}
\newtheorem{theoremn}{Theorem} 
\theoremstyle{myu}
{\upshape}
\newtheorem{x}[theoremn]{}

 \newtheoremstyle{mr}%
  {\upshape\item[ \indent{\theorem@headerfont ##2. \hspace{.2cm}}]}%
\makeatother
\theoremstyle{mr}
{\upshape}
\newtheorem{rf}[theoremn]{}

\title{The Exponential World}
\author{Garth Warner (Emeritus)\\
Department of Mathematics\\
University of Washington}
\date{}									

\titleformat{\chapter}[display]
{\normalfont\filcenter\huge\bfseries}{}{0pt}{\large}

\titleformat{\chapter}[display]
{\normalfont\filcenter\huge\bfseries}{}{0pt}{\large}

\usepackage[OT2,OT1]{fontenc} 
\newcommand\cyr
{
\renewcommand\rmdefault{wncyr} 
\renewcommand\sfdefault{wncyss} 
\renewcommand\encodingdefault{OT2} 
\normalfont
\selectfont
}

\DeclareTextFontCommand{\textcyr}{\cyr}

\makeindex 
\begin{document}

\maketitle                              

\titlespacing*{\chapter}{0pt}{-50pt}{40pt}
\setlength{\parskip}{0.1em}
\renewcommand{\thepage}{\roman{page}}
\include{__abstract}


\include{__tocX}
\pagenumbering{bychapter}
\setcounter{chapter}{0}
\include{_00_the_canonical_estimate}%
\include{_01_ordered_sets}%
\include{_02_real_numbers}

\include{_03_suprema}%
\include{_04_exponents_and_roots}%
\include{_05_exp_a_and_log_a}%
\include{_06_irrationality_of_sqrt_2}%
\include{_07_irrationality_theory_and_examples}%
\include{_08_irrationality_of_e}%
\include{_09_irrationality_of_e_a_div_b}%
\include{_10_irrationality_of_e_a_div_b_bis}%
\include{_11_irrationality_of_pi}%
\include{_12_irrationality_of_cos_x}%
\include{_13_irrationality_of_cosh_x}%
\include{_14_algebraic_and_transcendental_numbers}%
\include{_15_liouville_theory}%
\include{_16_the_mahler_classification}%
\include{_17_transcendence_of_e}%
\include{_18_symmetric_algebra}%
\include{_19_transcendence_of_pi}%
\include{_20_algebraic_in_dependence}%
\include{_21_the_lindeman_weierstrass_theorem}%
\include{_22_exceptional_sets}%
\include{_23_complex_logarithms_and_complex_powers}%
\include{_24_the_gelfond_schneider_theorem}%
\include{_25_interpolation_determinants}%
\include{_26_zero_estimates}%
\include{_27_gelfond_schneider_setting_stage}%
\include{_28_gelfond_schneider_execution}%
\include{_29_schneider_lang_criterion}%
\include{_30_schneider_lang_criteria}%
\include{_31_baker_statement}%
\include{_32_equivalences}%
\include{_33_baker_proof}%
\include{_34_estimates}%
\include{_35_matrices}%
\include{_36_the_six_exponentials_theorem}%
\include{_37_vector_spaces}%
\include{_38_vector_spaces_L}%
\include{_39_vector_spaces_LG}%
\include{_40_vector_spaces_Vmax_V_min}%
\include{_41}%
\include{_42_sharp_six_exp_theorem}%
\include{_43_strong_six_exp_theorem}%
\include{_44_the_four_exponentials_conjecture}%
\include{_45_strong_four_exp_conjecture_s4ec}%
\include{_46_transcendental_extensions}%
\include{_47_schanuels_conjecture}%
\include{_48_numerical_examples}%
\include{_49_the_zero_condition}%
\include{_50_property_matrix_abc_}%
\include{_51_vector_spaces_L_bis}%
\include{_52_on_the_equation_z_plus_exp_z}%
\include{_53_on_the_equation_P_z_exp_z}

\include{_54_zilber_fields}%
\include{_55_e_rings}%
\include{_56_shanuel_implies_shapiro}%
\include{_57_shapiro_conjecture_case_1}%
\include{_58_shapiro_conjecture_case_2}%
\include{_59_differential_algebra}%
\include{_60_formal_schanuel}%
\include{_61_an_arithmetic_criterion}%
\include{_62_real_numbers_bis}%
\include{_XX_Supp_trans_of_series}%
\include{_01_canonical_illustrations}%
\include{_02_the_role_of_the_cotangent}%
\include{_03_application_of_nesterenko}%
\include{_04_introduction_of_schc}%
\include{_05_introduction_of_schc_bis}%
\include{_06_consolidation}%
\include{_07_consideration_of_a_over_b}%
\include{_08_an_algebraic_series}%
\include{_XX_Supp_zeta_function_values}

\include{_01_Bernoulli_Numbers}

\include{_02_Zeta_2n}

\include{_03_Zeta_2}

\include{_04_Zeta_2_bis}

\include{_05_Zeta_2n_bis}

\include{_06_Zeta_3}

\include{_07_Conjugate_Bernoulli_Numbers}

\include{_08_Zeta_2n_1}

\include{__refs}
\setcounter{page}{1}
\renewcommand{\thepage}{Index-\arabic{page}}
\backmatter
\bibliography{}
\printindex
\end{document}

%% file: __abstract.tex
\chapter{
ABSTRACT}
\setlength\parindent{2em}
\setcounter{theoremn}{0}

\ \indent 

In this book there will be found an introduction to transcendental number theory, starting at the beginning and ending at the frontiers.  The emphasis is on the conceptual aspects of the subject, thus the effective theory has been more or less completely ignored, as has been the theory of \mE-functions and \mG-functions.  Still, a fair amount of ground is covered and while I take certain results without proof, this is done primarily so as not to get bogged down in technicalities, otherwise the exposition is detailed and little is left to the reader.


\chapter{
ACKNOWLEDGEMENT}
\setlength\parindent{2em}
\setcounter{theoremn}{0}

\ \indent 

My thanks to Judith Clare for a superb job of difficult technical typing.\\

Recently David Clark converted the typewritten manuscript to AMS-TeX.  
This was a monumental task and in so doing he made a number of constructive comments and useful suggestions which serve to enhance the exposition.  
His careful scrutiny of the manuscript has been invaluable. 
\\[3cm]

\[
\textbf{DEDICATION}
\]

This article is dedicated to the memory of Paul Sally.


%% file: __tocX.tex
\centerline{\textbf{\LARGE CONTENTS}}
\vspace{0.5cm}

\allowdisplaybreaks
\begin{align*}
\S0.  \qquad &\text{THE CANONICAL ESTIMATE} 
\\[10pt]
\S1.  \qquad &\text{ORDERED SETS} 
\\[10pt]
\S2.  \qquad &\text{REAL NUMBERS}  
\\[10pt]
\S3.  \qquad &\text{SUPREMA }
\\[10pt]
\S4.  \qquad &\text{EXPONENTS AND ROOTS}
\\[10pt]
\S5.  \qquad &\text{$\exp_a$ AND $\log_a$}
\\[10pt]
\S6.  \qquad &\text{IRRATIONALITY OF $\sqrt{2}$}
\\[10pt]
\S7.  \qquad &\text{IRRATIONALITY: THEORY AND EXAMPLES}
\\[10pt]
\S8.  \qquad &\text{IRRATIONALITY OF \be}
\\[10pt]
\S9.  \qquad &\text{IRRATIONALITY OF $\be^{\ba/\bb}$}
\\[10pt]
\S10.  \qquad &\text{IRRATIONALITY OF $\be^{\ba/\bb}$ (bis)}
\\[10pt]
\S11.  \qquad &\text{IRRATIONALITY OF $\boldsymbol{\pi}$}
\\[10pt]
\S12.  \qquad &\text{IRRATIONALITY OF \textbf{cos(x)}}
\\[10pt]
\S13.  \qquad &\text{IRRATIONALITY OF \textbf{cosh(x)}}
\\[10pt]
\S14.  \qquad &\text{ALGEBRAIC AND TRANSCENDENTAL NUMBERS}
\\[10pt]
\S15.  \qquad &\text{LIOUVILLE THEORY}
\\[10pt]
\S16.  \qquad &\text{THE MAHLER CLASSIFICATION}
\\[10pt]
\S17.  \qquad &\text{TRANSCENDENCE OF \be}
\\[10pt]
\S18.  \qquad &\text{SYMMETRIC ALGEBRA}
\\[10pt]
\S19.  \qquad &\text{TRANSCENDENCE OF $\boldsymbol{\pi}$}
\\[10pt]
\S20.  \qquad &\text{ALGEBRAIC (IN) DEPENDENCE}
\\[10pt]
\S21.  \qquad &\text{THE LINDEMANN-WEIERSTRASS THEOREM} 
\\[10pt]
\S22.  \qquad &\text{EXCEPTIONAL SETS}  
\\[10pt]
\S23.  \qquad &\text{COMPLEX LOGARITHMS AND COMPLEX POWERS}
\\[10pt]
\S24.  \qquad &\text{THE GELFOND-SCHNEIDER THEOREM}
\\[10pt]
\S25.  \qquad &\text{INTERPOLATION DETERMINANTS}
\\[10pt]
\S26.  \qquad &\text{ZERO ESTIMATES}
\\[10pt]
\S27.  \qquad &\text{GELFOND-SCHNEIDER: SETTING THE STAGE}
\\[10pt]
\S28.  \qquad &\text{GELFOND-SCHNEIDER: EXECUTION}
\\[10pt]
\S29.  \qquad &\text{THE SCHNEIDER-LANG CRITERION}
\\[10pt]
\S30.  \qquad &\text{SCHNEIDER-LANG CRITERIA}
\\[10pt]
\S31.  \qquad &\text{BAKER: STATEMENT} 
\\[10pt]
\S32.  \qquad &\text{EQUIVALENCES}  
\\[10pt]
\S33.  \qquad &\text{BAKER: PROOF}
\\[10pt]
\S34.  \qquad &\text{ESTIMATES}
\\[10pt]
\S35.  \qquad &\text{MATRICES}
\\[10pt]
\S36.  \qquad &\text{SIX EXPONENTIALS THEOREM}
\\[10pt]
\S37.  \qquad &\text{VECTOR SPACES}
\\[10pt]
\S38.  \qquad &\text{VECTOR SPACES: \mL}
\\[10pt]
\S39.  \qquad &\text{VECTOR SPACES: $L_G$}
\\[10pt]
\S40.  \qquad &\text{VECTOR SPACES: $\sV_{\max}$, $\sV_{\min}$}
\\[10pt]
\S41.  \qquad &\text{EXPONENTIALS (6 or 5)} 
\\[10pt]
\S42.  \qquad &\text{SHARP SIX EXPONENTIALS THEOREM}  
\\[10pt]
\S43.  \qquad &\text{STRONG SIX EXPONENTIALS THEOREM}
\\[10pt]
\S44.  \qquad &\text{FOUR EXPONENTIALS CONJECTURE (4EC)}
\\[10pt]
\S45.  \qquad &\text{STRONG FOUR EXPONENTIALS CONJECTURE (S4EC)}
\\[10pt]
\S46.  \qquad &\text{TRANSCENDENTAL EXTENSIONS}
\\[10pt]
\S47.  \qquad &\text{SCHANUEL'S CONJECTURE}
\\[10pt]
\S48.  \qquad &\text{NUMERICAL EXAMPLES}
\\[10pt]
\S49.  \qquad &\text{THE ZERO CONDITION}
\\[10pt]
\S50.  \qquad &\text{PROPERTY $\huge\bf{\binom{ A \ B}{C \ 0}}$}
\\[10pt]
\S51.  \qquad &\text{VECTOR SPACES: \mL (bis)} 
\\[10pt]
\S52.  \qquad &\text{ON THE EQUATION $z + e^z = 0$}  
\\[10pt]
\S53.  \qquad &\text{ON THE EQUATION $P(z,e^z) = 0$ }
\\[10pt]
\S54.  \qquad &\text{ZILBER FIELDS}
\\[10pt]
\S55.  \qquad &\text{\mE-RINGS}
\\[10pt]
\S56.  \qquad &\text{SCHANUEL $\implies$ SHAPIRO}
\\[10pt]
\S57.  \qquad &\text{SHAPIRO'S CONJECTURE: \ CASE 1}
\\[10pt]
\S58.  \qquad &\text{SHAPIRO'S CONJECTURE: \ CASE 2}
\\[10pt]
\S59.  \qquad &\text{DIFFERENTIAL ALGEBRA}
\\[10pt]
\S60.  \qquad &\text{FORMAL SCHANUEL}
\\[10pt]
\S61.  \qquad &\text{AN ARITHMETIC CRITERION} 
\\[10pt]
\S62.  \qquad &\text{REAL NUMBERS (bis)}  
\\[10pt]
\end{align*}
\vspace{0.5cm}


\centerline{\textbf{\large SUPPLEMENT I}}
\vspace{0.5cm}

\centerline{\textbf{TRANSCENDENCE OF SERIES}}
\vspace{0.25cm}

\allowdisplaybreaks
\begin{align*}
\S1.  \qquad &\text{CANONICAL ILLUSTRATIONS} \hspace{4.7cm}
\\[10pt]
\S2.  \qquad &\text{THE ROLE OF THE COTANGENT}  
\\[10pt]
\S3.  \qquad &\text{APPLICATION OF NESTERENKO}
\\[10pt]
\S4.  \qquad &\text{INTRODUCTION OF SCHC}
\\[10pt]
\S5.  \qquad &\text{INTRODUCTION OF SCHC (bis)}
\\[10pt]
\S6.  \qquad &\text{CONSOLIDATION}
\\[10pt]
\S7.  \qquad &\text{CONSIDERATION OF $\frac{A}{B}$}
\\[10pt]
\S8.  \qquad &\text{AN ALGEBRAIC SERIES}
\end{align*}
\vspace{0.5cm}


\centerline{\textbf{\large SUPPLEMENT II}}
\vspace{0.5cm}

\centerline{\textbf{ZETA FUNCTION VALUES}}
\vspace{0.25cm}

\allowdisplaybreaks
\begin{align*}
\S1.  \qquad &\text{BERNOULLI NUMBERS} \hspace{5.75cm}
\\[10pt]
\S2.  \qquad &\text{$\zeta(2n)$}  
\\[10pt]
\S3.  \qquad &\text{$\zeta(2)$}
\\[10pt]
\S4.  \qquad &\text{$\zeta(2)$ (bis)}
\\[10pt]
\S5.  \qquad &\text{$\zeta(2n)$ (bis)}
\\[10pt]
\S6.  \qquad &\text{$\zeta(3)$}
\\[10pt]
\S7.  \qquad &\text{CONJUGATE BERNOULLI NUMBERS}
\\[10pt]
\S8.  \qquad &\text{$\zeta(2n+1)$}
\end{align*}

%% file: _00_the_canonical_estimate.tex
\chapter{
$\boldsymbol{\S}$\textbf{0}.\quad  THE CANONICAL ESTIMATE}
\setlength\parindent{2em}
\setcounter{theoremn}{0}
\renewcommand{\thepage}{\S0-\arabic{page}}

\ \indent

\qquad {\small\bf THEOREM} \ 
Given a positive constant \mC, 
\[
\lim\limits_{n \ra \infty} \hsx \frac{C^n}{n!}
\ = \ 
0.
\]
\vspace{0.2cm}

PROOF \ 
Write
\[
n! 
\ = \ 
n^n e^{-n} \sqrt{n} \hsx \gamma_n \qquad \text{(Stirling's formula).}
\]
Here
\[
\frac{e}{\sqrt{2}} \leq \gamma_n \leq e 
\qquad \bigg( \implies \frac{\sqrt{2}}{e}  \geq \frac{1}{\gamma_n} \geq \frac{1}{e}\bigg).
\]
Choose $n \gg 0$ : $e \hsx C < n$ $-$then 
\begin{align*}
0 < \frac{C^n}{n!} \ 
&=\ \frac{C^n}{n^n e^{-n} \sqrt{n} \hsx \gamma_n}
\\[15pt]
&=\ \frac{\bigl(e C \bigr)^n}{n^n} \hsx \frac{1}{\sqrt{n} \hsx \gamma_n}
\\[15pt]
&\leq\ \bigg(\frac{e C}{n}\bigg)^n \hsx \frac{\sqrt{2}}{e} \hsx \frac{1}{\sqrt{n}}
\\[15pt]
&< \ \frac{\sqrt{2}}{e} \hsx \frac{1}{\sqrt{n}}
\\[15pt]
&\ra 0 \qquad (n \ra \infty).
\end{align*}


%% file: _01_ordered_sets.tex
\chapter{
$\boldsymbol{\S}$\textbf{1}.\quad  ORDERED SETS}
\setlength\parindent{2em}
\setcounter{theoremn}{0}
\renewcommand{\thepage}{\S1-\arabic{page}}

\ \indent 
Let \mX be a nonempty set.
\vspace{0.3cm}

\begin{x}{\small\bf DEFINITION} \ 
An 
\un{order}
\index{order (on a set \mX)} 
on \mX is a relation $<$ with the following properties.
\vspace{0.2cm}

\qquad \textbullet \quad 
\un{Trichotomy} \ 
Given $x, y \in X$, then one and only one of the statements 
\[
x < y, \quad x = y, \quad y < x
\]
is true.
\vspace{0.2cm}

\qquad \textbullet \quad 
\un{Transitivity} \ 
Given $x, y, z \in X$, if $x < y$ and $y < z$, then $x < z$.
\vspace{0.2cm}

\end{x}
\vspace{0.3cm}

\begin{x}{\small\bf \un{N.B.}} \ 
\vspace{0.2cm}

\qquad \textbullet  \quad
$y > x$ means $x < y$.
\vspace{0.2cm}

\qquad \textbullet \quad 
$x \leq y$ means $x < y$ or $x = y$.
\vspace{0.2cm}

\end{x}
\vspace{0.3cm}

\begin{x}{\small\bf DEFINITION} \ 
An 
\un{ordered set}
\index{ordered set} 
is a pair $(X, <)$, where \mX is a nonempty set equipped with an order $<$.
\end{x}
\vspace{0.3cm}

\begin{x}{\small\bf EXAMPLE} \ 
Take $X = \Q$ $-$then \mX is an ordered set if $p < q$ is defined to mean that $q - p$ is positive.
\end{x}
\vspace{0.3cm}

Let \mX be an ordered set, $S \subset X$ a nonempty subset.
\vspace{0.3cm}

\begin{x}{\small\bf NOTATION} \ 
\[
\tU(S) \ = \ \{x \in X \ : \ \forall \ s \in S, \ s \leq x\}.
\]
\end{x}
\vspace{0.3cm}

\begin{x}{\small\bf DEFINITION} \ 
\mS is 
\un{bounded above}
\index{bounded above} 
if $\tU(S) \neq \emptyset$, an element of $\tU(S)$ being called an 
\un{upper bound}
\index{upperbound}  
of \mS.
\end{x}
\vspace{0.3cm}
\begin{x}{\small\bf \un{N.B.}} \ 
The terms ``bounded below'' and ``lower bound'' are to be assigned the obvious interpretations, where now
\[
\tL(S) \ = \ \{x \in X \ : \ \forall \ s \in S, \ x \leq s\}.
\]

\end{x}
\vspace{0.3cm}

Let \mX be an ordered set, $S \subset X$ a nonempty subset such that $\tU(S) \neq \emptyset$.
\vspace{0.3cm}

\begin{x}{\small\bf DEFINITION} \ 
An element $x \in \tU(S)$ is a 
\un{least upper bound}
\index{least upper bound} 
of \mS if $y < x$ $\implies$ $y \notin \tU(S)$.
\end{x}
\vspace{0.3cm}

\begin{x}{\small\bf LEMMA} \ 
Least upper bounds are unique (if they exist at all) and one writes
\[
x \ = \ \lub S \quad \text{or} \quad x \ = \ \sup S \qquad \text{(``supremum'')}.
\]
\vspace{0.2cm}

[Note: \ 
The definition of ``greatest lower bound'' is analogous, such an element being denoted by
\[
x \ = \ \glb S \quad \text{or} \quad x \ = \ \inf S \qquad \text{(``infimum'')}.]
\]
\end{x}
\vspace{0.3cm}

\begin{x}{\small\bf EXAMPLE} \ 
Take $X = \Q$ and let $\tS = \bigg\{\ds\frac{1}{n}: \ n \in \N\bigg\}$ $-$then $\sup S = 1$ is in \mS but $\inf S = 0$ is not in \mS.
\end{x}
\vspace{0.3cm}

Let \mX be an ordered set.
\vspace{0.3cm}

\begin{x}{\small\bf DEFINITION} \ 
\mX has the 
\un{least upper bound}
\index{least upper bound propery} 
property if each nonempty subset $S \subset X$ which is bounded above has a least upper bound.
\end{x}
\vspace{0.3cm}

\begin{x}{\small\bf EXAMPLE} \ 
Take $X = \N$  $-$then \mX has the least upper bound property.
\end{x}
\vspace{0.3cm}

\begin{x}{\small\bf EXAMPLE} \ 
Take $X = \Q$ $-$then \mX does not have the least upper bound property.
\end{x}
\vspace{0.3cm}


[Assign to each rational $p > 0$ the rational
\[
q 
\ = \ 
p - \frac{p^2 - 2}{p + 2} 
\ = \ 
\frac{2p + 2}{p + 2}
\]
and note that
\[
q^2 - 2 
\ = \ 
\frac{2(p^2 - 2)}{(p + 2)^2}.
\]
Introduce
\[
\begin{cases}
\ A = \{p \in \Q \ : \ p > 0 \  \&\ p^2 < 2\}\\
\ B= \{p \in \Q \ : \ p > 0 \  \&\ p^2 > 2\}
\end{cases}
.
\]
Then
\[
\begin{cases}
\ p \in A \implies p < q \ \& \ q \in A\\
\ p \in B \implies q < p \ \& \ q \in B
\end{cases}
.
\]
Therefore
\[
\begin{cases}
\ \text{\mA has no largest element}\\
\ \text{\mB has no smalles element}
\end{cases}
.
\]
But
\[
\begin{cases}
\ \tU(A) = B\\
\ \tL(B) = A
\end{cases}
.
\]
So \mA does not have a least upper bound and \mB does not have a greatest lower bound.]
\vspace{0.3cm}

Let \mX be an ordered set.
\vspace{0.3cm}

\begin{x}{\small\bf LEMMA} \ 
Suppose that \mX has the least upper bound property.  
Let $S \subset X$
be nonempty and bounded below $-$then 
\[
\sup \tL(S) \ =  \ \inf  S.
\]
\vspace{0.2cm}

PROOF \ 
By hypothesis, $L(S) \neq \emptyset$ and 
\[
s \in S \implies s \in \tU (\tL(S)) \implies \tU (\tL(S)) \neq \emptyset.
\]
Therefore $\sup \tL(S)$ exists, call it $\lambda$.  
Given $s \in S$, there are three possibilities:
\[
s < \lambda, \quad s = \lambda, \quad \lambda < s.
\]
However $s < \lambda$ is untenable since it implies that 
\[
s \notin \tU (\tL(S)) \implies s \notin S.
\]
Accordingly
\[
s \in S \implies \lambda \leq s \implies \lambda \in \tL(S).
\]
If now $\lambda < \lambda^\prime$, then $\lambda^\prime \notin \tL(S)$ 
(for otherwise $\lambda^\prime \in \tL(S) \implies \lambda^\prime \leq \lambda$ by the very definition of $\lambda \ldots$), 
thus $\lambda = \inf S$.
\end{x}
\vspace{0.3cm}

\begin{x}{\small\bf DEFINITION} \ 
An 
\un{ordered field}
\index{ordered field} 
is an ordered set \mX which is also a field subject to the following conditions.
\vspace{0.2cm}

\qquad \textbullet \quad
If $y < z$, then $\forall \ x$, $x + y < x + z$. 
\vspace{0.2cm}

\qquad \textbullet \quad 
If $x > 0$ \& $y > 0$, then $x y > 0$.
\vspace{0.2cm}
\end{x}
\vspace{0.3cm}

\begin{x}{\small\bf EXAMPLE} \ 
Take $X = \Q$ $-$then \mX is an ordered field.
\end{x}
\vspace{0.3cm}


%% file: _02_real_numbers.tex
\chapter{
$\boldsymbol{\S}$\textbf{2}.\quad  REAL NUMBERS}
\setlength\parindent{2em}
\setcounter{theoremn}{0}
\renewcommand{\thepage}{\S2-\arabic{page}}

\ \indent 
The following result is the central theorem of existence.
\vspace{0.3cm}

\begin{x}{\small\bf THEOREM} \ 
There exists an ordered field $\R$ with the least upper bound property which contains $\Q$ as an ordered subfield.
\vspace{0.2cm}

[Note: \ 
Here there is an abuse of the language in that ``$\Q$'' is not necessarily the rationals but rather an isomorphic replica thereof.]
\end{x}
\vspace{0.3cm}

\begin{x}{\small\bf DEFINITION} \ 
The elements of $\R$ are called 
\un{real numbers}.
\index{real numbers}
\end{x}
\vspace{0.3cm}

\begin{x}{\small\bf \un{N.B.}} \ 
Suppose that $\R_1$ and $\R_2$ are two realizations of $\R$ $-$then there exists a unique order preserving field isomorphism 
$\phi:\R_1 \ra \R_2$ such that $\phi(\Q_1) = \Q_2$.
\end{x}
\vspace{0.3cm}

\begin{x}{\small\bf REMARK} \ 
There are three standard realizations of $\R$.
\vspace{0.2cm}

\qquad \textbullet \quad 
The set of infinite decimal expansions.
\vspace{0.2cm}

\qquad \textbullet \quad
The set of equivalence classes of Cauchy sequences of rational numbers.
\vspace{0.2cm}

\qquad \textbullet \quad 
The set of Dedekind cuts.

\vspace{0.2cm}

[Note: \ 
The fact that these models are actually ordered fields with the least upper bound property is not obvious, the actual verification involving a fair amount of tedious detail.]
\end{x}
\vspace{0.3cm}

\begin{x}{\small\bf REMARK} \ 
If \mS is a nonempty subset of $\R$ which is bounded below, then \mS has a greatest lower bound (cf. \S1, \#14).
\vspace{0.2cm}

[In fact, 
\[
\glb S 
\ = \ 
- \lub - S.]
\]
\end{x}
\vspace{0.3cm}

\begin{x}{\small\bf LEMMA} \ 
Let \mS be a nonempty subset of $\R$ which is bounded above $-$then
for each $\epsilon > 0$, there is an element $s \in S$ such that $s > \sup S - \epsilon$.
\vspace{0.2cm}

PROOF \ 
If there assertion were false, then for some $\epsilon > 0$ and for all $s \in S$, 
\[
\sup S - \epsilon \geq s.
\]
Accordingly, by definition of supremum,
\[
\sup S - \epsilon \geq \sup S,
\]
so $\epsilon \leq 0$, a contradiction.
\end{x}
\vspace{0.3cm}

\begin{x}{\small\bf LEMMA} \ 
Let \mS be a nonempty subset of $\R$ which is bounded above.  
Suppose that $\mu$ is an upper bound for \mS with the property that for each $\epsilon > 0$, there exists an element 
$s \in S$ such that $\mu - \epsilon < s$ $-$then $\mu = \sup S$.
\vspace{0.2cm}

PROOF \ 
If instead $\mu \neq \sup S$, then $\mu > \sup S$, hence $\mu - \sup S > 0$, thus for some $s \in S$, 
\[
\mu 
\ - \ 
(\mu \ - \ \sup S) 
\ = \ 
\sup S 
\ < \ 
s,
\]
a contradiction.
\end{x}
\vspace{0.3cm}

\begin{x}{\small\bf ARCHIMEDEAN PROPERTY} \ 
For every positive real $x$ and for every real $y$, there exists a natural number $n$ such that $n x > y$.
\vspace{0.2cm}

PROOF \ 
Suppose to the contrary that there exist real numbers $x > 0$ and $y$ such that $n x \leq y$ for every real number $n$.  
Let $S = \{n x: n \in \N\}$ $-$then \mS is bounded above (by $y$), hence has a supremum $\mu$, say.  
Because $\mu - x < \mu$ ($x$ is positive), there must be a natural number $n$ with the property that $n x > \mu - x$ 
(cf. \#6), so $(n + 1) x > \mu$.  
But $(n + 1) x$ belongs to \mS, thus the inequality $(n + 1) x > \mu$ contradicts the assumption that $\mu$ is, in particular, an upper bound for \mS.
\end{x}
\vspace{0.3cm}

\begin{x}{\small\bf COROLLARY} \ 
For every real number $x$, there exists a natural number $n$ such that $n > x$.
\end{x}
\vspace{0.3cm}

\begin{x}{\small\bf COROLLARY} \ 
For every real number $x$, there exists an integer $m$ such that $x > m$.
\vspace{0.2cm}

[Choose a natural number $n$ such that $n > -x$ (cf. \#9) $-$then $x > -n$, so we can take $m = -n$.]
\end{x}
\vspace{0.3cm}

\begin{x}{\small\bf COROLLARY} \ 
For every positive real number $x$, there exists a natural number $n$ such that $x > \ds\frac{1}{n}$.
\end{x}
\vspace{0.3cm}

\begin{x}{\small\bf EXAMPLE} \ 
Let $S = \bigg\{ \ds\frac{n}{n+1} : n \in \N \bigg\}$ $-$then $1 \in U(S)$ and we claim that $1 = \sup S$.  
Thus let $\mu = \sup S$ and suppose to the contrary that $\mu < 1$.  
Using \#11, choose a natural number $n > 1$ such that $\ds\frac{1}{n} < 1 - \mu$, hence
\[
\mu \ < \ 
1 - \frac{1}{n}  
\ = \ 
\frac{n-1}{n} ,
\]
which implies that $\mu$ is less than an element of \mS.
\end{x}
\vspace{0.3cm}

\begin{x}{\small\bf LEMMA} \ 
For every real number $x$, there exists an integer $m$ such that $x - 1 \leq m < x$.
\vspace{0.2cm}

PROOF \ 
Owing to \#9 and \#10, there exist integers $a$ and $b$ such that $a < x < b$.  
Let $m$ be the largest integer in the finite collection $a, a+1, \ldots, b$ such that $m < x$ $-$then $m + 1 \geq x$, 
hence $m \geq x - 1$.
\end{x}
\vspace{0.3cm}

\begin{x}{\small\bf DEFINITION} \ 
A nonempty subset \mS of $\R$ is said to be 
\un{dense}
\index{dense} 
in $\R$ if it has the following property: 
Between any two distinct real numbers there is an element of \mS.
\end{x}
\vspace{0.3cm}

\begin{x}{\small\bf THEOREM} \ 
$\Q$ is dense in $\R$.
\vspace{0.2cm}

PROOF \ 
Fix $x, y \in \R$ : $x < y$ $-$then $y - x > 0$, so there exists a natural number $n$ such that 
$y - x > \ds\frac{1}{n}$ (cf. \#11), i.e., such that $x < y - \ds\frac{1}{n}$.  
On the other hand, there exists an integer $m$ with the property that 
\[
n y - 1 \ \leq m \  < \  n y \qquad \text{(cf. \#13)},
\]
hence
\[
y - \frac{1}{n} \ \leq \  \frac{m}{n} \ < \ y
\]
from which 
\[
x < y \ - \  \frac{1}{n} \ \leq \ \frac{m}{n} \ < \ y.
\]
\end{x}
\vspace{0.3cm}

\begin{x}{\small\bf SCHOLIUM} \ 
If $x$ and $y$ are real numbers with $x < y$, then there exists an infinite set of rationals $q$ such that $x < q < y$.
\end{x}
\vspace{0.3cm}

The Archimedean Property is essentiallly ``additive'' in character; here is its ``multiplicative'' analog.
\vspace{0.3cm}

\begin{x}{\small\bf LEMMA} \ 
If $x > 1$ and $y$ are real numbers, then there exists a natural number $n$ such that $x^n > y$.
\vspace{0.2cm}

PROOF \ 
Proceeding by contradiction, suppose that there exist real numbers $x > 1$ and $y$ 
such that $x^n \leq y$ for every natural number $n$.  
Let $S = \{x^n : n \in \N\}$ $-$then \mS is bounded above (by $y$), hence has a supremum $\mu$, say.  
Because $x > 1$, $\mu$ is less than $\mu x$, hence $\mu / x < \mu$, so there must exsit an $n \in \N$ such that 
$\mu / x < x^n$.  
But then $\mu < x^{n+1}$ and, as $x^{n+1} \in S$, we have arrived at a contradiction.
\end{x}
\vspace{0.3cm}

\begin{x}{\small\bf EXAMPLE} \ 
Let $x > 0$ and $0 <r < 1$ be real numbers; let
\[
S = \bigg\{ \frac{x (1 - r^n)}{1 - r} : n \in \N\bigg\}.
\]
Then, in view of the relation
\[
\frac{x (1 - r^n)}{1 - r} 
\ = \ 
\frac{x}{1 - r} - \frac{x r^n }{1 - r}
\ < \ 
\frac{x}{1 - r} 
\qquad (n \in \N),
\]
it is clear that $\ds\frac{x}{1 - r}$ is an upper bound for \mS and we claim that 
\[
\frac{x}{1 - r}
\ = \ 
\sup S.
\]
To prove this, it suffices to show that if $\epsilon$ is any real number such that 
$0 < \epsilon < \ds\frac{x}{1 - r}$, then $\epsilon \notin U(S)$ (cf. \S1, \#8).  
So fix such an $\epsilon$ $-$then there exists a natural number $n$ such that 
\[
\frac{1}{r^n} 
>
\frac{x}{x - \epsilon (1 - r)}
\qquad (cf. \ \#17) \quad (0 < r < 1 \implies \frac{1}{r} > 1),
\]
thus
\[
r^n < \frac{x - \epsilon (1 - r)}{x} 
\ = \ 
1 - \epsilon \bigg(\frac{1 - r}{x}\bigg)
\]
or still, 
\[
\epsilon \ < \ \frac{x (1 - r^n)}{1 - r}
\implies 
\epsilon \notin U(S).
\]
\end{x}
\vspace{0.3cm}

\begin{x}{\small\bf DEFINITION} \ 
A real number $x$ is 
\un{irrational}
\index{irrational} 
if it is not rational.
\end{x}
\vspace{0.3cm}

\begin{x}{\small\bf NOTATION} \ 
$\PP$ is the subset of $\R$ whose elements are the irrational numbers.
\end{x}
\vspace{0.3cm}

\begin{x}{\small\bf \un{N.B.}} \ 
Therefore $\R = \PP \cup \Q$, where $\PP \cap \Q = \emptyset$.
\end{x}
\vspace{0.3cm}

\begin{x}{\small\bf LEMMA} \ 
Irrational numbers exist.
\vspace{0.2cm}

[In fact, $\R$ is not countable, hence $\PP$ is neither finite nor countable ($\Q$ being countable), hence $\PP \neq \emptyset$.]
\end{x}
\vspace{0.3cm}

\begin{x}{\small\bf THEOREM} \ 
$\PP$ is dense in $\R$.
\vspace{0.2cm}

PROOF \ 
Fix a positive irrational $p$ and fix $x, y \in \R$ : $x < y$.  
Using \#15, choose a nonzero rational $q$ such that 
\[
\frac{x}{p} \ < \ q \ < \  \frac{y}{p}.
\]
Then
\[
x \ < \ p  q \ < \ y
\]
and $p q \in \PP$.

\end{x}
\vspace{0.3cm}

\begin{x}{\small\bf \un{N.B.}} \ 
For the record, if $p \in \PP$, then $-p \in \PP$ and $\ds\frac{1}{p} \in \PP$.  
In addition, if $q \in \Q$ $(q \neq 0)$, then 
\[
p + q, \ p - q, \ p q, \ \frac{p}{q}
\]
are irrational.
\end{x}
\vspace{0.3cm}

\begin{x}{\small\bf DEFINITION} \ 
An element $x \in \R$ is 
\un{algebraic}
\index{algebraic} 
or 
\un{transcendental}
\index{transcendental} 
according to whether it is or is not a root of a nonzero polynomial in $\Z[X]$.
\end{x}
\vspace{0.3cm}

\begin{x}{\small\bf EXAMPLE} \ 
If $\ds\frac{a}{b}$ $(b \neq 0)$ is rational, then $\ds\frac{a}{b}$ is algebraic.
\vspace{0.2cm}

[Consider the polynomial $bX - a$.]
\end{x}
\vspace{0.3cm}

\begin{x}{\small\bf EXAMPLE} \ 
Let $r, s \in \Q$, $r > 0$ $-$then $r^s$ is algebraic.
\vspace{0.2cm}

[Write $s = \ds\frac{m}{n}$ $(m, n \in \Z, n > 0)$ and consider the polynomial $X^n - r^m$.]
\vspace{0.2cm}


[Note: \ 
Take $r = 2$, $s = \ds\frac{1}{2}$, hence $n = 2$ and 
$2^{^{\frac{1}{2}}} = \sqrt{2}$ 
is algebraic (but irrational 
(cf. \S6, \#2)).]
\end{x}
\vspace{0.3cm}

\begin{x}{\small\bf \un{N.B.}} \ 
It will be shown in due course that $e$ and $\pi$ are transcendental.  
However the status of $e + \pi$, $e - \pi$, $e \pi$, $e^e$, and $\pi^\pi$ is unknown.
\vspace{0.2cm}

[Note: \ 
$e^\pi$ is transcendental but whether this is true of $\pi^e$ remains an open question.]
\end{x}
\vspace{0.3cm}

\begin{x}{\small\bf EXAMPLE} \ 
Is $e + \pi$ irrational?  
Is $e \pi$ irrational?  
Answer: \ 
Nobody knows.  
But at least one of them must be irrational.  
To see this, consider the polynomial 
\[
X^2 - (e + \pi) X + e \pi.
\]
Its zeros are $e$ and $\pi$.  
So if both $e + \pi$ and $e \pi$ were rational, then $e$ and $\pi$ would be algebraic which they are not.
\end{x}
\vspace{0.3cm}

\begin{x}{\small\bf NOTATION} \ 
$\ov{\Q}$ is the subset of $\R$ whose elements are the algebraic numbers 
and $\T$ is the subset of $\R$ whose elements are the transcendental numbers.
\end{x}
\vspace{0.3cm}

\begin{x}{\small\bf \un{N.B.}} \ 
$\Q$ is a subset of $\ov{\Q}$ and $\T$ is a subset of $\PP$.
\end{x}
\vspace{0.3cm}

\begin{x}{\small\bf LEMMA} \ 
The cardinality of $\ov{\Q}$ is $\aleph_0$.
\end{x}
\vspace{0.3cm}

\begin{x}{\small\bf \un{N.B.}} \ 
Consequently, on purely abstract grounds, transcendental numbers exist.  
Historically, the first explicit transcendental number was constructed by Liouville, viz. 
\[
\sum\limits_{n = 1}^\infty \hsx
10^{-n!} 
\qquad \text{(cf. \S15, \#9)}.
\]
\end{x}
\vspace{0.3cm}
\begin{x}{\small\bf LEMMA} \ 
$\ov{\Q}$ is the algebraic closure of $\Q$ in $\R$ and 
\[
[\ov{\Q} : \Q] 
\ = \ 
\aleph_0.
\]
\end{x}
\vspace{0.3cm}

Being a field, $\ov{\Q}$ is closed under addition and multiplication.
\vspace{0.3cm}

\begin{x}{\small\bf LEMMA} \ 
If $x \neq 0$ is algebraic and $y$ is transcendental, then $x + y$ and $x y$ are transcendental.
\end{x}
\vspace{0.3cm}

\begin{x}{\small\bf EXAMPLE} \ 
$\sqrt{2} \hsx e$ and $\sqrt{2} + \pi$ are transcendental.
\end{x}
\vspace{0.3cm}

\begin{x}{\small\bf LEMMA} \ 
If $x \in \R$ is transcendental, then so is $x^2$.
\vspace{0.2cm}

[If $x^2$ were algebraic, then there would be a relation of the form
\[\
a_0 + a_2 x^2 + a_4 x^4 + \cdots + a_{2n} x^{2n} 
\ = \ 
0
\qquad (a_{2 k} \in \Q)
\]
or still, 
\[
a_0 + 0 x + a_2 x^2 + 0x^3 + a_4 x^4 + \cdots + a_{2n} x^{2n}
\ = \ 
0
\]
implying thereby that $x$ is algebraic.]
\end{x}
\vspace{0.3cm}

\begin{x}{\small\bf EXAMPLE} \ 
Not both $e \pi$ and $\ds\frac{\pi}{e}$ can be algebraic.
\vspace{0.2cm}

[In fact, 
\[
(e \pi) \bigl( \frac{\pi}{e}\bigr) = \pi^2.]
\]

\end{x}
\vspace{0.3cm}

\begin{x}{\small\bf \un{N.B.}} \ 
$\T$ is not closed under addition and multiplication.
\end{x}
\vspace{0.3cm}

\begin{x}{\small\bf CRITERION} \ 
Let $x$ and $y$ be real numbers.  
Suppose that $x \leq y + \epsilon$ for every $\epsilon > 0$ $-$then $x \leq y$.

\vspace{0.2cm}

PROOF \ 
Assume that $x > y$ and put $\epsilon = \ds\frac{1}{2} (x - y)$ $-$then $\epsilon > 0$.  
However
\[
y + \epsilon 
\ = \ 
\frac{1}{2} \hsx (x + y) 
\ < \ 
\frac{1}{2} \hsx (x + x) 
\ = \ 
x,
\]
contrary to the supposition that $y + \epsilon \geq x$ for every $\epsilon > 0$.
\end{x}
\vspace{0.3cm}

%% file: _03_suprema.tex
\chapter{
$\boldsymbol{\S}$\textbf{3}.\quad  SUPREMA}
\setlength\parindent{2em}
\setcounter{theoremn}{0}
\renewcommand{\thepage}{\S3-\arabic{page}}

\ \indent 
We shall record here some technicalities that will be of use in the sequel.
\vspace{0.3cm}

\begin{x}{\small\bf LEMMA} \ 
Let \mS be a nonempty subset of $\R$, \mT a nonempty subset of \mS.  
Suppose that \mS is bounded above $-$then \mT is also bounded above and $\sup T \leq \sup S$.
\vspace{0.2cm}

[This is obvious from the definitions.]
\end{x}
\vspace{0.3cm}

\begin{x}{\small\bf LEMMA}\ 
Let \mS and \mT be two nonempty subsets of $\R$, each being bounded above.  
Suppose further that given any $s \in S$ there is a $t \in T$ such that $s \leq t$ and that given any $t \in T$ there is an $s \in S$ such that 
$t \leq s$ $-$then $\sup S = \sup T$.
\vspace{0.2cm}

PROOF \ 
It suffices to rule out the other possibilities:
\[
\begin{cases}
\ \sup S < \sup T\\
\ \sup T < \sup S
\end{cases}
.
\]
If the first of these were true, then $\sup S \notin \tU(T)$, so there exists a $t \in T$ such that $\sup S < t \leq \sup T$.  
But, by hypothesis, there is an $s \in S$ such that $t \leq s$, hence $\sup S < s$, a contradiction.  
The second of these can be eliminated in the same way.
\end{x}
\vspace{0.3cm}

\begin{x}{\small\bf NOTATION} \ 
Given nonempty subsets \mS, \mT of $\R$, put
\[
S + T \ = \ \{s + t : s \in S, \ t \in T\}.
\]
\end{x}
\vspace{0.3cm}

\begin{x}{\small\bf LEMMA} \ 
Let \mS and \mT be nonempty subsets of $\R$, each being bounded above $-$then $S + T$ is bounded above and 
\[
\sup (S + T) \ = \ \sup S + \sup T.
\]
\vspace{0.2cm}

PROOF \ 
Let $r \in S + T$ $-$then there exist $s \in S$, $t \in T$ such that $r = s + t$ and so $r \leq \sup S + \sup T$.  
Since $r$ is an arbitrary element of $S + T$, it follows 
that $\sup S + \sup T$ is an upper bound for $S + T$, hence $\sup (S + T)$ exists and in fact
\[
\sup (S + T) \ \leq \  \sup S + \sup T.
\]
To reverse this, we shall employ \S2, \#40 and prove that 
\[
\sup S + \sup T \ \leq \ \sup (S + T) + \epsilon
\]
for every $\epsilon > 0$.  Thus fix $\epsilon > 0$ and choose $s \in S$, $t \in T$ such that 
\[
s \ > \ \sup S - \frac{\epsilon}{2}, \quad t \ > \ \sup T - \frac{\epsilon}{2} \qquad \text{(cf. \S2, \#6)}.
\]
Then
\[
s + t \ > \sup S + \sup T  - \epsilon
\]
or still, 
\allowdisplaybreaks
\begin{align*}
\sup S + \sup T \ 
&<\ \ s + t + \epsilon
\\[12pt]
&\leq 
\sup (S + T) + \epsilon.
\end{align*}
\end{x}
\vspace{0.3cm}

\begin{x}{\small\bf NOTATION} \ 
Given nonempty subsets \mS, \mT of $\R$, put
\[
S \cdot T \ = \ \{s t : s \in S, \ t \in T\}.
\]
\end{x}
\vspace{0.3cm}

\begin{x}{\small\bf LEMMA} \ 
Let \mS and \mT be nonempty subsets of $\R_{> 0}$, each being bounded above $-$then $S \cdot T$ is bounded above and 
\[
\sup (S \cdot T) \ = \ (\sup S) \cdot (\sup T).
\]
\vspace{0.2cm}

PROOF \ 
Note first that
\[
\sup S \ > \ 0 \quad \text{and} \quad \sup T \ > \ 0.
\]
This said, let $r \in S \cdot T$ $-$then there exist $s \in S$, $t \in T$ such that $r = s t$ and so $r \leq (\sup S) \cdot (\sup T)$.  
Since $r$ is an arbitrary element of $S \cdot T$, it follows that $(\sup S) \cdot (\sup T)$ is an upper bound for $S \cdot T$, hence $\sup (S \cdot T)$ exists and in fact
\[
\sup (S \cdot T) \ \leq \ (\sup S) \cdot (\sup T).
\]
To reverse this, we shall employ \S2, \#40 and prove that
\[
(\sup S) \cdot (\sup T)
\ \leq \ 
\sup (S \cdot T) + \epsilon
\]
for every $\epsilon > 0$.  Thus fix $\epsilon > 0$ and choose $s \in S$, $t \in T$ such that 
\[
s > \sup S - \frac{\epsilon}{\sup S + \sup T}, \quad
t > \sup T - \frac{\epsilon}{\sup S + \sup T} 
\qquad \text{(cf. \S2, \#6)}.
\]
Then 
\[
\sup S - s < \frac{\epsilon}{\sup S + \sup T}, \quad
\sup T - t < \frac{\epsilon}{\sup S + \sup T},
\]
from which 
\[
t ( \sup S - s) \ \leq \  \frac{\epsilon \cdot \sup T }{\sup S + \sup T} 
\]
and 
\[
\sup S ( \sup T - t) \ <  \ \frac{\epsilon \cdot \sup S }{\sup S + \sup T}.
\]
Therefore
\allowdisplaybreaks
\begin{align*}
(\sup S) \cdot (\sup T) - st \ 
&=\ 
\sup S (\sup T - t) + t( \sup S - s)
\\[12pt]
&<\ 
\frac{\epsilon \cdot \sup S }{\sup S + \sup T} + \frac{\epsilon \cdot \sup T }{\sup S + \sup T}
\\[12pt]
&=\ 
\epsilon, 
\end{align*}
i.e., 
\allowdisplaybreaks
\begin{align*}
(\sup S) \cdot (\sup T) \ 
&\leq \ 
st + \epsilon
\\[12pt]
&\leq \
\sup (S \cdot T) + \epsilon.
\end{align*}
\end{x}
\vspace{0.3cm}

\begin{x}{\small\bf REMARK} \ 
The assertion of \#6 may be false if we drop the assumption 
that \mS and \mT are nonempty subsets of $\R_{> 0}$.
\vspace{0.2cm}

[Take, e.g., $S = - \N$,  $T = - \N$, $-$then both \mS and \mT are bounded above but $S \cdot T$ is not.]
\end{x}
\vspace{0.3cm}


%% file: _04_exponents_and_roots.tex
\chapter{
$\boldsymbol{\S}$\textbf{4}.\quad  EXPONENTS AND ROOTS}
\setlength\parindent{2em}
\setcounter{theoremn}{0}
\renewcommand{\thepage}{\S4-\arabic{page}}

\ \indent 
Let $a > 0$ and $x$ be real numbers$-$then the primary objective of the present \S \ is to assign a meaning to the symbol $a^x$.
\vspace{0.3cm}

If $a$ is any real number and if $n$ is a natural number, then the power $a^n$ is defined inductively by the rule
\[
a^1 = a, \ \ a^{n+1} = a^n \cdot a.
\]
When $a \neq 0$, we define $a^0$ as 1; we do not define $0^0$.  
When $a \neq 0$, we define $a^{-n}$ as $\ds\frac{1}{a^n}$; we do not define $0^{-n}$.

\vspace{0.3cm}

\begin{x}{\small\bf LAWS OF EXPONENTS FOR INTEGRAL POWERS} \ 
Let $a$ and $b$ be nonzero real numbers; let $m$ and $n$ be integers.
\vspace{0.2cm}

\qquad (1)  \ 
$a^m \cdot a^n = a^{m + n}$;
\vspace{0.2cm}

\qquad (2)  \ 
$\bigl(a^m\bigr)^n = a^{m n}$;
\vspace{0.2cm}

\qquad (3)  \ 
$\ds\frac{a^m}{a^n} = a^{m - n}$; 
\vspace{0.2cm}

\qquad (4)  \ 
$(a b)^m = a^m b^m$; 
\vspace{0.2cm}

\qquad (5)  \ 
$\bigg(\ds\frac{a}{b}\bigg)^m = \ds\frac{a^m}{b^m}$;
\vspace{0.2cm}

\qquad (6) (i) \ 
If $n > 0$ and $a, \ b > 0$, then $a < b$ if and only if $a^n < b^n$.
\vspace{0.2cm}

\qquad (6)  (ii) \ 
If $n < 0$ and $a, \ b > 0$, then $a < b$ if and only if $a^n > b^n$.
\vspace{0.2cm}

\qquad (7) (i)   \ 
If $a > 1$, then $m < n$ if and only if $a^m < a^n$.
\vspace{0.2cm}

\qquad (7)  (ii)  \ 
If $0 < a < 1$, then $m < n$ if and only if $a^m > a^n$.

\end{x}
\vspace{0.3cm}

In order to define the symbol $a^r$ for rational $r$, it is first necessary to establish the existence and uniqueness of 
``$n^\nth$ roots''.
\vspace{0.3cm}

\begin{x}{\small\bf THEOREM} \ 
For every real $a > 0$ and every natural number $n$, there is one and only one real $x > 0$ such that $x^n = a$.
\vspace{0.2cm}

Uniqueness is immediate.  
For suppose that $x_1 > 0$, $x_2 > 0$ are such that $x_1^n = a$, $x_2^n = a$ $-$then these conditions imply that $x_1 = x_2$ 
(cf. \#1, 6(i)).
\vspace{0.2cm}

Turning to existence, let \mS be the set of all positive real numbers $s$ such that $s^n < a$.
\end{x}
\vspace{0.3cm}

\begin{x}{\small\bf LEMMA} \ 
\vspace{0.2cm}
\mS is nonempty and is bounded above.
\vspace{0.2cm}

PROOF \ 
To see that \mS is nonempty, observe that $\ds\frac{a}{1 +a}$ lies between 0 and 1, hence
\[
\frac{a^n}{(1 + a)^n} 
\ \leq \ 
\frac{a}{1 + a} 
\ < \ 
a
\implies 
\frac{a}{1 + a} \in S.
\]
In addition, $1 + a \in \tU(S)$.  
Indeed, if there exists $s \in S$ such that $s > 1 + a$ $(> 1)$, then $s^n > s > 1 + a > a$, a contradiction.

Let $\mu = \sup S$ $-$then we claim that $\mu^n = a$.  
To establish this, it suffices to eliminate the other possibilities: 
\[
\begin{cases}
\ \mu^n < a\\
\ \mu^n > a
\end{cases}
.
\]
\vspace{0.2cm}

\qquad $\un{\mu^n < a:}$ \ 
Since
\[
\frac{a - \mu^n}{(1 + \mu)^n - \mu^n}
\]
is a positive real number, one can choose a real number $\nu$ lying between 0 and 1 and such that
\[
\nu 
\ < \ 
\frac{a - \mu^n}{(1 + \mu)^n - \mu^n} 
\qquad \text{(e.g. quote \S2, \#15)}.
\]
Then
\allowdisplaybreaks
\begin{align*}
(\mu + \nu)^n \ 
&=\ 
\mu^n + \binom{n}{1} \mu^{n-1} \nu + \binom{n}{2} \mu^{n-2} \nu^2 + \cdots + \binom{n}{n} \nu^n
\\[12pt]
&\leq\ 
\mu^n + \nu\bigg[\binom{n}{1} \mu^{n-1} + \binom{n}{2} \mu^{n-2} + \cdots + \binom{n}{n}\bigg]
\\[12pt]
&=\ 
\mu^n  + \nu [(1 + \mu)^n - \mu^n] 
\\[12pt]
&<\ 
\mu^n + (a - \mu^n) 
\\[12pt]
&=\ 
a.
\end{align*}
Therefore $\mu + \nu \in S$, which contradicts the fact that $\mu$ is an upper bound for \mS.
\vspace{0.3cm}

\qquad $\un{\mu^n > a:}$ \ 
Choose a real number $\nu$ lying between 0 and 1 with the following properties:
\[
\nu < \mu 
\quad \text{and} \quad 
\nu < \frac{\mu^n - a}{(1 + \mu)^n - \mu^n}.
\]
Then for $s > \mu - \nu$, we have
\allowdisplaybreaks
\begin{align*}
s^n \ 
&\geq 
(\mu - \nu)^n \ 
\\[12pt]
&=\ 
\mu^n -\binom{n}{1}\mu^{n-1}\nu + \binom{n}{2} \mu^{n-2}\nu^2 - \cdots + (-1)^n \binom{n}{n}\nu^n
\\[12pt]
&=\ 
\mu^n - \nu\bigg[\binom{n}{1}\mu^{n-1} - \binom{n}{2} \mu^{n-2}\nu + \cdots - (-1)^n\binom{n}{n}\nu^{n-1}\bigg]
\\[12pt]
&\geq\ 
\mu^n - \nu\bigg[\binom{n}{1}\mu^{n-1} + \binom{n}{2} \mu^{n-2} + \cdots + \binom{n}{n}\bigg]
\\[12pt]
&=\ 
\mu^n - \nu[(1 + \mu)^n - \mu^n]
\\[12pt]
&>\ 
\mu^n  - (\mu^n- a)
\\[12pt]
&=\ 
a.
\end{align*}
Therefore $\mu - \nu$ is an upper bound for \mS, which contradicts the fact that $\mu$ is the supremum for \mS.  

Consequently
\[
\mu^n \ = \ a,
\]
as claimed.

Let $a > 0$ be a positive real number $-$then for each natural number $n$, the preceding theorem guarantees the existence and uniqueness of a real number $x > 0$ such that $x^n = a$.  
We write $\ds\sqrt[\leftroot{-1}\uproot{3}n]{a}$ for this $x$ and call $\ds\sqrt[\leftroot{-1}\uproot{3}n]{a}$ the 
\un{$n^\nth$ root}
\index{$n^\nth$ root} 
of a.
\vspace{0.2cm}

[Note: \ 
If $n = 1$, write a for $\ds\sqrt[\leftroot{-1}\uproot{3}1]{a}$; 
if $n = 2$, write $\ds\sqrt{a}$ for $\ds\sqrt[\leftroot{-1}\uproot{3}2]{a}$.]
\end{x}
\vspace{0.3cm}

\begin{x}{\small\bf EXAMPLE} \ 
$\sqrt{2}$ exists.
\end{x}
\vspace{0.3cm}

Suppose now that $a < 0$ is a negative real number $-$then for each odd natural number $n$, $\sqrt[\leftroot{-1}\uproot{3}n]{a}$ is taken to be the unique real $x < 0$ such that 
$-x = \sqrt[\leftroot{-1}\uproot{3}n]{-a}$ (e.g., $\sqrt[\leftroot{-1}\uproot{3}3]{-8} = -2)$.  
Since $n$ is odd, 
\[
x^n \ = \ (-(-x))^n \ = \  (-1)^n (-x)^n \ = \ - (-a) \ = \ a,
\]
thereby justifying the definition.
\vspace{0.2cm}

[Note: \ 
We do not define $\sqrt[\leftroot{-1}\uproot{3}n]{a}$ when $a < 0$ and $n$ is an even natural number.]
\vspace{0.3cm}

\begin{x}{\small\bf \un{N.B.}} \ 
Set $\sqrt[\leftroot{-1}\uproot{3}n]{0} = 0$ for all $n \in \N$.
\end{x}
\vspace{0.3cm}

Let $a > 0$ be a positive real number.  
Given a rational number $r$, let $\ds\frac{m}{n}$ be the representation of $r$ in lowest terms.
\vspace{0.3cm}

\begin{x}{\small\bf DEFINITION} \ 
\[
a^r
\ = \ 
\bigl(\sqrt[\leftroot{-1}\uproot{3}n]{a}\bigr)^m,
\]
the $m^\nth$ power of the $n^\nth$ root of $a$ (if $m = 1$, then $\ds a^{^{\frac{1}{n}}} = \sqrt[\leftroot{-1}\uproot{3}n]{a}$).
\vspace{0.2cm}

[Note: \ 
Regardless of the sign of $m$, it is clear that $a^r > 0$.]
\end{x}
\vspace{0.3cm}

\begin{x}{\small\bf LAWS OF EXPONENTS FOR RATIONAL POWERS} \ 
Let $a$ and $b$ be positive real numbers; let $r$ and $s$ be rational numbers.
\vspace{0.2cm}

\qquad (1)  \ 
$a^r \cdot a^s = a^{r + s}$;
\vspace{0.2cm}

\qquad (2)  \ 
$\bigl(a^r\bigr)^s = a^{r s}$;
\vspace{0.2cm}

\qquad (3)  \ 
$\ds\frac{a^r}{a^s} = a^{r - s}$; 
\vspace{0.2cm}

\qquad (4)  \ 
$(a b)^r = a^r b^r$; 
\vspace{0.2cm}

\qquad (5)  \ 
$\bigg(\ds\frac{a}{b}\bigg)^r = \ds\frac{a^r}{b^r}$;
\vspace{0.2cm}

\qquad (6) (i) \ 
If $r > 0$, then $a < b$ if and only if $a^r < b^r$.
\vspace{0.2cm}

\qquad (6)  (ii) \ 
If $r < 0$, then $a < b$ if and only if $a^r > b^r$.
\vspace{0.2cm}

\qquad (7) (i)   \ 
If $a > 1$, then $r < s$ if and only if $a^r < a^s$.
\vspace{0.2cm}

\qquad (7)  (ii)  \ 
If $0 < a < 1$, then $r < s$ if and only if $a^r > a^s$.
\end{x}
\vspace{0.3cm}

\begin{x}{\small\bf REMARK} \ 
If $p$ is a natural number, then
\[
(\sqrt[\leftroot{-1}\uproot{3}n]{a})^m 
\ = \ 
(\sqrt[\leftroot{-1}\uproot{3}np]{a})^{mp}.
\]
Therefore in the definition of the symbol $a^r$, it is not necessary to require that $r$ be reduced to lowest terms so, for example,
\[
a \ = \ a^1 \ = \ (\sqrt[\leftroot{-1}\uproot{3}n]{a})^n \qquad (n \in \N).
\]
\end{x}
\vspace{0.3cm}

\begin{x}{\small\bf LEMMA} \ 
Let $a > 0$, $a \neq 1$ $-$then
\[
\frac{a^r - 1}{r} 
\ < \ 
\frac{a^s - 1}{s}
\]
for all $r$, $s \in \Q - \{0\}$ with $r < s$.
\vspace{0.2cm}

PROOF \ 
Let us admit for the moment that the lemma is true when, in addition, $r$ and $s$ are nonzero integers with $r < s$.  
Proceeding to the general case, there is no loss of generality in supposing that $r = p/n$, $s = q/n$, where $n \in \N$, $p$ and 
$q \in \Z - \{0\}$, and $p < q$.  
It is then a question of proving that
\[
\frac{(a^{p/n} - 1)n}{p} 
\ < \ 
\frac{(a^{q/n} - 1)n}{q},
\]
or, equivalently, since $n > 0$, that
\[
\frac{a^{p/n} - 1}{p} 
\ < \ 
\frac{a^{q/n} - 1}{q}.
\]
Put $b = \sqrt[\leftroot{-1}\uproot{3}n]{a}$ $-$then, since we are granting temporarily the truth of the lemma in the integral case, it follows that
\[
\frac{b^p - 1}{p} 
\ < \ 
\frac{b^q - 1}{q},
\]
as desired.  
Turning now to the case when $r$ and $s$ are nonzero integers with $r < s$, it is enough to consider just three possibilities, namely  
(i) \ $0 < r < r+1 = s$; 
(ii) \ $r < r + 1 = s < 0$; 
(iii) \ $-1 = r < s = 1$.  
The first of these is the assertion that 
\[
\frac{a^r - 1}{r} 
\ < \ 
\frac{a^{r+1} - 1}{r+1}
\]
or still, upon multiplying both sides of the inequality by $r(r+1)$, that
\[
(r + 1)a^r - 1
\ < \ 
r a^{r+1},
\]
or still, that
\[
a^r - 1 
\ < \ 
r a^r (a - 1),
\]
or still, upon division by $a - 1 \neq 0$, that
\[
\begin{cases}
\ a^{r-1} + a^{r-2} + \cdots + a + 1 < r a^r \quad \text{if} \quad a > 1\\
\ a^{r-1} + a^{r-2} + \cdots + a + 1 > r a^r \quad \text{if} \quad 0 < a < 1
\end{cases}
.
\]
But these inequalities do in fact obtain (apply \#1, 7(i) and 7(ii)).  
The second case, $r < r + 1 = s < 0$, can be reduced to the first case by considering $-s$, $-r$, and $a^{-1}$.  
Finally, if $r = -1$ and $s = 1$, then the inequality to be established can be written $1 - a^{-1} < a - 1$ and this is certainly true for 
$a > 0$, $a \neq 1$.
\end{x}
\vspace{0.3cm}

Fix a real number $a > 1$.  Given a rational number $x$, let
\[
S \ = \ \{a^r: r \in \Q \quad \text{and} \quad r < x\}.
\]
\vspace{0.3cm}

\begin{x}{\small\bf SUBLEMMA} \ 
\mS is nonempty and has an upper bound \mM, say, thus \mS has a supremum.
\end{x}
\vspace{0.3cm}

\begin{x}{\small\bf LEMMA} \ 
$\sup S = a^x$.
\vspace{0.2cm}

PROOF \ 
Since $a^x \in \tU(S)$, is suffices to show that for each $\epsilon > 0$, there is a rational number $r < x$ such that $a^x - a^r < \epsilon$ (cf. \S2, \#7).  
Without yet committing ourselves, it can be assumed from the beginning that $0 < x - r < 1$, hence
\[
\frac{a^{x-r} - 1}{x - r} 
\ < \ 
a - 1 
\ < \ 
a + 1 
\qquad \text{(cf. \#9)},
\]
from which
\allowdisplaybreaks
\begin{align*}
a^x - a^r \ 
&=\ 
a^r \bigg[\frac{a^{x-r} - 1}{x - r}\bigg] (x - r)
\\[12pt]
&< \ M(a+1)(x-r),
\end{align*}
so if $r < x$ is chosen in such a way that
\[
0 \ <  \ x - r \ < \ \frac{1}{2} \min\bigg\{\frac{\epsilon}{M(a+1)},1\bigg\},
\]
then $a^x - a^r < \epsilon$.
\end{x}
\vspace{0.3cm}

Fix a real number $a > 1$.  Given a real number $x$, let 
\[
S \ = \ \{a^r : r \in \Q \quad \text{and} \quad r < x\}.
\]
\vspace{0.3cm}

\begin{x}{\small\bf SUBLEMMA} \ 
\mS is nonempty and bounded above.
\vspace{0.2cm}

[It is clear that \mS is nonempty (cf. \S2, \#10).  
On the other hand, if $n$ is any natural number $> x$ (cf. \S2, \#9), then 
\allowdisplaybreaks
\begin{align*}
r < x 
&\implies
r < n
\\[12pt]
&\implies
a^r  < a^n \qquad \text{(cf. \#7, 7(i))}
\\[12pt]
&\implies
a^n \in \tU(S) 
\\[12pt]
&\implies
\tU(S) \neq \emptyset.]
\end{align*}
\end{x}
\vspace{0.3cm}

\begin{x}{\small\bf DEFINITION} \ 
$a^x = \sup S$.
\vspace{0.2cm}

[Note: \ 
If $a = 1$, we define $a^x$ as 1.  
If $0 < a < 1$, then $1/a > 1$ and we define $a^x$ as $1 / (1/a)^x$.  
In all cases: $a^x > 0$.]
\end{x}
\vspace{0.3cm}

\begin{x}{\small\bf \un{N.B.}} \ 
Matters are consistent when restricted to rational $x$ (cf. \#11).
\end{x}
\vspace{0.3cm}

\begin{x}{\small\bf LAWS OF EXPONENTS FOR REAL POWERS} \ 
Let $a$ and $b$ be positive real numbers; let $x$ and $y$ be real numbers.
\vspace{0.2cm}

\qquad (1)  \ 
$a^x \cdot a^y = a^{x + y}$;
\vspace{0.2cm}

\qquad (2)  \ 
$\bigl(a^x\bigr)^y = a^{x y}$;
\vspace{0.2cm}

\qquad (3)  \ 
$\ds\frac{a^x}{a^y} = a^{x - y}$; 
\vspace{0.2cm}

\qquad (4)  \ 
$(a b)^x = a^x b^x$; 
\vspace{0.2cm}

\qquad (5)  \ 
$\bigg(\ds\frac{a}{b}\bigg)^x = \ds\frac{a^x}{b^x}$;
\vspace{0.2cm}

\qquad (6) (i) \ 
If $x > 0$, then $a < b$ if and only if $a^x < b^x$.
\vspace{0.2cm}

\qquad (6)  (ii) \ 
If $x < 0$, then $a < b$ if and only if $a^x> b^x$.
\vspace{0.2cm}

\qquad (7) (i)   \ 
If $a > 1$, then $x < y$ if and only if $a^x < a^y$.
\vspace{0.2cm}

\qquad (7)  (ii)  \ 
If $0 < a < 1$, then $x < y$ if and only if $a^x > a^y$.
\end{x}
\vspace{0.3cm}

The proof of this result is spelled out in the lines below.
\vspace{0.2cm}

[Note: \ 
We shall omit consideration of trivial, special cases (e.g., $1^x \cdot 1^y = 1^{x+y}$ etc.] 
\\[15pt]
\un{LAW 1:} \quad

\qquad \un{Case 1:} \ 
$a > 1$.  \ Let

\allowdisplaybreaks
\begin{align*}
S \ 
&=\ 
\{a^s : s \in \Q \ \text{and} \ s < x\}
\\[12pt]
T \ 
&=\ 
\{a^t : t \in \Q \ \text{and} \ t < y\}
\\[12pt]
U \ 
&=\ 
\{a^u : u \in \Q \ \text{and} \ u < x + y\},
\\[12pt]
\end{align*}
thus $a^x = \sup S$, $a^y = \sup T$, $a ^{x+y} = \sup U$.  
In addition,
\allowdisplaybreaks
\begin{align*}
a^x \cdot a^y \ 
&=\ 
(\sup S) \cdot (\sup T)
\\[12pt]
&=\ 
\sup (S \cdot T) \qquad \text{(cf. \S3, \#6)},
\end{align*}
and
\allowdisplaybreaks
\begin{align*}
S \cdot T \ 
&=\ 
\{a^s \cdot a^t : s, \ t\in \Q \ \text{and} \ s < x, t < y\}
\\[12pt]
&=\ 
\{a^{s + t}: s, \ t\in \Q \ \text{and} \ s < x, t < y\}.
\end{align*}
So, to prove that $a^x \cdot a^y = a^{x + y}$, it will be enough to prove that $\sup (S \cdot T) = \sup U$ and for this purpose, we shall employ \S3, \#2.  
Since $S \cdot T$ is a subset of \mU, it need only be shown that given any element $a^u$ ($u \in \Q$ and $u < x + y$) in \mU, there exist rational numbers $s$, $t$ with $s < x$, $t < y$ and such that $u < s + t$ (for then $a^u < a^{s +t} \in S \cdot T$).  
Noting that 
\[
\frac{u - x + y}{2}
\ < \ 
y, 
\quad 
\frac{u - y + x}{2}
\ < \ 
x, 
\]
choose rational numbers $s$ and $t$ such that
\[
\frac{u - y + x}{2} \ <  \ s \ < \ x, 
\quad 
\frac{u - x + y }{2} \ <  \ t \ < \ y
\qquad \text{(cf. \S2, \#15)}.
\]
Then
\[
u 
\ = \ 
\frac{u - y + x}{2} + \frac{u - x + y }{2}
\ < \ 
s + t.
\]

\qquad \un{Case 2:} \ 
$0 < a < 1$. \ 
We have
\allowdisplaybreaks
\begin{align*}
a^x \cdot a^y \ 
&=\ 
\frac{1}{(1/a)^x} \cdot \frac{1}{(1/a)^y} 
\\[12pt]
&=\ 
\frac{1}{(1/a)^x \cdot (1/a)^y}
\\[12pt]
&=\ 
\frac{1}{(1/a)^{x+y}}
\\[12pt]
&=\ 
a^{x+y}.
\end{align*}
\vspace{0.3cm}

A simple but importan consequence of LAW 1 is the fact that
\[
a^x \ = \ \frac{1}{a^{-x}}
\qquad (a > 0, \ x \in \R).
\]
Proof:
\[
1 \ = \ a^0 \ = \ a^{x-x} \ = \ a^x \cdot a^{-x} 
\implies
a^x \ = \ \frac{1}{a^{-x}}.
\]
\\[15pt]
\qquad \un{LAW 2:} \quad
\vspace{0.2cm}

\qquad \un{Case 1:} \ 
$y \in \Z$. \ 
Suppose first that $y \in \N$ and argue by induction.  
The assertion is trivial if $y = 1$.  
Assuming that the assertion is true for $y = n$, we have
\allowdisplaybreaks
\begin{align*}
(a^x)^{n+1} \ 
&=\ 
(a^x)^n \cdot a^x 
\qquad \text{(by definition)}
\\[12pt]
&=\ 
(a^{xn}) \cdot a^x 
\qquad \text{(by induction hypothesis)}
\\[12pt]
&=\ 
(a^x)^{n+1} 
\hspace{1.0cm} \text{(by LAW 1)}.
\\[12pt]
\end{align*}
It therefore follows that $(a^x)^y = a^{xy}$ for arbitrary $a > 0$, $x$ real, and $y$ a positive integer.  
The assertion is trivial if $y = 0$ and the reader can supply the details if $y$ is a negative integer.
\vspace{0.2cm}

\qquad \un{Case 2:} \ 
$y \in \Q$. \ 
Let $\ds\frac{m}{n}$ be the representation of $y$ in lowest terms.  By Case 1, $(a^x)^m = a^{xm}$.  Therefore
\allowdisplaybreaks
\begin{align*}
\bigg(a^x\bigg)^\frac{m}{n} \ 
&=\ 
\bigg(\big(a^x\big)^m\bigg)^\frac{1}{n}
\\[12pt]
&=\ 
\big(a^{xm}\big)^\frac{1}{n}
\\[12pt]
&=\ 
\big(a^{x \frac{m}{n} \cdot n}\big)^\frac{1}{n}
\\[12pt]
&=\ 
\bigg(\big(a^{x \frac{m}{n}}\big)^n\bigg)^\frac{1}{n} 
\qquad \text{(by Case 1)}
\\[12pt]
&=\ 
a^{x \frac{m}{n}}.
\end{align*}

\vspace{0.2cm}

\qquad \un{Case 3:} \ 
$a > 1$, $x > 0$, $y$ arbitrary.  \ 
Let
\allowdisplaybreaks
\begin{align*}
S \ 
&=\ 
\{(a^x)^s : s \in \Q \ \text{and} \ s < y\}
\\[12pt]
T \ 
&=\ 
\{a^t : t \in \Q \ \text{and} \ t < xy\},
\end{align*}
thus $(a^x)^y = \sup S$, $a^{xy} = \sup T$, the claim being that $\sup S = \sup T$.  
To this end, we shall utilize \S3, \#2.  
In view of Case 2, 
\[
S 
\ = \ 
\{a^{xs} : s \in \Q \quad \text{and} \quad s < y\}.
\]
Given $a^{xs} \in S$, choose a rational number $t$ such that $xs < t < xy$ $-$then $a^{xs} < a^t$ and $a^t \in T$.  
On the other hand, given $a^t \in T$, choose a rational number $s$ such that 
$\ds\frac{t}{x} < s < y$ $-$then $a^t < a^{xs}$ and $a^{xs} \in S$.

\vspace{0.2cm}

\qquad \un{Case 4:} \ 
$0 < a < 1$, $x > 0$, $y$ arbitrary.  \ 
Using LAW 4 below (whose proof does \un{not} depend on LAW 2), write
\[
(a^x)^y
\ = \ 
\bigg(\frac{1}{(1/a)^x}\bigg)^y 
\ = \ 
\frac{1}{((1/a)^x)^y}
\ = \ 
\frac{1}{(1/a)^{xy}}
\ = \ 
a^{xy}.
\]
\vspace{0.2cm}

\qquad \un{Case 5:} \ 
$0 < a$, $x < 0$, $y$ arbitrary. \ 
If $x < 0$, then $-x > 0$, hence
\[
(a^x)^y
\ = \ 
\bigg(\frac{1}{a^{-x}}\bigg)^y 
\ = \ 
\frac{1}{(a^{-x})^y}
\ = \ 
\frac{1}{a^{-xy}}
\ = \ 
a^{xy}.
\]
\\[15pt]
\qquad \un{LAW 3:} \quad
One need only observe that
\allowdisplaybreaks
\begin{align*}
a^x \ 
&=\ 
x^{x-y+y}
\\[12pt]
&=\ 
a^{x-y} \cdot a^y \qquad \text{(by LAW 1)},
\end{align*}
i.e., 
\[
\frac{a^x}{a^y} \ = \ a^{x-y}.
\]
\\[15pt]
\qquad \un{LAW 4:} \quad
\vspace{0.2cm}

\qquad \un{Case 1:} \ 
$a > 1$, $b > 1$.  \ 
Let
\allowdisplaybreaks
\begin{align*}
S \ 
&=\ 
\{a^s : s \in \Q \ \text{and} \ s < x\}
\\[12pt]
T \ 
&=\ 
\{b^t : t \in \Q \ \text{and} \ t < x\}
\\[12pt]
U \ 
&=\ 
\{(ab)^u : u \in \Q \ \text{and} \ u < x\},
\end{align*}
thus $a^x = \sup S$, $b^x = \sup T$, $(ab)^x = \sup U$.  
Meanwhile, 
\allowdisplaybreaks
\begin{align*}
a^x b^x \ 
&=\ 
(\sup S) \cdot (\sup T) 
\\[12pt]
&=\ 
(\sup S \cdot T) \qquad \text{(cf. \S3, \#6)}.
\end{align*}
So, to prove that $(ab)^x = a^x b^x$, it will be enough to prove that $\sup (S \cdot T) = \sup U$ and for this purpose, we shall employ \S3, \#2.  
Since \mU is a subset of $S \cdot T$, it suffices to go the other way.  
But a generic element of $S \cdot T$ is of the form $a^s b^t$, where $s$, $t \in \Q$ and $s < x$, $t < x$.  
And, assuming that $s \leq t$, we have
\[
a^s b^t \ \leq \ a^t b^t \ = \ (ab)^t \in U.
\]

\qquad \un{Case 2:} \ 
$0 < a < 1$, $0 < b < 1$.   
Since $0 < ab < 1$, from the definitions, 
\[
(ab)^x \ = \  \frac{1}{(1 / ab)^x}.
\]
Since $1/a > 1$, $1/b > 1$, it follows from the discussion in Case 1 that
\[
\bigg(\frac{1}{ab}\bigg)^x
\ = \  
\bigg(\frac{1}{a}\bigg)^x \bigg(\frac{1}{b}\bigg)^x.
\]
Therefore
\allowdisplaybreaks
\begin{align*}
(ab)^x \ 
&=\ 
 \frac{1}{(1 / ab)^x} 
 \\[12pt]
 &=\ 
 \frac{1}{(1/a)^x \cdot (1/b)^x} 
  \\[12pt]
 &=\ 
 a^x b^x.
\end{align*}
\vspace{0.2cm}

\qquad \un{Case 3:} \ 
$0 < a < 1$, $b > 1$.  
In this situation $1/a > 1$.  
Suppose first that $1 < 1/a \leq b$ $-$then $ab \geq 1$, so
\[
b^x \ = \ \bigg(ab \cdot \frac{1}{a}\bigg)^x \ = \ (ab)^x \bigg(\frac{1}{a}\bigg)^x,
\]
hence
\[
(ab)^x \ = \ b^x \frac{1}{(1/a)^x} \ = \ \frac{1}{(1/a)^x} b^x \ = \ a^x b^x.
\]
The other possibility is that $1 < b < 1/a$.  
Since in this situation both $1/ab$ and $b$ are greater than 1, we have
\[
\bigg(\frac{1}{a}\bigg)^x \ = \ \bigg(\frac{1}{ab} \cdot b\bigg)^x \ = \ \bigg(\frac{1}{ab}\bigg)^x b^x,
\]
so
\[
(ab)^x 
\ = \ 
\frac{1}{(1/ab)^x} 
\ = \ 
\frac{1}{(1/a)^x} b^x 
\ = \ 
a^x b^x.
\]

\vspace{0.2cm}
\qquad \un{Case 4:} \ 
$a > 1$, $0 < b < 1$. \ 
This is the same as Case 3 with the roles of $a$ and $b$ interchanged.
\vspace{0.3cm}

A simple but important consequence of LAW 4, used already in Case 4 of LAW 2 above, is the fact that
\[
\bigg(\frac{1}{a}\bigg)^x 
\ = \ 
\frac{1}{a^x} \qquad (a > 0, \ x \in \R).
\]
Proof:
\[
1^x \ = \  \bigg(a \cdot \frac{1}{a}\bigg)^x  \ = \ a^x \bigg(\frac{1}{a}\bigg)^x 
\implies 
\bigg(\frac{1}{a}\bigg)^x  \ = \ \frac{1}{a^x}.
\]
\\[15pt]
\qquad \un{LAW 5:} \quad
Write
\[
\bigg(\frac{a}{b}\bigg)^x
\ = \
 \bigg(a \cdot \frac{1}{b}\bigg)^x
\ = \
a^x \bigg(\frac{1}{b}\bigg)^x 
\ = \
a^x \frac{1}{b^x}
\ = \
 \frac{a^x}{b^x}.
\]
\\[15pt]
\qquad \un{LAW 6:} \quad
We shall consider (i), leaving (ii) for the reader, and of the two parts to (i), only the assertion 
$0< a < b\implies a^x < b^x$ will be dealt with explicitly.

Claim: \ If $c > 1$, $x > 0$, then $c^x > 1$.  
Granting the claim for the moment, note now that
\allowdisplaybreaks
\begin{align*}
0 < a < b 
&\implies 1 \ < \ \frac{b}{a} 
\\[12pt]
&\implies
1 < \bigg(\frac{b}{a}\bigg)^x \ = \ \frac{b^x}{a^x} \qquad \text{(by LAW 5)}
\\[12pt]
&\implies
a^x < b^x.
\end{align*}
Going back to the claim, fix a rational number $r$ such that $0 < r < x$ $-$then it will be enough to prove that $1 < c^r$.  
Since $1 < 2$ $\implies$ $r < 2r$ $\implies$ $c^r < c^{2r}$, it
follows that
\[
1 
\ = \ 
c^{r - r} 
\ < \ 
c^{2r - r} 
\ = \ c^r.
\]
\\[15pt]
\qquad \un{LAW 7:} \quad 
We shall consider (i), leaving (ii) for the reader, and of the two parts to (i), only the assertion 
$x < y \implies a^x < a^y$ will be dealt with explicitly.  
Choose $s \in \Q$ : $x < s < y$ $-$then
\[
r \in \Q \quad \text{and} \quad r < x \implies r < s \implies a^r < a^s \implies a^x \leq a^s.
\]
Choose $t \in \Q$ : $s < t < y$ $-$then $a^s < a^t$ and $a^t \leq a^y$, hence $a^x < a^y$.
\vspace{0.3cm}

\begin{x}{\small\bf LEMMA} \ 
Let $a > 0$, $a \neq 1$, $-$then 
\[
\frac{a^x - 1}{x}
\ < \ 
\frac{a^y - 1}{y}
\]
for all $x, y \in \R - \{0\}$ with $x < y$ (cf. \#9).
\end{x}
\vspace{0.3cm}

%% file: _05_exp_a_and_log_a.tex
\chapter{
$\boldsymbol{\S}$\textbf{5}.\quad  $\boldsymbol{\exp_a}$ AND $\boldsymbol{\log_a}$}
\setlength\parindent{2em}
\setcounter{theoremn}{0}
\renewcommand{\thepage}{\S5-\arabic{page}}

\ \indent 
Let $a \neq 1$ be a positive real number.
\vspace{0.3cm}

\begin{x}{\small\bf DEFINITION} \ 
The 
\un{exponential function to base $a$}
\index{exponential function to base $a$} 
is the function $\exp_a$ with domain $\R$ defined by the rule
\[
\exp_a (x) \ = \ a^x \qquad (x \in \R).
\]
\end{x}
\vspace{0.3cm}

\begin{x}{\small\bf LEMMA}  \ 
$\exp_a : \R \ra \R_{> 0}$ \ is injective (cf. \S4, \#15, 7(i) and 7(ii)).
\end{x}
\vspace{0.3cm}

\begin{x}{\small\bf LEMMA} \ 
$\exp_a : \R \ra \R_{> 0}$ \ is surjective.
\end{x}
\vspace{0.3cm}

This is not quite immediate and requires some preparation.
\vspace{0.3cm}

\begin{x}{\small\bf SUBLEMMA} \ 
Let $n > 1$ be a natural number and let $a \neq 1$ be a positive real number $-$then
\[
n (a^{1/n} - 1) \ < \ a - 1.
\]
\vspace{0.2cm}

PROOF \ 
In \S4, \#9, take $r = \ds\frac{1}{n}$, $s = 1$, then $r < s$ and 
\[
\frac{a^{1/n} - 1}{\frac{1}{n}} \ < \ \frac{a - 1}{1},
\]
i.e., 
\[
n (a^{1/n} - 1) \ < \ a - 1.
\]
\end{x}
\vspace{0.3cm}

To discuss \#3, distinguish two cases: $a > 1$ or $a < 1$.  
We shall work through the first of these, leaving the second to the reader.
\vspace{0.3cm}


\begin{x}{\small\bf SUBLEMMA} \ 
If $t > 1$ and 
\[
n \ > \ \frac{a - 1}{t - 1},
\]
then $\ds a^{1/n} < t$.
\vspace{0.2cm}

PROOF \ 
In fact, 
\[
a - 1 \ > \ n \big(a^{1/n} - 1\big) \ > \ \frac{a-1}{t-1} \big(a^{1/n} - 1\big)
\]
\qquad\qquad $\implies$
\[
1 \ > \ \frac{a^{1/n} - 1}{t - 1}
\]
\qquad\qquad $\implies$ 
\[
t - 1 \ > \ a^{1/n} - 1
\]
\qquad\qquad $\implies$ 
\[
t > a^{1/n}.
\]
\vspace{0.2cm}

Fix $y > 0$ $-$then the claim is that there is a real number $x$ such that $a^x = y$ ($x$ then being necessarily unique).  
So let
\[
S \ = \ \{w : a^w < y\}
\]
and put $x = \sup S$.
\vspace{0.2cm}

\qquad \textbullet \quad
$a^x < y$ is untenable.
\vspace{0.2cm}

[In \#5, take $t = \ds\frac{y}{a^x} > 1$ to get
\[
a^{1/n} \ < \ \frac{y}{a^x}
\]
for $n \gg 0$, thus
\[
a^{x + \frac{1}{n}} \ < \ y
\]
for $n \gg 0$.  
But then, for any such $n$, 
\[
x + \frac{1}{n} \in S
\]
which leads to the contradiction $x \geq x + \ds\frac{1}{n}$.]
\vspace{0.2cm}

\qquad \textbullet \quad
$a^x > y$ is untenable.
\vspace{0.2cm}

[In \#5, take $t = \ds\frac{a^x}{y} > 1$ to get
\[
a^{1/n} \ < \ \frac{a^x}{y}
\]
for $n \gg 0$, thus
\[
y \ < \ a^{x - \frac{1}{n}}
\]
for $n \gg 0$.  
Owing to \S2, \#6, for each $n \gg 0$, there exists $w_n \in S$: $w_n > x - \ds\frac{1}{n}$, hence
\allowdisplaybreaks
\begin{align*}
y \ 
&>\ 
a^{w_n} \
\\[12pt]
&>\ 
a^{x - \frac{1}{n}} \qquad \text{(cf. \S4, \#15, 7(i))}
\\[12pt]
&>\ 
y, 
\end{align*}
a contradiction.]
\vspace{0.2cm}

Therefore $a^x = y$, as contended.
\end{x}
\vspace{0.3cm}

\begin{x}{\small\bf SCHOLIUM} \ 
$\exp_a : \R \ra \R_{> 0}$ is bijective.

\end{x}
\vspace{0.3cm}

\begin{x}{\small\bf REMARK}\ 
There is another way to establish the surjectivity of $\exp_a$ if one is willing to introduce some machinery, the point being that the range of $\exp_a$  is an open subgroup of $\R_{> 0}$.  
One may then quote the following generality: \ A locally compact topological group is connected if and only if it has no proper open subgroups.
\end{x}
\vspace{0.3cm}

Since
\[
\exp_a : \R \ra \R_{> 0}
\]
is bijective, it admits an inverse.
\[
\exp_a^{-1} : \R_{> 0} \ra \R.
\]
\vspace{0.3cm}

\begin{x}{\small\bf NOTATION} \ 
Put
\[
\log_a  \ = \ \exp_a^{-1}.
\]
\end{x}
\vspace{0.3cm}

\begin{x}{\small\bf DEFINITION} \ 
The 
\un{logarithm function to a base $a$}
\index{logarithm function to a base $a$} 
is the function $\log_a$ defined by the rule
\[
\log_a (a^x) \ = \ x \qquad (x \in \R).
\]
\end{x}
\vspace{0.3cm}

\begin{x}{\small\bf LEMMA} \ 
Let $u$ and $v$ be positive real numbers $-$then
\[
\begin{cases}
\ \log_a (u v) \ = \ \log_a(u) + \log_a (v)\\[3pt]
\ \log_a \big(\frac{u}{v}\big) \ = \ \log_a(u) - \log_a (v)
\end{cases}
.
\]
\end{x}
\vspace{0.3cm}

\begin{x}{\small\bf LEMMA} \ 
Let $y$ be a positive real number, $r$ a real number $-$then
\[
\log_a (y^r) \ = \ r \hsx \log_a (y).
\]
\vspace{0.2cm}

PROOF \ 
Write $y = a^x$, thus
\allowdisplaybreaks
\begin{align*}
y^r \ 
&=\ 
(a^x)^r
\\[12pt]
&=\ 
a^{x r} \qquad \text{(cf. \S4, \#15, (2))}
\\[12pt]
&=\ 
a^{r x}
\end{align*}
\qquad\qquad $\implies$
\[
\log_a (y^r) \ = \ r x  \ = \ r \hsx \log_a (y).
\]
\end{x}
\vspace{0.3cm}

\begin{x}{\small\bf \un{N.B.}} \ 
Special cases:
\[
\log_a (1) \ = \ 0, \quad \log_a(a) \ = \ 1.
\]
\end{x}
\vspace{0.3cm}

\begin{x}{\small\bf LEMMA} \ 
Let $a \neq 1$, $b \neq 1$ be positive real numbers $-$then
\[
\log_a (b) \hsx  \log_b (a) \ = \ 1.
\]
\vspace{0.2cm}

PROOF \ 
Put
\[
x \ = \ \log_a (b), \quad y \ = \ \log_b (a), 
\]
so that
\[
a^x = b, \quad b^y = a, 
\]
hence
\[
a \ = \ b^y \ = \ \big(a^x\big)^y \ = \ a^{x y} \qquad \text{(cf. \S4, \#15, (2))}
\]
from which $x y = 1$.
\end{x}
\vspace{0.3cm}

\begin{x}{\small\bf DEFINITION} \ 
The 
\un{common logarithm}
\index{common logarithm} 
is $\log_{10}$.
\end{x}
\vspace{0.3cm}

\begin{x}{\small\bf EXAMPLE} \ 
$\log_{10} \hsx 2$ is irrational.
\vspace{0.2cm}

[Suppose that 
\[
\log_{10} \hsx 2 \ = \ \frac{a}{b},
\]
where $a$ and $b$ are positive integers $-$then
\[
2 \ = \ 10^{^{\frac{a}{b}}} 
\implies
2^b \ = \ 10^a \ = \ 2^a 5^a.
\]
But $2^b$ is not divisible by 5.]
\vspace{0.2cm}

[Note: \ 
It turns out that $\log_{10} \hsx 2$ is transcendental, a point that will be dealt with later on.]
\end{x}
\vspace{0.3cm}

There are irrational numbers $\alpha$, $\beta$ such that $\alpha^\beta$ is rational.
\vspace{0.3cm}

\begin{x}{\small\bf EXAMPLE} \ 
Take $\alpha = \sqrt{10}$ (cf. \S7, \#6), $\beta = 2 \log_{10} 2$ $-$then
\allowdisplaybreaks
\begin{align*}
(\sqrt{10})^{2 \log_{10} 2} \ 
&=\ 
\big(10^{^{\frac{1}{2}}}\big)^{2 \log_{10} 2}
\\[12pt]
&=\ 
10^{\log_{10} 2}
\\[12pt]
&=\ 
2.
\end{align*}
\end{x}
\vspace{0.3cm}

\[
\text{APPENDIX}
\]
\vspace{0.3cm}

Put
\[
E(x) \ = \ 
\sum\limits_{k=0}^\infty \hsx
\frac{x^k}{k!} \qquad (x \in \R).
\]
\vspace{0.2cm}

[Note: \ 
\[
E(1) \ \equiv \ e.]
\]
\vspace{0.3cm}

\qquad {\small\bf LEMMA} \ 
$E(x_1 + \cdots + x_n) \ = \ E(x_1) \cdots E(x_n)$.
\vspace{0.2cm}

[Note: \ 
\[
E(x) E(-x) \ = \ E(x - x) \ = \ E(0) \ = \ 1.]
\]
\vspace{0.2cm}

Take $x_1 = 1, \ldots, x_n = 1$ to get
\[
E(n) \ = \ e^n.
\]
If now $r = \ds\frac{m}{n}$ $(m, n \in \N)$, then 
\allowdisplaybreaks
\begin{align*}
(E(r))^n \ 
&=\ 
E(nr) = E(m) = e^m
\\[12pt]
&\implies 
E(r) = e^{^{\frac{m}{n}}} = e^r.
\end{align*}
And
\[
E(-r) \ = \ \frac{1}{E(r)} \ = \ \frac{1}{e^r} \ = \ e^{-r}.
\]
\vspace{0.2cm}

Summary:
\[
E(x) \ = \ e^x \qquad (x \in \Q).
\]

But now for any real $x$,
\[
e^x \ = \ \sup S,
\]
where
\[
S \ = \ \big\{e^r : r \in \Q \ \text{and} \ r < x \big\} \qquad \text{(cf. \S4, \#13)}.
\]
\vspace{0.3cm}

\qquad {\small\bf THEOREM} \ 
$\forall \ x \in \R$, 
\[
E(x) \ = \ e^x \qquad (= \exp_e(x)).
\]
\vspace{0.3cm}

\qquad {\small\bf REMARK} \ 
It can be shown that
\[
e \ = \ \sup\bigg\{\bigg(1 + \frac{1}{n}\bigg)^n \hsx : \hsx n \in \N\bigg\},
\]
a fact which is sometimes used as the definition of $e$.
\vspace{0.3cm}


%% file: _06_irrationality_of_sqrt_2.tex
\chapter{
$\boldsymbol{\S}$\textbf{6}.\quad  IRRATIONALITY OF  $\boldsymbol{\sqrt{2}}$}
\setlength\parindent{2em}
\setcounter{theoremn}{0}
\renewcommand{\thepage}{\S6-\arabic{page}}

\ \indent 
Recall that $\PP$ is the subset of $\R$ whose elements are irrational and, on abstract grounds, is uncountable, in particular, irrational numbers exist.  
Still, the problem of deciding whether a specific real number is irrational or not is generally difficult.
\vspace{0.3cm}

\begin{x}{\small\bf RAPPEL} \ 
$\sqrt{2}$ exists (cf. \S4, \#4).
\end{x}
\vspace{0.3cm}

\begin{x}{\small\bf THEOREM} \ 
$\sqrt{2}$ is irrational.
\end{x}
\vspace{0.3cm}

There are many proofs of this result.  
In what follows we shall give a representative sampling.
\vspace{0.3cm}

\qquad \un{First Proof:} \ 
Suppose that $\sqrt{2}$ is rational, say $\sqrt{2} = \ds\frac{x}{y}$, where $x$ and $y$ are positive integers and $\gcd(x,y) = 1$, 
$-$then $\ds\frac{x^2}{y^2} = 2$ or still, $x^2 = 2 y^2$, thus $2 \big | x^2$ and $x^2$ is even.  
But then $x$ must be even (otherwise, $x$ odd forces $x^2$ odd), so $x = 2n$ for some positive integer $n$.  And:
\allowdisplaybreaks
\begin{align*}
x^2 = 2 y^2 
&\implies
(2 n)^2 = 2 y^2
\\[12pt]
&\implies
 2n^2 = y^2 
 \\[12pt]
 &\implies
 2 \big | y^2 
  \\[12pt]
  &\implies
  2 \big | y.
\end{align*}
\vspace{0.3cm}
Therefore $\gcd(x,y) \neq 1$, a contradiction.
\vspace{0.3cm}

\qquad \un{Second Proof:} \ 
Suppose that $\sqrt{2}$ is rational, say $\sqrt{2} = \ds\frac{x}{y}$, where $x$ and $y$ are positive integers and $y$ is the smallest such $-$then $\ds\frac{x^2}{y^2} = 2$ or still, $x^2 = 2 y^2$.  
Next
\allowdisplaybreaks
\begin{align*}
y^2 \ 
&<\ 
2 y^2 = x^2 = (2y) y < (2y) (2y)
\\[12pt]
&\implies
y^2 < x^2 < (2y)^2
\\[12pt]
&\implies
y < x < 2y \qquad \text{(cf. \S4, \#1, 6(i))}.
\end{align*}
Put $u = x - y$, a positive intger:
\[
y + u 
\ = \ 
x < 2y 
\ = \ 
y + y 
\implies 
u < y.
\]
Put $v = 2y - x$, a positive integer:
\allowdisplaybreaks
\begin{align*}
v^2 - 2 u^2 \ 
&=\ 
(2y - x)^2 - 2(x - y)^2
\\[12pt]
&=\ 
4y^2 - 4yx + x^2 - 2(x^2 - 2xy + y^2)
\\[12pt]
&=\ 
4y^2 + x^2 - 2x^2 - 2y^2
\\[12pt]
&=\ 
(x^2 - 2y^2) - 2(x^2 - 2y^2)
\\[12pt]
&=\ 
(1 - 2)(x^2 - 2y^2)
\\[12pt]
&=\ 
(-1) (0)
\\[12pt]
&=\ 
0.
\end{align*}
\qquad\qquad $\implies$
\allowdisplaybreaks
\begin{align*}
v^2 \ = \ 2u^2
&\implies 
\frac{v^2}{u^2} \ = \ 2
\\[12pt]
&\implies  
\bigg(\frac{v^2}{u^2}\bigg)^{1/2}  \ = \ 2^{1/2} \ = \ \sqrt{2} 
\\[12pt]
&\implies 
\frac{v^{2(1/2)}}{u^{2(1/2)}} \ =  \ \sqrt{2} \qquad \text{(cf. \S4, \#7, 5)} 
\\[12pt]
&\implies 
\frac{v}{u} \ = \ \sqrt{2}.
\end{align*}
But now we have reached a contradiction: $u$ is less than $y$ whereas $y$ was the smallest positive integer with the property that 
$\ds\frac{x}{y} = \sqrt{2}$ for some positive integer $x$.  
\vspace{0.3cm}

\qquad \un{Third Proof:} \ 
Suppose that $\sqrt{2}$ is rational, say $\sqrt{2} = \ds\frac{x}{y}$, where $x$ and $y$ are positive integers.  
Write
\[
\sqrt{2} + 1 \ = \ \frac{1}{\sqrt{2} - 1},
\]
thus
\[
\frac{x}{y} + 1 
\ = \ 
\frac{1}{\frac{x}{y} - 1}
\ = \ 
\frac{y}{x - y}
\]
\qquad\qquad $\implies$
\[
\sqrt{2} 
\ = \ 
\frac{x}{y}
\ = \ 
\frac{y}{x - y} - 1
\ = \ 
\frac{2y - x}{x - y}
\ \equiv \ 
\frac{x_1}{y_1}.
\]
But
\[
1  \ < \ \sqrt{2}  \ < \ 2 
\implies
1  \ < \ \frac{x}{y}  \ < \ 2 
\implies 
y  \ < \ x  \ < \ 2y
\]
\[
\implies 
\begin{cases}
\ x_1 = 2y - x > 0\\
\ y_1 = x - y > 0
\end{cases}
\implies 
\begin{cases}
\ x_1 \in \N\\
\ y_1 \in \N
\end{cases}
.
\]
In addition
\[
2 y \ < \ 2x \ = \ x + x 
\implies 
2y - x  \ < \ x 
\implies
x_1 \  < \ x.
\]
Proceeding, there exist positive integers $x_2$ and $y_2$ such that
\[
\sqrt{2} 
\ = \ 
\frac{x_1}{y_1} 
\ = \ 
\frac{2 y_1 - x_1}{x_1 - y_1} 
\ \equiv \ 
\frac{x_2}{y_2}
\]
with $x_2 < x_1 < x$.  And so on, ad infinitum.  The supposition that $\sqrt{2}$ is irrational 
therefore leads to an infinite descending chain of natural numbers, an impossibility.
\vspace{0.3cm}

\qquad \un{Fourth Proof:} \ 
Suppose that $\sqrt{2}$ is rational, say $\sqrt{2} = \ds\frac{x}{y}$, where $x$ and $y$ are positive integers.  
Define sequences
\[
\begin{cases}
\ a_1, a_2, \ldots \\[3pt]
\ b_1, b_2, \ldots
\end{cases}
\]
of natural numbers recursively by 
\[
\begin{cases}
\ a_1 = 1,\ a_2 = 2, \ a_n = 2 a_{n-1} + a_{n-2} \qquad (n > 2) \\
\ b_1 = 1,\ b_2 = 3, \ b_n = 2 b_{n-1} + b_{n-2} \ \hspace{0.85cm} (n > 2)
\end{cases}
.
\]
Put 
\[
p_n(t) \ = \ a_n^2 t^2 - b_n^2 \qquad (n \geq 1).
\]
Then
\[
p_n(\sqrt{2}) \ = \ 2 a_n^2 - b_n^2
\]
is an integer and $\abs{p_n(\sqrt{2})} = 1$ (details below).  
On the other hand, 
\allowdisplaybreaks
\begin{align*}
1\ 
&=\ 
\abs{p_n(\sqrt{2})}
\\[12pt]
&=\ 
\abs{\bigl(a_n \sqrt{2} - b_n \bigr) \bigl(a_n \sqrt{2} + b_n\bigr)}
\\[12pt]
&=\ 
\abs{\bigl(a_n \frac{x}{y} - b_n \bigr) \bigl(a_n \frac{x}{y} + b_n\bigr)}
\\[12pt]
&=\ 
\abs{a_n x - b_n y} \hsx \bigg(\frac{a_n x + b_n y}{y^2}\bigg)
\end{align*}
\qquad\qquad $\implies$
\[
0 \ < \  \abs{a_n x - b_n y} 
\ = \ 
\frac{y^2}{a_n x + b_n y}.
\]
Since the sequence $\{a_n x + b_n y\}$ is strictly increasing, from some point on
\[
y^2 \ < \ a_n x + b_n y.
\]
I.e.: 
\[
n \gg 0 \implies \abs{a_n x  -b_n y} \ < \ 1.
\]
But there are no integers between 0 and 1.
\vspace{0.2cm}

[Inductively we claim that
\[
2 a_n^2 - b_n^2 \ = \ (-1)^{n+1} 
\quad \text{and} \quad
2 a_{n-1} a_n - b_{n-1}b_n \ = \ (-1)^n.
\]
These identities are certainly true when $n=1$ (take $a_0 = 0, b_0 = 1$).  
Assume therefore that they hold at level $n > 1$ $-$then at level $n + 1$:
\allowdisplaybreaks
\begin{align*}
2 a_{n+1}^2 - b_{n+1}^2 \ 
&=\ 
2(2 a_n + a_{n-1})^2 - (2 b_n + b_{n-1})^2
\\[12pt]
&=\ 
4(2 a_n^2 - b_n^2) + 4(2 a_{n-1} a_n - b_{n-1} b_n) + (2 a_{n-1}^2 - b_{n-1}^2)
\\[12pt]
&=\ 
4(-1)^{n+1} + 4(-1)^n + (-1)^n
\\[12pt]
&=\ 
(-1)^n
\\[12pt]
&=\ 
(-1)^{n+2}.
\end{align*}
And, analogously, 
\[
2 a_n a_{n+1} - b_n b_{n+1} \ = \ (-1)^{n+1}.  
\]
Finally
\allowdisplaybreaks
\begin{align*}
p_n(\sqrt{2}) \ 
&=\ 
2 a_n^2 - b_n^2 = (-1)^{n+1} 
\\[12pt]
&\implies \abs{p_n(\sqrt{2})} 
\\[12pt]
&=\ 
1.]
\end{align*}
\vspace{0.3cm}

\qquad \un{Fifth Proof:} \ 
Let \mS be the set of positive integers $n$ with the property 
that $n \sqrt{2}$ is a positive integer.  
If $\sqrt{2}$ were rational, then \mS would be nonempty, hence would have a smallest element, call it $k$.  
Now, from the definitions, 
\[
k \in S \implies (\sqrt{2} - 1) k \in \N.
\]
But 
\allowdisplaybreaks
\begin{align*}
((\sqrt{2} - 1) k) \sqrt{2} \ 
&=\ 
2k - k\sqrt{2}
\\[12pt]
&=\ 
(2 - \sqrt{2}) k
\end{align*}
is a positive integer, so $(\sqrt{2} - 1) k \in S$.  However
\[
(\sqrt{2} - 1) k \ < \ (2 - 1) k \ = \ k,
\]
which contradicts the assumption that $k$ is the smallest element of \mS.

%% file: _07_irrationality_theory_and_examples.tex
\chapter{
$\boldsymbol{\S}$\textbf{7}.\quad  IRRATIONALITY: \ THEORY AND EXAMPLES}
\setlength\parindent{2em}
\setcounter{theoremn}{0}
\renewcommand{\thepage}{\S7-\arabic{page}}

\ \indent 
For use below:
\vspace{0.3cm}

\begin{x}{\small\bf RAPPEL} \ 
Let $a$, $b$, $c$ be integers such that $a$, $b$ have no prime factors in common and $a \big | b^n c$ $(n \in \N)$ $-$then 
$a \big | c$.
\end{x}
\vspace{0.3cm}

The following result is the so-called ``rational roots test''.
\vspace{0.3cm}

\begin{x}{\small\bf THEOREM}\ 
Let 
\[
f(X) \ = \ a_0 + a_1 X + a_2 X^2 + \cdots + a_n X^n
\]
be a polynomial with integral coefficients.  
Suppose that it has a rational root $\ds\frac{p}{q}$ : $p$, $q \in \Z$ and $\gcd(p,q) = 1$  $-$then $p \big | a_0$ and $q \big | a_n$. 
\vspace{0.2cm}

PROOF \ 
Take $X = \ds\frac{p}{q}$ to get
\[
a_0 + a_1 \bigg(\frac{p}{q}\bigg) + a_2 \bigg(\frac{p}{q}\bigg)^2 + \cdots + a_n \bigg(\frac{p}{q}\bigg)^n  \ = \ 0
\]
so, after multiplying through by $q^n$, 
\allowdisplaybreaks
\begin{align*}
q^n a_0 \ 
&=\ 
- (a_1 p q^{n-1} + a_2 p^2 q^{n-2} + \cdots + a_n p^n)
\\[12pt]
&=\ 
- p(a_1 q^{n-1} + a_2 p q^{n-2} + \cdots + a_n p^{n-1}) \in \Z
\\[12pt]
&\qquad\qquad \implies
p \big | q^n a_0  
\\[12pt]
&\qquad\qquad \implies
p \big | a_0 \qquad \text{(cf. \#1)}.
\end{align*}
That $q \big| a_n$ can be established analogously.
\end{x}
\vspace{0.3cm}

\begin{x}{\small\bf \un{N.B.}} \ 
When specialized to the case where $a_n = 1$, the conclusion is
that if the polynomial 
\[
a_0 + a_1 X + a_2 X^2 + \cdots +  X^n
\]
has a rational root, then this root is an integer (which divides $a_0$).
\vspace{0.2cm}

[Consider a rational root $\ds\frac{p}{q}$ and take $q$ positive (in the event that $q$ were negative absorb the minus sign into $p$).    
From the above, $q$ divides $a_n = 1$, hence $q = 1$, hence $\ds\frac{p}{q} = \ds\frac{p}{1} = p$ (and $p \big | a_0$).]
\end{x}
\vspace{0.3cm}

\begin{x}{\small\bf EXAMPLE} \ 
If $p$ is a prime, then $\sqrt{p}$ is irrational.
\vspace{0.2cm}

[Consider the polynomial $X^2 - p$, thus $\bigl(\sqrt{p}\bigr)^2 - p = 0$, i.e., $\sqrt{p}$ is a root.  
Suppose that $\sqrt{p}$ was rational so for some $k \in \N$, 
\[
\sqrt{p} \ = \ k \implies p  \ = \ k^2.
\]
But $k^2$ has an even number of prime factors, from which it follows that the stated relation is impossible 
(or quote \#1: \ $a = p$, $b = k$, $n = 2$, $c = 1$, implying that $p \big |1$).]
\end{x}
\vspace{0.3cm}

Therefore in particular $\sqrt{2}$ and $\sqrt{3}$ are irrational but this does not automatically imply that $\sqrt{2} \ + \ \sqrt{3}$ is irrational (the sum of two irrationals may be either rational or irrational).
\vspace{0.3cm}

\begin{x}{\small\bf EXAMPLE} \ 
$\ \sqrt{2} \ + \ \sqrt{3}\ $ \ is  irrational.
\vspace{0.2cm}

[$\sqrt{2} \ + \ \sqrt{3}$ is a zero of the function 
\[
X^2 - 2X \sqrt{2} - 1,
\]
so $\sqrt{2} \ + \ \sqrt{3}\ $   is a root of the polynomial 
\[
(X^2 + 2X \sqrt{2} - 1) (X^2 - 2X\sqrt{2} - 1) \ = \ X^4 - 10X^2 + 1.
\]
From the above, the only possible rational roots of this polynomial are integers which divide 1, i.e., $\pm 1$.  
And $\sqrt{2} \ + \ \sqrt{3} \neq \pm 1$, thus $\sqrt{2} \ + \ \sqrt{3}\ $ is not among the possible roots of 
\[
X^4 - 10 X^2 + 1,
\]
thus is  irrational.]
\end{x}
\vspace{0.3cm}

\begin{x}{\small\bf EXAMPLE} \ 
Let $a$ and $n$ be positive integers $-$then $\sqrt[\leftroot{-1}\uproot{3}n]{a}$ is either irrational or a positive integer.  
And if $\sqrt[\leftroot{-1}\uproot{3}n]{a}$ is a positive integer, then $a$ is the $n^\nth$ power of a positive integer.
\vspace{0.2cm}

[Consider the polynomial $X^n - a$, hence ($\sqrt[\leftroot{-1}\uproot{3}n]{a})^n - a = a - a = 0$.  
There are now two possibilities, viz. either $\sqrt[\leftroot{-1}\uproot{3}n]{a}$ is irrational or else 
$\sqrt[\leftroot{-1}\uproot{3}n]{a}$ is rational in which case $\sqrt[\leftroot{-1}\uproot{3}n]{a} \equiv k$ is a positive integer 
(and $a = k^n$).]
\end{x}
\vspace{0.3cm}

\begin{x}{\small\bf REMARK} \ 
Consequently, if $a$ is a positive integer such that $\sqrt{a}$ is not a positive integer, then $\sqrt{a}$ is irrational (cf. \#4).
\vspace{0.2cm}

[Here is another proof.  Assume instead that $\sqrt{a}$ is rational, say $\sqrt{a} = \ds\frac{x}{y}$, where $x$ and $y$ are positive integers and $y$ is the smallest such:
\[
y \sqrt{a} \ = \  x \implies (y\sqrt{a}\hsx) \sqrt{a} \ = \ x \sqrt{a} \implies y a \ = \  x \sqrt{a}.
\]
Choose $n \in \N$ : $n < \sqrt{a} < n + 1$ $-$then
\allowdisplaybreaks
\begin{align*}
\sqrt{a}\ 
&=\ 
\frac{x}{y}
\\[12pt]
&=\ 
\frac{x(\sqrt{a} - n)}{y(\sqrt{a} - n)}
\\[12pt]
&=\ 
\frac{x \sqrt{a} - x n}{y \sqrt{a} - y n}
\\[12pt]
&=\ 
\frac{ya - xn}{x - yn}.
\end{align*}
The numerator and denominator of the fraction 
\[
\frac{ya - xn}{x - yn}
\]
are integers that, in fact, are positive:
\[
\begin{cases}
\ ya - xn = x \sqrt{a} - xn = x ( \sqrt{a} - n) > 0\\[3pt]
\ \ x - yn = y \sqrt{a} - yn = y(\sqrt{a} - n) > 0
\end{cases}
.
\]
And
\[
x - yn \ = \ y (\sqrt{a} - n) \ < \ y
\]
which contradicts the choice of $y$.]
\end{x}
\vspace{0.3cm}

\begin{x}{\small\bf THEOREM} \ 
Suppose that $a_1$, $a_2$, \ldots, $a_n$ are positive integers.  Assume: 
\[
\Sigma \ \equiv \ \sqrt{a_1} \ + \ \sqrt{a_2} \ + \cdots + \sqrt{a_n}
\]
is rational $-$then 
$\sqrt{a_1}, \hsx \sqrt{a_2}, \hsx \ldots , \sqrt{a_n}$ are rational.
\end{x}
\vspace{0.3cm}

\begin{x}{\small\bf APPLICATION} \ 
If for some $k$ $(1 \leq k \leq n)$, $\sqrt{a_k}$ is irrational, then
\[
\sqrt{a_1} \ + \ \sqrt{a_2} \ + \cdots + \sqrt{a_n}
\]
is irrational.
\end{x}
\vspace{0.3cm}

\begin{x}{\small\bf EXAMPLE} \ 
$\  \sqrt{2} +\sqrt{3}$ is irrational (cf. \#5).
\end{x}
\vspace{0.3cm}

\begin{x}{\small\bf EXAMPLE} \ 
$\  \sqrt{2} +\sqrt{3}  + \sqrt{5}$ is irrational.
\vspace{0.3cm}
\end{x}

Passing to the proof of \#8, it will be enough to show that $\sqrt{a_1}$ is rational.  For this purpose, introduce
\[
F(X;a_1) \ = \ \Pi (X - \sqrt{a_1} \hsx \pm \sqrt{a_2} \pm \cdots \pm \sqrt{a_n}\hsx),
\]
where the product ranges over all combinations of plus and minus signs, thus
\[
F(\Sigma; a_1) \ = \ 0.
\]
Next multiply out the expression defining $F(X;a_1)$ $-$then $\sqrt{a_1}$ appears to both even and odd powers but 
$\sqrt{a_2}\hsx, \ldots, \sqrt{a_n}$ appear only to even powers.  
Assemble the even powered terms in $\sqrt{a_1}$, call the result $G(X;a_1)$, and assemble the odd powered terms in $\sqrt{a_1}$, call the result $-\sqrt{a_1}\hsx H(X;a_1)$ $-$then 
\[
F(X;a_1) \ = \ G(X;a_1) -\sqrt{a_1}\hsx H(X;a_1)
\]
and $G(X;a_1)$, $H(X;a_1)$ are polynomials with integral coefficients.
\vspace{0.2cm}

E.g.: \ When $n = 2$, 
\allowdisplaybreaks
\begin{align*}
F(X) \ 
&=\ 
(X - \sqrt{a_1} + \sqrt{a_2}) (X - \sqrt{a_1} - \sqrt{a_2})
\\[12pt]
&=\ 
(X - \sqrt{a_1})^2 - (\sqrt{a_2})^2
\\[12pt]
&=\ 
(X^2 + (\sqrt{a_1}\hsx)^2 - (\sqrt{a_2}\hsx)^2) - \sqrt{a_1} \hsx (2X).
\end{align*}
Now evaluate the data at $X = \Sigma$:
\[
0 \ = \ F(\Sigma; a_1) \ = \ G(\Sigma;a_1) -\sqrt{a_1}\hsx H(\Sigma;a_1)
\]
\qquad\qquad  $\implies$
\[
\sqrt{a_1} \ = \ \frac{G(\Sigma;a_1)}{H(\Sigma;a_1)} \in \Q
\]
provided $\tH(\Sigma; a_1) \neq 0$.  
To check that this is so, write
\allowdisplaybreaks
\begin{align*}
F(\Sigma; a_1) - F(\Sigma; -a_1)\ 
&=\ 
0 - F(\Sigma; -a_1)
\\[12pt]
&=\ 
(G(\Sigma;a_1) - \sqrt{a_1}\hsx H(\Sigma;a_1)) - (G(\Sigma;a_1) + \sqrt{a_1}\hsx H(\Sigma;a_1))
\\[12pt]
&=\ 
- 2\hsx \sqrt{a_1}\hsx H(\Sigma;a_1)
\end{align*}
\qquad  $\implies$
\allowdisplaybreaks
\begin{align*}
H(\Sigma;a_1) \ 
&=\ 
\frac{1}{2 \sqrt{a_1}} \hsx  F(\Sigma ; -a_1)
\\[12pt]
&=\ 
\frac{1}{2 \sqrt{a_1}} \hsx  \prod (\Sigma + \sqrt{a_1} \pm \sqrt{a_2} \pm \cdots \pm \sqrt{a_n}\hsx\big)
\\[12pt]
&=\ 
\frac{1}{2 \sqrt{a_1}} \hsx  
\prod \hsx 
\big(2\sqrt{a_1} + (\sqrt{a_2} \pm \sqrt{a_2}\hsx) + \cdots + (\sqrt{a_n} \pm \sqrt{a_n}\hsx\big) 
\\[12pt]
&=\ 
\frac{1}{2 \sqrt{a_1}} \hsx 
\prod\limits_{S \subset \{\sqrt{a_2}, \ldots, \sqrt{a_n}\hsx\}}\hsx 
\big(2\sqrt{a_1} + 2 \hsx\sum\limits_{a_i \in S} \hsx \sqrt{a_i}\hsx)
\\[12pt]
&=\ 
\frac{1}{\sqrt{a_1}} \hsx 
\prod\limits_{S \subset \{\sqrt{a_2}, \ldots, \sqrt{a_n}\hsx\}} \hsx 
\big(\sqrt{a_1} + \sum\limits_{a_i \in S} \hsx \sqrt{a_i}\hsx\big).
\end{align*}
But
\[
\sqrt{a_1} \ + \  \sum\limits_{a_i \in S} \hsx \sqrt{a_i}
\]
is never zero.
\vspace{0.5cm}

\begin{x}{\small\bf THEOREM} \ 
Given $x \in \R$, there are infinitely many coprime solutions $p$, $q$ $(q > 0)$  to
\[
\abs{x - \frac{p}{q}} \ \leq \ \frac{1}{q}.
\]
\end{x}
\vspace{0.3cm}

One can say more is $x$ is irrational.
\vspace{0.3cm}
\begin{x}{\small\bf THEOREM} \ 
Given $x \in \PP$, there are infinitely many coprime solutions $p$, $q$ $(q > 0)$  to
\[
\abs{x - \frac{p}{q}} \ \leq \ \frac{1}{q^2}.
\]
\vspace{0.2cm}

[Note: \ 
This estimate can be sharpened to
\[
\abs{x - \frac{p}{q}} \ \leq \ \frac{1}{\sqrt{5} \hsx q^2}
\]
but $\ds\frac{1}{\sqrt{5}}$ cannot be replaced by a smaller real number unless some restriction is placed on $x$.  
To see this, take 
\[
x \ = \ \frac{\sqrt{5} - 1}{2}.
\]
Then it can be shown that there is a coprime sequence $\ds\frac{p_n}{q_n}$ $(q_n > 0)$ with the property that if 
\ $0 < C < \ds\frac{1}{\sqrt{5}}$, \ then
\[
\abs{x - \frac{p_n}{q_n}} \ > \ \frac{C}{q_n^2} \qquad \forall \ n \gg 0.]
\]
\end{x}
\vspace{0.3cm}

\begin{x}{\small\bf NOTATION} \ 
For any real number $r$, write
\[
\{r\} \ = \ r - [r], 
\]
the 
\un{fractional part}
\index{fractional part} 
of $r$.
\vspace{0.2cm}

[Note: \ 
$0 \leq \{r\} < 1.$]
\end{x}
\vspace{0.3cm}

\begin{x}{\small\bf BOX PRINCIPLE} \ 
If $n + 1$ objects are placed in $n$ boxes, then some box contains at least 2 objects.
\end{x}
\vspace{0.3cm}

\begin{x}
\begin{spacing}{1.75}
{\small\bf CONSTRUCTION} 
\ 
Let $n > 1$ be a positive integer and divide the interval 
$[0,1]$ into $n$ subintervals $\ds\bigg[\frac{j}{n}, \frac{j+1}{n}\bigg]\vspace{1cm}$ $(j = 0, 1, \ldots, n - 1)$.  
Assuming that $x$ is irrational, the $n+1$ numbers $0, \{x\}, \ldots, \{n x\}$ are distinct elements of $[0,1]$, hence by the Box Principle at least 2 of them must be in one of the subintervals $\ds\bigg[\frac{j}{n}, \frac{j+1}{n}\bigg]$ $(j = 0, 1, \ldots, n - 1)$.  
Arrange matters in such a way that $\{j_1 x\}$ and $\{j_2 x\}$ $(j_2 > j_1)$ are contained in one subinterval of width $\ds\frac{1}{n}$.  
Set
\end{spacing}
\[
p = [j_2 x] - [j_1 x], \quad q = j_2 - j_1 \geq 1 \qquad (q < n).
\]
Then
\[
\abs{\{j_2 x\} - \{j_1 x\}} \ < \ \frac{1}{n}
\]
\qquad  $\implies$
\[
\abs{(j_2 - j_1) x - ([j_2 x] - [j_1 x])} \ < \ \frac{1}{n}
\]
\qquad \qquad $\implies$
\[
\abs{qx - p} \ < \ \frac{1}{n}
\]
\qquad \qquad $\implies$
\[
\abs{x - \frac{p}{q}} \ < \ \frac{1}{nq} \ < \ \frac{1}{q^2}.
\]

Existence per \#13 is thereby established.  To conclude, it has to be ruled out that there is just a finite number of coprime solutions to
\[
\abs{x - \frac{p}{q}} \ \leq \ \frac{1}{q^2},
\]
say
\[
\frac{p_1}{q_1}, \frac{p_2}{q_2}, \ldots, \frac{p_k}{q_k}.
\]
Since $x$ is irrational, there exists a positive integer $m > 1$ such that
\[
\abs{x - \frac{p_i}{q_i}} \ > \ \frac{1}{m} \qquad (i = 1, 2, \ldots, k).
\]
In \#16, replace $n$ by $m^2$ and $\ds\frac{p}{q}$ by $\ds\frac{a}{b}$, thus
\[
\abs{x - \frac{a}{b}} \ < \ \frac{1}{m^2 b} \ < \ \frac{1}{b^2}.
\]
On the other hand,
\[
\frac{1}{m^2 b} \ < \ \frac{1}{m} \qquad (b \geq 1),
\]
so
\[
\abs{x - \frac{a}{b}} \ < \ \frac{1}{m}.
\]
But
\[
\frac{a}{b} \ = \ \frac{p_i}{q_i} \qquad (\exists \ i)
\]
which implies that
\[
\abs{x - \frac{a}{b}} \ > \ \frac{1}{m}.
\]
Contradiction.
\end{x}
\vspace{0.3cm}

\begin{x}{\small\bf THEOREM} \ 
Given $x = \ds\frac{a}{b} \in \Q$ $(a, b \in \Z, \ b > 0, \ \gcd(a,b) = 1$), for any coprime pair $(p,q)$ $(q > 0)$ with 
\[
\frac{a}{b} \ \neq \ \frac{p}{q}
\]
there follows
\[
\abs{\frac{a}{b} - \frac{p}{q}} \ \geq \ \frac{1}{bq}.
\]
\vspace{0.2cm}

PROOF\ 
\allowdisplaybreaks
\begin{align*}
\frac{a}{b} \neq \frac{p}{q} 
&\implies
aq - bp \neq 0
\\[12pt]
&\implies
\abs{a q - b p} \geq 1
\end{align*}
\qquad \qquad $\implies$
\allowdisplaybreaks
\begin{align*}
\abs{\frac{a}{b} - \frac{p}{q}} \ 
&=\ 
\abs{\frac{aq - bp}{bq}}
\\[12pt]
&=\ 
\frac{\abs{aq - bp}}{\abs{bq}}
\\[12pt]
&=\ 
\frac{\abs{aq - bp}}{bq}
\\[12pt]
&\geq\  
\frac{1}{b q}.
\end{align*}
\end{x}
\vspace{0.3cm}

\begin{x}{\small\bf CRITERION} \ 
Let $x \in \R$.  Assume: \ There exists a coprime sequence $p_n, q_n$ $(q_n > 0)$ such that 
$x \neq \ds\frac{p_n}{q_n}$ for all $n$ and $q_n x - p_n \ra 0$ as $n \ra \infty$ 
$-$then $x$ is irrational.
\vspace{0.2cm}

[Suppose instead that $x$ is rational, say $x = \ds\frac{a}{b}$ $(b > 0$, $\gcd(a,b) = 1)$, thus
\allowdisplaybreaks
\begin{align*}
\frac{\abs{q_n x - p_n}}{q_n}\ 
&=\ 
\abs{x - \frac{p_n}{q_n}}
\\[12pt]
&=\ 
\abs{\frac{a}{b} - \frac{p_n}{q_n}}
\\[12pt]
&\geq\ 
\frac{1}{b q_n}
\end{align*}
\qquad \qquad $\implies$
\[\abs{q_n x - p_n} \ \geq \ \frac{1}{b} \ > \ 0.
\]
But this is a contradiction since $q_n x - p_n \ra 0$ by hypothesis.
\end{x}
\vspace{0.3cm}

\begin{x}{\small\bf CRITERION} \ 
Let $x \in \R$.  Fix positive constants \mC and $\delta$.  
Assume: \ There are infinitely many coprime solutions $p$, $q$ $(q > 0)$  to
\[
\abs{x - \frac{p}{q}} \ < \ \frac{C}{q^{1 + \delta}}.
\]
Then $x$ is irrational.
\vspace{0.2cm}

[The contrapositive is the assertion that for a rational $x$ there are but finitely many coprime $p$, $q$ $(q > 0)$ satisfying the stated inequality.  Take $x$ as $\ds\frac{a}{b}$ per \#17, hence
\allowdisplaybreaks
\begin{align*}
\frac{C}{q^{1 + \delta}} \ 
&> \ 
\abs{x - \frac{p}{q}}
\\[12pt]
&=\ 
\abs{\frac{a}{b} - \frac{p}{q}}
\\[12pt]
&\geq\ 
\frac{1}{bq}
\end{align*}
\qquad \qquad $\implies$
\[
\frac{C}{q^\delta} \  > \  \frac{1}{b} 
\implies 
(C b)^{1/\delta} \ > \ q.
\]
Accordingly, there are but finitely many possibilities for $q$.  The same is true of $p$.  To see this, fix $p$ and $q$ subject to
\[
\abs{\frac{a}{b} - \frac{p}{q}} \ < \ \frac{C}{q^{1 + \delta}}
\]
and consider fractions of the form
\[
\frac{p + r}{q} \qquad (r \in \Z),
\]
where
\[
\abs{\frac{a}{b} - \frac{p+r}{q}} \ < \ \frac{C}{q^{1 + \delta}}.
\]
Then
\allowdisplaybreaks
\begin{align*}
\frac{\abs{r}}{q} \ 
&=\ 
\abs{\frac{r}{q} + \frac{p}{q} -  \frac{a}{b} - \frac{p}{q} + \frac{a}{b}}
\\[12pt]
&\leq\ 
\abs{\frac{p+r}{q} - \frac{a}{b}} + \abs{\frac{p}{q} - \frac{a}{b}}
\\[12pt]
&<\ 
\frac{2 C}{q^{1 + \delta}}
\end{align*}
\qquad \qquad $\implies$
\[
\abs{r} \ < \ \frac{2 C}{q^\delta} \ \leq \ 2C.
\]
Our contention is therefore manifest.]
\end{x}
\vspace{0.3cm}

\begin{x}{\small\bf APPLICATION} \ 
Let $x \in \R$.  
Assume: \ There is a $\delta > 0$ and a sequence $\ds\frac{p_n}{q_n}$ $(q_n > 0) \neq x$ of rational numbers such that 
\[
\abs{x - \frac{p_n}{q_n}} \ = \ O\big(q_n^{-(1 + \delta)}\big).
\]
Then $x$ is irrational.
\end{x}
\vspace{0.5cm}

\[
\text{APPENDIX}
\]
\vspace{0.5cm}

{\small\bf IRRATIONALITY CRITERIA} \ 
Let $x$ be a real number $-$then the following conditions are equivalent.
\vspace{0.2cm}

\qquad (i) \ 
$x$ is irrational.
\vspace{0.2cm}

\qquad (ii) \ 
$\forall \ \epsilon > 0$, $\exists \ \ds\frac{p}{q} \in \Q$ such that
\[
0 \ < \ \abs{x - \frac{p}{q}} \ < \ \frac{\epsilon}{q}.
\]
\vspace{0.2cm}

\qquad (iii) \ 
\vspace{0.2cm}
$\forall$ real number $Q > 1$, $\exists$ an integer $q$ in the range $1 \leq q < Q$ and a rational integer $p$ such that 
\[
0 \ < \ \abs{x - \frac{p}{q}} \ < \ \frac{1}{q Q}.
\]

\qquad (iv) \ 
$\exists$ infinitely many $\ds\frac{p}{q} \in \Q$ such that 
\[
\abs{x - \frac{p}{q}} \ < \ \frac{1}{\sqrt{5} \hsx q^2}.
\]


%% file: _08_irrationality_of_e.tex
\chapter{
$\boldsymbol{\S}$\textbf{8}.\quad  IRRATIONALITY OF $\boldsymbol{e}$}
\setlength\parindent{2em}
\setcounter{theoremn}{0}
\renewcommand{\thepage}{\S8-\arabic{page}}

\ \indent 
Recall that $e$ can be defined as
\[
\sup \bigg\{ \sum\limits_{k=0}^n \hsx \frac{1}{k!} : n \in \N\bigg\}
\]
or, equivalently, as 
\[
\sup \bigg\{ \bigg(1 + \frac{1}{n}\bigg)^n : n \in \N\bigg\}.
\]
\vspace{0.3cm}

\begin{x}{\small\bf \un{N.B.}}  \ 
\[
\sum\limits_{k=0}^n \hsx \frac{1}{k!} 
\ < \ 
\sum\limits_{k=0}^{n+1} \hsx \frac{1}{k!}
\quad \text{and} \quad 
\bigg(1 + \frac{1}{n}\bigg)^n 
\ < \ 
\bigg(1 + \frac{1}{n+1}\bigg)^{n+1}.
\]
\end{x}
\vspace{0.3cm}

\begin{x}{\small\bf SUBLEMMA} \ 
Let $0 < r < 1$ $-$then 
\[
\sum\limits_{n=0}^\infty \hsx r^n \ = \ \frac{1}{1 - r},
\]
so
\[\sum\limits_{n=1}^\infty \hsx r^n \ = \ \frac{r}{1 - r}.
\]
\end{x}
\vspace{0.3cm}

\begin{x}{\small\bf THEOREM} \ 
$e$ is irrational.
\vspace{0.2cm}

PROOF \ 
Suppose that $e$ is rational, say $e = \ds\frac{x}{y}$, where $x$ and $y$ are positive integers and $\gcd(x,y) = 1$.  
Since $2 < e < 3$, $y$ is $> 1$.  
Write
\[
e 
\ = \ 
\bigg(1 + \frac{1}{1!} + \cdots + \frac{1}{y!}\bigg) + \cdots \ .
\]
Then
\allowdisplaybreaks
\begin{align*}
y! e \ 
&=\ 
y! \frac{x}{y}
\\[12pt]
&=\ 
(y - 1)! x
\\[12pt]
&=\ 
(y! + \frac{y!}{1!} + \cdots + \frac{y!}{y!}) + R.
\end{align*}
Here
\[
R \ = \ 
y! \hsx
\bigg(\frac{1}{(y + 1)!} + \frac{1}{(y + 2)!} + \cdots \bigg)
\]
is a positive integer.  Continuing, 
\allowdisplaybreaks
\begin{align*}
y! \hsx
\bigg(\frac{1}{(y + 1)!} + \frac{1}{(y + 2)!} + \cdots \bigg)\ 
&=\ 
\frac{1}{y + 1} + \frac{1}{(y+1)(y+2)} + \cdots
\\[12pt]
&<\ 
\frac{1}{y + 1} + \frac{1}{(y+1)^2} + \cdots 
\\[12pt]
&=\ 
\sum\limits_{n=1}^\infty \hsx \frac{1}{(y+1)^n}
\\[12pt]
&=\ 
\frac{\frac{1}{y+1}}{1 - \frac{1}{y+1}}
\\[12pt]
&=\
\frac{1}{y}
\\[12pt]
&<\
1. 
\end{align*}
But this implies that \mR is less than 1, a contradiction.
\vspace{0.2cm}

[Note: \ 
The preceding is actually an instance of \  \S7, \#18.  \ 
Thus take \  $q_n = n!$, \ $p_n = q_n \hsx \ds\sum\limits_{k=0}^n \hsx \ds\frac{1}{k!}$ $-$then

\allowdisplaybreaks
\begin{align*}
q_n e - p_n \ 
&=\ 
q_n \bigg(e -  \sum\limits_{k=0}^n \hsx \frac{1}{k!}\bigg)
\\[12pt]
&=\ 
n! \bigg( \sum\limits_{k=n+1}^\infty \hsx \frac{1}{k!} \bigg)
\\[12pt]
&=\ 
\frac{1}{n+1} + \frac{1}{(n+1)(n+2)} + \cdots
\\[12pt]
&< \ 
\frac{1}{n} \qquad \text{(cf. supra)}
\\[12pt]
&\ra 
0 \qquad (n \ra \infty).]
\end{align*}
\vspace{0.3cm}

The foregoing argument can be extended to establish the irrationality of $e^2$.
\vspace{0.3cm}

Thus start as before by assuming that $e^2 = \ds\frac{x}{y}$, where $x$ and $y$ are positive integers and $\gcd(x,y) = 1$ 
$(y > 1)$, hence
\[
y e \ = \ \frac{x}{e}
\]
\qquad\qquad $\implies$
\[
y \bigg( \sum\limits_{k=0}^\infty \hsx \frac{1}{k!}\bigg)
\ = \ 
x\bigg(\sum\limits_{k=0}^\infty \hsx (-1)^k \hsx \frac{1}{k!}\bigg)
\]
\qquad\qquad $\implies$ $(\forall \ n \in \N)$
\allowdisplaybreaks
\begin{align*}
y \bigg( \sum\limits_{k=0}^n \hsx \frac{1}{k!} + \sum\limits_{k>n} \hsx \frac{1}{k!}\bigg) \ 
&=\ 
x \bigg( \sum\limits_{k=0}^n \hsx (-1)^k \frac{1}{k!} + \sum\limits_{k>n} \hsx (-1)^k  \frac{1}{k!}\bigg)
\\[12pt]
\end{align*}
\qquad\qquad $\implies$
\[
y \bigg(A_n +  \sum\limits_{k>n} \hsx \frac{1}{k!}\bigg) 
\ = \ 
x \bigg(B_n + \sum\limits_{k>n} \hsx (-1)^k  \frac{1}{k!}\bigg)
\qquad \text{where} \ 
\begin{cases}
\ A_n  = \sum\limits_{k=0}^n \hsx \frac{1}{k!}\\[11pt]
\ B_n  = \sum\limits_{k=0}^n \hsx (-1)^k \frac{1}{k!}
\end{cases}
.
\]
Now multiply both sides of the last relation by $n!$ to get
\[
y\bigg(C_n + n! \hsx 
\sum\limits_{k > n} \hsx \frac{1}{k!}\bigg)
\ = \ 
x\bigg(D_n + n! \hsx  
\sum\limits_{k > n} \hsx (-1)^k\frac{1}{k!}\bigg),
\]
\[
\begin{cases}
\ C_n = n! A_n\\
\ D_n = n! B_n
\end{cases}
\]
being integers.  Moving on, 
\allowdisplaybreaks
\begin{align*}
&y C_n + 
y\bigg(\frac{1}{n+1} + \frac{1}{(n+1)(n+2)} + \cdots \bigg) \ 
\\[12pt]
&\hspace{1.25cm}=\ 
x D_n + x(-1)^{n+1} 
\bigg(\frac{1}{n+1} - \frac{1}{(n+1)(n+2)} + \cdots \bigg) 
\end{align*}
or still, 
\allowdisplaybreaks
\begin{align*}
&y C_n  - x D_n \
\\[12pt]
&\hspace{1cm}=\ 
x(-1)^{n+1} \hsx 
\bigg(\frac{1}{n+1} - \frac{1}{(n+1)(n+2)} + \cdots \bigg) 
- 
y \bigg(\frac{1}{n+1} + \frac{1}{(n+1)(n+2)} + \cdots \bigg).
\end{align*}
Therefore
\allowdisplaybreaks
\begin{align*}
\abs{y C_n  - x D_n} \ 
&\leq\ 
x  
\abs{\frac{1}{n+1} - \frac{1}{(n+1)(n+2)} + \cdots} 
+ 
y 
\abs{\frac{1}{n+1} + \frac{1}{(n+1)(n+2)} + \cdots}
\\[12pt]
&\leq\
x\bigg(
\frac{1}{n+1} + \frac{1}{(n+1)(n+2)} + \cdots
\bigg)
+
y\bigg(
\frac{1}{n+1} + \frac{1}{(n+1)(n+2)} + \cdots
\bigg)
\\[12pt] 
&<\ 
x \frac{1}{n} + y \frac{1}{n} 
\\[12pt] 
&=\
\frac{x + y}{n}.
\end{align*}
Finally, for all $n \gg 0$, 
\[
\frac{x + y}{n} \ <\ 1.
\]
I.e.: \ For an infinite set of $n$, 
\[
\abs{y C_n - x D_n} \ = \ 0,
\]
or still, for an infinite set of $n$, 
\[
y C_n \ = \ x D_n,
\]
an impossibility.
\end{x}
\vspace{0.3cm}

\begin{x}{\small\bf DEFINITION} \ 
An irrational number $r$ is a 
\un{quadratic irrational}
\index{quadratic irrational} 
if there exist integers \mA, \mB, \mC not all zero such that
\[
A r^2 + Br + C = 0.
\]

[Note: \ 
A quadratic irrational is necessarily algebraic.]
\end{x}
\vspace{0.3cm}

\begin{x}{\small\bf EXAMPLE} \ 
$\sqrt{2}$ is a quadratic irrational.
\end{x}
\vspace{0.3cm}

\begin{x}{\small\bf THEOREM} \ 
$e$ is not a quadratic irrational.
\end{x}
\vspace{0.3cm}

The proof is detailed in the lines below.
\vspace{0.3cm}

To arrive at a contradiction, suppose that there are integers \mA, \mB, \mC not all zero such that
\[
A e^2 + Be + C = 0.
\]
\vspace{0.01cm}

\begin{x}{\small\bf \un{N.B.}} \ 
If $A = 0$, matters are clear.  
If $A \neq 0$ and if $B = 0$, matters are clear.  
If $A \neq 0$ and if $B \neq 0$ and if $C = 0$, matters are clear.  
One can accordingly assume from the beginning that $A \neq 0$, $B \neq 0$, $C \neq 0$.  
Moreover, we shall work instead with the equation
\[
A e + B + \frac{C}{e} = 0.
\]
\end{x}
\vspace{0.3cm}

\begin{x}{\small\bf SUBLEMMA} \ 
Given $n \in \N$, there is an integer $I_n$ such that 
\[
n! e 
\ = \ 
I_n + \frac{1}{n + \alpha_n},
\]
where $0 < \alpha_n < 1.$
\vspace{0.2cm}

PROOF \ 
Write
\[
n! \hsx e 
\ = \
\sum\limits_{k=0}^n \hsx 
\frac{n!}{k!} 
+ 
\sum\limits_{k=n+1}^\infty \hsx 
\frac{n!}{k!}.
\]
\vspace{0.01cm}

\allowdisplaybreaks
\begin{align*}
\text{\textbullet} \quad 
\sum\limits_{k=n+1}^\infty \hsx \frac{n!}{k!} 
&=\ 
\frac{n!}{(n+1)!} + \frac{n!}{(n+2)!} + \cdots \hspace{5cm}
\\[12pt]
&>\  
\frac{n!}{(n+1)!} 
\\[12pt]
&=\ 
\frac{1}{n + 1}.
\\[19pt]
\text{\textbullet} \quad 
\sum\limits_{k=n+1}^\infty \hsx \frac{n!}{k!} 
&=\ 
\frac{1}{n + 1} + \frac{1}{(n+1)(n+2)} + \cdots
\\[12pt]
&<\ 
\frac{1}{n+1} + \frac{1}{(n+1)^2} + \cdots
\\[12pt]
&=\ 
\frac{1}{n}.
\end{align*}
\vspace{0.2cm}

Therefore
\[
\frac{1}{n + 1}
\ < \ 
\sum\limits_{k=n+1}^\infty \hsx \frac{n!}{k!} 
\ < \ 
\frac{1}{n},
\]
from which
\[
\sum\limits_{k=n+1}^\infty \hsx \frac{n!}{k!} 
\ = \ 
\frac{1}{n + \alpha_n} 
\qquad (0 < \alpha_n < 1).
\]
To conclude, it remains only to set
\[
I_n 
\ = \ 
\sum\limits_{k=0}^n \hsx \frac{n!}{k!}.
\]
\end{x}
\vspace{0.3cm}

\begin{x}{\small\bf SUBLEMMA} \ 
Given $n \in \N$, there is an integer $J_n$ such that
\[
\frac{n!}{e} 
\ = \ J_n + \frac{(-1)^{n + 1}}{n + 1 + \beta_n},
\]
where $0 < \beta_n < 1$.
\vspace{0.2cm}

PROOF \ 
Write
\[
\frac{n!}{e} 
\ = \ 
\sum\limits_{k=0}^n \hsx (-1)^k \frac{n!}{k!} 
+ 
\sum\limits_{k=n+1}^\infty \hsx (-1)^k \frac{n!}{k!}.
\]
\vspace{0.0cm}

\allowdisplaybreaks
\begin{align*}
\text{\textbullet} \quad 
\sum\limits_{k=n+1}^\infty \hsx (-1)^k \frac{n!}{k!} \ 
&=\ 
\sum\limits_{\ell = 0}^\infty \hsx (-1)^{\ell + (n + 1)} \frac{n!}{(\ell + (n + 1))!} \hspace{3cm}
\\[12pt]
&=\ 
(-1)^{n+1} \hsx
\sum\limits_{\ell = 0}^\infty \hsx (-1)^{\ell} \frac{n!}{(\ell + (n + 1))!}
\\[12pt]
&\equiv\ (-1)^{n + 1} S.
\end{align*}

Put
\[
S_N \ = \ 
\sum\limits_{\ell = 0}^N \hsx 
(-1)^\ell \frac{n!}{(\ell + (n + 1))!}.
\]
Then
\[
S_N  \ <  \ S  \ <  \ S_{N+1} \qquad \text{($N$ odd)}.
\]
In particular $(N = 1)$:
\[
\frac{1}{n+1} - \frac{1}{(n+1)(n+2)}
\ < \ 
S
\ < \ 
\frac{1}{n+1} - \frac{1}{(n+1)(n+2)} + \frac{1}{(n+1)(n+2)(n+3)}.
\]
\vspace{0.0cm}
\allowdisplaybreaks
\begin{align*}
\text{\textbullet} \quad 
\frac{1}{n+1} - \frac{1}{(n+1)(n+2)}\ 
&=\ 
\frac{1}{n+1} \bigg(1 - \frac{1}{n+2}\bigg) \hspace{3.5cm}
\\[12pt]
&=\ 
\frac{1}{n+2}
\end{align*}
and
\allowdisplaybreaks
\begin{align*}
\frac{1}{n+1} \ - \ &\frac{1}{(n+1)(n+2)}  + \frac{1}{(n+1)(n+2)(n+3)} 
\\[12pt]
&\hspace{1cm}=\ 
\frac{1}{n+1} \hsx
\bigg(
1 - \frac{1}{n+2} +\frac{1}{(n+2)(n+3)}
\bigg) \hspace{2cm}
\\[12pt]
&\hspace{1cm}=\ 
\frac{1}{n+1} \hsx
\bigg(
1 + \frac{1}{n+2} \bigg(-1 + \frac{1}{n+3}\bigg)
\bigg)
\\[12pt]
&\hspace{1cm}=\ 
\frac{1}{n+1}\bigg(1 + \frac{1}{n+2} \bigg(\frac{-n - 3 + 1}{n+3}\bigg) \bigg) 
\\[12pt]
&\hspace{1cm}=\ 
\frac{1}{n+1}\bigg(1 + \frac{1}{n+2} \bigg(\frac{-n - 2}{n+3}\bigg) \bigg) 
\\[12pt]
&\hspace{1cm}=\ 
\frac{1}{n+1}\bigg(1 - \frac{1}{n+3}\bigg) 
\\[12pt]
&\hspace{1cm}<\ 
\frac{1}{n+1}.
\end{align*}
Therefore
\[
\frac{1}{n+2} < S < \frac{1}{n+1} 
\implies
S = \frac{1}{n + 1 + \beta_n} \qquad (0 < \beta_n < 1).
\]
And then
\[
\sum\limits_{k=n+1}^\infty \hsx 
(-1)^k \hsx \frac{n!}{k!} 
\ = \ 
(-1)^{n+1} \hsx S
\ = \ 
\frac{(-1)^{n+1} }{n + 1 + \beta_n}.
\]
To conclude, let
\[
J_n 
\ = \ 
\sum\limits_{k=0}^n \hsx (-1)^k \hsx \frac{n!}{k!} .
\]

Summary:
\[
\begin{cases}
\ n! e - I_n \ = \ O\bigg(\ds\frac{1}{n}\bigg)\\[15pt]
\ \ds\frac{n!}{e} - J_n \ = \ O\bigg(\frac{1}{n}\bigg) 
\end{cases}
.
\]
Return now to the equation
\[
A e + B + \frac{C}{e} = 0
\]
and consider
\allowdisplaybreaks
\begin{align*}
A (n! e - I_n) + C\bigl(\frac{n!}{e} - J_n\bigr)\ 
&=\ 
n! \bigl(A e + B  + \frac{C}{e}\bigr) -
\bigl( A I_n + Bn! + C J_n \bigr)
\\[12pt]
&=\
-\bigl(A I_n + B n! +  C J_n\bigl)
\\[12pt] 
&\equiv\  -K_n.
\end{align*}
Then $K_n$ is an integer.  But
\[
K_n 
\ = \ O\bigg(\frac{1}{n}\bigg).
\]
Therefore
\[
K_n \ = \ 0 \qquad (n \gg 0).
\]
\end{x}
\vspace{0.3cm}

\begin{x}{\small\bf SUBLEMMA} \ 
\[
K_{n+2} - (n + 1) (K_n + K_{n+1}) \ = \ 2 A.
\]
\vspace{0.2cm}

[
Use the relations
\[
\begin{cases}
\ I_{n+1} \ = \ 1 + (n + 1) I_n\\
\ J_{n+1} \ = \ (-1)^{n+1} + (n + 1)J_n
\end{cases}
.]
\]
\end{x}
\vspace{0.3cm}

Since $A \neq 0$, the relation figuring in \#10 is impossible for $n \gg 0$.  
And this contradiction closes out the proof of \#6.
\vspace{0.3cm}

\begin{x}{\small\bf SCHOLIUM} \ 
1, $e$, $e^2$ are linearly independent over $\Q$.
\end{x}
\vspace{0.5cm}

\[
\text{APPENDIX}
\]
\vspace{0.5cm}

\qquad {\small\bf EXAMPLE} 1 \ 
Suppose that $r$ is a nonzero rational number $-$then the number
\[
\sum\limits_{k=0}^\infty  \ 
\frac{r^k}{2^{k(k-1)/2}}
\ = \ 
1 + r + \frac{1}{2} r^2 + \frac{1}{8} r^3 + \cdots
\]
is irrational.
\vspace{0.3cm}

\qquad {\small\bf EXAMPLE} 2 \ 
Suppose that $r$ is a nonzero rational number subject to $0 < \abs{r} < 1$ $-$then the number
\[
\sum\limits_{k=0}^\infty \hsx 
r^{2^k}
\ = \ 
r + r^2 + r^4 + r^8 + \cdots
\]
is irrational.
\vspace{0.3cm}

\qquad {\small\bf EXAMPLE} 3 \ 
Suppose that \mM is an integer $\geq 2$ $-$then the number
\[
\sum\limits_{k=1}^\infty \hsx 
\frac{1}{M^{k^2}}
\]
is irrational.
\vspace{0.3cm}


%% file: _09_irrationality_of_e_a_div_b.tex
\chapter{
$\boldsymbol{\S}$\textbf{9}.\quad  IRRATIONALITY OF $\boldsymbol{e^{a/b}}$}
\setlength\parindent{2em}
\setcounter{theoremn}{0}
\renewcommand{\thepage}{\S9-\arabic{page}}

\ \indent 
Let $a/b$ be a nonzero rational number.
\vspace{0.3cm}

\begin{x}{\small\bf THEOREM} \ 
$e^{a/b}$ is irrational.
\vspace{0.2cm}

[Note: \ 
Special cases, namely $e$ and $e^2$ are irrational, as has been shown in \S8.]
\end{x}
\vspace{0.3cm}

\begin{x}{\small\bf LEMMA}  \ 
If $e^r$ is irrational for all integers $r \geq 1$, then $e^{a/b}$ is irrational for all nonzero rationals $a/b$.
\vspace{0.2cm}

PROOF \ 
Take $a \in \N$ and suppose that $e^{a/b}$ is rational, say $e^{a/b} = q \in \Q$ $-$then 
\[
e^a \ = \ \big(e^{a/b}\big)^b \ = \ q^b \in \Q.
\]

Working toward a contradiction, assume that for some $r \in \N$, $e^r$ is rational and choose a positive integer $m$ with the property that $me^r \in \N$.
\end{x}
\vspace{0.3cm}

The data in place, we shall now introduce the machinery that will be utilized to arrive at our objective.
\vspace{0.3cm}

\begin{x}{\small\bf NOTATION} \ 
Given $n \in \N$, let
\[
T_n(X) \ = \ \prod\limits_{j=n+1}^{2n} \hsx (X - j),
\]
an element of $\Z[X]$.
\end{x}
\vspace{0.3cm}

\begin{x}{\small\bf RAPPEL} \ 
\[
e^x \ = \ \sum\limits_{k=0}^\infty \hsx \frac{x^k}{k!}.
\]
\end{x}
\vspace{0.3cm}

Put
\[
\delta \ = \ x \hsx \frac{d}{dx}.
\]
\vspace{0.3cm}

\begin{x}{\small\bf SUBLEMMA} \ 
\[
T_n(\delta) x^k \ = \ T_n(k) x^k.
\]
\end{x}
\vspace{0.3cm}

\begin{x}{\small\bf LEMMA} \ 
\allowdisplaybreaks
\begin{align*}
T_n(\delta) e^x \ 
&=\ 
Q_n(x) e^x 
\\[12pt]
&=\
(x^n + \cdots) e^x
\\[12pt] 
&=\ 
P_n(x) + R_n(x),
\end{align*}
where
\[
P_n(x) 
\ = \ 
\sum\limits_{k=0}^n \hsx T_n(k) \frac{x^k}{k!} 
\ = \ 
(-1)^n \hsx \sum\limits_{k=0}^n \hsx \frac{(2n - k)!}{n!} \binom{n}{k} x^k
\]
and
\allowdisplaybreaks
\begin{align*}
R_n(x) \ 
&=\ 
\sum\limits_{k=n+1}^\infty \hsx T_n(k) \hsx \frac{x^k}{k!}
\\[12pt] 
&=\ 
\sum\limits_{k=2n+1}^\infty \hsx T_n(k) \hsx \frac{x^k}{k!}
\\[12pt] 
&=\ 
\sum\limits_{k=2n+1}^\infty \hsx \frac{(k - n - 1)!}{(k - 2n - 1)!} \frac{x^k}{k!}.
\end{align*}

\end{x}
\vspace{0.3cm}

\begin{x}{\small\bf \un{N.B.}} \ 
\[
\begin{cases}
\ Q_n(x) \in \Z[x]\\[3pt]
\ P_n(x) \in \Z[x]
\end{cases}
.
\]
Accordingly, at an $r \in \N$, 
\[
\begin{cases}
\ Q_n(r) \in \Z\\[3pt]
\ P_n(r) \in \Z
\end{cases}
.
\]
\end{x}
\vspace{0.3cm}

\begin{x}{\small\bf REMAINDER ESTIMATE} \ 
\allowdisplaybreaks
\begin{align*}
\abs{R_n(x)} \ 
&\leq\ 
\frac{n!}{(2n + 1)!} \hsx 
\sum\limits_{k=2n+1}^\infty \hsx 
\frac{\abs{x}^k}{(k - 2n - 1)!} 
\\[12pt]
&=\ 
\frac{n! \abs{x}^{2n + 1}}{(2n + 1)!} \hsx e^{\abs{x}}.
\end{align*}
\end{x}
\vspace{0.3cm}

Returning to the situation above, we claim that for sufficiently large $n$, 
\[
0 \ < \ m R_n(r) \ < \ 1.
\]
To see this, consider
\[
\frac{n! r^{2n+1}}{(2n + 1)!} \hsx e^r 
\ = \ 
\frac{n!}{(2n + 1)!} r^{2n} (re^r).
\]
Then
\allowdisplaybreaks
\begin{align*}
\frac{n!}{(2n + 1)!} r^{2n}\ 
&=\ 
\frac{n!}{n!} \cdot \frac{r^2}{n + 1} \cdot \frac{r^2}{n + 2} \cdots \frac{r^2}{n + n} \cdot \frac{1}{2n + 1}
\\[12pt]
&=\ 
\frac{r^2}{n + 1} \cdot \frac{r^2}{n + 2} \cdots \frac{r^2}{n + n} \cdot \frac{1}{2n + 1}.
\end{align*}
Choose $n \gg 0$:
\[
\frac{r^2}{n + 1} \ < \ 
1,
\]
thus
\[
\frac{n!}{(2n+1)!} \hsx r^{2n} \ < \ \frac{1}{2n + 1},
\]
from which the claim is immediate.

On the other hand,
\allowdisplaybreaks
\begin{align*}
m R_n (r) \ 
&=\ 
m(Q_n(r) e^r - P_n(r))
\\[12pt]
&=\ 
(m e^r) Q_n(r) - mP_n(r)
\\[12pt]
& 
\in \Z.
\end{align*}
\vspace{0.3cm}
But there are no integers between 0 and 1.

\begin{x}{\small\bf REMARK} \ 
It will be shown in due course that if $x\neq 0$ is algebraic, then $e^x$ is irrational, so e.g., $e^{\sqrt{2}}$ is irrational.
\end{x}
\vspace{0.5cm}

\[
\text{APPENDIX}
\]
\vspace{0.25cm}

\qquad \un{$0 \leq k \leq n$:} \ \ 
Here
\allowdisplaybreaks
\begin{align*}
(-1)^n \hsx \frac{(2n - k)!}{n!} \hsx \binom{n}{k} \ 
&=\ 
(-1)^n \hsx \frac{(2n - k)!}{n!} \hsx \frac{n!}{k! (n - k)!}
\\[12pt]
&=\ 
(-1)^n \hsx \frac{(2n - k)!}{(n - k)!} \hsx \frac{1}{k!}
\end{align*}
and the claim is that
\[
T_n(k) \ = \ (-1)^n \hsx \frac{(2n - k)!}{(n - k)!}.
\]
\vspace{0.2cm}

[
\textbullet \ $k = 0$:
\allowdisplaybreaks
\begin{align*}
T_n(0) \ 
&=\ 
\prod\limits_{j=n+1}^{2n} \hsx (0 - j)
\\[12pt]
&=\ 
- (n + 1) ( - (n + 2)) \cdots (- (2n))
\\[12pt]
&=\ 
(-1)^n (n + 1) (n + 2) \cdots (2n)
\\[12pt]
&=\ 
(-1)^n \hsx \frac{2n!}{n!}.
\end{align*}

\textbullet \ $k = 1$:
\allowdisplaybreaks
\begin{align*}
T_n(1) \ 
&=\ 
\prod\limits_{j=n+1}^{2n} \hsx (1 - j)
\\[12pt]
&=\ 
(1 - (n + 1)) (1 - (n + 2)) \cdots (1 - (2n))
\\[12pt]
&=\ 
(-n) \hsx (- n - 1) \cdots (-(2n -1))
\\[12pt]
&=\ 
(-1)^n \hsx(n) (n + 1) \cdots (2n - 1)
\\[12pt]
&=\ 
(-1)^n \hsx \frac{(2n - 1)!}{(n - 1)!}
\end{align*}

\vspace{0.3cm}
\[
\bcdot\bcdot\bcdot\bcdot\bcdot\bcdot\bcdot\bcdot\bcdot\bcdot\bcdot\bcdot\bcdot
\bcdot\bcdot\bcdot\bcdot\bcdot\bcdot\bcdot\bcdot\bcdot\bcdot\bcdot\bcdot\bcdot
\]
\vspace{0.3cm}

\textbullet \ $k = n$:
\allowdisplaybreaks
\begin{align*}
T_n(n) \ 
&=\ 
\prod\limits_{j=n+1}^{2n} \hsx (n - j)
\\[12pt]
&=\ 
(n - (n + 1)) (n - (n + 2)) \cdots (n - (2n))
\\[12pt]
&=\ 
(-1) (-2) \cdots (-n)
\\[12pt]
&=\ 
(-1)^n \hsx n!
\\[12pt]
&=\ 
(-1)^n \hsx \frac{(2n - n)!}{(n - n)!}.]
\end{align*}
\vspace{0.3cm}

\qquad \un{$2n + 1 \leq k < \infty$:} \ \ 
In
 this situation, the claim is that
\[
T_n(k) \ = \ \frac{(k - n - 1)!}{(k - 2n - 1)!}.
\]

[
\textbullet \ $k = 2n+1$:

\allowdisplaybreaks
\begin{align*}
T_n(2n + 1) \ 
&=\ 
\prod\limits_{j=n+1}^{2n} \hsx (2n + 1 - j)
\\[12pt]
&=\ 
(2n + 1 - (n + 1)) (2n + 1 - (n + 2)) \cdots (2n + 1 - 2n)
\\[12pt]
&=\ 
(n) (n -1) \cdots (1)
\\[12pt]
&=\ 
n!
\\[12pt]
&=\ 
\frac{(2n + 1 - n - 1)!}{(2n + 1 - 2n - 1)!}.
\end{align*}

\textbullet \ $k = 2n + 2$:

\allowdisplaybreaks
\begin{align*}
T_n(2n + 2) \ 
&=\ 
\prod\limits_{j=n+1}^{2n} \hsx (2n + 2 - j)
\\[12pt]
&=\ 
(2n + 2 - (n + 1)) (2n + 2 - (n + 2)) \cdots (2n + 2 - 2n)
\\[12pt]
&=\ 
(n + 1) (n) \cdots (2)
\\[12pt]
&=\ 
(n+1)!
\\[12pt]
&=\ 
\frac{(2n + 2 - n - 1)!}{(2n + 2 - 2n - 1)!}.  
\end{align*}

\vspace{0.3cm}
\[
\bcdot\bcdot\bcdot\bcdot\bcdot\bcdot\bcdot\bcdot\bcdot\bcdot\bcdot\bcdot\bcdot
\bcdot\bcdot\bcdot\bcdot\bcdot\bcdot\bcdot\bcdot\bcdot\bcdot\bcdot\bcdot\bcdot
\]
\vspace{0.3cm}

To prove the remainder estimate, one has to show that 
\[
\frac{(k - n - 1)!}{k!} \ \leq \ \frac{n!}{(2n + 1)!} \qquad (k \geq 2n + 1).
\]
Let $k = 2n + r$ $(r = 1, 2, \ldots)$ and take $r > 1$ $-$then 
\allowdisplaybreaks
\begin{align*}
\frac{(k - n - 1)!}{k!} \ 
&=\
\frac{(2n + r - n - 1)!}{(2n + r)!}
\\[12pt]
&=\
\frac{(n + r - 1)!}{(2n + r)!}
\\[12pt]
&=\
\frac{(n + r - 1)!}{(2n + 1)! (2n + 2) \cdots (2n + r)}.
\end{align*}
Cancelling the
\[
\frac{1}{(2n + 1)!},
\]
there remains the claim that
\[
\frac{(n + r - 1)!}{(2n + 2) \cdots (2n + r)} \ \leq \ n!.
\]
Write
\allowdisplaybreaks
\begin{align*}
(n + r - 1)! \ 
&=\ 
1 \cdot 2 \cdots (n - 1) (n + 1 - 1) (n + 2 - 1) \cdots (n + r - 1)
\\[12pt]
&=\ 
(n -1)! (n + 1 - 1) (n + 2 - 1) \cdots (n + r - 1).
\end{align*}
Cancelling the $(n-1)!$, matters thus reduce to
\[
\frac{(n + 1 - 1) (n + 2 - 1) \cdots (n + r - 1)}{(2n + 2) \cdots (2n + r)} \ \leq \ n
\]
or still, 
\[
\frac{ (n + 2 - 1) \cdots (n + r - 1)}{(2n + 2) \cdots (2n + r)} \ \leq \ 1,
\]
which is obvious.

%% file: _10_irrationality_of_e_a_div_b_bis.tex
\chapter{
$\boldsymbol{\S}$\textbf{10}.\quad  IRRATIONALITY OF $\boldsymbol{e^{a/b}}$ (bis)}
\setlength\parindent{2em}
\setcounter{theoremn}{0}
\renewcommand{\thepage}{\S10-\arabic{page}}

\ \indent 
There is another way to prove that $e^{a/b}$ is irrational ($a/b$ a nonzero rational number).  
Thus, proceeding as in \S9, suppose that for some $r \in \N$, $e^r$ is rational, say $e^r = \ds\frac{u}{v}$ $(u, v \in \Z, \ v > 0)$.

Let
\[
f(x) \ = \ \frac{x^n (1 - x)^n}{n!}.
\]
Then
\[
0 < x < 1 \implies 0 < f(x) < \frac{1}{n!}.
\]
\vspace{0.3cm}

\begin{x}{\small\bf LEMMA} \ 
\[
f^{(j)} (0) \in \Z \qquad (j = 1, 2, \ldots).
\]
\end{x}
\vspace{0.3cm}

\begin{x}{\small\bf \un{N.B.}} \ 
\[
f^{(j)} (1) \in \Z \qquad (j = 1, 2, \ldots).
\]
\vspace{0.01cm}

[This is because
\[
f(1 - x) \ = \ f(x).]
\]
\end{x}
\vspace{0.3cm}

Given $n \in \N$, put
\[
F(x) \ = \ 
r^{2n} f(x) - r^{2n-1} f^\prime(x) + r^{2n-2} f^{\prime\prime}(x) - \cdots - r f^{(2n-1)} (x) + f^{(2n)} (x), 
\]
and note that
\[
F(0), \ F(1) \in \Z.
\]

Obviously
\[
\frac{d}{dx} \big(e^{rx} F(x)\big) \ = \ 
e^{rx} (r F(x) + F^\prime (x)) 
\ = \ 
r^{2n + 1} e^{rx} f(x)
\]
\qquad\qquad $\implies$
\allowdisplaybreaks
\begin{align*}
vr^{2n+1} \hsx 
\int\limits_0^1 \hsx e^{rx} f(x) dx \ 
&=\ 
v \big(e^{r x} F(x) \bigg|_0^1 
\\[12pt]
&=\ 
v (e^r F(1)) - v F(0)
\\[12pt]
&=\ 
u F(1) - v F(0),
\end{align*}
an integer.  
On the other hand, 
\allowdisplaybreaks
\begin{align*}
0\ 
&<\ 
vr^{2n+1} \hsx 
\int\limits_0^1 \hsx e^{rx} f(x) dx
\\[12pt]
&<\ 
\frac{v r^{2n+1} e^r }{n!}
\\[12pt]
&=\ 
v r e^r \frac{(r^2)^n}{n!}
\\[12pt]
&<\  
1
\end{align*}
for $n \gg 0$ (cf. \S0), giving a contradiction.
\vspace{0.3cm}

This is a good place to insert an application.
\vspace{0.3cm}

\begin{x}{\small\bf DEFINITION} \ 
The 
\un{natural logarithm}
\index{natural logarithm} 
is $\log_e$.
\end{x}
\vspace{0.3cm}

\begin{x}{\small\bf NOTATION} \ 
Write $\elln$ in place of $\log_e$.
\end{x}
\vspace{0.3cm}

\begin{x}{\small\bf THEOREM} \ 
If $q \neq 1$ is rational and positive, then $\elln(q)$ is irrational.
\vspace{0.2cm}

PROOF \ 
Suppose that $\elln(q)$ is rational $-$then $e^{\elln(q)}$ is irrational.  
Meanwhile
\[
q \ = \ e^{\elln(q)}.
\]
\end{x}
\vspace{0.3cm}

\begin{x}{\small\bf SCHOLIUM} \ 
If $x \neq 1$ is a positive real number and if $\elln(x)$ is rational, then $x$ is irrational.
\end{x}
\vspace{0.3cm}

\[
\text{APPENDIX}
\]
\vspace{0.3cm}

Let $a \neq 1$, $b \neq 1$ be positive real numbers $-$then 
\allowdisplaybreaks
\begin{align*}
\log_a (b) \hsx \elln(a) \ 
&=\ 
\elln\big(a^{\log_a (b)} \big)
\\[12pt]
&=\ 
\elln(b),
\end{align*}
so
\[
\log_a (b) \ = \ \frac{\elln (b)}{\elln (a)}.
\]
\vspace{0.3cm}

{\small\bf  EXAMPLE}
\[
\log_3 9 
\ = \ 
\frac{\elln(9)}{\elln(3)}
\ = \ 
\frac{\elln(3^2)}{\elln(3)}
\ = \ 
2 \hsx \frac{\elln(3)}{\elln(3)}
\ = \ 
2.
\]


%% file: _11_irrationality_of_pi.tex
\chapter{
$\boldsymbol{\S}$\textbf{11}.\quad  IRRATIONALITY OF $\boldsymbol{\pi}$}
\setlength\parindent{2em}
\setcounter{theoremn}{0}
\renewcommand{\thepage}{\S11-\arabic{page}}

\ \indent 
There are many ways to introduce the number $\pi$.
\vspace{0.3cm}

\begin{x}{\small\bf DEFINITION} \ 
Geometrically, $\pi$ is the length of a semicircle of radius one, i.e., analytically, 
\[
\pi 
\ = \ 
\int\limits_{-1}^1 \hsx
\frac{dx}{\sqrt{1 - x^2}}.
\]
\end{x}
\vspace{0.3cm}

\begin{x}{\small\bf THEOREM} \ 
Consider the complex exponential function 
\[
\exp : \C \ra \C.
\]
Then $\pi$ is the unique positive real number with the property that
\[
\Ker (\exp) 
\ = \ 
2 \pi \hsx \sqrt{-1} \  \Z.
\]
\end{x}
\vspace{0.3cm}

\begin{x}{\small\bf THEOREM} \ 
$\pi$ is the unique positive real number such that $\ds\cos \hsx \frac{\pi}{2} = 0$ and $\cos x \neq 0$ for 
$0 \leq x < \ds\frac{\pi}{2}$.
\end{x}
\vspace{0.3cm}

\begin{x}{\small\bf THEOREM} \ 
$\pi$ is irrational.
\end{x}
\vspace{0.3cm}

We shall give four proofs of this result.
\vspace{0.5cm}

\un{First Proof:} \ 
Suppose that $\pi = \ds\frac{a}{b}$, where $a$ and $b$ are positive integers.   Introduce
\[
f(x) 
\ = \ 
\frac{x^n (a - b x)^n}{n!}
\]
and 
\[
F(x) 
\ = \ 
f(x) - f^{(2)}(x) + f^{(4)}(x)  - \cdots + (-1)^n f^{(2n)}(x),
\]
$n \in \N$ to be determined momentarily.  \ 
Note that 
$f^{(j)}(0) \in \Z$ \ $(j = 1, 2, \ldots)$, \ hence 
$f^{(j)}(\pi) \in \Z$ $(j = 1, 2, \ldots)$ 
(since $f(x) = f\bigg(\ds\frac{a}{b} - x\bigg) = f(\pi - x)$).  
Next
\allowdisplaybreaks
\begin{align*}
\frac{d}{dx} \bigl( F^\prime(x) \sin x - F(x) \cos x \bigr)\ 
&=\ 
F^{\prime\prime}(x)  \sin x + F(x) \sin x
\\[12pt]
&=\ 
f(x) \sin x 
\qquad \text{(since $F(x) + F^{\prime\prime}(x) = f(x)$)}.
\end{align*}
Therefore
\allowdisplaybreaks
\begin{align*}
\int\limits_0^\pi \hsx f(x) \sin x dx \ 
&=\ 
\big(
F^\prime(x) \sin x - F(x) \cos x \hsx \bigg|_0^\pi
\\[12pt]
&=\ 
F(\pi) + F(0).
\end{align*}
But $F(\pi) + F(0)$ is an integer.  On the other hand, 
\[
0 < f(x) \sin x < \frac{\pi^n a^n}{n!} \qquad (0 \leq x \leq \pi),
\]
so
\[
\int\limits_0^\pi \hsx f(x) \sin x dx \ < \ \pi \hsx \frac{\pi^n a^n}{n!}
\]
is positive and tends to zero as $n \ra \infty$ (cf. \S0).
\vspace{0.3cm}

\un{Second Proof:} \ 
This proof is a slightly more complicated variant of the preceding proof and has the merit that it establishes the stronger result that $\pi^2$ is irrational.  
Proceeding to the details, suppose that
$\pi^2 = \ds\frac{a}{b}$, where $a$ and $b$ are positive integers but this time introduce
\[
f(x) 
\ = \ 
\frac{x^n (1 - x)^n}{n!},
\]
a polynomial encountered earlier (cf. \S10).  
Put
\[
F(x) 
\ = \ 
b^n
(\pi^{2n} f(x) 
- 
\pi^{2n - 2} f^{(2)} (x) 
+ 
\pi^{2n - 4} f^{(4)} (x) 
- \cdots + 
(-1)^n f^{(2n)} (x))
\]
and note that
\[
F(0), \ F(1) \in \Z.
\]
Moreover
\allowdisplaybreaks
\begin{align*}
\frac{d}{dx} (F^\prime(x) \sin(\pi x) - \pi \hsx F(x) \cos(\pi x))\ 
&=\ 
(F^{(2)} (x) + \pi^2 F(x) ) \sin (\pi x)
\\[12pt]
&=\ 
b^n \hsx \pi^{2n + 2} \hsx f(x) \sin (\pi x)
\\[12pt]
&=\ 
\pi^2 a^n f(x) \sin(\pi x).
\end{align*}
Therefore
\allowdisplaybreaks
\begin{align*}
\pi a^n \hsx \int\limits_0^1 \hsx f(x) \sin (\pi x)\hsx dx \  
&=\ 
\bigg(
\frac{F^\prime(x) \sin (\pi x)}{\pi} - F(x) \cos (\pi x) \hsx \bigg|_0^1
\\[12pt]
&=\ 
F(1) + F(0),
\end{align*}
an integer.  
On the other hand, 
\[
0
\ < \ 
\pi a^n \hsx \int\limits_0^1 \hsx f(x) \sin(\pi x) dx 
\ < \ 
\frac{\pi a^n }{n!}
\ < \ 
1
\]
if $n \gg 0$, from which the usual contradiction.
\vspace{0.3cm}

\un{Third Proof:} \ 
Let
\[
I_n 
\ = \ 
\int\limits_{-1}^1 \hsx 
(1 - x^2)^n \hsx \cos\bigg(\frac{\pi x}{2}\bigg) \hsx dx \qquad (n = 0, 1, 2, \ldots).
\]
Then for $-1 < x < 1$,
\allowdisplaybreaks
\begin{align*}
0 \ 
&<\  
(1 - x^2)^n \cos\bigg(\frac{\pi x}{2}\bigg) < 1\ 
\\[12pt]
&\implies 0 < I_n < 2.
\end{align*}
In addition, there is a recurrence relation, viz.
\[
\frac{\pi^2}{4} I_n 
\ = \ 
2n (2n - 1) I_{n-1} - 4n (n - 1)I_{n-2} \qquad (n \geq 2),
\]
as can be seen by integration by parts (twice).  
Using this, it follows via induction that
\[
\bigg(\frac{\pi}{2}\bigg)^{2n+1} I_n 
\ = \ 
n! P_n,
\]
where $P_n$ is a polynomial in $\ds\frac{\pi^2}{4}$ with integral coefficients of degree $\bigg[\ds\frac{n}{2}\bigg]$:
\allowdisplaybreaks
\begin{align*}
\bigg(\frac{\pi}{2}\bigg)^{2n+3} I_{n+1} \ 
&=\ 
\bigg(\frac{\pi}{2}\bigg)^{2n+3}
\bigg(\frac{2}{\pi}\bigg)^2
(2(n+1) (2n + 1) I_n 
 - 
 4(n + 1)n I_{n-1})
\\[12pt]
&=\ 
\bigg(\frac{\pi}{2}\bigg)^{2n+1} 
(2(n+1)(2n+1)I_n 
 -  
 4(n + 1)n I_{n-1})
\\[12pt]
&=\ 
2(n+1) (2n + 1)
\bigg(\frac{\pi}{2}\bigg)^{2n+1} 
I_n
 - 
 4(n + 1) \hsx n \hsx
\bigg(\frac{\pi}{2}\bigg)^2
\bigg(\frac{\pi}{2}\bigg)^{2n-1} \hsx
I_{n-1}
\\[12pt]
&=\ 
2(n+1) (2n + 1) \hsx n! \hsx P_n - 4(n + 1) \hsx n \hsx
\bigg(\frac{\pi^2}{4}\bigg)
(n-1)! \hsx P_{n-1},
\end{align*}
the degree being that of the second term, i.e., 
\[
1 + \bigg[\frac{n - 1}{2}\bigg] 
\ = \ 
\bigg[1 + \frac{n - 1}{2}\bigg] 
\ = \ 
\bigg[\frac{n + 1}{2}\bigg].
\]
Suppose now that $\ds\frac{\pi^2}{4} = \frac{a}{b}$, where $a$ and $b$ are positive integers $-$then
\[
\bigg(\frac{\pi^2}{4}\bigg)^{2n+1} \hsx I_n^2
\ = \ 
(n!)^2 (P_n)^2
\]
\qquad\qquad $\implies$
\[
\bigg(\frac{a}{b}\bigg)^{2n+1}  I_n^2 
\ = \ 
(n!)^2 (P_n)^2
\]
\qquad\qquad $\implies$
\[
\frac{a^{2n+1}}{(n!)^2 } I_n^2 
\ = \ 
b^{2n+1} (P_n)^2.
\]
But $P_n$ is a polynomial in $\ds\frac{a}{b}$ with integral coefficients of degree $\ds\bigg[\frac{n}{2}\bigg]$, 
hence the degree of $(P_n)^2$ is $2\ds\bigg[\frac{n}{2}\bigg] < 2n + 1$, hence
$b^{2n+1}(P_n)^2$ is an integer.  
To get a contradiction, simply note that
\[
0 
\ < \ 
\frac{a^{2n+1}}{(n!)^2 } \hsx I_n^2 
\ < \ 
4a \hsx \frac{(a^2)^n}{n!} \ra 0 \qquad (n \ra \infty) \qquad \text{(cf. \S0)}.
\]
\vspace{0.3cm}

\un{Fourth Proof:} \ 
The machinery employed in \S9 can also be used to establish that $\pi$ is irrational.  
So assume once again that $\pi = \ds\frac{a}{b}$, where $a$ and $b$ are positive integers, and let 
$z_0 = \pi \hsx b \hsx \sqrt{-1} = a \hsx \sqrt{-1}$ $-$then
\allowdisplaybreaks
\begin{align*}
R_n(z_0) \ 
&=\ 
Q_n(a \hsx \sqrt{-1}) \hsx e^{\pi b \sqrt{-1}} - P_n(a \hsx \sqrt{-1})
\qquad \text{(cf. \S9, \#6)}
\\[12pt]
&=\ 
Q_n(a \hsx \sqrt{-1}) \hsx (e^{\pi \sqrt{-1}})^b - P_n(a \hsx \sqrt{-1})
\\[12pt]
&=\ 
Q_n(a \hsx \sqrt{-1}) \hsx (-1)^b - P_n(a \hsx \sqrt{-1}),
\end{align*}
an element of $\Z[\sqrt{-1}]$.
Replacing $x$ by $z_0$ in \S9, \#8 (a formal maneuver), it follows that
\[
R_n(z_0) 
\ = \ 
0 
\qquad (n \gg 0).
\]
Next
\allowdisplaybreaks
\begin{align*}
\Delta(x)\ 
&\equiv\ 
Q_n(x)R_{n+1}(x) 
 - 
 Q_{n+1}(x) R_n(x)
\\[12pt]
&=\ 
Q_n(x) ( Q_{n+1}(x) e^x - P_{n+1}(x)) 
- 
Q_{n+1}(x) ( Q_n(x) e^x - P_n(x))
\\[12pt]
&=\ 
-Q_n(x) P_{n+1}(x) + Q_{n+1}(x) P_n(x)
\\[12pt]
&=\ 
-(x^n + \cdots) \cdot (-1)^{n+1})(x^{n+1} + \cdots) 
\ + \ 
(x^{n+1} + \cdots) \cdot (-1)^n(x^n + \cdots)
\end{align*}
Consequently the term of highest degree of $\Delta(x)$ is $2 (-1)^n x^{2n+1}$.  
On the other hand, the smallest nonzero degree in the expression for $R_n$ is $2n+1$ 
thus the smallest possible degree of $\Delta(x)$ is $2n+1$.  
So
\[
\Delta(x) 
\ = \ 
2 (-1)^n x^{2n+1}.
\]
Therefore $\Delta(z_0) \neq 0$.  Meanwhile
\[
R_n(z_0) 
\ = \ 
R_{n+1}(z_0) 
\ = \ 
0 
\qquad (n \gg 0).
\]

%% file: _12_irrationality_of_cos_x.tex
\chapter{
$\boldsymbol{\S}$\textbf{12}.\quad  IRRATIONALITY OF $\boldsymbol{\cos (x)}$}
\setlength\parindent{2em}
\setcounter{theoremn}{0}
\renewcommand{\thepage}{\S12-\arabic{page}}

\ \indent 
Let $x$ be a nonzero rational number.
\vspace{0.3cm}

\begin{x}{\small\bf THEOREM} \ 
$\cos (x)$ is irrational.
\end{x}
\vspace{0.3cm}

\begin{x}{\small\bf APPLICATION}  \ 
$\pi$ is irrational.
\vspace{0.2cm}

[Suppose that $\pi$ is rational $-$then $\cos (\pi)$ is irrational.  
But $\cos (\pi) = -1 \ldots$ .]
\end{x}
\vspace{0.3cm}

\begin{x}{\small\bf LEMMA} \ 
Let $g(X) \in \Z[X]$ and put
\[
f(X) 
\ = \ 
\frac{X^n}{n!} \hsx g(X) \qquad (n \in \N).
\]
Then $\forall \ j \in \N$,
\[
f^{(j)} (0) \in \Z,
\]
and in addition,
\[
(n + 1) \big | f^{(j)} (0) 
\]
except perhaps for $j = n$ ($f^{(n)} (0) = g(0)$).
\end{x}
\vspace{0.3cm}

Let $a$, $b \in \N$ $(\gcd(a,b) = 1)$ and let $p > a$ be an odd prime.
\vspace{0.2cm}

Put
\[
f(X) 
\ = \ 
\frac{X^{p-1}}{(p - 1)!} g(X),
\]
where 
\[
g(X) 
\ = \ (a - bX)^{2p} (2a - bX)^{p-1}.
\]
Then \#3 is applicable (take $n = p - 1$), hence $\forall \ j \in \N$, 
\[
f^{(j)} (0) \in \Z,
\]
and in addition,
\[
p \big | f^{(j)} (0) 
\]
except perhaps for $j = p - 1$.
\vspace{0.3cm}

\qquad FACT \ 
\[
f^{(p-1)} (0)
\ = \ g(0) 
\ = \ 
a^{2p} (2a)^{p-1} 
\ = \ 
2^{p-1} a^{3p - 1}
\]
\qquad\qquad\qquad $\implies$
\[
p \not\bigl | f^{(p-1)} (0).
\]
\vspace{0.2cm}

\begin{x}{\small\bf LEMMA} \ 
Given a real number $r$, suppose that $\phi(X) \in \Z$ $[(r - X)^2]$, i.e., 
\[
\phi(X) \ = \ 
a_{2n} (r - X)^{2n} + a_{2n - 2} (r - X)^{2n -2} + \cdots + a_2 (r - X)^2 + a_0.
\]
Then for any positive odd integer $k$, $f^{(k)}(r) = 0$.
\vspace{0.5cm}

To ensure the applicability of \#4, take $r = \ds\frac{a}{b}$ and note that
\allowdisplaybreaks
\begin{align*}
f(X) \ 
&=\ 
\frac{(r - X)^{2p} (r^2 - (r - X)^2)^{p-1}}{(p - 1)!}
b^{3p - 1}
\\[12pt]
&\in \Z[(r - X)^2].
\end{align*}

Turning now to the proof of \#1, it suffices to establish that $\cos (x)$ $(x > 0)$ is irrational.  
This said, assume that $x = \ds\frac{a}{b}$, where $a$ $b \in \N$ $(\gcd(a,b) = 1)$.  
Working with $f(X)$ per supra ($p > a$ an odd prime), introduce
\[
F(X) 
\ = \ 
f(X) - f^{(2)} (X) + f^{(4)} (X)  - \cdots - f^{(4p - 2)} (X).
\]
Then
\[
F^{(2)} (X) + F(X) \ = \  f(X).
\]
Moreover 
\allowdisplaybreaks
\begin{align*}
\frac{d}{dX} \hsx (F^\prime (X) \sin (X) - F(X) \cos (X)) \ 
&=\ 
F^{(2)} (X) \sin (X) + F(X) \sin (X)
\\[12pt]
&=\ f(X) \sin (X)
\end{align*}
\qquad\qquad $\implies$
\[
\int\limits_0^x \hsx f(X) \sin (X) dX 
\ = \ 
F^\prime (x) \sin (x) - F(x) \cos (x) + F(0).
\]
\vspace{0.2cm}

From here, the procedure is to investigate the three terms on the right and see how the supposition that $\cos (x)$ is rational leads to a contradiction.
\vspace{0.5cm}

\qquad \textbullet \quad
$f^{(2j + 1)} (x) = 0 \implies F^\prime(x) = 0$.
\vspace{0.2cm}

\qquad \textbullet \quad 
$f^{(j)} (0) \in \Z \implies F(0) \in \Z$.
\vspace{0.2cm}

\qquad \textbullet \quad 
$p \bigl | f^{(j)} (0) \qquad (j \neq p - 1)$.
\vspace{0.2cm}

\qquad \textbullet \quad
$p \not \bigl | f^{(p-1)} (0)$.
\vspace{0.2cm}

\qquad \textbullet \quad 
$F(0) = q$ \quad\ ($\gcd(p,q) = 1$).
\vspace{0.2cm}

So far then
\[
\int\limits_0^x \hsx f(X) \sin (X) dX 
\ = \ 
-F(x) \cos (x) +q.
\]

Observe next that $f(X)$ can be viewed as a function of the variable $Y = x - X$:
\allowdisplaybreaks
\begin{align*}
f(X) \ 
&=\ 
h(Y)
\\[12pt]
&=\ 
\frac{Y^{2p} (x^2 - Y^2)^{p-1}}{(p - 1)!} 
b^{3p - 1}
\\[12pt]
&=\ 
\frac{Y^{p-1} Y^{p+1} (x^2 - Y^2)^{p-1}}{(p - 1)!} 
b^{3p - 1}
\\[12pt]
&=\ 
\frac{Y^{p-1}}{(p - 1)!} \hsx
(Y^{p+1} (x^2 - Y^2)^{p-1}) b^{3p - 1}.
\end{align*}
\vspace{0.2cm}

\qquad FACT \ 
$\forall \ j \in \N$, 
\[
f^{(j)} (x) 
\ = \ 
h^{(j)} (0).
\]

In view of \#3, the $h^{(j)}(0)$ are divisible by $p$ with the possible exception of $h^{(p-1)}(0)$.  
But here
\allowdisplaybreaks
\begin{align*}
h^{(p-1)} (0) \ 
&=\ 
(Y^{p+1} (x^2 - Y^2)^{p-1}) \bigg|_{Y = 0} \hsx b^{3p - 1}
\\[12pt]
&=\ 
0.
\end{align*}

Therefore
\[
F(x) \ = \ m p
\]
for some $m \in \Z$.

Assume henceforth that
\[
\cos (x) 
\ = \ \frac{c}{d} \qquad (c, d \in \Z, \ d > 0).
\]
Then
\[
\int\limits_0^x \hsx f(X) \sin (X) dX
\ = \ 
- m p \bigg(\frac{c}{d}\bigg) + q
\]
or still, 
\[
d \hsx \int\limits_0^x \hsx f(X) \sin (X) dX
\ = \ 
- m p c + d q.
\]
However for $0 < X < x$, 
\allowdisplaybreaks
\begin{align*}
0  \ 
&<\  f(X) 
\\[12pt]
&<\  
\frac{x^{2p} (x^2)^{p-1}}{(p-1)!} b^{3p - 1}
\\[12pt]
&=\ 
\frac{x^{4p - 2}}{(p-1)!} b^{3p - 1}
\end{align*}
\qquad $\implies$
\allowdisplaybreaks
\begin{align*}
\abs{d \hsx \int\limits_0^x \hsx f(X) \sin (X) dX} 
&=\ 
d \hsx \abs{\int\limits_0^x \hsx {f(X)} {\sin (X)} dX}
\\[12pt]
&\leq\ 
d \hsx \int\limits_0^x \hsx \abs{f(X)} \abs{\sin (X)} dX
\\[12pt]
&=\ 
d \hsx \int\limits_0^x \hsx f(X) \abs{\sin (X)} dX
\\[12pt]
&\leq\ 
d \hsx \int\limits_0^x \hsx f(X) dX
\\[12pt]
&< \ 
d x \frac{x^{4p-2}}{(p-1)!} b^{3p - 1}
\\[12pt]
&=\ 
d x^3 b^2 \frac{(x^4 b^3)^{p-1}}{(p-1)!}
\\[12pt]
&=\ 
\frac{K_1 K_2^{p-1}}{(p-1)!} ,
\end{align*}
where
\[
K_1 
\ = \ 
d x^3 b^2 
\quad \text{and} \quad 
K_2 = x^4 b^3.
\]
Since
\[
\lim\limits_{p \ra \infty} \hsx \frac{K_2^{p-1}}{(p-1)!} 
\ = \ 
0
\qquad \text{(cf. \S0)},
\]
it follows that
\[
\lim\limits_{p \ra \infty} \hsx d \hsx \int\limits_0^x \hsx f(X) \sin (X) dX 
\ = \ 
0.
\]

To arrive at a contradiction, choose $p \gg 0$:
\[
- m p c + d q \in \Z - \{0\}
\]
while simultaneously
\[
\abs{d \hsx \int\limits_0^x \hsx f(X) \sin (X) dX} 
\ < \ 
1.
\]
\end{x}
\vspace{0.3cm}

\begin{x}{\small\bf APPLICATION} \ 
The values of the trigonometric functions are irrational at any nonzero rational value of the argument.
\vspace{0.2cm}

[E.g.: \ 
If $\sin (x) \in \Q$ for some $0 \neq x \in \Q$, then 
\allowdisplaybreaks
\begin{align*}
\cos (2 x) \ 
&=\ 
1 - 2 \sin^2 (x) \in \Q
\\[12pt]
&\ldots \ .]
\end{align*}
\end{x}
\vspace{0.3cm}

\begin{x}{\small\bf \un{N.B.}}\ 
The squares of these numbers are irrational.
\vspace{0.2cm}

[E.g.: \ 
\[
\cos^2 (x) 
\ = \ 
\frac{1 + \cos (2x)}{2}.]
\]
\end{x}
\vspace{0.3cm}


%% file: _13_irrationality_of_cosh_x.tex
\chapter{
$\boldsymbol{\S}$\textbf{13}.\quad  IRRATIONALITY OF $\boldsymbol{\cosh (x)}$}
\setlength\parindent{2em}
\setcounter{theoremn}{0}
\renewcommand{\thepage}{\S13-\arabic{page}}

\ \indent 
Let $x$ be a nonzero rational number.
\vspace{0.3cm}

\begin{x}{\small\bf THEOREM} \ 
$\cosh (x)$ is irrational.
\end{x}
\vspace{0.1cm}

\begin{spacing}{1.75}
The proof is similar to that in the trigonometric case.  
Thus, as there, assume that $x = \ds\frac{a}{b}$, where $a$, $b \in \N$ $(\gcd(a,b) = 1)$ and define $f(X)$ as before.  
But this time let
\end{spacing}
\[
F(X) \ = \ f(X) + f^2(X) + f^4(X) + \cdots + f^{(4p-2)} (X).
\]
Then
\[
F(X) - F^{(2)} (X) \ = \ f(X).
\]
Moreover 
\allowdisplaybreaks
\begin{align*}
\frac{d}{dX} (F(X) \cosh(X) - F^\prime(X) \sinh(X)) 
&=\ 
F(X) \sinh(X) - F^{(2)} (X) \sinh(X) 
\\[12pt]
&=\ 
f(X) \sinh(X)
\end{align*}
\qquad\qquad $\implies$
\[
\int\limits_0^x \hsx f(X) \sinh(X) dX \ = \ F(x) \cosh(x) - F^\prime(x) \sinh(x) - F(0).
\]
Note that for $0 < X < x$, 
\[
f(X) \ > \ 0 \quad \text{and} \quad \sinh(X)  \ > \ 0,
\]
thus the integral on the left hand side is positive, a point that serves to simplify matters.

Proceeding, 
\[
F^\prime(x) \ = \ 0, \quad F(x) \in \Z, \quad \text{and} \quad F(0) \in \Z.
\]

Assume henceforth that
\[
\cosh(x) \ = \ \frac{c}{d} \qquad (c, d \in \Z, \ d > 0).
\]
Then
\[
\int\limits_0^x \hsx f(X) \sinh(X) dX \ = \ F(x) \frac{c}{d} - F(0)
\]
or still, 
\[
d \int\limits_0^x \hsx f(X) \sinh(X) dX \ = \ c F(x) - d F(0).
\]
The RHS is an integer while the LHS admits the estimate
\allowdisplaybreaks
\begin{align*}
0 \ 
&<\ 
d \int\limits_0^x \hsx f(X) \sinh(X) dX
\\[12pt]
&<\ 
dx \hsx \frac{x^{4p-2} b^{3p-1}}{(p-1)!} \cdot \frac{e^x - e^{-x}}{2}
\\[12pt]
&=\ 
\frac{dx^3 b^2 (e^x - e^{-x})}{2} \cdot \frac{(x^4 b^3)^{p-1}}{(p-1)!}
\end{align*}
which is \ $< 1$ if \ $p \gg 0$ \ (for this, $p$ could have been any positive integer).  \\
Contradiction.
\vspace{0.5cm}

\begin{x}{\small\bf APPLICATION}  \ 
The values of the hyperbolic functions are irrational at any nonzero rational value of the argument.
\vspace{0.2cm}

[Use the identities
\allowdisplaybreaks
\begin{align*}
\cosh(2X) \ 
&=\ 
1 + 2 \sinh^2 (X) 
\\[12pt]
&=\ 
\frac{1 + \tanh^2(X)}{1 - \tanh^2(X)} \hsx .]
\end{align*}
\end{x}
\vspace{0.3cm}


%% file: _14_algebraic_and_transcendental_numbers.tex
\chapter{
$\boldsymbol{\S}$\textbf{14}.\quad  ALGEBRAIC AND TRANSCENDENTAL NUMBERS}
\setlength\parindent{2em}
\setcounter{theoremn}{0}
\renewcommand{\thepage}{\S14-\arabic{page}}


\begin{x}{\small\bf DEFINITION} \ 
A complex number $x$ is said to be an 
\un{algebraic number}
\index{algebraic number} 
if it is the zero of a nonzero polynomial $P(X)$ in $\Z[X]$.
\end{x}
\vspace{0.3cm}

\begin{x}{\small\bf EXAMPLE} \ 
$\sqrt{-1}$ is algebraic (consider $P(X) = X^2 + 1$).
\end{x}
\vspace{0.3cm}

\begin{x}{\small\bf \un{N.B.}}  \ 
If $x$ is algebraic, then so is its complex conjugate $\ov{x}$ and its absolute value $\abs{x}$.
\end{x}
\vspace{0.3cm}

\begin{x}{\small\bf \un{N.B.}}  \ 
If $x = a + \sqrt{-1} \hsx b$ $(a, b \in \R)$, then $x$ is algebraic iff both $a$ and $b$ are algebraic.
\end{x}
\vspace{0.3cm}

\begin{x}{\small\bf NOTATION} \ 
$\Qbar$ is the algebraic closure of $\Q$ in $\C$.
\end{x}
\vspace{0.3cm}

\begin{x}{\small\bf LEMMA} \ 
$\Qbar$ is a countable subfield of $\C$.
\end{x}
\vspace{0.3cm}

\begin{x}{\small\bf LEMMA} \ 
Suppose that $x$ is an algebraic number $-$then there is a unique nonzero polynomial $f_x \in \Z[X]$ such that 
$f_x(x) = 0$, $f_x$ is irreducible in $\Q[X]$, the leading coefficient of $f_x$ is positive, and the coefficients of $f_x$ have greatest common divisor 1.
\vspace{0.2cm}

[Note: \ 
Spelled out, 
\[
f_x(X) \ = \ a_0 + a_1X + \cdots + a_n X^n \qquad (a_n > 0)
\]
with 
\[
\gcd(a_0, a_1, \ldots, a_n) \ = \ 1.]
\]
\end{x}
\vspace{0.3cm}

\begin{x}{\small\bf DEFINITION} \ 
The polynomial $f_x$ is called the 
\un{minimal polynomial}
\index{minimal polynomial (of an algebraic number)}
of $x$.
Its degree is the 
\un{degree $\td(x)$} 
\index{degree $\td(x)$}
of $x$, hence 
\[
\td(x) \ = \ [\Q(x) : \Q].
\]

[Note: \ 
The set of real algebraic numbers of fixed degree $n$ ($\geq 2$) is dense in $\R$.]
\end{x}
\vspace{0.3cm}

\begin{x}{\small\bf DEFINITION} \ 
The zeros of $f_x$ are called the 
\un{conjugates}
\index{conjugates (of an algebraic number)}
of $x$.
\vspace{0.2cm}

[Note: \ 
They too are, of course, algebraic.]
\end{x}
\vspace{0.3cm}

\begin{x}{\small\bf EXAMPLE} \ 
Take $x$ rational, say $x = \ds\frac{a}{b}$ $(a, b \in \Z, \ b > 0, \ \gcd(a,b) = 1)$ $-$then 
\[
f_x(X) \ = \ bX - a.
\]
\end{x}
\vspace{0.3cm}

\begin{x}{\small\bf DEFINITION} \ 
An algebraic number $x$ is said to be an 
\un{algebraic integer}
\index{algebraic integer} 
if its minimial polynomial $f_x$ has leading coefficient 1.
\end{x}
\vspace{0.3cm}

\begin{x}{\small\bf EXAMPLE} \ 
$\sqrt{5}$ is an algebraic integer (consider $X^2 - 5$) but $\sqrt{5} / 2$ is not an algebraic integer 
(consider $4X^2 - 5$).
\end{x}
\vspace{0.3cm}

\begin{x}{\small\bf EXAMPLE} \ 
The integers $\Z$ are algebraic integers and if $x$ is a rational number which is also an algebraic integer then $x \in \Z$.
\vspace{0.2cm}

[Note: \ 
Accordingly, a rational number which is not an integer is not an algebraic integer.]
\end{x}
\vspace{0.3cm}

\begin{x}{\small\bf LEMMA} \ 
Under the usual operations, the set of algebraic integers forms a ring.
\end{x}
\vspace{0.3cm}

\begin{x}{\small\bf LEMMA} \ 
If $x$ is an algebraic number, then $a_n x$ is an algebraic integer.
\vspace{0.2cm}

PROOF \ 
In fact, 
\[
f_x(x) \ = \ 0
\]
\qquad\qquad $\implies$
\[
a_nx^n + a_{n-1} x^{n-1} + \cdots + a_1 x + a_0 \ = \ 0
\]
\qquad\qquad $\implies$
\[
1 (a_n x)^n + a_{n-1}(a_n x)^{n-1} + \cdots + a_n^{n-2} a_1 (a_n x)  + a_n^{n-1}a_0 
\ = \ 
0.
\]
\end{x}
\vspace{0.3cm}

Given an algebraic number $x \in \Qbar$, let $D_x$ be the set of integers $n \in \Z$ such that $n x$ is an algebraic integer 
$-$then $D_x$ is a nonzero ideal of $\Z$.
\vspace{0.3cm}

\begin{x}{\small\bf \un{N.B.}} \ 
That $D_x$ is nonzero is implied by \#15.
\end{x}
\vspace{0.3cm}

\begin{x}{\small\bf DEFINITION} \ 
A positive element of $D_x$ is called a 
\un{denominator of $x$}.
\index{denominator (of an algebraic number}
\end{x}
\vspace{0.3cm}

\begin{x}{\small\bf DEFINITION} \ 
The positive generator $\td_x$ of $D_x$ is called 
\un{the denominator of $x$}.
\index{denominator (of an algebraic number}
\end{x}

\begin{x}{\small\bf \un{N.B.}} \ 
The $a_n$ of \#15 needn't be $\td_x$ (consider $4X^2 + 2X + 1$).
\end{x}
\vspace{0.3cm}

\begin{x}{\small\bf DEFINITION} \ 
A complex number $x$ is said to be a 
\un{transcendental number}
\index{transcendental number} 
if it is not an algebraic number.
\end{x}
\vspace{0.3cm}

Therefore the set of transcendental numbers is the complement of the field $\ov{\Q}$ in the field $\C$.
\vspace{0.3cm}


\begin{x}{\small\bf \un{N.B.}} \ 
In general, the sum or product of two transcendental numbers is not transcendental.  
However the sum of a transcendental number and an algebraic number is a transcendental number and the product of a transcendental number and a nonzero algebraic number is again a transcendental number.
\end{x}
\vspace{0.3cm}

\begin{x}{\small\bf EXAMPLE} \ 
$e$ is transcendental (cf. \S17, \#1) and $\pi$ is transcendental (cf \S19, \#1) but it is unknown whether $e + \pi$ and $e \pi$ are  transcendental (cf. \S2, \#29).
\end{x}
\vspace{0.3cm}

\[
\text{APPENDIX}
\]
\vspace{0.3cm}

Given an algebraic number $x \neq 0$, let $x_1 = x$, $x_2, \ldots, x_n$ $(n = d(x))$ be the conjugates of $x$ (cf. \#9) and put
\[
H(x) \ = \ \max\limits_{1 \leq j \leq n} \hsx \abs{x_j},
\]
the 
\un{house}
\index{house (of an algebraic number)}
of $x$.
\vspace{0.3cm}

{\small\bf LEMMA} \  
Let $T \in D_x$ $(T > 0)$ $-$then 
\[
\abs{x} \ \geq \ \frac{1}{T^n \hsx H(x)^{n-1}}\hsx.
\]


%% file: _15_liouville_theory.tex
\chapter{
$\boldsymbol{\S}$\textbf{15}.\quad  LIOUVILLE THEORY}
\setlength\parindent{2em}
\setcounter{theoremn}{0}
\renewcommand{\thepage}{\S15-\arabic{page}}


\begin{x}{\small\bf RAPPEL} \ 
(cf. \S7, \#17) \ 
Given $x = \ds\frac{a}{b} \in \Q$ ($a, \ b \in \Z$, $b > 0$, $\gcd(a,b) = 1$, for any coprime pair 
$(p,q)$ $(q > 0)$ with 
\[
\frac{a}{b} 
\ \neq \ 
\frac{p}{q}
\]
there follows
\[
\abs{\frac{a}{b}  - \frac{p}{q}} 
\ \geq \ 
\frac{1}{bq}.
\]
\end{x}
\vspace{0.3cm}

\begin{x}{\small\bf THEOREM} \ 
If $x$ is real and algebraic of degree $d(x) = n$ (cf. \S14, \#8), then there is a constant 
$C = C(x) > 0$ such that for any coprime pair $(p,q)$ $(q > 0)$, 
\[
\abs{x - \frac{p}{q}} 
\ > \ 
\frac{C}{q^n}.
\]

\vspace{0.2cm}

PROOF \ 
The case $d(x) = 1$ is \#1 above (choose $C = C(x) < \ds\frac{1}{b}$), so take $d(x) \geq 2$ and recall that
\[
f_x(X) 
\ = \ 
a_0 + a_1 X + \cdots + a_n X^n
\]
is the minimal polynomial of $x$.  
Let \mM be the maximum value of $\abs{f_x^\prime(X)}$ on $[x-1, x+1]$, let $\{y_1, \ldots, y_m\}$ 
$(m \leq n)$ be the distinct zeros of $f_x$ which are different from $x$, and then choose \mC:
\[
0 
\ < \ 
C 
\ < \ 
\min \bigg\{1, \frac{1}{M}, \abs{x - y_1}, \ldots, \abs{x - y_m}\bigg\}.
\]
To arrive at a contradiction, suppose that for some coprime pair $(p,q)$ $(q > 0)$
\[
\abs{x - \frac{p}{q}} 
\ \leq \ 
\frac{C}{q^n}
\]
or still, 
\[
\leq \ 
C 
\ < \ 
\min \{1,  \abs{x - y_1}, \ldots, \abs{x - y_m}\}.
\]
Of course, 
\[
\abs{x - \frac{p}{q}} 
\ > \ 
0,
\]
$x$ being irrational.   
And
\[
\abs{x - \frac{p}{q}} 
\ = \ 
\abs{\frac{p}{q} - x} 
\ < \ 
1 
\implies 
x - 1 
\ < \ 
\frac{p}{q} 
\ < \ 
x + 1.
\]
In addition
\allowdisplaybreaks
\begin{align*}
0 < \abs{x - \frac{p}{q}}  &< \abs{x - y_1}, \ldots, \abs{x - y_m}
\\[12pt]
\implies \frac{p}{q} 
&\neq y_k \qquad (k = 1, \ldots, m)
\\[12pt]
&\implies 
f_x\bigg(\frac{p}{q}\bigg) \neq 0.
\end{align*}
Owing to the mean value theorem, there is an $x_0$ between $\ds\frac{p}{q}$ and $x$ such that 
\[
\abs{f_x(x) - f_x\bigg(\frac{p}{q}\bigg)} 
\ = \ 
\abs{x - \frac{p}{q}} \hsx \abs{f_x^\prime(x_0)},
\]
i.e., 
\allowdisplaybreaks
\begin{align*}
\abs{f_x\bigg(\frac{p}{q}\bigg)} 
&=\ \abs{x - \frac{p}{q}} \hsx \abs{f_x^\prime(x_0)}
\\[12pt]
&\implies
\abs{f_x^\prime(x_0)} \neq 0
\\[12pt]
&\implies
\\[12pt]
&\qquad
\abs{x - \frac{p}{q}} = \frac{\abs{f_x\bigl(\frac{p}{q}\bigr)}}{ \abs{f_x^\prime(x_0)}}
\\[12pt]
&\qquad\qquad\quad\ \hsx  \geq \frac{f_x\bigl(\frac{p}{q}\bigr)}{M}.
\end{align*}
But
\allowdisplaybreaks
\begin{align*}
0 < \abs{f_x \bigg(\frac{p}{q}\bigg)} \ 
&=
\abs{\sum\limits_{j=0}^n \hsx a_j \bigg(\frac{p}{q}\bigg)^j}
\\[15pt]
&=\ 
\abs{\sum\limits_{j=0}^n \hsx a_j p^j q^{n-j}} \hsx \big/ q^n. 
\end{align*}
Since the numerator of this fraction is a positive integer, it follows that
\[
\abs{\sum\limits_{j=0}^n \hsx a_j p^j q^{n-j}} \ \geq \ 1,
\]
thus
\[
\abs{f_x \bigg(\frac{p}{q}\bigg)} \ \geq \ \frac{1}{q^n}.
\]
Finally
\allowdisplaybreaks
\begin{align*}
\abs{x - \frac{p}{q}}\ 
&\geq\ 
\frac{\abs{f_x \big(\frac{p}{q}\big)}}{M}
\\[15pt]
&\geq\ 
\frac{1}{Mq^n}
\\[15pt]
&>\ 
\frac{C}{q^n}.
\end{align*}
Contradiction.
\end{x}
\vspace{0.3cm}
\begin{x}{\small\bf REMARK} \ 
The preceding proof goes through if $f(X) \in \Z[X]$ has degree $n > 1$ and $x$ is an irrational root of $f(X)$.
\end{x}
\vspace{0.3cm}

\begin{x}{\small\bf DEFINITION} \ 
A real number $x$ is a 
\un{Liouville number}
\index{Liouville number} 
if for every positive integer $k$ there exist $p, \ q \in \Z$ $(q > 1$, $\gcd(p,q) = 1)$ such that 
\[
0 \ < \ \abs{x - \frac{p}{q}} \ < \  \frac{1}{q^k}.
\]
\end{x}
\vspace{0.3cm}
\begin{x}{\small\bf NOTATION} \ 
\bL is the subset of $\R$ whose elements are the Liouville numbers.
\end{x}
\vspace{0.3cm}

\begin{x}{\small\bf LEMMA} \ 
Every Liouville number is irrational.
\vspace{0.2cm}

\begin{spacing}{1.75}
PROOF \ 
Suppose instead that $x = \ds\frac{a}{b}$ $(a, b \in \Z$, $b > 0$, $\gcd(a,b) = 1$).  
Let $k$ be a positive integer: \ $2^{k-1} > b$ and take $p, q$: $\ds\frac{a}{b} \neq \ds\frac{p}{q}$ $-$then
\end{spacing}
\allowdisplaybreaks
\begin{align*}
\abs{x - \frac{p}{q}} \ 
&=\ 
\abs{\frac{a}{b} - \frac{p}{q}}
\\[12pt]
&=\ 
\frac{\abs{aq - bp}}{bq}
\\[12pt]
&\geq\ 
\frac{1}{bq}
\\[12pt]
&>\ 
\frac{1}{2^{k-1}q}
\\[12pt]
&\geq\ 
\frac{1}{q^{k-1}q} \qquad (q \geq 2)
\\[12pt]
&=\ 
\frac{1}{q^k}.
\end{align*}
So $x$ is not a Liouville number.
\vspace{0.2cm}

Therefore
\[
\bL \subset \PP.
\]
\end{x}
\vspace{0.3cm}

\begin{x}{\small\bf THEOREM} \ 
Every Liouville number is transcendental.
\vspace{0.2cm}

PROOF \ 
Assume that $x$ is an algebraic irrational number with $d(x) = n$, hence per \#2, for any coprime pair 
$(p,q)$ $(q > 0)$, 
\[
\abs{x - \frac{p}{q}} 
\ > \ 
\frac{C}{q^n}.
\]
Choose a positive integer $r: 2^r \geq \ds\frac{1}{C}$ and then, using the definition of Liouville number, 
choose $p, q$:
\[
0 \ < \ \abs{x - \frac{p}{q}}  \ < \ \frac{1}{q^{n+r}} \qquad (k \equiv n + r).
\]
But 
\[
\frac{1}{q^{n+r}} \ \leq \ \frac{1}{2^rq^n} \ \leq \ \frac{C}{q^n}
\]
\qquad\qquad $\implies$
\[
\abs{x - \frac{p}{q}}  \ < \ \frac{C}{q^n}.
\]
On the other hand, 
\[
\abs{x - \frac{p}{q}} 
\ > \ 
\frac{C}{q^n} \qquad \text{(cf. \#2)}.
\]
Contradiction.
\vspace{0.2cm}

Therefore
\[
\bL \subset \T \subset \PP.
\]
\end{x}
\vspace{0.3cm}

\begin{x}{\small\bf REMARK} \ 
Not every transcendental number is a Liouville number, e.g., $e$ and $\pi$ are transcendental but not in \bL.
\end{x}
\vspace{0.3cm}

\begin{x}{\small\bf EXAMPLE} \ 
Let $a$ be a positive integer $\geq 2$.  Put
\[
x 
\ = \ 
\sum\limits_{j=1}^\infty \hsx
\frac{1}{a^{j!}}.
\]
Then $x$ is a Liouville number.
\vspace{0.2cm}

[Define a sequence of rationals $\ds\frac{p_k}{q_k}$ $(k = 1, 2, \ldots)$ by the prescription
\[
\frac{p_k}{q_k} 
\ = \ 
\sum\limits_{j=1}^k \hsx
\frac{1}{a^{j!}}, 
\qquad q_k \ = \ a^{k!}.
\]
Then 
\[
\abs{x - \frac{p_k}{q_k}} 
\ = \ 
\sum\limits_{j=k+1}^\infty \hsx
\frac{1}{a^{j!}}.
\]
But 
\allowdisplaybreaks
\begin{align*}
\sum\limits_{j=k+1}^\infty \hsx\frac{1}{a^{j!}}\
&<\ 
\sum\limits_{j=(k+1)!}^\infty \hsx 
\frac{1}{a^j}
\\[15pt]
&=\
\frac{1}{a^{(k+1)!}}  \hsx 
\sum\limits_{j=0}^\infty \hsx\frac{1}{a^{j}}
\\[15pt]
&=\
\frac{1}{a^{(k+1)!}} 
\cdot
\frac{a}{a - 1}
\\[15pt]
&=\
\frac{1}{q_k^{k+1}} 
\cdot
\frac{a}{a - 1}
\\[15pt]
&\leq\
\frac{2}{q_k q_k^k}
\\[15pt]
&\leq\
\frac{1}{q_k^k} \qquad (q_k \geq 2).
\end{align*}
So, $\forall \ k \in \N$, 
\[
0 
\ < \ 
\abs{x - \frac{p_k}{q_k}} 
\ < \ 
\frac{1}{q_k^k}.
\]
Therefore $x$ is in \bL. (cf. \#4).]
\end{x}
\vspace{0.3cm}
\begin{x}{\small\bf \un{N.B.}} \ 
The preceding discussion can be generalized.  Thus fix an integer
$n \geq 2$ and a sequence of integers $m_j \in \{0, 1, 2, \ldots, n-1\}$ $(j = 1, 2 \ldots)$ such that 
$m_j \neq 0$ for infinitely many $j$.  Put
\[
x 
\ = \ 
\sum\limits_{j=1}^\infty \hsx
\frac{m_j}{n^{j!}}.
\]
Then $x$ is a Liouville number.
\vspace{0.2cm}

[Define a sequence of rationals $\ds\frac{p_k}{q_k}$ $(k = 1, 2, \ldots)$ by the prescription
\[
\frac{p_k}{q_k} 
\ = \ 
\sum\limits_{j=1}^k \hsx
\frac{m_j}{n^{j!}}, 
\qquad q_k \ = \ n^{k!}.
\]
Then 
\[
\abs{x - \frac{p_k}{q_k}} 
\ = \ 
\sum\limits_{j=k+1}^\infty \hsx
\frac{m_j}{n^{j!}}.
\]
But as above
\allowdisplaybreaks
\begin{align*}
\sum\limits_{j=k+1}^\infty \hsx\frac{m_j}{n^{j!}}\
&\leq\ 
\sum\limits_{j=k+1}^\infty \hsx 
\frac{n-1}{n^{j!}}
\\[15pt]
&<\ 
\sum\limits_{j=(k+1)!}^\infty \hsx 
\frac{n-1}{n^j}
\\[15pt]
&=\
\frac{n-1}{n^{(k+1)!}}  \hsx 
\sum\limits_{j=0}^\infty \hsx\frac{1}{n^{j}}
\\[15pt]
&=\
\frac{n-1}{n^{(k+1)!}}  \cdot  \frac{n}{n - 1}
\\[15pt]
&=\
\frac{n}{n^{(k+1)!}} 
\\[15pt]
&\leq\
\frac{n^{k!}}{n^{(k+1)!}} 
\\[15pt]
&=\
n^{k! - (k+1)!}
\\[15pt]
&=\
\bigl(n^{-k!}\bigr)^k
\\[15pt]
&=\
\bigl(q_k^{-1}\bigr)^k
\\[15pt]
&=\
\bigg(\frac{1}{q_k}\bigg)^k
\\[15pt]
&=\
\frac{1}{q_k^k}.
\end{align*}
So, $\forall \ k \in \N$, 
\[
0 
\ < \ 
\abs{x - \frac{p_k}{q_k}} 
\ < \ 
\frac{1}{q_k^k}.
\]
Therefore $x$ is in \bL (cf. \#4).]
\end{x}
\vspace{0.3cm}

\begin{x}{\small\bf EXAMPLE} \ 
Put
\[
x 
\ = \ 
\sum\limits_{j = 1}^\infty \hsx 
\frac{1}{2^{2^{^j}}}.
\]
Then $x$ is a Liouville number.
\end{x}
\vspace{0.3cm}

In \#10, it is traditional to take $n = 10$, hence $m_j \in \{0, 1, 2, \ldots, 9\}$ $(j = 1, 2, \ldots)$.
\vspace{0.3cm}

\begin{x}{\small\bf LEMMA} \ 
Put
\[
x 
\ = \ 
\sum\limits_{j = 1}^\infty \hsx 
m_j 10^{-j!}, 
\qquad 
y
\ = \ 
\sum\limits_{j = 1}^\infty \hsx 
n_j 10^{-j!}.
\]
Assume $m_j \neq n_j$ for some $j$ and let $k$ be the least index $j$ such that $m_j \neq n_j$ $-$then
$x \neq y$.
\vspace{0.2cm}

PROOF \ 
\allowdisplaybreaks
\begin{align*}
\abs{x - y} \ 
&=\ 
\abs{(m_k - n_k)10^{-k!} + \sum\limits_{j = k+1}^\infty \hsx (m_j - n_j)10^{-j!}\hsx} \ 
\\[15pt]
&\geq\
\abs{m_k - n_k} 10^{-k!} - 
\abs{
\sum\limits_{j = k+1}^\infty \hsx
(m_j - n_j)10^{-j!}
}
\\[15pt]
&\geq\
\abs{m_k - n_k} 10^{-k!} - 
\sum\limits_{j = k+1}^\infty \hsx
\abs{m_j - n_j}10^{-j!}
\\[15pt]
&\geq\
10^{-k!} \hsx -
\sum\limits_{j = k+1}^\infty \hsx
(9) 10^{-j!}
\\[15pt]
&>\
10^{-k!} \hsx -
\sum\limits_{j = (k+1)!}^\infty \hsx
(9) 10^{-j}
\\[15pt]
&=\ 
10^{-k!} \hsx - 
(9) \hsx
\big(
10^{-(k+1)!} + 10^{-(k+1)! - 1} + 10^{-(k+1)! - 2} + \cdots
\big)
\\[15pt]
&=\ 
10^{-k!} \hsx - (9) \hsx 10^{-(k+1)!} \hsx 
(1 + 10^{-1} + 10^{-2} + \cdots)
\\[15pt]
&=\ 
10^{-k!} \hsx - (9) \hsx 10^{-(k+1)!} \hsx \frac{1}{1 - \frac{1}{10}}
\\[15pt]
&=\ 
10^{-k!} \hsx - (9) \hsx 10^{-(k+1)!} \hsx \frac{1}{\frac{9}{10}}
\\[15pt]
&=\ 
10^{-k!} \hsx - 10^{-(k+1)!} \hsx (10)
\\[15pt]
&=\ 
10^{-k!} \hsx - 10^{-(k+1)! + 1}
\\[15pt]
&=\ 
10^{-k!} \hsx - 10^{-k! (k+1)} \hsx 10
\\[15pt]
&=\ 
10^{-k!} \hsx - 10^{-k! \hsx k - k!} 10
\\[15pt]
&=\ 
10^{-k!} \hsx - 
10^{-k!} 10^{-k! \hsx k} \hsx 10
\\[15pt]
&=\ 
10^{-k!} 
\big(
1 - 10^{-k! \hsx k} 10 
\big).
\end{align*}
And
\[
1 - 10^{-k! \hsx k} 10 
\ \geq \ 
0
\]
since
\[
\frac{1}{10} 
\ \geq \ 
\frac{1}{10^{k! \hsx k}}.
\]
\end{x}
\vspace{0.3cm}

\begin{x}{\small\bf SCHOLIUM} \ 
The set of Liouville numbers is uncountable.
\vspace{0.2cm}

[The Liouville numbers of the form
\[
\sum\limits_{j = 1}^\infty \hsx 
m_j 10^{-j!} 
\]
constitute an uncountable set (use a Cantor diagonalization argument).]
\end{x}
\vspace{0.3cm}

\begin{x}{\small\bf THEOREM} \ 
Suppose that $f(X) \in \Z[X]$ has degree $\geq 1$ and let $x \in \bL$ $-$then $f(x) \in \bL$.
\end{x}
\vspace{0.3cm}

To begin with:
\vspace{0.3cm}

\begin{x}{\small\bf LEMMA} \ 
If the degree of $f(X) \in \R[X]$ is $\geq 1$ and if $a \in \R$, then there is a polynomial $g(X) \in \R[X]$ such that 
\[
f(X) - f(a) 
\ = \ 
(X - a) g(X).
\]
\vspace{0.2cm}

PROOF \
Write
\[
f(X)
\ = \ 
\sum\limits_{j=0}^r \hsx C_j X^j.
\]
Then for $j \geq 1$, 
\allowdisplaybreaks
\begin{align*}
X^j - a^j \ 
&=\ 
(X- a) (X^{j-1} + aX^{j-2} + a^2 X^{j-3} + \cdots + a^{j-2} X + a^{j-1})
\\[12pt]
&=\ 
(X - a) g_j(X).
\end{align*}
Therefore
\allowdisplaybreaks
\begin{align*}
f(X) - f(a) \ 
&=\ 
C_0 + 
\sum\limits_{j=1}^r \hsx C_j X^j - C_0 - \sum\limits_{j=1}^r \hsx C_j a^j
\\[12pt]
&=\ 
\sum\limits_{j=1}^r \hsx C_j (X^j - a^j) 
\\[12pt]
&=\ 
\sum\limits_{j=1}^r \hsx C_j (X - a) g_j(X)
\\[12pt]
&=\ 
(X - a) \hsx
\sum\limits_{j=1}^r \hsx C_j g_j(X)
\\[12pt]
&\equiv\ 
(X - a) g(X).
\end{align*}
\end{x}
\vspace{0.3cm}

To set up the particulars for \#14, note first that 
$\{X:X \neq x \ \& \ f(X) = f(x)\}$ is a finite set (the degree of $f(X)$ being by assumption $\geq 1$).  
Fix $\delta > 0$ subject to 
\[
0 \ < \ \delta \ < \ 
\min\{\abs{X - x}: X \neq x \ \& \ f(X) = f(x) \}
\]
and put
\[
M \ = \ 
\max\{\abs{g(X)}: \abs{X - x} \leq \delta\}.
\]
Bearing in mind the definition figuring in \#4, let $k$ be a positive integer and choose a natural number 
$m > kr$ ($r$ the degree of $f$) such that 
\[
1 
\ < \ 
\delta 2^m 
\quad \text{and} \quad 
M 2^{k r} \ < \ 2^m.
\]
Next, determine $p, q \in \Z$ $(q > 1$, $\gcd(p,q) = 1)$: 
\[
0 
\ < \ 
\abs{x - \frac{p}{q}} 
\ < \ 
\frac{1}{q^m}.
\]
\vspace{0.2cm}

\qquad \un{Step 1:} \ 
\[
\abs{x - \frac{p}{q}} 
\ < \ 
\frac{1}{q^m}
\ \leq \ 
\frac{1}{2^m}
\ < \ 
\delta
\]
\qquad\qquad $\implies$
\[
\abs{g\bigg(\frac{p}{q}\bigg)} \leq M 
\quad \text{and} \quad 
f\bigg(\frac{p}{q}\bigg) \neq f(x).
\]
\vspace{0.2cm}

\qquad \un{Step 2:} \ 
\[
M 2^{k r}  \ <  \ 2^m 
\implies 
M  \ <  \ 2^{m - kr}
\]
\qquad\qquad $\implies$
\[
\abs{g\bigg(\frac{p}{q}\bigg)} \ \leq \ M  < 2^{m - kr}  \ \leq  \ q^{m - kr}.
\]
\vspace{0.2cm}

\qquad \un{Step 3:} \ 
\allowdisplaybreaks
\begin{align*}
0 \ <  \ \abs{f(x) - f\bigg(\frac{p}{q}\bigg)} \ 
&=\ 
\abs{x - \frac{p}{q}} \abs{g\bigg(\frac{p}{q}\bigg)}
\\[12pt]
&<\ 
\frac{1}{q^m} q^{m - kr}
\\[12pt]
&=\ 
\bigg(\frac{1}{q^r}\bigg)^k.
\end{align*}
\vspace{0.2cm}

\qquad \un{Step 4:} \ 
Write
\[
f(X) 
\ = \ 
\sum\limits_{j=0}^r \hsx
C_j X^j 
\qquad (C_j \in \Z).
\]
Then
\allowdisplaybreaks
\begin{align*}
f\bigg(\frac{p}{q}\bigg) \
&=\ 
\sum\limits_{j=0}^r \hsx
C_j \bigg(\frac{p}{q}\bigg)^j
\\[12pt]
&=\ 
\bigg(
\sum\limits_{j=0}^r \hsx
C_j p^j q^{r-j}\bigg) \big/ q^r
\\[12pt]
&=\ 
\frac{C}{q^r},
\end{align*}
where $C \in \Z$.
\vspace{0.2cm}

\qquad \un{Step 5:} \ 
\allowdisplaybreaks
\begin{align*}
0 \ 
&<\ 
\abs{f(x) - f\bigg(\frac{p}{q}\bigg)} 
\\[12pt]
&=\ \abs{f(x) - \frac{C}{q^r}}
\\[12pt]
&<\ 
\bigg(\frac{1}{q^r}\bigg)^k.
\end{align*}
To fullfill the requirements of \#4, it remains only to take 
\[
\begin{cases}
\ \text{``$p$''} \ = \  C\\
\ \text{``$q$''} \ = \  q^r
\end{cases}
.
\]
\vspace{0.3cm}

\begin{x}{\small\bf APPLICATION} \ 
If $a \neq 0$, $b \neq 0$ are integers and if $x \in \bL$, then 
\[
a + bx \in \bL.
\]
\vspace{0.2cm}

[Consider
\[
f(X) \ =  \ a + bX.]
\]
\end{x}
\vspace{0.3cm}

\begin{x}{\small\bf APPLICATION} \ 
If $x \in \bL$, then $\forall \ n \in \N$, $x^n \in \bL$.
\vspace{0.2cm}

[Consider
\[
f(X) \ =  \ X^n.]
\]
\end{x}
\vspace{0.3cm}


\begin{x}{\small\bf LEMMA} \ 
If $x$ is a Liouville number and if $r \in \Q$ is nonzero, then $r x \in \bL$.
\vspace{0.2cm}

PROOF \ 
Write $r = \ds\frac{a}{b}$ $(a, b \in \Z, \ b > 0$.  
Given a natural number $k$, choose a natural number $m > k$: 
\[
\abs{a} \hsx b^{k-1} \ < \ 2^{m-k}.
\]
Next, per the definition of \bL (cf. \#4), there exist $p, \ q \in \Z$ $(q > 1$, $\gcd(p,q) = 1$):
\[
0 \ < \ \abs{x - \frac{p}{q}}  \ < \ \frac{1}{q^m}.
\]
Therefore
\allowdisplaybreaks
\begin{align*}
0 \ 
&< \
\abs{rx - \frac{ap}{bq}}
\\[12pt]
&< \
\frac{\abs{r}}{q^m}
\\[12pt]
&< \
\hsx \frac{\abs{a}}{b q^m}
\\[12pt]
&< \
\frac{2^{m-k}}{b^{k-1}} \cdot \frac{1}{b q^m}
\\[12pt]
&\leq\ 
\frac{q^{m-k}}{b^{k-1}} \cdot \frac{1}{b q^m}
\\[12pt]
&=\ 
\frac{1}{(b q)^k}.
\end{align*}
\vspace{0.2cm}

[Note: \ 
The assertion may be false if $r$ is merely algebraic.  For example, consider
\[
\sqrt{3/2} \hspace{0.2cm}
\sum\limits_{j=1}^\infty \hsx
\frac{1}{10^{j!}}\hsx .]
\]
\end{x}
\vspace{0.3cm}


\begin{x}{\small\bf APPLICATION} \ 
Every interval $]a,b[$ $(a < b)$ contains a Liouville number.
\vspace{0.2cm}

[Take a positive Liouville number $x$ and consider
\[
\bigg]\frac{a}{x}, \frac{b}{x}\bigg[\ .
\]
Fix a nonzero rational number $r$:
\[
\frac{a}{x} \ < \ r \ < \ \frac{b}{x} \qquad \text{(cf. \S2, \#15)}.
\]
Then
\[
a \ < \ rx \ < \ b.]
\]
\end{x}
\vspace{0.3cm}

\begin{x}{\small\bf SCHOLIUM} \ 
\bL is a dense subset of $\R$ (cf. \S2, \#14).
\end{x}
\vspace{0.3cm}

\begin{x}{\small\bf THEOREM} \ 
Let $f(X) \in \Q[X]$ be nonconstant and suppose that $x \in \bL$ $-$then $f(x) \in \bL$.
\vspace{0.2cm}

PROOF \ 
Choose $n \in \N$:
\[
(n f) (X) \in \Z[X].
\]
Then
\[
(n f) (x) \in \bL \qquad \text{(cf. \#14)} 
\implies
\frac{1}{n} (n f) (x) \in \bL \qquad \text{(cf. \#18)},
\]
i.e., $f(x) \in \bL$.
\vspace{0.2cm}

[In particular, the sum of a rational number $\ds\frac{a}{b}$ and a Liouville number $x$ is again a Liouville number:
\[
\frac{a}{b} + x 
\ = \ 
\frac{1}{b} (a + bx).]
\]
\end{x}
\vspace{0.3cm}

\begin{x}{\small\bf THEOREM} \ 
The set of Liouville numbers in $[0,1]$ is a set of measure 0.
\vspace{0.2cm}

PROOF \ 
Fix $\epsilon > 0$.  Let $k$ be a positive integer such that
\[
4 \hsx 
\sum\limits_{q=2}^\infty \hsx 
\frac{1}{q^{k-1}} 
\ < \ 
\epsilon.
\]
That such a choice is possible can be seen by noting that 
\allowdisplaybreaks
\begin{align*}
4 \hsx \sum\limits_{q=2}^\infty \hsx \frac{1}{q^{k-1}} \ 
&=\ 
4 \hsx 
\bigg(
\frac{1}{2^{k-1}} + \frac{1}{3^{k-1}} + \cdots\bigg)
\\[15pt]
&=\ 
4 \hsx \cdot
\frac{1}{2^{k-3}} 
\bigg(
\frac{1}{2^2} + \frac{1}{3^2} + \cdots\bigg).
\end{align*}
This said, let $x$ be a Liouville number in $[0,1]$ and per \#4, write
\[
0
\ < \ 
\abs{x - \frac{p}{q}}
\ < \ 
\frac{1}{q^k}
\]
or still, 
\[
\frac{p}{q} - \frac{1}{q^k} 
\ < \
x 
\ < \
\frac{p}{q} + \frac{1}{q^k}.
\]
Put
\[
I_{p/q} 
\ = \ 
\bigg]
\frac{p}{q} - \frac{1}{q^k} , 
\frac{p}{q} + \frac{1}{q^k}\bigg[ \ ,
\]
an open interval of length
\[
\frac{p}{q} + \frac{1}{q^k} - \bigg( \frac{p}{q} - \frac{1}{q^k} \bigg) 
\ = \ 
\frac{2}{q^k}.
\]
Since $x \in [0,1]$ and $\ds\frac{1}{q^k} \leq \ds\frac{1}{2}$, it follows that 
\[
\frac{p}{q} \in \bigg] -\frac{1}{2}, \frac{3}{2} \bigg[\ ,
\]
i.e., 
\[
-\frac{1}{2} 
\ < \ 
\frac{p}{q}
\ < \ 
\frac{3}{2}
\implies
-\frac{q}{2} 
\ < \ 
p
\ < \ 
\frac{3q}{2}.
\]
Therefore the total number of $I_{p/q}$ is $\leq 2q$.\\
Put
\[
I(q) 
 \ = \ 
\bigcup\limits_{p/q} \hsx I_{p/q},
\]
a set of measure
\allowdisplaybreaks
\begin{align*}
\leq \ \sum\limits_{p/q} \hsx \frac{2}{q^k} \ 
&=\ 
\frac{2}{q^k} \hsx 
\sum\limits_{p/q} \hsx  1
\\[12pt]
&\leq \ 
\frac{2}{q^k} \cdot 2 q
\\[12pt]
&=\  
\frac{4 q}{q^k}.
\end{align*}
The set of Liouville numbers in $[0,1]$ is contained in
\[
\bigcup\limits_{q > 1} \hsx I(q),
\]
a set of measure
\[
\leq \ 
\sum\limits_{q=2}^\infty \hsx 
\frac{4q}{q^k} 
\ = \ 
4 \hsx 
\sum\limits_{q=2}^\infty \hsx 
\frac{1}{q^{k-1}} 
\ < \ 
\epsilon,
\]
from which the assertion.
\end{x}
\vspace{0.3cm}

\begin{x}{\small\bf APPLICATION} \ 
There are transcendental numbers that are not Liouville numbers.
\vspace{0.2cm}

[Let \mS be the set of algebraic numbers in $[0,1]$ and let \mT be the set of transcendental numbers 
in $[0,1]$ $-$then
\[
[0,1] 
\ = \ 
S \hsx \cup \hsx T, 
\qquad 
S \hsx \cap \hsx T 
\ = \ 
\emptyset.
\]
Since \mS is countable, it is of measure 0, hence \mT is of measure 1.]
\vspace{0.2cm}

[Note: \ 
Almost all transcendental numbers in $[0,1]$ are non-Liouville numbers.]
\end{x}
\vspace{0.3cm}

Working within $\R$, it follows that \bL is a set of measure 0.
\vspace{0.3cm}

\begin{x}{\small\bf NOTATION} \ 
Given $k \in \N$, put
\[
U_k 
\ = \ 
\bigcup\limits_{q \geq 2} \hsx 
\bigcup\limits_{p \in \Z} \hsx 
\bigg] \frac{p}{q} - \frac{1}{q^k}, \frac{p}{q} + \frac{1}{q^k} \bigg[ 
\ - \ 
\bigg\{\frac{p}{q}\bigg\}
\]
or still, 
\[
U_k 
\ = \ 
\bigcup\limits_{q \geq 2} \hsx 
\bigcup\limits_{p \in \Z} \hsx 
\bigg\{
x \in \R: \ 0 < \abs{x - \frac{p}{q}} \ < \ \frac{1}{q^k}
\bigg\}.
\]
\end{x}
\vspace{0.3cm}

\begin{x}{\small\bf LEMMA} \ 
$U_k$ is an open dense subset of $\R$.

\vspace{0.2cm}
[Each $\ds\frac{p}{q} \in \Q$ belongs to the closure of $U_k$.]
\end{x}
\vspace{0.3cm}

\begin{x}{\small\bf LEMMA} \ 
\[
\bL \ = \ 
\bigcap\limits_{k=1}^\infty \hsx  U_k.
\]
\end{x}
\vspace{0.3cm}

\begin{x}{\small\bf RAPPEL} \ 
A 
\un{$G_\delta$-subset}
\index{$G_\delta$-subset} 
of a topological space \mX is the countable intersection of open dense subsets of \mX.
\end{x}
\vspace{0.3cm}

Therefore \bL is a $G_\delta$-subset of $\R$.
\vspace{0.5cm}

\begin{x}{\small\bf RAPPEL} \ 
If \mX is a complete metric space and if $\{G_n\}$ is a sequence of open dense subsets of \mX, then 
\[
\bigcap\limits_{n=1}^\infty \hsx G_n
\]
is not empty and, in fact, is dense in \mX.
\end{x}
\vspace{0.3cm}

Therefore \bL is a dense subset of $\R$ (cf. \#20).
\vspace{0.5cm}

\begin{x}{\small\bf RAPPEL} \ 
If \mX is a complete metric space without isolated points and if \mS is a $G_\delta$-subset of \mX, then \mS is uncountable.
\end{x}
\vspace{0.3cm}

Therefore \bL is an uncountable subset of $\R$ (cf. \#13).
\vspace{0.5cm}

\begin{x}{\small\bf THEOREM} \ 
Every real number $x$ is the sum of two Liouville numbers:
\[
x 
\ = \ 
\alpha + \beta 
\qquad (\alpha, \hsx \beta  \in \bL).
\]
\end{x}
\vspace{0.3cm}

\begin{x}{\small\bf THEOREM} \ 
Every nonzero real number $x$ is the product of two Liouville numbers:
\[
x 
\ = \ 
\alpha \hsx \beta 
\qquad (\alpha, \hsx \beta  \in \bL).
\]
\end{x}
\vspace{0.75cm}

It will be enough to sketch the proof of \#30.
\vspace{0.5cm}

\qquad \un{Step 1:} \ 
Put
\[
\alpha 
\ = \ 
\sum\limits_{j=1}^\infty \hsx 10^{-j!}.
\]
Then
\[
0 \ = \ \alpha + (-1)\alpha, \ \ 
1 \ = \ \alpha + (1 + (-1)\alpha).
\]
Recalling \#21, these representations take care of the cases when $x = 0$, $x = 1$.  
But then matters follow if $x$ is any rational.
\vspace{0.2cm}

\qquad \un{Step 2:} \ 
Take $x$ irrational and introduce
\[
\alpha 
\ = \ 
\sum\limits_{j=1}^\infty \hsx \alpha_j 2^{-j}, 
\qquad 
\beta
\ = \ 
\sum\limits_{j=1}^\infty \hsx \beta_j 2^{-j}, 
\]
where for $k! \leq j < (k+1)!$, 
\[
\begin{cases}
\ \alpha_j = m_j \quad \quad \text{and} \quad \beta_j = 0 \qquad (k = 1, 3, 5, \ldots)\\
\ \alpha_j = \ 0 \hspace{.9cm} \text{and} \quad \beta_j = m_j   \hspace{.5cm} (k = 2, 4, 6, \ldots)
\end{cases}
.
\]
Then
\[
x 
\ = \ 
\alpha + \beta.
\]
\vspace{0.2cm}

\qquad \un{Step 3:} \ 
Assume that the series defining $\alpha$ is infinite $-$then in this case, $\alpha$ is a Liouville number.
\vspace{0.2cm}

[Break up the series
\[
\sum\limits_{j=1}^\infty \hsx
\alpha_j 2^{-j}
\]
as follows:
\allowdisplaybreaks
\begin{align*}
\sum\limits_{1! \leq j < 2!} \hsx
\alpha_j 2^{-j}
&+
\sum\limits_{2! \leq j < 3!} \hsx
\alpha_j 2^{-j}
+
\sum\limits_{3! \leq j < 4!} \hsx
\alpha_j 2^{-j}
+
\sum\limits_{4! \leq j < 5!} \hsx
\alpha_j 2^{-j}
+
\sum\limits_{5! \leq j < 6!} \hsx
\alpha_j 2^{-j} 
+ \cdots \
\\[12pt]
&=\ 
\sum\limits_{1! \leq j < 2!} \hsx
\alpha_j 2^{-j}
+
\sum\limits_{3! \leq j < 4!} \hsx
\alpha_j 2^{-j}
+
\sum\limits_{5! \leq j < 6!} \hsx
\alpha_j 2^{-j}
+ \cdots .
\end{align*}

Consider
\allowdisplaybreaks
\begin{align*}
0 
&<  \ 
\alpha - 
\sum\limits_{j=1}^{(2k)! - 1} \hsx \alpha_j 2^{-j}
\\[12pt]
&=\ 
\sum\limits_{j \geq (2k)!} \hsx 
\alpha_j 2^{-j}
\\[12pt]
&=\ 
\sum\limits_{(2k)! \leq j < (2k+1)!} \hsx
\alpha_j 2^{-j}
+
\sum\limits_{(2k+1)! \leq j < (2k+2)!} \hsx
\alpha_j 2^{-j}
\\[12pt]
&\hspace{4.5cm}
+
\sum\limits_{(2k+2)! \leq j < (2k+3)!} \hsx
\alpha_j 2^{-j}
+ \cdots 
\\[12pt]
&=\ 
0 
+ 
\sum\limits_{(2k+1)! \leq j < (2k+2)!} \hsx
\alpha_j 2^{-j}
+ 0 
+ \cdots 
\\[12pt]
&\leq\ 
\sum\limits_{j=(2k+1)!}^\infty \hsx
\alpha_j 2^{-j}
\\[12pt]
&=\ 
\frac{1}{2^{(2k+1)!}} \hsx 
\sum\limits_{j=0}^\infty \hsx
\frac{1}{2^j}
\\[12pt]
&=\ 
\frac{2}{2^{(2k+1)!}} 
\\[12pt]
&=\ 
2^{1 - (2k + 1)!}.
\end{align*}
Define a sequence of rationals $\ds\frac{p_k}{q_k}$ $(k = 1, 2, \ldots)$ by the prescription
\[
\frac{p_k}{q_k}
\ = \ 
\sum\limits_{j=1}^{(2k)! - 1} \hsx \alpha_j 2^{-j}, 
\qquad q_k \ = \ 2^{(2k)! - 1}.
\]
Then $p_k$ and $q_k$ are integers, $q_k > 1$, and 
\[
0 \ < \ \alpha - \frac{p_k}{q_k} \ < \ \frac{1}{q_k^k}.
\]
Therefore $\alpha$ is a Liouville number.]
\vspace{0.2cm}

[Note: \ 
Tacitly
\[
2^{1 - (2k+1)!} 
\ < \ 
2^{k - k(2k)!}.
\]
In fact, 
\allowdisplaybreaks
\begin{align*}
1 - (2k + 1)! + k(2k)!\ 
&=\ 
1 - (2k)!(2k+1) + k(2k)!
\\[12pt]
&=\ 
1 - (k + k)(2k)! - (2k)! + k(2k)!
\\[12pt]
&=\ 
1 - k(2k)! - k(2k)! - (2k)! + k(2k)!
\\[12pt]
&=\ 
1 - k(2k)! - (2k)! 
\\[12pt]
&<\ k.]
\end{align*}
\vspace{0.2cm}

\qquad \un{Step 4:} \ 
Assume that the series defining $\beta$  is infinite $-$then in this case, $\beta$ is a Liouville number.
\vspace{0.2cm}

\qquad \un{Step 5:} \ 
So if the series defining $\alpha$ and the series defining $\beta$ are infinite, we are done.
\vspace{0.2cm}

\qquad \un{Step 6:} \ 
If the series defining $\alpha$ is finite, then $\alpha$ is rational.  
If the series defining $\beta$  is infinite, then $\beta$  is a Liouville number, thus $x = \alpha + \beta$ is a Liouville number, 
thence $\ds\frac{x}{2}$ is a Liouville number and 
\[
x \ = \ \frac{x}{2} + \frac{x}{2}.
\]
\vspace{0.2cm}

\qquad \un{Step 7:} \ 
Reverse the roles of $\alpha$ and $\beta$  in the previous step.
\vspace{0.2cm}

\qquad \un{Step 8:} \ 
The case when both defining series are finite cannot occur (for then $\alpha$ and $\beta$ are rational, contradicting the assumption that 
$x = \alpha + \beta$ is irrational).
\vspace{0.5cm}

\begin{x}{\small\bf THEOREM} \ 
If $x$ is a Liouville number, then for any algebraic number $\alpha > 0$ $(\alpha \neq 1)$, the power $\alpha^x$ is transcendental.
\end{x}
\vspace{0.3cm}

It is a question of showing that $\alpha^x \neq \alpha^\prime$ for every algebraic \ $\alpha^\prime > 0$, \ i.e., that \ 
$\elln(\alpha^x) \neq \elln(\alpha^\prime)$, i.e., that $x \elln(\alpha) \neq \elln(\alpha^\prime)$, or still, that
\[
\abs{x \elln(\alpha) - \elln(\alpha^\prime)}
\ > \ 
0.
\]
If
\[
\frac{\elln(\alpha^\prime)}{\elln(\alpha)}
\]
were rational and if 
\[
\abs{x \elln(\alpha) - \elln(\alpha^\prime)} \ = \ 0,
\]
then it would follow that
\[
x \ = \ 
\frac{\elln(\alpha^\prime)}{\elln(\alpha)},
\]
which is impossible ($x$, being Liouville, is transcendental (cf. \#7)).  
So assume that
\[
\frac{\elln(\alpha^\prime)}{\elln(\alpha)}
\]
is irrational and write
\allowdisplaybreaks
\begin{align*}
\abs{x \elln(\alpha) - \elln(\alpha^\prime)} \ 
&=\ \abs{x \elln(\alpha) - \frac{p}{q}\elln(\alpha) + \frac{p}{q}\elln(\alpha) - \elln(\alpha^\prime)}
\\[12pt]
&=\ 
\abs{\bigg(x - \frac{p}{q}\bigg)\elln(\alpha) +\frac{p}{q}\elln(\alpha) - \elln(\alpha^\prime)}
\\[12pt]
&=\ 
\abs{\frac{p}{q}\elln(\alpha)  - \elln(\alpha^\prime) - - \bigg(x - \frac{p}{q}\bigg)\elln(\alpha)}
\\[15pt]
&\geq\ 
\abs{\frac{p}{q}\elln(\alpha)  - \elln(\alpha^\prime)} - 
\abs{-\bigg(x - \frac{p}{q}\bigg)\elln(\alpha)} 
\\[15pt]
&=\ 
\abs{\frac{p}{q}\elln(\alpha)  - \elln(\alpha^\prime)} - 
\abs{\bigl(x - \frac{p}{q}\bigr)\elln(\alpha)} 
\\[15pt]
&=\ 
\abs{\frac{p}{q}\elln(\alpha)  - \elln(\alpha^\prime)} - 
\abs{x - \frac{p}{q}} \hsx \abs{\elln(\alpha)}
\\[15pt]
&>\ 
\abs{\frac{p}{q}\elln(\alpha)  - \elln(\alpha^\prime)} - 
\frac{\abs{\elln(\alpha)}}{q^k}
\\[15pt]
&=\ 
\frac{\abs{p\elln(\alpha)  - q\elln(\alpha^\prime)}}{q} - 
\frac{\abs{\elln(\alpha)}}{q^k}
\\[15pt]
&=\ 
\frac{1}{q} \hsx
\bigg(
\abs{p\elln(\alpha)  - q\elln(\alpha^\prime)} - 
\frac{\abs{\elln(\alpha)}}{q^{k-1}}
\bigg)
\\[15pt]
\end{align*}
thereby reducing matters to the positivity of 
\[
\abs{p\elln(\alpha)  - q\elln(\alpha^\prime)} - 
\frac{\abs{\elln(\alpha)}}{q^{k-1}}.
\]
In any event, 
\[
\abs{p\elln(\alpha)  - q\elln(\alpha^\prime)} 
\]
is positive since otherwise
\[
\frac{p}{q}
\ = \ 
\frac{\elln(\alpha^\prime)}{\elln(\alpha)}
\]
contradicting the supposition that
\[
\frac{\elln(\alpha^\prime)}{\elln(\alpha)}
\]
is irrational.
\vspace{0.3cm}

\begin{x}{\small\bf LEMMA} \ 

\[
\abs{p \elln(\alpha) - q \elln(\alpha^\prime)} \ \geq \ 
\frac{1}{\max \{\abs{p},q\}^c},
\]
where $c > 0$ depends only on $\elln(\alpha)$ and $\elln(\alpha^\prime)$.
\vspace{0.2cm}

[This estimate will be established later on (cf. \S32, \#4).] 
\end{x}
\vspace{0.3cm}

Assume that $x \in [0,1]$, choose $k \gg 0$:
\[
\frac{\abs{\elln(\alpha)}}{q^{k-1-c}} 
\ < \ 
\frac{1}{2} \hsx\min\bigg\{ \bigg(\frac{2}{3}\bigg)^c, \frac{1}{2}\bigg\},
\]
and take $\abs{p} \neq 0$, hence 
\[
-\frac{q}{2} \ < \  p \ < \  \frac{3 q}{2} \qquad \text{(cf. \#22)}
\]
\qquad\qquad $\implies$
\[
0 < \abs{p} < \frac{3q}{2} 
\]
\qquad\qquad $\implies$
\[
\frac{1}{\abs{p}} > \frac{2}{3q}.
\]
There are now two possibilities:
\[
\abs{p \elln(\alpha) - q \elln(\alpha^\prime)} \ \geq \ 
\begin{cases}
\ \quad \text{$\ds\frac{1}{q^c}$}\\[11pt]
\ \quad \text{$\ds\frac{1}{\abs{p}^c}$}
\end{cases}
.
\]

\qquad \textbullet \quad
Work with $\ds\frac{1}{q^c}$ $-$then the issue is the positivity of 
\[
\frac{1}{q^c} - \frac{\abs{\elln(\alpha)}}{q^{k-1}}
\]
or still, the positivity of
\[
1 - \frac{\abs{\elln(\alpha)}}{q^{k-1-c}} 
\ > \ 
1 - \frac{1}{2} \cdot \frac{1}{2} 
\ = \ 
\frac{3}{4} 
\ > \ 0.
\]
\vspace{0.2cm}

\qquad \textbullet \quad 
Work with $\ds\frac{1}{\abs{p}^c}$ $-$then the issue is the positivity of 
\[
\frac{1}{\abs{p}^c}- \frac{\abs{\elln(\alpha)}}{q^{k-1}}
\]
or still, the positivity of
\[
\bigg( \frac{2}{3q}\bigg)^c - \frac{\abs{\elln(\alpha)}}{q^{k-1}}\ 
=\ 
\bigg( \frac{2}{3}\bigg)^c \hsx \frac{1}{q^c} \hsx - \frac{\abs{\elln(\alpha)}}{q^{k-1}}
\]
or still, the positivity of
\[
\bigg( \frac{2}{3}\bigg)^c  - 
\frac{\abs{\elln(\alpha)}}{q^{k-1-c}}
\ > \ 
\bigg( \frac{2}{3}\bigg)^c - 
\frac{1}{2} \hsx \bigg( \frac{2}{3}\bigg)^c 
\ = \ 
\frac{1}{2} \hsx \bigg(\frac{2}{3}\bigg)^c
\ > \ 
0.
\]
\vspace{0.2cm}

\vspace{0.3cm}

\begin{x}{\small\bf REMARK} \ 
Take $\alpha$ as above and assume that $x$ is positive $-$then 
\[
\elln(x \alpha) 
\quad \text{and} \quad 
x \elln(\alpha)
\]
are transcendental.
\end{x}
\vspace{0.3cm}

%% file: _16_the_mahler_classification.tex
\chapter{
$\boldsymbol{\S}$\textbf{16}.\quad  THE MAHLER CLASSIFICATION}
\setlength\parindent{2em}
\setcounter{theoremn}{0}
\renewcommand{\thepage}{\S16-\arabic{page}}

\ \indent 
What follows is a proofless summary of the relevant facts.
\vspace{0.3cm}

\begin{x}{\small\bf DEFINITION} \ 
Let $P(X) \in \C[X]$, say
\[
P(X) \ = \ a_0 + a_1 X + \cdots + a_n X^n.
\]
Then the 
\un{height}
\index{height}
of $P(X)$, denoted $\tH(P)$, is 
\[
\max\{\abs{a_0}, \abs{a_1}, \ldots, \abs{a_n}\}.
\]
\end{x}
\vspace{0.3cm}

\begin{x}{\small\bf NOTATION} \ 
Given a real number $x$, $w_n(x)$ $(n \in \N)$ is the supremum of the real numbers $w$ such that
\[
0 \ < \ \abs{P(x)} \ \leq \ \tH(P)^{-w}
\]
has infinitely many solutions $P(X) \in \Z[X]$ of degree at most $n$.
\end{x}
\vspace{0.3cm}

\begin{x}{\small\bf LEMMA} \ 
For any nonzero rational number $\ds\frac{a}{b}$, 
\[
w_n(x) \ = \ 
\begin{cases}
\ w_n\bigg(x + \ds\frac{a}{b}\bigg)\\[15pt]
\ w_n\bigg(\ds\frac{a}{b} x\bigg)
\end{cases}
.
\]
\end{x}
\vspace{0.3cm}

\begin{x}{\small\bf LEMMA} \ 
For any positive integer $n$, 
\[
0 \ \leq \ w_n(x) \ \leq \ \infty.
\]
\end{x}
\vspace{0.3cm}

\begin{x}{\small\bf \un{N.B.}}\ 
The sequence $\{w_n(x)\}$ is increasing: \ 
$w_1(x) \leq w_2(x) \leq \ldots$ and $w_n(x) \geq n$.  
\end{x}
\vspace{0.3cm}

\begin{x}{\small\bf MAIN PROBLEM} \ 
Suppose that $\{w_n\}$ is an increasing sequence of real numbers with $w_n \geq n$ $\forall \ n \in \N$.  
Does there exist a real number $x$ such that for all $n$, $w_n(x) = w_n$?
\end{x}
\vspace{0.3cm}

\begin{x}{\small\bf NOTATION}  \ 
Put
\[
w(x) \ = \ 
\lim\limits_{n \ra \infty} \hsx \sup \frac{w_n(x)}{n}.
\]
Therefore
\[
0 \ \leq \ w(x) \ \leq \ \infty.
\]
\vspace{0.2cm}

[Note: \ 
Real numbers with $0 < w(x) < 1$ do not exist.]
\end{x}
\vspace{0.3cm}

\begin{x}{\small\bf DEFINITION} \ 
A real number $x$ is an 
\vspace{0.2cm}

\qquad \textbullet \quad \un{\mA-number} if $w(x) = 0$; 
\index{\mA-number}
\vspace{0.2cm}

\qquad \textbullet \quad \un{\mS-number} if $0 < w(x) < \infty$; 
\index{\mS-number}
\vspace{0.2cm}

\qquad \textbullet \quad \un{\mT-number} if $w(x) = \infty \ \& \ \forall \ \ n \geq 1, \ w_n(x) < \infty$; 
\index{\mT-number}
\vspace{0.2cm}

\qquad \textbullet \quad \un{\mU-number} if $w(x) = \infty \ \& \ \forall \ n \gg 1, \ w_n(x) = \infty$. 
\index{\mU-number}
\vspace{0.2cm}
\end{x}
\vspace{0.3cm}

Write \mA, \mS, \mT, \mU for the corresponding sets (termed 
\un{Mahler classes}) 
\index{Mahler classes}
$-$then
\[
\R \ = \ 
A
\hsx \cup \hsx
S
\hsx \cup \hsx
T
\hsx \cup \hsx
U,
\]
a disjoint union.
\vspace{0.2cm}

[Note: \ 
The transcendentals $\T$ decompose as
\[
S
\hsx \cup \hsx
T
\hsx \cup \hsx
U\hsx .]
\]
\begin{x}{\small\bf THEOREM} \ 
The \mA-numbers are exactly the real algebraic numbers.
\end{x}
\vspace{0.3cm}

\begin{x}{\small\bf THEOREM} \ 
The Mahler classes \mS, \mT, \mU  are not empty.
\end{x}
\vspace{0.3cm}

\begin{x}{\small\bf REMARK} \ 
$A (=\Qbar)$ is a set of measure 0 (being countable).  
It can be shown that \mT and \mU are also sets of measure 0, hence almost all real numbers are \mS-numbers.
\end{x}
\vspace{0.3cm}

\begin{x}{\small\bf EXAMPLE} \ 
Suppose that $\alpha$ is a nozero algebraic number $-$then $e^\alpha$ is an \mS-number, thus in particular, $e$ is an \mS-number.
\end{x}
\vspace{0.3cm}

\begin{x}{\small\bf EXAMPLE} \ 
For any positive integer $d$, 
\[
\sum\limits_{j=1}^\infty \hsx 2^{-(d+1)j}
\]
is an \mS-number.
\end{x}
\vspace{0.3cm}

\begin{x}{\small\bf EXAMPLE} \ 
$\pi$ is not a \mU-number, so, being transcendental, is either an \mS-number or a \mT-number but no one knows which one.
\end{x}
\vspace{0.3cm}

\begin{x}{\small\bf \un{N.B.}} \ 
Exhibiting explicit \mT-numbers is complicated business.
\end{x}
\vspace{0.3cm}

\begin{x}{\small\bf DEFINITION} \ 
A \mU-number $x$ is a 
\un{$U_n$-number}
\index{$U_n$-number}
if $n$ is the smallest positive integer such that $w_n(x) = \infty$.
\end{x}
\vspace{0.3cm}

Write $U_n$ for the set of such.
\vspace{0.3cm}

\begin{x}{\small\bf THEOREM} \ 
Each $U_n$ is nonempty and
\[
U \ = \ \bigcup\limits_{n=1}^\infty \hsx U_n,
\]
a disjoint union.
\end{x}
\vspace{0.3cm}

\begin{x}{\small\bf EXAMPLE} \ 
$\forall \ n \in \N$, 
\[
\sqrt[\leftroot{-3}\uproot{3}n \hsx]{3/2} \hsx \cdot \hsx
\sum\limits_{j=1}^\infty \hsx 10^{-j!}
\]
is a $U_n$-number.
\end{x}
\vspace{0.3cm}

\begin{x}{\small\bf EXAMPLE} \ 
Let $m_j \in \{2, 4\}$ $(j = 1, 2, \ldots)$.  Put
\[
x \ = \ \big(3 + \hsx \sum\limits_{j=1}^\infty \hsx m_j \hsx 10^{-j!} \big) / 4.
\]
Then for all $n \geq 1$, the positive real $n^\nth$ root of $x$ is a $U_n$-number.
\end{x}
\vspace{0.3cm}

\begin{x}{\small\bf SCHOLIUM} \ 
$\forall \ n \geq 1$, $U_n$ is uncountable.
\end{x}
\vspace{0.3cm}

\begin{x}{\small\bf \un{N.B.}} \ 
$U_1 \ = \ \bL.$
\end{x}
\vspace{0.3cm}

\begin{x}{\small\bf DEFINITION} \ 
Two real numbers $x$ and $y$ are 
\un{algebraically dependent}
\index{algebraically dependent} 
if there is a nonzero polynomial $P(X,Y) \in \Z[X,Y]$ such that 
$P(x,y) = 0$ (cf. \S20, \#1).  
\vspace{0.2cm}

[Note: \ 
The denial is
\un{algebraically independent}.]
\index{algebraically independent} 
\end{x}
\vspace{0.3cm}

\begin{x}{\small\bf THEOREM} \ 
Algebraically dependent real numbers belong to the same Mahler class.
\end{x}
\vspace{0.3cm}

\begin{x}{\small\bf EXAMPLE} \ 
If $x$ is a \mU-number and $y$ is not a \mU-number, then $x$ and $y$ are algebraically independent.  
So, e.g., 
$\ds\sum\limits_{j=1}^\infty \hsx 10^{-j!}$ and $\pi$ are algebraically independent. 
\vspace{0.2cm}

[Note: \ 
$x + y$ is transcendental: \ Given
\[
\sum\limits_{j=0}^n \hsx
a_j (x + y)^j 
\ = \ 
0,
\]
consider
\[
P(X,Y) \ = \ 
\sum\limits_{j=0}^n \hsx a_j(X + Y)^j\hsx.]
\]
\end{x}
\vspace{0.3cm}

\begin{x}{\small\bf REMARK} \ 
In general, if $x$ and $y$ are transcendental numbers, then at least one of $x + y$ and $x y$ must be transcendental (cf. \S2, \#29).
\vspace{0.2cm}

[To see this, consider the polynomial
\[
X^2 - (x + y) X + xy.
\]
Its zeros are $x$ and $y$.  
So if both $x + y$ and $xy$ were algebraic, then $x$ and $y$ would be algebraic which they are not.]
\end{x}
\vspace{0.3cm}

\begin{x}{\small\bf EXAMPLE} \ 
It can be shown that the numbers $\pi$ and $e^\pi$ are algebraically independent but it is not known whether $e^\pi$ is or is not a \mU-number (recall that $\pi$ is not a \mU-number (cf. \#14)).
\end{x}
\vspace{0.3cm}


%% file: _17_transcendence_of_e.tex
\chapter{
$\boldsymbol{\S}$\textbf{17}.\quad  TRANSCENDENCE OF $\boldsymbol{e}$}
\setlength\parindent{2em}
\setcounter{theoremn}{0}
\renewcommand{\thepage}{\S17-\arabic{page}}

\ \indent 
We have seen that $e$ is irrational (cf. \S8, \#3) but more is true.
\vspace{0.3cm}

\begin{x}{\small\bf THEOREM} \ 
$e$ is transcendental.
\end{x}
\vspace{0.3cm}

\begin{x}{\small\bf SCHOLIUM} \ 
$\forall \ n \in \N$, $1, e, e^2, \ldots, e^n$ are linearly independent over $\Q$ (cf. \S8, \#11).
\end{x}
\vspace{0.3cm}

\begin{x}{\small\bf LEMMA} \ 
Given $f \in \R[X]$ of degree \mM, 
\[
e^x \hsx 
\int\limits_0^x \hsx 
f(t) \hsx e^{-t} \hsx dt 
\ = \ 
F(0) \hsx e^x 
\ - \ 
F(x),
\]
where
\[
F(x) 
\ = \
 \sum\limits_{\ell = 0}^M \hsx f^{(\ell)} (x).
\]
\vspace{0.2cm}

PROOF \ 
Integrate by parts to get
\[
\int\limits_0^x \hsx 
f(t) \hsx e^{-t} \hsx dt 
\ = \ 
f(0) - f(x) e^{-x}  + \int\limits_0^x \hsx f^\prime (t) e^{-t} dt.
\]
Then integrate this.
\vspace{0.2cm}

[Note: \ 
If $f$ has integer coefficients, then the same is true of \mF.]
\end{x}
\vspace{0.3cm}

Consider now a relation of the form
\[
a_0 + a_1 e + a_2 e^2 + \cdots + a_m e^m \ = \ 0,
\]
where $a_0 > 0$, $a_m \neq 0$ $(a_k \in \ \Z)$ $-$then from \#3, 
\[
F(0) e^k - F(k) \ = \ 
e^k \hsx \int\limits_0^k \hsx f(t) e^{-t} dt \qquad (k = 0, 1, \ldots, m),
\]
so
\[
F(0) \hsx 
\sum\limits_{k=0}^m \hsx 
a_k e^k - 
\sum\limits_{k=0}^m \hsx 
a_k F(k) 
\ = \ 
\sum\limits_{k=0}^m \hsx 
a_k e^k 
\int\limits_0^k \hsx
f(t) e^{-t} dt
\]
or still, 
\[
-\sum\limits_{k=0}^m \hsx 
a_k F(k) 
\ = \ 
\sum\limits_{k=0}^m \hsx 
a_k e^k 
\int\limits_0^k \hsx
f(t) e^{-t} dt,
\]
i.e., 
\[
-a_0 F(0) - 
\sum\limits_{k=1}^m \hsx 
a_k F(k) 
\ = \ 
\sum\limits_{k=0}^m \hsx 
a_k e^k 
\int\limits_0^k \hsx
f(t) e^{-t} dt.
\]

The polynomial $f$ is at our disposal and the trick is to choose it appropriately in order to reach a contradiction.  
One choice is to put
\[
g(X) 
\ = \ 
X^{n-1} (X - 1)^n \cdots (X - m)^n
\]
and let 
\[
f(X) 
\ = \ 
\frac{g(X)}{(n - 1)!},
\]
$n \in \N$ to be determined in due course.
\vspace{0.5cm}

{\small\bf FACTS}
\allowdisplaybreaks
\begin{align*}
& \deg f \ = \ (m+1) n - 1 \equiv M,\\[12pt]
& f^{(\ell)} (0) \ = \ 0 \qquad (0 \leq \ell \leq n - 2), \\[12pt]
& f^{(n - 1)} (0) \ = \ (-1)^{m n} \hsx (m!)^n, \\[12pt]
& n | f^{(\ell)} (0) \qquad (\forall \ \ell \neq n - 1).
\end{align*}
\vspace{0.2cm}

[Write
\allowdisplaybreaks
\begin{align*}
f(X) \ 
&=\ 
\frac{g(X)}{(n - 1)!} 
\\[12pt]
&= \ 
\frac{X^{n-1}}{(n - 1)!}  \hsx \bigl(b_0 + b_1 X + \cdots + b_{m n}X^{m n} \bigr)
\\[12pt]
&=\ 
\frac{1}{(n - 1)!}  \hsx  \bigl(b_0  X^{n-1}+ b_1 X^n + \cdots + b_{m n}X^{(m+1) n - 1} \bigr)
\\[12pt]
&=\ 
\frac{1}{(n - 1)!}  \hsx \sum\limits_{\ell = n - 1}^M \hsx c_\ell X^\ell  \qquad (c_{n-1} = b_0, \ c_n = b_1, \ldots).
\end{align*}
Then 
\[
\ell < n - 1 \implies f^{(\ell)} (0) = 0.
\]
And
\allowdisplaybreaks
\begin{align*}
\ell \geq n - 1 
&\implies \frac{f^{(\ell)}(0)}{\ell!} \ = \ \frac{c_\ell}{(n - 1)!}
\\[12pt]
&\implies  f^{(\ell)}(0) \ = \ \ell! \hsx \frac{c_\ell}{(n - 1)!} \in \Z.
\end{align*}
Therefore
\[
\ell \geq n \implies n | f^{(\ell)} (0)
\]
but
\allowdisplaybreaks
\begin{align*}
\ell \ = \ n - 1 \implies f^{(n - 1)} (0)
&=\ c_{n - 1} 
\\
&=\ b_0 
\\
&=\ (-1)^{m n} \hsx (m!)^n.]
\end{align*}
Consequently 
\allowdisplaybreaks
\begin{align*}
F(0) \ 
&=\ 
\sum\limits_{\ell=0}^M \hsx 
f^{(\ell)} (0)
\\[15pt]
&=\ 
\sum\limits_{\ell=n-1}^M \hsx 
f^{(\ell)} (0)
\\[15pt]
&=\ 
f^{(n-1)} (0) + f^{(n)} (0) + \cdots + f^{( (m + 1) n - 1)} (0)
\\[15pt]
&=\ 
 (-1)^{m n} \hsx (m!)^n + n C,
\end{align*}
C an integer.
\vspace{0.3cm}

The next step is to get a handle on the $F(k)$ $(1 \leq k \leq m)$.  
To this end, let
\allowdisplaybreaks
\begin{align*}
g_k(X) \ 
&=\ 
\frac{g(X)}{(X - k)^n}
\\[15pt]
&=\ 
X^{n - 1} \hsx 
\prod\limits_{\substack{\ell = 1\\ \ell \neq k}}^m \hsx 
(X - \ell)^m,
\end{align*}
a polynomial with integral coefficients.  
Using now the formula for differentiating a product, 
\[
g^{(j)}(X) 
\ = \ 
\sum\limits_{i=0}^j \hsx 
\binom{j}{i} ((X - k)^n)^{(i)} (g_k(X))^{(j-1)}.
\]
Due to the presence of the factor $X - k$, it follows that  
\[
g^{(j)} (k) 
\ = \ 
0 \qquad (j < n).
\]
On the other hand, if $j \geq n$, then 
\[
g^{(j)} (k) 
\ = \ 
\binom{j}{n} \hsx n! \hsx g_k^{(j-n)} (k).
\]
So, for all $j$, $g^{(j)}(k)$ is an integer divisible by $n!$, say
\[
g^{(j)} (k) 
\ = \ 
n! \hsx n_j (k).
\]
And then
\allowdisplaybreaks
\begin{align*}
F(k) \ 
&=\ 
\sum\limits_{\ell = 0}^M \hsx 
f^{(\ell)}(k)
\\[15pt]
&=\ 
\sum\limits_{\ell = n}^M \hsx 
f^{(\ell)}(k)
\\[15pt]
&=\ 
\sum\limits_{\ell = n}^M \hsx 
\frac{g^{(\ell)}(k)}{(n-1)!}
\\[15pt]
&=\ 
\sum\limits_{\ell = n}^M \hsx 
\frac{n! n_\ell(k)}{(n-1)!}
\\[15pt]
&=\ 
n \hsx \sum\limits_{\ell = n}^M \hsx n_\ell(k)
\\[15pt]
&=\ 
n n_k \qquad (n_k \in \Z).
\end{align*}
\vspace{0.1cm}

Take $n \gg 0$ ($n$ prime):
\[
 n > a_0 \quad \text{and} \quad \gcd(n,m!) = 1,
\]
hence
\[
n \ \not| \ a_0 F(0) \qquad \text{(cf. \S7, \#1)}.
\]
And this implies that
\allowdisplaybreaks
\begin{align*}
-a_0 F(0) - \sum\limits_{k = 1}^m \hsx a_k F(k) \ 
&=\ 
-a_0 F(0) -
\sum\limits_{k = 1}^m \hsx
a_k(n n_k) 
\\[15pt]
&=\ 
-a_0 F(0) - n \hsx
\bigg(
\sum\limits_{k = 1}^m \hsx
a_k n_k \bigg)
\\[15pt]
&\neq\ 
0.
\end{align*}
\vspace{0.1cm}

To recapitulate:
\[
-a_0 F(0) \hsx - \hsx \sum\limits_{k = 1}^m \hsx a_k F(k) 
\]
is a nonzero integer, thus
\[
\abs{\sum\limits_{k = 0}^m \hsx a_k F(k)} \ \geq \  1.
\]
Return now to 
\[
\sum\limits_{k = 0}^m \hsx a_k e^k \hsx 
\int\limits_0^k \hsx 
f(t) e^{-t} \hsx dt,
\]
an entity that depends on $n$ and which can be made arbitrarily small (leading thereby to the sought for contradiction).
\vspace{0.1cm}

To see this, note that
\[
\abs{f(x)} 
\ \leq \ 
\frac{m^M}{(n-1)!} \qquad (0 \leq x \leq m) \ (M = (m + 1) n - 1),
\]
so
\allowdisplaybreaks
\begin{align*}
\abs{\sum\limits_{k = 0}^m \hsx a_k e^k \hsx \int\limits_0^k \hsx f(t) e^{-t} \hsx dt} \ 
&\leq\ 
\frac{m^M}{(n-1)!} \hsx 
\sum\limits_{k = 0}^m \hsx
\abs{a_k} \hsx \int\limits_0^k \hsx e^{k - t} \hsx dt
\\[15pt]
&\leq\ 
\frac{m^{(m+1) n}}{(n-1)!} \hsx 
\sum\limits_{k = 0}^m \hsx
\abs{a_k} \hsx (e^k - 1)
\\[15pt]
&\leq\ 
\frac{m^{(m+1) n}}{(n-1)!} \hsx 
\sum\limits_{k = 0}^m \hsx
\abs{a_k} \hsx e^k
\\[15pt]
&\leq\ 
\frac{m^{(m+1) n}}{(n-1)!} \  
e^m \ 
\sum\limits_{k = 0}^m \hsx
\abs{a_k} \hsx
\\[15pt]
&=\ 
\frac{C^n}{(n-1)!}  \hsx
e^m \hsx
\sum\limits_{k = 0}^m \hsx
\abs{a_k} \hsx,
\end{align*}
where
\[
C 
\ = \ m^{m + 1}.
\]
But
\allowdisplaybreaks
\begin{align*}
\frac{C^n}{(n-1)!} \ 
&=\ 
C \cdot 
\frac{C^{n-1}}{(n-1)!}
\\[12pt]
&\ra 0 \qquad (n \ra \infty) \quad \text{(cf. \S0)}.
\end{align*}
\vspace{0.5cm}

Here is an application of \#1.
\vspace{0.5cm}

\begin{x}{\small\bf SCHOLIUM} \ 
Let $q$ be a nonzero rational number $-$then $e^q$ is transcendental (cf. \S9, \#1).
\vspace{0.2cm}

[Take $q > 0$ and suppose that $e^q$ is algebraic.  
Write $q = \ds\frac{a}{b}$ $(a, b > 0)$ $-$then $\ds\big( e^{\frac{a}{b}} \big)^b = e^a$ is algebraic, which implies that $e$ is algebraic (cf. \S2, \#37), a contradiction.]
\end{x}
\newpage
\[
\text{APPENDIX}
\]
\vspace{0.5cm}

Consider the transcendence status of the three examples figuring in the Appendix to \S8.  
\vspace{0.2cm}

\qquad \textbullet \quad 
Is the number
\[
\sum\limits_{k = 0}^\infty \hsx \frac{r^k}{2^{k(k - 1)/2}}
\]
transcendental? \ Ans: \ Unknown.
\vspace{0.2cm}

\qquad \textbullet \quad 
Is the number
\[
\sum\limits_{k = 0}^\infty \hsx r^{2^k}
\]
transcendental? \ Ans: \ Yes. 
\vspace{0.2cm}

\qquad \textbullet \quad 
Is the number
\[
\sum\limits_{k = 1}^\infty \hsx \frac{1}{M^{k^2}}
\]
transcendental? \ Ans: \  Yes.

%% file: _18_symmetric_algebra.tex
\chapter{
$\boldsymbol{\S}$\textbf{18}.\quad  SYMMETRIC ALGEBRA}
\setlength\parindent{2em}
\setcounter{theoremn}{0}
\renewcommand{\thepage}{\S18-\arabic{page}}


\begin{x}{\small\bf RAPPEL} \ 
Let \mA be a commutative ring with unit $-$then a polynomial 
\[
P(X_1, \ldots, X_n) \in A[X_1, \ldots, X_n]
\]
is 
\un{symmetric}
\index{symmetric polynomial}
if for any permutation $\sigma$ of $\{1, \ldots, n\}$, 
\[
P(X_{\sigma(1)}, \ldots, X_{\sigma(n)}) 
\ = \ 
P(X_1, \ldots, X_n).
\]
\end{x}
\vspace{0.3cm}

\begin{x}{\small\bf DEFINITION}  \ 
The 
\un{elementary symmetric polynomials}
\index{elementary symmetric polynomials} 
$s_1, s_2, \ldots, s_n$ in $n$ variables $x_1, x_2, \ldots, x_n$ appear as coefficients in the monic polynomial of degree $n$ and roots 
$x_1, x_2, \ldots, x_n$:
\[
(X -x_1) (X - x_2) \ldots (X - x_n) 
\ = \ 
X^n - s_1 X^{n-1} + \cdots + (-1)^n s_n.
\]
\vspace{0.2cm}

Explicated:
\begin{align*}
s_1 \ 
&=\ 
x_1 + x_2 + \cdots + x_n
\\[12pt]
s_2 \ 
&=\ 
x_1 x_2 + x_1 x_3 + \cdots + 
x_2 x_3 + x_2 x_4 + \cdots +x_{n-1} x_n
\\[12pt]
\vdots
\\[12pt]
s_n \ 
&=\ 
x_1 x_2 \cdots x_n.
\end{align*}
\end{x}
\vspace{0.3cm}

\begin{x}{\small\bf THEOREM} \ 
Every symmetric polynomial can be written as a polynomial in the elementary symmetric polynomials: \ 
If $P \in A[X_1, \ldots, X_n]$ is symmetric, then there exists a polynomial $F \in A[s_1, \ldots, s_n]$ such that 
\[
P 
\ = \ 
F(s_1, \ldots, s_n).
\]
E.g.:
\allowdisplaybreaks
\begin{align*}
P(X_1, X_2) \ 
&=\ 
3(X_1 X_2)^3 - ((X_1 + X_2)^2 - 2X_1 X_2)
\\[12pt]
&=\ 
3s_2^2 - s_1^2 - 2 s_2
\\[12pt]
&\equiv\ 
F(s_1, s_2).
\end{align*}
\end{x}
\vspace{0.3cm}

\begin{x}{\small\bf LEMMA} \ 
Let $\alpha$ be an algebraic number, let $d = \deg \alpha$ $(\equiv d(\alpha))$, let
$\alpha_1, \ldots, \alpha_d$ $(\alpha = \alpha_1)$ be the zeros of $f_\alpha$ (cf. \S14, \#7), and let
\[
F 
\ = \ 
F(X;\alpha_1, \ldots, \alpha_d) \in \Q[X;\alpha_1, \ldots, \alpha_d].
\]
Assume: \  As a polynomial in $\alpha_1, \ldots, \alpha_d$ with coefficients in $\Q[X]$, \mF is symmetric $-$then
\[
F 
\ = \ 
F(X) \in \Q[X].
\]

PROOF \ 
Write
\allowdisplaybreaks
\begin{align*}
f_\alpha(z) \ 
&=\ 
a_0 + a_1 z + \cdots + a_d z^d \qquad (a_0, a_1, \ldots, a_d \in \Z)
\\[12pt]
&=\ 
a_d (z - \alpha_1) (z - \alpha_2) \cdots (z - \alpha_d)
\\[12pt]
&=\ 
a_d (z^d - (\alpha_1 + a_2 + \cdots + \alpha_d) z^{d-1} 
\\[12pt]
& \qquad\qquad
+ (\alpha_1 \alpha_2 + \alpha_1 \alpha_3 + \cdots + \alpha_{d-1} \alpha_d)z^{d-2} 
\\[12pt]
& \qquad\qquad\qquad\qquad
+ \cdots + (-1)^d (\alpha_1 \alpha_2 \cdots \alpha_d)),
\end{align*}
from which
\allowdisplaybreaks
\begin{align*}
s_1 \ 
&=\ 
\alpha_1 + \alpha_2 + \cdots + \alpha_d 
\ = \ 
-\frac{a_{d-1}}{a_d}
\\[12pt]
s_2 \ 
&=\ 
\alpha_1 \alpha_2 + \alpha_1 \alpha_3 + \cdots + \alpha_{d-1} \alpha_d
\ = \ 
\frac{a_{d-2}}{a_d}
\\[12pt]
\vdots
\\[12pt]
s_d \ 
&=\ 
\alpha_1 \alpha_2 \cdots \alpha_d
\ = \ 
(-1)^d \hsx \frac{a_0}{a_d},
\end{align*}
implying thereby that the elementary symmetric polynomials in the $\alpha_1, \alpha_2, \ldots, \alpha_d$ are rational numbers.
Turning now to \mF, being a symmetric polynomial in $\alpha_1, \alpha_2, \ldots, \alpha_d$, it can be written as a polynomial in the elementary symmetric polynomials $s_1, s_2, \ldots, s_d$ with coefficients in $\Q[X]$.  
But $s_1, s_2, \ldots, s_d \in \Q[X]$, 
hence
\[
F 
\ = \ 
F(X) \in \Q[X].
\]
\end{x}
\vspace{0.3cm}

\begin{x}{\small\bf \un{N.B.}} \ 
Suppose that $\alpha$ is an algebraic integer and let
\[
F 
\ = \ 
F(X;\alpha_1, \ldots, \alpha_d) \in \Z[X;\alpha_1, \ldots, \alpha_d].
\]
Assume: \ 
As a polynomial in $\alpha_1, \ldots, \alpha_d$ with coefficeints in $\Z[X]$, \mF is symmetric $-$then
\[
F 
\ = \ 
F(X) \in \Z[X].
\]
\end{x}
\vspace{0.3cm}


%% file: _19_transcendence_of_pi.tex
\chapter{
$\boldsymbol{\S}$\textbf{19}.\quad  THE TRANSCENDENCE OF $\boldsymbol{\pi}$}
\setlength\parindent{2em}
\setcounter{theoremn}{0}
\renewcommand{\thepage}{\S19-\arabic{page}}

\ \indent 
Here is the objective:
\vspace{0.3cm}

\begin{x}{\small\bf THEOREM} \ 
$\pi$ is transcendental.
\end{x}
\vspace{0.1cm}

Suppose that $\pi$ is algebraic $-$then $\alpha \equiv \pi \hsx \sqrt{-1}$ is algebraic.  
Agreeing to use the notation of \S18, \#4, in view of the relation $e^{\pi \hsx \sqrt{-1}} + 1 = 0$, it follows that 
\[
(1 + e^{\alpha_1}) \hsx (1 + e^{\alpha_2}) \hsx \cdots \hsx (1 + e^{\alpha_d}) 
\ = \ 
0
\]
or still, upon expanding the product, 
\[
\sum\limits_{\epsilon_1 = 0}^1 \hsx
\sum\limits_{\epsilon_2 = 0}^1 \hsx
\cdots
\sum\limits_{\epsilon_d = 0}^1 \hsx
e^{\epsilon_1 \alpha_1 + \epsilon_2 \alpha_2 + \cdots \epsilon_d \alpha_d }
\ = \ 
0.
\]
\vspace{0.3cm}

\begin{x}{\small\bf EXAMPLE}  \ 
Take $\epsilon_1 = 1, \ \epsilon_2 = \cdots = \epsilon_d = 0$ $-$then 
\[
\epsilon_1 \alpha_1 + \epsilon_2 \alpha_2 + \cdots + \epsilon_d \alpha_d \ \neq \ 0.
\]
Take $\epsilon_1 = \epsilon_2 = \cdots = \epsilon_d = 0$ $-$then 
\[
\epsilon_1 \alpha_1 + \epsilon_2 \alpha_2 + \cdots + \epsilon_d \alpha_d \ = \ 0.
\]
\end{x}
\vspace{0.3cm}

Denoting the exponents by $\beta_k$, rewrite matters in the form
\[
1 + \sum\limits_{k=1}^{2^d - 1} \hsx
e^{\beta_k} 
\ = \ 
0,
\]
where things have been arranged so that the nonzero $\beta_k$ are placed first:
\[
\beta_1 \neq 0, \ 
\beta_2 \neq 0, \ 
\ldots, 
\beta_r \neq 0, 
\  0, \ldots, 0.
\]

Put 
\[
A = 1 + (2^d - 1) - r.
\]
Then $A \geq 1$ and 
\[
1 + \sum\limits_{k=1}^{2^d - 1} \hsx
e^{\beta_k} 
\ = \ 
A + e^{\beta_1} + e^{\beta_2} + \cdots + e^{\beta_r} 
\ = \ 0.
\]
\vspace{0.3cm}

\begin{x}{\small\bf LEMMA} \ 
The nonzero numbers $\beta_1, \ldots, \beta_r$ are the set of roots of a polynomial $\phi(X) \in \Z[X]$ of degree $r$ (hence are algebraic).
\vspace{0.2cm}

PROOF \ 
Let 
\[
\psi(X) 
\ = \ 
\prod\limits_{\epsilon_1 = 0}^1 \hsx
\prod\limits_{\epsilon_2 = 0}^1 \hsx
\cdots
\prod\limits_{\epsilon_d = 0}^1 \hsx
(X - (\epsilon_1 \alpha_1 + \epsilon_2 \alpha_2 + \cdots + \epsilon_d \alpha_d)).
\]
Viewed as a polynomial in $\alpha_1, \alpha_2, \ldots, \alpha_d$ with coefficients in $\Q[X]$, it is symmetric.  
Therefore $\psi(X)$ is in $\Q[X]$ (cf. \S18, \#4).  
On the other hand, the roots of $\psi(X)$ are the $\beta_k$ $(1 \leq k \leq r)$ and 0 with multiplicity 
\mA ($r + A = r + 2^d - r = 2^d$, the degree of $\psi(X))$, thus the roots of the polynomial
\[
X^{-A} \psi(X)
\]
are $\beta_1, \ldots, \beta_r$.  
Denoting by $m$ the least common denominator of the coefficients of this polynomial, take
\allowdisplaybreaks
\begin{align*}
\phi(X) \ 
&=\ 
m X^{-A} \psi(X)
\\[12pt]
&=\ 
C_r X^r + \cdots + C_1 X + C_0
\\[12pt]
&\in \Z[X] \quad (C_r > 0, \ C_0 \neq 0).
\end{align*}
\end{x}
\vspace{0.3cm}

\begin{x}{\small\bf RAPPEL} \ 
Given $f \in \R[X]$ of degree \mM, 
\[
e^x \hsx
\int\limits_0^x \hsx 
f(t) e^{-t} \hsx dt 
\ = \ 
F(0) e^x - F(x) 
\qquad \text{(cf. \S17, \#3)}.
\]
\vspace{0.2cm}

[Note: \ 
Complex $x$ are admitted in which case the integral $\ds\int\limits_0^x$ 
is calculated along the line segment joining 0 and $x$.]
\vspace{0.3cm}

Feed into this relation $x = \beta_1, \ldots, x = \beta_r$ to get:
\[
\begin{cases}
\ e^{\beta_1} \hsx \ds\int\limits_0^{\beta_1} \hsx f(t) e^{-t} \hsx dt \ = \ F(0) \hsx e^{\beta_1} - F(\beta_1)\\[15pt]
\hspace{3.15cm} \vdots\\[12pt]
\ e^{\beta_r} \hsx \ds\int\limits_0^{\beta_r} \hsx f(t) e^{-t} \hsx dt \ = \ F(0) \hsx e^{\beta_r} - F(\beta_r)
\end{cases}
.
\]
But
\[
A + e^{\beta_1} + \cdots + e^{\beta_r} \ = \ 0.
\]
Therefore
\[
-A F(0) \ - \ 
\sum\limits_{k=1}^r \hsx 
F(\beta_k) 
\ = \ 
\sum\limits_{k=1}^r \hsx 
e^{\beta_k} \hsx
\int\limits_0^{\beta_k} \hsx f(t) e^{-t} \hsx dt.
\]

Just as in the proof of the transcendence of $e$, the modus operandi at this juncture is to choose $f$ judiciously so as to bring about a contradiction.  
To this end, let
\[
f(X) 
\ = \ 
\frac{1}{(n-1)!} \hsx
(C_r)^{n r - 1} X^{n - 1} (\phi(X))^n
\]
or still, 
\allowdisplaybreaks
\begin{align*}
f(X) \ 
&=\ 
\frac{1}{(n-1)!} \hsx
(C_r)^{n r - 1} X^{n - 1} (C_r(X - \beta_1) \cdots(X - \beta_r))^n
\\[15pt]
&=\ 
\frac{1}{(n-1)!} \hsx
(C_r)^{n(r + 1) - 1} X^{n - 1} (X - \beta_1)^n \cdots(X - \beta_r)^n,
\end{align*}
$n \in \N$ a ``large'' natural number to be held in abeyance for the moment.
\vspace{0.5cm}

\qquad {\small\bf FACTS}
\allowdisplaybreaks
\begin{align*}
&\deg f \ = \ n (r + 1) - 1 \ \equiv \ M, 
\\[12pt]
&
f^{(\ell)} (0)  \ = \ 0 \quad (0 \leq \ell \leq n - 2),
\\[12pt]
&
f^{n-1} (0) \ = \ (C_r)^{n r - 1} C_0^n, 
\\[12pt]
&
n | f^{(\ell)} (0) \quad (\forall \ \ell \neq n - 1).
\end{align*}
Consequently
\allowdisplaybreaks
\begin{align*}
F(0)\ 
&=\ 
\sum\limits_{\ell=0}^M \hsx
f^{(\ell)} (0)
\\[12pt]
&=\ 
\sum\limits_{\ell=n-1}^M \hsx
f^{(\ell)} (0)
\\[12pt]
&=\ 
f^{n-1} (0) + f^{n} (0) + \cdots + f^{(n(r+1) - 1)} (0)
\\[12pt]
&=\ 
(C_r)^{n r - 1} C_0^n + nC,
\end{align*}
\mC an integer.

Moving on, from the definitions, 
\[
F(\beta_k) 
\ = \ 
\sum\limits_{\ell=0}^M \hsx
f^{(\ell)} (\beta_k).
\]
And $\beta_k$ is a root of $f(X)$ of multiplicity $n$, thus
\[
f^{(\ell)} (\beta_k) 
\ = \ 
0 
\qquad (0 \leq \ell \leq n - 1, 1 \leq k \leq r),
\]
leaving
\[
F(\beta_k) 
\ = \ 
\sum\limits_{\ell=n}^M \hsx
f^{(\ell)} (\beta_k).
\]
\end{x}
\vspace{0.3cm}
\begin{x}{\small\bf LEMMA} \ 
If $p(X) \in \Z[X]$, then $\forall \ \ell \in \N$, all the coefficients of the $\ell^{\nth}$ derivative 
$p^{(\ell)}(X)$ are divisible by $\ell !$.
\vspace{0.2cm}

PROOF \ 
Since differentiation is a linear operation, it suffices to check this on the powers $X^k$, 
restricting ourselves to when $1 \leq \ell \leq k$, 
in which case the $\ell^\nth$ derviative of $X^k$ is equal to
\[
\ell ! \hsx \binom{k}{\ell} X^{k - \ell}
\]
and the binomial coefficient $\ds\binom{k}{\ell}$ is a positive integer.
\vspace{0.2cm}

It therefore follows that for $\ell \geq n$, the coefficients of $f^{(\ell)}(X)$ are  integers divisible by $n (C_r)^{n r - 1}$.
\vspace{0.2cm}

[In detail, the polynomial
\[
X^{n-1} (\phi(X))^n \in \Z[X] \qquad \text{(cf. \#3)}
\]
and its $\ell^\nth$ derviative has all coefficients divisible by $\ell !$, so for $\ell \geq n$, its $\ell^\nth$ 
derviative has all coefficients divisible by $n!$ ($\ell ! = n!$ $(n + 1) \ldots \ell$).  
If $\ell \geq n$ and if generically, $n! \hsx W$ $(W \in \Z$) is a coefficient of 
\[
(X^{n-1} (\phi(X))^n)^{(\ell)} , 
\]
then
\[
\frac{1}{(n - 1)!} \hsx 
(C_r)^{n r - 1} \hsx n! \hsx W 
\ = \ 
n (C_r)^{n r - 1} W
\]
is a coefficient of $f^{(\ell)}(X)$.]
\end{x}
\vspace{0.3cm}

\begin{x}{\small\bf LEMMA} \ 
Let $P(X_1, \ldots, X_r)$ be a polynomial with integer coefficients of degree $s \leq t$ symmetric in the $X_k$ $-$then
\[
C_r^t \hsx P(\beta_1, \ldots, \beta_r)
\]
is an integer.
\vspace{0.2cm}

PROOF \ 
The algebraic numbers $C_r \beta_1, \ldots, C_r \beta_r$ are the roots of the monic polynomial
\[
(C_r)^{r - 1} \phi\bigg(\frac{X}{C_r}\bigg)
\ = \ 
X^r + C_{r-1}X^{r-1} + C_r C_{r-2}X^{r-2} + \cdots + C_r^{r-1}C_0,
\]
thus the elementary symmetric polynomials per $C_r \beta_1, \ldots, C_r \beta_r$ are integers, since
\[
s_1  \ = \ -\frac{C_{r-1}}{1}, \ 
s_2  \ = \ \frac{C_r C_{r-2}}{1}, 
\ldots, 
s_r  \ = \ (-1)^r \frac{C_r^{r-1} C_0}{1}.
\]
If $p(X_1, \ldots, X_r)$ is a homogeneous symmetric polynomial of degree $s \leq t$ with integer coefficients, then 
\[
C_r^s \hsx p(\beta_1, \ldots, \beta_r) 
\ = \ 
p(C_r \beta_1, \ldots, C_r \beta_r).
\]
But the right hand side can be written as a polynomial with integer coefficients in the elementary symmetric polynomials per 
$C_r \beta_1, \ldots, C_r \beta_r$, hence 
\[
C_r^s \hsx p(\beta_1, \ldots, \beta_r) 
\]
is an integer, hence a fortiori
\[
C_r^t \hsx p(\beta_1, \ldots, \beta_r) 
\]
is an integer.  
To treat the general case, simply separate the polynomial \mP into a sum of homogeneous polynomials $p$.
\end{x}
\vspace{0.3cm}

Fix $\ell$: $n \leq \ell \leq M$ and pass to 
\[
\sum\limits_{k=1}^r \hsx
f^{(\ell)} (\beta_k)
\]
or still, in suggestive notation,
\[
n (C_r)^{n r - 1} \hsx 
\sum\limits_{k = 1}^r \hsx 
g_\ell (\beta_k).
\]
\vspace{0.3cm}

\begin{x}{\small\bf \un{N.B.}} \ 
The degree of $f^{(n)}$ is 
\[
M - n 
\ = \ 
(n (r + 1) - 1) - n 
\ = \ n r - 1,
\]
so the degree of $f^{(\ell)}$ $(n \leq \ell \leq M)$ is $\leq nr - 1$.
\vspace{0.2cm}

Applying \#6 to
\[
(C_r)^{n r - 1} \hsx 
\sum\limits_{k = 1}^r \hsx 
g_\ell (\beta_k),
\]
legal since the sum is symmetric in the $\beta_k$, we conclude that
\[
\sum\limits_{k = 1}^r \hsx 
f^{(\ell)} (\beta_k) 
\ = \ 
n N_\ell,
\]
$N_\ell$ an integer.

Therefore
\allowdisplaybreaks
\begin{align*}
\sum\limits_{k = 1}^r \hsx F(\beta_k) \ 
&=\ 
\sum\limits_{k = 1}^r \hsx
\sum\limits_{\ell = n}^M \hsx 
f^{(\ell)}(\beta_k)
\\[15pt]
&=\ 
\sum\limits_{\ell = n}^M \hsx 
\sum\limits_{k = 1}^r \hsx
f^{(\ell)}(\beta_k)
\\[15pt]
&=\ 
n \hsx 
\sum\limits_{\ell = n}^M \hsx 
N_\ell
\\[15pt]
&\equiv n B.
\end{align*}
Now assemble what has been established thus far:
\allowdisplaybreaks
\begin{align*}
A F(0) + \sum\limits_{k = 1}^r \hsx F(\beta_k) \
&=\ 
A((C_r)^{n r - 1} C_0^n + n C) + nB
\\[12pt]
&=\ 
A(C_r)^{n r - 1} C_0^n + n(A C + B).
\end{align*}
Choose $n \gg 0$ ($n$ prime):
\[
n > A \quad \& \quad \gcd(n, C_r C_0) \ = \ 1.
\]
Then
\[
A (C_r)^{n r - 1} C_0^n + n(A C + B)
\]
is an integer not divisible by $n$, hence in particular is nonzero, hence
\[
\abs{A F(0) \ + \ \sum\limits_{k = 1}^r \hsx F(\beta_k)} \ \geq \ 1.
\]

It remains to estimate
\[
\sum\limits_{k = 1}^r \hsx
e^{\beta_k} \hsx 
\int\limits_0^{\beta_k} \hsx
f(t) \hsx e^{-t} \hsx dt.
\]
Suppose that
\[
\abs{\beta_k} \leq R \qquad (k = 1, \ldots, r)
\]
and put
\[
T 
\ = \ 
\max\limits_{\abs{z} \leq R} \hsx
\abs{(C_r)^r \phi(z)} \qquad (C_r \geq 1 \implies \frac{1}{C_r} \leq 1).
\]
Then
\allowdisplaybreaks
\begin{align*}
\max\limits_{\abs{z} \leq R} \hsx \abs{f(z)} \ 
&\leq\ 
\max\limits_{\abs{z} \leq R} \hsx
\frac{1}{(n-1)!}\hsx
\abs{(C_r)^{nr - 1} z^{n-1} (\phi(z))^n}
\\[15pt]
&\leq\ 
\max\limits_{\abs{z} \leq R} \hsx
\frac{1}{(n-1)!}\hsx
\abs{z}^{n-1}
\frac{1}{C_r} \hsx 
\abs{(C_r^r)^n (\phi(z))^n}
\\[15pt]
&\leq\ 
\frac{R^{n-1}}{(n-1)!}\hsx
\max\limits_{\abs{z} \leq R} \hsx
\abs{(C_r^r)^n (\phi(z))^n}
\\[15pt]
&\leq\ 
\frac{R^{n-1} T^n}{(n-1)!}\hsx.
\end{align*}
Consequently, for all $n$ per supra
\allowdisplaybreaks
\begin{align*}
\abs{\sum\limits_{k = 1}^r \hsx e^{\beta_k} \hsx \int_0^{\beta_k} \hsx f(t) \hsx e^{-t} \hsx dt}\ 
&\leq\ 
\sum\limits_{k = 1}^r \hsx
\abs{e^{\beta_k} \hsx \int_0^{\beta_k} \hsx f(t) \hsx e^{-t} \hsx dt}
\\[15pt]
&\leq\ 
\sum\limits_{k = 1}^r \hsx
\abs{\int_0^{\beta_k} \hsx\abs{f(t)}\hsx \abs{e^{(\beta_k - t)}}\hsx dt}
\\[15pt]
&\leq\ 
\frac{R^{n-1} T^n}{(n-1)!}\ 
\sum\limits_{k = 1}^r \hsx
\abs{\int_0^{\beta_k} \abs{e^{(\beta_k - t)}}\hsx dt}
\\[15pt]
&\leq\ 
\frac{R^{n-1} T^n}{(n-1)!}\ 
e^R \ 
\sum\limits_{k = 1}^r \hsx
\abs{\int_0^{\beta_k} \hsx dt}
\\[15pt]
&\leq\ 
\frac{R^{n-1} T^n}{(n-1)!} \ 
e^R (r R)
\\[15pt]
&=\ 
r e^R \  
\frac{(R T)^n}{(n-1)!}\hsx
\\[15pt]
&=\ 
r e^R \   (RT) \hsx \frac{(R T)^{n-1}}{(n-1)!},
\end{align*}
which leads to a contradiction in the usual way (cf. \S0).
\end{x}
\vspace{0.3cm}


%% file: _20_algebraic_in_dependence.tex
\chapter{
$\boldsymbol{\S}$\textbf{20}.\quad  ALGEBRAIC (IN)DEPENDENCE }
\setlength\parindent{2em}
\setcounter{theoremn}{0}
\renewcommand{\thepage}{\S20-\arabic{page}}


\begin{x}{\small\bf TERMINOLOGY} \ 
Let $\LL$ be a field, $\K \subset \LL$ a subfield.
\vspace{0.3cm}

\qquad \textbullet \quad 
A finite subset $S = \{\alpha_1, \ldots, \alpha_n\} \subset \LL$ is 
\un{algebraically dependent over $\K$}
\index{algebraically dependent over $\K$} 
if there is a nonzero polynomial $P \in \K[X_1, \ldots, X_n]$ such that 
\[
P(\alpha_1, \ldots, \alpha_n) \ = \ 0.
\]

\qquad \textbullet \quad 
A finite subset $S = \{\alpha_1, \ldots, \alpha_n\} \subset \LL$ is 
\un{algebraically independent over $\K$}
\index{algebraically dependent over $\K$} 
if there is no nonzero polynomial $P \in \K[X_1, \ldots, X_n]$ such that 
\[
P(\alpha_1, \ldots, \alpha_n) \ = \ 0.
\]
\end{x}
\vspace{0.1cm}

\begin{x}{\small\bf \un{N.B.}} \ 
Take $S = \{\alpha\}$, a one element set $-$then by definition, $\alpha$ is 
\un{algebraic over $\K$} 
\index{algebraic over $\K$} 
if \mS is algebraically dependent over $\K$
and $\alpha$ is 
\un{transcendental over $\K$} 
\index{transcendental over $\K$} 
if \mS is algebraically independent over $\K$ 
i.e., $\alpha \in S$ is algebraic or transcendental over $\K$ according to whether it is or is not a root of a polynomial in $\K[X]$ (cf. \S2, \#25).
\end{x}
\vspace{0.3cm}

\begin{x}{\small\bf LEMMA} \ 
Suppose that \mS is algebraically independent over $\K$ $-$then the elements of \mS are transcendental over $\K$.
\end{x}
\vspace{0.3cm}

The setup for us is when 
\[
\LL \ = \ \C \quad \text{and} \quad \K \ = \ \Q,
\]
in which case one can work either with polynomials \mP in $\Q[X_1, \ldots, X_n]$ or in $\Z[X_1, \ldots, X_n]$.
\vspace{0.2cm}

[Note: \ 
Here, of course, ``algebraic'' means algebraic over $\Q$ and ``transcendental'' means transcendental over $\Q$ and to say that the complex numbers 
$x_1, \ldots, x_n$ are 
algebraically dependent or algebraically independent means that the set $\{x_1, \ldots, x_n\}$ is algebraically dependent over $\Q$ or algebraically independent over $\Q$.]
\vspace{0.3cm}

\begin{x}{\small\bf REMARK} \ 
A complex number $x$ is transcendental iff the numbers $1,\ x, \ x^2, \ldots$ are linearly independent over $\Q$.  
And, in general, the complex numbers $x_1, \ldots, x_n$ are algebraically independent over $\Q$ iff the powers
\[
x_1^{k_1} \cdots x_n^{k_n} \qquad (k_i \in \Z,k_i \geq 0)
\]
are linearly independent over $\Q$.
\end{x}
\vspace{0.3cm}

\begin{x}{\small\bf LEMMA} \ 
Suppose that $S \subset \C$ is algebraically independent over $\Q$ $-$then the elements of \mS are transcendental over $\Q$ (cf. \#3).
\vspace{0.2cm}

[Note: \ 
If any of the elements in \mS are algebraic over $\Q$, then \mS is algebraically dependent over $\Q$.]
\end{x}
\vspace{0.3cm}

\begin{x}{\small\bf REMARK} \ 
It can happen that all the elements of \mS are transcendental over $\Q$, yet \mS is not algebraically independent over $\Q$.
\vspace{0.2cm}

[The real numbers $\sqrt{\pi}$ and $2 \pi + 1$ are transcendental but $\{\sqrt{\pi}, 2 \pi + 1\}$ is not algebraically independent over $\Q$.  
Thus consider
\[
P(X,Y) 
\ = \ 
2 X^2 - Y + 1.
\]
Then
\[
P(\sqrt{\pi}, 2 \pi + 1) 
\ = \ 
0.]
\]
\end{x}
\vspace{0.3cm}

\begin{x}{\small\bf LEMMA} \ 
If $\{x_1, \ldots, x_n\}$ is algebraically independent over $\Q$, then $\{x_1, \ldots, x_n\}$ is algebraically independent over $\ov{\Q}$ and for any nonconstant polynomial 
$P \in \ov{\Q}[X_1, \ldots, X_n]$, the number
\[
P(x_1, \ldots, x_n)
\]
is transcendental.
\end{x}
\vspace{0.3cm}

\begin{x}{\small\bf EXAMPLE} \ 
The numbers $e^{\sqrt{2}}$, $e^{3\sqrt{2}}$ are algebraically dependent.
\vspace{0.2cm}

[Consider
\[
P(X_1, X_2) 
\ = \ 
X_1^3 - X_2.
\]
Then
\allowdisplaybreaks
\begin{align*}
P(e^{\sqrt{2}},e^{3\sqrt{2}}) \ 
&=\ 
\bigl(e^{\sqrt{2}}\bigr)^3 - e^{3\sqrt{2}}
\\[12pt]
&=\ 
e^{3\sqrt{2}} - e^{3\sqrt{2}}
\\[12pt]
&=\ 
0.]
\end{align*}
\end{x}
\vspace{0.3cm}

\begin{x}{\small\bf EXAMPLE} \ 
Let $a$ and $b$ be relatively prime natural numbers $> 1$ $-$then the Liouville numbers (cf. \S15, \#9)
\[
x 
\ = \ 
\sum\limits_{j=1}^\infty \hsx
\frac{1}{a^{j !}} 
\quad \text{and} \quad 
y
\ = \ 
\sum\limits_{j=1}^\infty \hsx
\frac{1}{b^{j !}} 
\]
are algebraically independent over $\Q$.
\end{x}
\vspace{0.3cm}

\begin{x}{\small\bf EXAMPLE} \ 
Nesterenko proved in 1996 that the numbers $\pi$, $e^\pi$ are algebraically independent over $\Q$.
\end{x}
\vspace{0.3cm}

\begin{x}{\small\bf REMARK} \ 
The question of whether the numbers $e$, $\pi$ are algebraically dependent over $\Q$ or algebraically independent over $\Q$ is open.
\end{x}
\vspace{0.3cm}

In addition to numbers, one must also deal with functions.
\vspace{0.3cm}

\begin{x}{\small\bf DEFINITION} \ 
A function $f(z)$ of a complex variable $z$ is 
\un{algebraic}
\index{algebraic function}
if there is a nonzero polynomial $P \in \C[X,Y]$ such that $\forall \ z$, 
\[
P(z, f(z)) 
\ = \ 0.
\]
\end{x}
\vspace{0.3cm}

\begin{x}{\small\bf THEOREM} \ 
An entire function is algebraic iff it is a polynomial.
\end{x}
\vspace{0.3cm}

\begin{x}{\small\bf DEFINITION} \ 
An \hfill entire \hfill function \hfill which \hfill is \hfill not \hfill algebraic \hfill is \hfill said \hfill to \hfill be\\ 
\un{transcendental}.
\index{transcendental}
\end{x}
\vspace{0.3cm}

\begin{x}{\small\bf EXAMPLE} \ 
$e^z$,  $\cos z$,  $\sin z$  are  transcendental,  as  is  the  function
\[
z \ra \ds\int\limits_0^z \hsx e^{-t^2} \hsx dt.
\]
\end{x}
\vspace{0.3cm}

\begin{x}{\small\bf DEFINITION} \ 
\ A collection of entire functions $f_1, \ldots, f_n$ is said to be \ 
\un{algebraically dependent over $\C$}
\index{algebraically dependent over $\C$} 
if there is a nonzero polynomial $P \in \C[X_1, \ldots, X_n]$ such that $P(f_1, \ldots, f_n)$ is the zero function.
\end{x}
\vspace{0.3cm}

\begin{x}{\small\bf DEFINITION} \ 
\ A collection of entire functions $f_1, \ldots, f_n$ is said to be \ 
\un{algebraically independent over $\C$}
\index{algebraically independent over $\C$} 
if for any nonzero polynomial $P \in \C[X_1, \ldots, X_n]$, the function $P(f_1, \ldots, f_n)$ is not the zero function.
\end{x}
\vspace{0.3cm}

\begin{x}{\small\bf EXAMPLE} \ 
Let $I(z) = z$ be the identity function $-$then an entire function $f$ is algebraic (transcendental) iff $I$ and $f$ are algebraically dependent (independent) over $\C$.
\end{x}
\vspace{0.3cm}

\begin{x}{\small\bf EXAMPLE} \ 
$\sin z$ and $\cos z$ are algebraically dependent over $\C$.
\vspace{0.2cm}

[Consider
\[
P(X,Y) 
\ = \ X^2 + Y^2 - 1.
\]
Then
\[
P(\sin z, \cos z) 
\ = \ 
(\sin z)^2 + (\cos z)^2 - 1 
\ = \ 
1 - 1 
\ = \ 0.]
\]
\end{x}
\vspace{0.3cm}

\begin{x}{\small\bf EXAMPLE} \ 
Take
\[
f_1(z) 
\ = \ 
e^z, 
\quad 
f_2(z) 
\ = \ 
e^{\frac{1}{2} z}.
\]
Then the functions $f_1$, $f_2$ are algebraically dependent over $\C$.
\vspace{0.2cm}

[Consider 
\[
P(X_1, X_2) 
\ = \ 
X_2^6 - X_1 X_2^4 + X_1^2 X_2^2 - X_1^3.
\]
Then
\allowdisplaybreaks
\begin{align*}
P\big(e^z, e^{\frac{1}{2} z}\big) \
&=\ 
e^{3z} - e^{3z} + e^{3z} - e^{3z}
\\
&=\ 0.]
\end{align*}
\end{x}
\vspace{0.3cm}

\begin{x}{\small\bf EXAMPLE} \ 
Take
\[
f_1(z) 
\ = \ e^z, 
\quad
f_2(z) 
\ = \ e^{\sqrt{-1} \hsx z}.
\qquad 
\]
Then the functions $f_1$, $f_2$ are algebraically independent over $\C$ (cf. \#26 infra).

\end{x}
\vspace{0.3cm}

\begin{x}{\small\bf EXAMPLE} \ 
The functions $1, z, z^2, \ldots, z^n$ are linearly independent over $\C$ and the functions 
$z, e^z, e^{z^2}, \ldots, e^{z^n}$ are algebraically independent over $\C$.
\end{x}
\vspace{0.3cm}

\begin{x}{\small\bf LEMMA} \ 
Let $\lambda_1, \ldots, \lambda_n$ be distinct complex numbers $-$then the entire functions
\[
e^{\lambda_1 z}, \ldots, e^{\lambda_n z} 
\]
are linearly independent over $\C(z)$.
\vspace{0.3cm}

PROOF \ 
The case $n = 1$ is trivial.  
Proceed from here by induction, assuming that the statement is true at level $n - 1$ $(n > 1)$ and consider the dependence relation
\[
F_1 e^{\lambda_1 z} + \cdots + F_n e^{\lambda_n z} 
\ = \ 
0,
\]
where $F_1, \ldots, F_n$ are nonzero elements of $\C(z)$, the objective being to derive a contradiction from this.  
Divide by $F_n$:
\allowdisplaybreaks
\begin{align*}
\frac{F_1}{F_n} e^{\lambda_1 z} + \cdots + \frac{F_n}{F_n}  e^{\lambda_n z} \ 
&\equiv \ 
G_1 e^{\lambda_1 z} + \cdots + 1  e^{\lambda_n z} \qquad (G_n = 1)
\\[15pt]
&=\ 
0
\end{align*}
or still, 
\[
e^{\lambda_n z} \bigl( G_1 e^{(\lambda_1 - \lambda_n) z} + \cdots + 1  e^{0 z} \bigr) \ = \ 0
\]
or still, 
\[
G_1 e^{(\lambda_1 - \lambda_n) z} + \cdots + 1  e^{0 z} \ = \ 0
\]
or still, 
\[
G_1 e^{\sigma_1 z} + \cdots + 1 e^{\sigma_n z} 
\ = \ 
0,
\]
where
\[
\sigma_1 
\ = \ 
\lambda_1 - \lambda_n 
\ \neq \ 
0 , \ldots, \sigma_n 
\ = \ 0.
\]
Now differentiate:
\[
(G_1^\prime + \sigma_1 G_1) e^{\sigma_1 z} + \cdots + (G_{n-1}^\prime + \sigma_{n-1} G_{n-1}) e^{\sigma_{n-1} z} 
\ = \ 
0,
\]
thereby leading to a dependence relation at level $n-1$ with distinct exponents $\sigma_1, \ldots, \sigma_{n-1}$, so
\[
G_1^\prime + \sigma_1 G_1 
\ = \ 
0, \ldots, G_{n-1}^\prime + \sigma_{n-1} G_{n-1} 
\ = \ 
0.
\]
But each of these coefficients is nonzero, hence the purported dependence relation
\[
F_1 e^{\lambda_1 z} + \cdots +  F_n e^{\lambda_n z} 
\ = \ 
0
\]
has led to a contradiction.
\end{x}
\vspace{0.3cm}
\begin{x}{\small\bf APPLICATION} \ 
Let $\lambda_1, \ldots, \lambda_n$ be distinct complex numbers $-$then the entire function
\[
c_1 e^{\lambda_1 z} + \cdots + c_n e^{\lambda_n z} \qquad (c_1, \ldots, c_n \in \C)
\]
is not identically zero if the $c_i$ are not all zero.
\end{x}
\vspace{0.3cm}

\begin{x}{\small\bf LEMMA} \ 
Let $\lambda_1, \ldots, \lambda_n$ be distinct complex numbers which are linearly independent over $\Q$ $-$then the entire functions
\[
e^{\lambda_1 z}, \ldots, e^{\lambda_n z}
\]
are algebraically independent over $\C$.
\vspace{0.2cm}

PROOF
Let 
\[
P(X_1, \ldots, X_n) 
\ \in \ 
\C[X_1, \ldots, X_n]
\]
be a nonzero polynomial $-$then the claim is that 
\[
f(z) 
\ = \ 
P(e^{\lambda_1 z}, \ldots, e^{\lambda_n z})
\]
is not identically zero.  
To this end, write
\[
P(X_1, \ldots, X_n) 
\ = \ 
\sum\limits_{(k_1, \ldots, k_n)} \hsx 
a_{k_1, \ldots, k_n} \hsx X_1^{k_1} \cdots X_n^{k_n},
\]
where the $a_{k_1, \ldots, k_n} \in \C$ and not all of them are zero, thus
\[
f(z) 
\ = \ 
\sum\limits_{(k_1, \ldots, k_n)} \hsx 
a_{k_1, \ldots, k_n} \hsx
\exp((k_1 \lambda_1 + \cdots + k_n \lambda_n) z).  
\]
But, due to our assumption on $\lambda_1, \ldots, \lambda_n$, the complex numbers
\[
k_1 \lambda_1 + \cdots + k_n \lambda_n
\]
are distinct:
\[
k_1 \lambda_1 + \cdots + k_n \lambda_n 
\ = \ 
\ell_1 \lambda_1 + \cdots + \ell_n \lambda_n
\]
\qquad\qquad $\implies$
\[
(k_1 - \ell_1) \lambda_1 + \cdots + (k_n - \ell_n) \lambda_n 
\ = \ 
0
\]
\qquad\qquad $\implies$
\[
(k_1 - \ell_1)
\ = \ 
0
, \ldots, 
(k_n - \ell_n)
\ = \ 
0.
\]
To conclude that $f(z)$ is not identically zero, it remains only to quote \#24.
\end{x}
\vspace{0.3cm}

\begin{x}{\small\bf EXAMPLE} \ 
Take $\lambda_1 = 1$, $\lambda_2 = \beta \notin \Q$ $-$then $e^z$, $e^{\beta z}$ are algebraically independent over $\C$ 
(take $\beta = \sqrt{-1}$ to recover \#21).
\end{x}
\vspace{0.3cm}

%% file: _21_the_lindeman_weierstrass_theorem.tex
\chapter{
$\boldsymbol{\S}$\textbf{21}.\quad  THE LINDEMANN-WEIERSTRASS THEOREM}
\setlength\parindent{2em}
\setcounter{theoremn}{0}
\renewcommand{\thepage}{\S21-\arabic{page}}

\ \indent 
This is the following statement.
\vspace{0.3cm}

\begin{x}{\small\bf THEOREM} \ 
Let $\alpha_0, \alpha_1, \ldots, \alpha_t$ be distinct algebraic numbers 
$-$then 
$e^{\alpha_0}, e^{\alpha_1}, \ldots, e^{\alpha_t}$ are linearly independent over 
$\ov{\Q}$, i.e., if $b_0, b_1, \ldots, b_t$ are algebraic numbers not all zero, then
\[
b_0 e^{\alpha_0} + b_1 e^{\alpha_1} + \cdots + b_t e^{\alpha_t} \ \neq \ 0.
\]

[It is a corollary that
\[
b_0 e^{\alpha_0} + b_1 e^{\alpha_1} + \cdots + b_t e^{\alpha_t} 
\]
is a transcendental number provided 
$\alpha_i \neq 0$ $\forall \ i = 0, 1, \ldots, t$.  
For suppose it was algebraic, say
\[
b_0 e^{\alpha_0} + b_1 e^{\alpha_1} + \cdots + b_t e^{\alpha_t}   
\ = \ c \ (\in \Qbar) 
\ = \ 
c e^0.
\]
Then
\[
b_0 e^{\alpha_0} + b_1 e^{\alpha_1} + \cdots + b_t e^{\alpha_t}  - ce^0 = 0.
\]
But $\alpha_0, \alpha_1, \ldots, \alpha_t, 0$ 
are distinct,  from which the obvious contradiction.]
\vspace{0.2cm}

[Note: \ 
Some condition on the data is necessary as can be seen by taking 
\[
\alpha_0 = 0, \ \alpha_1 = 1, \ b_0 = 1, \ b_1 = 0.]
\]

\end{x}
\vspace{0.3cm}


\begin{x}{\small\bf \un{N.B.}} \ 
We are working here in the complex domain, hence $\sqrt{-1}$ is algebraic (consider $X^2 + 1 = 0$) and 
$\ov{\Q}$, computed in $\C$, is a field.
\end{x}
\vspace{0.3cm}

\begin{x}{\small\bf LEMMA} \ 
Suppose that $a$ and $b$ are real $-$then $a + \sqrt{-1} \hsx b$ is algebraic iff $a$ and $b$ are algebraic (cf. \S14, \#4).
\vspace{0.2cm}

PROOF \ 
If $a$ and $b$ are algebraic, then the combination $a + \sqrt{-1} \hsx b$ is algebraic ($\ov{\Q}$ being a field).  
Conversely, if $a + \sqrt{-1} \hsx b$ is algebraic, then $p(a + \sqrt{-1} \hsx b) = 0$, where $p(X)$ is a polynomial with rational 
coefficients, thus also $p(a - \sqrt{-1} \hsx b) = 0$.
Therefore
\[
\begin{cases}
\ (a + \sqrt{-1} \hsx b) + (a - \sqrt{-1} \hsx b) = 2 a \in \ov{\Q}\\[11pt]
\ (a + \sqrt{-1} \hsx b) - (a - \sqrt{-1} \hsx b) = 2 \sqrt{-1} \hsx b \in \ov{\Q}
\end{cases}
\implies
\begin{cases}
\ \ds\frac{1}{2} \hsx (2 a) = a \in \ov{\Q}\\[15pt]
\ -\ds\frac{\sqrt{-1}}{2} \bigl(2 \sqrt{-1} \hsx b \bigr) = b \in \ov{\Q}
\end{cases}
,
\]
i.e., $a$ and $b$ are algebraic.]
\end{x}
\vspace{0.3cm}

Before tackling the proof of the theorem, we shall consider some applications and examples.
\vspace{0.3cm}

\begin{x}{\small\bf LEMMA} \ 
If $\alpha$ is a nonzero algebraic number, then $e^\alpha$ is transcendental (Hermite-Lindemann).
\vspace{0.2cm}

[A nontrivial relation of the form
\[
q_0 + q_1 e^\alpha + \cdots + q_n e^{n \alpha} 
\ = \ 0 \qquad (q_k \in \Q)
\]
is impossible. 
Alternatively, if $e^\alpha$ were algebraic then take in \#1 
$\alpha_0 = 0$, 
$\alpha_1 = \alpha$, 
$b_0 = e^\alpha$, 
$b_1 = -1$
to get
\[
(e^\alpha) e^0 + (-1) e^\alpha
\ = \ 
0.
\]
Contradiction.]
\vspace{0.2cm}

[Note: \ 
Consequently, if $\alpha$ is a nonzero complex number, then at least one of the numbers $\alpha$ or $e^\alpha$ is transcendental.]
\end{x}
\vspace{0.3cm}

In particular: \ $e$ is transcendental (cf. \S17, \#1).  
And if $a, \ b \in \N$, then $e^a \neq b$.
\vspace{0.3cm}

\begin{x}{\small\bf EXAMPLE} \ 
$e^{\sqrt{2}}$ is transcendental.
\end{x}
\vspace{0.3cm}

\begin{x}{\small\bf EXAMPLE} \ 
$\pi$ is transcendental (cf. \S19, \#1).
\vspace{0.2cm}

[For if $\pi$ were algebraic, then $\pi \sqrt{-1}$ would be algebraic, hence $e^{\pi \sqrt{-1}}$ would be transcendental (cf. \#4), 
contrary to the fact that $1 + e^{\pi \sqrt{-1}} = 0$.]
\end{x}
\vspace{0.3cm}

\begin{x}{\small\bf EXAMPLE} \ 
Let $\alpha$ be a real nonzero algebraic number $-$then $\cos (\alpha)$ is transcendental (cf. \S12, \#1).
\vspace{0.2cm}

[Suppose instead that $\cos (\alpha) \equiv \beta$ was algebraic.  
Write
\[
\cos(\alpha) 
\ = \ 
\frac{e^{\alpha \sqrt{-1}} +e^{-\alpha\sqrt{-1}}}{2 \sqrt{-1}}
\ = \ 
\frac{e^{\alpha \sqrt{-1}}}{2 \sqrt{-1}} + \frac{e^{-\alpha \sqrt{-1}}}{2 \sqrt{-1}} 
\]
or still, 
\[
\bigg(-\frac{\sqrt{-1}}{2} \bigg) e^{\sqrt{-1} \hsx \alpha} 
+ 
\bigg(-\frac{\sqrt{-1}}{2} \bigg) e^{-\sqrt{-1} \hsx \alpha} 
+ 
(-\beta)e^0
\ = \ 
0,
\]
a contradiction (cf. \#1) ($\sqrt{-1} \hsx \alpha$ and $-\sqrt{-1} \hsx \alpha$ are obviously distinct).]
\vspace{0.2cm}

[Note: \ 
Consider the unique real fixed point of the cosine function, thus 
$\cos(x) = x = 0.739085\ldots$ $-$then $x$ is transcendental.  
For suppose that $x$ is algebraic $-$then $\cos(x)$ would be transcendental.  
But $\cos(x) = x$.]
\end{x}
\vspace{0.3cm}

The story for $\sin(\alpha)$ is analogous, as are the stories for 
\[
\begin{cases}
\ \cosh(\alpha)\\
\ \sinh(\alpha)
\end{cases}
.
\]
\vspace{0.3cm}
\begin{x}{\small\bf EXAMPLE} \ 
Let $\alpha$ be a real nonzero algebraic number $-$then $\tan(\alpha)$ is transcendental.
\vspace{0.2cm}

[Assuming the opposite, write
\[
\tan(\alpha) 
\ = \ 
\frac{e^{\alpha \sqrt{-1}} - e^{-\alpha\sqrt{-1}}} {\sqrt{-1} (e^{\alpha \sqrt{-1}} + e^{-\alpha\sqrt{-1}})}
\ \equiv \ 
\beta
\]
\qquad\qquad $\implies$
\[
(1 - \beta \sqrt{-1}) e^{\alpha \sqrt{-1}} - (1 + \beta\sqrt{-1}) e^{-\alpha \sqrt{-1}} 
\ = \ 0
\]
and note that $1 - \beta \sqrt{-1}$ and $1 + \beta \sqrt{-1}$ cannot simultaneously be zero.]
\end{x}
\vspace{0.3cm}

\begin{x}{\small\bf EXAMPLE} \ 
Let $\alpha \neq 1$ be a positive algebraic number $-$then $\elln(\alpha)$ is transcendental.
\vspace{0.2cm}

[If $\elln(\alpha)$ were algebraic, then $e^{\elln(\alpha)}$ would be transcendental (cf. \#4).  
But $e^{\elln(\alpha)} = \alpha \ldots$ .]
\end{x}
\vspace{0.3cm}

\begin{x}{\small\bf LEMMA} \ 
Let $\alpha$ be a nonreal algebraic number $-$then
\[
\begin{cases}
\ \Rex(e^\alpha)\\
\ \Imx(e^\alpha)
\end{cases}
\]
are transcendental.
\vspace{0.2cm}

PROOF \ 
Write $\alpha = a + \sqrt{-1} \hsx b$ $(b \neq 0)$ $-$then $a$ and $b$ are algebraic (cf. \#3).  Moreover, by definition, 
\[
e^\alpha 
\ = \ 
e^{a + \sqrt{-1} \hsx b}
\ = \ 
e^a (\cos b + \sqrt{-1} \hsx \sin b)
\]
and the claim is that
\[
\begin{cases}
\ e^a \cos b\\
\ e^a \sin b
\end{cases}
\]
are transcendental.  
To deal with the first of these, proceed by contradiction and assume that $e^a \cos b \equiv \beta$ is algebraic, thus 
$\beta \neq 0$ (the zeros of the cosine are transcendental).  
Next 
\allowdisplaybreaks
\begin{align*}
e^{a + \sqrt{-1} \hsx b} + e^{a - \sqrt{-1} \hsx b} \ 
&=\ 
e^a ( e^{\sqrt{-1} \hsx b} + e^{-\sqrt{-1} \hsx b}) 
\\[12pt]
&=\ 
e^a (\cos b + \sqrt{-1} \hsx \sin b + \cos (-b) + \sqrt{-1} \hsx \sin (-b))
\\[12pt]
&=\ 
2 e^a \cos b 
\\[12pt]
&=\ 
2 \beta
\end{align*}
\qquad\qquad $\implies$
\[
2 \beta e^0 - e^{a + \sqrt{-1} \hsx b} - e^{a - \sqrt{-1} \hsx b} 
\ = \ 
0.
\]
Owing to \#1, the algebraic numbers $0$, $a + \sqrt{-1} \hsx b$, $a - \sqrt{-1} \hsx b$ are not distinct, hence $b = 0$.  
On the other hand, $\alpha$ is not real, so $b \neq 0$.
\end{x}
\vspace{0.3cm}

\begin{x}{\small\bf \un{N.B.}} \ 
If in \#10, $\alpha$ was real, then matters are covered by \#4.
\end{x}
\vspace{0.3cm}

\begin{x}{\small\bf THEOREM} \ 
Suppose that $\beta_1, \ldots, \beta_r$ are nonzero algebraic numbers which are linearly independent over $\Q$ $-$then 
the transcendental numbers $e^{\beta_1}, \ldots, e^{\beta_r}$ are algebraically independent over $\Q$. 

\vspace{0.2cm}

PROOF \ 
Assume instead that for some nonzero polynomial
\[
P(X_1, \ldots, X_r) \in \Q[X_1, \ldots, X_r],
\]
say 
\[
P(X_1, \ldots, X_r) \
\ = \ 
\sum\limits_{k_1, \ldots, k_r} \hsx 
a_{k_1, \ldots, k_r} \hsx X_1^{k_1} \cdots X_r^{k_r},
\]
we have
\[
P(e^{\beta_1}, \ldots, e^{\beta_r}) 
\ = \ 0
\]
or still, 
\[
\sum\limits_{k_1, \ldots, k_r} \hsx 
a_{k_1, \ldots, k_r} \hsx e^{k_1 \beta_1 + \cdots + k_r \beta_r}
\ = \ 
0,
\]
where the $a_{k_1, \ldots, k_r} \in \Q$ and not all of them are zero.  
To settle the issue and arrive at a contradiction, it suffices to check that the exponents 
\[
k_1 \beta_1 + \cdots + k_r \beta_r
\]
are distinct (since then one can quote \#1).  
So suppose that
\[
(k_1, \ldots, k_r) \ \neq \ (\ell_1, \ldots, \ell_r)
\]
with 
\[
k_1 \beta_1 + \cdots + k_r \beta_r
\ = \ 
\ell_1 \beta_1 + \cdots + \ell_r \beta_r,
\] 
thus
\[
(k_1 - \ell_1)\beta_1 + \cdots + (k_r - \ell_r)\beta_r \ = \ 0,
\]
a nontrivial dependence relation over $\Q$.

\end{x}
\vspace{0.3cm}

\begin{x}{\small\bf EXAMPLE} \ 
The transcendental numbers $e$, $e^{\sqrt{2}}$  are algebraically independent over $\Q$. 
\vspace{0.2cm}

[For it is clear that the algebraic numbers 1, $\sqrt{2}$ are linearly independent over $\Q$.]
\end{x}
\vspace{0.3cm}

\begin{x}{\small\bf THEOREM} \ 
Suppose that $\beta_1, \ldots, \beta_r$ are nonzero algebraic numbers for which the transcendental numbers
 $e^{\beta_1}, \ldots, e^{\beta_r}$ are algebraically independent over $\Q$ $-$then $\beta_1, \ldots, \beta_r$ 
 are linearly independent over $\Q$.
\vspace{0.2cm}

PROOF \ 
Consider a nontrivial dependence relation over $\Q$:
\[
b_1 \beta_1 + \cdots + b_r \beta_r \ = \ 0.
\]
Clear the denominators and take the $b_k$ integral $-$then not all of them are zero and 
\[
1 \ = \ e^0 
\ = \ 
e^{b_1 \beta_1 + \cdots b_r \beta_r}.
\]
Define
\[
P(X_1, \ldots, X_r) \in \Q[X_1, \ldots, X_r]
\]
by the prescription
\[
P(X_1, \ldots, X_r)  
\ = \ 
X_1^{b_1} \cdots X_r^{b_r} - 1.
\]
Then
\allowdisplaybreaks
\begin{align*}
P(e^{\beta_1}, \ldots, e^{\beta_r}) \ 
&=\ 
e^{b_1 \beta_1} \cdots  e^{b_r \beta_r} - 1
\\[12pt]
&=\ 
e^{b_1 \beta_1 + \cdots + b_r \beta_r} - 1
\\[12pt]
&=\ 
1 - 1
\\[12pt]
&=\ 
0.
\end{align*}
But  $e^{\beta_1}, \ldots, e^{\beta_r}$ are algebraically independent over $\Q$.
Therefore
\[
P(X_1, \ldots, X_r) 
\ \equiv \ 
0 
\implies 
b_1 = 0, \ldots, b_r = 0,
\]
a contradiction.
\end{x}
\vspace{0.3cm}

\begin{x}{\small\bf SCHOLIUM} \ 
Nonzero algebraic numbers $\beta_1, \ldots, \beta_r$ are linearly independent over $\Q$ iff the transcendental numbers 
$e^{\beta_1}, \ldots, e^{\beta_r}$ algebraically independent over $\Q$.
\end{x}
\vspace{0.3cm}

\begin{x}{\small\bf LEMMA} \ 
Let $\alpha$ be an algebraic number whose real and imaginary parts are both nonzero $-$then the transcendental numbers 
$\Rex (e^\alpha)$, $\Imx(e^\alpha)$ are algebraically independent over $\Q$ (cf. \#10).
\end{x}
\vspace{0.3cm}

We need a preliminary.

\begin{x}{\small\bf SUBLEMMA} \ 
Let $x$ and $y$ be nonzero real numbers $-$then $x$ and $y$ are algebraically dependent over $\Q$ iff 
$x + \sqrt{-1} \hsx y$ and $x - \sqrt{-1} \hsx y$ are algebraically dependent over $\Q$.
\vspace{0.2cm}

PROOF \  
To deal with one direction, assume that there exists a nonzero polynomial
\[
P(X,Y) \ = \ 
\sum\limits_{m, n} \hsx 
a_{m n} \hsx X^m Y^n \ \in \Q[X,Y]
\]
such that
\[
P(x,y) \ = \ 0.
\]
Let
\[
\begin{cases}
\ \alpha = x + \sqrt{-1} \hsx y\\[11pt]
\ \ov{\alpha} = x - \sqrt{-1} \hsx y
\end{cases}
\implies
\begin{cases}
\ x = \ds\frac{\alpha + \ov{\alpha}}{2}\\[15pt]
\ y = \ds\frac{\alpha - \ov{\alpha}}{2 \sqrt{-1}}
\end{cases}
.
\]
Then
\[
\sum\limits_{m, n} \hsx 
a_{m n} \hsx 
\bigg(\frac{1}{2}\bigg)^{m + n} \bigl(-\sqrt{-1}\bigr)^n \hsx 
(\alpha + \ov{\alpha})^m (\alpha - \ov{\alpha})^n 
\ = \ 0.
\]
Introduce
\allowdisplaybreaks
\begin{align*}
Q(X,Y) \ 
&=\ 
\sum\limits_{m, n} \hsx 
a_{m n} \hsx \bigg(\frac{1}{2}\bigg)^{m + n} \bigl(-\sqrt{-1}\bigr)^n X^m Y^n
\\[15pt]
\ov{Q}(X,Y) 
&=\ 
\sum\limits_{m, n} \hsx 
a_{m n} \hsx \bigg(\frac{1}{2}\bigg)^{m + n} \bigl(\sqrt{-1}\bigr)^n X^m Y^n.
\end{align*}
Thus
\[
Q, \ \ov{Q} \in \C[X,Y]
\]
but
\[
Q \hsx \ov{Q} \in \Q[X,Y].
\]
Put now
\[
P^+(X,Y)  \ = \ 
Q(X + Y, X - Y) \hsx \ov{Q} (X + Y, X - Y).
\]
Then
\[
Q(\alpha + \ov{\alpha}, \alpha - \ov{\alpha}) \ = \ 0,
\]
so
\[
P^+ (\alpha, \ov{\alpha}) \ = \ 0,
\]
thereby establishing that $\alpha$ and $\ov{\alpha}$ are algebraically dependent over $\Q$.
\vspace{0.3cm}

Passing to the proof of \#16, write 
$\alpha = a + \sqrt{-1} \hsx b$ (thus $a \neq 0$, $b \neq 0$ are algebraic (cf. \#3)) $-$then 
$e^a \cos b$ and $e^a \sin b$ are algebraically dependent over $\Q$ iff 
\[
e^\alpha  \ = \ e^a \cos b + \sqrt{-1} \hsx e^a \sin b 
\quad \text{and} \quad 
e^{\ov{\alpha}}  \ = \ e^a \cos b - \sqrt{-1} \hsx e^a \sin b 
\]
are algebraically dependent over $\Q$ (cf. \#17), i.e., iff  $\alpha$ and $\ov{\alpha}$ are linearly dependent over $\Q$ 
(cf. \#15), i.e., iff $a = 0$ or $b = 0$, which cannot be.
\vspace{0.3cm}

We shall conclude this \S \ with an indication of the steps leading up to a proof of \#1.  
So let as there 
$b_0, b_1, \ldots, b_t$ be algebraic numbers not all zero but with 
\[
b_0 e^{\alpha_0} + b_1 e^{\alpha_1} + \cdots + b_t e^{\alpha_t}  
\ = \ 0.
\]

\qquad \un{Step 1:} \ 
By discarding terms whose coefficients are zero and rearranging the notation, it can be assumed that no coefficient is zero and 
\[
b_1 e^{\alpha_1} + \cdots + b_t e^{\alpha_t} 
\ = \ 
0.
\]
Consider the Taylor series expansion
\[
b_1 e^{\alpha_1 z} + \cdots + b_t e^{\alpha_t z} 
\ = \ 
\sum\limits_{n = 0}^\infty \hsx
\frac{u_n}{n!} \hsx z^n.
\]
\vspace{0.2cm}

\qquad \un{Step 2:} \  
$\forall \ n = 0, 1, \ldots$, 
\[
u_n \ = \ 
\sum\limits_{i = 1}^t \hsx
b_i \alpha_i^n.
\]
Define $a_1, \ldots, a_t$ by writing 
\[
(X - \alpha_1) \cdots (X - \alpha_t) 
\ = \ 
X^t - a_1 X^{t - 1} - \cdots - a_t.
\]
\vspace{0.2cm}

\qquad \un{Step 3:} \ 
$\forall \ n = 0, 1, \ldots$, 
\[
\alpha_i^{t+n} 
\ = \ 
a_1 \alpha_i^{t + n - 1} + \cdots + a_t \alpha_i^n \qquad (i = 1, \ldots, t).
\]
\vspace{0.2cm}

\qquad \un{Step 4:} \ 
$\forall \ n = 0, 1, \ldots$, 
\[
u_{n+t} 
\ = \ 
a_1 u_{n+t-1} + \cdots + a_t u_n.
\]
\vspace{0.2cm}

\qquad \un{Step 5:} \ 
It suffices to treat the case in which the $u_n \in \Q$ $(n = 0, 1, \ldots)$ and the $a_i \in \Q$ $(i = 1, \ldots, t)$.
\vspace{0.2cm}

[Consider the product
\[
\prod\limits_\sigma \hsx 
\big(
\sigma(b_1) e^{\sigma(\alpha_1) z} + \cdots + \sigma(b_t) e^{\sigma(\alpha_t) z}
\big),
\]
where
\[
\sigma \in \Gal(\Q(b_1, \ldots, b_t, \alpha_1, \ldots, \alpha_t) / \Q).
\]
This expression is still 0 (one of its factors is zero) and upon expanding has the form 
\[
\sum\limits_i \hsx 
b_i^\prime \hsx e^{\alpha_i^\prime z}.
\]
Since the sets $\{b_i^\prime\}$, $\{\alpha_i^\prime\}$ are Galois stable, the numbers $u_n^\prime$ and $a_i^\prime$ are rational.]
\vspace{0.5cm}

\qquad \un{Step 6:} \ 
Upon clearing denominators if necessary, it can be assumed that $u_0, \ldots, u_{t-1} \in \Z$, thus using Step 4 recursively, 
$\forall \ n \geq 0$, 
\[
d^n u_n \in \Z,
\]
where $d$ is a common denominator of the $a_i$ $(i = 1, \ldots, t)$.
\vspace{0.2cm}

[So, if $d = 1$, then the $u_n$ are integers.]
\vspace{0.5cm}

\qquad \un{Step 7:} \ 
Put
\[
A 
\ = \ 
\max\{1, \abs{\alpha_1}, \ldots, \abs{\alpha_t}\}.
\]
Then there exists a positive constant \mC such that $\forall \ n \geq 0$, 
\[
\abs{u_n} 
\ \leq \ 
C A^n \qquad \text{(use Step 2).}
\]

Recall now that the assumption is that
\[
b_1 e^{\alpha_1} + \cdots + b_t e^{\alpha_t} 
\ = \ 
0,
\]
hence 
\[
\sum\limits_{n = 0}^\infty \hsx 
\frac{u_n}{n!} 
\ = \ 
0.
\]

Given $k \in \N$, put
\[
v_k
\ = \ 
k! \hsx 
\sum\limits_{n = 0}^k \hsx 
\frac{u_n}{n!}  \qquad (v_0 \equiv u_0).
\]
\vspace{0.5cm}

\qquad \un{Step 8:} \ 
$\forall \ k$: $A < k + 1$, 

\allowdisplaybreaks
\begin{align*}
\abs{v_k} \ 
&=\ 
k! \hsx
\abs{\sum\limits_{n=0}^k \hsx \frac{u_n}{n!}}
\\[15pt]
&=\ 
k! \hsx
\abs{\sum\limits_{n=k+1}^\infty \hsx \frac{u_n}{n!}}
\\[15pt]
&\leq\ 
k! \hsx
\sum\limits_{n=k+1}^\infty \hsx \frac{\abs{u_n}}{n!}
\\[15pt]
&\leq\ 
C \hsx k! \hsx
 \sum\limits_{n=k+1}^\infty \hsx 
 \frac{A^n}{n!}
\\[15pt]
&=\ 
C \hsx k! \hsx
\bigg(
\frac{A^{k+1}}{(k+1)!} + \frac{A^{k+2}}{(k+2)!} + \cdots
\bigg)
\\[15pt]
&=\ 
C \hsx
\bigg(
\frac{A^{k+1}}{k+1} + \frac{A^{k+2}}{(k+1)(k+2)} + \cdots
\bigg)
\\[15pt]
&\leq\ 
C \hsx
\bigg(
\frac{A^{k+1}}{k+1} + \frac{A^{k+2}}{(k+1)^2} + \cdots
\bigg)
\\[15pt]
&=\
C A^k \hsx
\bigg(
\frac{A}{k+1} + \frac{A^2}{(k+1)^2} + \cdots
\bigg)
\\[15pt]
&=\
C A^k \hsx
\bigg(
\frac{\frac{A}{k + 1}}{1 - \frac{A}{k + 1}}
\bigg)
\hspace{2cm}
\bigg(
\frac{A}{k + 1} 
< 1
\bigg)
\quad \text{(cf. \S8, \#2)} 
\\[15pt]
&=\ 
C A^k \hsx
\bigg(
\frac{A}{k + 1 - A}
\bigg)
\\[15pt]
&=\ 
C \hsx 
\frac{A^{k+1}}{k + 1 - A}.
\end{align*}
\vspace{0.5cm}

\qquad \un{Step 9:} \ 
$\forall \ k$: $2A < k + 1$, 
\[
0 < k + 1 - 2A
\]
\qquad\qquad $\implies$
\[
k + 1 < 2( k + 1) - 2A
\]
\qquad\qquad $\implies$
\[
\frac{1}{k + 1 - A} \ < \ \frac{2}{k +1}.
\]
\vspace{0.2cm}

To recapitulate: \ $\forall \ k$: $2A < k + 1$, 
\allowdisplaybreaks
\begin{align*}
\abs{v_k} \ 
&\leq \  
C \hsx 
\frac{A^{k+1}}{k + 1 - A}
\\[12pt]
&<\ 
2 C \hsx 
\frac{A^{k+1}}{k + 1}.
\end{align*}
\vspace{0.2cm}

[Note: \ 
If $d = 1$, then the $v_k \in \Z$ (cf. Step 6) and if in addition, $A = 1$, then $\forall \ k \gg 0$, $v_k = 0$ 
(thus $\ds\sum\limits_{k=0}^\infty \hsx v_k X^k$ is a polynomial) and we would have a contradiction but, of course, in general $d > 1$ and $A > 1$.]
\vspace{0.5cm}

\qquad \un{Step 10:} \ 
Define $v_k(n)$ by the stipulation
\[
\sum\limits_{k = 0}^\infty \hsx v_k(n)X^k 
\ = \ 
(1 - a_1 X - \cdots - a_t X^t)^n \ 
\sum\limits_{k = 0}^\infty \hsx v_k X^k.
\]
Then $\forall \ n \geq 0$, 
\[
v_k(n+1) 
\ = \ 
v_k(n) - a_1 v_{k-1}(n) - \cdots - a_t v_{k-t} (n) \qquad (k \geq t).
\]
\vspace{0.2cm}

\qquad \un{Step 11:} \ 
Let 
\[
T 
\ = \ 
1 + \abs{a_1} + \cdots + \abs{a_t}.
\]
Then $\forall \ k \geq n t$, 
\[
\abs{v_k(n)} \leq (2 C) A^k T^n.
\]
Moreover
\[
d^k v_k(n) \in \Z
\]
and 
\[
n! \quad \text{divides} \quad d^k v_k(n).
\]
\vspace{0.2cm}

\qquad \un{Step 12:} \ 
If $k \geq n t$ and if $v_k(n) \neq 0$, then
\allowdisplaybreaks
\begin{align*}
n! \ 
&\leq\ 
\abs{d^k v_k(n)}
\\[12pt]
&=\ 
d^k\abs{v_k(n)}
\\[12pt]
&\leq\ 
d^k (2 C) A^k T^n
\\[12pt]
&=\ 
(2 C) (d A)^k T^n.
\end{align*}
So, if
\[
n! > (2 C) (d A)^k T^n
\]
and if $k \geq n t$, then $v_k(n) = 0$.
\vspace{0.5cm}

\qquad \un{Step 13:} \ 
Choose $n_0$ so large that $\forall \ n \geq n_0$, 
\[
n! > (2 C) (d A)^{10 n t} T^n.
\]
\vspace{0.2cm}

\qquad \un{Step 14:} \ 
\[
v_k(n) = 0  \ \forall \ n \geq n_0, \ n t \leq k \leq 10 \hsx n\hsx t.
\]
In particular: 
\[
v_k(n_0) = 0 \quad \text{if} \quad n_o t \leq k \leq 10 \hsx n_0 \hsx t.
\]
\vspace{0.2cm}

\qquad \un{Step 15:} \ 
\[
v_k(n) \ = \ 0 
\quad \text{if} \quad
n_0 \leq n \leq k / 10t,
\]
thus 
\[
v_k(n_0) \ = \ 0 
\quad \text{if} \quad
10 n_0 t \leq k.
\]
\vspace{0.2cm}

\qquad \un{Step 16:} \ 
$\forall \  k \ge n_0 \hsx t$, 
\[
v_k(n_0) \ = \ 0.
\]
\vspace{0.3cm}

Recall now the definition of $v_k(n)$, viz.
\[
\sum\limits_{k=0}^\infty \hsx 
v_k(n) X^k 
\ = \ 
(1 - a_1 X - \cdots - a_t X^t)^n \hsx
\sum\limits_{k=0}^\infty \hsx 
v_k X^k.
\]
Take $n = n_0$ $-$then in view of Step 16, 
\[
\sum\limits_{k=0}^\infty \hsx 
v_k(n_0) X^k \in \Q[X].  
\]
Therefore
\[
\sum\limits_{k=0}^\infty \hsx 
v_k X^k 
\in \Q(X),
\]
i.e., 
\[
\sum\limits_{k=0}^\infty \hsx 
v_k X^k
\]
is a rational function.  
\vspace{0.3cm}

To finish this sketch, let
\[
v(X) 
\ = \ 
\sum\limits_{k=0}^\infty \hsx 
v_k X^k.
\]
Then from the definitions
\[
\frac{v_k}{k!} 
\ - \ 
\frac{v_{k-1}}{(k-1)!} 
\ = \ 
\frac{u_k}{k!} 
\]
\qquad\qquad $\implies$
\[
v_k - k v_{k-1} \ = \ u_k
\]
\qquad\qquad $\implies$
\allowdisplaybreaks
\begin{align*}
\sum\limits_{k=0}^\infty \hsx 
(v_k - k v_{k-1}) X^k \ 
&=\ 
\sum\limits_{k=0}^\infty \hsx 
u_k X^k
\\[15pt]
&=\ 
\sum\limits_{n=0}^\infty \hsx 
u_n X^n
\\[15pt]
&=\ 
\sum\limits_{n=0}^\infty \hsx 
\bigg(
\sum\limits_{i=1}^t \hsx 
b_i \alpha_i^n\bigg) X^n 
\qquad \text{(cf. Step 2)}
\\[15pt]
&=\ 
\sum\limits_{i=1}^t \hsx 
b_i \hsx
\bigg(
\sum\limits_{n=0}^\infty \hsx 
\alpha_i^n X^n \bigg)
\\[15pt]
&=\ 
\sum\limits_{i=1}^t \hsx 
\frac{b_i}{1 - \alpha_i X}.
\end{align*}
On the other hand,
\allowdisplaybreaks
\begin{align*}
\sum\limits_{k=0}^\infty \hsx (v_k - k v_{k-1}) X^k \ 
&=\ 
v(X) - X \frac{d}{d X} (X v(X))
\\[12pt]
&=\ 
(1 - X) v(X) - X^2 \frac{d}{d X} v(X).
\end{align*}
Accordingly, if
\[
L 
\ \equiv \ 
-X^2 \frac{d}{dX} + (1 - X), 
\]
then $v(X)$ satisfies the differential equation 
\[
L v(X) 
\ = \ 
\sum\limits_{i=1}^t \hsx 
\frac{b_i}{1 - \alpha_i X}.
\]
And $v(X)$ is a rational function, thus the order of the nonzero poles of $L v(X)$ is at least 2.  
But the poles of the rational function
\[
\sum\limits_{i=1}^t \hsx 
\frac{b_i}{1 - \alpha_i X}
\]
are at the $\ds\frac{1}{\alpha_i}$ and are simple.  
Contradiction.
\end{x}
\vspace{0.3cm}

%% file: _22_exceptional_sets.tex
\chapter{
$\boldsymbol{\S}$\textbf{22}.\quad  EXCEPTIONAL SETS}
\setlength\parindent{2em}
\setcounter{theoremn}{0}
\renewcommand{\thepage}{\S22-\arabic{page}}

\ \indent 
Is it true that ``in general'' a transcendental function takes transcendental values at algebraic points?
\vspace{0.3cm}

\begin{x}{\small\bf DEFINITION} \ 
The 
\un{exceptional set}
\index{exceptional set}\index{$E_f$}
$E_f$ of an entire function $f$ is the set of algebraic numbers $\alpha$ such that $f(\alpha)$ is algebraic:
\[
E_f \ = \ \{\alpha \in \Qbar : f(\alpha) \in \Qbar\}.
\]
\end{x}
\vspace{0.3cm}

\begin{x}{\small\bf EXAMPLE} \ 
Take $f(z) = e^z$ $-$then $E_f = \{0\}$ (cf. \S21, \#4).
\end{x}
\vspace{0.3cm}

\begin{x}{\small\bf DEFINITION} \ 
A subset \mS of $\Qbar$ is 
\un{exceptional}
\index{exceptional} 
if there exists a transcendental function $f$ such that $E_f = S$.
\end{x}
\vspace{0.3cm}

\begin{x}{\small\bf EXAMPLE} \ 
An arbitrary finite subset
\[
\{\alpha_1, \ldots, \alpha_n\} \subset \Qbar
\]
is exceptional.
\vspace{0.2cm}

[Consider
\[
f(z) \ = \ e^{(z - \alpha_1) \cdots (z - \alpha_n)}.
\]
If $\alpha \in \Qbar$ and if $\alpha \neq \alpha_i$ $(i = 1, \ldots, n)$, then
\[
(\alpha - \alpha_1) \ldots (\alpha - \alpha_n) \in \Qbar
\]
is nonzero, hence $f(\alpha)$ is transcendental. (cf. \S21, \#4).]
\end{x}
\vspace{0.3cm}

\begin{x}{\small\bf EXAMPLE} \ 
Take 
\[
f(z) \ = \ e^z + e^{z+1}.
\]
Then $E_f = \emptyset$.
\vspace{0.2cm}

[First, $f(0) = 1 + e$ is not algebraic (since $e$ is transcendental) (cf. \S17, \#1).  
Suppose therefore that $\alpha$ is a nonzero algebraic number.
In \S21, \#1, take
\[
\alpha_0 = \alpha, \quad
\alpha_1 = \alpha + 1, \quad
b_0 = 1, \quad 
b_1 = 1,
\]
thus
\[
e^\alpha + e^{\alpha + 1}
\]
is transcendental.]
\end{x}
\vspace{0.3cm}

\begin{x}{\small\bf THEOREM} \ 
Given any subset $S \subset \Qbar$, there exists a transcendental function $f$ such that $E_f = S$.
\end{x}
\vspace{0.3cm}

\begin{x}{\small\bf \un{N.B.}} \ 
It was proved in 1895 by St\"ackel that there exists a transcendental function $f$ such that $E_f = \Qbar$.
\end{x}
\vspace{0.3cm}

\begin{x}{\small\bf DEFINITION} \ 
The exceptional set $E_f(\tmul)$ with multiplicities of an entire function $f$ is the subset of $\Qbar \times \Z_{\geq 0}$ consisting of those points $(\alpha, n)$ such that $f^{(n)}(\alpha) \in \Qbar$.
\vspace{0.2cm}

[Note: \ 
Here $f^{(n)}$ is the $n^\nth$ derivative of $f$.]
\end{x}
\vspace{0.3cm}

\begin{x}{\small\bf THEOREM} \ 
Given any subset $S \hsx \subset \hsx \Qbar \times \Z_{\geq 0}$, there exists a transcendental function $f$ such that 
$E_f(\tmul) = S$.
\end{x}
\vspace{0.3cm}


%% file: _23_complex_logarithms_and_complex_powers.tex
\chapter{
$\boldsymbol{\S}$\textbf{23}.\quad  COMPLEX LOGARITHMS AND COMPLEX POWERS}
\setlength\parindent{2em}
\setcounter{theoremn}{0}
\renewcommand{\thepage}{\S23-\arabic{page}}


\begin{x}{\small\bf DEFINITION} \ 
Given a complex number $z \neq 0$, a 
\un{logarithm of $z$}
\index{logarithm of $z$} 
is a complex number $w$ such that $e^w = z$, denoted $\log z$.
\vspace{0.2cm}

[Note: \ 
$\log 0$ is left undefined (there is no complex number $w$ such that $e^w = 0$).]
\end{x}
\vspace{0.3cm}

Therefore
\[
\log z 
\ = \ 
\elln(\abs{z}) + \sqrt{-1} \hsx \arg z,
\]
where $\elln(\abs{z})$ is the natural logarithm of $\abs{z}$ (cf. \S10, \#3 \& \#4) and arg $z$ is given all admissible values.  
Since the latter differ by multiples of $2 \pi$, it follows that the various determinations of log $z$ differ by multiples of $2 \pi \sqrt{-1}$.
\vspace{0.3cm}

\begin{x}{\small\bf DEFINITION}\ 
The 
\un{principal determination}
\index{principal determination} 
of the logarithm corresponds to the choice
\[
-\pi
\ < \ 
\Arg z 
\ \leq \ 
\pi,
\]
so
\[
-\pi
\ < \ 
\Imx (\log z)
\ \leq \ 
\pi
\]
and one signifies this by writing $\Log z$, thus $\restr{\Log}{\R_{>0}} = \elln$.
\end{x}
\vspace{0.3cm}

\begin{x}{\small\bf EXAMPLE} \ 
\[
\Log (-3 \sqrt{-1}) 
\ = \ 
\elln(3) - \frac{\pi \sqrt{-1}}{2}.
\]
\end{x}
\vspace{0.3cm}

\begin{x}{\small\bf \un{N.B.}} \ 
The restriction of the exponential function to the horizontal strip \mS consisting of all complex numbers 
$x + \sqrt{-1} \hsx y$ $(-\pi < y \leq \pi)$ has an inverse: 
$\restr{\exp}{S}$ maps \mS bijectively to $\C^\times = \C - \{0\}$ and the inverse of this restriction is
$\Log: \C^\times \ra S$, hence
\[
\begin{cases}
\ \Log \circ \restr{\exp}{S} = \id_S\\
\ \exp \circ \hsx \Log =  \id_{\C^\times}.
\end{cases}
\]
\vspace{0.2cm}

[Note: \ 
Log is discontinuous at each negative real number but is continuous everywhere else on $\C^\times$.]
\end{x}
\vspace{0.3cm}

\begin{x}{\small\bf REMARK} \ 
It is always true that
\[
\Log (z_1 z_2) 
\ \equiv \ 
\Log z_1 + \Log z_2 
\qquad (\modx 2 \pi \sqrt{-1})
\]
but the relation
\[
\Log (z_1 z_2) 
\ = \ 
\Log z_1 + \Log z_2 
\]
can fail.  E.g.:
\allowdisplaybreaks
\begin{align*}
\Log((-1)\sqrt{-1}) \ 
&=\ 
\Log \bigl(-\sqrt{-1}\bigr)
\\[12pt]
&=\ 
\elln\abs{-\sqrt{-1}} - \frac{\pi \sqrt{-1}}{2}
\\[12pt]
&=\ 
\elln(1) - \frac{\pi \sqrt{-1}}{2}
\\[12pt]
&=\ 
-\frac{\pi \sqrt{-1}}{2}
\end{align*}
while
\allowdisplaybreaks
\begin{align*}
\Log(-1) + \Log (\sqrt{-1}) \ 
&=\ 
(\elln(1) + \pi \sqrt{-1} + \bigl(\elln(1) + \frac{\pi \sqrt{-1}}{2}\bigr)
\\[12pt]
&=\ 
\frac{3 \pi \sqrt{-1}}{2} 
\\[12pt]
&\neq\ 
- \frac{\pi \sqrt{-1}}{2}.
\end{align*}
\end{x}
\vspace{0.3cm}

\begin{x}{\small\bf LEMMA} \ 
\[
\Log z 
\ = \ 
\int\limits_1^z \hsx 
\frac{dt}{t}
\qquad \bigl(\abs{\arg z} < \pi),
\]
the integral being taken along the line segment $[1,z]$.
\end{x}
\vspace{0.3cm}

\begin{x}{\small\bf LEMMA} \ 
\[
\Log z 
\ = \ 
\sum\limits_{n=1}^\infty \hsx
\frac{(-1)^{n-1}}{n} \hsx (z - 1)^n 
\qquad (\abs{z-1} < 1).
\]
\end{x}
\vspace{0.3cm}

\begin{x}{\small\bf DEFINITION} \ 
Let \mD be an open simply connected region in the complex plane that does not contain 0 $-$then a 
\un{branch of log $z$}
\index{branch of log $z$} 
is a continuous function \mL with domain \mD such that $L(z)$ is a logarithm of $z$ for each $z$ in \mD:
\[
e^{L(z)} \ = \ z.
\]
\end{x}
\vspace{0.3cm}

\begin{x}{\small\bf EXAMPLE} \ 
Take $D = \C - \R_{\leq 0}$ $-$then the restriction of Log to \mD is a branch of log $z$.
\end{x}
\vspace{0.3cm}

\begin{x}{\small\bf CONSTRUCTION} \ 
A branch of $\log z$ with domain \mD can be obtained by first fixing a point $a$ in \mD, then choosing a logarithm $b$ of $a$, 
and then defining \mL by the prescription
\[
L(z) 
\ = \ 
b + \int\limits_a^z \hsx \frac{dw}{w}.
\]
Here the integration is along any path in \mD that connects $a$ and $z$.
\end{x}
\vspace{0.3cm}

\begin{x}{\small\bf LEMMA} \ 
$L(z)$ is holomorphic in \mD, its derivative being $\ds\frac{1}{z}$.
\vspace{0.2cm}

[Note: \ 
Different choices of $b$ will in general lead to different functions.]
\end{x}
\vspace{0.3cm}

\begin{x}{\small\bf RAPPEL} \ 
If $\alpha$ is a nonzero algebraic number, then $e^\alpha$ is transcendental (cf. \S21, \#4) (Hermite-Lindemann).
\end{x}
\vspace{0.3cm}


\begin{x}{\small\bf EXAMPLE} \ 
(cf. \S21, \#9) \ 
Let $\alpha$ be a nonzero algebraic number $-$then $\Log \alpha$ is transcendental.
\vspace{0.2cm}

[The point is that $e^{\Log \alpha} = \alpha$.] 
\end{x}
\vspace{0.3cm}

Let $a$ be a complex number with $a \neq 0$, $\neq e$.
\vspace{0.3cm}

\begin{x}{\small\bf DEFINITION} \ 
The 
\un{principal power of a}
\index{principal power of a} 
is the holomorphic function
\[
z \ra a^z 
\ = \ 
e^{z \Log a}.
\]
\end{x}
\vspace{0.3cm}

\begin{x}{\small\bf DEFINITION} \ 
The 
\un{$k^\nth$ associate of $a^z$}
\index{$k^\nth$ associate of $a^z$} 
$(k \in \Z)$ is the holomorphic function
\allowdisplaybreaks
\begin{align*}
z \ra 
&e^{z \bigl(\Log a + 2 k \pi \sqrt{-1}\bigr)}
\\[12pt]
&=\ a^z\bigl(e^{2 k \pi \sqrt{-1} \hsx z}\bigr).
\end{align*}
\end{x}
\vspace{0.3cm}

\begin{x}{\small\bf \un{N.B.}}  \ 
The reason for excluding $e$ is that we want $e^z$ to remain single valued and to mean the power series
\[
\sum\limits_{n=0}^\infty \hsx
\frac{z^n}{n!}.
\]
\end{x}
\vspace{0.3cm}

\begin{x}{\small\bf EXAMPLE} \ 
\[
1^z 
\ = \ 
e^{z \Log 1} 
\ = \ 
e^{z (\elln(1) + \sqrt{-1} \hsx 0)} 
\ = \ 
1^0 
\ = \ 
1
\]
and its $k^\nth$ associate is
\[
e^{z (\Log 1 + 2 k \pi \sqrt{-1})}
\ = \
e^{2 k \pi \sqrt{-1}\hsx z}.
\]
\end{x}
\vspace{0.3cm}

\begin{x}{\small\bf EXAMPLE} \ 
Take $a = \sqrt{-1}$ and take $z = -2\sqrt{-1}$ $-$then with this data,
\allowdisplaybreaks
\begin{align*}
\sqrt{-1}^{\hsx -2 \sqrt{-1}} \ 
&=\ 
e^{-2 \sqrt{-1} \hsx\hsx \Log(\sqrt{-1})}
\\[12pt]
&=\ 
e^{-2 \sqrt{-1} \bigl(\sqrt{-1} \hsx \frac{\pi}{2}\bigr)}.
\end{align*}
Therefore the associates of $\sqrt{-1}^{\hsx -2 \hsx \sqrt{-1}}$ are the
\[
e^{-2 \sqrt{-1} \bigl(\sqrt{-1} \hsx \frac{\pi}{2} + 2 k \pi \sqrt{-1}\bigr)}
\ = \ 
e^{\pi (4k + 1)} 
\qquad (k \in \Z).
\]
\end{x}
\vspace{0.3cm}

\begin{x}{\small\bf EXAMPLE} \ 
Let $n$ be a positive integer and write
\[
a 
\ = \ 
\abs{a} \hsx 
e^{\sqrt{-1} \hsx \theta} 
\qquad (-\pi \leq \theta \leq \pi).
\]
Then
\allowdisplaybreaks
\begin{align*}
a^{\frac{1}{n}} \ 
&=\ 
e^{\frac{1}{n} \Log a}
\\[15pt]
&=\ 
e^{\frac{1}{n} \big(\elln(\abs{a}) + \sqrt{-1}\hsx \theta \big)}
\\[15pt]
&=\ 
e^{\frac{1}{n} \elln(\abs{a})} \hsx e^{\frac{1}{n} \sqrt{-1}\hsx \theta}
\\[15pt]
&=\ 
e^{\elln(\abs{a}^{\frac{1}{n}})} \hsx e^{\frac{1}{n} \sqrt{-1}\hsx \theta}
\\[15pt]
&=\ 
\abs{a}^{\frac{1}{n}} \hsx e^{\frac{1}{n} \sqrt{-1}\hsx \theta}.
\end{align*}
Therefore the associates of $\ds a^{\frac{1}{n}}$ are the 
\[
\abs{a}^{\frac{1}{n}} \hsx e^{\frac{1}{n} \sqrt{-1}\hsx (\theta + 2 k \pi)}
\qquad (k \in \Z).
\]
And there are only $n$ different values for 
\[
\exp\bigg(\frac{1}{n} 2 k \pi \sqrt{-1}\bigg).
\]
\end{x}
\vspace{0.3cm}

The \hfill laws \hfill of \hfill exponents \hfill spelled \hfill out \hfill in \hfill \S4 \hfill over\hfill  $\R$ \hfill do \hfill not \hfill hold \hfill without \hfill qualification \\
over $\C$.

\qquad \textbullet \quad 
In general, $\bigl(a^b \bigr)^c$ has more values than $a^{b c}$.
\vspace{0.2cm}

\qquad \textbullet \quad 
In general, $a^b a^c$ has more values than $a^{b + c}$.


%% file: _24_the_gelfond_schneider_theorem.tex
\chapter{
$\boldsymbol{\S}$\textbf{24}.\quad  THE GELFOND-SCHNEIDER THEOREM}
\setlength\parindent{2em}
\setcounter{theoremn}{0}
\renewcommand{\thepage}{\S24-\arabic{page}}

\ \indent 

This is the following statement.
\vspace{0.3cm}

\begin{x}{\small\bf THEOREM} \ 
If $\alpha \neq 0$, 1 is algebraic and if $\beta \notin \Q$ is algebraic, then $\alpha^\beta$ is transcendental.
\vspace{0.2cm}

[Note: \ 
Here $\alpha^\beta$ is the principal power (cf. \S23, \#14): 
\[
\alpha^\beta \ = \ e^{\beta \Log \alpha}.
\]
Nevertheless it can be shown that the Gelfond-Schneider theorem goes through if the principal power $\alpha^\beta$ is replaced by any of its associates (cf. \S31, \#16).]
\end{x}
\vspace{0.3cm}

Special Cases:
\vspace{0.3cm}

\begin{x}{\small\bf EXAMPLE} \ 
$2^{\sqrt{2}}$ is transcendental.
\end{x}
\vspace{0.3cm}

\begin{x}{\small\bf EXAMPLE} \ 
$\sqrt{3}^{\sqrt{2}}$ is transcendental.
\end{x}
\vspace{0.3cm}

\begin{x}{\small\bf EXAMPLE} \ 
$\sqrt{-1}^{\sqrt{-1}} = e^{^{-\frac{\pi}{2}}}$ is transcendental.
\end{x}
\vspace{0.3cm}

\begin{x}{\small\bf EXAMPLE} \ 
$e^\pi$ is transcendental.
\vspace{0.2cm}

[Starting from the fact that
\[
e^{\pi \hsx \frac{\sqrt{-1}}{2}} \ = \ \sqrt{-1}
\]
and using the principal determination of the logarithm:
\allowdisplaybreaks
\begin{align*}
\Log \sqrt{-1} \ 
&=\ 
\elln(\abs{\sqrt{-1}}) + \sqrt{-1} \hsx \Arg \sqrt{-1} 
\\[12pt]
&=\ 
\elln(1) + \sqrt{-1} \hsx  \frac{\pi}{2}
\\[12pt]
&=\ 
\sqrt{-1} \hsx \frac{\pi}{2}
\end{align*}
\qquad\qquad $\implies$
\[
\pi \ = \ 
- 2 \sqrt{-1} \hsx \Log \sqrt{-1}
\]
\qquad\qquad $\implies$
\[
e^\pi 
\ = \ 
e^{-2 \sqrt{-1} \hsx \Log \sqrt{-1}} 
\ = \ 
\sqrt{-1}^{-2 \sqrt{-1}}
\qquad \text{(cf. \S23, \#18)}
\]
and the entity on the right is transcendental.]
\end{x}
\vspace{0.3cm}

\begin{x}{\small\bf EXAMPLE} \ 
Take $f(z) = 2^z$, thus $2^z = e^{z \Log 2} = e^{z \elln(2)}$.\\[12pt]
If $\alpha = 2$ in Gelfond-Schneider and if $z \notin \Q$ is algebraic, then $2^z$ is transcendental.  \\[12pt]
On the other hand, the $\ds 2^{1/n}$ $(n \in \N)$ are algebraic, as are the $\ds\big(2^{1/n}\big)^m$ $(m \in \Z)$.  
Therefore the exceptional set $E_f$ of $f$ is equal to $\Q$.
\vspace{0.2cm}

[Note: \ 
$f^\prime(z) = 2^z \elln(2)$, so
\[
E_f \hsx \cap \hsx E_{f^\prime} \ = \ \emptyset
\]
since $\elln(2)$ is transcendental (cf. \S21, \#9).]
\end{x}
\vspace{0.3cm}

\begin{x}{\small\bf EXAMPLE} \ 
Take $f(z) = e^{\pi \sqrt{-1} \hsx z}$ $-$then 
\[
e^{\pi \sqrt{-1} \hsx z} \ = \ (-1)^z,
\]
the principal power.  In fact, 
\allowdisplaybreaks
\begin{align*}
(-1)^z \ 
&=\ 
e^{z \Log - 1}
\\[12pt]
&=\ 
e^{z (\elln(\abs{-1}) + \pi \sqrt{-1})}
\\[12pt]
&=\ 
e^{\pi \sqrt{-1} \hsx z}.
\end{align*}
So, if $\alpha = -1$ in Gelfond-Schneider and if $z \notin \Q$ is algebraic, then $e^{\pi \sqrt{-1} \hsx z}$ is
transcendental.  
As for what happens if $z \in \Q$, write
\[
e^{\pi \sqrt{-1} \hsx z} \ = \cos(\pi z) + \sqrt{-1} \sin(\pi z)
\]
and quote the wellknown fact that the trigonometric functions cos and sin are algebraic numbers at arguments that are rational multiplies of $\pi$.  
Therefore the exceptional set $E_f$ of $f$ is equal to $\Q$.
\end{x}
\vspace{0.3cm}

\begin{x}{\small\bf THEOREM} \ 
Given nonzero complex numbers $a$ and $b$ with $a \notin \Q$, then at least one of $a$, $e^b$, $e^{ab}$ is transcendental.
\end{x}
\vspace{0.3cm}

\begin{x}{\small\bf \un{N.B.}} \ 
$\#8 \Leftrightarrow \#1$.
\vspace{0.2cm}

[To see that $\#8 \implies \#1$, take $a = \beta$, $b = \Log \alpha$ $-$then at least one of the following numbers is transcendental: \ 
$\beta$, $\ds e^{\Log \alpha} = \alpha$, or $\ds e^{\beta \Log \alpha} = \ds{\alpha^\beta}$.  
But the first two of these are algebraic, hence $\ds\alpha^\beta$ must be transcendental.  
That $\#1 \implies \#8$ is analogous.]
\end{x}
\vspace{0.3cm}

\begin{x}{\small\bf EXAMPLE} \ 
Let $\alpha$, $\beta$ be algebraic numbers not equal to 0 or 1.  
Suppose that 
\[
\frac{\Log \beta}{\Log \alpha} \notin \Q.
\]
Then
\[
\frac{\Log \beta}{\Log \alpha}
\]
is transcendental.
\vspace{0.2cm}

[In \#8, take
\[
a \ = \ \frac{\Log \beta}{\Log \alpha} 
\quad \text{and} \quad 
b = \Log \alpha.
\]
Then at least one of the following numbers is transcendental:
\[
\frac{\Log \beta}{\Log \alpha}, \ 
e^{\Log \alpha} = \alpha, \ 
e^{\ds\frac{\Log \beta}{\Log \alpha} \hsx \Log \alpha} = e^{\Log \beta} = \beta.]
\]
\vspace{0.2cm}

[Note: \ 
If $\Log \alpha$ and $\Log \beta$ are linearly independent over $\Q$, then 
\[
\frac{\Log \beta}{\Log \alpha} \notin \Q, 
\]
as can be seen by putting
\[
\gamma \ = \ \frac{\Log \beta}{\Log \alpha}
\]
and considering the dependence relation
\[
\gamma \hsx \Log \alpha - \Log \beta \ = \ 0.
\]
Consequently
\[
\frac{\Log \beta}{\Log \alpha}
\]
is transcendental, thus for any nonzero algebraic numbers $\mu$ and $\nu$, 
\[
\mu \hsx \Log \alpha + \nu \hsx \Log \beta \ \neq \ 0, 
\]
i.e., $\Log \alpha$ and $\Log \beta$ are linearly independent over $\Qbar$.]
\end{x}
\vspace{0.3cm}

\begin{x}{\small\bf EXAMPLE} \ 
Let $r$ be a positive rational number.
Write (see the Appendix to \S10)
\[
\log_{10} (r) \ = \ \frac{\elln(r)}{\elln(10)}.
\]
Therefore, if $\log_{10} (r)$ is not rational, then by the above it must be transcendental (cf. \S5, \#15).
\end{x}
\vspace{0.3cm}

Question: \ For what pairs $(\beta,t)$ $(\beta \in \Qbar$, $\beta \neq 0$ and $t \in \R^\times)$ 
is $e^{t \beta}$ algebraic?
\vspace{0.3cm}

\begin{x}{\small\bf EXAMPLE} \ 
Take $\beta \in \Qbar \cap \R$ $(\beta \neq 0)$ and
\[
t \ = \ \frac{\elln(2)}{\beta}.
\]
Then
\[
e^{t \beta} 
\ = \ 
e^{\elln(2)}
\ = \ 
2.
\]
\end{x}
\vspace{0.3cm}

\begin{x}{\small\bf EXAMPLE} \ 
Take $\beta \in \Qbar \cap \sqrt{-1} \hsx \R$ $(\beta \neq 0)$ and
\[
t \ = \ \frac{\sqrt{-1} \hsx \pi}{\beta}.
\]
Then
\[
e^{t \beta} 
\ = \ 
e^{\sqrt{-1} \hsx \pi}
\ = \ 
-1.
\]
\end{x}
\vspace{0.3cm}

\begin{x}{\small\bf THEOREM} \ 
Let $\beta \in \Qbar$ $(\beta \neq 0)$ and $t \in \R^\times$.  
Assume: \ $\beta \notin \R \hsx \cup \hsx \sqrt{-1} \hsx \R$ $-$then $e^{t \beta}$ is transcendental.
\vspace{0.2cm}

PROOF \ 
Put $\alpha = e^{t \beta}$ $-$then the complex conjugate $\ov{\alpha}$ of $\alpha$ is 
$e^{t \ov{\beta}} = e^{\ov{\beta}/ \beta}$.  
The algebraic number $\ov{\beta}/\beta$ is not real 
(for $\abs{\ov{\beta}/\beta} = 1$ but $\ov{\beta}/\beta \neq \pm 1$), hence is not rational.  
In \#8, take
\[
a \ = \ \ov{\beta}/\beta, \quad
b = t \beta,
\]
leading thereby to
\[
\ov{\beta}/\beta, \quad 
e^{t \beta} \ = \ \alpha, \quad
e^{t \ov{\beta}} \ = \ \ov{\alpha}.
\]
Since $\ov{\beta}/\beta$ is algebraic, either $\alpha$ or $\ov{\alpha}$ must be transcendental.   
But $\alpha$ is transcendental iff $\ov{\alpha}$  is transcendental.
\end{x}
\vspace{0.3cm}

It remains to give a proof of Gelfond-Schneider, a task that requires some preliminaries.


%% file: _25_interpolation_determinants.tex
\chapter{
$\boldsymbol{\S}$\textbf{25}.\quad  INTERPOLATION DETERMINANTS}
\setlength\parindent{2em}
\setcounter{theoremn}{0}
\renewcommand{\thepage}{\S25-\arabic{page}}


\begin{x}{\small\bf NOTATION} \ 
Given $w \in \C$, $R \in \R_{\geq 0}$, let
\[
\begin{cases}
\ D(R,w) = \{z \in \C: \abs{z - w} < R\}\\[8pt]
\ \ov{D}(R,w) = \{z \in \C: \abs{z - w} \leq R\} 
\end{cases}
. 
\]
\vspace{0.2cm}

[Note: \ 
Write
\[
\begin{cases}
\ D(R)\\
\ \ov{D}(R)
\end{cases}
\]
if $w = 0$.]
\end{x}
\vspace{0.3cm}

\begin{x}{\small\bf NOTATION} \ 
Let $\abs{f}_R$ stand for the maximum value of $\abs{f(z)}$ when $\abs{z} = R$.
\end{x}
\vspace{0.3cm}

\begin{x}{\small\bf RAPPEL} \ 
If $f(z)$ is a function holomorphic in $D(R)$ and continuous in $\ov{D}(R)$, then
\[
\abs{f(z)} \ \leq \ \abs{f}_R
\]
for every $z \in \ov{D}(R)$.
\end{x}
\vspace{0.3cm}

\begin{x}{\small\bf LEMMA} \ 
Let \mT be a nonnegative integer, let $r$ and \mR be positive real numbers subject to $0 < r \leq R$, and let $F(z)$ be a function of one complex variable holomorphic in $D(R)$ and continuous in $\ov{D}(R)$.  
Assume: \ \mF has a zero of multiplicity at least \mT at 0 $-$then
\[
\abs{F}_r 
\ \leq \ 
\bigg(\frac{R}{r}\bigg)^{-T} \hsx \abs{F}_R.
\]
\vspace{0.2cm}

PROOF \ 
Put
\[
G(z) \ = \ z^{-T} F(z).
\]
Then
\[
\abs{G}_r \ \leq \ \abs{G}_R
\]
or still, 
\[
r^{-T} \abs{F}_r \ \leq \ R^{-T} \abs{F}_R
\]
or still, 
\[
\abs{F}_r \ \leq \ \bigg(\frac{R}{r}\bigg)^{-T} \abs{F}_R.
\]
\end{x}
\vspace{0.3cm}

\begin{x}{\small\bf THEOREM} \ 
Let $r$ and \mR be positive real numbers subject to $0 < r \leq R$, let $f_1(z), \ldots, f_L(z)$ be functions of one complex variable which are holomorphic in $D(R)$ and continuous in $\ov{D}(R)$, and let $\zeta_1, \ldots, \zeta_L$ belong to the disc 
$\abs{z} \leq r$.  Put
\[
\Delta \ = \ \det
\begin{pmatrix}
f_1(\zeta_1) &\cdots\cdots &f_L(\zeta_1)\\
\vdots &&\vdots\\
f_1(\zeta_L) &\cdots\cdots &f_L(\zeta_L)\\
\end{pmatrix}
.
\]
Then
\[
\abs{\Delta} 
\ \leq \ 
\bigg(\frac{R}{r}\bigg)^{-L(L-1)/2} \hsx
L! \hsx\hsx 
\prod\limits_{j=1}^L \hsx 
\abs{f_j}_R.
\]
\vspace{0.2cm}

PROOF \ 
Let $F(z)$ be the determinant of the $L \times L$ matrix
\[
(f_j(\zeta_i z))_{1 \leq j, i \leq L} \quad (\implies F(1) = \Delta).
\]
Since the $\zeta_i$ satisfy $\abs{\zeta_i} \leq r$, the functions $f_j(\zeta_i z)$ are holomorphic in $D(R/r)$ and continuous in 
$\ov{D}(R/r)$.  
And since the determinant is a sum of products of the 
$f_j(\zeta_i z)$, the determinant $F(z)$ itself is holomorphic in $D(R/r)$ and continuous in $\ov{D}(R/r)$.  
The claim then is that $F(z)$ vanishes at 0 with multiplicity at least $L(L-1)/2$.  
To see this, put
\[
K \ = \ L(L-1) / 2
\]
and consider the expansion
\[
f_j(\zeta_i z) 
\ = \ 
\sum\limits_{k=0}^{K-1} \hsx 
a_k(j) \zeta_i^k z^k + z^K g_{i j} (z), 
\]
where $a_k(j) \in \C$ and $g_{i j}(z)$ is holomorphic in $D(R/r)$ and continuous in $\ov{D}(R/r)$.  
Since the determinant is linear in its columns, one can view $F(z)$ as $z^K$ times a function holomorphic in $D(R/r)$ plus terms involving the factor 
\[
z^{n_1 + n_2 + \cdots +  n_L} \hsx \det\big(\zeta_i^{n_j}\big),
\]
i.e., 
\[
z^{n_1 + n_2 + \cdots +  n_L} \  \det
\begin{pmatrix}
\zeta_1^{n_1} &\cdots\cdots &\zeta_1^{n_L}\\
\vdots &&\vdots\\
\zeta_L^{n_1} &\cdots\cdots &\zeta_L^{n_L}\\
\end{pmatrix}
,
\]\\
where $n_1, n_2, \ldots, n_L \in \Z_{\geq 0}$ and $n_j \in \{0, 1, \ldots, K-1\}$.  
The determinant vanishes if two of the $n_j$ are identical, so the nonzero terms satisfy
\[
n_1 + n_2 + \cdots + n_L 
\ \geq \ 
0 + 1 + \cdots + (L - 1) 
\ = \ 
\frac{L(L-1)}{2}.
\]
Take now in \#4
\[
T \ = \ L(L-1)/2
\]
and replace $r$ by 1 and \mR by $R/r$, hence
\allowdisplaybreaks
\begin{align*}
\abs{\Delta} \ 
&=\ 
\abs{F(1)} 
\\[12pt]
&\leq\ 
\abs{F}_1
\\[12pt]
&\leq\ 
\bigg(\frac{R}{r}\bigg)^{-L(L-1)/2} \abs{F}_{R/r}.
\end{align*}
It remains to bound $\abs{F}_{R/r}$.  
From its very definition, the determinant of an $L \times L$ matrix is the sum of $L!$ products, where each product consists of \mL entries such that for each row and column only one entry is a part of a product.  
Since $\abs{z} = R/r$ $\implies$ $\abs{\zeta_i z} \leq R$, for each column index $j$, 
\[
\abs{f_j(\zeta_i z)} \ \leq \ \abs{f_j}_R \qquad (i = 1, 2, \ldots, L).
\]
Therefore
\[
\abs{F}_{R/r}
\ \leq \
L! \hsx 
\prod\limits_{j=1}^L \hsx 
\abs{f_j}_R.
\]
So finally
\[
\abs{\Delta} \ \leq \ 
\bigg(\frac{R}{r}\bigg)^{-L(L-1)/2} \hsx
L! \hsx 
\prod\limits_{j=1}^L \hsx 
\abs{f_j}_R.
\]
\end{x}
\vspace{0.3cm}

\begin{x}{\small\bf REMARK} \ 
The derivatives of $F(z)$ can be calculated via an application of the product rule, viz:
\[
\bigg(\frac{d}{dz}\bigg)^k F(z) 
\ = \ 
\sum\limits_{\kappa_1 + \cdots + \kappa_L = k} \hsx
\frac{k!}{\kappa_1! \cdots \kappa_L!} \hsx
\det\bigg(\bigg(\frac{d}{dz}\bigg)^{\kappa_i} \hsx f_j(\zeta_j z)\bigg)_{1 \leq j, \ i \leq L}.
\]
\end{x}
\vspace{0.3cm}

The foregoing can be generalized by incorporating derivatives.
\vspace{0.3cm}

\begin{x}{\small\bf THEOREM} \ 
Let $r$ and \mR be positive real numbers subject to $0 < r \leq R$, 
let $\sigma_1, \ldots, \sigma_L$ be nonnegative integers, let $f_1, \ldots, f_L$ be entire functions, and let 
$\zeta_1, \ldots, \zeta_L$ belong to the disc $\abs{z} \leq r$.  
Put
\[
\Delta \ = \ 
\det\bigg(\bigg(\frac{d}{dz}\bigg)^{\sigma_i} \hsx f_j(\zeta_i)\big)_{1 \leq j, i \leq L}.
\]
Then
\[
\abs{\Delta} 
\ \leq \ 
\bigg(\frac{R}{r}\bigg)^{-L(L-1)/2 + \sigma_1 + \cdots + \sigma_L} L! \hsx
\prod\limits_{j=1}^L \hsx
\max\limits_{1 \leq i \leq L} \hsx 
\sup\limits_{\abs{z} = R} \hsx
\abs{\bigg(\frac{d}{dz}\bigg)^{\sigma_i} \hsx f_j(z)}.
\]
\end{x}
\vspace{0.3cm}

\newpage
\[
\text{APPENDIX}
\]
\vspace{0.5cm}

Suppose that $1 \leq j \leq p_k \ (\in \N)$, $1 \leq k \leq \ell$, $1 \leq i \leq n$ $-$then
\[
\frac{d^{i-1}}{dz^{i-1}} \hsx \big(z^{j-1} e^{w_k z}\big)\bigg|_{z = 0}
\ = \ 
\frac{d^{j-1}}{dz^{j-1}} \hsx \big(z^{i-1}\big)\bigg|_{z = w_k},
\]
their common value being
\[
\begin{cases}
\ \ds\frac{(i-1)!}{(i-j)!} \hsx w_k^{i-j}\qquad \text{if} \quad i \geq j\\[12pt]
\ \hspace{0.5cm} 0 \hspace{2.2cm} \text{if} \quad i < j
\end{cases}
.
\]


%% file: _26_zero_estimates.tex
\chapter{
$\boldsymbol{\S}$\textbf{26}.\quad  ZERO ESTIMATES}
\setlength\parindent{2em}
\setcounter{theoremn}{0}
\renewcommand{\thepage}{\S26-\arabic{page}}


\begin{x}{\small\bf LEMMA} \ 
Let $P_1, \ldots, P_n$ be nonzero polynomials in $\R[X]$ of degrees $d_1, \ldots, d_n$ and let $w_1, \ldots, w_n$ be distinct real numbers $-$then
\[
F(x)
\ = \ 
\sum\limits_{j=1}^n \hsx
P_j(x) e^{w_j x}
\]
has at most
\[
d_1 + \cdots + d_n + n - 1
\]
real zeros counting multiplicities.
\end{x}
\vspace{0.3cm}

To begin with:
\vspace{0.3cm}

\begin{x}{\small\bf SUBLEMMA}\ 
If a continuously differentiable function \mF of a real variable $x$ has at least \mN real zeros counting multiplicities (\mN a positive integer), then its derivative $F^\prime$ has at least $N - 1$ real zeros counting multiplicities.
\vspace{0.2cm}

PROOF \ 
Let $x_1, \ldots, x_k$ $(k \geq 1)$ be distinct real zeros of \mF arranged in increasing order: $x_1 < \cdots < x_k$ 
with $n_1$ the multiplicity of $x_1, \ldots, n_k$ the multiplicity of $x_k$ and 
$n_1 + \cdots + n_k \geq N$ $-$then $x_i$ is a zero of $F^\prime$ of multiplicity $\geq n_i - 1$ $(1 \leq i \leq k)$.  
Owing to Rolle's theorem, $F^\prime$ has at least one zero in the open interval 
$]x_i, x_{i+1}[$ $(1 \leq i \leq k)$, so all told, $F^\prime$ has at least
\allowdisplaybreaks
\begin{align*}
(n_1 - 1) + \cdots + (n_k - 1) + (k - 1) \ 
&\geq
N - k + (k - 1) 
\\[12pt]
&=\ 
N - 1
\end{align*}
real zeros counting multiplicities.
\end{x}
\vspace{0.3cm}

Passing to the proof of \#1, upon multiplying through by $e^{-w_n x}$, it can be 
assumed that $w_n = 0$ and $w_j \neq 0$ for $j = 1, \ldots, n - 1$.  
Put
\[
D \ = \ d_1 + \cdots + d_n + n
\]
and proceed from here by induction on \mD, matters being clear if $D = 1$ (since $n = 1$ and $d_1 = 0$) so in this case there are at most $D - 1 = 0$ real zeros.  
Suppose now that the lemma holds if $k = 2, \ldots, D-1$ and consider the situation at level $k = D$.  
Take the first derivative of $F(x)$:
\[
F^\prime(x) 
\ = \ 
\sum\limits_{j=1}^{n-1} \hsx
(w_j P_j(x) + \frac{d}{dx} P_j(x)) + \frac{d}{dx} P_n(x).
\]
Then
\[
w_j P_j(x) + \frac{d}{dx} P_j(x)
\]
is a polynomial of degree $d_j$ whereas $\ds\frac{d}{dx} P_n(x)$ is a polynomial of degree $d_n - 1$.  
It therefore follows from the induction hypothesis that $F^\prime(x)$ has at most
\[
d_1 + \cdots + d_{n-1} + d_n - 1 + n - 1 
\ = \ 
d_1 + \cdots + d_n + n - 2
\]
real zeros counting multiplicities.  
Let \mN be a postive integer such that \mF has at least \mN real zeros counting multiplicities, hence by \#2,
\[
N  - 1
\ \leq \  
d_1 + \cdots + d_n + n - 2
\]
\qquad\qquad $\implies$
\[
N 
\ \leq \ 
d_1 + \cdots + d_n + n - 1.
\]
\vspace{0.3cm}

\begin{x}{\small\bf REMARK} \ 
Let $d_1, \ldots, d_n$ be nonnegative integers and let $w_1, \ldots, w_n$ be distinct real numbers.  
Fix distinct real numbers $x_1, \ldots, x_N$, where
\[
N \ = \ d_1 + \cdots + d_n + n - 1.
\]
Then there are polynomials $P_1, \ldots, P_n$ in $\R[X]$ of degrees $d_1, \ldots, d_n$ such that the function 
\[
F(x) 
\ = \ 
\sum\limits_{j=1}^n \hsx
P_j(x) e^{w_j x}
\]
has a simple zero at each point $x_1, \ldots, x_N$ and no other zeros.
\vspace{0.2cm}

[Note: \ 
This can be generalized by dropping the requirement that the $x_1, \ldots, x_N$ be distinct and incorporating multiplicities.]

\vspace{0.3cm}
\end{x}

\begin{x}{\small\bf \un{N.B.}}  \ 
The upper bound in \#1 is thus the best possible.
\end{x}
\vspace{0.3cm}

There is also an estimate in the complex domain.
\vspace{0.3cm}

\begin{x}{\small\bf LEMMA} \ 
Let $P_1, \ldots, P_n$ be nonzero polynomials in $\C[X]$ of degrees $d_1, \ldots, d_n$ and let $w_1, \ldots, w_n$ be distinct complex numbers.  Put
\[
\Omega
\ = \ 
\max\{\abs{w_1}, \ldots, \abs{w_n}\}.
\]
Then the number of zeros counting multiplicities of
\[
F(z) 
\ = \ 
\sum\limits_{j=1}^n \hsx
P_j(z) e^{w_j z}
\]
in the disc $\abs{z} \leq R$ is at most
\[
3(d_1 + \cdots + d_n + n - 1) + 4 R \Omega.
\]
\end{x}
\vspace{0.3cm}

\begin{x}{\small\bf NOTATION} \ 
If $f(z)$ is a function  continuous in $\ov{D}(R,w)$, put
\[
M(R, w, f) 
\ = \ 
\max\limits_{z \in \ov{D}(R,w)} \abs{f(z)}.
\]
\vspace{0.2cm}

[Note: \ 
Write
\[
M(R,f)
\]
if $w = 0$.]
\end{x}
\vspace{0.3cm}

\begin{x}{\small\bf NOTATION}\ 
If $f(z)$ is a function holomorphic in $D(R,w)$ and continuous in $\ov{D}(R,w)$, denote by
\[
N(R,w,f) 
\]
the number of zeros counting multiplicities of $f(z)$ in $\ov{D}(R,w)$.
\vspace{0.2cm}

[Note: \ 
Write
\[
N(R,f)
\]
if $w = 0$.]
\end{x}
\vspace{0.3cm}

\begin{x}{\small\bf RAPPEL} \ 
(Jensen) \ 
Let $R > 0$, $s > 1$ $-$then
\[
\int\limits_0^{s R} \hsx
\frac{N(r, w, f)}{r} dr 
\ = \ 
\frac{1}{2 \pi} \hsx 
\int\limits_0^{2 \pi} \hsx
\elln\bigg(\abs{f(w + s R e^{\sqrt{-1} \hsx \theta})}\bigg) d \theta 
- 
\elln(\abs{f(w)}).
\]
\end{x}
\vspace{0.3cm}

\begin{x}{\small\bf SUBLEMMA} \ 
Let \mR, $s$, $t \in \R_{> 0}$, $s > 1$, and  let $f \not\equiv 0$ be holomorphic in 
$D((st + s + t) R)$ and continuous in $\ov{D}((st + s + t)R)$ $-$then
\[
N(R,f) 
\ \leq \ 
\frac{1}{\elln(s)} \hsx 
\elln\bigg(
\frac{M((st + s + t)R, f)}{M(tR,f)}
\bigg).
\]
\vspace{0.2cm}

PROOF \ 
Choose $w \in \ov{D}(t R)$: $\abs{f(w)} = M(t R, f)$ (cf. \S25, \#3) $-$then $\abs{w} = t R$.  
So
\allowdisplaybreaks
\begin{align*}
z \in \ &\ov{D}(R)
\\[12pt]
&\implies
\abs{z - w} \ \leq \ \abs{z} + \abs{w}
\\[12pt]
&\hspace{2.4cm}
\leq \ R + tR 
\\[12pt]
&\hspace{2.4cm}
=\  (1 + t) R
\\[12pt]
&\implies
\\[12pt]
&\hspace{2cm}
\text{\textbullet} \quad
\ov{D}(R) \ \subset \ \ov{D}((1 + t)R,w)
\end{align*}
and
\allowdisplaybreaks
\begin{align*}
z \hsx \in \  &\ov{D}((st + s)R,w)
\\[12pt]
&\implies
\\[12pt]
&\hspace{2cm}
\abs{z} 
\ =\  \abs{z - w + w}
\\[12pt]
&\hspace{2.6cm}
\leq \ 
\abs{z - w} + \abs{w}
\\[12pt]
&\hspace{2.6cm}
\leq \ 
(st + s) R + tR
\\[12pt]
&\hspace{2.6cm}
=\ (st + s + t)R
\\[12pt]
&\implies
\\[12pt]
&\hspace{2.4cm}
\text{\textbullet} \quad
\ov{D}((st + s)R,w) \ \subset \ \ov{D}((st + s + t)R).  
\end{align*}
Next
\allowdisplaybreaks
\begin{align*}
N(R,w,f)\ 
&=\ 
\frac{1}{\elln(s)} 
\int\limits_R^{sR} \hsx \
\frac{N(R,w,f)}{r} dr 
\\[12pt]
&\leq \ 
\frac{1}{\elln(s)} 
\int\limits_0^{sR} \hsx \
\frac{N(r,w,f)}{r} dr 
\\[12pt]
&= \ 
\frac{1}{\elln(s)} 
\bigg[
\frac{1}{2 \pi} \hsx 
\int\limits_0^{2 \pi} \hsx 
\elln \bigg(\abs{f\big(w + sR e^{\sqrt{-1} \hsx \theta}\big)}\bigg) d \theta 
- 
\elln(\abs{f(w)}) 
\bigg]
\\[12pt]
&= \ 
\frac{1}{\elln(s)} 
\bigg[
\frac{1}{2 \pi} \hsx 
\int\limits_0^{2 \pi} \hsx 
\elln \bigg(\abs{f\big(w + sR e^{\sqrt{-1} \hsx \theta}\big)}\bigg) d \theta 
- 
\frac{1}{2 \pi} \hsx 
\int\limits_0^{2 \pi} \hsx 
\elln(\abs{f(w)}) d \theta 
\bigg]
\\[12pt]
&= \ 
\frac{1}{\elln(s)} 
\bigg[
\frac{1}{2 \pi} \hsx 
\int\limits_0^{2 \pi} \hsx 
\elln\bigg(
\frac{|f\big(w + sR e^{\sqrt{-1} \hsx \theta}\big)|}{\abs{f(w)}}
\bigg)
d \theta
\bigg].
\end{align*}
Take
\[
z \ = \ w + s R e^{\sqrt{-1} \hsx \theta}.
\]
Then
\allowdisplaybreaks
\begin{align*}
\abs{z - w} \ 
&=\ 
\abs{w + s R e^{\sqrt{-1} \hsx \theta} - w}
\\[12pt]
&=\ 
\abs{s R e^{\sqrt{-1} \hsx \theta}}
\\[12pt]
&=\ 
sR.
\end{align*}
Therefore
\[
N(R,w,f) 
\ \leq \ 
\frac{1}{\elln(s)} \hsx 
M\bigg(sR, w, \elln\bigg(\frac{\abs{f}}{\abs{f(w)}}\bigg) \bigg). 
\]
Spelled out
\allowdisplaybreaks
\begin{align*}
N(R,w,f) \ 
&=\ 
\abs{N(R,w,f)}
\\[12pt]
&\leq\ 
\frac{1}{\elln(s)} \hsx 
\frac{1}{2 \pi} \hsx
\abs{
\int\limits_0^{2 \pi} \hsx
\elln \bigg(\frac{|f(w + s R e^{\sqrt{-1} \hsx \theta})|}{\abs{f(w)}}\bigg)d \theta 
}
\\[12pt]
&\leq\ 
\frac{1}{\elln(s)} \hsx 
\frac{1}{2 \pi} \hsx
\int\limits_0^{2 \pi} \hsx
\abs{\elln \bigg(\frac{|f(w + s R e^{\sqrt{-1} \hsx \theta})|}{\abs{f(w)}}\bigg)}d \theta 
\\[12pt]
&\leq\ 
\frac{1}{\elln(s)} \hsx 
\frac{1}{2 \pi} \hsx
\int\limits_0^{2 \pi} \hsx
M\bigg(sR, w, \elln\bigg(\frac{\abs{f}}{\abs{f(w)}}\bigg) \bigg) d \theta 
\\[12pt]
&=\ 
\frac{1}{\elln(s)} \hsx  
M\bigg(sR, w, \elln\bigg(\frac{\abs{f}}{\abs{f(w)}}\bigg) \bigg).
\end{align*}
Finally
\[
N(R,f) 
\ \leq \ 
N((1 + t) R, w, f)
\]
since
\[
\ov{D}(R) \ \subset \ \ov{D}((1 + t)R, w).
\]
And working in the above with $(1 + t)R$ rather than just \mR, it follows that
\[
\tN((1 + t)R, w, f)
\]
is majorized by
\[
\frac{1}{\elln(s)} \hsx 
M\bigg(s(1 + t)R,  w, \elln\bigg(\frac{\abs{f}}{\abs{f(w)}}\bigg)\bigg)
\]
or still, by
\[
\frac{1}{\elln(s)} \hsx 
M\bigg((st + s)R,  w, \elln\bigg(\frac{\abs{f}}{M(t R,f)}\bigg)\bigg)
\]
which in turn is
\[
\leq \ \frac{1}{\elln(s)} \hsx 
M\bigg((st + s + t)R, \elln\bigg(\frac{\abs{f}}{M(t R,f)}\bigg)\bigg)
\]
because
\[
\ov{D}((st + s)R,w) \ \subset \ \ov{D}((st + s + t)R).
\]
Accordingly
\allowdisplaybreaks
\begin{align*}
N(R,f) \
&\leq\ 
\frac{1}{\elln(s)} \hsx 
M\bigg((st + s + t)R, \elln\bigg(\frac{\abs{f}}{M(t R,f)}\bigg)\bigg)
\\[12pt]
&\leq\ 
\frac{1}{\elln(s)} \hsx 
\elln\bigg(\frac{M((st + s + t)R,f)}{M(t R,f)}\bigg).
\end{align*}

Keep to the notation and assumptions of \#5 and set for simplicity
\[
D \ = \ 
\sum\limits_{j=1}^n \hsx d_j + n.
\]

\qquad {\small\bf FACT} \ 
Let \mR, $\gamma \in \R_{> 0}$, $\gamma > 1$ $-$then
\[
M(\gamma R, F) 
\ \leq \ 
\frac{\gamma^D - 1}{\gamma -1} \hsx 
e^{R \Omega (\gamma + 1)} M(R,F).
\]
[This technicality is dispensed with in the Appendix to this \S.] 
\vspace{0.3cm}

With this preparation, let us take up the proof of \#5.  
In the preceding, work with $tR$ rather than \mR, hence
\[
M(\gamma t R, F) 
\ \leq \ 
\frac{\gamma^D - 1}{\gamma -1} \hsx 
e^{t R \Omega (\gamma + 1)} M(t R,F).
\]
Now specialize and take
\[
\gamma \ = \ (st + s + t) / t.
\]
Then
\allowdisplaybreaks
\begin{align*}
\frac{\gamma^D - 1}{\gamma -1}\ 
&\leq\ 
\frac{1}{\gamma - 1}  \gamma^D
\\[12pt]
&=\ 
\frac{t}{st + s} \gamma^D
\\[12pt]
&=\ 
\frac{t}{s(t+1)}
\bigg(\frac{st + s + t}{t}\bigg)^D
\\[12pt]
&=\ 
\frac{t}{s(t+1)}
\bigg(\frac{st + s + t}{t}\bigg)^{D-1}
\bigg(\frac{st + s + t}{t}\bigg)
\\[12pt]
&=\ 
\frac{1}{s(t+1)}
\bigg(\frac{st + s + t}{t}\bigg)^{D-1}
(s (t + 1) + t)
\\[12pt]
&=\ 
\bigg(1 + \frac{t}{s(t+1)}
\bigg(\frac{st + s + t}{t}\bigg)^{D-1}
\\[12pt]
&\leq\ 
\bigg(1 + \frac{1}{s}\bigg) 
\bigg(\frac{st + s + t}{t}\bigg)^{D-1}.
\end{align*}
Therefore
\[
M((st + s + t)R, F) 
\ \leq \ 
\bigg(1 + \frac{1}{s}\bigg) 
\bigg(\frac{st + s + t}{t}\bigg)^{D-1}
e^{(st + s + 2t)R \Omega}
M(t R, F)
\]
\qquad\qquad $\implies$
\[
\frac{M((st + s + t)R, F)}{M(t R, F)} 
\ \leq \ 
\bigg(1 + \frac{1}{s}\bigg) 
\bigg(\frac{st + s + t}{t}\bigg)^{D-1}
e^{(st + s + 2t)R \Omega}
\]
\qquad\qquad $\implies$
\[
N(R,F) 
\ \leq \ 
\frac{1}{\elln(s)} \hsx
\elln\bigg(
\frac{M((st + s + t)R, F)}{M(t R, F)} 
\bigg)
\qquad \text{(cf. \#9)}
\]
\qquad\qquad $\implies$
\allowdisplaybreaks
\begin{align*}
N(R,F) \ 
&\leq \ 
\frac{1}{\elln(s)} \hsx
\bigg[
\elln\bigg(1 + \frac{1}{s}\bigg) + (D - 1) \elln\bigg(\frac{st + s + t}{t}\bigg) + (st + s + 2t) R \Omega
\bigg]
\\[12pt]
&\leq \ 
\frac{1}{\elln(s)} \hsx
\bigg[
\frac{1}{s} + (D - 1) \elln\bigg(\frac{st + s + t}{t}\bigg) + (st + s + 2t) R \Omega
\bigg].
\end{align*}
Into this relation insert $s = 5$ and $t = \ds\frac{1}{5}$.  
Toss the ``$\ds\frac{1}{s}$'' and note that
\[
\frac{\elln(31)}{\elln(5)} < 2.2 
\quad \text{and} \quad 
\frac{32}{5 \elln(5)} < 3.9
\]
giving
\[
N(R,F) \ \leq \ 3 (D-1) + 4 R \Omega,
\]
the assertion of \#5.
\end{x}
\vspace{0.3cm}

\begin{x}{\small\bf \un{N.B.}} \ 
One can replace the origin by any complex number $w$ and, upon consideration of $F(z - w)$, conclude that still
\[
N(R, w, F) \ \leq \ 3(D-1) + 4 R \Omega. 
\]
\end{x}
\vspace{0.3cm}

\[
\text{APPENDIX}
\]
\vspace{0.5cm}

Recall the setup of \#5.  
Thus, as there, let $P_1, \ldots, P_n$ be nonzero polynomials in $\C[X]$ of degrees $d_1, \ldots, d_n$ and let 
$w_1, \ldots, w_n$ be distinct complex numbers.  
Put
\[
\Omega
\ = \ 
\max\{\abs{w_1}, \ldots, \abs{w_n}\}, 
\quad
D 
\ = \ 
\sum\limits_{j=1}^n \hsx
d_j + n,
\]
and form
\[
F(z) 
\ = \ 
\sum\limits_{j=1}^n \hsx
P_j(z) e^{w_j z}.
\]
\vspace{0.3cm}

{\small\bf PREFACT} \ 
Fix a point $z_0 \in \C$ $-$then
\[
\abs{F(z_0)}
\ \leq \ 
e^{(\abs{z_0} + 1) \Omega}
\bigg(\sum\limits_{k=0}^{D-1} \hsx \abs{z_0}^k\bigg) \hsx
\max\limits_{j = 1, \ldots, D} \hsx \abs{\frac{F^{(j-1)}(0)}{(j-1)!}}.
\]
\vspace{0.3cm}

{\small\bf FACT} \ 
Let $R$, $\gamma \in \R_{> 0}$, $\gamma > 1$ $-$then
\[
M(\gamma R, F) 
\ \leq \ 
\frac{\gamma^D - 1}{\gamma -1} \hsx 
e^{R \Omega (\gamma+1)} \hsx
M(R, F).
\]

PROOF \ 
Choose $z_0 \ (\abs{z_0} = \gamma)$:
\[
\abs{F(z_0 R)}
\ = \ 
\max\limits_{\abs{z} \leq \gamma R} \hsx \abs{F(z)}.
\]
Consider 
\[
G(z) 
\ = \ 
F(z R) 
\ = \ 
\sum\limits_{j=1}^n \hsx
P_j(z R) e^{w_j R z}.
\]
Then by the above applied to \mG (hence now it is a question of $w_j R$ rather than $w_j$ and it is also a question of 
$R \Omega$ rather than $\Omega$) we have
\[
\abs{G(z_0)}
\ < \ 
e^{(\gamma+1) R \Omega} \hsx
\bigg(\sum\limits_{k=0}^{D-1} \hsx \gamma^k \bigg) \hsx
\max\limits_{j = 1, \ldots, D} \hsx \abs{\frac{R^{j-1} F^{(j-1)}(0)}{(j-1)!}}.
\]
But
\[
\sum\limits_{k=0}^{D-1} \hsx \gamma^k 
\ = \ 
\frac{\gamma^D - 1}{\gamma -1}
\]
and, thanks to Cauchy's inequality,
\[
\max\limits_{j = 1, \ldots, D} \hsx \abs{\frac{R^{j-1} F^{(j-1)}(0)}{(j-1)!}}
\ \leq \ 
\max\limits_{\abs{z} \leq R} \hsx \abs{F(z)}.
\]
Therefore
\allowdisplaybreaks
\begin{align*}
M(\gamma R, F) \ 
&=\ 
\max\limits_{\abs{z} \leq \gamma R} \hsx \abs{F(z)}
\\[10pt]
&=\ 
\abs{F(z_0 R)}
\\[10pt]
&=\ 
\abs{G(z_0)}
\\[10pt]
&\leq\ 
\frac{\gamma^D - 1}{\gamma -1} \hsx 
e^{R \Omega (\gamma+1)} \hsx
\max\limits_{\abs{z} \leq R} \hsx \abs{F(z)}
\\[10pt]
&=\ 
\frac{\gamma^D - 1}{\gamma -1} \hsx
e^{R \Omega (\gamma+1)} \hsx
M(R,F).
\end{align*}
\vspace{0.2cm}

{\small\bf REMARK} \ 
The estimate figuring in \#5 can be sharpened to
\[
N(R,F) 
\ \leq \ 
2(D - 1) + \frac{4}{\pi} R \Omega.
\]

%% file: _27_gelfond_schneider_setting_stage.tex
\chapter{
$\boldsymbol{\S}$\textbf{27}.\quad  GELFOND-SCHNEIDER: \ SETTING THE STAGE}
\setlength\parindent{2em}
\setcounter{theoremn}{0}
\renewcommand{\thepage}{\S27-\arabic{page}}

\ \indent 
Recall the claim:
\vspace{0.3cm}

\begin{x}{\small\bf THEOREM} \ 
If $\alpha \neq 0, 1$ is algebraic and if $\beta \notin \Q$ is algebraic, then $\alpha^\beta$ is transcendental.

\vspace{0.2cm}

[Note: \ 
Here $\alpha^\beta$ is the principal power (cf. \S23, \#14):
\[
\alpha^\beta \ = \ e^{\beta \hsx \Log \alpha}.]
\]
\end{x}
\vspace{0.3cm}

Methodology: \ Assume that $\alpha \neq 0, 1$ is algebraic, that $\beta$ is algebraic, and that $\alpha^\beta$ is algebraic $-$then 
the theorem will follow if it can be shown that $\beta \in \Q$.
\vspace{0.3cm}

\begin{x}{\small\bf NOTATION}  \ 
Given a positive odd integer $N \gg 0$, put
\[
L \ = \ N^8, 
\quad
S \ = \ \frac{1}{2} (N^4 - 1),
\]
and 
\[
\begin{cases}
\ L_0 \ = \ N^6 - 1\\
\ L_1 \ = \ N^2 - 1
\end{cases}
.
\]
\vspace{0.2cm}

[Note: \ 
Restricting \mN to be odd guarantees that \mS is an integer.]
\end{x}
\vspace{0.3cm}

\begin{x}{\small\bf LEMMA} \ 
\[
L \ = \ (L_0 + 1)(L_1 + 1) \ = \ (2S + 1)^2.
\]
\vspace{0.2cm}

PROOF \ 
\[
\begin{cases}
\ L_0  + 1\ = \ N^6 \\[8pt]
\ L_1 + 1\ = \ N^2
\end{cases}
\implies (L_0 + 1) (L_1 + 1) \ = \  N^8.
\]
And
\[
(2S + 1)^2 \ = \ (N^4 - 1 + 1)^2 \ = \ N^8.
\]
\end{x}
\vspace{0.3cm}

During the ensuing analysis, there will emerge a positive absolute constant \mC.
\vspace{0.5cm}

\begin{x}{\small\bf LEMMA} \ 
Given $C \in \R_{> 0}$, $\exists \ N_0(C) \gg 0$ : $\forall \ N > N_0(C)$, 
\[
C \hsx L_0 \hsx \elln(S) \ \leq \ L
\quad \text{and} \quad
C \hsx L_1 \hsx S  \ \leq \ L.
\]
\end{x}
\vspace{0.3cm}

\begin{x}{\small\bf \un{N.B.}}  \ 
Therefore
\allowdisplaybreaks
\begin{align*}
C \hsx L \hsx(L_0 \hsx \elln(S) + L_1 S)\ 
&=\ 
L(C \hsx L_0 \elln(S)  + C \hsx L_1  \hsx S)
\\[12pt]
&\leq\ 
L(L) + L(L)
\\[12pt]
&=\ 
2 L^2.
\end{align*}
\end{x}
\vspace{0.3cm}

\qquad \textbullet \quad 
Choose an ordering of the integral pairs $(s_1, s_2)$ with $\abs{s_1} \leq S$ and $\abs{s_2} \leq S$,  i.e., 
$(s_1, s_2) \in \Z^2$ and $-S \leq s_1, s_2, \leq S$.
\vspace{0.2cm}

[Note: \ 
There are $S + (S + 1)$ choices for $s_1$ and $S + (S + 1)$ choices for $s_2$, hence there are all told
\[
(2S + 1) \times (2S + 1) \ = \ (2S + 1)^2 \ = \ L
\]
integral pairs $(s_1, s_2)$.]
\vspace{0.2cm}

\qquad \textbullet \quad 
Choose an ordering of the integral pairs
\[
(u,v) \in \{0, \ldots, L_0\} \times \{0, \ldots, L_1\}.
\]
\vspace{0.2cm}

[Note: \ 
There are $L_0 + 1$ choices for $u$ and $L_1 + 1$ choices for $v$, hence there are all told
\[
(L_0 + 1) (L_1 + 1) \ = \  L \ \ (= (2S + 1)^2)
\]
choices for $(u,v)$.]
\vspace{0.3cm}

\begin{x}{\small\bf NOTATION} \ 
Introduce an $L \times L$ matrix $\sM$ via the prescription 
\[
\sM 
\ = \ \big((s_1(i) + s_2(i) \beta)^{u(j)}\big(\alpha^{s_1(i) + s_2(i) \beta}\big)^{v(j)}\big)
\]
and let 
\[
\Delta \ = \ \det (\sM).
\]

[Note: \ 
$j$ is the column index and $i$ is the row index.]
\end{x}
\vspace{0.3cm}

\begin{x}{\small\bf \un{N.B.}} \ 
The orderings for the columns and rows has not been explicated but a change in these orderings simply changes matters by a factor 
$\pm 1$, which has no effect on the absolute value $\abs{\Delta}$ of $\Delta$.
\end{x}
\vspace{0.3cm}

Define a function of one complex variable $z$ by
\[
f_j(z) 
\ = \ 
z^{u(j)} \alpha^{v(j) z} \qquad (1 \leq j \leq L)
\]
and put
\[
\zeta_i 
\ = \ 
s_1(i) + s_2(i) \beta \qquad (1 \leq i \leq L).
\]
\vspace{0.3cm}

\begin{x}{\small\bf SUBLEMMA} \ 
$\forall$ complex numbers $z_1$, $z_2$, 
\[
\abs{e^{z_1 z_2}} \ = \ e^{\Rex(z_1 z_2)} 
\ \leq \ 
e^{\abs{z_1 z_2}} \ = \ e^{\abs{z_1}}e^{\abs{z_2}}.
\]
\end{x}
\vspace{0.3cm}

\begin{x}{\small\bf LEMMA} \ 
$\forall \ R \in \R_{> 0}$, 
\[
M(R,f_j) \ \leq \ R^{u(j)} e^{v(j) R \abs{\Log \alpha}}.
\]
\vspace{0.2cm}

PROOF \ 
For by definition,
\[
\alpha^{v(j) z} \ = \ \exp(v(j) \hsx z \hsx \Log \alpha).
\]
Therefore
\allowdisplaybreaks
\begin{align*}
\elln(M(R, f_j)) \ 
&\leq\ 
u(j) \elln(R) + v(j) R \abs{\Log \alpha}
\\[12pt]
&\leq\ 
L_0 \elln(R) + L_1 R \abs{\Log \alpha}.
\end{align*}
\end{x}
\vspace{0.3cm}

\begin{x}{\small\bf RAPPEL} \ 
In the notation of \S25, \#5, 
\[
\abs{\Delta}
\ \leq\ 
\bigg(\frac{R}{r}\bigg)^{-L(L-1)/2} \hsx
L! \hsx
\prod\limits_{j=1}^L \hsx
\abs{f_j}_R.
\]
\vspace{0.2cm}

[Note: \ 
The symbols $\abs{f_j}_R$ and $M(R, f_j)$ mean one and the same thing.]
\end{x}
\vspace{0.3cm}

In the case at hand, 
\[
\Delta \ = \ \det(f_j, (\zeta_i)),
\]
thus the foregoing generality is applicable.
\vspace{0.3cm}

\qquad \textbullet \quad 
Take $r = S(1 + \abs{\beta})$ and note that
\allowdisplaybreaks
\begin{align*}
\abs{\zeta_i} \ 
&=\ 
\abs{s_1(i) + s_2(i) \beta}
\\[12pt]
&\leq\ 
\abs{s_1(i)} + \abs{s_2(i) \beta}
\\[12pt]
&\leq\ 
S + S \abs{\beta}
\\[12pt]
&=\ 
S(1 + \abs{\beta}).
\end{align*}

\qquad \textbullet \quad 
Take $R = e^2 r$ and note that
\allowdisplaybreaks
\begin{align*}
\bigg(\frac{R}{r}\bigg)^{-L(L-1)/2} \ 
&=\ 
\bigg(\frac{e^2 r}{r}\bigg)^{-L(L-1)/2} 
\\[12pt]
&=\ 
e^{-L(L - 1)}.
\end{align*}
\begin{x}{\small\bf LEMMA} \ 
\[
\abs{\Delta}
\ \leq\ 
e^{-L(L - 1)} \hsx
L! \hsx
\prod\limits_{j=1}^L \hsx
M(R, f_j),
\]
where
\[
R \ = \ e^2 S (1 + \abs{\beta}).
\]
\end{x}
\vspace{0.3cm}

\begin{x}{\small\bf LEMMA} \ 
\[
\elln(\abs{\Delta}) \ \leq \  -\frac{L^2}{2}.
\]
\vspace{0.2cm}

PROOF \ 
Starting with \#11,
\allowdisplaybreaks
\begin{align*}
\elln(\abs{\Delta}) \ 
&\leq\ 
-L(L - 1) + \elln(L!) + \sum\limits_{j=1}^L  \hsx \elln(M(R, f_j))
\\[12pt]
&\leq\ 
- L^2 + L + L \elln(L) + L \max\limits_{1 \leq j \leq L} \elln(M(R, f_j))
\\[12pt]
&\leq\ 
-L^2 + L(1 + \elln(L) + L_0 \elln(R) + L_1 R \abs{\Log \alpha}).
\end{align*}

\allowdisplaybreaks
\begin{align*}
\text{\textbullet} \quad 
L L_0 \elln(R) \ 
&=\ 
L L_0 \elln(e^2 S (L + \abs{\beta})) \hspace{2cm}
\\[12pt]
&=\ 
L L_0 (\elln(e^2) + \elln(S) + \elln (1 + \abs{\beta}))
\\[12pt]
&=\ 
L L_0 \elln(e^2) + L L_0 \elln (1 + \abs{\beta}) + L L_0 \elln(S)
\\[12pt]
&\leq\ 
C_1 L L_0 \elln(S).
\end{align*}

\allowdisplaybreaks
\begin{align*}
\text{\textbullet} \quad
L L_1 R \abs{\Log \alpha} \ 
&=\ 
L L_1 e^2 S(1 + \abs{\beta}) \abs{\Log \alpha} \hspace{2.4cm}
\\[12pt]
&=\ 
e^2 (1 + \abs{\beta}) \abs{\Log \alpha} L L_1 S
\\[12pt]
&\leq\ 
C_2 L L_1 S.
\end{align*}

Therefore
\allowdisplaybreaks
\begin{align*}
-L^2 + L(1 + &\elln(L) + L_0 \elln(R) + L_1 R \abs{\Log \alpha})
\\[12pt]
&\leq
-L^2 + L(1 + \elln(L)) + C_1 L L_0 \elln(S) + C_2 L L_1 S
\\[12pt]
&\leq
-L^2 + C_3(L L_0 \elln(S) + L L_1S) + C_1 L L_0 \elln(S) + C_2 L L_1 S
\\[12pt]
&\leq
-L^2 + C_4(L L_0 \elln(S) + L L_1 S),
\end{align*}
the positive absolute constant $C_4$ being independent of $N \gg 0$.  
Take now $C \geq 4 C_4$ and unravel the data:
\allowdisplaybreaks
\begin{align*}
\elln(\abs{\Delta})
&\leq
-L^2 + C_4(L L_0 \elln(S) + L L_1 S)
\\[12pt]
&\leq
- L^2 + \frac{C}{4} (L L_0 \elln(S) + L L_1 S)
\\[12pt]
&=\
-L^2 + \frac{1}{4} C L (L_0 \elln(S) + L_1 S)
\\[12pt]
&\leq
-L^2 + \frac{1}{4}  (2 L^2) \qquad \text{(cf. \#5)}
\\[12pt]
&=\
-L^2 + \frac{L^2}{2} 
\\[12pt]
&=\
-\frac{L^2}{2},
\end{align*}
thereby completing the proof.
\end{x}
\vspace{0.3cm}

\begin{x}{\small\bf LEMMA} \ 
\[
\elln(\abs{\Delta}) \ \geq \ - \frac{L^2}{3}
\]
if $\Delta \neq 0$.
\end{x}
\vspace{0.3cm}

\begin{x}{\small\bf \un{N.B.}}  \ 
Granted this, we have a contradiction: $\ds\frac{1}{3} \geq \ds\frac{1}{2}$.  
Thus the conclusion is that
\[
\Delta \ = \ \det(\sM) \ = \ 0.
\]
\end{x}
\vspace{0.3cm}

Bearing in mind that for \#13, $\Delta \neq 0$, fix $T \in \N$ such that $T \alpha$, $T \beta$, and $T \alpha^\beta$ are algebraic integers (recall that $\forall \ x \in \Qbar$, $D_x$ is a nonzero ideal of $\Z$ (cf. \S14)) $-$then $T\struty^{L_0 + 2 L_1 S}$ times any element of the matrix 
$\sM$ is an algebraic integer.  
Moreover the algebraic integer
\[
T\strutx^{L(L_0 + 2 L_1 S)} \hsx \Delta
\]
is a zero of a monic polynomial of degree $d$, where $d$ is at most the product of the degrees of the minimal polynomials of 
$\alpha$, $\beta$, and $\alpha^\beta$.
\vspace{0.3cm}

\begin{x}{\small\bf SUBLEMMA} \ 
\[
H(\Delta) 
\ \leq \ 
L! \hsx S^{L_0 L} 
(1 + H(\beta))^{L_0 L} 
(1 + H(\alpha))^{L_1 L S} 
(1 + H(\alpha^\beta))^{L_1 L S}.
\]
\vspace{0.2cm}

[Note: \ 
The house of an algebraic number $x \neq 0$ is, by definition, the maximum of the absolute values of $x$ and its conjugates 
(see the Appendix to \#14, in particular the result formulated there, to be used infra).]
\end{x}
\vspace{0.3cm}

On the other hand, 
\[
\Delta \ \neq \ 0 
\quad \text{and} \quad
T\strutx^{L(L_0 + 2 L_1 S)} 
\in 
D_\Delta, 
\]
hence
\allowdisplaybreaks
\begin{align*}
\abs{\Delta} \ 
&\geq \ 
T^{-dL(L_0 + 2L_1 S)} H(\Delta)^{1 - d}
\\[12pt]
&\geq \ 
T^{-dL(L_0 + 2L_1 S)} H(\Delta)^{-d}
\end{align*}
\qquad\qquad $\implies$
\allowdisplaybreaks
\begin{align*}
\abs{\Delta} \ 
&\geq \ 
T^{-dL(L_0 + 2L_1 S)} 
(L!)^{-d_S^{-d L_0 L}}
\\[12pt]
&\hspace{1cm}
\times 
(1 + H(\beta))^{-d L_0 L}
(1 + H(\alpha))^{-d L_1 L S}
(1 + H(\alpha^\beta))^{-d L_1 L S}
\end{align*}
\qquad\qquad $\implies$
\allowdisplaybreaks
\begin{align*}
\elln(\abs{\Delta}) \ 
&\geq \ 
-d L (L_0 + 2 L_1 S) \elln(T) - d L \elln(L) - d L_0 L \elln(S)
\\[12pt]
&\hspace{1cm}
- d L_0 L \elln(1 + H(\beta)) - d L_1 L S \elln(1 + H(\alpha)) - dL_1 L S \elln(1 + H(\alpha^\beta))
\end{align*}
\qquad\qquad $\implies$
\allowdisplaybreaks
\begin{align*}
\elln(\abs{\Delta}) \ 
&\geq \ 
-K_1 L(L_0 + \elln(L) + L_0 \elln(S) + L_1 S)
\end{align*}
\qquad\qquad $\implies$
\allowdisplaybreaks
\begin{align*}
\elln(\abs{\Delta}) \ 
&\geq \ 
-K_2 L (L_0 \elln(S) + L_1S),
\end{align*}
the positive absolute constant $K_2$ being independent of $N \gg 0$.  
Take now $C \geq 6 K_2$ $-$then
\allowdisplaybreaks
\begin{align*}
\elln(\abs{\Delta})
&\geq 
-\frac{C}{6} L(L_0 \elln(S) + L_1 S)
\\[12pt]
&=\ 
\frac{1}{6}(-C L(L_0 \elln(S) + L_1 S))
\\[12pt]
&\geq 
\frac{1}{6} (-2 L^2) \qquad \text{(cf. \#5)}
\\[12pt]
&=\ 
-\frac{L^2}{3},
\end{align*}
the assertion of \#13.

%% file: _28_gelfond_schneider_execution.tex
\chapter{
$\boldsymbol{\S}$\textbf{28}.\quad  GELFOND-SCHNEIDER: \ EXECUTION}
\setlength\parindent{2em}
\setcounter{theoremn}{0}
\renewcommand{\thepage}{\S28-\arabic{page}}

\ \indent 

Under the assumption that $\alpha \neq 0, 1$ is algebraic, that $\beta$ is algebraic, and that $\alpha^\beta$ is algebraic, the central conclusion of \S27 is that
\[
\Delta \ = \ \det(f_j (\zeta_i)) \ = \ 0,
\]
the goal being to show that $\beta \in \Q$.
\vspace{0.2cm}

Proceeding, assume momentarily that $\alpha$, $\beta$, $\alpha^\beta \in \Qbar \hsx \cap \hsx \R$ $(\alpha > 0)$, 
hence all data is real 
and the columns of the matrix $(f_j(\zeta_i))$ are linearly dependent over $\R$, thus there exist real numbers 
$b_1, \ldots, b_L$ not all zero such that
\[
\sum\limits_{j=1}^L \hsx
b_j f_j(\zeta_i) 
\ = \ 0 \qquad (1 \leq i \leq L).
\]
But
\[
f_j(\zeta_i) \ = \ \zeta_i^{u(j)} \hsx \alpha^{v(j) \zeta_i},
\]
so
\[
\sum\limits_{j=1}^L \hsx
b_j \hsx \zeta_i^{u(j)} \hsx \alpha^{v(j) \zeta_i}
\ = \ 0 \qquad (1 \leq i \leq L)
\]
or still,
\[
\sum\limits_{v=0}^{L_1} \hsx
\bigg(
\sum\limits_{u=0}^{L_0} \hsx
b_{(L_0 + 1)v + u + 1} \hsx \zeta_i^u
\bigg)^{v \zeta_i}
\ = \ 0.
\]
Introduce
\[
a_v(t) 
\ = \ 
\sum\limits_{u=0}^{L_0} \hsx
b_{(L_0 + 1)v + u + 1} \hsx t^u,
\]
where $t \in \R$, and consider
\[
\sum\limits_{v=0}^{L_1} \hsx
a_v(t)  e^{w_v t} \qquad (w_v = v \hsx \Log \alpha).
\]
Since
\[0
\ = \ 
\sum\limits_{v=0}^{L_1} \hsx
a_v(\zeta_i) e^{w_v \zeta_i} 
\qquad (\zeta_i = s_1(i) + s_2(i) \beta),
\]
it follows that each of the \mL values of $\zeta_i$ is a zero of
\[
A(t) 
\ \equiv \ 
\sum\limits_{v=0}^{L_1} \hsx
a_v(t) e^{w_v t}.
\]
\vspace{0.2cm}

At this point, \#1 of \S26 is applicable:
\vspace{0.2cm}

\qquad \textbullet \quad 
The degree of $a_v(t)$ is $\leq L_0$.
\vspace{0.2cm}

\qquad \textbullet \quad 
The $w_v$ are distinct real numbers.
\vspace{0.2cm}

\qquad \textbullet \quad 
The sum defining $A(t)$ consists of $L_1 + 1$ polynomials.
\vspace{0.2cm}

Accordingly $A(t)$ has at most
\[
L_0(L_1 + 1) + (L_1 + 1) - 1
\]
real zeros counting multiplicities.  
And:
\allowdisplaybreaks
\begin{align*}
L_0(L_1 + 1) + (L_1 + 1) - 1 \ 
&=\ 
L_0 L_1 + L_0 + L_1  + 1 - 1
\\[12pt]
&=\ 
(L_0 + 1)(L_1 + 1) - 1
\\[12pt]
&=\ 
L - 1 \qquad \text{(cf. \S27, \#3)}
\\[12pt]
&<\
L.
\end{align*}
Consequently two of the $\zeta_i$ must be the same, so
\[
s_1(i) + s_2(i) \beta
\ = \ 
s_1(i^\prime) + s_2(i^\prime) \beta
\]
for some $i$, $i^\prime$ with $1 \leq i < i^\prime \leq L$.  
However, since the pairs $(s_1(i), s_2(i))$ and 
$(s_1(i^\prime), s_2(i^\prime))$ are distinct, either
\[
\beta
\ = \ 
\frac{s_1(i) - s_1(i^\prime)}{s_2(i^\prime) - s_2(i)}
\quad \text{if} \quad
s_2(i^\prime) \neq s_2(i)
\]
or
\[
\frac{1}{\beta}
\ = \ 
\frac{s_2(i^\prime) - s_2(i)}{s_1(i) - s_1(i^\prime)}
\quad \text{if} \quad
s_1(i) \neq s_1(i^\prime).
\]
in any event, $\beta$ is rational \ldots \ .
\vspace{0.3cm}

To discuss the general case, it is necessary to elaborate on what has been said in \S27.
\vspace{0.5cm}

\qquad \qquad \un{Step 1:} \quad
Redefine \mS and replace $\ds\frac{1}{2} (N^4 - 1)$ by $2N^4$ $-$then
\[
\frac{S}{2} 
\ = \ 
N^4 
\implies
\frac{S^2}{4} 
\ = \ 
N^8 = L.
\]
And
\allowdisplaybreaks
\begin{align*}
(2S + 1)^2 \ 
&=\ 
4S^2 + 4S + 1
\\[12pt]
&=\ 
16N^8 + 8N^4 + 1
\\[12pt]
&>\ 
16N^8
\\[12pt]
&=\ 
16L
\\[12pt]
&>\ 
L.
\end{align*}

\qquad \qquad \un{Step 2:} \quad
Define the $(2S + 1)^2 \times L$ matrix $\sM$  as in \S27 and note that all the $L \times L$ submatrices of $\sM$ have determinant zero, as can be gleaned from the argumentation 
used there.

\qquad \qquad \un{Step 3:} \quad
The columns of the matrix $\sM$ are linearly dependent over $\C$, thus there exist complex numbers 
$b_1, \ldots, b_L$ not all zero such that 
\[
\sum\limits_{j=1}^L \hsx b_j f_j (\zeta_i) 
\ = \ 
0
\qquad (i \in \{1, \ldots, (2S + 1)^2\}).
\]

\qquad \qquad \un{Step 4:} \quad
Introduce as before
\[
A(t) 
\ \equiv \ 
\sum\limits_{v=0}^{L_1} \hsx
a_v(t) e^{w_v t} 
\]
and observe that
\[
A(\zeta_i) = 0 \qquad (i \in \{1, \ldots, (2S + 1)^2\}).
\]

Owing to \S26, \#5, 
\[
N(R,A) \ \leq \ 3 (D - 1) + 4 R \Omega
\]
or better, its improvement
\[
N(R,A) \ \leq \ 2 (D - 1) + \frac{4}{\pi} R \Omega,
\]
as noted in the Appendix to \S26.  Here
\[
D \ \leq \ L_0 (L_1 + 1) + (L_1 + 1) \ = \ L.
\]
And
\[
\zeta_i \ = \ s_1(i) + s_2(i) \beta,
\]
where a priori $\beta$ is complex and $\abs{s_1}$, $\abs{s_2} \leq S$, the choice
\[
R \ = \  S(1 + \abs{\beta})
\]
ensures that the disc of radius \mR centered at the origin contains all the points $\zeta_i$.  
In addition
\allowdisplaybreaks
\begin{align*}
\Omega \ 
&=\ 
\max\limits_{v= 0, \ldots, L_1} \abs{w_v} 
\\[12pt]
&=\ 
\max\limits_{v= 0, \ldots, L_1} \abs{v \Log \alpha}
\\[12pt]
&=\ 
L_1 \abs{\Log \alpha}.
\end{align*}
Therefore
\[
N(R,A) 
\ \leq \ 
2(L - 1) + \frac{4}{\pi} S(1 + \abs{\beta}) L_1 \abs{\Log \alpha}
\]
or still, 
\[
N(R,A) 
\ \leq \ 
2(L - 1) + K S L_1, 
\]
where
\[
K \ = \ \frac{4}{\pi}  (1 + \abs{\beta}) \abs{\Log \alpha}.
\]
But:
\allowdisplaybreaks
\begin{align*}
&\text{\textbullet} \quad
2(L - 1) < 2 L \ = \ 2 \bigg(\frac{S^2}{4}\bigg) = \frac{S^2}{2} \hspace{6.5cm}
\\[12pt]
&\text{\textbullet} \quad
K S L_1 = K(2N^4) (N^2 - 1)
\\[12pt]
&\hspace{1.75cm}
< \ 2 K N^6
\\[12pt]
&\hspace{1.75cm}
< \ N^8 \qquad (N \gg 0)
\\[12pt]
&\hspace{1.75cm}
= \ \frac{S^2}{4}
\\[12pt]
&\implies
\\[12pt]
&\hspace{.75cm}
N(R,A) \ < \  \frac{S^2}{2} + \frac{S^2}{4}
\\[12pt]
&\hspace{2.3cm}
= \ \frac{3}{4}S^2
\\[12pt]
&\hspace{2.3cm}
< \ (2S + 1)^2.
\end{align*}
Since \mA admits $(2 S + 1)^2$ zeros $\zeta_i$, two of them must be the same, forcing in the end the rationality of $\beta$.


%% file: _29_schneider_lang_criterion.tex
\chapter{
$\boldsymbol{\S}$\textbf{29}.\quad  THE SCHNEIDER-LANG CRITERION}
\setlength\parindent{2em}
\setcounter{theoremn}{0}
\renewcommand{\thepage}{\S29-\arabic{page}}

\ \indent 
Fix an algebraic number field $\K$.
\vspace{0.2cm}

[Note: \ Therefore $\K$ is a subfield of $\C$ which, when considered as a vector space over $\Q$, is finite dimensional, the dimension being denoted 
$[\K:\Q]$ and called the 
\un{degree} 
\index{degree (of an algebraic number field)}
of $\K$ over $\Q$.]
\vspace{0.3cm}

\begin{spacing}{1.50}
\begin{x}{\small\bf THEOREM} \ 
Let $f_1$, $f_2$ be entire functions of finite strict orders $\leq \beta_1$, $\leq \beta_2$.  
Assume: \ $f_1$, $f_2$ are algebraically independent over $\C$ and that the derivatives 
$\ds\frac{d}{dz} f_1$, $\ds\frac{d}{dz} f_2$ belong to the ring $\K[f_1, f_2]$ 
(i.e., can be written as polynomials in $f_1$, $f_2$) $-$then the set
\[
S \ = \ \{w \in \C : f_1(w), \ f_2(w) \in \K\}
\]
is finite.
\end{x}
\end{spacing}
\vspace{0.3cm}

There are two ``canonical'' examples that illustrate this criterion.
\vspace{0.3cm}

\begin{x}{\small\bf APPLICATION} \ 
Schneider-Lang $\implies$ Hermite-Lindemann. 
\vspace{0.3cm}

I.e.: If $\alpha$ is a nonzero algebraic number, then $e^\alpha$ is transcendental (cf. \S21, \#4).
\vspace{0.3cm}

[Suppose instead that $e^\alpha$ is algebraic, 
let $\K = \Q(\alpha, e^\alpha)$, and take 
$f_1(z) = z$ $(\rho_1 = 0)$, $f_2(z) = e^z$ $(\rho_2 = 1)$ 
(which are algebraically independent over $\C$ (cf. \S20, \#18)).  
Since it is clear that
\[
\frac{d}{dz} z, \ 
\frac{d}{dz} e^z 
\in \K[f_1(z), f_2(z)],
\]
the assumptions of \#1 are satisfied.  
On the other hand, $\forall \ n \in \N$, 
\[
f_1(n \alpha) = n \alpha \in \K, \quad 
f_2(n \alpha) = e^{n \alpha} \in \K,
\]
an infinite set of conditions, from which a contradiction.]
\end{x}
\vspace{0.3cm}

\begin{x}{\small\bf APPLICATION} \ 
Schneider-Lang $\implies$ Gelfond-Schneider. 
\vspace{0.2cm}

I.e.: If $\alpha \neq 0, 1$ is algebraic and if $\beta \in \Q$ is algebraic, 
then $\alpha^\beta$ is transcendental (cf. \S24, \#1).
\vspace{0.2cm}

[Suppose instead that $\alpha^\beta$ is algebraic, 
let $\K = \Q(\alpha, \beta, \alpha^\beta)$, and take 
$f_1(z) = e^z$ $(\rho_1 = 1)$, $f_2(z) = e^{\beta z}$ $(\rho_2 = 1)$ $-$then 
$f_1(z)$, $f_2(z)$ are algebraically independent over $\C$ $(\beta \notin \Q)$ (cf. \S20, \#26).  
Moreover
\[
\frac{d}{dz} f_1 = f_1, \quad \frac{d}{dz} f_2 = \beta f_2,
\]
so $\K[f_1(z),f_2(z)]$ is closed under differentiation, thus in view of \#1 there are but finitely many points $w \in \C$ such that 
$f_1(w) \in \K$ and $f_2(w) \in \K$.  
But for all $k = 1, 2, \ldots,$
\[
f_1(k \hsx  \Log \alpha) = \alpha^k \in \K
\quad \text{and} \quad
f_2(k \hsx \Log \alpha) = (\alpha^\beta)^k \in \K,
\]
an infinite set of conditions, from which a contradiction.]
\end{x}
\vspace{0.3cm}

\begin{x}{\small\bf REMARK} \ 
The objective is to show that the set \mS figuring in \#1 is finite.  
In fact, it will turn out that the cardinality of \mS is bounded by
\[
(\rho_1 + \rho_2) \hsx [\K:\Q].
\]
As for the proof, we shall not provide all the details but will say enough to
render the whole affair believable.
\end{x}
\vspace{0.3cm}

Let $N \gg 0$ be a positive integer.
\vspace{0.3cm}

\begin{x}{\small\bf NOTATION} \ 
Put
\[
R_1 \ (= R_1(N))  
\ = \ 
\bigg[N\strutx^{\frac{\raisebox{.1cm}{$\scriptstyle{\rho_2}$}}{\rho_1 + \rho_2}} (\elln(N))^{1/2}\bigg]
\]
and
\[
R_2 \ (= R_2(N))  
\ = \ 
\bigg[N\strutx^{\frac{\raisebox{.1cm}{$\scriptstyle{\rho_1}$}}{\rho_1 + \rho_2}} (\elln(N))^{1/2}\bigg].
\]
\end{x}
\vspace{0.3cm}

\begin{x}{\small\bf \un{N.B.}} \ 
\begin{align*}
R_1 \hsx R_2 \ 
&\leq \ 
N\strutx^{\frac{\rho_2}{\rho_1 + \rho_2}} (\elln(N))^{1/2} 
\ 
N\strutx^{\frac{\rho_1}{\rho_1 + \rho_2}} (\elln(N))^{1/2} 
\\[12pt]
&= \ N \elln(N).
\end{align*}
Therefore
\[
(R_1 + 1)(R_2 + 1)
\ \geq \ 
N \elln(N).
\]

[Note: \ 
If $C \in \R_{> 0}$, then
\[
N \elln(N) + CN \ \leq \ 2 N \elln(N)
\]
provided \mN is large enough: 
\[
N \gg 0 
\implies
\frac{N}{N \elln(N)} 
\ < \ 
\frac{1}{C}.]
\]
\end{x}
\vspace{0.3cm}

Let $w_1, \ldots, w_r$ be elements of \mS.
\vspace{0.3cm}

\begin{x}{\small\bf SUBLEMMA}\ 
There exists a nonzero polynomial $P_N \in \Z[X_1, X_2]$ whose degree 
w.r.t. $X_1$ is $\leq R_1$ and whose degree w.r.t. $X_2$ is $\leq R_2$ such that the function
\[
F_N \ = \ P_N(f_1, f_2)
\]
has the property that
\[
\frac{d^n}{d z^n} F_N(w_j) \ = \ 0 \qquad (n = 0, \ldots, N-1; \ j = 1, \ldots, r).
\]

[Note: \ 
Explicated, there are integers
\[
C_{\lambda_1, \lambda_2} : \ 
\begin{cases}
\ 0 \leq \lambda_1 \leq R_1\\
\ 0 \leq \lambda_2 \leq R_2
\end{cases}
\]
with
\[
F_N 
\ = \ 
\sum\limits_{\lambda_1=0}^{R_1} \hsx
\sum\limits_{\lambda_2=0}^{R_2} \hsx
C_{\lambda_1, \lambda_2} \hsx
f_1^{\lambda_1} f_2^{\lambda_2}.
\]
Moreover
\[
0 
\ < \ 
\max\limits_{\lambda_1, \lambda_2} \hsx \abs{C_{\lambda_1, \lambda_2}} 
\ < \ 
e^{3 \hsx [\K:\Q] \hsx r N}.]
\]
\end{x}
\vspace{0.3cm}

Bearing in mind that, by assumption, $f_1(z)$, $f_2(z)$ are algebraically independent over $\C$, let \mM be the smallest positive integer with the property that for some $j_0 $ : $1 \leq j_0 \leq r$, 
\[
\gamma_{\raisebox{-.1cm}{$\scriptstyle N$}} \ 
\ \equiv \ 
\frac{d^M}{d z^M} \hsx F_N(w_{j_0}) 
\ \neq \ 
0.
\]
\vspace{0.3cm}

\begin{x}{\small\bf \un{N.B.}} \ 
$\gamma_N \in \K$ is an algebraic number.  In addition
\[
\frac{d^m}{d z^m} F_N(w_j) = 0 : \ 
\begin{cases}
\ 1 \leq j \leq r\\
\ 0 \leq m \leq M - 1
\end{cases}
,
\]
hence $N \leq M$.
\end{x}
\vspace{0.3cm}

\begin{x}{\small\bf NOTATION} \ 
Put
\[
R \ = \  M\strutx^{\frac{1}{\rho_1 + \rho_2}}.
\]
\end{x}
\vspace{0.3cm}

Ultimately, all relevant data depends on $N \gg 0$.  
This said, choose $N \gg 0$ so as to force $M \gg 0$:
\[
\abs{w_j} \ < \ \frac{R}{2} \qquad (j = 1, \ldots, r).
\]
\vspace{0.3cm}

\begin{x}{\small\bf LEMMA} \ 
If $\abs{z} = R$, then $\forall \ j = 1, \ldots, r$,
\[
\frac{1}{\abs{z - w_j}} \ \leq \ \frac{2}{R}.
\]

PROOF \ 
\[
\abs{z - w_j} 
\ \geq \ 
\abs{\abs{z} - \abs{w_j}}
\]
\qquad\qquad $\implies$
\allowdisplaybreaks
\begin{align*}
\frac{1}{\abs{z - w_j}} \ 
&\leq \ 
\frac{1}{\abs{\abs{z} - \abs{w_j}}}
\\[12pt]
&=\  
\frac{1}{\abs{R - \abs{w_j}}}.
\end{align*}
But
\allowdisplaybreaks
\begin{align*}
\abs{w_j} \ < \ \frac{R}{2} 
&
\implies - \abs{w_j} \ > \ -\frac{R}{2}
\\[12pt]
&
\implies R - \abs{w_j} > R - \frac{R}{2}  \ = \ \frac{R}{2} 
\\[12pt]
&
\implies  \frac{1}{\abs{R - \abs{w_j}}} \ < \  \frac{2}{R}.
\end{align*}

The function
\[
G_N(z) 
\ = \ 
F_N(z) \hsx 
\prod\limits_{j=1}^r \hsx
(z - w_j)^{-M}
\]
is entire and
\[
\gamma_N
\ = \ 
M! \hsx G_N(w_{j_0}) \hsx
\prod\limits_{j\neq j_0} \hsx
(w_{j_0} - w_j)^M.
\]
To estimate $\abs{\gamma_N}$, write
\[
\abs{\gamma_N}
\ \leq \ 
M! \hsx 
\prod\limits_{j \neq j_0} \hsx
\abs{w_{j_0} - w_j}^M 
\hsx\cdot\hsx
\sup\limits_{\abs{z} = R} \  
\prod\limits_{j=1}^r \hsx
\abs{z - w_j}^{-M} 
\hsx\cdot\hsx
\abs{F_N}_R
\]

\qquad \textbullet \quad
$M! \leq M^M$
\vspace{0.2cm}

\qquad \textbullet \quad
$\ds\prod\limits_{j\neq j_0} \hsx \abs{w_{j_0} - w_j}^M \ \equiv \  C^M \qquad (C \in \R_{> 0})$
\vspace{0.2cm}

\qquad \textbullet \quad
$\ds\frac{1}{\abs{z - w_j}^M} \ \leq \ \bigg(\ds\frac{2}{R}\bigg)^{M}$
\vspace{0.75cm}

\qquad\qquad $\implies$

\[
\sup\limits_{\abs{z} = R} \ 
\prod\limits_{j=1}^r \hsx
\abs{z - w_j}^{-M} 
\ \leq \ 
\bigg(\frac{2}{R}\bigg)^{r M}
\]
\vspace{0.2cm}


\allowdisplaybreaks
\begin{align*}
\text{\textbullet}
\hspace{1cm}\abs{F_N}_R \ 
&=\ 
\abs{
\sum\limits_{\lambda_1=0}^{R_1} \hsx
\sum\limits_{\lambda_2=0}^{R_2} \hsx
C_{\lambda_1, \lambda_2} \hsx
f_1^{\lambda_1} f_2^{\lambda_2}
}_R
\\[12pt]
&\leq\ 
\sum\limits_{\lambda_1=0}^{R_1} \hsx
\sum\limits_{\lambda_2=0}^{R_2} \hsx
\abs{C_{\lambda_1, \lambda_2}} \hsx 
\abs{f_1^{\lambda_1} f_2^{\lambda_2}}_R
\\[12pt]
&\leq\ 
(R_1 + 1)(R_2 + 1) \hsx
\max\limits_{\lambda_1, \lambda_2} 
\abs{C_{\lambda_1, \lambda_2}} \hsx 
\abs{f_1^{\lambda_1} f_2^{\lambda_2}}_R
\\[12pt]
&\leq\ 
(R_1 + 1)(R_2 + 1) \hsx
e^{3 \hsx [\K:\Q]  \hsx r N} \hsx
\big(
\abs{f_1}_R + 1\big)^{R_1} \hsx 
\big(
\abs{f_2}_R + 1\big)^{R_2}
\end{align*}

\[
\text{\textbullet} \hspace{1cm} \abs{z} \leq R \quad \   \implies \quad 
\begin{cases}
\ \abs{f_1(z)} \leq K_1 R\strutx^{\rho_1}\\[8pt]
\ \abs{f_2(z)} \leq K_2 R\strutx^{\rho_2}
\end{cases}
\hspace{1cm} (\exists \ K_1, K_2 \in \R_{> 0})
\]
\vspace{0.2cm}
$\implies$ 
\allowdisplaybreaks
\begin{align*}
\big(\abs{f_1}_R + 1\big)^{R_1} \hsx \big(\abs{f_2}_R + 1\big)^{R_2} \ 
&\leq \ 
\big(K_1 R\strutx^{\rho_1} + 1\big)^{R_1} \hsx \big(K_2 R\strutx^{\rho_2} + 1\big)^{R_2} 
\\[12pt]
&\leq \ 
K\big(R_1 R\strutx^{\rho_1} + R_2 R\strutx^{\rho_2}\big).
\end{align*}
The next step is to use these majorants and derive an estimate for $\elln(\abs{\gamma_N})$.
\vspace{0.2cm}

\quad {\small\bf FACT}
For $N \gg 0$,
\[
\elln(\abs{\gamma_N}) 
\ \leq \ 
\bigg(1 - \frac{r}{\rho_1 + \rho_2}\bigg) \hsx M \elln(M) 
+ 
M(\elln(M))^{3/4}.
\]
\end{x}
\vspace{0.3cm}

\begin{x}{\small\bf LEMMA} \ 
Let $x \in \K$ be a nonzero algebraic number $-$then
\[
\elln(\abs{x}) + [\K:\Q] \hsx \elln(d_x) + ([\K:\Q] - 1) \hsx \elln(H(x)) \ \geq \ 0.
\]

[Here $d_x$ is the denominator of $x$ and $H(x)$ is the house of $x$ (cf. \S14).]
\vspace{0.3cm}

Take $x = \gamma_N$ in \#11.
\vspace{0.5cm}

\quad {\small\bf FACT} \ 
$\elln(d_{\gamma_N}) \leq M(\elln(M))^{1/2}$.
\vspace{0.5cm}

\quad {\small\bf FACT} \ 
$\elln(H(\gamma_N)) \leq M \elln(M) + M(\elln(M))^{1/2}$.
\vspace{0.2cm}

Therefore
\[
\elln(\abs{\gamma_N}) 
+ 
[\K:\Q]  M(\elln(M))^{1/2} 
+\big([\K:\Q] - 1\big) \big(M \elln(M) + M(\elln(M))^{1/2}\big)
\ \geq \ 
0
\]
or still, 
\allowdisplaybreaks
\begin{align*}
&\bigg(1 - \frac{r}{\rho_1 + \rho_2}\bigg) M \elln(M) 
+M(\elln(M))^{3/4}
+ 
[\K:\Q]  M(\elln(M))^{1/2} 
\\[12pt]
&
\hspace{3cm}
+\big([\K:\Q] - 1\big) \big(M \elln(M) + M(\elln(M))^{1/2}\big)
\\[12pt]
&
\hspace{1.5cm}
\geq \ 0
\end{align*}
or still,
\allowdisplaybreaks
\begin{align*} 
&\bigg([\K:\Q] - \frac{r}{\rho_1 + \rho_2} \bigg) M \elln(M) + M (\elln(M))^{3/4} 
+ [\K:\Q] M \big(\elln(M)\big)^{1/2} 
\\[12pt]
&
\hspace{3cm}
+ [\K:\Q] M \big(\elln(M)^{1/2}\big) -M \big(\elln(M)^{1/2}\big)
\\[12pt]
&
\hspace{1.5cm}
\geq \ 0
\end{align*}
or still, 
\[
\bigg([\K:\Q] - \frac{r}{\rho_1 + \rho_2} \bigg) M \elln(M) 
+ M (\elln(M))^{3/4}
+ (2[\K:\Q] - 1) M (\elln(M))^{1/2}
\ \geq \ 
0
\]
or still, 
\[
\bigg([\K:\Q] - \frac{r}{\rho_1 + \rho_2} \bigg) M \elln(M)
\ \geq \ 
-M (\elln(M))^{3/4} - (2[\K:\Q] - 1) M (\elln(M))^{1/2}
\]
or still, 
\[
\bigg(\frac{r}{\rho_1 + \rho_2} - [\K:\Q]\bigg) M \elln(M)
\ \leq \ 
M (\elln(M))^{3/4} + (2[\K:\Q] - 1) M (\elln(M))^{1/2}
\]
or still, 
\[
\bigg(\frac{r}{\rho_1 + \rho_2} - [\K:\Q]\bigg) \elln(M)
\ \leq \ 
(\elln(M))^{3/4} + (2[\K:\Q] - 1) (\elln(M))^{1/2}
\]
or still, 
\[
\bigg(\frac{r}{\rho_1 + \rho_2} - [\K:\Q]\bigg)
\ \leq \ 
(\elln(M))^{-1/4} + (2[\K:\Q] - 1) (\elln(M))^{-1/2}.
\]
But $N \ra \infty$ $\implies$ $M \ra \infty$, hence
\[
\frac{r}{\rho_1 + \rho_2} - [\K:\Q]
\ \leq \ 
0
\]
\qquad\qquad $\implies$
\[
\frac{r}{\rho_1 + \rho_2} 
\ \leq \ 
[\K:\Q]
\]
\qquad\qquad $\implies$
\[
r
\ \leq \ 
(\rho_1 + \rho_2) [\K:\Q],
\]
from which the claimed bound on \mS (cf. \#4).
\end{x}
\vspace{0.3cm}

\begin{x}{\small\bf EXAMPLE} \ 
Take $\K = \Q$, $f_1(z) = z$, $f_2(z) = e^z$ $-$then
\[
S \ = \ \{w \in \C : w, e^w \in \Q\}.
\]
But
\[
w \in \Q \quad (w \neq 0) 
\implies 
e^w \in \PP \qquad \text{(cf. \S 9, \#1)},
\]
so $S = \{0\}$, a set of cardinality 1.  
On the other hand, 
\[
\rho_1 = 0, \ \rho_2 = 1 
\implies
\rho_1 + \rho_2 = 1,
\]
thus in this case, the estimate
\[
(\rho_1 + \rho_2) [\K:\Q]
\]
is the best possible.
\end{x}
\vspace{0.3cm}

\[
\text{APPENDIX}
\]
\vspace{0.3cm}

We shall indicate the derivation of the estimate
\[
\elln(\abs{\gamma_N}) 
\ \leq \ 
\bigg(1- \frac{r}{\rho_1 + \rho_2} \bigg) M \elln(M) + M (\elln(M))^{3/4}.
\]
First of all, the term
\[
M (\elln(M))^{3/4}
\]
results from the discussion of $\abs{F_N}_R$, hence can be set aside.  
As for
\[
\bigg(1- \frac{r}{\rho_1 + \rho_2} \bigg) M \elln(M),
\]
note that
\vspace{0.2cm}
\allowdisplaybreaks
\begin{align*}
&\text{\textbullet} \quad
\elln(M!) \ \leq \ M \elln(M) 
\\[12pt]
&\text{\textbullet} \quad
\elln(C^M) \ \leq \ M \elln(C)
\\[12pt]
&\text{\textbullet} \quad
\elln\bigg(\frac{2}{R}\bigg)^{r M}\ 
=\ 
\elln\big(2^{r M}\big)- \elln\bigg(M\strutx^{\frac{r M}{\rho_1 + \rho_2}}\bigg)
\\[12pt]
&
\hspace{2.3cm}
=
M r \elln(2) - \frac{r}{\rho_1 + \rho_2}M \elln(M).
\end{align*}
One must then add these terms.  
But since $N \gg 0$ $\implies$ $M \gg 0$, one can ignore 
\[
M \elln(C) \quad \text{and} \quad M r \elln(2), 
\]
leaving 
\[
M \hsx \elln(M) - \frac{r}{\rho_1 + \rho_2} M \elln(M) 
\ = \ 
\bigg(1- \frac{r}{\rho_1 + \rho_2} \bigg) M \elln(M) .
\]


%% file: _30_schneider_lang_criteria.tex
\chapter{
$\boldsymbol{\S}$\textbf{30}.\quad  SCHNEIDER-LANG CRITERIA}
\setlength\parindent{2em}
\setcounter{theoremn}{0}
\renewcommand{\thepage}{\S30-\arabic{page}}

\ \indent 
There are extensions and variants of the Schneider-Lang criterion (cf. \S29, \#1), e.g., work with meromorphic functions (i.e., quotients of two entire functions) or raise the variables from 1 to $n$ (i.e., replace $\C$ by $\C^n$).
\vspace{0.3cm}

Fix an algebaic number field $\K$.
\vspace{0.3cm}

\begin{x}{\small\bf RAPPEL} \ 
A meromorphic function is said to be of finite strict order $\leq \rho$ if it is the quotient of two entire functions each of finte strict order $\leq \rho$.
\end{x}
\vspace{0.3cm}

\begin{spacing}{1.50}
\begin{x}{\small\bf THEOREM} \ 
Let $f_1, f_2, \ldots, f_n$ $(n \geq 2)$ be meromorphic functions such that $f_1$, $f_2$ are of finite strict orders $\leq \rho_1$, $\leq \rho_2$.  
Assume: \ $f_1$, $f_2$ are algebraically independent over $\C$ and that the derivative $\ds\frac{d}{dz}$ maps the ring $\K[f_1, f_2, \ldots, f_n]$ into itself $-$then the set \mS of $w \in \C$ which are not among the singularities of $f_1, f_2, \ldots, f_n$ but such that 
\[
f_i(w) \in \K \qquad (1 \leq i \leq n)
\]
is finite and in fact the cardinality of \mS is bounded by
\[
(\rho_1 + \rho_2) \hsx [\K:\Q].
\]
\end{x}
\end{spacing}
[The argument is a straight forward extension of that used to establish the Schneider-Lang criterion.  
Thus let $w_1, \ldots, w_r$ be elements of \mS which are not among the singularities of $f_1, f_2, \ldots, f_n$ but such that
\[
f_i(w_j) \in \K \qquad (1 \leq i \leq n; \ 1 \leq j \leq r).
\]
Choose entire functions $g_1$, $g_2$ of finite strict orders $\leq \rho_1$, $\leq \rho_2$, with the property
that $g_1 f_1$, $g_2 f_2$ are entire and
\[
\begin{cases}
\ g_1(w_j) \neq 0 \qquad (1 \leq j \leq r)\\
\ g_2(w_j) \neq 0 \qquad (1 \leq j \leq r)
\end{cases}
.
\]
Define $F_N$ as in \S29, \#7 and form
\[
g_1^{R_1} g_2^{R_2} F_N,
\]
an entire function admitting $w_1, \ldots, w_r$ as zeros of order at least equal to \mM.  
Put
\[
G_N (z)
\ = \ 
g_1(z)^{R_1} \hsx g_2(z) ^{R_2} \hsx F_N (z) \hsx
\prod\limits_{j=1}^r \hsx
(z - w_j)^{-M},
\]
take \mR as in \S29, \#9, and note that
\[
\gamma_N 
\ = \ 
M! \hsx G_N(w_{j_0}) \hsx g_1(w_{j_0}) ^{-R_1}\hsx g_2(w_{j_0})^{-R_2} \hsx
\prod\limits_{j\neq j_0} \hsx
(w_{j_0} - w_j)^M.
\]
Proceed from this point as before.]
\vspace{0.5cm}

There are also versions of Schneider-Lang where $\C$ is replaced by $\C^n$.
\vspace{0.3cm}

\begin{spacing}{1.50}
To set matters up, fix an algebraic number field $\K$ and suppose that $f_1, \ldots, f_m$ are entire functions of the complex variables 
$z_1, \ldots, z_n$ with $m \geq n + 1$.  
Assume: \ $f_1, \ldots, f_{n+1}$ are algebraically independent over $\C$ of finite strict orders $\leq \rho_1, \ldots, \leq \rho_{n+1}$ and that the partial deriviatives $\ds\frac{\partial}{\partial z_i}$ $(1 \leq i \leq n)$ map the ring $\K[f_1, \ldots, f_m]$ into itself.  
Denote by \mS the set of $w \in \C^n$ such that
\[
f_k(w) \in \K \qquad (1 \leq k \leq m).
\]
\end{spacing}

\begin{x}{\small\bf REMARK} \ 
It can be shown that \mS is contained in an algebraic hyper-surface of degree at most 
\[
n(\rho_1 + \cdots + \rho_{n+1}) \hsx [\K:\Q].
\]
\vspace{0.2cm}

[Note: \ 
This means that \mS is the set of zeros of a nonzero polynomial in $\C[X_1, \ldots, X_n]$, its degree being the minimum of the degrees of the nonzero polynomials which annihilate \mS.]
\end{x}
\vspace{0.3cm}

\begin{x}{\small\bf THEOREM} \ 
Let $e_1, \ldots, e_n$ be a basis for $\C^n$ over $\C$ and let $S_1, \ldots, S_n$ be subsets of $\C$.  
Suppose further that 
\[
S \supset 
\{s_1 e_1 + \cdots + s_n e_n: (s_1, \ldots, s_n) \in S_1 \times \cdots \times S_n\}.
\]
I.e.: \ $\forall \ (s_1, \ldots, s_n) \in S_1 \times \cdots \times S_n$:
\[
f_k(s_1 e_1 + \cdots + s_n e_n) \in \K \qquad (1 \leq k \leq m).
\]
Then
\[
\min\limits_{1 \leq i \leq n} \card S_i 
\ \leq \ 
n(\rho_1 + \cdots + \rho_{n+1}) \hsx [\K:\Q].
\]
\vspace{0.2cm}

[Note: \ 
Take $n = 1$ to recover the Schneider-Lang criterion.]
\end{x}
\vspace{0.3cm}

\begin{x}{\small\bf \un{N.B.}} \ 
Therefore the set \mS cannot contain a product $S_1 \times \cdots \times S_n$, where each $S_i$ is infinite. 
\end{x}
\vspace{0.3cm}

Let $\Gamma$ be an additive subgroup of $\C^n$ which contains a basis for $\C^n$ over $\C$ $-$then the points of $\Gamma$ are linearly independent over the complex numbers and this allows one to change coordinates so as to render $\Gamma$ a product:
\[
\Gamma \ \approx \  S_1 \times \cdots \times S_n.
\]
Consider the values
\[
f_k(\zeta_1, \ldots, \zeta_n)
\qquad (1 \leq k \leq m),
\]
where
\[
(\zeta_1, \ldots, \zeta_n) \in \Gamma.
\]
Then the set \mS cannot contain $\Gamma$ (cf. \#5).
\vspace{0.3cm}

\begin{x}{\small\bf EXAMPLE} \ 
It is shown in \S31, \#13 that
\[
\int\limits_0^1 \hsx 
\frac{1}{1 + x^3} dx 
\ = \ 
\frac{1}{3} \bigg(\elln(2) + \frac{\pi}{\sqrt{3}}\bigg)
\]
is transcendental.  
Here is another approach.  
Suppose that
\[
\frac{1}{3} \bigg(\elln(2) + \frac{\pi}{\sqrt{3}}\bigg)
\]
is algebraic $-$then
\allowdisplaybreaks
\begin{align*}
\alpha \ 
&\equiv\ 
3 \sqrt{3} \hsx \sqrt{-1} \hsx \cdot \frac{1}{3} \bigg(\elln(2) + \frac{\pi}{\sqrt{3}}\bigg)
\\[12pt]
&=\ 
\sqrt{3} \hsx \sqrt{-1} \hsx \elln(2) + 3 \sqrt{-1} \hsx \pi
\end{align*}
is algebraic.  Work in $\C^2$ with the functions
\[
f_1(z_1, z_2) \ = \ \exp(z_1), \ 
f_2(z_1, z_2) \ = \ \exp(z_2), \ 
f_3(z_1, z_2) \ = \ z_1 + \sqrt{3} \hsx \sqrt{-1} \hsx z_2
\]
and let $\K = \Q(\sqrt{3} \hsx \sqrt{-1}, \alpha)$.  
Denote by $\Gamma$ the additive subgroup of $\C^2$ generated by the points
\[
u \ = \ 
(3 \pi \sqrt{-1}, \elln(2)), \ 
v \ = \ 
(-3\elln(2), 3 \pi \sqrt{-1}) 
\]
\qquad\qquad $\implies$ 
\[
\Gamma \ = \ \Z u + \Z v.
\]
Then these points are linearly independent over $\C$ since their determinant
\[
\begin{pmatrix}
 3 \pi \sqrt{-1} &&\elln(2)\\
 \\
 -3 \elln(2) &&3 \pi \sqrt{-1}
\end{pmatrix}
\ = \ 
-9 \pi^2 + 3 \big(\elln(2)\big)^2 \ \neq \ 0.
\]
The claim now is that $S \subset \Gamma$, a contradiction.  It is trivial that
\[
f_1(\Gamma) \subset \K, \ 
f_2(\Gamma) \subset \K.
\]
As for $f_3$, we have
\allowdisplaybreaks
\begin{align*}
f_3(3 \pi \sqrt{-1}, \elln(2))\
&=\ 
3 \pi \sqrt{-1} + \sqrt{3} \hsx \sqrt{-1} \hsx \elln(2)
\\[12pt]
&=\ 
\sqrt{3} \hsx \sqrt{-1} \hsx \elln(2) + 3 \sqrt{-1} \hsx \pi
\\[12pt]
&=\ 
\alpha
\end{align*}
and
\allowdisplaybreaks
\begin{align*}
f_3(-3 \hsx \elln(2), 3 \pi \sqrt{-1}) \ 
&=\ 
-3 \elln(2) +  \sqrt{3} \hsx \sqrt{-1} \hsx 3 \pi \sqrt{-1} 
\\[12pt]
&=\ 
-3 \elln(2) - 3 \sqrt{3} \hsx \pi.
\end{align*}
By construction, $\sqrt{3}\hsx \sqrt{-1} \in \K$.  With this in mind, consider
\allowdisplaybreaks
\begin{align*}
\sqrt{3} \hsx \sqrt{-1} \hsx (-3 \elln(2) - 3 \sqrt{3} \hsx \pi)\ 
&=\ 
-3(\sqrt{3} \hsx \sqrt{-1} \hsx \elln(2) + 3 \sqrt{-1} \hsx \pi)
\\[12pt]
&=\ 
-3 \alpha
\end{align*}
or still, 
\allowdisplaybreaks
\begin{align*}
-3 \elln(2) - 3 \sqrt{3} \hsx \pi \ 
&=\ 
\frac{-3}{\sqrt{3} \hsx \sqrt{-1}} \alpha
\\
&\in \K.
\end{align*}
\end{x}
\vspace{0.3cm}
\begin{x}{\small\bf NOTATION}\ 
Given
\[
\begin{cases}
\ \bz = (z_1, \ldots, z_n)\\[7pt]
\ \bw = (w_1, \ldots, w_n)
\end{cases}
\]
in $\C^n$, write
\[
\bzw \ = \ z_1 w_1 + \cdots z_n w_n.
\]
\end{x}
\vspace{0.3cm}

Let $d_0$, $d_1$, and $n$ be integers with 
\[
0 \ \leq \ d_0 \ \leq \ n \ < \ d_0 + d_1.
\]
\vspace{0.1cm}

\begin{x}{\small\bf \un{N.B.}}  \ 
The role of $m$ above is played at this juncture by
\[
d \ \equiv \ d_0 + d_1 \ > \ n 
\implies 
n + 1 \ \leq \ d.
\]
\end{x}
\vspace{0.3cm}

Let $\bx_1, \ldots, \bx_{d_1}$ be $\Q$-linearly independent elements of $\Qbar^{\hspace{.05cm} n}$ 
and let $\by_1, \ldots, \by_n$ be a basis for $\C^n$ over $\C$.  Write
\[
\by_j 
\ = \ 
(y_{1 j}, \ldots y_{n  j}) \qquad (1 \leq j \leq n)
\]
and call $\Gamma$ the additive subgroup of $\C^n$ generated by the $\by_j$.
\vspace{0.3cm}

\begin{x}{\small\bf THEOREM} \ 
At least one of the following numbers
\[
y_{h j} \quad (1 \leq h \leq d_0), 
\qquad
e^{\bx_i \by_j} \quad (1 \leq i \leq d_1, \ 1 \leq j \leq n)
\]
is transcendental.
\vspace{0.2cm}

PROOF \ 
Consider the functions
\[
f_h(\bz) \ = \ z_h \qquad (1 \leq h \leq d_0), \ 
f_{d_0 + i} (\bz) \ = \ e^{\bx_i \bz} \qquad (1 \leq i \leq d_1).
\]
The condition on the ``finite strict orders'' is certainly satisfied and since 
$\bx_1, \ldots, \bx_{d_1}$ are linearly independent over $\Q$, the functions $f_1, \ldots, f_d$ are algebraically independent over the field 
$\Q(z_1, \ldots, z_n)$.  
Moreover
\[
\frac{\partial}{\partial z_j} f_h
\ = \ 
\delta_{h j} 
\ = \ 
\begin{cases}
\ 0 \quad \text{if} \quad h \neq j\\
\ 1 \quad \text{if} \quad h = j
\end{cases}
\qquad (1 \leq h \leq d_0)
\]
and
\[
\frac{\partial}{\partial z_j} f_{d_0 + i}
\ = \ x_{j i} \hsx f_{d_0 + i}
\qquad (1 \leq i \leq d_1),
\]
where $\bx_i = (x_{1 i}, \ldots, x_{n i})$ $(1 \leq i \leq d_1)$.  
Therefore the partial derivative requirement is satisfied.  
Now let $\K$ be the field generated over $\Q$ by the $(d_0 + 2 d_1) n$ numbers
\[
x_{j i}, \ 
f_h(\by_j) \ = \ y_{h j}, \ 
f_{d_0 + i} (\by_j) \ = \ e^{\bx_i \by_j},
\]
the range of the parameters being
\[
1 \leq h \leq d_0, \ 1 \leq i \leq d_1, \ 1 \leq j \leq n.
\]
To arrive at a contradiction, assume that these numbers are algebraic, hence that $\K$ is an algebraic number field.  
Take a typical point
\[
Y \ \equiv \ s_1 \by_1 + \cdots + s_n \by_n \qquad (\bss = (s_1, \ldots, s_n) \in \Z^n)
\]
on $\Gamma$ $-$then
\[
f_1(Y) \in \K, \ 
\ldots, \ 
f_d(Y) \in \K.
\]
I.e.: \ $\Gamma \subset S$, an impossibility (cf. supra).  
Accordingly the supposition that $\K$
is an algebraic number field is false.  
Since the $x_{j i}$ are algebraic (by hypothesis), it follows that at least one of the following numbers
\[
y_{h j} \quad (1 \leq h \leq d_0), 
\qquad
e^{\bx_i \by_j} \quad (1 \leq i \leq d_1, \ 1 \leq j \leq n)
\]
is transcendental.
\end{x}
\vspace{0.3cm}

\begin{x}{\small\bf APPLICATION} \ 
Take $d_0 = 0$, so $d = d_1 > n$ (formally, this just means to ignore in the above anything involving $d_0$), hence $y_{h j}$ is no longer part of the theory and the conclusion is that at least one of the
\[
e^{\bx_i \by_j} \qquad (1 \leq i \leq d, \ 1 \leq j \leq n)
\]
is transcendental, hence at least one of the
\[
\bx_i \by_j \qquad (1 \leq i \leq d, \ 1 \leq j \leq n)
\] 
does not belong to $\fL$ (cf. \S31, \#1).

\vspace{0.2cm}

[Note: \ 
It suffices for the analysis that the set $\{\by_1, \ldots, \by_d\}$ contain a basis for $\C^n$ over $\C$.]
\end{x}
\vspace{0.3cm}

\begin{x}{\small\bf EXAMPLE} \ 
Let $\lambda_1, \lambda_2$, $\lambda_3$ be elements of $\fL$ and assume that
\[
\lambda_1 
+ 
\sqrt[\leftroot{0}\uproot{3}3]{2} \hsx \lambda_2 
+
\sqrt[\leftroot{0}\uproot{3}3]{4} \hsx \lambda_3 
\ = \ 0.
\]
Then
\[
1, \ \sqrt[\leftroot{0}\uproot{3}3]{2}, \ \sqrt[\leftroot{0}\uproot{3}3]{4}
\]
belong to $\Qbar$ and we claim that
\[
\lambda_1 \ = \ 0, \ 
\lambda_2 \ = \ 0, \ 
\lambda_3 \ = \ 0.
\]
To see this, start by multiplying the given relation by $\sqrt[\leftroot{-1}\uproot{3}3]{2}$ and $\sqrt[\leftroot{-1}\uproot{3}3]{4}$:
\[
2 \lambda_3  
+ 
\sqrt[\leftroot{-1}\uproot{3}3]{2} \hsx \lambda_1 
+ 
\sqrt[\leftroot{-1}\uproot{3}3]{4}\hsx \lambda_2
\ = \ 
0
\quad \text{and} \quad 
2 \lambda_2  
+ 
2  \sqrt[\leftroot{-1}\uproot{3}3]{2} \hsx \lambda_3 
+ 
\sqrt[\leftroot{-1}\uproot{3}3]{4}\hsx \lambda_1
\ = \ 
0.
\]
Put
\[
\begin{cases}
\ \bx_1 = (1, 0), \ \bx_2 = (0, 1), \ \bx_3 = (\sqrt[\leftroot{-1}\uproot{3}3]{2},\sqrt[\leftroot{-1}\uproot{3}3]{4})\\[7pt]
\ \by_1 = (\lambda_2,\lambda_3), \ \by_2 = (\lambda_1,\lambda_2), \ \by_3 = (2 \lambda_3,\lambda_1)
\end{cases}
.
\]
\\
Here $d = 3$, $n = 2$ and
\[
\begin{matrix}
{\bx_1 \by_1 = \lambda_2,\ \hspace{.2cm} }  
&&{\bx_1 \by_2 = \lambda_1, \hspace{.45cm} }
&&{\bx_1 \by_3 = 2 \lambda_3, \hspace{.4cm} }
\\
{\bx_2 \by_1 = \lambda_3, \hspace{.3cm}}
&&{\bx_2 \by_2 = \lambda_2, \hspace{.5cm} }
&&{\bx_2 \by_3 =  \lambda_1, \hspace{.6cm} }
\\
\bx_3 \by_1 = -\lambda_1, 
&&\bx_3 \by_2 = -2 \lambda_3, 
&&{\bx_3 \by_3 = -2 \lambda_2. \hspace{.1cm} }
\end{matrix}
\]
Moreover if $\lambda_1 \neq 0$, $\lambda_2 \neq 0$, $\lambda_3 \neq 0$, then the matrix
\[
\begin{pmatrix}
\lambda_2 &\lambda_1 &2\lambda_3\\[8pt]
\lambda_3 &\lambda_2 &\lambda_1
\end{pmatrix}
\]
has rank 2, thus $\{\by_1, \by_2, \by_3\}$ contains a basis for $\C^2$ over $\C$.  
Therefore this data realizes the setup of \#10, hence at least one of the
\[
\bx_i \by_j  \qquad (1 \leq i \leq 3, \ 1 \leq j \leq 3)
\]
does not belong to $\fL$, an impossibility.  
Since the supposition that $\lambda_1 \neq 0$, $\lambda_2 \neq 0$, $\lambda_3 \neq 0$ has led to a contradiction, at least one of the 
$\lambda_1$, $\lambda_2$, $\lambda_3$ is 0, say $\lambda_1 = 0$, leaving $\lambda_2$ and $\lambda_3$: 
\[
\sqrt[\leftroot{-1}\uproot{3}3]{2} \lambda_2 
+
\sqrt[\leftroot{-1}\uproot{3}3]{4} \lambda_3 
\ = \ 
0.
\]
Obviously
\[
\begin{cases}
\ \lambda_2 = 0 \implies \lambda_3 = 0\\[7pt]
\ \lambda_3 = 0 \implies \lambda_2 = 0
\end{cases}
.
\]
If now both $\lambda_2$ and $\lambda_3$ are nonzero, then on general grounds (cf. \S24, \#10), the ratio 
$\lambda_2 / \lambda_3$ is either rational or transcendental.  But $\lambda_2 / \lambda_3$ is not rational but is algebraic \ldots \ . 
\end{x}
\vspace{0.3cm}

\begin{x}{\small\bf APPLICATION} \ 
Take $d_0 = 1$, $d_1 = n$ $(\implies d = 1 + n)$.  
Work this time with $\bx_1, \ldots, \bx_n$ $\Q$-linearly independent elements of $\Qbar^{\hspace{0.05cm} n}$ and $\by_1, \ldots, \by_n$ as a basis for $\C^n$ over 
$\C$.  Write
\[
\by_j \ = \ (y_{1 j}, \ldots, y_{n j}) \qquad (1 \leq j \leq n)
\]
and assume that the numbers
\[
y_{1 j} \qquad (1 \leq j \leq n) \quad (h = 1) 
\]
are algebraic $-$then the conclusion is that at least one of the
\[
e^{\bx_i \by_j} \qquad (1 \leq i \leq n, \ 1 \leq j \leq n)
\]
is transcendental, hence at least one of the
\[
\bx_i \by_j \qquad (1 \leq i \leq n, \ 1 \leq j \leq n)
\]
does not belong to $\fL$.
\vspace{0.2cm}

[Note: \ 
This is a literal transcription of \#9 to the current setting.  For later use, observe that the symbol $d$ does not appear in any of the formulas.  
Because of this, one can replace $n$ by $d$ throughout, so now at least one of the 
\[
\bx_i \by_j \qquad (1 \leq i \leq d, \ 1 \leq j \leq d)
\]
does not belong to $\fL$.]
\end{x}
\vspace{0.3cm}


%% file: _31_baker_statement.tex
\chapter{
$\boldsymbol{\S}$\textbf{31}.\quad  BAKER: \ STATEMENT}
\setlength\parindent{2em}
\setcounter{theoremn}{0}
\renewcommand{\thepage}{\S31-\arabic{page}}


\begin{x}{\small\bf NOTATION} \ 
Put
\[
\fL \ = \ \{\lambda \in \C : e^\lambda \in \Qbar^\times\}
\]
or still, 
\[
\fL \ = \ \exp^{-1} \big(\Qbar^\times\big).
\]
\end{x}
\vspace{0.3cm}
\index{$\fL$}

\begin{x}{\small\bf LEMMA} \ 
$\fL$ is a $\Q$-vector space
\end{x}
\vspace{0.3cm}

\begin{x}{\small\bf LEMMA} \ 
$\Qbar \cap \fL \ = \ \{0\} \qquad \text{(cf. \S21, \#4)}$.
\end{x}
\vspace{0.3cm}

\begin{x}{\small\bf \un{N.B.}} \ 
Therefore every nonzero element of $\fL$ is transcendental.
\end{x}
\vspace{0.3cm}

\begin{x}{\small\bf THEOREM} \ 
The following assertions are equivalent.
\vspace{0.2cm}

\qquad\qquad \textbullet \quad
If $\alpha$ is a nonzero algebraic number, then $e^\alpha$ is transcendental (Hermite-Lindemann).
\vspace{0.2cm}

\qquad\qquad \textbullet \quad
If $\lambda \in \fL$ is nonzero, then 1, $\lambda$ are $\Qbar$-linearly independent.
\vspace{0.2cm}

\qquad\qquad \textbullet \quad
If $a$ is a nonzero complex number, then at least one of the two numbers $a$, $e^a$ is transcendental.
\end{x}
\vspace{0.3cm}

\begin{x}{\small\bf THEOREM} \ 
The following assertions are equivalent.
\vspace{0.2cm}

\qquad\qquad \textbullet \quad
If $\alpha \neq 0, 1$ is algebraic and if $\beta \notin \Q$ is algebraic, then $\alpha^\beta$ is transcendental (Gelfond-Schneider).
\vspace{0.2cm}

\qquad\qquad \textbullet \quad
If $\lambda_1 \in \fL$, $\lambda_2 \in \fL$ are nonzero $\Q$-linearly independent, then $\lambda_1$, $\lambda_2$ are $\Qbar$-linearly independent.
\vspace{0.2cm}

\qquad\qquad \textbullet \quad
If $a$, $b$ are nonzero complex numbers with $a \notin \Q$, then at least one of the three numbers $a$, $e^b$, $e^{ab}$ is transcendental.
\end{x}
\vspace{0.3cm}
\begin{x}{\small\bf REMARK}\ 
$\fL$ is not a $\Qbar$-vector space.
\end{x}
\vspace{0.3cm}

Items 5 and 6 serve to motivate the central result which is due to Baker.
\vspace{0.3cm}

\begin{x}{\small\bf THEOREM} \ 
If $\lambda_1 \in \fL, \ldots, \lambda_n \in \fL$ are nonzero and $\Q$-linearly independent, then 
$1, \lambda_1, \ldots, \lambda_n$ are $\Qbar$-linearly independent.
\end{x}
\vspace{0.3cm}

\begin{x}{\small\bf \un{N.B.}} \ 
This is the so-called ``inhomogeneous case''.
\index{inhomogeneous case}  
Dropping the ``1'' gives the ``homogeneous case''.
\index{homogeneous case}    
I.e.: \ If $\lambda_1, \in \fL, \ldots, \lambda_n \in \fL$ are nonzero and $\Q$-linearly independent, then 
$\lambda_1, \ldots, \lambda_n$ are $\Qbar$-linearly independent.
\end{x}
\vspace{0.3cm}

We shall postpone the proof of \#8 until \S33 and simply assume its validity for the remainder of this \S.
\vspace{0.3cm}

\begin{x}{\small\bf SCHOLIUM} \ 
If $\lambda_1 \in \fL, \ldots, \lambda_n \in \fL$ are nonzero and $\Q$-linearly independent, then 
\[
\beta_0 + \beta_1 \lambda_1 + \cdots + \beta_n \lambda_n \  \neq \  0
\]
for every tuple $(\beta_0, \beta_1, \ldots, \beta_n)$ of algebraic numbers different from $(0,0, \ldots, 0)$.
\end{x}
\vspace{0.3cm}

\begin{x}{\small\bf LEMMA} \ 
Every nonzero linear combination
\[
\beta_1 \lambda_1 + \cdots + \beta_n \lambda_n \qquad (\lambda_1 \in \fL, \ldots, \lambda_n \in \fL)
\]
with algebraic coefficients is transcendental.
\vspace{0.2cm}

PROOF \ 
Argue by induction on $n$, starting with $n = 1$, the validity in this case being ensured by \#4.  
Proceeding, suppose first that $\lambda_1, \ldots, \lambda_n$ are nonzero and $\Q$-linearly independent and suppose that 
\[
\beta_1 \lambda_1 + \cdots + \beta_n \lambda_n \ \equiv \ -\beta_0
\]
is algebraic, hence
\[
\beta_0 + \beta_1 \lambda_1 + \cdots + \beta_n \lambda_n \ = \ 0
\]
\qquad\qquad $\implies$
\[
\beta_1 = 0, \ldots, \beta_n = 0,
\]
contradicting the assumption that
\[
\beta_1 \lambda_1 + \cdots + \beta_n \lambda_n \ \neq \ 0.
\]
If now instead there exist rationals $q_1, \ldots, q_n$ such that
\[
q_1 \lambda_1 + \cdots + q_n \lambda_n \ =\  0
\]
with $q_n \neq 0$, then
\allowdisplaybreaks
\begin{align*}
&
q_n (\beta_1 \lambda_1 + \cdots + \beta_n \lambda_n)
\\[12pt]
&\hspace{1.5cm}=\ 
q_n \beta_1 \lambda_1 + \cdots + q_n \beta_n \lambda_n
\\[12pt]
&\hspace{1.5cm}=\ 
q_n \beta_1 \lambda_1 + \cdots + q_n \beta_n \lambda_n 
-\beta_n(q_1 \lambda_1 + \cdots + q_n \lambda_n)
\\[12pt]
&\hspace{1.5cm}=\ 
(q_n \beta_1 - q_1 \beta_n) \lambda_1  + \cdots + (q_n \beta_n - q_n \beta_n) \lambda_n
\\[12pt]
&\hspace{1.5cm}=\ 
(q_n \beta_1 - q_1 \beta_n) \lambda_1 + \cdots + (q_n \beta_{n-1} - q_{n-1} \beta_n) \lambda_{n-1},
\end{align*}
a number which, by the induction hypothesis, is transcendental.
\end{x}
\vspace{0.3cm}

\begin{x}{\small\bf APPLICATION} \ 
If $\alpha$, $\beta$, are nonzero algebraic numbers, then
\[
\beta \pi + \Log \alpha
\]
is transcendental.

\vspace{0.2cm}

[In \#11, take
\[
\begin{cases}
\ \lambda_1 = 2 \pi \sqrt{-1} \quad \big(e^{\lambda_1} = 1\big), \  \lambda_2 = \Log \alpha\\[8pt]
\ \beta_1 = \sqrt{-1} \hsx \beta, \ \beta_2 = -2
\end{cases}
.
\]
Then
\[
\sqrt{-1} \hsx \beta (2 \pi \sqrt{-1}) + (-2) \Log \alpha
\]
is transcendental, i.e., 
\[
-\beta 2 \pi + (-2) \Log \alpha
\]
is transcendental, i.e., 
\[
-\frac{1}{2} \big(-\beta 2 \pi + (-2) \Log \alpha \big)
\]
is transcendental, i.e., 
\[
\beta  \pi + \Log \alpha
\]
is transcendental.
\vspace{0.2cm}

[Note: \ 
Take $\alpha = 1$, $\beta = 1$ and conclude that $\pi$ is transcendental (cf. \S19, \#1).  
On the other hand, if $\alpha \neq 1$, then $\Log \alpha$ is transcendental (cf. \#4).]
\end{x}
\vspace{0.3cm}

\begin{x}{\small\bf EXAMPLE} \ 
Put
\[
I \ = \ \int\limits_0^1 \hsx \frac{1}{1 + x^3} dx.
\]
Then
\[
I \ = \  \frac{1}{3} \bigg(\elln(2) + \frac{\pi}{\sqrt{3}}\bigg)
\]
is transcendental.
\end{x}
\vspace{0.3cm}

\begin{x}{\small\bf LEMMA} \ 
If $\alpha_1, \ldots, \alpha_n$ and $\beta_0, \beta_1, \ldots, \beta_n$ are nonzero algebraic numbers, then
\[
e^{\beta_0} \alpha_1^{\beta_1} \cdots \alpha_n^{\beta_n} \qquad \text{(principal powers)}
\]
is transcendental.
\vspace{0.2cm}

PROOF \ 
Suppose that
\[
\alpha_{n+1} \ \equiv \ 
e^{\beta_0} \alpha_1^{\beta_1} \cdots \alpha_n^{\beta_n}
\]
were algebraic.  
Take Log's $-$then for some $k \in \Z$, 
\allowdisplaybreaks
\begin{align*}
\Log \alpha_{n+1} \ 
&=\ 
\Log (e^{\beta_0} \alpha_1^{\beta_1} \cdots \alpha_n^{\beta_n})
\\[12pt]
&=\ 
\beta_0 + \beta_1 \Log \alpha_1 + \cdots + \beta_n \Log \alpha_n + 2 \pi \sqrt{-1} \hsx k 
\qquad \text{(cf. \S23, \#5)}.
\end{align*}
But
\allowdisplaybreaks
\begin{align*}
\Log -1 \ 
&=\ 
\elln(\abs{-1}) + \pi \sqrt{-1}
\\[12pt]
&=\ 
\pi \sqrt{-1}.
\end{align*}
Therefore
\[
\Log \alpha_{n+1}
\ = \ 
\beta_0 + \beta_1 \Log \alpha_1 + \cdots + \beta_n \Log \alpha_n + 2k \hsx \Log - 1
\]
or still, 
\[
\beta_1 \Log \alpha_1 + \cdots + \beta_n \Log \alpha_n + 2k \hsx  \Log - 1 - \Log \alpha_{n+1} 
\ = \ 
-\beta_0.
\]
But the RHS is algebraic and nonzero, thus so is the LHS, which contradicts \#11.
\end{x}
\vspace{0.3cm}

\begin{x}{\small\bf EXAMPLE} \ 
$e^{\sqrt{2}} \hsx 2^{\sqrt{3}}$ is transcendental.
\end{x}
\vspace{0.3cm}

\begin{x}{\small\bf EXAMPLE} \ 
Consider
\[
e^{\pi \alpha + \beta} \qquad(\alpha, \ \beta \in \Qbar, \ \alpha \neq 0, \ \beta \neq 0).
\]
Then
\[
e^{\pi \alpha} 
\ = \ 
(-1)^{- \sqrt{-1} \hsx \alpha}
\ = \ 
e^{- \sqrt{-1} \hsx\hsx \alpha \hsx \Log -1}.
\]
In the preceding, take
\[
\alpha_1 = -1, \ 
\beta_0 = \beta, \ 
\beta_1 = - \sqrt{-1} \hsx \alpha.
\]
Then
\[
e^{\beta_0} \alpha_1^{\beta_1} 
\ = \ 
e^\beta (-1)^{- \sqrt{-1} \hsx \alpha}
\ = \ 
e^\beta e^{\pi \alpha} 
\ = \ 
e^{\pi \alpha + \beta} 
\]
is transcendental.
\vspace{0.2cm}

[Note: \ 
Take $\alpha = 2 \sqrt{-1}$ and conclude that $e^\beta$ is transcendental  (cf. \S21, \#4).]
\end{x}
\vspace{0.3cm}

\begin{x}{\small\bf LEMMA} \ 
If $\alpha_1 \neq 0, 1, \ldots, \alpha_n \neq 0, 1$ are algebraic numbers and if $\beta_1, \ldots, \beta_n$ are algebraic numbers  with 
$1, \beta_1, \ldots, \beta_n$ $\Q$-linearly independent, then
\[
\alpha_1^{\beta_1} \cdots \alpha_n^{\beta_n} \qquad \text{(principal powers)}
\]
is transcendental.
\vspace{0.2cm}

PROOF \ 
Suppose that
\[
\alpha_{n+1} 
\ \equiv \ 
\alpha_1^{\beta_1} \cdots \alpha_n^{\beta_n} 
\]
was algebraic.  Write
\allowdisplaybreaks
\begin{align*}
\alpha_1^{\beta_1} \cdots \alpha_n^{\beta_n} \ 
&=\ 
e^{\beta_1 \Log \alpha_1} \cdots e^{\beta_n \Log \alpha_n}
\\[12pt]
&=\ 
e^{\beta_1 \Log \alpha_1 \hsx + \cdots + \hsx \beta_n \Log \alpha_n}
\\[12pt]
&=\ 
e^\Lambda
\end{align*}
if
\[
\Lambda
\ = \ 
\beta_1 \Log \alpha_1 + \cdots + \beta_n \Log \alpha_n.
\]
Then
\[
e^\Lambda
\ = \ 
\alpha_{n+1} 
\implies 
\Lambda \in \fL.
\]
Put
\[
\lambda_1 = \Log \alpha_1, \ldots, 
\hsx
\lambda_n = \Log \alpha_n, 
\hsx
\lambda_{n+1} = \Lambda
\]
to get
\[
\beta_1 \lambda_1 + \cdots + \beta_n \lambda_n  + 1(-\lambda_{n+1}) \ = \ 0.
\]
On the other hand, thanks to the assumption that $1, \beta_1, \ldots, \beta_n$ are $\Q$-linearly independent, the entity
\[
\beta_1 \lambda_1 + \cdots + \beta_n \lambda_n + 1 (-\lambda_{n+1})
\]
is nonzero (cf. \S32, \#3 (ii)).  Contradiction.
\end{x}
\vspace{0.3cm}

\begin{x}{\small\bf REMARK} \ 
Consider Gelfond-Schneider (cf. \#6).  
Here $\alpha^\beta = e^{\beta \hsx \Log \alpha}$ is the principal power.
Pass to its $k^\nth$ associate:
\[
\alpha^\beta \big(e^{2 k \pi \sqrt{-1} \hsx \beta} \big) \qquad (k \in \Z) \quad \text{(cf. \S23, \#15)}
\]
and write
\allowdisplaybreaks
\begin{align*}
e^{2 k \pi \sqrt{-1} \hsx \beta} \ 
&=\ 
e^{\pi (2 k \sqrt{-1}\hsx \beta)}
\\[12pt]
&=\ 
(-1)^{- \sqrt{-1}\hsx (2 k \sqrt{-1}\hsx \beta)} \qquad \text{(cf. \#16)}
\\[12pt]
&=\ 
(-1)^{2 k \beta}.
\end{align*}
Therefore
\[
\alpha^\beta \big(e^{2 k \pi \sqrt{-1} \hsx \beta} \big) 
\ = \ 
\alpha^\beta  (-1)^{2 k \beta}
\]
is transcendental.
\end{x}
\newpage

\[
\text{APPENDIX}
\]
\vspace{0.5cm}

For the record,
\[
e^{\Log z} \ = \ z
\]
but
\[
\Log e^z \ \equiv \ z \qquad (\modx 2 \pi \sqrt{-1}).
\]
\vspace{0.2cm}

{\small\bf EXAMPLE} \ 
Consider $\alpha^\beta$ $-$then $\exists\ k \in \Z$:
\allowdisplaybreaks
\begin{align*}
\Log \alpha^\beta \ 
&= \ 
\Log e^{\beta \hsx \Log \alpha}
\\[12pt]
&= \ 
\beta \hsx \Log \alpha + 2 \pi \sqrt{-1} \hsx k
\end{align*}
and
\allowdisplaybreaks
\begin{align*}
e^{\beta \hsx \Log \alpha \hsx + \hsx 2 \pi \sqrt{-1} \hsx k}\ 
&= \ 
e^{\beta \hsx \Log \alpha } e^{2 \pi \sqrt{-1} \hsx k}
\\[12pt]
&= \ 
\alpha^\beta \cdot 1
\\[12pt]
&= \
\alpha^\beta.
\end{align*}

%% file: _32_equivalences.tex
\chapter{
$\boldsymbol{\S}$\textbf{32}.\quad  EQUIVALENCES}
\setlength\parindent{2em}
\setcounter{theoremn}{0}
\renewcommand{\thepage}{\S32-\arabic{page}}

\ \indent 
In this \S, we shall formulate various statements that are equivalent to inhomogeneous Baker or homogeneous Baker.
\vspace{0.2cm}

\begin{x}{\small\bf THEOREM} \ 
The following assertions are equivalent.
\vspace{0.2cm}

\qquad\qquad (i) \quad
If $\lambda_1 \in \fL, \ldots, \lambda_n \in \fL$ are 
nonzero 
\ifcomments
\textit{\textcolor{red}{this is redundant - note a fact pointed out later on in the text}} 
\fi
and $\Q$-linearly independent, then 
$1, \lambda_1, \ldots, \lambda_n$ are $\Qbar$-linearly independent (inhomogeneous Baker).
\vspace{0.2cm}

\qquad\qquad (ii) \quad
If $\lambda_1 \in \fL, \ldots, \lambda_{n-1} \in \fL$ are 
nonzero 
\ifcomments
\textit{\textcolor{red}{this is redundant}} 
\fi 
and $\Q$-linearly independent and if 
$\beta_0, \beta_1, \ldots, \beta_{n-1}$ are algebraic numbers such that
\[
\beta_0 + \beta_1 \lambda_1 + \cdots+ \beta_{n-1} \lambda_{n-1}
\]
is an element of $\fL$, then $\beta_0 = 0$ and $\beta_1, \ldots, \beta_{n-1}$ are rational.
\vspace{0.2cm}

\qquad\qquad (iii) \quad
If $\lambda_1 \in \fL, \ldots, \lambda_{n-1}\in \fL$ are 
nonzero 
\ifcomments
\textit{\textcolor{red}{this is redundant}} 
\fi 
and $\Q$-linearly independent and if 
$\beta_0, \beta_1, \ldots, \beta_{n-1}$ are algebraic numbers such that
\[
\beta_0 + \beta_1 \lambda_1 + \cdots+ \beta_{n-1} \lambda_{n-1}
\]
is an element of $\fL$, then $\beta_0 = 0$ and $\beta_1, \ldots, \beta_{n-1}$ are $\Q$-linearly dependent.
\vspace{0.5cm}

The proof proceeds according to the scheme: 
\[
(ii) \implies (iii), \ 
(i) \implies (ii), \
(iii) \implies (i).
\]

\qquad\qquad $(ii) \implies (iii)$: \quad
Obvious.
\vspace{0.2cm}

\qquad\qquad $(i) \implies (ii)$: \quad
Fix the data per the assumption:
\[
\beta_0 + \beta_1 \lambda_1 + \cdots+ \beta_{n-1} \lambda_{n-1} \in \fL.
\]
Then there exists $\lambda_n \in \fL$:
\[
\beta_0 + \beta_1 \lambda_1 + \cdots+ \beta_{n-1} \lambda_{n-1} - \lambda_n  \ = \ 0.
\]
Therefore $1, \lambda_1, \ldots, \lambda_n$ are $\Qbar$-linearly dependent.  
But $\lambda_1, \ldots, \lambda_{n-1}$ are $\Q$-linearly independent, so by (i), there are rational numbers 
$q_1, \ldots, q_{n-1}$ not all zero such that
\[
\lambda_n 
\ = \ 
q_1 \lambda_1 + \cdots+ q_{n-1} \lambda_{n-1},
\]
hence
\[
\beta_0 + \beta_1 \lambda_1 + \cdots+ \beta_{n-1} \lambda_{n-1}
\ - \ 
(q_1 \lambda_1 + \cdots+ q_{n-1} \lambda_{n-1}) 
\ = \ 
0
\]
or still, 
\[
\beta_0 + (\beta_1 - q_1)\lambda_1 + \cdots+ (\beta_{n-1} - q_{n-1})\lambda_{n-1} \ = \ 0.
\]
Finally, appealing to (i) once again, it follows that $\beta_0 = 0$ and $\beta_i = q_i$ $(1 \leq i \leq n-1)$, thus 
$\beta_1, \ldots, \beta_{n-1}$ are rational.
\vspace{0.2cm}

\qquad\qquad $(iii) \implies (i)$: \quad
Denote by $\sP(\fL)$ the set of finite nonempty subsets \mS of $\fL$ subject to:
\vspace{0.2cm}

\qquad 1. \quad
The elements of \mS are $\Q$-linearly independent.
\vspace{0.2cm}

\qquad 2. \quad
The elements of $S \cup \{1\}$ are $\Qbar$-linearly dependent.\\[8pt]
Then the claim is that $\sP(\fL) = \emptyset$, which will do it.  
Suppose instead that $\sP(\fL) \neq \emptyset$ $-$then
\[
n \ \equiv \ \inf \{\card S : S \in \sP(\fL)\}
\]
is $\geq 1$.  
Fix an element $S = \{\lambda_1, \ldots, \lambda_n\} \in \sP(\fL)$ at which the inf is attained $-$then the 
$\lambda_i \ (1 \leq i \leq n)$ are $\Q$-linearly independent and by definition of $\sP(\fL)$ there exist algebraic numbers 
$\beta_0, \beta_1, \ldots, \beta_n$ with $\beta_1, \ldots, \beta_n$ not all zero: 
\[
\beta_0 + \beta_1 \lambda_1 + \cdots+ \beta_n \lambda_n \ = \ 0.
\]
Assume now without loss of generality that $\beta_n \neq 0$, so
\[
\frac{\beta_0}{-\beta_n} 
 + 
 \frac{\beta_1}{-\beta_n} \lambda_1
 + \cdots + 
 \frac{\beta_n}{-\beta_n} \lambda_n
\ = \ 
0.
\]
Adjusting the notation, one can suppose from the beginning that $\beta_n = -1$ and work with
\[
\beta_0 + \beta_1 \lambda_1 + \cdots+ (-1) \lambda_n \ = \ 0,
\]
hence
\[
\beta_0 + \beta_1 \lambda_1 + \cdots+ \beta_{n-1} \lambda_{n-1} \ = \ \lambda_n \in \fL.
\]
Therefore $\beta_0 = 0$ and $\beta_1, \ldots, \beta_{n-1}$ are $\Q$-linearly dependent (cf. (iii)), thus there exist rational numbers $q_1, \ldots, q_{n-1}$ not all zero such that
\[
q_1 \beta_1 + \cdots+ q_{n-1} \beta_{n-1} \ = \ 0.
\]
Choose
\[
q_k \in \{q_1, \ldots, q_{n-1}\} : q_k \neq 0,  \quad \beta_k \neq 0
\]
\qquad\qquad $\implies$
\[
\beta_k 
\ = \ 
\sum\limits_{\substack{i = 1 \\i \neq k}}^{n-1} \hsx
\bigg(- \frac{q_i}{q_k}\bigg) \beta_i
\]
implying thereby that not all the $\beta_i$ $(i \neq k)$ are zero.  
Meanwhile, since $\beta_0 = 0$, 
\[
\beta_1 \lambda_1 + \cdots+ \beta_n \lambda_n \ = \ 0 \qquad (\beta_n = -1)
\]
\qquad\qquad $\implies$
\allowdisplaybreaks
\begin{align*}
0\ 
&=\ 
\sum\limits_{\substack{i = 1 \\i \neq k}}^{n} \hsx
\lambda_i \beta_i + \lambda_k \beta_k
\\[12pt]
&=\ 
\sum\limits_{\substack{i = 1 \\i \neq k}}^{n} \hsx
\lambda_i \beta_i 
-
\lambda_k \hsx
\sum\limits_{\substack{i = 1 \\i \neq k}}^{n-1} \hsx
\frac{q_i}{q_k} \beta_i
\\[12pt]
&=\ 
-\lambda_n + 
\sum\limits_{\substack{i = 1 \\i \neq k}}^{n-1} \hsx
\bigg(\lambda_i  - \lambda_k \frac{q_i}{q_k} \bigg)\beta_i.
\end{align*}
Put
\[
\begin{cases}
\ \gamma_i \ = \ \lambda_i - \lambda_k \hsx \ds\frac{q_i}{q_k} \qquad (i < n, \ i \neq k) \\[8pt]
\ \gamma_i \ = \ \lambda_n \hspace{2cm} (i = n)
\end{cases}
.
\]
Then the $\gamma_i \in \fL$ $(i \neq k)$ are $\Q$-linearly independent (see infra) and
\[
\sum\limits_{\substack{i = 1 \\i \neq k}}^{n} \hsx
\gamma_i \beta_i \ = \ 0.
\]
Because the $\beta_i$ $(i \neq k)$ are not all zero, we have reached a contradiction to the minimality of $n$.
\vspace{0.2cm}

[Note: \ 
To check that the $\gamma_i$ $(i \neq k)$ are $\Q$-linearly independent, consider a dependence relation

\begin{align*}
0 \ 
&=\ 
\sum\limits_{\substack{i = 1 \\i \neq k}}^{n} \hsx
C_i \gamma_i 
\qquad (C_i \in \Q).
\\[12pt]
&=\
C_n \lambda_n + 
\sum\limits_{\substack{i = 1 \\i \neq k}}^{n-1} \hsx
C_i\bigg(\lambda_i - \lambda_k \frac{q_i}{q_k} \bigg) 
\\[12pt]
&=\ 
\sum\limits_{\substack{i = 1 \\i \neq k}}^{n} \hsx
C_i \lambda_i 
- 
\sum\limits_{\substack{i = 1 \\i \neq k}}^{n-1} \hsx
\lambda_k \hsx C_i \hsx \frac{q_i}{q_k} 
\\[12pt]
&=\ 
\sum\limits_{\substack{i = 1 \\i \neq k}}^{n} \hsx
C_i \lambda_i - C \lambda_k,
\end{align*}
where
\[
C 
\ = \ 
\sum\limits_{\substack{i = 1 \\i \neq k}}^{n-1} \hsx
C_i \hsx \frac{q_i}{q_k} \  \in \ \Q.
\]
But the $\lambda_i$ $(1 \leq i \leq n)$ are $\Q$-linearly independent (by hypothesis), so $C_i = 0$ $(i \neq k)$ 
(and $C = 0$).] 
\end{x}
\vspace{0.3cm}

\begin{x}{\small\bf \un{N.B.}} \ 
The proof that we shall give of Baker in \S33 does not go through 
items (ii) or (iii).
\end{x}

\begin{x}{\small\bf THEOREM} \ 
The following assertions are equivalent.
\vspace{0.2cm}

\qquad (i) \quad
If $\lambda_1 \in \fL, \ldots, \lambda_n \in \fL$ are nonzero and $\Q$-linearly independent, then 
$\lambda_1, \ldots, \lambda_n$ are $\Qbar$-linearly independent (homogeneous Baker).
\vspace{0.2cm}

\qquad (ii) \quad
If $\lambda_1 \in \fL, \ldots, \lambda_n \in \fL$ are nonzero and if 
$\beta_1, \ldots, \beta_n$ are $\Q$-linearly independent elements of $\Qbar$, then
\[
\beta_1 \lambda_1 + \cdots+ \beta_n \lambda_n \ \neq \ 0.
\]

\qquad (iii) \quad
If $\lambda_1 \in \fL, \ldots, \lambda_n \in \fL$ are nonzero and $\Q$-linearly independent and if 
$\beta_1, \ldots, \beta_n$ are $\Q$-linearly independent elements of $\Qbar$, then
\[
\beta_1 \lambda_1 + \cdots+ \beta_n \lambda_n \ \neq \ 0.
\]

The proof preceeds according to the scheme: 
\[
(i) \implies (iii), \ 
(ii) \implies (i), \
(iii) \implies (ii).
\]

\qquad $(i) \implies (iii)$: \quad
Obvious.
\ifcomments
\\ \textit{\textcolor{red}{this is trivial - why cite (cf. \S31, \#10)? - Also this is consistent with the similar proof supra.}} 
\fi
\vspace{0.2cm}

\qquad $(ii) \implies (i)$: \quad
Assume $\lambda_1 \in \fL, \ldots, \lambda_n \in \fL$ are $\Q$-linearly independent and that 
\[
\beta_1 \lambda_1 + \cdots+ \beta_n \lambda_n \ = \ 0 \qquad (\beta_j \in \Qbar, \ \ 1 \leq j \leq n).
\]
Observe that since (ii) is in force, $\beta_1, \ldots, \beta_n$ are not $\Q$-linearly independent, so
let $\gamma_1, \ldots, \gamma_m$ $(m < n)$ be a basis for the $\Q$-span of $\{\beta_1, \ldots, \beta_n\}$, thus
\[
\beta_i 
\ = \ 
\sum\limits_{j=1}^m \hsx
c_{i j} \gamma_j \qquad (1 \leq i \leq n \quad \text{with} \quad c_{i j} \in \Q).
\]
Then
\allowdisplaybreaks
\begin{align*}
0\ 
&=\ 
\beta_1 \lambda_1 + \cdots + \beta_n \lambda_n
\\[12pt]
&=\ 
\bigg(
\sum\limits_{j=1}^m \hsx 
c_{1 j} \gamma_j\bigg) \lambda_1 
+ \cdots + 
\bigg(
\sum\limits_{j=1}^m \hsx 
c_{n j} \gamma_j\bigg) \lambda_n
\\[12pt]
&=\ 
\sum\limits_{j=1}^m \hsx 
\gamma_j 
\bigg(
\sum\limits_{i=1}^n \hsx 
c_{i j} \lambda_i
\bigg)
\\[12pt]
&=\ 
\sum\limits_{j=1}^m \hsx 
\gamma_j 
\hsx
\lambda_j^\prime
\qquad 
\big(
\text{where} \quad
\lambda_j^\prime
\equiv
\sum\limits_{i=1}^n \hsx 
c_{i j} \lambda_i \ \in \fL \ \text{(cf. \S31, \#2)}
\big).
\end{align*}
In view of (ii) at least one and hence all of the $\lambda_j^\prime$ $(j = 1, \ldots, m)$ must be zero.  
Therefore $\forall  \ j = 1, \ldots, m$
\[
\lambda_{j}^\prime 
\ = \ 
c_{1 j} \lambda_1 + \cdots + c_{n j} \lambda_n
\ = \
0
\]
But $\lambda_1 \in \fL, \ldots, \lambda_n \in \fL$ are $\Q$-linearly independent.  Therefore 
\[
c_{1 j} \ = \ \cdots \ = \ c_{n j} = 0 \qquad j = 1, \ldots, m.
\]
And this implies that $\beta_1 = 0, \ldots, \beta_n = 0$, i.e., that the $\lambda_1, \ldots, \lambda_n$ are $\Qbar$-linearly independent. 
\vspace{0.3cm}

\qquad $(iii) \implies (ii)$: \quad
If
\[
 \beta_1 \lambda_1 + \cdots+ \beta_n \lambda_n \ = \ 0,
\]
where $\beta_1, \ldots, \beta_n$ are $\Q$-linearly independent elements of $\Qbar$, then it will be shown that 
\[
\lambda_1 \ = \ 0, \ 
\ldots, \ 
\lambda_n \ = \ 0,
\]
from which the result.  
Renumbering the data if necessary, assume that $\lambda_1, \ldots, \lambda_m$ $(0 \leq m \leq n)$ is a basis for the 
$\Q$-span of $\{\lambda_1, \ldots, \lambda_n\}$:
\[
\lambda_i 
\ = \ 
\sum\limits_{j=1}^m \hsx
c_{i j} \lambda_j 
\qquad (m+1 \leq i \leq n),
\]
where the $c_{i j} \in \Q$.  
Then
\[
0 \ = \ 
\sum\limits_{j=1}^m \hsx
\gamma_j\lambda_j 
\qquad \big(\gamma_j = \beta_j + \sum\limits_{i = m+1}^n c_{i j} \beta_i\big).
\]
Now apply (iii) (with $n$ replaced by $m$): $\lambda_1, \ldots, \lambda_m$ are $\Q$-linearly independent, hence 
$\gamma_1, \ldots, \gamma_m$ are $\Q$-linearly dependent.  
However $\beta_1, \ldots, \beta_n$ are $\Q$-linearly independent, so the only possibility is $m = 0$, implying that
\[
\lambda_1 = 0, \ 
\ldots, 
\lambda_n = 0.
\]
\vspace{0.2cm}

[Note: \ 
If $C_j \in \Q$ $(1 \leq j \leq m)$, then
\allowdisplaybreaks
\begin{align*}
\sum\limits_{j=1}^m \hsx
C_j \gamma_j
&=\ 
\sum\limits_{j=1}^m \hsx
C_j \big(\beta_j + 
\sum\limits_{i=m+1}^n \hsx
c_{i j} \beta_i\big)
\\[12pt]
&=\ 
\sum\limits_{j=1}^m \hsx
C_j \beta_j 
+
\sum\limits_{i=m+1}^n \hsx
\bigg(
\sum\limits_{j=1}^m \hsx
c_{i j} C_j \bigg) \beta_i.]
\end{align*}
\end{x}
\vspace{0.3cm}

\begin{x}{\small\bf REMARK} \ 
One can add a fourth condition, viz.
\vspace{0.2cm}

\qquad\qquad (iv) \quad
If $\lambda_1, \ldots, \lambda_{n+1}$ are nonzero elements of $\fL$ such that $\lambda_1, \ldots, \lambda_n$ are 
$\Qbar$-linearly independent and if $\beta_1, \ldots, \beta_n$ are elements of $\Qbar$ such that
\[
\beta_1 \lambda_1 + \cdots+ \beta_n \lambda_n \ = \ \lambda_{n+1},
\]
then $\beta_1, \ldots, \beta_n$ are rational.
\vspace{0.2cm}

[Note: \ 
Suppose that homogeneous Baker is in force.  Consider item (ii) of \#1 $-$then the crux is to prove that $\beta_0 = 0$.]
\end{x}
\vspace{0.3cm}

\begin{x}{\small\bf \un{N.B.}} \ 
Consider the arrow of inclusion:
\[
\fL \ra \C.
\]
Then it lifts to an arrow
\[
\fL \otimes_{\Q} \Qbar \hsx  \ra \hsx \C
\]
which remains injective iff item (iv) supra is in force.
\end{x}
\vspace{0.3cm}

\begin{x}{\small\bf LEMMA} \ 
Baker's inhomogeneous theorem is equivalent to the conjunction of \S31, \#11 and \S31, \#16.
\end{x}
\vspace{0.3cm}

\begin{x}{\small\bf LEMMA} \ 
Baker's homogeneous theorem is equivalent to \S31, \#11.
\end{x}
\vspace{0.3cm}

\begin{x}{\small\bf \un{N.B.}} \ 
\[
\S31, \ \#11 \Leftrightarrow \S31, \ \#14.
\]
\end{x}
\vspace{0.3cm}


%% file: _33_baker_proof.tex
\chapter{
$\boldsymbol{\S}$\textbf{33}.\quad  BAKER: \ PROOF}
\setlength\parindent{2em}
\setcounter{theoremn}{0}
\renewcommand{\thepage}{\S33-\arabic{page}}

\ \indent 

Our objective is to establish that if $\lambda_1 \in \fL, \ldots, \lambda_n \in \fL$ are nonzero and $\Q$-linearly independent, then 
$1, \lambda_1, \ldots, \lambda_n$ are $\Qbar$-linearly independent (cf. \S31, \#8).  
I.e.: \  If $\gamma_0, \gamma_1, \ldots, \gamma_n$ are algebraic numbers and if
\[
\gamma_0 + \gamma_1 \lambda_1 + \cdots + \gamma_n \lambda_n \ = \ 0,
\]
then
\[
\gamma_0 \ = \ 0, \ 
\gamma_1 \ = \ 0, \ 
\ldots, 
\gamma_n \ = \ 0.
\]
\vspace{0.3cm}

\begin{x}{\small\bf THEOREM} \ 
Let $\K$ be an algebraic number field of degree $d$ over $\Q$, let $\{\beta_1, \ldots, \beta_d\}$ be a basis of the 
$\Q$-vector space $\K$, and let $\lambda_1, \ldots, \lambda_d$ be elements of $\fL$.  
Assume:
\[
\beta_1 \lambda_1 + \cdots + \beta_d \lambda_d \hsx \in \hsx \Qbar.
\]
Then
\[
\lambda_1 \ = \ 0, \ 
\ldots, 
\lambda_d \ = \ 0.
\]

\end{x}
\vspace{0.3cm}

\begin{x}{\small\bf REMARK}\ 
Granted Baker's theorem (in its inhomogeneous version), it follows that \#11 of \S31 is in force.  
So, if
\[
\beta_1 \lambda_1 + \cdots + \beta_d \lambda_d 
\]
is nonzero, then
\[
\beta_1 \lambda_1 + \cdots + \beta_d \lambda_d 
\]
must be transcendental.  
On the other hand, under the assumption that it is algebraic, it must be zero:
\[
\beta_1 \lambda_1 + \cdots + \beta_d \lambda_d \ = \ 0.
\]
Still, this does not imply that
\[
\lambda_1 \ = \ 0, \ 
\ldots, 
\lambda_d \ = \ 0.
\]
\vspace{0.1cm}

The foregoing result can be used to give a quick proof of Baker's inhomogeneous theorem.  
So suppose that
\[
\gamma_0 + \gamma_1 \lambda_1 + \cdots + \gamma_n \lambda_n \ = \ 0.
\]
Put $\K = \Q(\gamma_1, \ldots, \gamma_n)$, choose a basis $\{\beta_1, \ldots, \beta_d\}$ for the 
$\Q$-vector space $\K$, and write
\[
\gamma_j 
\ = \ 
\sum\limits_{i=1}^d \hsx
c_{j i} \beta_i \qquad (1 \leq j \leq n)
\]
with $c_{j i} \in \Q$ $-$then
\allowdisplaybreaks
\begin{align*}
-\gamma_0 \ (\in \Qbar)\ 
&=\ 
\sum\limits_{j=1}^n \hsx
\gamma_j \lambda_j
\\[12pt]
&=\ 
\sum\limits_{j=1}^n \hsx
\bigg(
\sum\limits_{i=1}^d \hsx
c_{j i} \beta_i
\bigg) \lambda_j
\\[12pt]
&=\ 
\sum\limits_{i=1}^d \hsx
\beta_i 
\sum\limits_{j=1}^n \hsx
c_{j i} \lambda_j 
\\[12pt]
&=\ 
\sum\limits_{i=1}^d \hsx
\beta_i \lambda_i^\prime,
\end{align*}
where
\[
\lambda_i^\prime
\ = \ 
\sum\limits_{j=1}^n \hsx
c_{j i} \lambda_j 
\in \fL.
\]
Owing to \#1, 
\[
\lambda_1^\prime \ = \ 0, \ 
\ldots, 
\lambda_d^\prime \ = \ 0.
\]
But $\lambda_1, \ldots, \lambda_n$ are nonzero and $\Q$-linearly independent, thus the relations
\[
\sum\limits_{j=1}^n \hsx
c_{j i} \lambda_j 
\ = \ 0
\]
imply that
\[
c_{j i} \ = \ 0 \qquad (1 \leq i \leq d, 1 \leq j \leq n),
\]
hence
\allowdisplaybreaks
\begin{align*}
\gamma_1 \ = \ 0, \ &\ldots, \gamma_n \ = \ 0
\\[12pt]
&\implies \gamma_0 \ = \ 0.
\end{align*}
\end{x}
\vspace{0.3cm}

\begin{x}{\small\bf RAPPEL}\ 
Let $\K$ be an algebraic number field $-$then the 
\un{trace}
\index{trace (of an algebraic number field)}
$\K \ra \Q$ is the $\Q$-linear map
\[
\gamma \ra \sum\limits_\sigma \hsx \gamma^\sigma,
\]
where $\sigma$ runs over the set of  complex embeddings of $\K$ (a set of cardinality $[\K:\Q]$) and $\gamma^\sigma$ is the image of 
$\gamma$ under $\sigma$. 
\end{x}
\vspace{0.3cm}

\begin{x}{\small\bf NOTATION} \ 
Let $\K$ be an algebraic number field, let $\{\beta_1, \ldots, \beta_d\}$ be a basis for the $\Q$-vector space $\K$, and let 
$\sigma_1: \K \ra \C, \ldots, \sigma_d: \K \ra \C$ be the complex embeddings of $\K$ (label matters so that $\sigma_1$ is the arrow 
$\K \ra \C$ of inclusion).
\end{x}
\vspace{0.3cm}

\begin{x}{\small\bf LEMMA} \ 
\[
\det(\tr(\beta_i \beta_j))_{1 \leq i, j \leq d} 
\ = \ 
(\det B)^2,
\]
where
\[
B \ = \ \big(\beta_k^{\sigma_i}\big)_{1 \leq i, k \leq d}
\]
is nonsingular.
\end{x}
\vspace{0.3cm}

We shall now take up the proof of \#1.
\vspace{0.3cm}

\begin{x}{\small\bf NOTATION} \ 
Put
\[
\Lambda_i \ = \ 
\sum\limits_{k=1}^d \hsx
\beta_k^{\sigma_i} \hsx \lambda_k 
\qquad (1 \leq i \leq d).
\]
\end{x}
\vspace{0.3cm}

\qquad \un{Case 1:} \quad
At least one but not all of the $\Lambda_i$ vanish.
\vspace{0.2cm}

[Arrange the notation so that
\[
\Lambda_1 \neq 0, \ldots, \Lambda_n \neq 0, \Lambda_{n+1} = 0, \ldots, \Lambda_d = 0.
\]
\vspace{0.2cm}

\qquad \textbullet \quad
Define $\bx_i \in \Qbar^{\hsx n}$ by
\[
\bx_i \ = \ \big(\beta_i^{\sigma_1}, \ldots, \beta_i^{\sigma_n}\big) \qquad (1 \leq i \leq d).
\]
If $q_1, \ldots, q_d$ are rational numbers such that
\[
q_1 \bx_1 + \cdots + q_d \bx_d \ = \ (0, \ldots, 0),
\]
I.e., 
\[
q_1 \ 
\begin{pmatrix}
\hsx \beta_1^{\sigma_1} \hsx\\
\vdots\\
\hsx \beta_1^{\sigma_b} \hsx
\end{pmatrix} 
\ + \cdots + \ 
q_d \ 
\begin{pmatrix}
\hsx \beta_d^{\sigma_1}\hsx \\
\vdots\\
\hsx \beta_d^{\sigma_b} \hsx
\end{pmatrix} 
\ = \ 
\begin{pmatrix}
\hsx 0_{\hsx } \\
\vdots\\
\hsx 0_{ \hsx}
\end{pmatrix}
.
\]
So in particular
\begin{align*}
0 \ 
&=\ 
q_1 \beta_1^{\sigma_i} + \cdots + q_d \beta_d^{\sigma_i} \qquad (1 \leq i \leq d)
\\[12pt]
&=\ 
\big(q_1 \beta_1\big)^{\sigma_i} + \cdots + \big(q_d \beta_d\big)^{\sigma_i}
\\[12pt]
&=\ 
\big(q_1 \beta_1 + \cdots + q_d \beta_d\big)^{\sigma_i}
\end{align*}
\qquad\qquad $\implies$
\[
\sum\limits_{i=1}^d \hsx
q_i \beta_i
\ = \ 
0
\]
\qquad\qquad $\implies$
\[
q_1 = 0, \ldots, q_d = 0.
\]
Therefore the $\bx_1, \ldots, \bx_d$ are $\Q$-linearly independent elements of $\Qbar^n$.
\vspace{0.2cm}

\qquad \textbullet \quad
Define $\by_j \in \C^n$ by
\[
\by_j \ = \ \big(\beta_j^{\sigma_1} \Lambda_1, \ldots, \beta_j^{\sigma_n} \Lambda_n\big) \
\qquad (1 \leq j \leq d).
\]
Since the matrix
\[
B 
\ = \ 
\big(\beta_k^{\sigma_i}\big)_{1 \hsx \leq \hsx i, k \hsx \leq \hsx d}
\]
has rank $d$, the $d \times n$ matrix
\[
B_n 
\ = \ 
\big(\beta_k^{\sigma_i}\big)_{1 \leq k \hsx \leq \hsx d, \ 1 \leq i \hsx \leq \hsx n}
\]
has rank $n$ (its $n$ columns are independent in $\K^d$).  
The product of $B_n$ by the $n \times n$ diagonal matrix 
\[
\diag (\Lambda_1, \ldots, \Lambda_n)
\]
is the $d \times n$ matrix whose row vectors are $\by_1, \ldots, \by_d$:
\[
\begin{pmatrix}
\beta_1^{\sigma_1} \Lambda_1 &\cdots &\beta_1^{\sigma_n} \Lambda_n\\
\vdots & &\vdots\\
\beta_d^{\sigma_1} \Lambda_1 &\cdots &\beta_d^{\sigma_n} \Lambda_n
\end{pmatrix}
\ = \ 
\begin{pmatrix}
\beta_1^{\sigma_1}  &\cdots &\beta_1^{\sigma_n} \\
\vdots & &\vdots\\
\beta_d^{\sigma_1} &\cdots &\beta_d^{\sigma_n}
\end{pmatrix}
\times
\begin{pmatrix}
\Lambda_1 &\cdots &0\\
\vdots & &\vdots\\
0 &\cdots &\Lambda_n
\end{pmatrix}
.
\]
Therefore the set $\{\by_1, \ldots, \by_d\}$ contains a basis for $\C^n$ over $\C$.]
\vspace{0.3cm}

The preceding considerations set the stage for an application of \S30, \#10,
hence at least one of the
\[
\bx_i \by_j \qquad (1 \leq i \leq d, \ 1 \leq j \leq n)
\]
does not belong to $\fL$, which, however is false.  
To see this, recall that
\[
\Lambda_{n+1} = 0, \ldots, \Lambda_d = 0,
\]
and write
\allowdisplaybreaks
\begin{align*}
\bx_i \by_j \ 
&=\ 
\sum\limits_{m=1}^n \hsx
\beta_i^{\sigma_m} \beta_j^{\sigma_m} \Lambda_m
\\[12pt]
&=\ 
\sum\limits_{m=1}^d \hsx
\beta_i^{\sigma_m} \beta_j^{\sigma_m} \Lambda_m
\\[12pt]
&=\ 
\sum\limits_{m=1}^d \hsx
\beta_i^{\sigma_m} \beta_j^{\sigma_m}
\sum\limits_{k=1}^d \hsx 
\beta_k^{\sigma_m} \lambda_k
\\[12pt]
&=\ 
\sum\limits_{k=1}^d \hsx 
c_{i j k} \lambda_k, 
\end{align*}
where
\allowdisplaybreaks
\begin{align*}
c_{i j k}
&=\ 
\sum\limits_{m=1}^d \hsx
\beta_i^{\sigma_m} \beta_j^{\sigma_m}  \beta_k^{\sigma_m} 
\\[12pt]
&=\ 
\tr(\beta_i \beta_j \beta_k) \in \Q.
\end{align*}
But $\fL$ is a $\Q$-vector space (cf. \S31, \#2).  
Consequently
\[
\bx_i \hsx \by_j \in \fL,
\]
a contradiction.
\vspace{0.3cm}

\qquad \un{Case 2:} \quad
None of the $\Lambda_i$ vanish.
\vspace{0.2cm}

[To begin with
\[
\Lambda_1 
\ = \ 
\sum\limits_{k=1}^d \hsx
\beta_k^{\sigma_1} \lambda_k
\ = \ 
\sum\limits_{k=1}^d \hsx
\beta_k\lambda_k \hsx
\in \hsx \Qbar
\]
by hypothesis.

\qquad \textbullet \quad
Define $\bx_k \in \C^d$ by 
\[
\bx_k \ = \ 
\big(\beta_k^{\sigma_1}, \ldots, \beta_k^{\sigma_d}\big) \qquad (1 \leq k \leq d).
\]
Since the matrix
\[
B 
\ = \ 
\big(\beta_k^{\sigma_i}\big)_{1 \hsx \leq \hsx i, k \hsx \leq \hsx d}
\]
is nonsingular, $\bx_1, \ldots, \bx_d$ are $\Q$-linearly independent elements of $\Qbar^{\hsx d}$ .
\vspace{0.2cm}

\qquad \textbullet \quad
Define $\by_j\in \C^d$ by 
\[
\by_j 
\ = \ \big(\beta_j^{\sigma_1}\Lambda_1, \ldots, \beta_j^{\sigma_d}\Lambda_d\big) \qquad (1 \leq j \leq d).
\]
Since \mB has rank $d$ and since
\[
\begin{pmatrix}
\beta_1^{\sigma_1} \Lambda_1 &\cdots &\beta_1^{\sigma_d} \Lambda_d\\
\vdots & &\vdots\\
\beta_d^{\sigma_1} \Lambda_1 &\cdots &\beta_d^{\sigma_d} \Lambda_d
\end{pmatrix}
\ = \ 
\begin{pmatrix}
\beta_1^{\sigma_1}  &\cdots &\beta_1^{\sigma_d} \\
\vdots & &\vdots\\
\beta_d^{\sigma_1} &\cdots &\beta_d^{\sigma_d}
\end{pmatrix}
\times
\begin{pmatrix}
\Lambda_1 &\cdots &0\\
\vdots & &\vdots\\
0 &\cdots &\Lambda_d
\end{pmatrix}
,
\]

it follows that $\by_1, \ldots, \by_d$ is a basis for $\C^d$ over $\C$.  In addition, 
\[
y_{1 j} 
\ = \ 
\beta_j^{\sigma_1} \Lambda_1 
\ = \ 
\beta_j \Lambda_1 \in \Qbar.]
\]

Therefore the assumptions of \S30, \#12 are satisfied, hence at least one of the
\[
\bx_i \hsx \by_j \qquad (1 \leq i \leq d, \ 1 \leq j \leq d)
\]
doen not belong to $\fL$.  On the other hand,
\[
\bx_i \hsx \by_j 
\ = \ 
\sum\limits_{k=1}^d \hsx
\tr(\beta_i \beta_j \beta_k) \lambda_k \in \fL
\]
and we again have a contradiction.
\vspace{0.3cm}

\qquad \un{Case 3:} \quad
All of the $\Lambda_i$ vanish.  Consider the system:
\[
\begin{matrix}
\Lambda_1 : \beta_1^{\sigma_1} \lambda_1  &+ \cdots + &\beta_d^{\sigma_1} \lambda_d = \ 0\\
\vdots { \qquad \qquad  }    & &   \qquad \qquad  \vdots\\
\Lambda_d : \beta_1^{\sigma_d} \lambda_1  &+ \cdots + &\beta_d^{\sigma_d} \lambda_d = \ 0
\end{matrix}
\quad.
\]
Its matrix is the transpose of \mB, thus is nonsingular, thus
\[
\lambda_1 = 0, \ldots, \lambda_d = 0,
\]
as desired.


%% file: _34_estimates.tex
\chapter{
$\boldsymbol{\S}$\textbf{34}.\quad  ESTIMATES}
\setlength\parindent{2em}
\setcounter{theoremn}{0}
\renewcommand{\thepage}{\S34-\arabic{page}}

\ \indent 
Given algebraic numbers $\alpha_1 \neq 0, 1, \ldots, \alpha_n \neq 0, 1$ and nonzero integers $b_1, \ldots, b_n$, put
\[
\Lambda \ = \ 
b_1 \hsx \Log \alpha_1 + \cdots b_n \hsx \Log \alpha_n.
\]
Then for the applications, it is important to estimate $\abs{\Lambda}$ from below.
\vspace{0.5cm}

\begin{x}{\small\bf NOTATION} \ 
Put
\[
B \ = \ \max\{2, \abs{b_1}, \ldots, \abs{b_n}\}.
\]
\end{x}
\vspace{0.3cm}

\begin{x}{\small\bf THEOREM} \ 
Assume: \ $\Lambda \neq 0$ $-$then
\[
\abs{\Lambda} \ \geq \ B^{-C},
\]
where $C > 0$ is a constant depending only on $n$ and $\alpha_1, \ldots, \alpha_n$.
\end{x}
\vspace{0.3cm}

\begin{x}{\small\bf REMARK} \ 
The reason for introducing the ``2'' is to accommodate the case when all the $b_i$ are $\pm 1$ since then
\[
\max\{\abs{b_1}, \ldots, \abs{b_n}\} \ = \ 1 
\quad \text{and} \quad 
1^C \ = \ 1.
\]
\end{x}
\vspace{0.3cm}

\begin{x}{\small\bf EXAMPLE} \ 
Suppose that $\ds\frac{p}{q}$ is a nonzero rational number with $q \geq 2$.  
Let $\alpha > 0$ $(\alpha \ne 1)$, $\alpha^\prime > 0$ $(\alpha^\prime \ne 1)$ be algebraic numbers $-$then
\[
\abs{p \hsx \elln(\alpha) - q \hsx \elln(\alpha^\prime)} 
\ \geq \ 
\frac{1}{\max\{\abs{p}, q\}^c}
\qquad \text{(cf. \S15, \#33)},
\]
where $c > 0$ depends only on $\elln(\alpha)$ and  $\elln(\alpha^\prime)$.
\vspace{0.3cm}

[Note: \ 
In the context of \S15, \#32, it is automatic that $\alpha^\prime \neq 1$.  For if
$\alpha^x = \alpha^\prime = 1$, then
\[
\elln(\alpha^x ) \ = \ \elln(1) 
\implies
x \elln(\alpha) \ = \ 0 
\implies 
\elln(\alpha) \ = \ 0 
\implies
\alpha \ = \ 1,
\]
which was ruled out at the beginning.]
\end{x}
\vspace{0.3cm}

Obviously
\allowdisplaybreaks
\begin{align*}
e^\Lambda \ 
&=\ 
\exp(b_1 \Log \alpha_1 + \cdots + b_n \Log \alpha_n)
\\[12pt]
&=\ 
\alpha_1^{b_1}, \ldots, \alpha_n^{b_n}.
\end{align*}
\vspace{0.3cm}

\begin{x}{\small\bf THEOREM} \ 
Assume
\[
\alpha_1^{b_1} \cdots \alpha_n^{b_n} \ \neq \ 1.
\]
Then
\[
\abs{\alpha_1^{b_1} \cdots \alpha_n^{b_n} - 1} \ \geq \ B^{\hsx -C},
\]
where $C > 0$ is a constant depending only on $n$ and $\alpha_1, \ldots, \alpha_n$.
\end{x}
\vspace{0.3cm}

Some elementary preliminaries are needed in order to make the transition from \#2 to \#5.
\vspace{0.2cm}

[Note: \ 
The ``C'' in \#5 is not the ``C'' in \#2.]
\vspace{0.3cm}

\begin{x}{\small\bf RAPPEL} \ 
\[
\Log z 
\ = \ 
\sum\limits_{n=1}^\infty \hsx
\frac{(-1)^{n-1}}{n} \hsx (z - 1)^n \qquad (\abs{z - 1} < 1) \quad \text{(cf. \S23, \#7)}.
\]
Put $z = 1 + w$, hence
\[
\Log (1 + w) 
\ = \ 
\sum\limits_{n=1}^\infty \hsx
\frac{(-1)^{n-1}}{n} \hsx w^n \qquad (\abs{w} < 1).
\]
\end{x}
\vspace{0.3cm}
\begin{x}{\small\bf LEMMA} \ 
\[
\abs{w} \ \leq \ \frac{1}{2}
\implies
\abs{\Log (1 + w)}  \ \leq \ 2 \abs{w}.
\]
\end{x}
\vspace{0.3cm}

Passing to the proof of \#5, put $w = \alpha_1^{b_1} \cdots \alpha_n^{b_n}$   
$-$then there are two possibilities.
\vspace{0.3cm}

\qquad \textbullet \quad
$\abs{w} > \ds\frac{1}{2}$. \ \ By definition,
\allowdisplaybreaks
\begin{align*}
&B = \max\{2, \abs{b_1}, \ldots, \abs{b_n}\}
\\[12pt]
\implies &B \geq 2 
\\[12pt]
\implies &\frac{1}{B} \leq \frac{1}{2} 
\\[12pt]
\implies &\frac{1}{B} < \abs{w},
\end{align*}
so $C = 1$ will work.
\vspace{0.2cm}

\qquad \textbullet \quad 
$\abs{w} \leq \ds\frac{1}{2}$. \ \ To begin with, for some $k \in \Z$, 
\allowdisplaybreaks
\begin{align*}
\Log (1 + w) \ 
&= 
\Log (\alpha_1^{b_1} \cdots \alpha_n^{b_n})
\\[12pt]
&= 
\Log \alpha_1^{b_1} + \cdots + \Log \alpha_n^{b_n} + 2 \pi \sqrt{-1} \hsx k \qquad \text{(cf. \S23, \#5)}
\\[12pt]
&= 
b_1 \Log \alpha_1 + \cdots + b_n \Log \alpha_n + 2 \pi \sqrt{-1} \hsx k.
\end{align*}
But
\allowdisplaybreaks
\begin{align*}
\Log \hsx {-1} \ 
&=\ 
\elln(\abs{-1}) + \pi \sqrt{-1}
\\[12pt]
&=\ 
\pi \sqrt{-1}.
\end{align*}
Therefore
\[
\Log (1 + w) \ = \ 
b_1 \Log \alpha_1 + \cdots + b_n \Log \alpha_n + 2 k \Log -1. 
\]
The right hand side has the form needed for an application of \#2 (ignore $2k \Log -1$ if $k = 0$), thus setting
\[
B_0 \ = \ 
\max\{2, \abs{b_1}, \ldots, \abs{b_n}, \abs{2k}\},
\]
it follows that
\[
\abs{\Log (1 + w)} \ \geq \ B_0^{\hsx -C_0}
\]
for some $C_0 > 0$.  Now estimate $\abs{2 \pi \sqrt{-1} \hsx k}$: 
\allowdisplaybreaks
\begin{align*}
\abs{2 \pi \sqrt{-1} \hsx k}\ 
&\leq\ 
\abs{\Log (1 + w)} + \sum\limits_{i=1}^n \hsx \abs{b_i} \hsx \abs{\Log \alpha_i}
\\[12pt]
&\leq\ 
2 \abs{w} + \sum\limits_{i=1}^n \hsx \abs{b_i} \hsx \abs{\Log \alpha_i}
\\[12pt]
&\leq\ 
1 + B \hsx \sum\limits_{i=1}^n \hsx \abs{\Log \alpha_i}
\\[12pt]
&\leq\ 
B \big(1 + \sum\limits_{i=1}^n \hsx \abs{\Log \alpha_i}\big)
\end{align*}
\qquad\qquad $\implies$
\allowdisplaybreaks
\begin{align*}
\abs{2 k}\ 
&\leq\ 
B \big(1 + \sum\limits_{i=1}^n \hsx \abs{\Log \alpha_i}\big) / \pi
\\[12pt]
&\leq\ 
B \big(1 + \sum\limits_{i=1}^n \hsx \abs{\Log \alpha_i}\big)
\\[12pt]
&\equiv\ 
C_1 B \qquad (C_1 > 1)
\end{align*}
\qquad\qquad $\implies$
\allowdisplaybreaks
\begin{align*}
B_0
&=\ 
\max \{B,\abs{2k}\}
\\[12pt]
&\leq\ 
\max \{B, C_1 B\}
\\[12pt]
&=\ 
C_1 B
\end{align*}
\qquad\qquad $\implies$
\allowdisplaybreaks
\begin{align*}
2 \abs{w}
&\geq\ 
\abs{\Log (1 + w)}
\\[12pt]
&\geq\ 
B_0^{\hsx -C_0}
\\[12pt]
&>\ 
(C_1 B)^{\hsx -C_0}
\end{align*}

\qquad\qquad $\implies$
\[
\abs{w} \ \geq \ \frac{1}{2} \hsx (C_1 B)^{\hsx -C_0} .
\]
Write
\[
2(C_1 B)^{C_0} \ = \ 2(C_1)^{C_0} B^{C_0}.
\]
Choose \mD:
\[
2(C_1)^{C_0} \ \leq \ B^D.
\]
Then
\begin{align*}
2(C_1)^{C_0} B^{C_0}  \ 
&\leq \ 
B^D \hsx B^{C_0}
\\[12pt]
&= \ B^{D + C_0}.
\end{align*}
Let $C = C + C_0$ to conclude that
\[
\frac{1}{2} \hsx (C_1 B)^{\hsx -C_0} \ \geq \ B^{\hsx -C},
\]
so 
\[
\abs{w}  \ \geq \ B^{\hsx -C}
\]
thereby completing the proof of \#5.
\vspace{0.3cm}

Under the circumstances, one can go beyond \#5.
\vspace{0.3cm}

\begin{x}{\small\bf THEOREM} \ 
Let
\[
\begin{cases}
\ a_1, \ldots, a_n\\
\ b_1, \ldots, b_n
\end{cases}
\]
be nonzero integers.  Assume: 
\[
a_1 \geq 2, \ldots, a_n \geq 2
\]
and
\[
a_1^{b_1} \cdots  a_n^{b_n} \ \neq \ 1.
\]
Then
\[
\abs{a_1^{b_1} \cdots a_n^{b_n} - 1} \ \geq \ \exp(-C(n) \hsx \elln(B) \hsx \elln(\alpha_1) \cdots \elln(a_n)),
\]
where $C(n) > 0$ is a constant depending only on $n$.
\end{x}
\vspace{0.3cm}

\begin{x}{\small\bf REMARK} \ 
According to Waldschmidt, an admissible value for $C(n)$ is
\[
2^{26 n} \hsx n^{3n}.
\]

\qquad {\small\bf FACT} \ 
If $\abs{b_1} \geq 2$, $\abs{b_2} \geq 2$, then
\[
\abs{b_1 \elln(2) + b_2 \elln(3)} 
\ \geq \ 
B^{-13.3}.
\]
\end{x}
\vspace{0.3cm}
\newpage

\[
\textbf{APPENDIX}
\]
\vspace{0.5cm}

{\small\bf DEFINITION} \ 
Complex numbers $\alpha_1, \ldots, \alpha_n$ are 
\un{multiplicatively independent}
\index{multiplicatively independent} 
if none are zero and if for any relation
\[
\alpha_1^{a_1} \cdots \alpha_n^{a_n} \ = \ 1,
\]
where $(a_1, \ldots, a_n) \in \Z^n$, there follows
\[
a_1 = 0, \ldots, a_n = 0.
\]
\vspace{0.3cm}

{\small\bf LEMMA} \ 
Suppose that $\alpha_1, \ldots, \alpha_n$  are multiplicatively independent $-$then for any choice $(\lambda_1, \ldots, \lambda_n) \in \C^n$ with 
$e^{\lambda_i} = \alpha_i$ $(1 \leq i \leq n)$, the $n + 1$ complex numbers 
$2 \pi \sqrt{-1}$, $\lambda_1, \ldots, \lambda_n$ are $\Q$-linearly independent.
\vspace{0.3cm}

Suppose given algebraic numbers $\alpha_1 \neq 0, 1, \ldots, \alpha_n \neq 0, 1$ and assume that they are multiplicatively independent, hence that
\[
\alpha_1^{k_1} \cdots \alpha_n^{k_n} \ \neq \ 1
\]
if the exponents are not all zero.
\vspace{0.3cm}

Turning to \#2, it can be shown that if
\[
\abs{\Lambda} \ < \ B^{\hsx -C}
\]
for a sufficiently large positive constant \mC depending only on $n$ and $\alpha_1, \ldots, \alpha_n$, then $\alpha_1, \ldots, \alpha_n$ must be 
multiplicatively dependent \ldots \hsx .


%% file: _35_matrices.tex
\chapter{
$\boldsymbol{\S}$\textbf{35}.\quad  MATRICES}
\setlength\parindent{2em}
\setcounter{theoremn}{0}
\renewcommand{\thepage}{\S35-\arabic{page}}

\ \indent 

Let \mA be an $m \times n$ matrix with entries in the complex numbers ($m$ rows and $n$ columns).
\vspace{0.3cm}

\begin{x}{\small\bf DEFINITION} \ 
The 
\un{column space}
\index{column space} 
of \mA is the vector space spanned by its columns and the 
\un{column rank}
\index{column rank} 
of \mA is the dimension of the column space of \mA.
\end{x}
\vspace{0.3cm}

\begin{x}{\small\bf DEFINITION} \ 
The 
\un{row space}
\index{row space} 
of \mA is the vector space spanned by its rows and the 
\un{row rank}
\index{row rank} 
of \mA is the dimension of the row space of \mA.
\end{x}
\vspace{0.3cm}

\begin{x}{\small\bf THEOREM} \ 
The column rank of \mA equals the row rank of \mA.
\end{x}
\vspace{0.3cm}

Therefore the number of linearly independent columns of \mA equals the number of linearly independent rows of \mA, 
their common value being the 
\un{rank of \mA}: $\rank A$.
\index{rank of \mA}
\vspace{0.2cm}

[Note: \ 
Only a zero matrix has rank 0.]
\vspace{0.3cm}

\begin{x}{\small\bf EXAMPLE} \ 
\[
\rank \ 
\begin{pmatrix}
 \hsm 1 &\hsm 0 &\hsx 1\\
-2 &-3 &\hsx1 \\
\hsm 3 &\hsm 3 &\hsx 0
\end{pmatrix}
\ =  \ 2.
\]
\vspace{0.2cm}

[The first two rows are linearly independent, so the rank is at least 2 but the three rows in total are linearly dependent (the third is equal to the second subtracted from the first), thus the rank is less than 3.]
\end{x}
\vspace{0.3cm}

\begin{x}{\small\bf \un{N.B.}} \ 
Denote by $A^\Tee$ the transpose of \mA $-$then
\[
\rank A \ = \ \rank A^\Tee.
\]
\end{x}
\vspace{0.3cm}

\begin{x}{\small\bf EXAMPLE} \ 
\[
\rank \ 
\begin{pmatrix}
\hsm 1 &\hsm1 &0 &\hsm 2\\
-1 &-1 &0 &-2
\end{pmatrix}
\ =  1.
\]
In fact, there are nonzero columns so the rank is positive.  
On the other hand
\[
\rank \ 
\begin{pmatrix}
1 &-1\\
1 &-1\\
0 &\hsm 0\\
2 &-2
\end{pmatrix}
\ =  1.
\]
\end{x}
\vspace{0.3cm}

\begin{x}{\small\bf LEMMA} \ 
The rank of \mA is the smallest integer $k$ such that \mA can be factored as a product $A = BC$, where \mB is an 
$m \times k$ matrix and \mC is a $k \times n$ matrix.  
\end{x}
\vspace{0.3cm}

\begin{x}{\small\bf LEMMA} \ 
The rank of \mA is the largest integer $r$ for which there exists a nonsingular $r \times r$ submatrix of \mA.
\vspace{0.2cm}

[Note: \ 
A nonsingular 
\un{$r$-minor}
\index{$r$-minor} 
is an $r \times r$ submatrix with nonzero determinant.] 
\end{x}
\vspace{0.3cm}

\begin{x}{\small\bf LEMMA} \ 
The rank of \mA is the smallest integer $k$ such that \mA can be written as a sum of $k$ rank 1 matrices.
\vspace{0.2cm}

[Note: \ 
A matrix has rank 1 if it can be written as a nonzero product $C  R$ of a column vector \mC and a row vector \mR:
\[
C \ = \ 
\begin{pmatrix}
c_1\\
\vdots\\
c_m
\end{pmatrix}
, \quad R \ = \ 
\begin{pmatrix}
r_1 &\cdots &r_n
\end{pmatrix}
\]
\qquad $\implies$
\[
C  R \ = \ 
\begin{pmatrix}
c_1 r_1 &\cdots &c_1 r_n\\[8pt]
\vdots &&\vdots \\[8pt] 
c_m r_1 &\cdots &c_m r_n
\end{pmatrix}
.]
\]
\end{x}
\vspace{0.3cm}

\begin{x}{\small\bf } \ 
Take \mA as in \#6 $-$then
\[
A \ = \ 
\begin{pmatrix}
\hsm 1\\
-1
\end{pmatrix}
\ 
\begin{pmatrix}
1 &1 &0 &2\\
\end{pmatrix}
\]
has rank 1.
\end{x}
\vspace{0.3cm}

\begin{x}{\small\bf LEMMA} \ 
The rank of \mA is $\leq \min\{m, n\}$.
\end{x}
\vspace{0.3cm}

\begin{x}{\small\bf DEFINITION} \ 
If
\[
\rank A \ = \ \min\{m, n\},
\]
then \mA is said to have 
\un{full rank};
\index{full rank} 
otherwise \mA is 
\un{rank deficient}.
\index{rank deficient}
\end{x}
\vspace{0.3cm}

\begin{x}{\small\bf LEMMA} \ 
If \mA is a square matrix (i.e., if $m = n$), then \mA is invertible iff \mA has rank $n$, thus is full rank.
\end{x}
\vspace{0.3cm}

\begin{x}{\small\bf LEMMA} \ 
If \mB is an $n \times k$ matrix, then 
\[
\rank A \hsx B 
\ \leq \ 
\min\{\rank A, \rank B\}
\]
and if $\rank B = n$, then 
\[
\rank A \hsx B 
\ = \ 
\rank A.
\]
\end{x}
\vspace{0.3cm}

\begin{x}{\small\bf LEMMA} \ 
The rank of \mA is equal to $r$ iff there exists an invertible
$m \times m$ matrix \mX and an invertible $n \times n$ matrix \mY such that
\[
X A Y \ = \ 
\begin{pmatrix}
I_r &0\\
0 &0
\end{pmatrix}
,
\]
where $I_r$ is the $r \times r$ identity matrix.
\end{x}
\vspace{0.3cm}

\begin{x}{\small\bf NOTATION} \ 
$\ov{A}$ is the complex conjugate of \mA and $A^*$ is the conjugate transpose of \mA.
\end{x}
\vspace{0.3cm}

\begin{x}{\small\bf LEMMA} \ 
\allowdisplaybreaks
\begin{align*}
\rank A\ 
&=\ 
\rank \ov{A}
\\[12pt]
&=\ 
\rank A^*
\\[12pt]
&=\ 
\rank A^* A
\\[12pt]
&=\ 
\rank A A^*.
\end{align*}
\end{x}
\vspace{0.3cm}

Attached to \mA is the linear map
\[
f_A : \C^n \ra \C^m
\]
defined by
\[
f_A(x) \ = \  A x.
\]
\vspace{0.3cm}

\begin{x}{\small\bf LEMMA} \ 
The rank of \mA equals the dimension of the image of $f_A$.
\end{x}
\vspace{0.3cm}

\begin{x}{\small\bf LEMMA} \ 
\vspace{0.2cm}

\qquad \textbullet \quad 
$f_A$ is injective iff $\rank A = n$.
\vspace{0.2cm}

\qquad \textbullet \quad 
$f_A$ is surjective iff $\rank A = m$.
\vspace{0.2cm}

\end{x}
\vspace{0.3cm}

\[
\textbf{APPENDIX}
\]
\vspace{0.5cm}

\qquad {\small\bf SIEGEL'S LEMMA} \ 
Let \\[2pt]
\[
\begin{cases}
\ c_{1 \hsy 1} x_1 + c_{1 \hsy 2} x_2 + \cdots + c_{1 \hsy n} x_n \ = \ 0\\[5pt]
\ \quad \vdots \\[5pt] 
\ c_{m \hsy 1} x_1 + c_{m \hsy 2} x_2 + \cdots + c_{m \hsy n} x_n \ = \ 0
\end{cases}
\]
\\[-2pt]
be a system of $m$ linear equations in $n > m$ unknowns $x_1, x_2, \ldots, x_n$, where $c_{i \hsy j} \in \Z$ not all 0.  
Let $C \in \N$:
\[
\max\limits_{i, j} \hsx \abs{c_{i \hsy j}} \ \leq \ C.
\]
There there exists a nontrivial integral solution $\{x_j \hsx \in \hsx \Z\}_{j = 1}^n$ with 
\[
\abs{x_j} 
\ \leq \ 
(n C)^{m / (n - m)}.
\]


%% file: _36_the_six_exponentials_theorem.tex
\chapter{
$\boldsymbol{\S}$\textbf{36}.\quad  SIX EXPONENTIALS THEOREM}
\setlength\parindent{2em}
\setcounter{theoremn}{0}
\renewcommand{\thepage}{\S36-\arabic{page}}

\ \indent 
This is the following statement.
\vspace{0.5cm}

\begin{x}{\small\bf THEOREM} \ 
Suppose given $\Q$-linearly independent complex numbers
\[
\{x_1, \ldots, x_m\} \ \text{and} \ \{y_1, \ldots, y_n\}.
\]
Assume: 
\[
m \hsy n \ > \  m + n.
\]
Then at least one of the numbers
\[
\exp(x_i \hsy y_j) \qquad (1 \leq i \leq m, 1 \leq j \leq n)
\]
is transcendental.
\vspace{0.2cm}

[As regards the proof, one can extend the ideas used in the proof of Gelfond-Schneider but we shall omit the details opting instead for a ``geometric argument'' later on (cf. \S41, \#1).]
\end{x}
\vspace{0.3cm}

Special Cases: \ $m = 3$, $n = 2$ or $m = 2$, $n = 3$.
\vspace{0.3cm}

\begin{x}{\small\bf EXAMPLE} \ 
Take
\[
x_1 = 1, \ x_2 = e, \ y_1 = e, \ y_2 = e^2, \ y_3 = e^3,
\]
where \S17, \#2 has been silently invoked $-$then the six exponentials are 
\[
e^e, \ e^{e^2}, \ e^{e^3}, \ e^{e^2}, \ e^{e^3}, \ e^{e^4},
\]
thus at least one of the numbers
\[
e^e, \ e^{e^2}, \ e^{e^3},  \ e^{e^4}
\]
is transcendental.
\end{x}
\vspace{0.3cm}

\begin{x}{\small\bf EXAMPLE} \ 
Take
\[
x_1 = 1, \ x_2 = \pi, \ y_1 = \elln(2), \ y_2 = \pi \elln(2), \ y_3 = \pi^2 \elln(2).
\]
Then the six exponentials are 
\[
2, \ 2^\pi, \ 2^{\pi^2}, \ 2^\pi, \ 2^{\pi^2}, \ 2^{\pi^3},
\]
thus at least one of the numbers
\[
2^\pi, \ 2^{\pi^2}, \ 2^{\pi^3}
\]
is transcendental.
\vspace{0.2cm}

[Note: \ 
Consider a dependence relation
\[
q_1 \elln(2) \hsx + \hsx q_2 \pi \elln(2) \hsx + \hsx q_3 \pi^2 \elln(2) \ = \ 0
\]
where $q_1, q_2, q_3 \in \Q$ $-$then 
\[
q_1 \hsx + \hsx q_2 \pi \hsx + \hsx q_3 \pi^2 \ = \ 0
\]
\qquad\qquad\qquad $\implies$
\[
q_1 = 0, \ q_2 = 0, \ q_3 = 0,
\]
$\pi$ being transcendental (cf. \S19, \#1).]
\end{x}
\vspace{0.3cm}

\begin{x}{\small\bf REMARK} \ 
It is unknown whether one of the numbers 
\[
2^\pi, \quad 2^{\pi^2}
\]
is transcendental.
\end{x}
\vspace{0.3cm}

\begin{x}{\small\bf EXAMPLE} \ 
Fix $t \in \R$, $t \notin \Q$.  Take
\[
x_1 = 1, \  x_2 = t, \ y_1 = \elln(2), \ y_2 = \elln(3), \ y_3 = \elln(5).
\]
Then the six exponentials are
\[
2, \ 3, \ 5, \ 2^t, \ 3^t, \ 5^t, 
\]
thus at least one of the numbers 
\[
 2^t, \ 3^t, \ 5^t
\]
is transcendental.
\vspace{0.2cm}

[Note: \ 
$\elln(2)$, $\elln(3)$, $\elln(5)$ are $\Q$-linearly independent.  
To see this, consider a dependence relation
\[
q_1 \elln(2) + q_2 \elln(3) + q_3 \elln(5) \ = \ 0,
\]
where $q_1, q_2, q_3 \in \Q$.  
Write
\[
q_1 \ = \ \frac{m_1}{n_1}, \ \ 
q_2 \ = \ \frac{m_2}{n_2}, \ \ 
q_3 \ = \ \frac{m_3}{n_3}.
\]
Here
\[
n_1 \neq 0, \ n_2 \neq 0, \ n_3 \neq 0
\]
and the claim is that 
\[
m_1 = 0, \ m_2 = 0 , \ m_3 = 0.
\]
Clear the denominators and exponentiate to get
\[
2^{m_1 n_2 n_3} \hsx
3^{n_1 m_2 n_3} \hsx
5^{n_1 n_2 m_3} 
\ = \ 1
\]
\qquad\qquad\qquad $\implies$
\[
m_1 n_2 n_3 \ = \ 0, \ \ 
n_1 m_2 n_3 \ = \ 0, \ \ 
n_1 n_2 m_3 \ = \ 0,
\]
so
\[
m_1 = 0, \ m_2 = 0, \ m_3 = 0 \hsx.]
\]
\end{x}
\vspace{0.3cm}

\begin{x}{\small\bf DEFINITION} \ 
Let
\[
E_\infty 
\ = \ 
\{t \in \R : 2^t, 3^t, 5^t, \ldots \in \N\}.
\]
Then $E_\infty = \N$.
\vspace{0.2cm}

[Introduce
\allowdisplaybreaks
\begin{align*}
E_1 \ 
&=\ 
\{t \in \R : 2^t  \in \N\}
\\[12pt]
E_2 \ 
&=\ 
\{t \in \R : 2^t, 3^t  \in \N\}
\\[12pt]
E_3 \ 
&=\ 
\{t \in \R : 2^t, 3^t, 5^t \in \N\}.
\end{align*}
Then
\[
\N \hsx \subset  \hsx E_\infty \hsx \subset \hsx E_3 \hsx \subset \hsx E_2 \hsx \subset \hsx E_1.
\]

\noindent
Actually, we shall prove a stronger result, namely that $E_3 \ = \ \N$.  In fact, 
\allowdisplaybreaks
\begin{align*}
\N \ 
&\subset \ 
E_\infty
\\[12pt]
&\subset \ 
E_3
\\[12pt]
&=\ 
(E_3\hsx \cap \Q)  \ \amalg \hsx (E_3 \hsx \cap \hsx (\R - \Q))
\\[12pt]
&=\ 
E_3\hsx \cap \Q  \qquad \text{(cf. \#5)}
\\[12pt]
&\subset \ 
E_1\hsx \cap \Q 
\\[12pt]
&= \ 
\N \qquad \text{(cf. \#7)}.
\end{align*}
So, 
\[
\N \ = \ E_\infty \ = \ E_3.
\]
\\[-10pt]


[Note: \ 
True or False: \ $E_2 = \N$ \ \  (cf. \S44, \#6).]
\end{x}
\vspace{0.3cm}

\begin{x}{\small\bf \un{N.B.}} \ 
By definition,
\[
E_1 \ = \ \{t \in \R : 2^t  \in \N\}.
\]
And 
\allowdisplaybreaks
\begin{align*}
2^t = n \in \N \  
&\implies 
\elln\bigl(2^t\bigr) \hsx = \hsx \elln(n)
\\[12pt]
&\implies 
t \hsx = \hsx \frac{\elln(n)}{\elln(2)}.
\end{align*}
If $t \in \Q$, say $t = \ds\frac{p}{q}$ $-$then 

\allowdisplaybreaks
\begin{align*}
2^{\frac{p}{q}} \ = \ n \ 
&\implies 
2^p \ = \ n^q
\\[12pt]
&\implies 
n \ = \ 2^k \qquad (\exists \ k \in \N)
\\[12pt]
&\implies 
p = k \hsy q
\\[12pt]
&\implies 
t \ = \ \frac{p}{q} \ = \ k \in \N.
\end{align*}
Therefore 
\[
E_1 \hsx \cap \hsx \Q \ = \ \N.
\]


\end{x}
\vspace{0.3cm}

\begin{x}{\small\bf EXAMPLE} \ 
Let $x_1$, $x_2$ be two elements of $\R \cup \sqrt{-1} \hsx \R$ which are $\Q$-linearly independent.  
Let $y_1$, $y_2$ be two complex numbers subject to $y_1$, $y_2$, $\ov{y_2}$ being $\Q$-linear independent 
$-$then at least one of the numbers
\[
e^{x_1 y_1}, \ e^{x_1 y_2}, \ e^{x_2 y_1}, \ e^{x_2 y_2}
\]
is transcendental.
\vspace{0.2cm}

[Taking $y_3 = \ov{y_2}$, \#1 is applicable so it is a matter of eliminating $e^{x_1 y_3}$, $e^{x_2 y_3}$ from consideration.  
E.g.: 
\vspace{0.2cm}

(1) \ 
Suppose $x_1 \in \R$ $-$then
\[
e^{x_1 y_3} 
\ = \ 
e^{x_1 \ov{y}_2}
\ = \ 
e^{\ov{\ov{x}}_1 \hsx \ov{y}_2}
\ = \ 
e^{\ov{\ov{x}_1 y_2}}
\ = \ 
e^{\ov{x_1 y_2}}.
\]
But $e^{\ov{x_1 y_2}}$ is transcendental iff $e^{{x_1 y_2}}$ is transcendental.
\vspace{0.2cm}

(2) \ 
Suppose $x_1 \in \sqrt{-1} \hsx \R$ $-$then
\[
e^{x_1 y_3} 
\ = \ 
e^{x_1 \ov{y}_2}
\ = \ 
e^{\ov{\ov{x}}_1 \ov{y}_2}
\ = \ 
e^{\ov{\ov{x}_1 y_2}}
\ = \ 
e^{\ov{-x_1 y_2}}.
\]
But 
\[
e^{{-x_1 y_2}}
\ = \ 
\frac{1}{e^{x_1 y_2}}
\]
is transcendental iff $e^{{x_1 y_2}}$ is transcendental.
Meanwhile $e^{{-x_1 y_2}}$  is transcendental iff $\ds e^{\ov{-x_1 y_2}}$ is transcendental.]
\vspace{0.2cm}

[Note: \ 
$\alpha$ transcendental $\Leftrightarrow$ $\ov{\alpha}$ transcendental 
and $\alpha$ transcendental iff $\ds\frac{1}{\alpha}$ transcendental.]
\end{x}
\vspace{0.3cm}

\begin{x}{\small\bf LEMMA} \ 
Consider a nonzero $m \times n$ matrix
\[
M \ = \ 
\begin{pmatrix}
\lambda_{1 1} &\ldots &\lambda_{1 n}\\[8pt]
\vdots & &\vdots \\[8pt]
\lambda_{m 1} &\ldots &\lambda_{m n}
\end{pmatrix}
,
\]
where $\lambda_{i j} \in \fL$.  Assume: 
\vspace{0.2cm}

\qquad \textbullet \quad The $m$ rows
\[
[\lambda_{1 1}, \ldots, \lambda_{1 n}], 
\ldots, 
[\lambda_{m 1}, \ldots, \lambda_{m n}]
\]
are $\Q$-linearly independent in $\C^n$.
\vspace{0.2cm}

\qquad \textbullet \quad The $n$ columns
\[
\begin{pmatrix}
\lambda_{1 1}\\[8pt]
\vdots\\[8pt]
\lambda_{m 1}
\end{pmatrix}
, \ldots, 
\begin{pmatrix}
\lambda_{1 n}\\[8pt]
\vdots\\[8pt]
\lambda_{m n}
\end{pmatrix}
\]
are $\Q$-linearly independent in $\C^m$.
\\[8pt]
Then
\[
m \hsy n \  > \  m + n
\]
implies that the rank of \mM is $\geq 2$.
\vspace{0.2cm}

PROOF \ 
To get a contradiction, suppose that 
\[
\rank M \ = \ 1.
\]
Write (cf. \S35, \#9)
\[
\lambda_{i j} \ = \ x_i y_j.
\]
The point then is to check that the conditions of \#1 are satisifed, i.e., that 
\[
\begin{cases}
\ x_1, \ldots, x_m\\[8pt]
\ y_1, \ldots, y_n
\end{cases}
\text{are $\Q$-linearly independent.}
\]
For then the conclusion is that there is a pair $(x_i,y_j)$ such that 
\[
\exp(x_i y_j)
\]
is transcendental.  
But
\[
\exp(x_i y_j) 
\ = \ 
\exp(\lambda_{i j}) \in \ov{\Q}^\times,
\]
a contradiction.  
So consider the dependence relations
\[
\begin{cases}
\ q_1 x_1 + \cdots + q_m x_m \ = \ 0\\[8pt]
\ p_1 y_1 + \cdots + p_n y_n \ = \ 0
\end{cases}
\qquad (q_i \in \Q, \ p_j \in \Q)
\]
and for the sake of argument, set down a generic rational dependence relation for 
the columns:
\[
A_1
\begin{pmatrix}
x_1 y_1\\
\vdots\\
x_m y_1
\end{pmatrix}
\ + \cdots + \ 
A_n 
\begin{pmatrix}
x_1 y_n\\
\vdots\\
x_m y_n
\end{pmatrix}
\ = \ 
\begin{pmatrix}
0\\
\vdots\\
0
\end{pmatrix}
\ \in \C^m
\]
\qquad\qquad $\implies$
\[
\begin{cases}
\ A_1 x_1 y_1 + \cdots + A_n x_1 y_n \ = \ 0\\
\hspace{2cm} \vdots\\
\ A_1 x_m y_1 + \cdots + A_n x_m y_n \ = \ 0
\end{cases}
.
\]
We have
\[
p_1 y_1 + \cdots + p_n y_n \ = \ 0
\]
\qquad\qquad\qquad $\implies$
\begin{align*}
p_1 x_1 y_1 \ +
&\cdots + p_n x_1 y_n  \ = \ 0
\\[8pt]
&\quad \vdots
\\[8pt]
p_1 x_m  y_1 + &\cdots + p_n x_m  y_n \ = \ 0 
\end{align*}
Take now 
\[
A_1 = p_1, \ldots, \ A_n \ = \ p_n.
\]
Since by hypothesis, the columns are $\Q$-linearly independent in $\C^m$, it follows that 
$A_1 = 0, \ldots, A_n = 0$, or still, $p_1 = 0, \ldots, p_n = 0$ \hsx.]
\end{x}
\vspace{0.3cm}

\begin{x}{\small\bf SCHOLIUM} \ 
Take $m = 2$, $n = 3$, and consider a nonzero $2 \times 3$ matrix \mM with entries in $\fL$:
\[
M \ = \ 
\begin{pmatrix}
\lambda_{1 1} &\lambda_{1 2} &\lambda_{1 3}\\[8pt]
\lambda_{2 1} &\lambda_{2 2} &\lambda_{2 3}
\end{pmatrix}
\]
Suppose that its rows are $\Q$-linearly independent and its columns are $\Q$-linearly independent $-$then in view of \#9, the rank of \mM is $\geq 2$.  
However, on general grounds (cf. \S35, \#11), the rank of \mM is $\leq \min(2,3) = 2$.  
Therefore
\[
\rank M = 2, 
\]
hence \mM has full rank (cf. \S35, \#12). 
\end{x}
\vspace{0.3cm}

\begin{x}{\small\bf \un{N.B.}} \  \ 
We have seen above that $\#1 \implies \#9$.  
The converse is also true: $\#9 \implies \#1$.
\vspace{0.2cm}

[To begin with, the assumption that 
\[
\{x_1, \ldots, x_m \} \quad \text{and} \quad \{y_1, \ldots, y_n\}
\]
are $\Q$-linearly independent implies the $\Q$-linear independence of the rows and columns of \mM.  
E.g.: \ To deal with the columns, note that there is at least one $x_i \neq 0$, say $x_1 \neq 0$, thus from 
\[
A_1 x_1 y_1 + \cdots + A_n x_1 y_n = 0
\]
there follows
\[
A_1 y_1 + \cdots + A_n y_n = 0
\]
\qquad\qquad\qquad $\implies$
\[
A_1 = 0, \ldots, A_n = 0.
\]
Put $\lambda_{i j} = x_i y_j$ and suppose that $\forall \ i, \ j$ : $\lambda_{i j} \in \fL$ $-$then the rank of 
\[
M \ = \ 
\begin{pmatrix}
\lambda_{1 1} &\cdots &\lambda_{1 n}\\
\vdots &&\vdots\\
\lambda_{m 1} &\cdots &\lambda_{m n}
\end{pmatrix}
\]
is $\geq 2$ (bear in mind that $m n > m + n$).  
But this is false: \ $\rank M = 1$.  
Consequently \ $\exists \ i, \ j$ :$\lambda_{i j} \notin \fL$, so 
\[
\exp(\lambda_{i j}) 
\ = \ 
\exp(x_i y_j)
\]
is transcendental.
\end{x}
\vspace{0.75cm}

\[
\textbf{APPENDIX}
\]
\vspace{0.5cm}

{\small\bf QUESTION} \ 
If $m n / (m + n)$ is large, can one find a lower bound for the rank of \mM which is $> 2$?  
Without additional conditions, the answer is ``no''.  
To see this, consider 
\[
M_m \ = \ 
\begin{pmatrix}
\elln(2) &\elln(3) &\ldots &&\elln(p_m)\\
\elln(3)\\
\vdots &&0 \\
\elln(p_m)
\end{pmatrix}
,
\]
where $p_m$ is the $m^\nth$ prime $-$then rank $M_m = 2$ for each $m > 2$
(here $m = n$ and $m^2 > 2m \implies m > 2$).  
Therefore the mere $\Q$-linear indpendence of the rows and the columns does not suffice.
\vspace{0.5cm}

{\small\bf CRITERION} \ 
Let
\[
M \ = \ 
\begin{pmatrix}
\lambda_{1 1} &\cdots &\lambda_{1 n}\\[8pt]
\vdots &&\vdots\\
\lambda_{m 1} &\cdots &\lambda_{m n}
\end{pmatrix}
\]
be an $m \times n$ matrix with terms in $\fL$.  
Assume: 
\[
\begin{cases}
\ \forall \ (t_1, \ldots, t_m) \in \Z^m \ - \ \{(0, \ldots, 0)\}\\\\[8pt]
\ \forall \ (s_1, \ldots, s_n) \in \Z^n \ - \ \{(0, \ldots, 0)\}
\end{cases}
,
\]
the sum 
\[
\sum\limits_{i=1}^m \hsx\hsx
\sum\limits_{j=1}^n \hsx
t_i \hsy s_j \hsy \lambda_{i j} \ \neq \ 0.  
\]
Then the rank of \mM is 
\[
\geq \  \frac{m n}{m + n}.
\]
\vspace{0.2cm}

[Note: \ 
\[
 \ \lambda_{i j} \neq 0 \qquad (\forall \ i, j) \hsx.]
\]
\vspace{0.5cm}

{\small\bf EXAMPLE} \ 
Take $m = d > 1$, $n = d > 1$ $-$then 
\[
\frac{m n}{m + n} 
\ = \ 
\frac{d^2}{2d} 
\ = \ 
\frac{d}{2}.
\]
\vspace{0.5cm}

{\small\bf LEMMA} \ 
Under these circumstances, the rows and columns are $\Q$-linearly independent.
\vspace{0.2cm}

PROOF \ 
Consider
\[
A_1
\begin{pmatrix}
\lambda_{1 1}\\[8pt]
\vdots\\[8pt]
\lambda_{m 1}
\end{pmatrix}
+ \cdots + 
A_n
\begin{pmatrix}
\lambda_{1 n}\\[8pt]
\vdots\\[8pt]
\lambda_{m n}
\end{pmatrix}
,
\]
where without loss of generality, the $A_j \in \Z$ are not all zero $-$then the claim is that this expression is $\neq 0$.  
To be specific, assume $A_1 \neq 0$ and tailor the expression
\[
\sum\limits_{i=1}^m \hsx\hsx
\sum\limits_{j=1}^n \hsx
t_i \hsy s_j \hsy \lambda_{i j}
\]
as follows: \ Choose
\[
t_1 = 1, \ t_2 = 0, \ldots, \ t_m = 0
\]
to get 
\[
\sum\limits_{j=1}^n \hsx 
s_j \hsx \lambda_{1 j} 
\ = \ 
s_1 \lambda_{1 1} + s_2 \lambda_{1 2} + \cdots + s_n \lambda_{1 n} \neq 0.
\]
Take
\[
s_1 \ = \ A_1, \ 
s_2 \ = \ A_2, \ 
\ldots, \ 
s_n \ = \  A_n, 
\]
hence
\[
A_1 \lambda_{1 1} + A_2 \lambda_{1 2} + \cdots + A_n \lambda_{1 n} \  \neq \  0.
\]

Assume in addition that
\[
m \hsy n \ > \  m + n.
\]
Then what has been said above implies \#9 which in turn implies \#1 (cf. \#11).
\vspace{0.5cm}

{\small\bf EXAMPLE} \ 
Take $m = d > 1$, $n = d > 1$ $-$then the foregoing says that the rank of \mM is $\geq \ds\frac{d}{2}$.  
On the other hand, the theory also says that the rank of \mM is $\geq 2$ (cf. \#9).  
To check consistency, note that
\[
m n > m + n 
\quad \text{becomes} \quad 
d^2 > 2d 
\implies 
d > 2 
\implies
\frac{d}{2} > 1.
\]
\vspace{0.1cm}

\un{Case 1:} \quad
$d = 2r$ $(r = 1, 2, \ldots)$ $-$then
\[
1 < \frac{d}{2} 
\ = \ 
r 
\implies 
r \geq 2
\]
\qquad\qquad\qquad $\implies$
\[
2 \leq r \leq \rank M.
\]
\vspace{0.1cm}

\un{Case 2:} \quad 
$d = 2r + 1$ $(r = 1, 2, \ldots)$ $-$then
\vspace{0.2cm}

\qquad \un{$r = 1$:} \ 
Here
\[
\frac{d}{2} 
\ = \ 
\frac{3}{2} \ \leq \ \rank M.
\]
But rank \mM is a positive integer, so $\rank M \geq 2$.
\vspace{0.2cm}

\qquad \un{$r > 1$:} \ 
Simply write
\[
2 \  \leq \ r \  \leq \ \frac{2 r + 1}{2}
\ = \ 
\frac{d}{2} \  \leq \  \rank M.
\]
Therefore matters are in fact consistent.

%% file: _37_vector_spaces.tex
\chapter{
$\boldsymbol{\S}$\textbf{37}.\quad  VECTOR SPACES}
\setlength\parindent{2em}
\setcounter{theoremn}{0}
\renewcommand{\thepage}{\S37-\arabic{page}}

\ \indent 
Let $\K$ be a field, $\bk \subset \K$ a subfield.
\vspace{0.5cm}

\begin{x}{\small\bf \un{N.B.}} \ 
Typically
\[
\K \ = \ \C, \ \ 
\bk = \Qbar \ \text{or} \ \Q.
\]
\end{x}
\vspace{0.3cm}

\begin{x}{\small\bf LEMMA} \ 
Let $\sV \subset \K^d$ be a $\K$-vector subspace $-$then the following conditions are equivalent.
\\[-11pt]

\qquad (i) \ 
$\sV$ has a basis whose elements belong to $\bk^d$.
\vspace{0.2cm}

\qquad (ii) \ 
$\sV$ is the intersection of hyperplanes defined by linear forms with coefficients in \bk.
\vspace{0.2cm}

[Note: \ Such a subspace $\sV$ is said to be 
\un{rational over \bk}.]
\index{rational over \bk}
\end{x}
\vspace{0.3cm}

\begin{x}{\small\bf DEFINITION} \ 
Let $\sV$ be a $\K$-vector subspace $-$then a 
\un{\bk-structure on $\sV$}
\index{\bk-structure on $\sV$} 
is a \bk-vector subspace $\sV^\prime$ of $\sV$ such that any basis of 
$\sV^\prime$ over \bk is a basis of $\sV$ over $\K$.
\end{x}
\vspace{0.3cm}

\begin{x}{\small\bf LEMMA} \ 
Let $\sV \subset \K^d$ be a $\K$-vector subspace $-$then $\sV \cap \bk^d$ is a \bk-structure on $\sV$ iff $\sV$ is rational over \bk.
\end{x}
\vspace{0.3cm}

\begin{x}{\small\bf EXAMPLE} \ 
\\

\qquad \textbullet \quad 
$\Q^d$ is a $\Q$-structure on $\C^d$.
\vspace{0.2cm}

\qquad \textbullet \quad
$\Qbar^{\hsx d}$ is a $\Qbar$-structure on $\C^d$.
\end{x}
\vspace{0.3cm}

\begin{x}{\small\bf DEFINITION} \ 
Given $\K$-vector subspaces
\[
\begin{cases}
\ \sV_1 \subset \K\strutz^{d_1} \\[8pt]
\ \sV_2 \subset \K\strutz^{d_2} 
\end{cases}
\]
endowed with \bk-structures
\[
\begin{cases}
\ \sV_1^\prime \subset \sV_1 \\[8pt]
\ \sV_2^\prime \subset \sV_2
\end{cases}
,
\]
\end{x}
a $\K$-linear map $f:\sV_1 \ra \sV_2$ is 
\un{rational over \bk}
\index{rational over \bk ($\K$-linear map)}
if $f(\sV_1^\prime) \subset \sV_2^\prime$.
\vspace{0.3cm}

\begin{x}{\small\bf EXAMPLE}  \ 
Take $\sV_1 = \C\strutz^{d_1}$, $\sV_2 = \C\strutz^{d_2}$ to arrive at the notion of a $\C$-linear map 
$f:\C\strutz^{d_1}\ra \C\strutz^{d_2}$ which is rational over $\Q$ (or $\Qbar$).
\end{x}
\vspace{0.5cm}

\[
\textbf{APPENDIX}
\]
\vspace{0.5cm}

\qquad{\small\bf NOTATION} \ 
Let $\be_1, \ldots, \be_d$ be the canonical basis for $\K\strutz^d$.
\vspace{0.3cm}

Let $\sV \subset \K\strutz^d$ be a $\K$-vector subspace of dimension $n$.  
Consider the following properties.
\vspace{0.2cm}

\qquad (1) \quad 
If $\pi_{\sV}:\K^d \ra \K^d / \sV$ is the canonical projection, then $(\pi_{\sV}(\be_1), \ldots, \pi_{\sV}(\be_{d-n}))$ is a basis for $\K^d / \sV$.
\vspace{0.2cm}

\qquad (2) \quad 
Given $\bz = (z_1, \ldots, z_d) \in \sV$, the conditions
\[
z_{d - n + 1} = \cdots = z_d \ = \ 0 
\implies
\bz = \boldsymbol{0}.
\]

\qquad (3) \ 
The restriction to $\sV$ of the projection $\K\strutz^d \ra \K\strutz^n$ of the last $n$ coordinates is injective.
\vspace{0.2cm}

\qquad (4) \ 
$\sV$ is the intersection of $d - n$ hyperplanes defined by the equations
\[
z_j 
\ = \ 
\sum\limits_{i = d - n + 1}^d \hsx
a_{i j} z_i \qquad (1 \leq j \leq d - n).
\]

\vspace{0.2cm}

{\small\bf FACT} \ 
Properties (1), (2), (3), (4) are equivalent.

%% file: _38_vector_spaces_L.tex
\chapter{
$\boldsymbol{\S}$\textbf{38}.\quad  VECTOR SPACES : \ $\fL$}
\setlength\parindent{2em}
\setcounter{theoremn}{0}
\renewcommand{\thepage}{\S38-\arabic{page}}

\ \indent 
Recall that in \S32, \#3, various conditions were formulated which are equivalent to homogeneous Baker.  
What follows is a supplement to that list.
\vspace{0.3cm}

\begin{x}{\small\bf THEOREM} \ 
The following assertions are equivalent to homogeneous Baker.
\vspace{0.3cm}

\qquad (i) \ 
Let $\sV \subset \C^d$ be a $\C$-vector subspace rational over $\Qbar$ with $\sV  \hsx \cap \hsx \Q^d = \{0\}$ $-$then
$\sV\hsx \cap \hsx \fL^d = \{0\}$.
\vspace{0.3cm}

\qquad (ii) \ 
Let $\sV \subset \C^d$ be a $\C$-vector subspace rational over $\Qbar$ $-$then there exists a $\C$-vector subspace $\sV_0$ of 
$\C^d$ rational over $\Q$ and contained in $\sV$ such that
\[
\sV\hsx \cap \hsx \fL^d 
\ = \ 
\sV_0\hsx \cap \hsx \fL^d.
\]

[E.g.: \ 
To see that $(ii) \implies (i)$, note that if $\sV \hsx \cap \hsx \Q^d = \{0\}$, then the only $\C$-vector subspace $\sV_0$ of $\C^d$ rational over $\Q$ and contained in $\sV$ is $\{0\}$, hence
\[
\sV\hsx \cap \hsx \fL^d 
\ = \ 
\sV_0\hsx \cap \hsx \fL^d
\ = \ 
\{0\}\hsx \cap \hsx \fL^d
\ = \ 
\{0\} .]
\]
\end{x}
\vspace{0.3cm}

\begin{x}{\small\bf REMARK}  \ 
One can replace item (ii) by a weaker assertion, viz.: \ If $\sV \subset \C^d$ is a 
$\C$-vector subspace rational over $\Qbar$, then
\[
\sV\hsx \cap \hsx \fL^d 
\ = \ \bigcup\limits_{\sV_0} \hspace{.1cm} \sV_0\hsx \cap \hsx \fL^d,
\]
where $\sV_0$ ranges over the $\C$-vector subspaces of $\C^d$ rational over $\Q$ and contained in $\sV$.
\end{x}
\vspace{0.3cm}

\begin{x}{\small\bf THEOREM} \ 
Let $\sV \subset \C^d$ be a $\C$-vector subspace $-$then the  $\Q$-vector space $\sV \subset L^d$ is finite dimensional iff 
$\sV \hsx \cap \hsx \Q^d = \{0\}$.

The implication
\[
\dim_{\Q} (\sV\hsx \cap \hsx \fL^d)  < \infty
\implies
V\hsx \cap \hsx \Q^d \ = \ \{0\},
\]
i.e., 
\[
V \cap \Q^d \ \neq \ \{0\}
\implies
\dim_{\Q} (\sV\hsx \cap \hsx \fL^d)  \ = \  \infty
\]
is straightforward: \ Take
\[
\bq = (q_1, \ldots, q_d) \ \neq \ 0
\]
in $\sV \cap \Q^d$ $-$then $\forall \ \lambda \in \fL$, 
\[
(q_1 \lambda, \ldots, q_d \lambda) \in \sV \hsx \cap \hsx \fL^d 
\implies
\dim_{\Q} (\sV\hsx \cap \hsx \fL^d) = \infty.
\]

As for the converse, i.e., 
\[
\sV \cap \Q^d = \{0\} 
\implies
\dim_{\Q} (\sV\hsx \cap \hsx \fL^d) < \infty,
\]
it is not so easy to establish.  
However there is one situation when matters are immediate.  
For suppose that $\sV \hsx \cap \hsx \Q^d = \{0\}$ \ AND \ in addition that $\sV$ is rational over 
$\Qbar$ $-$then $\sV \hsx \cap \hsx \fL^d = \{0\}$ 
(cf. \#1 (i)).
\end{x}
\vspace{0.3cm}

\begin{x}{\small\bf \un{N.B.}}  \ 
If $\sV$ is not rational over $\Qbar$ but $\sV \cap \Q^d = \{0\}$, then
\[
\dim_{\Q} (\sV\hsx \cap \hsx \fL^d) 
\]
may very well be positive (but, of course, finite) (cf. \#7).
\end{x}
\vspace{0.3cm}

\begin{x}{\small\bf THEOREM} \ 
Let $\sV \subset \C^d$ be a $\C$-vector subspace such that $\sV \hsx \cap \hsx \Q^d = \{0\}$ $-$then
\[
\dim_{\Q} (\sV\hsx \cap \hsx \fL^d) \ \leq \  n (n+1),
\]
where
\[
n \ = \ \dim_{\C} (\sV).
\]
\end{x}
\vspace{0.3cm}

\begin{x}{\small\bf EXAMPLE} \ 
Take for $\sV$ a complex line in $\C^d$, hence $n = 1$.  
Suppose that $\sV$ contains three $\Q$-linearly independent points of $\fL^d$ $-$then $\sV$ contains a nonzero point of $\Q^d$.
\vspace{0.2cm}

[In fact, if $\sV \cap \Q^d = \{0\}$, then
\[
\dim_{\Q} (\sV\hsx \cap \hsx \fL^d) 
\ \leq \ 
n(n+1) 
\ = \ 
1 (1 + 1) 
\ = \ 
2.
\]
But the assumption implies that
\[
\dim_{\Q} (\sV\hsx \cap \hsx \fL^d) 
\ \geq \ 
3.
\]
Therefore $\sV \cap \Q^d \neq \{0\}$ \hsx.]
\end{x}
\vspace{0.3cm}

It is conjectured that $n(n+1)$ in \#5 can be replaced by $n(n+1)/2$ but this remains to be seen.
\vspace{0.3cm}

\begin{x}{\small\bf EXAMPLE} \ 
Fix nonzero $\Q$-linearly independent elements $\lambda_1, \ldots, \lambda_{n+1}$ of $\fL$ and define $\sV$ by the equations
\[
\lambda_1 z_1 + \cdots + \lambda_{n+1} z_{n+1} 
\ = \ 
0, 
\quad 
z_{n+2} \ = \ \cdots \ = \ z_d \ = \ 0.
\]
Then $\sV \cap \Q^d = \{0\}$ and $\sV\hsx \cap \hsx \fL^d$ contains the $n(n+1)/2$ points
\[
w_{i j} 
\ = \ 
(w_{i j 1}, \ldots, w_{i j d}) \in \C^d \qquad (1 \leq i < j \leq d),
\]
where
\[
w_{i j k} \ = \ \lambda_j \ (k = i), 
\quad 
w_{i j k} \ = \ -\lambda_i \ (k = j), 
\]
and $w_{i j k} = 0$ otherwise $(1 \leq k \leq d)$.  
And these points are $\Q$-linearly independent,
hence
\[
\dim_{\Q} (\sV\hsx \cap \hsx \fL^d) 
\ \geq \ 
n (n + 1) / 2.
\]
\end{x}
\vspace{0.3cm}

\begin{x}{\small\bf RAPPEL} \ 
Let \mX be a vector space, $S \subset X$ a nonempty subset $-$then the 
\un{span}
\index{span} 
$\langle S \rangle$
\index{$\langle S \rangle$} 
of \mS is the intersection of all subspaces containing \mS or still, the set of all finite linear combinations of the elements of \mS.
\end{x}
\vspace{0.3cm}

\begin{x}{\small\bf NOTATION} \ 
Given a $\C$-vector subspace  $\sV \subset \C^d$, put
\[
t 
\ = \ 
\dim_{\C} \hsx \langle \sV  \hsx \cap \hsx \Qbar^{\hsx d} \rangle,
\]
the dimension of the $\C$-vector space spanned by $\sV \hsx \cap \hsx \Qbar^{\hsx d}$.
\end{x}
\vspace{0.3cm}

\begin{x}{\small\bf \un{N.B.}} \ 
For the record,
\[
0 \ \leq \ t \ \leq \ n \ < \ d,
\]
it being assumed that $\sV \neq \C^d$.
\end{x}
\vspace{0.3cm}

\begin{x}{\small\bf THEOREM} \ 
Let $\sV \subset \C^d$ be a $\C$-vector subspace such that $\sV \cap \Q^d = \{0\}$ $-$then
\allowdisplaybreaks
\begin{align*}
\dim_{\Q} (\sV\hsx \cap \hsx \fL^d) \ 
&\leq \ 
d (n - t) 
\\[12pt]
&\leq \ 
d (d - 1 - t), 
\end{align*}
where
\[
n \ = \ \dim_{\C} (\sV).
\]
\end{x}
\vspace{0.3cm}

\begin{x}{\small\bf REMARK} \ 
Sometimes this estimate is better than the one provided by \#5 but it can also be worse.
\vspace{0.2cm}

\qquad \textbullet \quad
Suppose that
\[
n \ = \ \dim_{\C} (\sV) \ = \ d - 1, 
\quad 
t = n.
\]
Then
\allowdisplaybreaks
\begin{align*}
d(n - t) \ 
&=\ 
d(d - 1 - t) 
\\[12pt]
&=\ 
d(d - 1 - (d - 1))
\\[12pt]
&=\ 
0
\\[12pt]
\implies 
&\dim_{\Q} (\sV\hsx \cap \hsx \fL^d) = 0 
\\[12pt]
\implies 
&\sV\hsx \cap \hsx \fL^d = \{0\}
\end{align*}
in accordance with expectation ($\sV$ being rational over $\Qbar$).  
As for \#5, it just gives
\[
\dim_{\Q} (\sV\hsx \cap \hsx \fL^d) \ \leq \ (d-1) (d).
\]
\vspace{0.2cm}

\qquad \textbullet \quad 
Suppose that
\[
n = \dim_{\C} (\sV) = 1, \quad t = 0.
\]
Then
\[
d(n - t) \ = \ d(1-0) \ = \ d,
\]
whereas
\[
n(n+1) \ = \ 2
\]
which is less than $d$ if $d \geq 3$.

\end{x}
\vspace{0.3cm}

\begin{x}{\small\bf EXAMPLE} \ 
Let $\sV \subset \C^3$ be the hyperplane defined by the equation
\[
\sqrt{2} \hsx z_1 + e z_2 + z_3 \ = \ 0.
\]
Then $\sqrt{2}$, $e$, 1 are $\Q$-linearly independent.  
To check this, consider a rational dependence relation
\[
q_1 \sqrt{2} + q_2 e + q_3 \ = \ 0.
\]
\vspace{0.02cm}

\allowdisplaybreaks
\begin{align*}
\text{\un{Case 1:}} \quad q_1 = 0 
&\implies q_2 e + q_3 = 0 \hspace{6cm}
\\[12pt]
&\implies q_2 = 0, \ q_3 = 0.
\\[18pt]
\text{\un{Case 2:}} \quad q_1 \neq 0 
&\implies \sqrt{2} + \ds\frac{q_2}{q_1} e + \ds\frac{q_3}{q_1} \ = \ 0
\\[12pt]
&\implies
\frac{q_2}{q_1} e = -\sqrt{2} - \frac{q_3}{q_1}
\\[12pt]
&\implies
e = \frac{q_1}{q_2} \bigg(-\sqrt{2} - \frac{q_3}{q_1}\bigg).
\end{align*}
I.e.: $e$ is algebraic which it isn't.  
Consequently, $\sV \cap \Q^3 = \{0\}$.  
Since here
\[
d = 3, \ n = 2, \ t = 1, 
\]
it therefore follows from \#11 that
\[
\dim_{\Q} (\sV \hsx \cap \hsx \fL^3) \ \leq \ 3(2 - 1) \ = \ 3.
\]

[Note: \ 
There are three possibilities for $t$: 0, 1, 2. 
But
\[
(1, 0, -\sqrt{2}) \hsx \in \hsx \sV \hsx \cap \hsx \Qbar^{\hsx 3}
\]
which implies that $t \geq 1$.  
And $t = 2$ is impossible ($\sV$ is not rational over $\Qbar$), thus $t = 1$ \hsx .]
\end{x}
\vspace{0.3cm}

It has been observed above that \hsy \#1(i) \hsy is a particular instance of \hsx \#11 \hsx (cf. \#12 (first \textbullet)).  
To repeat:
\vspace{0.3cm}

\begin{spacing}{1.55}
\begin{x}{\small\bf THEOREM} \ 
Let $\sV \subset \C^d$ be a $\C$-vector subspace rational over \hsx $\Qbar$  \hsx with   \hsx$\sV \hsx  \cap \hsx \Q^d = \{0\}$ 
$-$then 
$\sV\hsx \cap \hsx \fL^d = \{0\}$.
\end{x}
\end{spacing}
\vspace{0.3cm}

\begin{spacing}{1.55}
\begin{x}{\small\bf APPLICATION} \ 
Here is one version of Gelfond-Schneider: \ 
Let $\lambda_1 \in \fL$, $\lambda_2 \in \fL$, let $\beta \in \Qbar$, $\beta \notin \Q$, 
and suppose that $\lambda_2 = \beta \lambda_1$ $-$then the claim is that $\lambda_1 = \lambda_2 = 0$.  
To establish this, work in $\C^2$ and let $\sV \subset \C^2$ be the complex line $\C(1, \beta)$ $-$then 
$\sV \  \cap \  \Q^2 = \{0\}$ ($(z, z\beta) = (q_1, q_2)$ $\implies$ 
$z = q_1$ $\implies$ $q_1 \beta = q_2$ $\implies$ 
$\beta = q_2 /  q_1$ if $q_1 \neq 0$).  
Moreover $\sV$ is rational over $\Q$ ($\sV$ being defined by the equation 
$z_2 = \beta z_1$).  
The assumptions of \#14 are therefore satisfied, hence $\sV \cap \fL^2 = \{0\}$.  
But $(\lambda_1, \lambda_2) \in \sV \hsx \cap \hsx \fL^2$, thus $\lambda_1 = \lambda_2 = 0$, as contended.
\end{x}
\end{spacing}
\vspace{0.3cm}

\begin{x}{\small\bf APPLICATION} \ 
Let $\beta_1 \neq 0, \ldots, \beta_d \neq 0$ be algebraic numbers.  
Denote by $\sV \subset \C^d$ the hyperplane defined by the equation
\[
\beta_1 z_1 + \cdots + \beta_d z_d \ = \ 0.
\]
Then $\sV$ is rational over $\Qbar$. \   
Assume: \ $\sV \cap \Q^d = \{0\}$  \ $-$then $\sV\hsx \cap \hsx \fL^d \ =  \ \{0\}$ \quad (cf. \#14).  
\\[2pt]
Next $\beta_1, \ldots, \beta_d$ are $\Q$-linearly independent:
\[
q_1 \beta_1 + \cdots + q_d \beta_d = 0 
\implies 
(q_1, \ldots, q_d) \in \sV \cap \Q^d = \{0\}.
\]
To exploit this, take nonzero $\lambda_1 \in \fL, \ldots, \lambda_n \in \fL$ and consider
\[
\beta_1 \lambda_1 + \cdots + \beta_d \lambda_d,
\]
which we claim is nonzero.  
For otherwise
\[
(\lambda_1, \ldots, \lambda_d) \in \sV\hsx \cap \hsx \fL^d \ = \ \{0\}.
\]
Now quote \S32, \#3(ii) to see that this setup implies homogeneous Baker.  
\vspace{0.2cm}

[Note: \ 
In \S32, \#3(ii), the supposition is that $\beta_1, \ldots, \beta_d$ are $\Q$-linearly independent (replace $n$ by $d$).  
This implies that $\sV \cap \Q^d = \{0\}$.  
Proof:
\[
(z_1, \ldots, z_d) \ = \ (q_1, \ldots, q_d) \in \sV \cap \Q^d
\]
\qquad\qquad $\implies$
\[
\beta_1 z_1 + \cdots + \beta_d z_d \ = \ 0
\]
\qquad\qquad $\implies$
\[
\beta_1 q_1 + \cdots + \beta_d q_d \ = \ 0 \hsx .]
\]
\end{x}
\vspace{0.3cm}


%% file: _39_vector_spaces_LG.tex
\chapter{
$\boldsymbol{\S}$\textbf{39}.\quad  VECTOR SPACES: \ $\fL_G$}
\setlength\parindent{2em}
\setcounter{theoremn}{0}
\renewcommand{\thepage}{\S39-\arabic{page}}

\ \indent 
It will be useful to generalize the considerations in \S38 as this provides a convenient forum for certain important applications.
\vspace{0.5cm}

\begin{x}{\small\bf NOTATION} \ 
Let $d_0 \geq 0$, $d_1 \geq 1$ be integers and let $d = d_0 + d_1$.  
Put
\[
\begin{cases}
\ G_0 \ = \ \C \times \cdots \times \C \qquad \hspace{0.53cm} \text{($d_0$ factors)}\\[8pt]
\ G_1 \ = \ \C^\times \times \cdots \times \C^\times \qquad \text{($d_1$ factors)}
\end{cases}
\]
and set
\[
G \ = \ G_0 \times G_1.
\]
\end{x}
\vspace{0.3cm}

\begin{x}{\small\bf NOTATION} \ 
\[
\fL_G \ = \ \Qbar^{\hsx d_0} \times \fL^{d_1}.
\]
\vspace{0.2cm}

[Note: \ 
Accordingly an element $\fL_G$ is a $d_0 + d_1$ tuple
\[
(\beta_1, \ldots, \beta_{d_0}, \lambda_1, \ldots, \lambda_{d_1}),
\]
where $\beta_1, \ldots, \beta_{d_0}$ are algebraic numbers i.e., are in $\Qbar$ and $\lambda_1, \ldots, \lambda_{d_1}$ are logarithms of algebraic numbers, i.e., are in $\fL$.]
\end{x}
\vspace{0.3cm}

\begin{x}{\small\bf \un{N.B.}}  \ 
The choice $d_0 = 0$ puts us back into the setting of \S38.
\end{x}
\vspace{0.3cm}

\begin{x}{\small\bf LEMMA} \ 
$\fL_G$ is a $\Q$-vector subspace of $\C^d$.
\end{x}
\vspace{0.3cm}

\begin{x}{\small\bf LEMMA} \ 
Let $\sV \subset \C^d$ be a $\C$-vector subspace.
\vspace{0.5cm}

\qquad \textbullet \quad
If $\sV \cap (\{0\} \times \Q^{d_1}) \neq \{0\}$, then
\[
\dim_{\Q} (\sV\hsx \cap \hsx \fL_G) \ = \ \infty.
\]

[Take \ 
\[
\bq \ = \ (0, \ldots, 0, q_1, \ldots, q_{d_1}) \ \neq \ \bzero
\]
in $\sV \cap (\{0\} \times \Q^{d_1})$ $-$then $\forall \ \lambda \in \fL$, 
\[
(0, \ldots, 0, q_1 \lambda, \ldots, q_{d_1} \lambda) \in \sV\hsx \cap \hsx \fL_G 
\implies 
\dim_{\Q} (\sV\hsx \cap \hsx \fL_G) \ = \ \infty.]
\]
\\
[-20pt]

\qquad \textbullet \quad
If $\sV \cap (\Qbar^{\hsx d_0}\times \{0\}) \neq \{0\}$, then
\[
\dim_{\Q} (\sV\hsx \cap \hsx \fL_G) \ = \ \infty.
\]
\\[-20pt]

[Take \
\[
\bbeta \ = \  (\beta_1, \ldots, \beta_{d_0}, 0, \ldots, 0) \ \neq \ \boldsymbol{0}
\]
in $\sV \cap (\Qbar^{\hsx d_0} \times \{0\})$ $-$then $\forall \ \gamma \in \Qbar$, 
\[
(\beta_1 \gamma, \ldots,  \beta_{d_0}\gamma, 0, \ldots, 0) \in \sV\hsx \cap \hsx \fL_G 
\implies 
\dim_{\Q} (\sV\hsx \cap \hsx \fL_G) \ = \ \infty. \hsx]
\]

\end{x}
\vspace{0.3cm}

\begin{x}{\small\bf SCHOLIUM} \ 
If
\[
\dim_{\Q} (\sV\hsx \cap \hsx \fL_G) \ < \ \infty,
\]
then
\[
\sV \cap (\{0\} \times \Q^{d_1}) \ = \ \{0\}
\quad \text{and} \quad
\sV \cap (\Qbar^{\hsx d_0}\times \{0\}) = \{0\}.
\]
\end{x}
\vspace{0.3cm}

\begin{x}{\small\bf DEFINITION} \ 
The relations
\[
\sV \cap (\{0\} \times \Q^{d_1}) \ = \ \{0\}
\quad \text{and} \quad
\sV \cap (\Qbar^{\hsx d_0}\times \{0\}) = \{0\}
\]
are the 
\un{canonical conditions}.
\index{canonical conditions}
\end{x}
\vspace{0.3cm}


\begin{x}{\small\bf THEOREM} \ 
Let $\sV \subset \C^d$ be a $\C$-vector subspace for which the canonical conditions are in force $-$then
\[
\dim_{\Q} (\sV\hsx \cap \hsx \fL_G) \ < \ \infty
\]
and, in fact, 
\[
\dim_{\Q} (\sV\hsx \cap \hsx \fL_G) \ \leq \ d_1(n - t).
\]

[Note: \ 
As in \S38, 
\[
n \ = \ \dim_{\Q} (\sV) 
\quad \text{and} \quad
t \ = \ \dim_{\C}  \langle \sV \cap \Qbar^{\hsx d} \rangle \hsx .]
\]
\end{x}
\vspace{0.3cm}

\begin{x}{\small\bf REMARK} \ 
Taking $d_0 = 0$ recovers \S38, 
\#11.  
As for the proof, it will be omitted since it depends on the so-called ``linear subgroup theorem'' which we shall not stop to formulate.]
\end{x}
\vspace{0.3cm}

\begin{x}{\small\bf APPLICATION} \ 
Homogeneous Baker
\index{Homogeneous Baker} 
is the assertion that if $\lambda_1 \in \fL, \ldots, \lambda_d \in \fL$ are nonzero and $\Q$-linearly independent, then 
$\lambda_1, \ldots, \lambda_d$ are $\Qbar$-linearly independent.
\vspace{0.2cm}

[
Suppose that $\lambda_1, \ldots, \lambda_d$ are $\Qbar$-linearly dependent, say
\[
\beta_1 \lambda_1 + \cdots + \beta_{d-1} \lambda_{d-1} \ = \ \lambda_d, 
\]
where $\beta_1, \ldots, \beta_{d-1}$ are algebraic.  It can be assumed in addtion that $\lambda_1, \ldots, \lambda_{d-1}$ are 
$\Qbar$-linearly independent.  
Take now for $\sV$ the hyperplane in $\C^d$ defined by the equation
\[
\lambda_1 z_1 + \cdots + \lambda_{d-1} z_{d-1} = z_d.
\]
Explicate the parameters: \ 
$d_0 = n = d-1$, $d_1 = 1$ (so $d \equiv d_0 + d_1 = n + 1 = (d - 1) + 1 = d \ldots$), $t = 0$.  
The definitions imply that the canonical conditions are in 
force, thus by \#8, 
\[
\dim_{\Q} (\sV\hsx \cap \hsx \fL_G) 
\ \leq \ 
d_1 (n - t) \ = \ 1(d - 1 - 0) \ = \ d - 1.
\]
On the other hand,
\[
\sV\hsx \cap \hsx \fL_G
\ = \ 
\sV \cap (\Qbar^{\hsx d-1} \times \fL)
\]
contains $d$ \hsy $\Q$-linearly independent points $\zeta_1, \ldots, \zeta_d$, namely
\[
\zeta_i 
\ = \ 
(\delta_{i \hsy1}, \ldots, \delta_{i \hsy (d-1)}, \lambda_i) \qquad (1 \leq i \leq d-1)
\]
and
\[
\zeta_d 
\ = \ 
(\beta_1, \ldots, \beta_{d-1}, \lambda_d) \hsx .]
\]

[Note: \ 
Take a point in $\sV \cap \Qbar^{\hsx d}$, say $(\beta_1, \ldots, \beta_d)$, subject to
\[
\lambda_1 \beta_1 + \cdots + \lambda_{d-1} \beta_{d-1} \ = \ \beta_d.  
\]
Argue that necessarily $\beta_d = 0$ (cf. \#14), hence 
$\beta_1 = 0, \ldots, \beta_{d-1} = 0$ $(\lambda_1, \ldots, \lambda_{d-1}$ 
are $\Qbar$-linearly independent), hence 
$\sV \hsx \cap \hsx  \Qbar^{\hsx d} = \{0\}$, hence $t = 0$.]
\end{x}
\vspace{0.3cm}

\begin{x}{\small\bf APPLICATION} \ 
Inhomogeneous Baker 
is the assertion that if $\lambda_1 \in \fL, \ldots, \lambda_d \in \fL$ are nonzero and $\Q$-linearly independent, then 
$1, \lambda_1, \ldots, \lambda_d$ are $\Qbar$-linearly independent.
\vspace{0.2cm}

[
Suppose that $1, \lambda_1, \ldots, \lambda_d$ are $\Qbar$-linearly dependent, say
\[
\beta_0 + \beta_1 \lambda_1 + \cdots + \beta_{d-1} \lambda_{d-1} \ = \ \lambda_d,
\]
where $\beta_0, \beta_1, \ldots, \beta_{d-1}$ are algebraic.  
It can be assumed in addition that $\lambda_1 , \ldots, \lambda_{d}$ are $\Q$-linearly independent 
and  $1, \lambda_1, \ldots, \lambda_{d-1}$ are $\Qbar$-linearly independent.  
Take now
for $\sV$ the hyperplane in $\C^{d+1}$ defined by the equation
\[
z_0  + \lambda_1 z_1 + \cdots + \lambda_{d-1} z_{d-1} \ = \ z_d.
\]
Explicate the parameters: \ $d_0 = n = d$, $d_1 = 1$ (the role of $d$ in the theory is played in this situation by 
$d + 1$: $d_0 + d_1 = d+1$, $t \geq 1$ (since $(1, 0, \ldots, 0, 1) \in \sV$).  
The definitions imply that the canonical conditions are in force, thus by \#8
\[
\dim_{\Q} (\sV\hsx \cap \hsx \fL_G) 
\ \leq \ 
d_1 (n - t) 
\ = \ 
1(d-t) 
\ \leq \ 
d - 1.
\]
On the other hand, 
\[
\sV\hsx \cap \hsx \fL_G 
\ = \ 
\sV \hsx \cap \hsx \big(\Qbar^{\hsx d} \times \fL)
\]
contains $d$ \hsy $\Q$-linearly independent points $\zeta_1, \ldots, \zeta_d$, namely
\[
\zeta_i 
\ = \ 
(0, \delta_{i \hsy 1}, \ldots, \delta_{i \hsy (d-1)}, \lambda_i) \qquad (1 \leq i \leq d - 1)
\]
and 
\[
\zeta_d
\ = \ 
(\beta_0, \beta_1, \ldots, \beta_{d-1}, \lambda_d) \hsx .]
\]

[Note: \ 
\[
t \geq 1 
\implies
- t \leq -1 
\implies
d - t \leq d -1.
\]
Also, on general grounds, $\beta_0 = 0$ (cf. \#14).]
\end{x}
\vspace{0.3cm}

\begin{x}{\small\bf THEOREM} \ 
Let $\sV \subset \C^d$ be a $\C$-vector subspace rational over $\Qbar$ and for which the canonical conditions are in force $-$then 
$\sV\hsx \cap \hsx \fL_G = \{0\}$.
\vspace{0.2cm}

PROOF \ 
In \#8, take $t = n$ to get
\[
\dim_{\Q} \hsx (\sV\hsx \cap \hsx \fL_G) \  = \  \{0\}.
\]
\end{x}
\vspace{0.3cm}

\begin{x}{\small\bf APPLICATION} \ 
\vspace{0.2cm}

\qquad \textbullet \quad
If $\alpha$ is a nonzero algebraic number, then $e^\alpha$ is transcendental (cf. \S21, \#4).
\vspace{0.3cm}

\qquad \textbullet \ \textbullet \quad
If $\beta$ is an algebraic number such that $e^\beta$ is algebraic, then $\beta = 0$.
\vspace{0.3cm}

Claim: \ 
\textbullet \ \textbullet $\implies$ \textbullet \quad
For if $e^\alpha$ was not transcendental, then it would be algebraic, hence that $\alpha = 0$, contradiction.
\vspace{0.2cm}

To establish \textbullet \ \textbullet, take $d_0 = 1$, $d_1 = 1$, so that $d = 1 + 1 = 2$ and 
$\fL_G = \Qbar \times \fL$.    
The complex line $\sV = \C(1,1)$ in $\C^2$ is rational over $\Qbar$ and contains $(\beta,\beta) \in \fL_G$.  
Moreover it is clear that the canonical conditions hold.  
Therefore
\[
\sV\hsx \cap \hsx \fL_G \ = \ \{0\} \quad \text{(cf. \#12)} \implies \beta = 0.
\]
\end{x}
\vspace{0.3cm}

\begin{x}{\small\bf APPLICATION} \ 
Suppose given a relation
\[
\beta_0 \ + \beta_1 \lambda_1 + \cdots + \beta_d \lambda_d \ = \ 0,
\]
where $\beta_0, \beta_1, \ldots, \beta_d$ are algebraic and 
$\lambda_1 \in \fL, \ldots, \lambda_d \in \fL$ $-$then $\beta_0 = 0$.
\vspace{0.2cm}

[Argue by contradiction and assume that $\beta_0 \neq 0$ with $d$ minimal, thus $\beta_1, \ldots, \beta_d$ are $\Q$-linearly independent and $\lambda_1, \ldots, \lambda_d$ are $\Qbar$-linearly independent.  
Let $\sV \subset \C^{d+1}$ be the hyperplane defined by the equation
\[
\beta_0 z_0 + \beta_1 z_1 + \cdots + \beta_d z_d \ = \ 0.
\]
Then $\sV$ is rational over $\Qbar$ and the canonical conditions are satisfied.  But
\[
(1, \lambda_1, \ldots , \lambda_d) \in \sV
\]
and
\[
(1, \lambda_1, \ldots , \lambda_d)  \in \fL_G \ = \ \Qbar \times \fL^d \qquad (d_0 = 1, \ d_1 = d).
\]
Meanwhile
\[
\sV\hsx \cap \hsx \fL_G \ = \ \{0\} \qquad \text{(cf. \#12)}.]
\]
\end{x}
\vspace{0.3cm}


\begin{x}{\small\bf SCHOLIUM} \ 
Suppose given a relation
\[
\beta_1 \lambda_1 + \cdots + \beta_d \lambda_d \ = \ 0,
\]
where $\beta_1, \ldots, \beta_d$ are algebraic and $\lambda_1 \in \fL, \ldots, \lambda_d \in \fL$.
\vspace{0.3cm}

\qquad \textbullet \quad
If $(\beta_1, \ldots, \beta_d) \neq (0, \ldots, 0)$, then $\lambda_1, \ldots, \lambda_d$ are $\Q$-linearly dependent.
\vspace{0.3cm}

\qquad \textbullet \quad
If $(\lambda_1, \ldots, \lambda_d) \neq (0, \ldots, 0)$, then $\beta_1, \ldots, \beta_d$ are $\Q$-linearly dependent.

\vspace{0.2cm}

\end{x}
\vspace{0.3cm}

\begin{x}{\small\bf \un{N.B.}} \ 
Recall that every nonzero linear combination
\[
\beta_1 \lambda_1 + \cdots + \beta_d \lambda_d 
\]
is transcendental (cf. \S31, \#11).
\end{x}
\vspace{0.3cm}

\begin{x}{\small\bf LEMMA} \ 
Suppose that $\lambda_1, \ldots, \lambda_d$ are nonzero elements of $\fL$ and $\beta_1, \ldots, \beta_d$ are nonzero elements of 
$\Qbar$.  Assume: 
\[
\beta_1 \lambda_1 + \cdots + \beta_d \lambda_d \ = \ 0.
\]
Then there exist nonzero integers $k_1, \ldots, k_d$ such that
\[
k_1 \beta_1 + \cdots + k_d \beta_d \ = \ 0.
\]
\end{x}
\vspace{0.3cm}


%% file: _40_vector_spaces_Vmax_V_min.tex
\chapter{
$\boldsymbol{\S}$\textbf{40}.\quad  VECTOR SPACES: \ $V_{\max}, \ V_{\min}$}
\setlength\parindent{2em}
\setcounter{theoremn}{0}
\renewcommand{\thepage}{\S40-\arabic{page}}


\begin{x}{\small\bf CONSTRUCTION} \ 
Let $\sV \subset \C^d$ be a $\C$-vector subspace $-$then $\sV$ contains a unique maximal subspace $\sV_{\max}$ of the form 
$W_0 \times W_1$, where $W_0$ is a subspace of $\C^{d_0}$ rational over $\Qbar$ and 
$W_1$ is a subspace of $\C^{d_1}$ rational over $\Q$.
\end{x}
\vspace{0.3cm}

\begin{x}{\small\bf LEMMA} \ 
$W_0$ is the subspace of $\C^{d_0}$ spanned by 
\[
\sV \cap \big(\Qbar^{\hsx d_0} \times \{0\}\big)
\]
and $W_1$ is the subspace of $\C^{d_1}$ spanned by
\[
\sV \cap \big(\{0\} \times \Q^{d_1}\big).
\]
\end{x}
\vspace{0.3cm}

\begin{x}{\small\bf RAPPEL} \ 
(cf. \S39, \#7) \ The relations
\[
\sV \cap \big(\{0\} \times \Q^{d_1}\big) \ = \ \{0\} 
\quad \text{and} \quad 
\sV \cap \big(\Qbar^{\hsx d_0} \times \{0\}\big) \ = \ \{0\} 
\]
are the 
\un{canonical conditions}.
\index{canonical conditions}
\end{x}
\vspace{0.3cm}

\begin{x}{\small\bf \un{N.B.}} \ 
$\sV_{\max} = \{0\}$ iff the canonical conditions are in force.
\end{x}
\vspace{0.3cm}

\begin{x}{\small\bf THEOREM} \ 
Let $\sV \subset \C^d$ be a $\C$-vector subspace.  
Assume: \ $\sV$ is rational over $\Qbar$ $-$then
\[
\sV \hsx \cap \hsx \fL_G 
\ = \ 
\sV_{\max} \hsx \cap \hsx  \fL_G.
\]

PROOF \ 
Trivially, 
\[
\sV_{\max} \hsx \cap \hsx \fL_G \hsx \subset \hsx \sV \hsx \cap \hsx \fL_G.
\]
This said, if first the canonical conditions hold, then $\sV \hsx \cap \hsx \fL_G = 0$ (cf. \S39 \# 12).
But also $\sV_{\max} = \{0\}$ (cf. \#4), hence $\sV_{\max} \hsx \cap \hsx \fL_G = 0$.  
Proceeding in general, write
\[
\sV_{\max} \ = \ W_0 \times W_1,
\]
put
\[
d_0^{\hsy \prime} 
\ = \ 
\dim_{\C} \bigg(\frac{\C^d}{W_0}\bigg) , \quad
d_1^{\hsy \prime}
\ = \ 
\dim_{\C} \bigg(\frac{\C^d}{W_1}\bigg),
\]
and introduce
\[
\begin{cases}
\ G_0^{\hsy \prime} = \C \times \cdots \times \C \qquad \hspace{0.5cm} \text{($d_0^{\hsy \prime}$ factors)}\\[8pt]
\ G_1^{\hsy \prime} = \C^\times \times \cdots \times \C^\times \qquad \text{($d_1^{\hsy \prime}$ factors)}
\end{cases}
.
\]
Let $\C^{d_0} \ra \C^{d_0^{\hsy \prime}}$ be a surjective linear map, rational over $\Qbar$, with kernel $W_0$ and 
let $\C^{d_1} \ra \C^{d_1^{\hsy \prime}}$ be a surjective linear map, rational over $\Q$, with kernel $W_1$.  
Denote by  $\phi$ their product
\[
\C^{d_0} \times \C^{d_1} \ra \C^{d_0^{\hsy \prime}} \times \C^{d_1^{\hsy \prime}}.
\]
Then the kernel of $\phi$ is $\sV_{\max}$ and $\phi(\fL_G) = \fL_{G^{\prime}}$.  
Moreover the canonical conditions hold for the subspace $\sV^\prime = \phi(\sV)$ of 
$\C^{d_0^{\hsy \prime}} \times \C^{d_1^{\hsy \prime}}$, hence $\sV^\prime \hsx \cap \hsx \fL_{G^\prime} = \{0\}$.  
Therefore
\[
\sV \cap \fL_G \hsx \subset \hsx \phi^{-1} (\sV^\prime \hsx \cap \hsx \fL_{G^\prime}) 
\ = \ \Ker \phi 
\ = \ 
\sV_{\max}
\]
\qquad\qquad $\implies$
\[
\sV \hsx \cap \hsx \fL_G \hsx \subset \hsx \sV_{\max} \hsx \cap \hsx \fL_G.
\]
\end{x}
\vspace{0.3cm}

\begin{x}{\small\bf CONSTRUCTION} \ 
Let $\sV \subset \C^d$ be a $\C$-vector subspace $-$then $\sV$ is contained in a unique miminal subspace $\sV_{\min}$ of the form 
$W_0 \times W_1$, where $W_0$ is a subspace of $\C^{d_0}$ rational over $\Qbar$ and $W_1$ is a subspace of 
$\C^{d_1}$ rational over $\Q$.
\end{x}
\vspace{0.3cm}


\begin{x}{\small\bf LEMMA} \ 
$W_0$ is the intersection of all hyperplanes of $\C^{d_0}$ rational over $\Qbar$ which contain the projection of $\sV$ onto $\C^{d_0}$ and $W_1$ is the intersection of all hyperplanes of $\C^{d_1}$ rational over $\Q$ which contain the projection of $\sV$ onto $\C^{d_1}$.
\end{x}
\vspace{0.3cm}

\begin{x}{\small\bf \un{N.B.}} \ 
$\sV_{\min} = \C^d$ means that $W_0 = \C^{d_0}$ and $W_1 = \C^{d_1}$.
\end{x}
\vspace{0.3cm}

\[
\textbf{APPENDIX}
\]
\vspace{0.5cm}

{\small\bf FACT} \ 
Let $\sV \subset \C^d$ be a $\C$-vector subspace.  
Assume: \ The canonical conditions are in force $-$then there exists a hyperplane $\sH \subset \C^d$ containing $\sV$ and for which the canonical conditions are also in force.

%% file: _41.tex
\chapter{
$\boldsymbol{\S}$\textbf{41}.\quad  EXPONENTIALS (6 or 5)}
\setlength\parindent{2em}
\setcounter{theoremn}{0}
\renewcommand{\thepage}{\S41-\arabic{page}}

\ \indent 

Specialized to the case when $m = 2$, $n = 3$, the six exponentials theorem is the following statement (cf. \S36, \#1):
\vspace{0.3cm}

\begin{x}{\small\bf THEOREM} \ 
Let $\{x_1, x_2\}$  and $\{y_1, y_2, y_3\}$  be two $\Q$-linearly independent sets of complex numbers $-$then at least one of the six numbers 
\[
e^{x_1 y_1}, \ 
e^{x_1 y_2}, \ 
e^{x_1 y_3}, \ 
e^{x_2 y_1}, \ 
e^{x_2 y_2}, \ 
e^{x_2 y_3}
\]
is transcendental.
\vspace{0.2cm}

\begin{spacing}{1.60}
PROOF \ 
To arrive at a contradiction, assume that the six numbers $x_i y_j$ $(i = 1, 2, \ j = 1, 2, 3)$ all belong to $\fL$ (the vectors in a linearly independent set are nonzero, thus $x_i \neq 0$ $(i = 1, 2)$, $y_j \neq 0$ $(j = 1, 2, 3)$, so $x_i y_j \neq 0$).  
Work in $\C^2$ and take for $\sV$ the complex line $\C \bx = \C\{x_1, x_2\}$ $-$then $\sV \hsx \cap \hsx \Q^2 = \{0\}$.  
For suppose that 
\end{spacing}
\[
z \bx \ = \ (z x_1, z x_2) \in \sV \cap \Q^2 \qquad (z \in \C, \ z \neq 0).
\]
Then
\[
\begin{cases}
\ z x_1 = q_1\\
\ z x_2 = q_2 
\end{cases}
\qquad (q_1, q_2 \in \Q)
\]
and the claim is that $q_1 = 0$, $q_2 = 0$.  
Consider the four possibilities.
\vspace{0.2cm}

\qquad \textbullet \quad
$q_1 \neq 0$, $q_2 \neq 0$ $\implies$
\[
\frac{1}{z} \ = \ \frac{x_1}{q_1}, \quad 
\frac{1}{z} \ = \ \frac{x_2}{q_2}
\]
\allowdisplaybreaks
\begin{align*}
&\implies 
q_2 x_1 = q_1 x_2 
\\[12pt]
&\implies 
q_2 x_1 - q_1 x_2 = 0
\\[12pt]
&\implies
q_1 = 0, \ q_2 = 0,
\end{align*}
$\{x_1, x_2\}$ being $\Q$-linearly independent.  
\vspace{0.3cm}

\qquad \textbullet \quad
$q_1 \neq 0$, $q_2 = 0 \implies z x_2 = 0 \implies x_2 = 0$.
\vspace{0.2cm}

\qquad \textbullet \quad
$q_1 = 0$, $q_2 \neq 0 \implies z x_1 = 0 \implies x_1 = 0$.
\\[5pt]
Therefore these three possibilities are untenable, leaving \ $q_1 \hsx = \hsx 0$, $q_2 \hsx = \hsx 0$, as claimed.  
\\[5pt]
Next, $\sV \hsx \cap \hsx \fL^2$ contains the points
\[
y_1 \bx, \ 
y_2 \bx, \ 
y_3 \bx
\]
which are $\Q$-linearly independent.  
To see this, consider a rational dependence relation
\[
q_1 y_1 \bx + q_2 y_2 \bx + q_3 y_3 \bx = \bzero,
\]
i.e., 
\[
\begin{cases}
\ q_1 x_1 y_1 + q_2 x_1 y_2  + q_3 x_1 y_3  = 0 \\[8pt]
\ q_1 x_2 y_1 + q_2 x_2 y_2  + q_3 x_2 y_3  = 0
\end{cases}
.
\]
Dividing the first of these relations by $x_1 \neq 0$ (or the second of these relations by $x_2 \neq 0)$ gives
\[
q_1 y_1 + q_2 y_2 + q_3 y_3 \ = \ 0
\]
\qquad\qquad $\implies$
\[
q_1 = 0, \ 
q_2 = 0, \ 
q_3 = 0,
\]
$\{y_1, y_2, y_3\}$  being $\Q$-linearly independent.  
Therefore
\[
3 \ \leq \ \dim_{\Q} (\sV \hsx \cap \hsx \fL^2).
\]
On the other hand (cf. \S38, \#5), 
\[
\dim_{\Q} (\sV \hsx \cap \hsx \fL^2) \ \leq \ 1(1 + 1) = 2.  
\]
Contradiction.
\end{x}
\vspace{0.3cm}

The next result is known as the five exponentials theorem.
\vspace{0.3cm}

\begin{x}{\small\bf THEOREM} \ 
Let $\{x_1, x_2\}$  and $\{y_1, y_2\}$  be two $\Q$-linearly independent sets of complex numbers.  
Let further $\gamma$ be a nonzero algebraic number $-$then at least one of the five numbers 
\[
e^{x_1 y_1}, \ 
e^{x_1 y_2}, \ 
e^{x_2 y_1}, \ 
e^{x_2 y_2}, \ 
e^{\gamma x_1 / x_2}
\]
is transcendental.
\vspace{0.2cm}

PROOF \ 
With \S39, \#8 in mind, take $d_0 = 1$, $d_1 = 2$ $(\implies d = 3)$ and let $\sV$ be the hyperplane in $\C^3$ defined by the equation
\[
\gamma x_1 z_1 - x_2 z_2 + x_1 z_3 \ = \ 0 \qquad (\implies n = 2).
\]
Note that
\[
(1, 0, -\gamma) \in \sV \hsx \cap \hsx \Qbar^{\hsx 3},
\]
\begin{spacing}{1.6}
\noindent hence $t \geq 1$.  
If both $x_1$, $x_2$ are algebraic, then $\gamma x_1 / x_2 \neq 0$ is algebraic, so $e^{\gamma x_1 / x_2}$ is transcendental 
(cf. \S39, \#13).  
It can therefore be assumed that either $x_1$ or $x_2$ is transcendental, thus $\sV$ is not rational over $\Qbar$, thus 
$t \neq 2 \implies t = 1$.  
Moving on, since $x_1$, $x_2$ are $\Q$-linearly independent and $\gamma \neq 0$, it follows that the
canonical conditions are in force.  
Consequently
\end{spacing}
\[
\dim_{\Q} (\sV \hsx \cap \hsx  \fL_G) \ \leq \ d_1(n - t) \ = \ 2(2 - 1) \ = \ 2.
\]
On the other hand, $\sV$ contains the $\Q$-linearly independent points
\[
(1, \gamma x_1 / x_2, 0), \ 
(0, x_1 y_1, x_2 y_1), \ 
(0, x_1 y_2, x_2 y_2),
\]
so at least one of these does not belong to
\[
\fL_G \ = \ \Qbar \times \fL^2 \ = \ \Qbar \times \fL \times \fL.
\]
E.g.: \ Suppose that
\[
(0, x_1 y_1, x_2 y_1) \hsx \notin \hsx \Qbar \times \fL \times \fL.
\]
Then
\[
x_1 y_1 \notin \fL
\quad \text{or} \quad 
x_2 y_1 \notin \fL
\quad \text{(or both)}
\]
\qquad\qquad $\implies$
\[
e^{x_1 y_1} \ \text{transcendental or} \ 
e^{x_2 y_1} \ \text{transcendental (or both).} \ 
\]
\end{x}
\vspace{0.3cm}

\begin{x}{\small\bf EXAMPLE} \ 
Suppose that $\lambda_1 \in \fL$, $\lambda_2 \in \fL$.  
Assume: \ $\{\lambda_1, \lambda_2\}$ is $\Q$-linearly independent.  
Let $w \in \C$ $(w \notin \Q)$ and let $\beta \in \Qbar$ $(\beta \neq 0)$ $-$then at least one of the three numbers
\[
e^{w \lambda_1}, \ 
e^{w \lambda_2}, \ 
e^{\beta w}
\]
is transcendental.
\vspace{0.2cm}

[In \#2, take $x_1 = w$ $(\notin \Q)$, $x_2 = 1$, $y_1 = \lambda_1$, $y_2 = \lambda_2$ $-$then at least one of 
\[
e^{w \lambda_1}, \ 
e^{w \lambda_2}, \ 
e^{\lambda_1}, \ 
e^{\lambda_2}, \ 
e^{\beta w}
\]
is transcendental or still, at least one of 
\[
e^{w \lambda_1}, \ 
e^{w \lambda_2}, \ 
e^{\beta w}
\]
is transcendental.]
\vspace{0.2cm}

[Note: \ 
Put
\[
\begin{cases}
\ \alpha_1 = e^{\lambda_1} \\
\ \alpha_2 = e^{\lambda_2}
\end{cases}
.
\]
Then at least one of 
\[
\alpha_1^w, \ \alpha_2^w, \ e^{\beta w}
\]
is transcendental.]
\end{x}
\vspace{0.3cm}

\begin{x}{\small\bf EXAMPLE} \ 
Fix $\lambda \neq 0$ in $\fL$.  
Let $w \in \C$ $(w \notin \Q)$ and let $\beta \in \Qbar$ $(\beta \neq 0)$ $-$then at least one of the three numbers
\[
e^{w^2 \lambda}, \ 
e^{w \lambda}, \ 
e^{\beta w}
\]
is transcendental.
\vspace{0.2cm}

[In \#2, take $x_1 = w  \ (\notin \Q)$, $x_2 = 1$, $y_1 = w \lambda$, $y_2 = \lambda$ $-$then at least one of 
\[
e^{w^2 \lambda}, \ 
e^{w \lambda}, \ 
e^{w \lambda}, \ 
e^{\lambda}, \ 
e^{\beta w}
\]
is transcendental or still, at least one of 
\[
e^{w^2 \lambda}, \ 
e^{w \lambda}, \ 
e^{\beta w}
\]
is transcendental.]
\vspace{0.2cm}

[Note: \
Put $\alpha = e^\lambda$ $-$then at least one of 
\[
\alpha^{w^2}, \ 
\alpha^{w}, \ 
e^{\beta w}
\]
is transcendental.]
\end{x}
\vspace{0.3cm}

\begin{x}{\small\bf EXAMPLE} \ 
Let $\lambda_0 \in \fL$ $(\lambda_0 \neq 0)$, $\lambda_1 \in \fL$, $\lambda_2 \in \fL$, $\beta \in \Qbar$ $(\beta \neq 0)$, 
$\gamma = \ds\frac{1}{\beta}$.
Assume: \ $\{\lambda_1, \lambda_2\}$ is $\Q$-linearly independent $-$then at least one of the two numbers
\[
e^{\beta \lambda_0 \lambda_1}, \ e^{\beta \lambda_0 \lambda_2}
\]
is transcendental.
\vspace{0.2cm}

[In \#2, take $x_1 = \lambda_0 \beta \ (\notin \Q)$, $x_2 = 1$, $y_1 = \lambda_1$, $y_2 = \lambda_2$, hence at least one of
\[
e^{\beta \lambda_0   \lambda_1}, \ 
e^{\beta \lambda_0   \lambda_2}, \ 
e^{\lambda_1}, \ 
e^{\lambda_2}, \ 
e^{\frac{1}{\beta}\lambda_0 \beta} \ = \ e^{\lambda_0}
\]
is transcendental or still, at least one of
\[
e^{\beta \lambda_0 \lambda_1}, \ 
e^{\beta \lambda_0 \lambda_2}
\]
is transcendental.]
\vspace{0.2cm}

[Note: \
$\lambda_0 \beta$ is not rational (for if it were, then $\lambda_0$ would be algebraic whereas it is transcendental).]
\end{x}
\vspace{0.3cm}

\begin{x}{\small\bf EXAMPLE} \ 
Let $\lambda_0$, $\lambda_1$ be nonzero elements of $\fL$ and let $\beta \in \Qbar$ $(\beta \neq 0)$ $-$then at least one of the two numbers 
\[
e^{\beta \lambda_0 \lambda_1}, \ 
e^{(\beta \lambda_0)^2 \lambda_1}
\]
is transcendental.
\vspace{0.2cm}

[To illustrate, take $\beta = 1$, $\lambda_0 = \elln(2)$, $\lambda_1 = \elln(2)$ $-$then at least one of
\[
2^{\elln(2)}, \ 
2^{(\elln(2))^2}
\]
is transcendental.]
\end{x}
\vspace{0.3cm}

\begin{x}{\small\bf REMARK} \ 
Is it true that
\[
\text{five exponentials $\implies$ six exponentials?}
\]
In the literature, it is asserted that this is the case but no proof has been offered.
\vspace{0.2cm}


[To see the difficulty, in \#2, take $\gamma = 1$, and consider
\[
\begin{cases}
\ e^{x_1 y_1}, \ e^{x_1 y_2}, \ e^{x_2 y_1}, \ e^{x_2 y_2}, \ e^{x_1 / x_2} \\[8pt]
\ e^{x_1 y_3}, \ e^{x_2 y_3}, \ e^{x_1 y_1}, \ e^{x_2 y_1}, \ e^{x_1 / x_2} 
\end{cases}
.
\]
If $\ds{e^{x_1 / x_2}}$ is algebraic, then we are done since one of the exponentials in the first row preceding $\ds{e^{x_1 / x_2}}$ or in the second row preceding $\ds{e^{x_1 / x_2}}$ must be transcendental.  
However, if $\ds{e^{x_1 / x_2}}$ is transcendental, then it is conceivable that the first four exponentials in both rows are algebraic 
\ldots \ .]
\end{x}
\vspace{0.3cm}


%% file: _42_sharp_six_exp_theorem.tex
\chapter{
$\boldsymbol{\S}$\textbf{42}.\quad  SHARP SIX EXPONENTIALS THEOREM}
\setlength\parindent{2em}
\setcounter{theoremn}{0}
\renewcommand{\thepage}{\S42-\arabic{page}}

\ \indent 
This is the following statement.
\vspace{0.3cm}

\begin{x}{\small\bf THEOREM} \ 
Let $\{x_1, x_2\}$ and $\{y_1, y_2, y_3\}$ be two $\Q$-linearly independent sets of complex numbers.  
Let further $\beta_{i j}$ $(i = 1, 2, \ j = 1, 2, 3)$ be algebraic numbers.  \\
Assume: \ The six numbers
\[
e^{x_i y_j - \beta_{i j}}
\]
are algebraic, hence that the $\lambda_{i j} = x_i y_j - \beta_{i j}$ are in $\fL$ $-$then
\[
x_i y_j \ = \ \beta_{i j} \qquad (i = 1, 2, \ j = 1, 2, 3).
\]
\\[-1.1cm]

PROOF \ 
With \S39, \#8 in mind, take $d_0 = 2$, $d_1 = 2$ $(\implies d = 4)$ and let $\sV \subset \C^4$ be the hyperplane defined by the equation 
\[
x_2 (z_1 + z_3) \ = \ x_1 (z_2+ z_4) \qquad (\implies n = 3).
\]
Note that
\[
\begin{cases}
\ (1, 0, -1, 0) \in \sV \hsx \cap \hsx \Qbar^{\hsx 4}\\[8pt]
\ (0, -1, 0, 1) \in \sV \hsx \cap \hsx \Qbar^{\hsx 4}
\end{cases}
\implies t \geq 2.
\]
Note in addition that for $j = 1, 2, 3$, 
\[
\eta_j \ \equiv \ 
\big(\beta_{1 j}, \beta_{2 j}, \lambda_{1 j}, \lambda_{2 j}\big) \hsx \in \hsx \sV \hsx \cap \hsx  \fL_G 
\ = \ 
\sV \hsx \cap \hsx \big(\Qbar^{\hsx 2} \times \fL^2\big).
\]
\ifdmcx

\qquad \textbullet \quad The $\eta_j$ are $\Q$-linearly independent (see below).\\[2pt]

\qquad \textbullet \quad The canonical conditions are not satisfied (see below).\\[2pt]
\fi
\ifdmc
\qquad \textbullet \quad The $\eta_j$ are $\Q$-linearly independent.\\[5pt]

[Consider a dependence relation over $\Q$:
\[
q_1 \eta_1 + q_2 \eta_2 + q_3 \eta_3 \ = \ (0, 0, 0, 0)
\] 
So
\allowdisplaybreaks
\begin{align*}
q_1(\beta_{1 1}, \beta_{2 1}, x_1 y_1 &- \beta_{1 1}, x_2 y_1 - \beta_{2 1}) 
+ q_2(\beta_{1 2}, \beta_{2 2}, x_1 y_2 - \beta_{1 2}, x_2 y_2 - \beta_{2 2})
\\[12pt]
&
\hspace{2cm}
+ q_3(\beta_{1 3}, \beta_{2 3}, x_1 y_3 - \beta_{1 3}, x_2 y_3 - \beta_{2 3})
\\[12pt]
\hspace{3cm}
&= \ (0, 0, 0, 0)
\end{align*}

\hspace{1.5cm} $\implies$
\[
q_1 \beta_{1 1} + q_2 \beta_{1 2}  + q_3 \beta_{1 3}  \ = \ 0
\]

\hspace{1.5cm}  $\implies$
\allowdisplaybreaks
\begin{align*}
q_1(x_1 y_1 - \beta_{1 1}) \  + \  &q_2(x_1 y_2 - \beta_{1 2}) + q_3(x_1 y_3 - \beta_{1 3})
\\[12pt]
&=\ 
q_1 x_1 y_1 + q_2 x_1 y_2 + q_3 x_1 y_3 
\ - \ 
q_1 \beta_{1 1} + q_2 \beta_{1 2}  + q_3 \beta_{1 3}
\\[12pt]
&=\ 
q_1 x_1 y_1 + q_2 x_1 y_2 + q_3 x_1 y_3
\\[12pt]
&=\ 
0
\end{align*}
or still, upon dividing by $x_1 \neq 0$, 
\[
q_1 y_1 + q_2 y_2 + q_3 y_3  \ = \ 0
\]
\qquad\qquad\qquad $\implies$
\[
q_1 = 0, \ 
q_2 = 0, \ 
q_3 = 0.]
\]
\fi
Therefore
\[
\sV \hsx \cap \hsx \big(\Qbar^{\hsx 2} \times \{0\} \big) \ \neq \ \{0\},
\]
say
\[
(z_1, z_2, z_3, z_4) \in \sV \hsx \cap \hsx \big(\Qbar^{\hsx 2} \times \{0\} \big)
\]
\qquad\qquad $\implies$
\[
z_1 \in \Qbar,\ 
z_2 \in \Qbar\ 
\&\ 
z_3 = 0, \ 
z_4 = 0.
\]
And
\[
x_2(z_1 + z_3) - x_1(z_2 + z_4) \ = \ 0
\]
\qquad\qquad $\implies$
\[
x_2(z_1) - x_1 (z_2) \ = \ 0
\]
\qquad\qquad $\implies$
\[
\frac{x_2}{x_1} (z_1) \ = \ z_2.
\]
\ifdmcx
\\[-.1cm]

\qquad \textbullet \quad Neither $z_1$ nor $z_2$ can be zero (see below).  
\\[-.1cm]

\qquad \textbullet \quad Thus
\[
\frac{x_2}{x_1} \ = \ \frac{z_2}{z_1} 
\]
is an algebraic number not in $\Q$ (see below). 
\fi
\\[-.1cm]

\noindent Now put $\gamma = \ds\frac{x_2}{x_1}$ and write
\[
\lambda_{2 j} + \beta_{2 j}  
\ = \ 
\gamma(\lambda_{1 j} + \beta_{1 j}) 
\qquad (j = 1, 2, 3)
\]
or still, 
\[
\gamma \lambda_{1 j}  - \lambda_{2 j} 
\ = \ 
\beta_{2 j}  - \gamma \beta_{1 j}
\qquad (j = 1, 2, 3).
\]
\ifdmcx
The entity $\beta_{2 j}  - \gamma \beta_{1 j}$ is an algebraic number.
\\[-.1cm]

\qquad \textbullet \quad 
 Then on general grounds (see below)
\[
\beta_{2 j}  - \gamma \beta_{1 j} \ = \ 0.
\]
Therefore
\fi
\[
\gamma \lambda_{1 j}  - \lambda_{2 j} 
\ = \ 
0 
\implies
\gamma \lambda_{1 j}  \ = \  \lambda_{2 j}.
\]
To finish the proof, make the claim that
\[
\begin{cases}
\ \lambda_{1 j} = 0\\[8pt]
\ \lambda_{2 j} = 0
\end{cases}
\quad(j = 1, 2, 3).
\]
To argue this, assume that $\lambda_{1 j} \neq 0$, so 
\ifdmcx
\\[-.1cm]

\qquad \textbullet \quad 
$\gamma = \ds\frac{\lambda_{2 j}}{\lambda_{1 j}}$ is transcendental (see below) (recall that $\gamma \notin \Q)$.  \\
\\[-.5cm]

\noindent Accordingly
\begin{align*}
&\gamma_{1 j}  \ = \ 0
\\[5pt]
\implies \qquad
&\gamma \hsx 0 - \lambda_{2 j} \ = \ 0 
\\[5pt]
\implies \qquad
&\lambda_{2 j}  \ = \ 0.
\end{align*}
\fi

[Note: \ 
Details$-$
\vspace{.25cm}

\qquad \textbullet \ 
Consider a dependence relation over $\Q$:
\[
q_1 \eta_1 + q_2 \eta_2 + q_3 \eta_3 \ = \ (0, 0, 0, 0)
\]
which, when unraveled, becomes
\allowdisplaybreaks
\begin{align*}
q_1(\beta_{1 1}, \beta_{2 1}, x_1 y_1 &- \beta_{1 1}, x_2 y_1 - \beta_{2 1}) 
+ q_2(\beta_{1 2}, \beta_{2 2}, x_1 y_2 - \beta_{1 2}, x_2 y_2 - \beta_{2 2})
\\[12pt]
&
\hspace{2cm}
+ q_3(\beta_{1 3}, \beta_{2 3}, x_1 y_3 - \beta_{1 3}, x_2 y_3 - \beta_{2 3})
\\[12pt]
\hspace{3cm}
&= \ (0, 0, 0, 0)
\end{align*}

\hspace{1.5cm} $\implies$
\[
q_1 \beta_{1 1} + q_2 \beta_{1 2}  + q_3 \beta_{1 3}  \ = \ 0
\]

\hspace{1.5cm}  $\implies$
\allowdisplaybreaks
\begin{align*}
q_1(x_1 y_1 - \beta_{1 1}) \  + \  &q_2(x_1 y_2 - \beta_{1 2}) + q_3(x_1 y_3 - \beta_{1 3})
\\[12pt]
&=\ 
q_1 x_1 y_1 + q_2 x_1 y_2 + q_3 x_1 y_3
\\[12pt]
&=\ 
0
\end{align*}
or still, upon dividing by $x_1 \neq 0$, 
\[
q_1 y_1 + q_2 y_2 + q_3 y_3  \ = \ 0
\]
\qquad\qquad\qquad $\implies$
\[
q_1 = 0, \ 
q_2 = 0, \ 
q_3 = 0.
\]

\qquad \textbullet \quad
Suppose that the canonical conditions were satisfied $-$then
\allowdisplaybreaks
\begin{align*}
\dim_{\Q} (\sV \cap \fL_G) \ 
&\leq \ d_1(n - t)
\\[12pt]
&=\ 
2(3 - t).
\end{align*}
There are two possibilities for $t$:
\[
\begin{cases}
\ t = 2 \implies 2(3 - 2) = 2\\[8pt]
\ t = 3 \implies 2(3 - 3) = 0
\end{cases}
.
\]
But
\[
\dim_{\Q} (\sV \hsx \cap \hsx \fL_G) \ \geq \ 3, 
\]
$\eta_1$, $\eta_2$, $\eta_3$ being three $\Q$-linearly independent points of $\sV \hsx \cap \hsx \fL_G$.
\vspace{0.3cm}

\qquad \textbullet \quad
The formula
\[
x_2(z_1) - x_1(z_2) \ = \ 0
\]
is a $\Qbar$ dependence relation per $\{x_1, x_2\}$.  
Claim: \ $z_1 \neq 0$, $z_2 \neq 0$.  
E.g.: \ Suppose $z_1 = 0$, hence $x_1(z_2) = 0$ $\implies$ $z_2 = 0$ $(x_1 \neq 0)$.
\vspace{0.4cm}

\qquad \textbullet \quad
$\ds\frac{x_2}{x_1}$ is a nonzero algebraic number and $\ds\frac{x_2}{x_1} \notin \Q$.  
For if $\ds\frac{x_2}{x_1} \in \Q$, we could write
\[
x_2 - \bigg(\frac{x_2}{x_1}\bigg) x_1 \ = \ 0
\]
and thereby contradict the $\Q$-linear independence of $x_1$, $x_2$.
\\[-.1cm]

\qquad \textbullet \quad
If
\[
\beta_0 + \beta_1 \lambda_1 + \cdots + \beta_d \lambda_d \ = \ 0,
\]
where $\beta_0, \beta_1, \ldots, \beta_d$ are algebraic and  
$\lambda_1 \in \fL, \ldots, \lambda_d \in \fL$, then $\beta_0 = 0$ 
(cf. \S39, \#14).
\\[-.1cm]

\qquad \textbullet \quad
The quotient $\ds\frac{u}{v}$ of two nonzero elements of $\fL$ is either rational or transcendental.
\end{x}
\vspace{0.3cm}

\begin{x}{\small\bf IMPLICATION} \ 
\[
\text{sharp six exponentials $\implies$ six exponentials }.
\]

[Take $\beta_{i j} = 0$, so $\forall \ i$, $\forall \ j$, $x_i y_j = 0$, which is false 
$(\forall \ i$, $x_i \neq 0$, $\forall \ j$, $y_j \neq 0$).  
The supposition that the six numbers
\[
e^{x_i y_j}
\]
are algebraic is therefore contradictory, thus at least one of the 
\[
e^{x_i y_j}
\]
is transcendental.]
\end{x}
\vspace{0.3cm}

\begin{x}{\small\bf IMPLICATION} \ 
\[
\text{sharp six exponentials $\implies$ five exponentials}.
\]

[Explicate the parameters in \S41, \#2:
\[
e^{x_1 y_1}, \ 
e^{x_1 y_2}, \ 
e^{x_2 y_1}, \ 
e^{x_2 y_2}, \
e^{\gamma x_1 / x_2}.
\]
Put
\[
y_3 \ = \ \gamma / x_2,
\]
let
\[
\beta_{1 1} \ = \ \beta_{1 2} \ = \ \beta_{1 3}  \ = \ \beta_{2 1} \ = \ \beta_{2 2} \ = \ 0,
\]
and let
\[
\beta_{2 3} \ = \ \gamma.
\]
To incorporate the denial of \S41, \#2, assume that the six numbers
\[
e^{x_i y_j \hsx - \hsx \beta_{i j}}
\]
are algebraic.  Note that
\[
e^{x_1 y_3 \hsx - \hsx \beta_{1 3}}
\ = \ 
e^{x_1 y_3 \hsx -  \hsx 0}
\ = \ 
e^{\gamma x_1 / x_2}
\]
and
\[
e^{x_2 y_3 \hsx - \hsx  \beta_{2 3}} \ = \ e^{\gamma \hsx - \hsx \gamma} \ = \ 1.
\]
Now apply \#1: 
\[
x_i y_j \ = \ \beta_{i j} \qquad (i = 1, 2, \ j = 1, 2, 3),
\]
so
\[
x_1 y_1 \ = \ \beta_{1 1} \ = \ 0,\ 
x_1 y_2 \ = \ \beta_{1 2} \ = \ 0,\ 
x_2 y_1 \ = \ \beta_{2 1} \ = \ 0,\ 
x_2 y_2 \ = \ \beta_{2 2} \ = \ 0, 
\]
so we have our contradiction.  Of course
\[
x_1 y_3 \ = \ \beta_{1 3} \ = \ 0, \quad
x_2 y_3 \ = \ \beta_{2 3} \ = \ \gamma
\]
but these formulas do not figure in the deduction and are merely part of the formalism.
\vspace{0.2cm}

[Note: \ 
There is a potential gap in the argument, viz. why is $\{y_1, y_2, y_3\}$ a $\Q$-linearly independent set?  
Thus consider a rational dependence relation
\[
q_1 y_1 + q_2 y_2 + q_3 \gamma / x_1 \ = \ 0.
\]
Multiply  through by $x_1$:
\[
q_1 x_1 y_1 + q_2 x_1 y_2 + q_3 \gamma \ = \ 0.
\]
Since
\[
e^{x_1 y_1} \in \Qbar, \quad 
e^{x_1 y_2} \in \Qbar,
\]
it follows that
\[
\lambda_1 \ \equiv \ x_1 y_1 \in \fL, \ 
\lambda_2 \ \equiv \ x_1 y_2 \in \fL
\]
and our relation reads
\[
q_3 \gamma + q_1 \lambda_1 + q_2 \lambda_2 \ = \ 0.
\]
But $\{x_1, x_2\}$ is a $\Q$-linearly independent set, $\lambda_1 \in \fL$, $\lambda_2 \in \fL$ are nonzero and 
$\Q$-linearly independent, hence with 
\[
\beta_0 \ = \ q_3 \gamma,\ 
\beta_1 \ = \ q_1,\ 
\beta_2 \ = \ q_2,
\]
we have
\[
\beta_0 + \beta_1 \lambda_1 + \beta_2 \lambda_2 \ = \ 0.
\]
Therefore $\beta_0 = 0$ (cf. \S39, \#14)
\begin{align*}
\implies &q_3 \hsx = \hsx 0\ 
\\[12pt]
\implies &q_1 \hsx = \hsx 0, \ q_2 \hsx = \hsx 0.]
\end{align*}
\end{x}
\vspace{0.3cm}


%% file: _43_strong_six_exp_theorem.tex
\chapter{
$\boldsymbol{\S}$\textbf{43}.\quad  STRONG SIX EXPONENTIALS THEOREM}
\setlength\parindent{2em}
\setcounter{theoremn}{0}
\renewcommand{\thepage}{\S43-\arabic{page}}

\ \indent 

Denote by $\fL^*$ the $\Qbar$-vector space spanned by 1 and $\fL$ in $\C$, thus
\[
\fL^* \ 
= \ 
\{\beta_0 + \beta_1 \lambda_1 + \cdots + \beta_n \lambda_n \ : 
n \geq 0, (\beta_0, \beta_1, \ldots, \beta_n) \in \Qbar^{\hsx n+1}, (\lambda_1, \ldots, \lambda_n) \in \fL^n\}.
\]
\\[-1.4cm]

[Note: \ 
$\fL^*$, like $\fL$, is stable under complex conjugation.]
\vspace{0.5cm}

\begin{x}{\small\bf THEOREM} \ 
Let $\{x_1, x_2\}$  and $\{y_1, y_2, y_3\}$  be two $\Qbar$-linearly independent sets of complex numbers $-$then 
\[
\{x_1 y_1, x_1 y_2, x_1 y_3, x_2 y_1, x_2 y_2, x_2 y_3\} \not\subset \fL^*,
\]
i.e., $\exists \ i \in \{1,2\}$, $\exists \ j \in \{1,2,3\}$:
\[
x_i y_j \notin \fL^*,
\]
hence $e^{x_i y_j}$ is transcendental.
\end{x}
\vspace{0.3cm}

This result, due to Damien Roy, is the strong six exponentials theorem (proof omitted).
\\[-.5cm]

[Note: \ 
The reason for the appelation ``strong'' as compared with the six exponentials theorem per se is that one of the 
$x_i y_j$ $(1 \leq i \leq 2, 1 \leq j \leq 3)$ is not in $\fL$ but even more, viz. it is not in $\fL^*$.]
\vspace{0.3cm}

\begin{x}{\small\bf STRONG CONDITION X} \ 
Suppose that $\lambda_0 \in \fL^*$, $\lambda_1 \in \fL^*$, $\lambda_2 \in \fL^*$, $\lambda_3 \in \fL^*$.  
Assume: \ $\{\lambda_0, \lambda_1\}$ is $\Qbar$-linearly independent and 
$\{\lambda_0, \lambda_2, \lambda_3\}$ is $\Qbar$-linearly independent 
$-$then
\[
\bigg\{\frac{\lambda_1 \lambda_2}{\lambda_0}, \frac{\lambda_1 \lambda_3}{\lambda_0}\bigg\} \not\subset \fL^*.
\]
\vspace{0.2cm}

PROOF \ 
In \#1, take
\[
x_1 = 1, \ 
x_2 = \frac{\lambda_1}{\lambda_0}, \ 
y_1 = \lambda_0, \ 
y_2 = \lambda_2, \ 
y_3 = \lambda_3.
\]
Then
\[
\bigg\{\lambda_0, \lambda_2, \lambda_3,\lambda_1, \frac{\lambda_1 \lambda_2}{\lambda_0}, \frac{\lambda_1 \lambda_3}{\lambda_0}\bigg\} 
\not\subset \fL^*.
\]
But by hypothesis,
\[
\{\lambda_0, \lambda_2, \lambda_3,\lambda_1 \} \subset \fL^*.
\]
Therefore
\[
\bigg\{\frac{\lambda_1 \lambda_2}{\lambda_0}, \frac{\lambda_1 \lambda_3}{\lambda_0}\bigg\} \not\subset \fL^*.
\]
\end{x}
\vspace{0.3cm}

\begin{x}{\small\bf THEOREM} \ 
The strong condition X implies the strong six exponentials theorem.
\vspace{0.2cm}

PROOF \ 
To devise a contradiction, assume that the six products $x_i y_j$ $(1 \leq i \leq 2, 1 \leq j \leq 3)$ are in $\fL^*$.  
Apply strong condition X as follows: \ Take
\[
\lambda_0 = x_1 y_1, \ 
\lambda_1 = x_2 y_1, \ 
\lambda_2 = x_1 y_2, \ 
\lambda_3 = x_1 y_3.
\]
Then $\{\lambda_0, \lambda_1\}$ is $\Qbar$-linearly independent, as is $\{\lambda_0, \lambda_2, \lambda_3\}$.  
Consequently either 
\[
\frac{\lambda_1 \lambda_2}{\lambda_0} \notin \fL^* 
\quad \text{or} \quad
\frac{\lambda_1 \lambda_3}{\lambda_0} \notin \fL^* 
\qquad \text{(or both)}.
\]
But
\[
\begin {cases}
\ \ds\frac{\lambda_1 \lambda_2}{\lambda_0}  
\ = \ 
\ds\frac{x_2 y_1 x_1 y_2}{x_1 y_1}
\ = \ 
x_2 y_2 \hsx \in \hsx \fL^*\\[15pt]
\ds\frac{\lambda_1 \lambda_3}{\lambda_0}
\ = \ 
\ds\frac{x_2 y_1 x_1 y_3}{x_1 y_1}
\ = \ 
x_2 y_3 \hsx \in \hsx \fL^* 
\end {cases}
.
\]
Contradiction.
\end{x}
\vspace{0.3cm}

\begin{x}{\small\bf LEMMA} \ 
Suppose that $\lambda_1$, $\lambda_2 \in \fL^*$ $(\lambda_2 \neq 0)$.  
Assume: \ $\{1, \lambda_1, 1 / \lambda_2\}$ is $\Qbar$-linearly indedendent $-$then
\[
\{\lambda_1 \lambda_2, 1 / \lambda_2\} \not\subset \fL^*.
\]
\vspace{0.2cm}

PROOF \ 
If $1 / \lambda_2 \notin \fL^*$, then we are done.  
Otherwise, apply strong condition X to the family $\{1 / \lambda_2, 1, \lambda_1, 1\}$ and conclude that 
\[
\{\lambda_1 \lambda_2, \lambda_2\} \not\subset \fL^*,
\]
hence $\lambda_1 \lambda_2 \notin \fL^*$.
\end{x}
\vspace{0.3cm}

\begin{spacing}{1.45}
\begin{x}{\small\bf SCHOLIUM} \ 
Suppose that $\lambda \in \fL^*$ $(\lambda \neq 0)$ is transcendental $-$then
\[
\{\lambda^2, 1 / \lambda\} \not\subset \fL^*.
\]

[In \#4, take $\lambda_1 = \lambda$, $\lambda_2 = \lambda$ $-$then the isssue is whether 
$\{1, \lambda, 1 / \lambda\}$ is $\Qbar$-linearly independent.  So consider a dependence relation
\[
r + s \lambda + t(1 / \lambda) \ = \ 0,
\]
where $r$, $s$, $t \in \Qbar$.  
Multiply by $\lambda$ to get
\[
r \lambda + s \lambda^2 + t \ = \ 0.
\]
Since $\lambda$ is transcendental, it follows that $\{\lambda, \lambda^2, 1\}$ is algebraically independent 
over $\Q$, hence is algebraically independent over $\Qbar$ (cf. \S20, \#7), hence is $\Qbar$-linearly independent, 
hence $r = 0$, $s = 0$, $t = 0$.]
\end{x}
\end{spacing}
\vspace{0.3cm}

\begin{x}{\small\bf APPLICATION} \ 
Take $\lambda = \pi \sqrt{-1}$ $-$then $\lambda \in \fL \subset \fL^*$ and 
\[
\big\{-\pi^2, 1 / \pi \sqrt{-1} \big\} \not\subset \fL^*.
\]
Therefore
\[
\pi^2 \notin \fL^* 
\quad \text{or} \quad 
1 / \pi \notin \fL^* 
\qquad \text{(or both)}
\]
which implies that either
\[
 \text{$e^{\pi^2}$ is transcendental or $e^{1 / \pi}$ is transcendental (or both).}
\]
\end{x}
\vspace{0.3cm}

\begin{spacing}{1.45}
\begin{x}{\small\bf SUBLEMMA} \ 
Let $x_1$, $x_2$, $y_1$, $y_2$ be complex numbers and let $\gamma$ be a nonzero algebraic number.  
Suppose that $\{x_1, x_2\}$ is $\Qbar$-linearly independent and $\{y_1, y_2, \gamma / x_1\}$ is $\Qbar$-linearly independent.  
Assume: 
\[
\gamma x_2 / x_1 \in \fL^*.
\]
Then
\[
\{x_1 y_1, x_1 y_2, x_2 y_1, x_2 y_2\} \not\subset \fL^*.
\]

PROOF \ 
Apply \#1 to 
\[
\{x_1, x_2\}
\quad \text{and} \quad 
\{y_1, y_2, \gamma / x_1\}
\]
which leads to
\[
\{x_1 y_1, x_1 y_2, x_1(\gamma / x_1), x_2 y_1, x_2 y_2, x_2(\gamma / x_1)\}.
\]
Of course, 
\[
x_1(\gamma / x_1) \ = \ \gamma \in \fL^*
\]
and by hypothesis,
\[
x_2(\gamma / x_1) \ = \ \gamma x_2 / x_1 \in \fL^*,
\]
leaving
\[
\{x_1 y_1, x_1 y_2, x_2 y_1, x_2 y_2\}.
\]
\end{x}
\end{spacing}
\vspace{0.3cm}

\begin{spacing}{1.45}
\begin{x}{\small\bf LEMMA} \ 
Let $x_1$, $x_2$, $y_1$, $y_2$ be complex numbers and let $\gamma$ be a nonzero algebraic number.  
Suppose that $\{x_1, x_2\}$  is $\Q$-linearly independent and $\{y_1, y_2\}$  is $\Q$-linearly independent.  
Assume:
\[
\gamma x_2 / x_1 \in \fL^*.
\]
Then
\[
\{x_1 y_1, x_1 y_2, x_2 y_1, x_2 y_2\} \not\subset \fL.
\]
\vspace{0.2cm}

PROOF \ 
Assume instead that
\[
\{x_1 y_1, x_1 y_2, x_2 y_1, x_2 y_2\} \subset \fL.
\]

\qquad \textbullet \quad
$\{x_1 y_1, x_2 y_1\}$  is $\Q$-linearly independent, hence is $\Qbar$-linearly independent (Gelfond-Schneider) 
(for $x_1 y_1 \in \fL$, $x_2 y_1 \in \fL$), hence $\{x_1, x_2\}$  is $\Qbar$-linearly independent.
\vspace{0.3cm}

\qquad \textbullet \quad
$\{x_1 y_1, x_1 y_2\}$  is $\Q$-linearly independent, hence $\{1, x_1 y_1, x_1 y_2\}$ is $\Qbar$-linearly independent 
(inhomogeneous Baker) (for $x_1 y_1 \in \fL$, $x_1 y_2 \in \fL$), hence $\{\gamma / x_1, y_1, y_2\}$ is $\Qbar$-linearly independent.
\vspace{0.2cm}

Therefore (cf. \#7)
\[
\gamma x_2 / x_1 \notin \fL^*.
\]
\vspace{0.2cm}

[Note: \ 
To check that $\{\gamma / x_1, y_1, y_2\}$ is $\Qbar$-linearly independent, write
\[
r(\gamma / x_1) + s y_1 + t y_2 \ = \ 0,
\]
where $r$, $s$, $t \in \Qbar$ $-$then
\allowdisplaybreaks
\begin{align*}
&r \gamma  + s x_1 y_1 + t x_1 y_2 \ = \ 0
\\
\implies \hspace{1cm}
&
\\
&r \gamma = 0, \ s = 0, \ t = 0.
\end{align*}
But $\gamma \in \Qbar$ is nonzero, so $r = 0$.]
\end{x}
\end{spacing}
\vspace{0.3cm}

\begin{x}{\small\bf \un{N.B.}} \ 
The strong six exponentials theorem intervenes in \#8 via an application of \#7.
\end{x}
\vspace{0.3cm}

\begin{x}{\small\bf RAPPEL} \ 
Let $\{x_1, x_2\}$ and $\{y_1, y_2\}$ be two $\Q$-linearly independent sets of complex numbers.  
Let further $\gamma$ be a nonzero algebraic number $-$then at least one of the five numbers
\[
e^{x_1 y_1}, \ 
e^{x_1 y_2}, \ 
e^{x_2 y_1}, \ 
e^{x_2 y_2}, \ 
e^{\gamma x_2 / x_2} 
\]
is transcendental.
\vspace{0.2cm}

[This is the five exponentials theorem (cf. \S41, \#2) (switch the roles of $x_1$ and $x_2$).]
\end{x}
\vspace{0.3cm}

\begin{x}{\small\bf IMPLICATION} \ 
\[
\text{strong six exponentials} 
\ \implies \ 
\text{five exponentials.}
\]

[The claim is that at least one of the five numbers
\[
e^{x_1 y_1}, \ 
e^{x_1 / y_2}, \  
e^{x_2 y_1}, \ 
e^{x_2 y_2}, \ 
e^{\gamma x_2 / x_1} 
\]
is transcendental.
\vspace{0.3cm}

\qquad \textbullet \quad \un{Case 1:} \ 
$ \gamma x_2 / x_1 \notin \fL^*$ $-$then
\[
e^{\gamma x_2 / x_1}
\]
is transcendental.
\vspace{0.3cm}

\qquad \textbullet \quad \un{Case 2:} \ 
$ \gamma x_2 / x_1 \in \fL^*$ $-$then
\[
\{
x_1 y_1, \ 
x_1 y_2, \  
x_2 y_1, \ 
x_2 y_2
\} \not\subset \fL \qquad \text{(cf. \#8)},
\]
i.e., $\exists \ i \in \{1, 2\}$, $\exists \ j \in \{1, 2, 3\}$:
\[
x_i y_j \notin \fL,
\]
hence $e^{x_i y_j}$ is transcendental.]
\end{x}
\vspace{0.3cm}

\begin{x}{\small\bf REMARK} \ 
Refer to \S41, \#7.  
Make the assumption that $x_2/ x_1 \in \fL^*$ $-$then for some pair $(i,j)$ : $x_i y_j \notin \fL$, implying thereby that $e^{x_i y_j}$ is transcendental, as desired.
\end{x}
\vspace{0.3cm}

\begin{x}{\small\bf RAPPEL} \ 
Let $\{x_1, x_2\}$ and $\{y_1, y_2, y_3\}$ be two $\Q$-linearly independent sets of complex numbers $-$then
\[
\{
x_1 y_1, \ 
x_1 y_2, \ 
x_1 y_3, \ 
x_2 y_1, \ 
x_2 y_2, \ 
x_2 y_3
\} \not\subset \fL.
\]

[This is the six exponentials theorem.]
\end{x}
\vspace{0.3cm}

\begin{x}{\small\bf CONDITION X} \ 
Suppose that $\lambda_0 \in \fL$, $\lambda_1 \in \fL$, $\lambda_2 \in \fL$, $\lambda_3 \in \fL$.  
Assume: \ 
$\{\lambda_0, \lambda_1\}$ is $\Q$-linearly independent and 
$\{\lambda_0, \lambda_2, \lambda_3\}$ is $\Q$-linearly independent $-$then
\[
\bigg\{\frac{\lambda_1 \lambda_2}{\lambda_0} \hsx , \frac{\lambda_1 \lambda_3}{\lambda_0} \bigg\} 
\not\subset \fL.
\]

[In \#2, replace $\Qbar$ by $\Q$ and $\fL^*$ by $\fL$.]

\end{x}
\vspace{0.3cm}

Imitating the proof that the strong six exponentials theorem is equivalent to strong condition X, it follows that the six exponentials theorem is equivalent to condition X.
\vspace{0.3cm}

\begin{spacing}{1.45}
\begin{x}{\small\bf IMPLICATION} \ 
\[
\text{strong six exponentials} 
\ \implies \ 
\text{six exponentials.}
\]
\\[-1.5cm]

[Start with the data for condition X $-$then thanks to homogeneous Baker, $\{\lambda_0, \lambda_1\}$ is $\Qbar$-linearly independent and $\{\lambda_0, \lambda_2, \lambda_3\}$ is $\Qbar$-linearly independent, the setup for strong condition X, hence (cf. \#2),
\[
\bigg\{\frac{\lambda_1 \lambda_2}{\lambda_0} , \frac{\lambda_1 \lambda_3}{\lambda_0} \bigg\} 
\not\subset 
\fL^*
\]
\qquad\qquad $\implies$
\[
\bigg\{\frac{\lambda_1 \lambda_2}{\lambda_0} , \frac{\lambda_1 \lambda_3}{\lambda_0} \bigg\} 
\not\subset  
\fL.]
\]
\end{x}
\end{spacing}
\vspace{0.3cm}

\[
\textbf{APPENDIX}
\]
\vspace{0.5cm}

It was established in \S36 that the six exponentials theorem is equivalent to the following statement.
\vspace{0.5cm}

{\small\bf SCHOLIUM} \ 
Consider a nonzero $2 \times 3$ matrix $\sM$ with entries in $\fL$:
\[
\sM \ = \ 
\begin{pmatrix}
\lambda_{1 1} &\lambda_{1 2} &\lambda_{1 3}\\
\\
\lambda_{2 1} &\lambda_{2 2} &\lambda_{2 3}
\end{pmatrix}
.
\]
Suppose that its rows are $\Q$-linearly independent and its columns are $\Q$-linearly independent $-$then
\[
\rank \sM \ = \ 2.
\]
\vspace{0.1cm}

Analogously, the strong exponentials theorem is equivalent to the following statement. 
\vspace{0.3cm}

{\small\bf SCHOLIUM} \ 
Consider a nonzero $2 \times 3$ matrix $\sM$ with entries in $\fL^*$:
\[
\sM \ = \ 
\begin{pmatrix}
\lambda_{1 1} &\lambda_{1 2} &\lambda_{1 3}\\
\\
\lambda_{2 1} &\lambda_{2 2} &\lambda_{2 3}
\end{pmatrix}
.
\]
Suppose that its rows are $\Qbar$-linearly independent and its columns are $\Qbar$-linearly independent $-$then
\[
\rank \sM \ = \ 2.
\]
\vspace{0.3cm}

{\small\bf \un{N.B.}} \ 
Once again, 
\[
\text{strong six exponentials \ $\implies$ \ six exponentials.}
\]

[Start with
\[
\sM \ = \ 
\begin{pmatrix}
\lambda_{1 1} &\lambda_{1 2} &\lambda_{1 3}\\
\\
\lambda_{2 1} &\lambda_{2 2} &\lambda_{2 3}
\end{pmatrix}
\qquad (\lambda_{i j} \in \fL).
\]\\[8pt]
Then the assumption of the $\Q$-linear independence of its rows and columns implies the $\Qbar$-linear independence of its rows and columns (homogeneous Baker).]
\vspace{0.75cm}

Finally, the sharp six exponentials theorem is equivalent to the following statement.
\vspace{0.75cm}

{\small\bf SCHOLIUM} \ 
Consider a nonzero $2 \times 3$ matrix $\sM$ with entries in $\Qbar + \fL$:
\[
\sM \ = \ 
\begin{pmatrix}
\lambda_{1 1} &\lambda_{1 2} &\lambda_{1 3}\\
\\
\lambda_{2 1} &\lambda_{2 2} &\lambda_{2 3}
\end{pmatrix}
.
\]
Suppose that its rows are $\Qbar$-linearly independent and its columns are $\Qbar$-linearly independent $-$then
\[
\rank \sM \ = \ 2.
\]
\vspace{0.3cm}

{\small\bf REMARK} \ 
Consequently
\[
\text{strong six exponentials \ $\implies$ \ sharp six exponentials.}
\]

To help keep it all straight, make a chart of the various implications:

\[
\begin{tikzcd}
&{\text{strong 6 exponentials}}  
\ar[bend right]{ddddl}
\ar[dddd,bend left=90,looseness=1.5,swap]
\ar[Rightarrow]{dd}
\arrow[Rightarrow, dd]
\\
\\
&{\text{sharp 6 exponentials}}
\arrow[Rightarrow, dd]
\\
\\
{\text{5 exponentials}} \ar{r}[swap]{?} &{\text{6 exponentials}}
\end{tikzcd}
.
\]

\vspace{0.3cm}

%% file: _44_the_four_exponentials_conjecture.tex
\chapter{
$\boldsymbol{\S}$\textbf{44}.\quad  FOUR EXPONENTIALS CONJECTURE (4EC)}
\setlength\parindent{2em}
\setcounter{theoremn}{0}
\renewcommand{\thepage}{\S44-\arabic{page}}

\ \indent 
This is the following statement
\vspace{0.5cm}

\begin{x}{\small\bf CONJECTURE} \ 
Let $\{x_1, x_2\}$ and $\{y_1, y_2\}$ be two $\Q$-linearly independent sets of complex numbers $-$then 
\[
\{x_1 y_1, x_1 y_2, x_2 y_1, x_2 y_2 \} \not\subset \fL, 
\]
thus at least one of the numbers
\[
e^{x_1 y_1}, \ 
e^{x_1 y_2}, \ 
e^{x_2 y_1}, \ 
e^{x_2 y_2}
\]
is transcendental.
\end{x}
\vspace{0.3cm}

In terms of matrices (see the Appendix to \S43):
\vspace{0.5cm}

\begin{x}{\small\bf CONJECTURE}  \ 
Consider a $2 \times 2$ matrix $\sM$ with entries in $\fL$:
\[
\sM \ = \ 
\begin{pmatrix}
\lambda_{1 1} &\lambda_{1 2}\\[8pt]
\lambda_{2 1} &\lambda_{2 2}
\end{pmatrix}
.
\]
Suppose that its rows are $\Q$-linearly independent and its columns are $\Q$-linearly independent $-$then 
\[
\rank \sM \ = \ 2.
\]
\end{x}
\vspace{0.3cm}

\begin{x}{\small\bf EXAMPLE} \ 
Consider the matrix
\[
\begin{pmatrix}
1 &\pi\\[8pt]
\pi &\pi^2
\end{pmatrix}
.
\]
Its determinant is 0 and its rank is 1.  
This is not a contradiction since $\pi$, $\pi^2 \notin \fL$.
\vspace{0.2cm}

[Note: \ 
Still, its rows and columns are $\Q$-linearly independent.]

\end{x}
\vspace{0.3cm}

\begin{x}{\small\bf LEMMA} \ \ 
$\#1 \Leftrightarrow \#2$.
\end{x}
\vspace{0.3cm}

\begin{x}{\small\bf REMARK} \ 
The four exponentials conjecture is a long outstanding open problem in transcendence theory.
\end{x}
\vspace{0.3cm}

\begin{x}{\small\bf EXAMPLE} \  (Admit 4EC) \ 
Use the notation of \S36, \#6.  Introduce as there
\[
E_2 \ = \ \{t \in \R : 2^t, \ 3^t \in \N\}.
\]
Then
\[
E_2 \ = \ \N.
\]
\vspace{0.2cm}

[Given $t \in \R$, $t \notin \Q$, take in \#1
\[
\begin{cases}
\ x_1 = 1\\
\ x_2 = t
\end{cases}
,
\qquad\qquad
\begin{cases}
\ y_1 = \elln(2)\\
y_2 = \elln(3)
\end{cases}
.
\]
Then the four exponentials are 
\[
2, \ 
3, \ 
2^t, \ 
3^t
\]
and either 
\[
2^t \quad \text{or} \quad 3^t \quad \text{(or both)}
\]
is (are) transcendental.  Therefore 
\[
(\R - \Q) \hsx \cap \hsx E_2 \ = \ \emptyset.
\]
But
\[
E_2 \hsx \cap \hsx \Q \ = \ \N.
\]
And
\allowdisplaybreaks
\begin{align*}
E_2 \hsx \cap \hsx \Q \ 
&=\ 
E_2 \hsx \cap \hsx (\Q \hsx \cup \hsx (\R - \Q))
\\[12pt]
&=\ 
E_2 \hsx \cap \hsx \R
\\[12pt]
&=\ 
E_2.]
\end{align*}
\end{x}
\vspace{0.3cm}

\begin{x}{\small\bf EXAMPLE} \  (Admit 4EC) \ 
Let $\lambda \in \fL$, $\lambda \notin \R$ $-$then $e^{\abs{\lambda}}$ is transcendental.
\vspace{0.3cm}

[In \#1, take
\[
\begin{cases}
\ x_1 = 1\\
\ x_2 = \abs{\lambda} / \lambda
\end{cases}
,\qquad\qquad
\begin{cases}
\ y_1 = \lambda\\
y_2 = \abs{\lambda}
\end{cases}
.
\]
Then the four exponentials are 
\[
e^\lambda, \ 
e^{\abs{\lambda}}, \ 
e^{\abs{\lambda}}, \ 
e^{\abs{\lambda}^2} / \lambda.
\]
Here $e^\lambda \in \ov{\Q}$.  
And
\allowdisplaybreaks
\begin{align*}
\abs{\lambda}^2 \ = \ \lambda \ov{\lambda} 
&\implies
{\abs{\lambda}^2} / \lambda \ = \ \frac{\lambda \ov{\lambda}}{\lambda} \ = \ \ov{\lambda}
\\[12pt]
&\implies
e^{\abs{\lambda}^2} / \lambda \ = \ e^{\ov{\lambda}} \hsx \in \hsx  \Qbar.
\end{align*}
Therefore $e^{\abs{\lambda}}$ is transcendental.]
\vspace{0.2cm}

[Note: \ 
One should check that $\{x_1, x_2\}$ and $\{y_1, y_2\}$  are $\Q$-linearly independent. \\
E.g.: \ Suppose that 
\[
p y_1 + q y_2 \ = \ 0 \qquad (p, q \in \Q)
\]
or still, if $\lambda = a + \sqrt{-1}\hsx b$ $(b \neq 0)$, 
\[
p (a + \sqrt{-1} \hsx b) + q \sqrt{a^2 + b^2} \ = \ 0
\]
\qquad\qquad $\implies$
\[
\begin{cases}
\ p a + q \sqrt{a^2 + b^2} \ = \ 0\\
\ p b \ = \ 0
\end{cases}
\]
\qquad\qquad $\implies$
\[
p \ = \ 0 
\implies
q \sqrt{a^2 + b^2} \ = \ 0
\implies
q \ = \ 0.]
\]
\end{x}
\vspace{0.3cm}

\begin{x}{\small\bf EXAMPLE} \  (Admit 4EC) \ 
In \#1, take 
\[
\begin{cases}
\ x_1  \ = \  1\\
\ x_2 \ = \ \sqrt{2}
\end{cases}
,\qquad\qquad 
\begin{cases}
\ y_1 \ = \ \sqrt{-1} \hsx \pi \\
\ y_2 \ = \ \sqrt{-1} \hsx \pi \hsx \sqrt{2}
\end{cases}
.
\]
Then the four exponentials are 
\[
e^{\sqrt{-1} \hsx \pi}, \ 
e^{\sqrt{-1} \hsx \pi \hsx \sqrt{2}}, \ 
e^{\sqrt{-1} \hsx \pi \hsx \sqrt{2}}, \ 
e^{2 \hsx \sqrt{-1} \hsx \pi}.
\]
The first of these is $-1$, the fourth is $+1$, leaving
\[
e^{\sqrt{-1} \hsx \pi \hsx \sqrt{2}},
\]
which must therefore be transcendental (a consequence already of Gelfond-Schneider: 
\[
e^{\sqrt{-1} \hsx \pi \hsx \sqrt{2}}
 \ = \ 
e^{\sqrt{2} \hsx \Log -1} 
 \ = \ 
 (-1)^{\sqrt{2}}).
\]
\end{x}
\vspace{0.3cm}

\begin{x}{\small\bf EXAMPLE} \  (Admit 4EC) \ 
Let $\lambda \in \fL - \{0\}$ and let $w \in \C - \Q$ (a complex irrational number) $-$then at least one of the two numbers 
\[
e^{\lambda w}, \quad e^{\lambda /  w}
\]
is transcendental.
\vspace{0.2cm}

[In \#1, take 
\[
\begin{cases}
\ x_1  \ = \  \lambda\\
\ x_2 \ = \ w \lambda
\end{cases}
,\qquad\qquad 
\begin{cases}
\ y_1 \ = \ 1\\
\ y_2 \ = \ 1 / w
\end{cases}
.
\]
Then the four exponentials are
\[
e^\lambda \in \ov{\Q}, \ 
e^{\lambda / w}, \ 
e^{w \lambda}, \ 
e^\lambda \in \ov{\Q}.
\]
\vspace{0.2cm}

\begin{spacing}{1.55}
[Note: \ 
There are circumstances when 4EC need not be invoked.  
E.g.: \ Consider the situation when $w \in \ov{\Q} - \Q$.  
In view of \S24, \#8, one of the numbers $w, \ e^\lambda$, and $e^{w \lambda}$ is transcendental.  
But $w$ is algebraic (by hypothesis), $e^\lambda$ is algebraic (by definition), thus $e^{w \lambda}$ is transcendental.]
\end{spacing}
\end{x}
\vspace{0.3cm}

\begin{x}{\small\bf EXAMPLE} \  (Admit 4EC) \ 
Let $w \in \C - \Q$ $-$then
\[
\exp(2 \hsx \pi \hsx \sqrt{-1} \hsx w) 
\quad \text{and} \quad 
\exp(-2 \hsx \pi \hsx \sqrt{-1} / w)
\]
are not simultaneously algebraic.
\vspace{0.2cm}

[Modify \#9 in the obvious way.]
\end{x}
\vspace{0.3cm}

\begin{x}{\small\bf EXAMPLE} \  (Admit 4EC) \ 
Let $\alpha_1$, $\alpha_2$ be positive algebraic numbers different from 1 $-$then $\pi^2$ and 
$\elln(\alpha_1) \hsx \elln(\alpha_2)$ are $\Q$-linearly independent.
\vspace{0.2cm}

[Proceed by contradiction and assume that $\pi^2$ and $\elln(\alpha_1) \hsx \elln(\alpha_2)$ are $\Q$-linearly dependent, 
say for $n, \ m \in \Z$ nonzero, 
\[
n\hsx(\elln(\alpha_1)) \hsx (\elln(\alpha_2)) 
\ = \ 
4 m \pi^2.
\]
Put 
\[
\beta_1 
\ = \ 
\alpha_1^n, 
\quad
\beta_2 
\ = \ 
\exp\bigg( \frac{1}{m} \elln(\alpha_2)\bigg).
\]
Then $\beta_1$, $\beta_2$ are algebraic, nonzero, and $\abs{\beta_1} \neq 1$, $\abs{\beta_2} \neq 1$.  
Moreover
\allowdisplaybreaks
\begin{align*}
\elln(\beta_1) \hsx \elln(\beta_2) \ 
&=\ 
(n\hsx\elln(\alpha_1)) \hsx \bigg( \frac{1}{m} \elln(\alpha_2)\bigg)
\\[12pt]
&=\ 
\frac{n}{m} \hsx \elln(\alpha_1) \hsx \elln(\alpha_2)
\\[12pt]
&=\ 
\frac{n}{m} \hsx \frac{4m}{n} \hsx \pi^2
\\[12pt]
&=\ 
4 \hsx \pi^2.
\end{align*}
Let now
\[
w 
\ = \ 
\elln(\beta_1) \hsx  / \hsx 2 \hsx \pi \hsx \sqrt{-1},
\]
so
\[
\elln(\beta_1) 
\ = \ 
2 \hsx \pi \hsx \sqrt{-1} \hsx w.
\]
Then
\allowdisplaybreaks
\begin{align*}
\elln(\beta_2) \ 
&=\ 
\frac{4 \pi^2}{\elln(\beta_1)}
\\[12pt]
&=\ 
-2 \hsx \pi \hsx \sqrt{-1} / w.
\end{align*}
Since
\[
\begin{cases}
\ \exp(2 \hsx \pi \hsx \sqrt{-1} \hsx w) \ = \ \beta_1\\[8pt]
\ \exp(-2 \hsx \pi \hsx \sqrt{-1} / w) \ = \ \beta_2
\end{cases}
,
\]
it follows that
\[
\exp(2 \hsx \pi \hsx \sqrt{-1} \hsx w) 
\quad \text{and} \quad 
\exp(-2 \hsx \pi \hsx \sqrt{-1} / w)
\]
are algebraic, which contradicts \#10.]
\vspace{0.2cm}

[Note: \ 
In the literature, this result is known as 
\un{Bertrand's conjecture}.]
\index{Bertrand's conjecture}
\end{x}
\vspace{0.3cm}

\begin{x}{\small\bf EXAMPLE} \  (Admit 4EC) \ 
Let $w \in \C - \Q$.  
Assume: \ $\abs{w}^2 \in \Q$ $-$then
\[
\exp(2 \pi\hsx \sqrt{-1} \hsx w)
\]
is transcendental.
\vspace{0.2cm}
%
%

[Assume $\exp(2 \pi\hsx \sqrt{-1} \hsx w)$ is not transcendental.  
and write $q = \abs{w}^2 \in \Q$.
\\[2pt]
So
\begin{align*}
&2 \pi\hsx \sqrt{-1} \hsx w \in \fL
\\[1pt]
\implies \quad
&\ov{2 \pi\hsx \sqrt{-1} \hsx w} \ = \  -2 \pi\hsx \sqrt{-1} \hsx \ov{w} \in \fL
\end {align*}
And
\begin{align*}
&2 \pi\hsx \sqrt{-1}  \in \fL
\\[9pt]
\implies \quad
&q \cdot 2 \pi\hsx \sqrt{-1}  \in \fL
\\[9pt]
\implies \quad
&\ov{2 \pi\hsx \sqrt{-1}} \ = \  -2 \pi\hsx \sqrt{-1}  \in \fL
\noindent 
\end{align*}
Then
\[
\det
\begin{pmatrix}
2 \pi\hsx \sqrt{-1} \hsx w && q \cdot 2 \pi\hsx \sqrt{-1}\\[8pt]
2 \pi\hsx \sqrt{-1} && -2 \pi\hsx \sqrt{-1} \hsx \ov{w}
\end{pmatrix}
\ = \ 
4 \pi^2 w \ov{w} - 4 \pi^2q 
\ = \ 0
\]
contradicts 4EC.]\\[2pt]

[Note: The rows and columns of the matrix are $\Q$-linearly independent.  E.g. Suppose
\[
m \cdot 2 \pi\hsx \sqrt{-1} \hsx w + n \cdot q \cdot 2 \pi\hsx \sqrt{-1} = 0.
\]
Then
\[
m \cdot w + n \cdot q = 0 \qquad (m, \ n \in \Z \quad \text{($\Z$ is sufficient)})
\]
\qquad $\implies$
\[
m = 0, \quad n = 0, \qquad (w \notin \Q).
\]
\end{x}
\vspace{0.3cm}

\begin{x}{\small\bf EXAMPLE} \  (Admit 4EC) \ 
Let $w \in \C$.  
Assume: \ $\abs{w} \in \Q$ and $\exp(2 \hsx \pi \hsx \sqrt{-1} \hsx w)$ algebraic $-$then $w \in \Q$.
\vspace{0.2cm}

[In fact, 
\[
\abs{w} \in \Q 
\implies 
\abs{w}^2 \in \Q,
\]
so if $w \in \C - \Q$, then 
\[
\exp(2 \hsx \pi \hsx \sqrt{-1} \hsx w)
\]
is transcendental (cf. \#12).]
\end{x}
\vspace{0.3cm}

\begin{x}{\small\bf REMARK} \  (Admit 4EC) \ 
The 
\un{Diaz curve}
\index{Diaz curve} 
is the set of points
\[
\exp(2 \hsx \pi \hsx \sqrt{-1} \hsx w) \qquad (\abs{w} \ = \ 1).
\]
If $w = \pm 1$, then 
\[
\exp(2 \hsx \pi \hsx \sqrt{-1} \hsx w)
\]
is algebraic.  Otherwise
\[
\exp(2 \hsx \pi \hsx \sqrt{-1} \hsx w)
\]
is transcendental.
\vspace{0.3cm}
\end{x}

Here is one situation where the 4EC can be verified.
\vspace{0.3cm}

\begin{x}{\small\bf THEOREM} \ 
Suppose that $x_1$, $x_2$ are elements of $\R \hsx \cup \hsx \sqrt{-1} \hsx \R$ which are $\Q$-linearly independent and suppose that $y$ is a nonreal complex number with irrational real part $-$then at least one of the numbers
\[
e^{x_1}, \ 
e^{x_1 y}, \ 
e^{x_2}, \ 
e^{x_2 y}
\]
is transcendental.
\vspace{0.2cm}

[Note: \ 
In the notation of \#1, $y_1 = 1$, $y_2 = y$. ]
\vspace{0.2cm}

Proceed in steps.
\vspace{0.4cm}

\qquad \textbullet \  
The set $\{1, \hsy y, \hsy \ov{y}\}$ is $\Q$-linearly independent.
\vspace{0.2cm}

[Consider a rational dependence relation
\[ 
a + by + c\ov{y} \ = \ 0.
\]
Then
\[
\begin{cases}
\ a + (b + c) \hsy \Rex  y \ = \ 0\\[8pt]
\ \quad (b - c) \hsy \Imx y \ = \ 0
\end{cases}
.
\]
Since $y$ is nonreal, $\Imx y \ne 0$, hence 
\[
b - c \ = \ 0 
\implies 
b \ = \ c 
\implies
a + 2 b (\Rex y) \ = \ 0 
\implies
a \ = \ 0, \ b \ = \ 0.]
\]

\qquad \textbullet \quad 
Apply the six exponentials theorem to $\{x_1, x_2\}$ and $\{1, \hsy y,  \hsy \ov{y}\}$ (cf. \S41, \#1).  
\\[-.2cm]

\noindent Therefore at least one of the six numbers
\[
e^{x_1}, \ 
e^{x_1 y}, \ 
e^{x_1 \ov{y}}, \ 
e^{x_2}, \ 
e^{x_2 y}, \
e^{x_2 \ov{y}}
\]
is transcendental.
\vspace{0.4cm}

\qquad \textbullet \quad
By hypothesis, 
\[
\ov{x}_1 \ = \ \varepsilon_1 x_1, 
\quad 
\ov{x}_2 \ = \ \varepsilon_2 x_2
\quad
(\varepsilon_1, \ \varepsilon_2 \in \{1, -1\}),
\]
so 
\[
e^{x_1 \ov{y}} \ = \ e^{\ov{\varepsilon_1 x_1 y}},
\quad 
e^{x_2 \ov{y}} \ = \ e^{\ov{\varepsilon_2 x_2 y}}.
\]
Therefore at least one of the numbers
\[
e^{x_1}, \ 
e^{x_1 y}, \ 
e^{x_2}, \ 
e^{x_2 y}
\]
is transcendental.
\vspace{0.2cm}

[Note: 
If $e^{x_1 y}$ (or $e^{x_2 y}$) were algebraic, then the same would be true of 
$e^{x_1 \ov{y}}$ (or $e^{x_2 \ov{y}}$).]
\end{x}
\vspace{0.3cm}

%% file: _45_strong_four_exp_conjecture_s4ec.tex
\chapter{
$\boldsymbol{\S}$\textbf{45}.\quad  STRONG FOUR EXPONENTIALS CONJECTURE (S4EC)}
\setlength\parindent{2em}
\setcounter{theoremn}{0}
\renewcommand{\thepage}{\S45-\arabic{page}}

\ \indent 
This is the following statement.
\vspace{0.5cm}

\begin{x}{\small\bf CONJECTURE} \ 
Let $\{x_1, x_2\}$ and $\{y_1, y_2\}$  be two $\Qbar$-linearly independent sets of complex numbers $-$then
\[
\{x_1 y_1, x_1 y_2, x_2 y_1, x_2 y_2\} \not\subset \fL^*.
\]
\end{x}
\vspace{0.3cm}

In terms of matrices (cf. \S44, \#2):
\vspace{0.3cm}

\begin{x}{\small\bf CONJECTURE} \ 
Consider a nonzero $2 \times 2$ matrix $\sM$ with entries in $\fL^*$:
\[
\sM \ = \ 
\begin{pmatrix}
\lambda_{1 1} &&\lambda_{1 2}\\
\\
\lambda_{2 1} &&\lambda_{2 2}
\end{pmatrix}
.
\]
Suppose that its rows are $\Qbar$-linearly independent and its columns are  $\Qbar$-linearly independent $-$then
\[
\rank \sM \ = \ 2.
\]
\end{x}
\vspace{0.3cm}

\begin{x}{\small\bf IMPLICATION} \ 
\[
\text{strong four exponentials $\implies$ four exponentials }.
\]
\end{x}
\vspace{0.3cm}

\begin{x}{\small\bf CONDITION} \ 
PQ \quad
Let $\lambda_0$, $\lambda_1$, $\lambda_2 \in \fL^* - \{0\}$.  
Assume: 
\[
\lambda_1 / \lambda_0 \notin \Qbar
\quad \text{and} \quad 
\lambda_2 / \lambda_0 \notin \Qbar.
\]
Then
\[
(\lambda_1 \lambda_2) / \lambda_0 \notin \fL^*.
\]
\end{x}
\vspace{0.3cm}

\begin{x}{\small\bf LEMMA} \ 
\[
\text{S4EC $\Leftrightarrow$ PQ}.
\]
\vspace{0.2cm}

PROOF 
\vspace{0.2cm} 

\qquad \textbullet \quad 
S4EC $\implies$ PQ.
\vspace{0.2cm}

[In \#1, take
\[
\begin{cases}
\ x_1 = \lambda_0\\[8pt]
\ x_2 = \lambda_2
\end{cases}
, \quad
\begin{cases}
\ y_1 = 1\\
\ y_2 = \lambda_1 / \lambda_0
\end{cases}
\]
to arrive at
\[
\lambda_0, \lambda_1, \lambda_2, (\lambda_1 \lambda_2) / \lambda_0.
\]
But $\lambda_0$, $\lambda_1$, $\lambda_2 \in \fL^* - \{0\}$, thus it must be the case that
\[
(\lambda_1 \lambda_2) / \lambda_0 \notin \fL^*.]
\]

\qquad \textbullet \quad 
PQ $\implies$ S4EC.
\vspace{0.2cm}

[Start wtih $\{x_1, x_2\}$ and $\{y_1, y_2\}$  $\Qbar$-linearly independent sets of complex numbers.  
Assume that
\[
x_1 y_1, \ 
x_1 y_2, \ 
x_2 y_2
\]
are in $\fL^*$ and then claim that $x_2 y_1 \notin \fL^*$.  
Put 
\[
\lambda_0  \hsx = \hsx x_1 y_2, \quad 
\lambda_1 \hsx = \hsx x_1 y_1, \quad
\lambda_2 \hsx = \hsx x_2 y_2
\]
which, by hypothesis, are in $\fL^* - \{0\}$.  
Since
\[
\lambda_1 / \lambda_0 = y_1 / y_2 \notin \Qbar, \ 
\lambda_2 / \lambda_0 = x_2 / x_1 \notin \Qbar,
\]
it follows that
\[
(\lambda_1 \lambda_2) / \lambda_0 \ = \ x_2 y_1 \notin \fL^*.]
\]
\end{x}
\vspace{0.3cm}

\begin{x}{\small\bf APPLICATION} \ 
(Admit S4EC) \ 
Let $\lambda_1$, $\lambda_2 \in \fL^* - \Qbar$ \  $-$then $\lambda_1 \lambda_2 \notin \fL^*$.
\vspace{0.2cm}

[In \#4 above, take $\lambda_0 = 1$.]
\end{x}
\vspace{0.3cm}

\begin{x}{\small\bf \un{N.B.}} \ 
So in particular, if $\lambda_1$, $\lambda_2 \in \fL - \{0\}$, then $\lambda_1 \lambda_2 \notin \fL^*$, hence
\[
\lambda_1 \lambda_2 \notin \Qbar
\quad \text{and} \quad
\lambda_1 \lambda_2 \notin \fL.
\]
\\[-1.2cm]

[Note: \ Bear in mind that $\fL \hsx \cap \hsx \Qbar = \{0\}$.]
\end{x}
\vspace{0.3cm}

\begin{x}{\small\bf EXAMPLE} \  
(Admit S4EC) \ \ 
$e^{\pi^2}$ is transcendental (cf. \S43, \#6).
\vspace{0.2cm}

[In \#7, take
\[
\lambda_1 \ = \ \lambda_2 \ \equiv \ \lambda \ = \ \pi \sqrt{-1}.
\]
Then
\[
\lambda^2 \ = \ -\pi^2 \notin \fL^* 
\implies
\pi^2 \notin \fL^*.
\]
Therefore $e^{\pi^2}$ is transcendental.]
\end{x}
\vspace{0.3cm}

\begin{x}{\small\bf THEOREM} \ 
(Admit S4EC) \ \ 
If $\lambda \in \fL$ is nonzero, then $\abs{\lambda}$ is transcendental.  
\vspace{0.2cm}

PROOF \ 
In \#7, take $\lambda_1 = \lambda$, $\lambda_2 = \ov{\lambda}$, thus
\[
\lambda_1 \lambda_2 
\ = \ 
\lambda \hsx \ov{\lambda} 
\ = \ 
\abs{\lambda}^2 \notin \fL^*,
\]
thus $\abs{\lambda}^2$ is transcendental, thus $\abs{\lambda}$ is transcendental (if $\abs{\lambda}$ were algebraic, then 
 $\abs{\lambda}^2$ would be algebraic).
\end{x}
\vspace{0.3cm}

\begin{x}{\small\bf EXAMPLE} \ 
(Admit S4EC) \ 
Take
\[
\lambda \ = \ \elln(2) + \sqrt{-1} \hsx \pi.
\]
Then $\lambda \in \fL$ and
\[
\abs{\lambda} \ = \ \big(\elln(2)^2 + \pi^2\big)^{1/2}
\]
is transcendental.
\end{x}
\vspace{0.3cm}

\begin{x}{\small\bf THEOREM} \ 
(Admit S4EC) \ 
Let $w \in \C - \{0\}$.  
Assume: \ $\abs{w}$ is algebraic $-$then $e^w$ is transcendental (cf. \S44, \#7).
\vspace{0.3cm}

[In \#1, take
\[
\begin{cases}
\ x_1 = 1\\[8pt]
\ x_2 = e^w
\end{cases}
, \quad
\begin{cases}
\ y_1 = 1\\
\ y_2 = e^{\ov{w}}
\end{cases}
.
\]
Then
\[
x_1 y_1 = 1,\ 
x_1 y_2 = e^{\ov{w}}, \ 
x_2 y_1 = e^w, \ 
x_2 y_2 = e^w e^{\ov{w}}.
\]
\\[-1cm]

\qquad \textbullet \quad
 $\{x_1, x_2\}$,  $\{y_1, y_2\}$ are $\Qbar$-linearly independent.]
\vspace{0.3cm}

[To deal with $\{x_1, x_2\}$, suppose that
\[
\alpha + \beta e^w \ = \ 0 \qquad (\alpha, \beta \in \Qbar).
\]  
Then $\beta = 0$ 
\[
\implies \alpha = 0. \hspace{3.8cm}
\]
And $\beta \neq 0$
\allowdisplaybreaks
\begin{align*}
&\implies
e^w = -\frac{\alpha}{\beta} \in \Qbar - \{0\}
\\[12pt]
&\implies
w \in \fL 
\\[12pt]
&\implies
\abs{w} \quad \text{transcendental (cf. \#9)},
\end{align*}
contrary to the assumption that $\abs{w}$ is algebraic.  
Therefore $\beta$ must be zero, as must $\alpha$.]
\vspace{0.5cm}

Consider now the relation
\[
\{1, e^{\ov{w}}, e^w, e^w e^{\ov{w}}\} \not\subset \fL^*.
\]
If $e^w$ was algebraic, then the same would be true of $e^{\ov{w}}$ and $e^w e^{\ov{w}}$, an impossibility.
\vspace{0.2cm}

[Note: \ 
One can proceed without S4EC when 
\[
w \in \R \cup \sqrt{-1} \hsx \R \qquad (w \neq 0).
\]
For in this situation,
\[
\begin{cases}
\ \abs{w} = \pm w \quad (w \in \R)\\[8pt]
\ \abs{w} = \pm \sqrt{-1} \hsx w \quad (w \in \sqrt{-1} \hsx \R)
\end{cases}
.
\]
Therefore
\[
w \in \Qbar - \{0\} \implies e^w 
\quad 
\text{transcendental (Hermite-Lindemann (\S21, \#4))}.]
\]
\end{x}
\vspace{0.3cm}

\begin{x}{\small\bf LEMMA} \ 
(Admit S4EC) \ 
Let $\lambda \in \fL^*$.  
Assume: \ $\{\lambda, \ov{\lambda}\}$ is $\Qbar$-linearly independent 
$-$then $\abs{\lambda} \notin \fL^*$.  
\vspace{0.2cm}

PROOF \ \ 
We shall utilize condition PQ.  \ \ 
To this end, note that \ $\{\lambda, \abs{\lambda}\}$ \ is also $\Qbar$-linearly independent: 
\begin{align*}
\abs{\lambda} 
\ = \ 
\alpha \lambda 
\ 
(\alpha \in \Qbar) 
\ 
&\implies \ 
\abs{\lambda}^2 
\ = \ 
\alpha^2 \lambda^2 
\\[12pt]
&\implies \ 
\lambda \hsx \ov{\lambda} 
\ = \ 
\alpha^2 \lambda^2 
\\[12pt]
&\implies \ 
\ov{\lambda} = \alpha^2 \lambda.
\end{align*}
Supposing that $\abs{\lambda} \notin \fL^*$, take in \#4
\[
\lambda_0 = \lambda, \ 
\lambda_1 = \lambda_2 = \abs{\lambda}.
\]
Then
\[
\lambda_1/  \lambda_0 \notin \Qbar
\quad \text{and} \quad
\lambda_2 / \lambda_0 \notin \Qbar
\]
\qquad\qquad $\implies$
\[
(\lambda_1 \lambda_2) / \lambda_0 \notin \fL^*.
\]
On the other hand, 
\[
(\lambda_1 \lambda_2) / \lambda_0 \ = \ \ov{\lambda} \in \fL^*.
\]
Contradiction.
\end{x}
\vspace{0.3cm}

\begin{x}{\small\bf LEMMA} \ 
(Admit S4EC) \ 
\vspace{0.5cm}

\qquad \textbullet \quad 
If $\lambda \in \fL^* - \Qbar$, then the quotient $1 / \lambda$ is not in $\fL^*$.
\vspace{0.3cm}

\qquad \textbullet \quad 
If $\lambda_1$, $\lambda_2 \in \fL^* - \Qbar$, then the product $\lambda_1 \lambda_2$ is not in $\fL^*$.
\end{x}
\vspace{0.3cm}

\[
\textbf{APPENDIX}
\]
\vspace{0.5cm}

Let $\lambda \in \fL - \{0\}$ and let $w \in \C - \{0\}$ with $\abs{w} \in \Qbar$.  
Assume: \ $e^{\lambda w}$ is algebraic $-$then either $w \in \Q$ or else $w \lambda / \ov{\lambda} \in \Q$.
\vspace{0.2cm}

[Note: \ 
Tacitly S4EC is in force.]


%% file: _46_transcendental_extensions.tex
\chapter{
$\boldsymbol{\S}$\textbf{46}.\quad  TRANSCENDENTAL EXTENSIONS}
\setlength\parindent{2em}
\setcounter{theoremn}{0}
\renewcommand{\thepage}{\S46-\arabic{page}}


\begin{x}{\small\bf NOTATION} \ 
Let $\K$ be a field $-$then the field $\K(X_1, \ldots, X_n)$ of rational functions in $X_1, \ldots, X_n$ is the quotient field of the polynomial ring $\K[X_1, \ldots, X_n]$, hence consists of all quotients
\[
f(X_1, \ldots, X_n) / g(X_1, \ldots, X_n)
\]
of polynomials in $X_1, \ldots, X_n$ with $g \neq 0$.
\end{x}
\vspace{0.3cm}

Let $\LL$ be a field, $\K \subset \LL$ a subfield.
\vspace{0.3cm}

\begin{x}{\small\bf NOTATION} \ 
Fix a subset $S \subset \LL$.
\vspace{0.3cm}

\qquad \textbullet \quad
The ring $\K[S]$ generated by $\K$ and \mS is the intersection of all subrings of $\LL$ that contain $\K$ and \mS.
\vspace{0.2cm}

\qquad \textbullet \quad
The field $\K(S)$ generated by $\K$ and \mS is the intersection of all subfields of $\LL$ that contain $\K$ and \mS.
\vspace{0.2cm}

[Note: \ If $S = \{\alpha_1, \ldots, \alpha_n\}$ is finite, write
\[
\K[S] \ = \ \K[\alpha_1, \ldots, \alpha_n]
\]
and 
\[
\K(S)\ = \ \K(\alpha_1, \ldots, \alpha_n).]
\]
\end{x}
\vspace{0.3cm}

\begin{x}{\small\bf \un{N.B.}} \ 
If \mS is finite, then the field $\K(S)$ is said to be a 
\un{finitely generated} \un{extension}
\index{finitely generated extension} 
of $\K$.
\vspace{0.2cm}

[Note: 
\allowdisplaybreaks
\begin{align*}
&\text{finite extension $\implies$ finitely generated extension}
\\[12pt]
&\text{finitely generated extension $\centernot\implies$ finite extension}.]
\end{align*}
\end{x}
\vspace{0.3cm}

\begin{x}{\small\bf LEMMA} \ 
$\K(S)$ is the set of all elements of $\LL$ that can be expressed as quotients of finite linear combinations with coefficients in $\K$ of finite products of elements of \mS.
\end{x}
\vspace{0.3cm}

\begin{spacing}{1.45}
\begin{x}{\small\bf TERMINOLOGY} \ 
Let $\LL$ be a field, $\K \subset \LL$ a subfield.
\vspace{0.3cm}

\qquad \textbullet \quad
A finite subset $S = \{\alpha_1, \ldots, \alpha_n\} \subset \LL$ is 
\un{algebraically dependent over $\K$}
\index{algebraically dependent over $\K$} 
if there is a nonzero polynomial $P \in \K[X_1, \ldots, X_n]$ such that
\[
P(\alpha_1, \ldots, \alpha_n) \ = \ 0.
\]

\qquad \textbullet \quad
A finite subset $S = \{\alpha_1, \ldots, \alpha_n\} \subset \LL$ is 
\un{algebraically independent over $\K$}
\index{algebraically independent over $\K$} 
if there is no nonzero polynomial $P \in \K[X_1, \ldots, X_n]$ such that
\[
P(\alpha_1, \ldots, \alpha_n) \ = \ 0.
\]
\end{x}
\vspace{0.3cm}

\begin{x}{\small\bf EXAMPLE} \ 
Take $\LL = \K(X_1, \ldots, X_n)$, the field of rational functions in $X_1, \ldots, X_n$ $-$then $\{X_1, \ldots, X_n\}$ is algebraically independent over $\K$.
\vspace{0.2cm}

[Note: \ 
Suppose that $r_1, \ldots, r_n$ are positive integers $-$then $\{X_1^{r_1}, \ldots, X_n^{r_n}\}$ is algebraically independent over $\K$.]
\end{x}
\vspace{0.3cm}

\begin{x}{\small\bf EXAMPLE} \ 
Working still with $\LL = \K(X_1, \ldots, X_n)$, let $A = [a_{i j}]$ be an $n \times n$ matrix with coefficients in $\K$.  
Put $f_j = \sum\limits_i \hsx a_{i j} X_i$ $-$then $\{f_1, \ldots, f_n\}$ is algebraically independent over $\K$ iff $\det A \neq 0$.
\end{x}
\vspace{0.3cm}

\begin{x}{\small\bf \un{N.B.}} \ 
Take $S = \emptyset$, the empty set $-$then it is deemed to be algebraically independent over $\K$.
\end{x}
\vspace{0.3cm}

\begin{x}{\small\bf LEMMA} \ 
If $\alpha_1, \ldots, \alpha_n \in \LL$ are algebraically independent over $\K$, then
$\K[\alpha_1, \ldots, \alpha_n]$ and $\K[X_1, \ldots, X_n]$ are $\K$-isomorphic rings, hence $\K(\alpha_1, \ldots, \alpha_n)$ and 
$\K(X_1, \ldots, X_n)$ are $\K$-isomorphic fields.
\vspace{0.2cm}

[Note: \ 
The property is characteristic in that if $\K(\alpha_1, \ldots, \alpha_n)$ and $\K(X_1, \ldots, X_n)$  are $\K$-isomorphic fields, then 
$\{\alpha_1, \ldots, \alpha_n\}$ is
algebraically independent over $\K$.]
\end{x}
\vspace{0.3cm}

\begin{x}{\small\bf REMARK} \ 
The algebraic independence of $\alpha_1, \ldots, \alpha_n \in \LL$ over $\K$ is equivalent to the requirement that for each $i$, 
$\alpha_i$ is transcendental over $\K(\alpha_1, \ldots, \alpha_{i-1})$.
\end{x}
\vspace{0.3cm}

\begin{x}{\small\bf DEFINITION} \ 
A subset \mS of $\LL$ is a 
\un{transcendence basis}\index{transcendence basis} 
for $\LL/\K$ if \mS is algebraically independent over $\K$ and if $\LL$ is algebraic over $\K(S)$.
\vspace{0.2cm}

[Note: \ 
A priori, \mS is infinite, the convention being that \mS is algebraically independent over $\K$ if every finite subset of \mS is algebraically independent over $\K$.] 
\end{x}
\vspace{0.3cm}

\begin{x}{\small\bf EXAMPLE} \ 
In the setup of \#6, $\{X_1^{r_1}, \ldots, X_n^{r_n}\}$ is algebraically independent over $\K$.  
So, to establish that  $\{X_1^{r_1}, \ldots, X_n^{r_n}\}$ is a transcendence basis for $\LL / \K$, it has to be shown that $\LL$ is algebraic over $\K(X_1^{r_1}, \ldots, X_n^{r_n})$.  
But for each $i$, the element $X_i$ is a zero of the polynomial $T^{r_i} - X_i^{r_i} \in \LL[T]$.
\end{x}
\vspace{0.3cm}

\begin{x}{\small\bf \un{N.B.}} \ 
If $S = \emptyset$ is a transcendence basis for $\LL / \K$, then $\LL/ \K$ is algebraic (and conversely).
\end{x}
\vspace{0.3cm}

\begin{x}{\small\bf THEOREM} \ 
There exists a transcendence basis for $\LL / \K$.
\end{x}
\vspace{0.3cm}

\begin{x}{\small\bf REMARK} \ 
If $S_1 \subset S_2 \subset \LL$, if $S_1$ is algebraically independent over $\K$, if $\LL/\K(S_2)$ is algebraic, then there exists a transcendence basis \mX for $\LL/\K$ with $S_1 \subset X \subset S_2$.
\end{x}
\vspace{0.3cm}

\begin{x}{\small\bf THEOREM} \ 
If $S_1 \subset \LL$, $S_2 \subset \LL$ are transcendence bases for $\LL/\K$, then 
\[
\card S_1 \ = \ \card S_2.
\]
\end{x}
\vspace{0.3cm}

\begin{x}{\small\bf DEFINITION} \ 
The 
\un{transcendence degree}
\index{transcendence degree} 
\[
\trdeg_{\K} (\LL / \K) 
\]
is the cardinality of any transcendence basis for $\LL/\K$.
\end{x}
\vspace{0.3cm}

\begin{x}{\small\bf \un{N.B.}} \ 
If
\[
\trdeg_{\K} (\LL / \K)  \ = \ 0,
\]
then $\LL/\K$ is algebraic (and conversely).
\end{x}
\vspace{0.3cm}

\begin{x}{\small\bf EXAMPLE} \ 
Take $\K = \Q$, $\LL = \C$ $-$then
\[
\trdeg_{\Q} (\C/\Q) \ = \ \mfc.
\]
\end{x}
\vspace{0.3cm}

\begin{x}{\small\bf THEOREM} \ 
Let $\bk \subset \K \subset \LL$ be fields $-$then 
\[
\trdeg_{\bk}(\LL/\bk) \ = \ \trdeg_{\K}(\LL/\K) + \trdeg_{\bk}(\K/\bk).
\]
\end{x}
\vspace{0.3cm}

The situtation when $\LL$ is a finitely generated extension of $\K$ occupies center stage.
\vspace{0.3cm}

\begin{x}{\small\bf SCHOLIUM} \ 
Let $\LL = \K(\alpha_1, \ldots, \alpha_n)$ $-$then a maximal algebraically 
independent subset of the set $\{\alpha_1, \ldots, \alpha_n\}$ is a transcendence basis for $\LL/\K$ and 
\[
\trdeg_{\K} (\LL/\K) \ \leq \ n.
\]
Assuming that $S = \{\alpha_1, \ldots, \alpha_m\}$, it follows that $\LL$ is a finite extension of $\K(\alpha_1, \ldots, \alpha_m)$ and if this is separable (which is always the case in characteristic 0), then
\[
\LL \ = \ \K(\alpha_1, \ldots, \alpha_m,\beta)
\]
for some $\beta$ in $\LL$ (primitive element).
\vspace{0.2cm}

[Note: \ 
The extension $\LL/\K$ can be broken up into a series of subextensions, viz. 
let $\K_i = \K(\alpha_1, \ldots, \alpha_i)$ (put $\K_0 = \K$) $-$then
\[
\K \ = \ \K_0 \hsx \subset \hsx \K_1 \hsx \subset \hsx \K_2 \hsx \subset \hsx \cdots \hsx \subset \hsx \K_n \ = \ \LL,
\]
where $\K_{i+1} = \K_i(\alpha_{i+1})$.]
\end{x}
\vspace{0.3cm}

\begin{x}{\small\bf LEMMA} \ 
Let $\LL$ be a field, $\K \subset \LL$ a subfield.  
Let \mS be a subset of $\LL$ with the property that each $\alpha \in S$ is algebraic over $\K$ $-$then $\K(S)$ is algebraic over $\K$ and
\[
\text{\mS finite $\implies [\K(S):\K]$ \ finite}.
\]
\end{x}
\vspace{0.3cm}

\begin{x}{\small\bf EXAMPLE} \ 
Take $\K = \Q$ and consider $\Q(\sqrt{2},\pi)$ $-$then it is clear that $\{\sqrt{2}\}$ is not algebraically independent, nor is 
$\{\sqrt{2}, \pi\}$, which leaves $\{\pi\}$, the claim being that it is a transcendence basis for $\Q(\sqrt{2},\pi)/\Q$ 
(per the theory spelled out in \#21).  
To check this, in \#22 take $\K = \Q(\pi)$, $\LL = \Q(\sqrt{2},\pi)$, $S = \{\sqrt{2}, \pi\}$.  
\vspace{0.3cm}

\qquad \textbullet \quad
$\sqrt{2}$ is algebraic over $\Q(\pi)$: \ Work with $X^2 - 2 \in \Q(\pi) [X]$.
\vspace{0.2cm}

\qquad \textbullet \quad
$\pi$ is algebraic over $\Q(\pi)$: \ Work with $X - \pi \in \Q(\pi) [X]$.
\vspace{0.2cm}

Therefore $\Q(\pi)(\sqrt{2}, \pi)$ is algebraic over $\Q(\pi)$.\\
\noindent And
\[
\trdegQ \Q(\sqrt{2}, \pi) \ = \ 1.
\]
\end{x}
\vspace{0.3cm}

\begin{x}{\small\bf REMARK} \ 
The transcendence degree
\[
\trdegQ \Q(\pi, e) 
\]
is either 1 or 2 but whether it is 1 or whether it is 2 is unknown since it is not known if $\pi$ and $e$ are algebraically independent or not.
\end{x}
\end{spacing}
\vspace{0.3cm}

\begin{x}{\small\bf RATIONAL RECAPITULATION} \ 
Let \mM and \mN be finite subsets of $\C$.
\vspace{0.3cm}

\qquad \textbullet \quad
If $N \subset \Qbar$, then
\[
\trdegQ \Q(M \cup N) 
\ = \ 
\trdegQ \Q(M).
\]
Therefore algebraic numbers do not contribute to the transcendence degree.
\vspace{0.3cm}

\qquad \textbullet \quad
If $N \subset M$, then
\[
\trdegQ \Q(M \cup N) 
\ = \ 
\trdegQ \Q(M).
\]
Therefore only distinct numbers can contribute to the transcendence degree.
\vspace{0.3cm}

\qquad \textbullet \quad
If the transcendence degree
\[
\trdegQ \Q(M) 
\]
of the field $\Q(M)$ is $\card M$, then \mM is algebraically independent over $\Q$ and conversely.
\vspace{0.3cm}

\qquad \textbullet \quad
If $M = \{m\}$, then the transcendence degree
\[
\trdegQ \Q(m) 
\]
of the field $\Q(m)$ is 0 if $m$ is algebraic and 1 if $m$ is transcendental.
\vspace{0.3cm}

\qquad \textbullet \quad
$\Q \hsx \ldots \hsx \Qbar$:
\[
\trdegQ \Q(M) \ = \ \trdegQbar \Qbar(M).
\]
\vspace{0.2cm}
\end{x}
\vspace{0.3cm}

\begin{spacing}{1.5}
\begin{x}{\small\bf LEMMA} \ 
Suppose that $\alpha_1, \ldots, \alpha_n$ are algebraically independent over $\K$ $-$then so are 
$\alpha_1^{p_1 / q_1}, \ldots, \alpha_n^{p_n / q_n}$ for nonzero rational numbers $p_1/q_1, \ldots, p_n/q_n$.
\vspace{0.2cm}

PROOF \ 
The transcendence degree of $\K(\alpha_1, \ldots, \alpha_n)$ over $\K$ is $n$ (cf. \#9), whereas
\[
\K\big(\alpha_1^{1/q_1}, \ldots, \alpha_n^{1/q_n}\big)
\]
is algebraic over $\K(\alpha_1, \ldots, \alpha_n)$ since $\bigg(\alpha_j^{1/q_j}\bigg)^{q_j} = \alpha_j$.  
Therefore the transcendence degree of 
\[
\K\big(\alpha_1^{1/q_1}, \ldots, \alpha_n^{1/q_n}\big)
\]
over $\K$ is also $n$.  
The numbers $\big\{\alpha_1^{1/q_1}, \ldots, \alpha_n^{1/q_n}\big\}$ are algebraically independent over $\K$, thus the same is true of the numbers $\big\{\alpha_1^{p_1 / q_1}, \ldots, \alpha_n^{p_n / q_n}\big\}$ (cf. \#6).
\end{x}
\end{spacing}
\vspace{0.3cm}

\begin{x}{\small\bf LEMMA} \ 
Suppose that $\alpha_1, \ldots, \alpha_n$ are algebraically independent over $\K$.  Let
\[
\frac{A[X_1, \ldots, X_n]}{B[X_1, \ldots, X_n]}
\]
be two nonzero polynomials whose quotient is not in $\K$ $-$then
\[
\frac{A(\alpha_1, \ldots, \alpha_n)}{B(\alpha_1, \ldots, \alpha_n)}
\]
is not in $\K$.
\vspace{0.2cm}

PROOF \ 
If the ratio was equal to some $\alpha \in \K$, then
\[
A(\alpha_1, \ldots, \alpha_n) \ - \ \alpha B(\alpha_1, \ldots, \alpha_n) \ = \ 0,
\]
which contradicts the algebraic independence of the $\alpha_j$'s.
\end{x}
\vspace{0.3cm}


%% file: _47_schanuels_conjecture.tex
\chapter{
$\boldsymbol{\S}$\textbf{47}.\quad  SCHANUEL'S CONJECTURE (SCHC)}
\setlength\parindent{2em}
\setcounter{theoremn}{0}
\renewcommand{\thepage}{\S47-\arabic{page}}

\ \indent 

This is the following statement.
\vspace{0.5cm}

\begin{x}{\small\bf CONJECTURE} \ 
Suppose that $x_1, \ldots, x_n$ are $\Q$-linearly independent complex numbers $-$then among the $2n$ numbers
\[
x_1, \ldots, x_n, e^{x_1}, \ldots, e^{x_n}, 
\]
at least $n$ are algebraically independent over $\Q$, i.e., 
\[
\trdegQ \Q(x_1, \ldots, x_n, e^{x_1}, \ldots, e^{x_n}) \ \geq \ n \qquad \text{(cf. \S46, \#21)}.
\]
\end{x}
\vspace{0.3cm}

This conjecture has many consequences, some of which are delineated below.
\vspace{0.3cm}

\begin{x}{\small\bf LEMMA} \ 
The set of $n$-tuples $(x_1, \ldots, x_n)$  in $\C^n$ such that the $2n$ numbers 
\[
x_1, \ldots, x_n, e^{x_1}, \ldots, e^{x_n}
\]
are algebraically independent over $\Q$ is a $G_\delta$-subset of $\C^n$ and its complement is a set of Lebesgue measure 0.
\end{x}
\vspace{0.3cm}

\begin{x}{\small\bf \un{N.B.}} \ 
The transcendence degree can be as small as $n$ (cf. \#6).
\end{x}
\vspace{0.3cm}

\begin{x}{\small\bf THEOREM} \ 
Take $n = 1$ and consider $x$, $e^x$ $(x \neq 0)$ $-$then at least one of $x$, $e^x$ is transcendental (cf. \S31, \#5), thus
\[
\trdegQ \Q(x, e^x) \ \geq \ 1,
\]
which is Schanuel in the simplest situation.
\end{x}
\vspace{0.3cm}


\begin{x}{\small\bf \un{N.B.}} \ 
Take $n = 2$ and consider  $x_1$, $x_2$, $e^{x_1}$, $e^{x_2}$ $-$then the claim is that
\[
\trdegQ \Q(x_1, x_2 , e^{x_1}, e^{x_2}) 
\ \geq \ 
2
\]
but this has never been verified in general.
\vspace{0.2cm}

[Note: \ 
Let $w_1$, $w_2$ be two nonzero complex numbers $-$then SCHC implies that
\[
\trdegQ \Q(w_1 w_2 , e^{w_1}, e^{w_2}) 
\ \geq \ 
1.]
\]
\end{x}
\vspace{0.3cm}

\begin{x}{\small\bf THEOREM} \ 
Suppose that $x_1, \ldots, x_n$ are $\Q$-linearly independent algebraic numbers $-$then the transcendental numbers 
$e^{x_1}, \ldots, e^{x_n}$ are algebraically independent over $\Q$ (cf. \S21, \#12), so 
\[
\trdegQ \Q(x_1, \ldots, x_n, e^{x_1}, \ldots, e^{x_n})
\ \geq \ 
n,
\]
thereby settling Schanuel in the particular case when $x_1, \ldots, x_n$ are algebraic.
\end{x}
\vspace{0.3cm}

\begin{x}{\small\bf THEOREM}  \ 
(Admit SCHC) \ 
Let $\lambda_1, \ldots, \lambda_n$ be $\Q$-linearly independent elements of $\fL$ (thus transcendental (cf. \S31, \#4)) $-$then 
$e^{\lambda_1}, \ldots, e^{\lambda_n}$ are algebraic numbers, hence
\allowdisplaybreaks
\begin{align*}
\trdegQ \Q(\lambda_1, \ldots, \lambda_n, e^{\lambda_1}, \ldots, e^{\lambda_n})
&=\ 
\trdegQ \Q(\lambda_1, \ldots, \lambda_n)
\\[12pt]
&\leq \ n.
\end{align*}
On the other hand, by Schanuel, 
\[
\trdegQ \Q(\lambda_1, \ldots, \lambda_n, e^{\lambda_1}, \ldots, e^{\lambda_n})
\ \geq \ 
n.
\]
Therefore
\[
\trdegQ \Q(\lambda_1, \ldots, \lambda_n) 
\ = \ n,
\]
which implies that $\{\lambda_1, \ldots, \lambda_n\}$ is algebraically independent over $\Q$ (cf. \S46, \#9).
\end{x}
\vspace{0.3cm}

\begin{x}{\small\bf EXAMPLE} \ 
It is not true in general that
\[
\text{linear independence}
\implies
\text{algebraic independence}.
\]
Thus, e.g., $\{1, \sqrt{2}, \sqrt{3}, \sqrt{6}\}$ is linearly independent over $\Q$ but is not algebraically independent over $\Q$ as can be seen by noting that if
\[
P(X_1, X_2, X_3, X_4) 
\ = \ 
X_2 X_3 - X_4,
\]
then
\[P(1, \sqrt{2}, \sqrt{3}, \sqrt{6}) \ = \ 0.
\]
\end{x}
\vspace{0.3cm}

\begin{spacing}{1.45}
\begin{x}{\small\bf IMPLICATION} \ 
\[
\text{Schanuel}
\implies
\text{inhomogeneous Baker}.
\]

[If $\lambda_1 \in \fL, \ldots, \lambda_n \in \fL$ are $\Q$-linearly independent, then $\lambda_1, \ldots, \lambda_n$ are $\Q$-algebraically independent (cf. \#7) or still, $\lambda_1, \ldots, \lambda_n$ are $\Qbar$-algebraically independent (cf. \S20, \#7), hence $1, \lambda_1, \ldots, \lambda_n$ are $\Qbar$-linearly independent.  
Proof: \ Given $\gamma, \gamma_1, \ldots, \gamma_n$ algebraic and 
\[
\gamma + \gamma_1 \lambda_1 + \cdots + \gamma_n \lambda_n \ = \ 0,
\]
work with
\[
P(X_1, \ldots, X_n) 
\ = \ 
\gamma + \gamma_1 X_1 + \cdots + \gamma_n X_n.]
\]
\end{x}
\vspace{0.3cm}

\begin{x}{\small\bf THEOREM}  \ 
(Admit SCHC) \ 
Suppose given elements $\lambda_1, \ldots, \lambda_n$ in $\fL$ 
and elements $\alpha_1, \ldots, \alpha_m$ in $\Qbar$.  
Assume: \ $\lambda_1, \ldots, \lambda_n$ are $\Q$-linearly independent and 
$\alpha_1, \ldots, \alpha_m$ are $\Q$-linearly independent $-$then
\[
\trdegQ \Q(\lambda_1, \ldots, \lambda_n, e^{\alpha_1}, \ldots, e^{\alpha_m}) 
\ = \ 
m + n,
\]
thus
\[
\{\lambda_1, \ldots, \lambda_n, e^{\alpha_1}, \ldots, e^{\alpha_m}\}
\]
is algebraically independent over $\Q$ (cf. \S46, \#9).
\vspace{0.2cm}

PROOF \ 
Define $\beta_j$: $j  = 1, \ldots, m + n$ by $\beta_j = \lambda_j$ for $j = 1, \ldots, n$ and 
$\beta_{j+n} = \alpha_j$ for $j = 1, \ldots,  m$.  Claim:
\[
\beta_1, \ldots, \beta_{m+n}
\]
is $\Q$-linearly independent.  For suppose that
\[
q_1 \beta_1 + \cdots + q_{m+n} \beta_{m+n} \ = \ 0
\]
is a rational dependence relation, hence
\[
q_1 \lambda_1 + \cdots + q_n \lambda_n + q_{n+1} \alpha_1 + \cdots + q_{m+n} \alpha_m \ = \ 0.
\]
From the definitions,
\[
q_{n+1} \alpha_1 + \cdots + q_{m+n} \alpha_m
\]
is an algebraic number, i.e., is in $\Qbar$.  Accordingly, thanks to inhomogeneous Baker,
\[
q_1 = 0, \ldots, q_n = 0, 
\quad \text{and} \quad
q_{n+1} \alpha_1 + \cdots + q_{m+n} \alpha_m \ = \ 0.
\]
But $\alpha_1, \ldots, \alpha_m$ are $\Q$-linearly independent.  Therefore
\[
q_{n+1} = 0, \ldots, q_{m+n} = 0, 
\]
hence the claim.  
Now apply Schanuel: \ The transcendence degree over $\Q$ of
\[
\Q\big(\beta_1, \ldots, \beta_{m+n}, e^{\beta_1}, \ldots, e^{\beta_{m+n}}\big)
\]
is $\hsx \geq \hsx m + n$.  
To cut this down, note that
\[
\beta_{1+n} = \alpha_1, \ldots, \beta_{m+n} = \alpha_m
\]
are algebraic, as are
\[
e^{\beta_1} = e^{\lambda_1}, \ldots, e^{\beta_n} = e^{\lambda_n}.
\]
So we are left with
\[
\trdegQ  \Q(\lambda_1, \ldots, \lambda_n, e^{\alpha_1}, \ldots, e^{\alpha_m})
\ \geq \ 
m + n,
\]
which suffices.
\end{x}
\end{spacing}
\vspace{0.3cm}

\begin{x}{\small\bf THEOREM} \ 
(Admit SCHC) \ \ 
If  $\alpha \neq 0, \ 1$ is algebraic and if $1, \beta_1, \ldots, \beta_n \in \Qbar$ are linearly independent over $\Q$, then the numbers $\Log \alpha$ and 
\[
\alpha^{\beta_1}, \ldots, \alpha^{\beta_n} \qquad \text{(principal powers)}
\]
are algebraically independent over $\Q$, hence are transcendental (cf. \S31, \#17).
\vspace{0.3cm}

PROOF \ 
To begin with,
\[
\beta_1 \hsx \Log \alpha, \ldots, \beta_n \hsx \Log \alpha, \Log \alpha
\]
are $\Q$-linearly independent, thus the transcendence degree of the field
\[
\Q(\beta_1 \hsx \Log \alpha, \ldots, \beta_n \hsx \Log \alpha, \Log \alpha, \alpha^{\beta_1}, \ldots, \alpha^{\beta_n}, \alpha)
\]
is $\geq n + 1$ (quote Schanuel).  But
\[
\beta_1 \ = \ (\beta_1 \hsx \Log \alpha) (\Log \alpha)^{-1}, \ldots
\]
\qquad\qquad $\implies$
\allowdisplaybreaks
\begin{align*}
\Q(\beta_1 \hsx \Log \alpha,& \ldots, \beta_n \hsx \Log \alpha, \Log \alpha, \alpha^{\beta_1}, \ldots, \alpha^{\beta_n}, \alpha)
\\[12pt]
&=\ 
\Q(\beta_1, \ldots, \beta_n, \Log \alpha, \alpha^{\beta_1}, \ldots, \alpha^{\beta_n}, \alpha)
\end{align*}
\qquad\qquad $\implies$
\allowdisplaybreaks
\begin{align*}
\trdegQ \Q(\beta_1, \ldots, \beta_n&,\Log \alpha, \alpha^{\beta_1}, \ldots, \alpha^{\beta_n}, \alpha)
\\[12pt]
&=\ 
\trdegQ \Q(\Log \alpha, \alpha^{\beta_1}, \ldots, \alpha^{\beta_n})
\\[12pt]
&\geq \ 
n + 1
\end{align*}
\qquad\qquad $\implies$
\[
\trdegQ \Q(\Log \alpha, \alpha^{\beta_1}, \ldots, \alpha^{\beta_n})
\ = \ 
n + 1,
\]
from which the algebraic independence over $\Q$ of $\Log \alpha$ and 
\[
\alpha^{\beta_1}, \ldots, \alpha^{\beta_n}.
\]
\end{x}
\vspace{0.3cm}

\begin{x}{\small\bf \un{N.B.}} \ 
In \#11, take $n = 1$ and assume that $\beta \notin \Q$ $-$then $\Log \alpha$ and $\alpha^\beta$ are algebraically independent over $\Q$.
\end{x}
\vspace{0.3cm}

\begin{x}{\small\bf THEOREM} \ 
(Admit SCHC) \ \ 
If $\alpha \neq 0, \hsy 1$ is algebraic and if $\beta \in \Qbar$ has degree $d \geq 2$, then 
\[
\trdegQ  \Q(\Log \alpha, \alpha^\beta, \ldots, \alpha^{\beta^{d-1}})
\ = \ 
d.
\]

\vspace{0.2cm}

PROOF \ 
First of all, $1, \beta, \ldots, \beta^{d-1}$ are linearly independent over $\Q$.  
In fact, the minimal polynomial of $\beta$ has degree $d \geq 2$, whereas a rational dependence relation
\[
q + q_1 \beta + \cdots + q_{d-1} \beta^{d-1} \ = \ 0
\]
leads to a contradiction upon consideration of
\[
P(X_0, X_1, \ldots, X_{d-1}) 
\ = \ 
X_0 + q_1 X_1 + \cdots + q_{d-1} X^{d-1}.
\]
So, applying \#11, the numbers $\Log \alpha$ and
\[
\alpha^\beta, \ldots, \alpha^{\beta^{d - 1}}
\qquad
\text{(principal powers)}
\]
are algebraically independent over $\Q$, from which the result.
\vspace{0.2cm}

[Note: \ 
It is not necessary to appeal to SCHC when $d = 2$ or $d = 3$ as these special cases have been resolved.  
For a case in point, take
\[
d = 3, \ \alpha = 2, \ \beta = 2^{1/3}.
\]
Then
\[
\elln(2), \ 2^{2^{1/3}}, \ 2^{2^{2/3}}
\]
are algebraically independent over $\Q$.] 
\end{x}
\vspace{0.3cm}

\begin{x}{\small\bf REMARK} \ 
It can be shown that unconditionally
\[
\trdegQ \Q\big(\alpha^\beta, \ldots, \alpha^{\beta^{d-1}}\big)
\ \geq \ 
\bigg[\frac{d+1}{2}\bigg],
\]
the symbol on the right standing for the greatest integer less than or equal to $\ds\frac{d+1}{2}$.

\end{x}
\vspace{0.3cm}

\begin{x}{\small\bf THEOREM} \ 
(Admit SCHC) \ 
If $x_1, \ldots, x_n$ are complex numbers linearly independent over $\Q$ and if $y$ is a transcendental number, then
\[
\trdegQ \Q\big(e^{x_1}, \ldots, e^{x_n}, e^{x_1 y}, \ldots, e^{x_n y} \big)
\ \geq \ 
n - 1.
\]
\\[-1.2cm]

PROOF \ 
Order the numbers $x_1, \ldots, x_n$ in such a way that a basis for the $\Q$-vector space generated by
\[
\{x_1, \ldots, x_n, x_1 y, \ldots, x_n y\}
\]
is
\[
\{x_1, \ldots, x_n, x_1 y, \ldots, x_m y\} \qquad (0 \leq m \leq n).
\]
Claim:
\[
\trdegQ \Q(x_1, \ldots, x_n,  y)
\ \leq \ 
m + 1.
\]
For $y$ is transcendental (by hypothesis), so there is a transcendence basis for
\[
\Q(x_1, \ldots, x_n,  y)
\]
which is
\[
\{x_{i_1}, \ldots, x_{i_k},  y\}
\]
with
\[
1 \leq i_1 < i_2 < \cdots < i_k \leq n.
\]
Then
\[
x_1, \ldots, x_n, x_{i_1} y, \ldots, x_{i_k} y
\]
are $\Q$-linearly independent, thus
\[
k + n \leq m + n 
\implies 
k \leq m 
\implies 
k + 1 \leq m + 1,
\]
which establishes the claim.  Next, invoking SCHC, 
\[
\trdegQ 
\Q(x_1, \ldots, x_n, x_1 y, \ldots, x_m y, e^{x_1}, \ldots, e^{x_n}, e^{x_1 y}, \ldots, e^{x_m y}) 
\ \geq \ 
n + m
\]
\qquad $\implies$
\[
\trdegQ 
\Q(x_1, \ldots, x_n, x_1 y, \ldots, x_n y, e^{x_1}, \ldots, e^{x_n}, e^{x_1 y}, \ldots, e^{x_n y}) 
\ \geq \ 
n + m.
\]
Taking into account the claim, it follows that at least $n-1$ of the numbers $e^{x_i}$, $e^{x_i y}$ $(i = 1, \ldots, n)$ are algebraically independent.
\end{x}
\vspace{0.3cm}

\begin{x}{\small\bf \un{N.B.}} \ 
Specialized to the case $n=2$, the upshot is that at least one 
of the numbers
\[
e^{x_1}, \ e^{x_2}, \ e^{x_1 y}, \ e^{x_2 y}
\]
is transcendental.
\end{x}
\vspace{0.3cm}

\begin{x}{\small\bf IMPLICATION} \ 
\[
\text{SCHC} \implies \text{4EC}.
\]
\end{x}
\vspace{0.3cm}

\begin{x}{\small\bf RAPPEL} \ 
(4EC) \ 
Let $\{x_1, x_2\}$ and $\{y_1, y_2\}$ be two $\Q$-linearly independent sets of complex numbers $-$then
\[
\{x_1 y_1, x_1 y_2, x_2 y_1, x_2 y_2\} \subset \fL,
\]
thus at least one of the numbers
\[
e^{x_1 y_1}, \ e^{x_1 y_2}, \ e^{x_2 y_1}, \ e^{x_2 y_2}
\]
it transcendental.

When dealing with 4EC, there is a little trick that can be used to advantage, viz. let
\[
w_1 = x_1 y_1, \ 
w_2 = x_2 y_1, \ 
z_1 = y_2 / y_1, \ 
z_2 = 1.
\]
Then
\[
w_1 z_1 = x_1 y_2, \ 
w_1 z_2 = x_1 y_1, \ 
w_2 z_1 = x_2 y_2, \ 
w_2 z_2 = x_2 y_1.
\]
So the list
\[
e^{x_1 y_1},\ e^{x_1 y_2},\ e^{x_2 y_1},\ e^{x_2 y_2}
\]
becomes the list
\[
e^{w_1 z_2},\ e^{w_1 z_1},\ e^{w_2 z_2},\ e^{w_2 z_1},
\]
i.e., the list
\[
e^{w_1},\ e^{w_1 z_1},\ e^{w_2},\ e^{w_2 z_1},
\]
i.e., the list
\[
e^{w_1},\ e^{w_2},\ e^{w_1 y},\ e^{w_2 y},
\]
where
\[
y \ = \ z_1.
\]

In order to utilize \#16, it is necessary that $y$ be transcendental.
\vspace{0.2cm}

\qquad\qquad \un{Case 1:} \quad
$y \notin \fL^*$ $-$then $y$ is transcendental (otherwise, $y$ would be algebraic, while $\Qbar \hsx \subset \hsx \fL^*$).
\vspace{0.2cm}

\qquad\qquad \un{Case 2:} \quad
$y \in \fL^*$ $-$then \#16 need not be applicable but in view of \S43, \#8, 
\[
\{x_1 y_1, x_1 y_2, x_2 y_1, x_2 y_2\} \hsx \not\subset \hsx \fL,
\]
thus at least one of the numbers
\[
e^{x_1 y_1}, \ e^{x_1 y_2}, \ e^{x_2 y_1}, \ e^{x_2 y_2}
\]
it transcendental.
\vspace{0.2cm}

[Note:  \ 
In the reference to \S43, \#8, take $\gamma = 1$ and replace $x_2 / x_1$ by $y_2 / y_1$ (as is certainly permissible).]
\end{x}
\vspace{0.3cm}

\begin{x}{\small\bf RAPPEL} \ 
(Admit S4EC) \ \ 
Let $w \in \C - \{0\}$.  
Assume: \ $\abs{w}$ is algebraic $-$then $e^w$ is transcendental (cf. \S45, \#11).
\vspace{0.2cm}

[Drop S4EC, impose instead SCHC, and bear in mind that the crux is when 
$w \notin \R \hsx \cup \hsx \sqrt{-1} \hsx \R$, thus $w$, $\ov{w}$ are $\Q$-linearly independent, so
\[
\trdegQ \Q(w, \ov{w}, e^w, e^{\ov{w}})
\ \geq \ 
2.
\]
If $e^w$ was algebraic, then $e^{\ov{w}} = \ov{e^{\raisebox{.02cm}{$\scriptstyle{w}$}}}$ would be too, reducing matters to 
\[
\trdegQ \Q(w, \ov{w}) 
\ \geq \ 
2,
\]
which is false since $\abs{w} \in \Qbar \implies \abs{w}^2 \in \Qbar = w \hsx \ov{w} \in \Qbar$.]
\end{x}
\vspace{0.3cm}

\begin{x}{\small\bf NOTATION}  \ 
Write
\[
\bx \ = \ (x_1, \ldots, x_n)
\]
and
\[
e^{\bx} \ = \ (e^{x_1}, \ldots, e^{x_n}).
\]
\end{x}
\vspace{0.3cm}

\begin{x}{\small\bf \un{N.B.}}\ 
SCHC can thus be abbreviated to
\[
\trdegQ \Q(\bx, e^{\bx}) \ \geq \ n.
\]
\end{x}
\vspace{0.3cm}

Consider a $\Q$-linear combination
\[
x_{n+1} \ = \ q_1 x_1 + \cdots + q_n x_n.
\]
Let \mM be a nonzero integer such that $M q_k$ is an integer for all $k = 1, \ldots, n$ and assume without loss of generality that 
\[
M q_1, \ldots, M q_t
\]
are nonnegative and
\[
M q_{t+1} , \ldots, M q_n
\]
are negative for some $0 \leq t \leq n$.  Let
\[
P(X_1, \ldots, X_{n+1}) 
\ = \ 
\prod\limits_{k=1}^t \hsx
X_k^{M q_k} 
\ - \ 
X_{n+1}^M \hsx 
\prod\limits_{k=t+1}^n \hsx 
X_k^{-M q_k}.
\]
Then
\allowdisplaybreaks
\begin{align*}
P(e^{x_1}&, \ldots, e^{x_{n+1}})
\\[12pt]
&=\ 
\prod\limits_{k=1}^t \hsx 
e^{x_k M q_k } - e^{x_{n+1} M} \hsx
\prod\limits_{k=t+1}^n \hsx 
e^{-x_k M q_k}
\\[12pt]
&=\ 
\prod\limits_{k=1}^t \hsx 
e^{M(q_k x_k)} - e^{M(q_1 x_1 \hsx + \cdots + \hsx q_n x_n)} \hsx
\prod\limits_{k=t+1}^n \hsx 
e^{-M(q_k x_k)}
\\[12pt]
&=\ 
\exp \bigg(\sum\limits_{k=1}^t \hsx M q_k x_k \bigg) 
-
\exp \bigg(M
\bigg(\sum\limits_{k=1}^t \hsx q_k x_k + \sum\limits_{k=t+1}^n\hsx q_k x_k\bigg) \bigg)
\exp \bigg(-\sum\limits_{k=t+1}^n \hsx M q_k x_k \bigg)
\\[12pt]
&=\ 
\exp\bigg(\sum\limits_{k=1}^t \hsx M q_k x_k \bigg)
\bigg(1 -
\exp \bigg(\sum\limits_{k=t+1}^n \hsx M q_k x_k \bigg)
\exp \bigg(-\sum\limits_{k=t+1}^n \hsx M q_k x_k \bigg)\bigg)
\\[12pt]
&=\ 
\exp \bigg(\sum\limits_{k=1}^t \hsx M q_k x_k \bigg) 
\bigg(1 -
\exp \bigg(\sum\limits_{k=t+1}^n \hsx M q_k x_k 
-
\sum\limits_{k=t+1}^n \hsx M q_k x_k \bigg)\bigg)
\\[12pt]
&=\ 
\exp \bigg(\sum\limits_{k=1}^t \hsx M q_k x_k \bigg) (1 - 1)
\\[12pt]
&=\ 
0.
\end{align*}

\begin{x}{\small\bf SCHOLIUM} \ 
The collection
\[
 e^{x_1}, \ldots, e^{x_n}, e^{x_{n+1}}
\]
is $\Q$-algebraically dependent.
\vspace{0.3cm}

So adding $x_{n+1}$, $e^{x_{n+1}}$ to
\[
\Q(x_1, \ldots, x_n, e^{x_1}, \ldots, e^{x_n})
\]
does not change the transcendence degree.
\end{x}
\vspace{0.3cm}

\begin{x}{\small\bf NOTATION} \ 
Given complex numbers $x_1, \ldots, x_n$, let
\[
\lindim_{\Q} \ \bx
\]
denote the linear dimension of the vector space over $\Q$ spanned by $x_1, \ldots, x_n$.
\end{x}
\vspace{0.3cm}

\begin{x}{\small\bf CONJECTURE} \ 
(SCHC) \ 
$\forall \ \bx$, 
\[
\trdegQ \Q(\bx, e^{\bx}) \ \geq \ \lindim_{\Q}\  \bx.
\]
\end{x}
\vspace{0.3cm}

To say that \bx is a counterexample to SCHC means that $x_1, \ldots, x_n$ are linearly independent over $\Q$ but
\[
\trdegQ \Q(\bx, e^{\bx}) \ < \ n.
\]
\vspace{0.3cm}

\begin{x}{\small\bf LEMMA} \ 
If there is a counterexample to SCHC, then there is a dense subset of $\C^n$ comprised of counterexamples.
\vspace{0.2cm}

PROOF \ 
If \bx is a counterexample to SCHC, then for any nonzero $q_1, \ldots, q_n$ in $\Q$, $q_1 x_1, \ldots, q_n x_n$ is also a counterexample.
\end{x}
\vspace{0.3cm}

\begin{x}{\small\bf NOTATION} \ 
Given \bx, put
\[
\delta(\bx) \ = \ \trdegQ \Q(\bx, e^{\bx}) - \lindim_{\Q} \ \bx,
\]
the 
\un{predimension}
\index{predimension} 
of \bx.
\end{x}
\vspace{0.3cm}
\begin{x}{\small\bf REMARK} \ 
SCHC is thus the claim that $\forall \ \bx$, 
\[
\delta(\bx) \ \geq \ 0,
\]
so a counterexample to Schanuel is an \bx with 
\[
\delta(\bx) \ < \ 0.
\]
If
\[
\delta(\bx) \ < \ -1,
\]
then for any complex number \mC, 
\[
\delta(\bx C) 
\ \leq \ 
\delta(\bx) + 1 
\ < \ 
0,
\]
leading therefore to continuum-many counterexamples..
\end{x}
\vspace{0.3cm}

\begin{x}{\small\bf LEMMA} \ 
$\forall \ n \in \N$, the set $X_n \subset \C^n$ of $n$-tuples which do not satisfy Schanuel's condition is first category and of Lebesgue measure 0.
\end{x}
\vspace{0.3cm}


\[
\textbf{APPENDIX}
\]
\vspace{0.5cm}

{\small\bf THEOREM}
(Admit SCHC) \ 
Let $\alpha \neq 1$ be a positive algebraic number and let $\beta$ be a positive irrational number.  
Assume: 
\[
\alpha^{\alpha^\beta} \ = \ \beta.
\]
Then $\beta$ is transcendental.
\vspace{0.5cm}

PROOF \ 
Suppose to the contrary that $\beta$ is algebraic, so by Gelfond-Schneider, $\alpha^\beta$ is transcendental.  
Claim: \ $1, \ \beta, \ \alpha^\beta$ are $\Q$-linearly independent.  
For suppose that
\[
r + s \beta + t \alpha^\beta \ = \ 0
\]
is a rational dependence relation:
\allowdisplaybreaks
\begin{align*}
r + s \beta  \in \ 
&\Qbar, \quad t\alpha^\beta \notin \Qbar \quad (\text{if} \ t \neq 0)
\\[12pt]
&\implies t = 0 
\\[12pt]
&\implies r, \ s = 0 \quad (\beta \in \PP).  
\end{align*}
Now multiply $1, \ \beta, \ \alpha^\beta$ by $\elln(\alpha) \neq 1$, hence
\[
\elln(\alpha), \ \beta \elln(\alpha), \ \alpha^\beta \elln(\alpha)
\]
are also $\Q$-linearly independent, hence by SCHC,
\[
\trdegQ 
\Q (\elln(\alpha), \ \beta \elln(\alpha), \ \alpha^\beta \elln(\alpha), \alpha, \alpha^\beta, \alpha^{\alpha^\beta}) 
\ \geq \ 3,
\]
i.e., 
\[
\trdegQ \Q (\elln(\alpha), \ \beta \elln(\alpha), \ \alpha^\beta \elln(\alpha),\alpha^\beta) 
\ \geq \ 3,
\]
i.e., 
\[
\trdegQbar \Qbar (\elln(\alpha), \ \beta \elln(\alpha), \ \alpha^\beta \elln(\alpha),\alpha^\beta) 
\ \geq \ 3.
\]
But
\begin{align*}
\trdegQbar \Qbar(\elln(\alpha),  \beta \elln(\alpha),  \alpha^\beta \elln(\alpha), \alpha^\beta) \ 
&= \ 
\trdegQbar \Qbar(\elln(\alpha), \alpha^\beta)
\\[12pt]
&\leq \ 
2.
\end{align*}
Contradiction.

%% file: _48_numerical_examples.tex
\chapter{
$\boldsymbol{\S}$\textbf{48}.\quad  SCHC: NUMERICAL EXAMPLES}
\setlength\parindent{2em}
\setcounter{theoremn}{0}
\renewcommand{\thepage}{\S48-\arabic{page}}

\ \indent 
Unless stipulated to the contrary, throughout the \S \ SCHC is in force.
\vspace{0.5cm}

\begin{x}{\small\bf EXAMPLE} \ 
The numbers $e$ and $e^e$ are algebraically independent over $\Q$. 
\vspace{0.2cm}

[Take $x_1 = 1$, $x_2 = e$ $-$then
\[
\trdegQ \Q(1, e, e^1, e^e) 
\ \geq \ 
2,
\]
i.e., 
\[
\trdegQ \Q(e, e^e) 
\ \geq \ 
2.]
\]
\end{x}
\vspace{0.3cm}

\begin{x}{\small\bf EXAMPLE}  \ \ 
The numbers \hsy $\elln(2)$ \hsy and \hsy $2^{\elln(2)}$ \hsy are algebraically independent over $\Q$. 
\vspace{0.2cm}

[Take $x_1 = \elln(2)$, $x_2 = (\elln(2))^2$ $-$then
\[
\trdegQ \Q(\elln(2), (\elln(2))^2, 2, 2^{\elln(2)}) 
\ \geq \ 
2,
\]
i.e., 
\[
\trdegQ \Q(\elln(2),\ 2^{\elln(2)}) 
\ \geq \ 
2.]
\]
\end{x}
\vspace{0.3cm}

\begin{x}{\small\bf EXAMPLE} \ \ 
The numbers \hsy $\elln(2)$ \hsy and \hsy $\elln(3)$ \hsy are algebraically independent over $\Q$. 
\vspace{0.2cm}

[Take $x_1 = \elln(2)$, $x_2 = \elln(3)$ $-$then
\[
\trdegQ \Q(\elln(2), \elln(3), 2, 3) 
\ \geq \ 
2,
\]
i.e., 
\[
\trdegQ \Q(\elln(2), \elln(3)) 
\ \geq \ 
2.]
\]
\vspace{0.2cm}

[Note: \ 
Recall that $\ds\frac{\elln(3)}{\elln(2)}$ is transcendental (cf. \S24, \#10), hence irrational.]
\end{x}
\vspace{0.3cm}

\begin{spacing}{1.45}
\begin{x}{\small\bf EXAMPLE} \ 
The numbers $e$ and $\pi$ are algebraically independent over $\Q$. 
\vspace{0.2cm}

[Take $x_1 = 1$, $x_2 = \sqrt{-1} \hsx \pi$ $-$then
\[
\trdegQ \Q(1, \sqrt{-1} \hsx \pi, e^1, e^{\sqrt{-1} \hsx \pi} = -1) 
\ \geq \ 
2,
\]
i.e., 
\[
\trdegQ \Q(\sqrt{-1} \hsx \pi, e)
\ \geq \ 
2.
\]
Therefore $e$ and $\sqrt{-1} \hsx \pi$ are algebraically independent over $\Q$.  
Suppose now that $e$ and $\pi$ are algebraically dependent over $\Q$, so there exists $P(X,Y) \in \Q[X,Y]$ nonzero such that 
$P(e, \pi) = 0$.  
Let $G(X,Y) = P(X, - \sqrt{-1} \hsx Y)$ and $H(X,Y) = \ov{P(X, - \sqrt{-1} \hsx Y)}$ $-$then 
\[
G(e, \sqrt{-1} \hsx \pi) 
\ = \ 
P(e, (-\sqrt{-1}) \sqrt{-1} \hsx \pi) 
\ = \ 
P(e, \pi) 
\ = \ 
0
\]
and
\[
H(e, \sqrt{-1} \hsx \pi) 
\ = \ 
\ov{P(e, (-\sqrt{-1}) \sqrt{-1} \hsx \pi) }
\ = \ 
\ov{P(e, \pi) }
\ = \ 
\ov{0}
\ = \ 
0.
\]
Consequently
\[
(G + H) (e, \sqrt{-1} \hsx \pi) 
\ = \ 
0.
\]
But $G + H$ is a nonzero polynomial with rational coefficients, thereby contradicting the algebraic independence over 
$\Q$ of $e$ and $\sqrt{-1} \hsx \pi$.]
\end{x}
\end{spacing}
\vspace{0.2cm}

[Three applications: 
\vspace{0.2cm}

\qquad \textbullet \quad 
$e + \pi$ is transcendental.  

\vspace{0.2cm}
[Suppose $e + \pi = \alpha \in \Qbar$.  Form
\[
P(X,Y) \ = \  X + Y - \alpha,
\]
an element of $\Qbar[X,Y]$ $-$then 
\[
P(e, \pi) \ = \ e + \pi - \alpha \ = \ 0.
\]
Contradiction.]
\vspace{0.2cm}

\qquad \textbullet \quad 
$e \hsx \pi$ is transcendental. 

\vspace{0.2cm}
[Suppose $e \hsx \pi = \alpha \in \Qbar$.  Form
\[
P(X,Y) \ = \  X Y - \alpha,
\]
an element of $\Qbar[X,Y]$ $-$then 
\[
P(e, \pi) \ = \ e \pi - \alpha \ = \ 0.
\]
Contradiction.]
\vspace{0.2cm}

\qquad \textbullet \quad 
$e / \pi$ is transcendental (hence $\pi / e$ is too).
\vspace{0.2cm}

[Suppose $e / \pi = \alpha \in \Qbar$.  Form
\[
P(X,Y) \ = \  X - \alpha Y,
\]
an element of $\Qbar[X,Y]$ $-$then 
\begin{align*}
P(e, \pi) \ 
&= \ e - \alpha  \pi 
\\[12pt]
&=\ 
\alpha \pi - \alpha \pi
\\[12pt]
&= \ 
0.
\end{align*}
Contradiction.]]
\vspace{0.3cm}

\begin{x}{\small\bf REMARK} \ 
It can be shown that unconditionally at least one of the following statements is true.
\vspace{0.2cm}

\qquad \textbullet \quad
The number $\ds e^{\pi^2}$ is transcendental.
\vspace{0.2cm}

\qquad \textbullet \quad
The numbers $e$ and $\pi$ are algebraically independent over $\Q$.
\vspace{0.2cm}

[Note: \ It is unknown whether $\ds e^{\pi^2}$ is even irrational.]
\end{x}
\vspace{0.3cm}

\begin{x}{\small\bf EXAMPLE} \ \ 
The numbers \hsy $e, \ \elln(2)$, \hsy and \hsy $\pi$  \hsy are algebraically independent over $\Q$.   
\vspace{0.2cm}

[Take $x_1 = 1$, $x_2 = \elln(2)$, $x_3 = \sqrt{-1} \hsx \pi$ to arrive at
\[
\trdegQ \Q(1, \elln(2), \sqrt{-1} \hsx \pi, e, 2, -1) \ \geq \  3.]
\]
\\[-1.2cm]

[Note: \ 
The numbers $1,\ \elln(2),\  \sqrt{-1} \hsx \pi$ are $\Q$-linearly independent (because $\elln(2)$ is irrational (cf. \S10, \#5).]
\end{x}
\vspace{0.3cm}
\begin{x}{\small\bf LEMMA} \ 
The eight numbers
\[
1,\ 
\sqrt{-1}\hsx \pi,\ 
\pi^2,\ 
e,\ 
e^2,\ 
\elln(2),\ 
2^{1/3}\elln(2),\ 
4^{1/3}\elln(2)
\]
are $\Q$-linearly independent.
\vspace{0.2cm}

PROOF \ 
The numbers $\sqrt{-1}\hsx \pi,\ e,\ \elln(2)$ are algebraically independent over $\Q$, hence are algebraically independent over $\Qbar$ 
(cf. \S20, \#7).  
Consider now a rational dependence relation
\[
A + B \sqrt{-1} \hsx \pi + C\pi^2 + De + F e^2 + G \elln(2) + H 2^{1/3} \elln(2) + K4^{1/3} \elln(2) 
\ = \ 
0.
\]
Define a polynomial $P \in \Qbar[X,Y,Z]$ by the prescription
\[
P(X,Y,Z) 
\ = \ 
A + BX - CX^2 + DY + FY^2 + GZ + H 2^{1/3} Z + K4^{1/3} Z.
\]
Then
\allowdisplaybreaks
\begin{align*}
P(\sqrt{-1} \hsx \pi, e, \elln(2)) \ 
&=\ 
A + B \sqrt{-1} \hsx \pi + C\pi^2 + De + F e^2 + G \elln(2)
\\[12pt]
&\qquad\qquad\qquad 
+ H 2^{1/3} \elln(2) + K4^{1/3} \elln(2)
\\[12pt]
&=\ 
0.
\end{align*}
Therefore
\[
A \ = \ B \ = \  C \ = \  D \ = \  F \ = \  G \ = \  H \ = \ K \ = \ 0.
\]
\end{x}
\vspace{0.3cm}

\begin{x}{\small\bf APPLICATION} \ 
The eight numbers
\[
e,\ 
\pi,\ 
e^e,\ 
e^{e^2},\ 
e^{\pi^2},\ 
2^{2^{1/3}},\ 
2^{2^{2/3}},\ 
\elln(2)
\]
are algebraically independent over $\Q$.
\vspace{0.2cm}

[Consider
\allowdisplaybreaks
\begin{align*}
&
1,\ 
\sqrt{-1}\hsx \pi,\ 
\pi^2,\ 
e,\ 
e^2,\ 
\elln(2),\ 
2^{1/3}\elln(2),\ 
4^{1/3}\elln(2),\  
\\[12pt]
&
e,\ 
-1,\ 
e^{\pi^2},\ 
e^e,\ 
e^{e^2},\ 
2,\ 
2^{2^{1/3}},\ 
2^{2^{2/3}}.]
\end{align*}
\end{x}
\vspace{0.3cm}

The next objective is \#14 infra, the verification of which proceeds in a series of steps.
\vspace{0.3cm}

\begin{x}{\small\bf LEMMA} \ 
Suppose that $x_1, \ldots, x_n$ is an algebraically independent set of positive real numbers $-$then $x_1, \ldots, x_n$ is multiplicatively independent (cf. \S34, Appendix).
\end{x}
\vspace{0.3cm}

\begin{x}{\small\bf EXAMPLE} \ 
The numbers $2,\ 3,\ \pi$, and $\elln(2)$ are multiplicatively independent: 
\[
2^a 3^b \pi^c (\elln(2))^d \ = \ 1 \qquad (a, b, c, d \in \Z)
\]
\[
\implies
a = b = c =d = 0.
\]

[The numbers $\pi$ and $\elln(2)$ are algebraically independent over $\Q$ (cf. \#6).  
This said, suppose that
\[
2^a 3^b \pi^c (\elln(2))^d \ = \ 1 \qquad (a, b, c, d \in \Z),
\] 
take for the sake of argument $c \geq 0$, $d \geq 0$, and introduce the polynomial
\[
P(X,Y) 
\ = \ 
2^a 3^b X^c Y^d - 1.
\] 
Then
\[
P(\pi, \elln(2)) \ = \ 2^a 3^b \pi^c (\elln(2))^d - 1
\]
\begin{align*}
&\implies 
c = 0,\ 
d = 0 
\\[12pt]
&\implies 
2^a 3^b - 1 = 0 
\\[12pt]
&\implies 
a = 0, \ b= 0.]
\end{align*}
\end{x}
\vspace{0.3cm}

\begin{x}{\small\bf LEMMA} \ 
Suppose that $x_1, \ldots, x_n$ is a multiplicatively independent set of positive real numbers $-$then the set 
$\elln(x_1), \ldots, \elln(x_n)$ is $\Q$-linearly independent.
\end{x}
\vspace{0.3cm}

\begin{x}{\small\bf EXAMPLE} \ 
The numbers
\[
\elln(\pi),\ 
\elln(2),\
\elln(3),\
\elln(\elln(2))
\]
are $\Q$-linearly independent (cf. \#10).
\end{x}
\vspace{0.3cm}

Therefore the numbers
\[
\sqrt{-1} \pi,\  
\elln(\pi),\ 
\elln(2),\
\elln(3),\
\elln(\elln(2))
\]
are $\Q$-linearly independent (consider real and imaginary parts).
\vspace{0.1cm}

Now use SCHC to arrive at
\[
\trdegQ 
\Q(
\sqrt{-1} \pi,\  
\elln(\pi),\ 
\elln(2),\
\elln(3),\
\elln(\elln(2)),\ 
-1,\ 
\pi,\ 
2,\ 
3,\ 
\elln(2)) 
\ \geq \ 5,
\]
from which the conclusion that
\[
 \pi,\  
\elln(\pi),\ 
\elln(2),\
\elln(3),\
\elln(\elln(2))
\]
are algebraically independent over $\Q$.
\\[-.1cm]

Next the numbers
\[
1,\ 
\sqrt{-1} \pi,\  
\elln(\pi),\ 
\elln(2),\
\elln(3),\
\elln(\elln(2))
\]
are $\Q$-linearly independent, thus invoking SCHC once again gives
\[
\trdegQ 
\Q(
1,\ 
\sqrt{-1} \pi,\ 
\elln(\pi),\ 
\elln(2),\
\elln(3),\
\elln(\elln(2)),\ 
e,\ 
-1,\ 
\pi,\ 
2,\ 
3,\ 
\elln(2)) 
\ \geq \ 6,
\]
so
\[
e,\ 
\pi,\ 
\elln(\pi),\ 
\elln(2),\
\elln(3),\
\elln(\elln(2))
\]
are algebraically independent over $\Q$.
\\[.2cm]

\begin{x}{\small\bf LEMMA} \ 
The seventeen numbers
\allowdisplaybreaks
\begin{align*}
& 
1,\ 
\sqrt{-1} \pi,\ 
\pi,\ 
\elln(\pi),\ 
e,\ 
e \hsx \elln(\pi),\ 
\pi \hsx\elln(\pi),\
\elln(2),\  
\\[12pt]
&
\pi \hsx\elln(2),\
e \hsx \elln(2),\
\sqrt{-1} \hsx \elln(2),\
\sqrt{-1},\
\sqrt{-1} \hsx\elln(\pi),\
\elln(3),\ 
\\[12pt]
&
\elln(\elln(2)),\ 
(\elln(3)) \hsx (\elln(\elln(2))),\ 
\sqrt{2} \hsx \elln(2)
\end{align*}
\end{x}
are $\Q$-linearly independent (cf. \#7).
\vspace{0.5cm}

\begin{x}{\small\bf THEOREM} \ 
(Waldschmidt's menagerie) \ (Admit SCHC) \ 
The seventeen numbers
\allowdisplaybreaks
\begin{align*}
& 
\pi,\ 
\elln(\pi),\ 
e,\ 
\elln(2),\ 
\elln(3),\ 
\elln(\elln(2)),\ 
e^\pi,\ 
e^e,\ 
\\[12pt]
&
\pi^e,\ 
\pi^\pi,\ 
2^\pi,\ 
2^e,\ 
2^{\sqrt{-1}},\ 
e^{\sqrt{-1}},\  
\pi^{\sqrt{-1}},\ 
\big(\elln(2)\big)^{\elln(3)},\ 
\sqrt{2}^{\sqrt{2}}
\end{align*}
are algebraically independent over $\Q$.

\end{x}
\vspace{0.3cm}

\begin{x}{\small\bf REMARK} \ 
$e^\pi$ is transcendental (unconditionally) (cf. \S20, \#10) but it is not even known whether 
$e^e$, $\pi^\pi$, and $\pi^e$ are irrational, let alone transcendental.
\end{x}
\vspace{0.3cm}

\begin{x}{\small\bf MISCELLANEA} \ 
(Admit SCHC) \ 
\vspace{0.5cm}

\qquad \textbullet \quad
$\ds{\sqrt{2}^{\sqrt{2}^{\sqrt{2}}}}$ \hspace{1.22cm} is transcendental.
\vspace{0.2cm}

\qquad \textbullet \quad
$\ds{\sqrt{-1}^{\sqrt{-1}^{\sqrt{-1}}}}$ \hspace{.5cm} is transcendental.
\vspace{0.2cm}

\qquad \textbullet \quad
$\ds\sqrt{-1}^{\hsx \textstyle e^\pi}$  \hspace{1.25cm} is transcendental.
\end{x}


%% file: _49_the_zero_condition.tex
\chapter{
$\boldsymbol{\S}$\textbf{49}.\quad  THE ZERO CONDITION}
\setlength\parindent{2em}
\setcounter{theoremn}{0}
\renewcommand{\thepage}{\S49-\arabic{page}}

\ \indent 
To begin with:
\vspace{0.3cm}

\begin{x}{\small\bf THE FUNDAMENTAL CONJECTURE} \ 
(FDC) \ 
Let $\lambda_1, \ldots, \lambda_d$ be elements of $\fL$ which are linearly independent over $\Q$ $-$then $\lambda_1, \ldots, \lambda_d$ are algebraically independent over $\Q$, hence are algebraically independent over $\Qbar$ (cf. \S20, \#7).
\vspace{0.2cm}

[Note: \ 
To appreciate how far away this conjecture lies, there is no known example of a $\Q$-linearly independent pair 
$\{\lambda_1, \lambda_2\}$ which is algebraically independent over $\Q$.]
\end{x}
\vspace{0.3cm}

\begin{x}{\small\bf \un{N.B.}} \ 
Recall that the fundemental conjecture is implied by SCHC (cf. \S47, \#7).
\end{x}
\vspace{0.3cm}

\begin{x}{\small\bf NOTATION} \ 
Fix $P \in \Q[X_1, \ldots, X_d]$, put
\[
Z(P) \ = \ \{\bx \in \C^d : P(\bx) = 0\}.
\]
\end{x}
\vspace{0.3cm}

\begin{x}{\small\bf DEFINITION} \ 
A nonzero polynomial $P \in \Q[X_1, \ldots, X_d]$ is said to satisfy the 
\un{zero condition}
\index{zero condition} 
if 
\[
Z(P) \hsx \cap \hsx \fL^d 
\ = \ 
\bigcup\limits_{\sV} \hsx
\sV \hsx \cap \hsx \fL^d,
\]
where $\sV$ ranges over the $\C$-vector subspaces of $\C^d$ rational over $\Q$ and contained in $Z(P)$.
\end{x}
\vspace{0.3cm}

\begin{x}{\small\bf EXAMPLE} \ 
Suppose that
\[
P(X_1, \ldots, X_d)  
\ = \ 
C_1 X_1 + \cdots + C_d X_d,
\]
where $C_1, \ldots, C_d \in \Q$ $-$then \mP satisfies the zero condition.
\end{x}
\vspace{0.3cm}


\begin{x}{\small\bf LEMMA} \ 
If every nonzero $P \in \Q[X_1, \ldots, X_d]$ satisfies the zero condition, then the fundamental conjecture is in force.
\vspace{0.2cm}

PROOF \ 
To get a contradiction, assume that $\lambda_1, \ldots, \lambda_d$ are linearly independent over $\Q$ but not algebraically independent over $\Q$, hence there exists a nonzero polynomial \mP in  $\Q[X_1, \ldots, X_d]$ such that 
$P(\lambda_1, \ldots, \lambda_d) = 0$, hence there is a $\C$-vector subspace $\sV$ of $\C^d$ rational over $\Q$ and contained in $Z(P)$ with 
\[
\blambda \ = \ (\lambda_1, \ldots, \lambda_d) 
\hsx \in \hsx \sV \hsx \cap \hsx \fL^d.
\]
Using the rationality of $\sV$ over $\Q$, write $\sV$ as the intersection of hyperplanes defined by linear forms with coefficients in $\Q$ (cf. \S37, \#2).  
Denoting by 
\[
\{(z_1, \ldots, z_d) \in \C^d : \beta_1 z_1 + \cdots + \beta_d z_d = 0 \qquad (\beta_1, \ldots, \beta_d \ \text{in} \ \Q)\}
\]
a typical such hyperplane, we then have
\[
\beta_1 \lambda_1 + \cdots +  \beta_d \lambda_d = 0,
\]
thus
\[
\beta_1 = 0, \ 
\ldots, \ 
\beta_d = 0
\]
and so $\sV = \{0\}$.  But
\[
(\lambda_1, \ldots, \lambda_d) \in \sV \hsx \cap \hsx  \fL^d 
\ = \ 
\{0\} \hsx  \cap\hsx  \fL^d 
\ = \ 
(0, \ldots, 0).
\]
\end{x}
\vspace{0.3cm}

\begin{x}{\small\bf REMARK} \ 
It is also true that the fundamental conjecture implies that every nonzero polynomial $P \in \Q[X_1, \ldots, X_d]$ satisfies the zero condition.
\end{x}
\vspace{0.3cm}

Our objective now will be to establish the four exponentials conjecture modulo yet another conjecture.

\vspace{0.2cm}

[Note: \ 
It was shown already in \S47, \#17 that 
\[
\text{SCHC} \implies \text{4EC}.]
\]
\vspace{0.1cm}

\begin{x}{\small\bf CONJECTURE} \ 
Work in $\C^4$ and define $P \in \Q[X_1, X_2, X_3, X_4]$ by
\[
P(X_1, X_2, X_3, X_4) \ = \ X_1 X_4 - X_2 X_3.
\]
Then \mP satisfies the zero condition.
\end{x}
\vspace{0.3cm}

\begin{x}{\small\bf CONJECTURE} \ 
Consider a $2 \times 2$ matrix \mM with entries in $\fL$:
\[
M \ = \ 
\begin{pmatrix}
\lambda_{1 1}  &&& \lambda_{1 2} \\[12pt]
\lambda_{2 1}  &&& \lambda_{2 2}
\end{pmatrix}
.
\]
Suppose that its rows are $\Q$-linearly independent and its columns are $\Q$-linearly independent $-$then 
\[
\rank M \ = \ 2 \qquad \text{(cf. \S44, \#2)}.
\]
\end{x}
\vspace{0.3cm}

\begin{x}{\small\bf \un{N.B.}}  \ 
The claim now is that
\[
\#8 \implies \#9.
\]
\end{x}
\vspace{0.3cm}

Here is another way to phrase it: \ If
\[
M \ = \ 
\begin{pmatrix}
\lambda_{1 1}  &&& \lambda_{1 2}  \\[12pt]
\lambda_{2 1}  &&& \lambda_{2 2}
\end{pmatrix}
\]
 is a $2 \times 2$ matrix with entries in $\fL$ and if 
 \[
\rank M \ = \ 1,
\]
then either its rows are $\Q$-linearly dependent or its columns are $\Q$-linearly dependent.
\vspace{0.3cm}

\begin{x}{\small\bf \un{N.B.}}  \ 
The condition
 \[
\rank M \ = \ 1
\]
implies that
\allowdisplaybreaks
\begin{align*}
\det M \ 
&=\ 
\lambda_{1 1} \lambda_{2 2} - \lambda_{1 2} \lambda_{2 1}
\\[12pt]
&=\ 
0.
\end{align*}
\end{x}
\vspace{0.3cm}

Per \#8, take for \mP the polynomial
\[
P(X_1, X_2, X_3, X_4) \ = \ X_1 X_4 - X_2 X_3.
\]
Substitute in
\[
X_1 \ = \ \lambda_{1 1}, \ 
X_4 \ = \ \lambda_{2 2}, \ 
X_2 \ = \ \lambda_{1 2}, \ 
X_3 \ = \ \lambda_{2 1}, \ 
\]
thus
\allowdisplaybreaks
\begin{align*}
P (\lambda_{1 1}, \lambda_{1 2}, \lambda_{2 1},\lambda_{2 2})  \ 
&=\ 
\lambda_{1 1} \lambda_{2 2} - \lambda_{1 2} \lambda_{2 1}
\\[12pt]
&=\ 
0
\end{align*}
and so
\[
(\lambda_{1 1}, \lambda_{1 2}, \lambda_{2 1},\lambda_{2 2})\in Z(P) \cap \fL^4.
\]
But
\[
Z(P) \hsx  \cap \hsx \fL^4 \ = \ \bigcup\limits_{\sV} \hsx \sV \cap \hsx \fL^4.
\]
Choose $\sV$: \ A $\C$-vector subspace of $\C^4$ rational over $\Q$ and contained in $Z(P)$ with 
\[
(\lambda_{1 1}, \lambda_{1 2}, \lambda_{2 1},\lambda_{2 2}) \in \sV \hsx \cap \hsx \fL^4.
\]

\begin{x}{\small\bf LEMMA} \ 
$\exists \ (a:b) \in \PP^1(\Q)$ such that $\sV$ is included either in the plane 
\[
*_1 : \{(z_1, z_2, z_3, z_4) \in \C^4 : a z_1 = b z_2, \ a z_3 = b z_4\}
\]
or in the plane 
\[
*_2 : \{(z_1, z_2, z_3, z_4) \in \C^4 : a z_1 = b z_3, \ a z_2 = b z_4\}.
\]
\vspace{0.2cm}

[Note: \ 
See the Appendix for the verification.]
\end{x}
\vspace{0.3cm}

\begin{x}{\small\bf \un{N.B.}} \ 
$(a:b)$ is the class of $(a,b)$ in the projective line $\PP^1(\Q)$.
\vspace{0.5cm}

Return to
\[
M = \ 
\begin{pmatrix}
\lambda_{1 1} && \lambda_{1 2}\\
\\
\lambda_{2 1} && \lambda_{2 2}
\end{pmatrix}
.
\]

\qquad \textbullet \quad
Assume $*_1$ in \#12 and work with the columns of \mM:
\[
\begin{pmatrix}
\lambda_{1 1}\\[12pt]
\lambda_{2 1}
\end{pmatrix}
, \quad 
\begin{pmatrix}
\lambda_{1 2}\\[12pt]
\lambda_{2 2}
\end{pmatrix}
.
\]
Then
\[
\begin{cases}
\ a \lambda_{1 1} = b \lambda_{1 2}\\[8pt]
\ a \lambda_{2 1} = b \lambda_{2 2}
\end{cases}
.
\]
Form now
\[
-a
\begin{pmatrix}
\lambda_{1 1}\\[12pt]
\lambda_{2 1}
\end{pmatrix}
\ + \ 
b
\begin{pmatrix}
\lambda_{1 2}\\
\\
\lambda_{2 2}
\end{pmatrix}
\]
or still,
\allowdisplaybreaks
\begin{align*}
\begin{pmatrix}
-a\lambda_{1 1}  \ + \  b\lambda_{1 2} \\[12pt]
-a\lambda_{2 1}  \ + \  b\lambda_{2 2}
\end{pmatrix}
&=\ 
\begin{pmatrix}
-b\lambda_{1 2}  \ + \  b\lambda_{1 2} \\[12pt]
-b\lambda_{2 2}  \ + \  b\lambda_{2 2}
\end{pmatrix}
\\[12pt]
&=\ 
\begin{pmatrix}
\ 0_{\ } \ \\
\ 0_{\ } \ 
\end{pmatrix}
.
\end{align*}
Since $(a:b) \in \PP^1(\C)$, the columns of \mM are linearly dependent and the four exponentials conjecture is thereby established.
\vspace{0.2cm}

\qquad \textbullet \quad
Assume $*_2$ in \#12 and work with the rows of \mM:
\[
[\lambda_{1 1} \hspace{.35cm} \lambda_{1 2}], \ 
[\lambda_{2 1}\hspace{.35cm}  \lambda_{2 2}] .
\]
This time
\[
\begin{cases}
\ a \lambda_{1 1} = b \lambda_{2 1}\\[8pt]
\ a \lambda_{1 2} = b \lambda_{2 2}
\end{cases}
\]
and one can consider
\[
-a \hsx [\lambda_{1 1} \hspace{.35cm} \lambda_{1 2}] + b \hsx [\lambda_{2 1} \hspace{.35cm} \lambda_{2 2}]. 
\]
\end{x}
\vspace{0.3cm}

\begin{spacing}{1.5}
It is not necessary to utilize \#8 in order to arrive at a restricted but unconditional result, the idea being to reduce the elements 
$\blambda$ in $Z(P) \hsx \cap \hsx \fL^4$ for which there is a $\sV$ : \ 
A $\C$-vector subspace of $\C^4$ rational over $\Q$ and contained in $Z(P)$ with $\blambda \in \sV \hsx \cap \hsx \fL^4$.
\end{spacing}
\vspace{0.3cm}

\begin{x}{\small\bf THEOREM} \ 
Take a 
\[
\blambda \ = \ (\lambda_1, \lambda_2, \lambda_3,\lambda_4) \in Z(P) \hsx \cap \hsx \fL^4 .
\]
Then either $\blambda \in \sV$ for some $\sV$ per supra or else
\[
\trdegQ \Q(\lambda_1, \lambda_2, \lambda_3,\lambda_4) \ \geq \ 2.
\]
\end{x}
\vspace{0.3cm}

\begin{x}{\small\bf SCHOLIUM} \ 
The statement of the four exponentials conjecture holds true for the set of those
\[
\blambda \ = \ (\lambda_1, \lambda_2, \lambda_3,\lambda_4) \in Z(P) \hsx \cap\hsx  \fL^4
\]
with the property that 
\[
\trdegQ \Q(\lambda_1, \lambda_2, \lambda_3,\lambda_4) \ = \ 1.
\]
\vspace{0.2cm}

[Note: \ 
The point, of course, is that for this set of $\blambda$, \#12 is applicable.]
\end{x}
\vspace{0.3cm}

\begin{x}{\small\bf \un{N.B.}} \ 
The $\lambda_i$ $(i = 1, 2, 3, 4)$ are transcendental (if not zero).
\end{x}
\vspace{0.3cm}

\[
\textbf{APPENDIX}
\]
\vspace{0.3cm}

The issue is the validity of \#12.  Write
\[
\begin{cases}
\ *_1 \ = \ W_1(a:b)\\[8pt]
\ *_2 \ = \ W_2(a:b)
\end{cases}
\]
and note that
\[
\begin{cases}
\ W_1(a:b)\\[8pt]
\ W_2(a:b)
\end{cases}
\subset Z(P).
\]
Matters are trivial if $\sV$ is contained in 
\[
W_1(0:1)
\quad \text{or} \quad
W_1(1:0) 
\quad \text{or} \quad
W_2(0:1)
\quad \text{or} \quad
W_2(1:0).
\]
Assume, therefore, that there exists $\bv = (w,x,y,z) \in \sV$ such that   
$w \hsy x \hsy y \hsy z \neq 0$.
Since $w z = x y$, we have $(x:w) = (z:y)$ and $(y:w) = (z:x)$, the claim then being that the supposition
\[
\sV \not\subset W_1(x:w) \quad \text{and} \quad \sV \not\subset W_2(y:w)
\]
leads to a contradiction.  
Choose $\bv^\prime = (w^\prime ,x^\prime ,y^\prime ,z^\prime )$ in $\sV$ which does not belong to
\[
W_1(0:1) \cup
W_1(1:0) \cup
W_2(0:1) \cup
W_2(1:0) \cup
W_1(x:w) \cup
W_2(y:w).
\]
Accordingly
\[
w^\prime x^\prime y^\prime z^\prime \ \neq \ 0.
\]
Moreover
\[u \bv + u^\prime \bv^\prime \in \sV
\]
for all $(u, u^\prime) \in \C^2$, hence
\[
P(u \bv + u^\prime \bv^\prime) \ = \ 0
\]
or still, 
\[
P((u w, u x, u y, u z) + (u^\prime  w^\prime , u^\prime  x^\prime , u^\prime  y^\prime , u^\prime  z^\prime )) 
\ = \ 0
\]
or still, 
\[
P(u w + u^\prime  w^\prime , u x + u^\prime  x^\prime , u y + u^\prime  y^\prime , u z + u^\prime  z^\prime ) 
\ = \ 0
\]
or still, 
\[
(u w + u^\prime  w^\prime )(u z + u^\prime  z^\prime ) - (u x + u^\prime  x^\prime ) (u y + u^\prime  y^\prime ) 
\ = \ 0
\]
or still, 
\[
(w z - x y) u^2 + (w^\prime z - x y^\prime - x^\prime y + w z^\prime) u u^\prime 
+ (w^\prime z^\prime  - x^\prime  y^\prime ){u^\prime}^2 \ = \ 0
\]
\qquad\qquad $\implies$
\[
w z = x y, \ 
w^\prime z^\prime = x^\prime y^\prime, \ 
w^\prime z + w z^\prime = x y^\prime + x^\prime y,
\]
$(u, u^\prime) \in \C^2$ being arbitrary.  Therefore
\allowdisplaybreaks
\begin{align*}
(y z^\prime - y^\prime z) (x z^\prime - x^\prime z)\ 
&=\ 
z z^\prime (w^\prime z - x y^\prime - x^\prime y + w z^\prime) 
\\[12pt]
&=\ 
0.
\end{align*}
So at least one of the numbers
\[
y z^\prime - y^\prime z, \quad
x z^\prime - x^\prime z 
\]
must vanish.
\vspace{0.2cm}

\qquad \textbullet \quad
$y z^\prime - y^\prime z = 0$\\

$\implies$
\[
\frac{w}{x} 
\ = \ 
\frac{y}{z} 
\ = \ 
\frac{y^\prime}{z^\prime} 
\ = \ 
\frac{w^\prime}{x^\prime} 
\implies
\bv^\prime \in W_1(x:w),
\]
a contradiction.
\vspace{0.2cm}

\qquad \textbullet \quad
$x z^\prime - x^\prime z = 0$\\

$\implies$
\[
\frac{w}{y} 
\ = \ 
\frac{x}{z} 
\ = \ 
\frac{x^\prime}{z^\prime} 
\ = \ 
\frac{w^\prime}{y^\prime} 
\]
\qquad $\implies$
\[
\bv^\prime \in W_2(y:w),
\]
a contradiction.
\vspace{0.2cm}

Since $\sV$ is rational over $\Q$ (by hypothesis), there is a basis $\be_1, \ldots, \be_d$ for $\sV \  (d  \leq 2)$ with 
\[
\be_i \ = \ (e_{i 1}, e_{i 2}, e_{i 3},e_{i 4}) \in \Q^4.
\]
If $\sV$ is included in $W_1(a:b)$ for some $(a:b) \in \PP^1(\C)$, then the system of equations
\[
u e_{i 1} \ = \ u^\prime e_{i 2}, \quad 
u e_{i 3} \ = \ u^\prime e_{i 4}\ 
\qquad (i = 1, \ldots, d)
\]
\begin{spacing}{1.6}
\noindent has a nontrival solution $(u, u^\prime) \in \C^2$, thus it has a nontrivial solution $(u, u^\prime) \in \Q^2$.  
Consequently $\sV$ is included in $W_1(a:b)$ for some $(a:b) \in \PP^1(\Q)$.  
The story for $W_2(a:b)$ is analogous.
\end{spacing}

%% file: _50_property_matrix_abc_.tex
\chapter{
$\boldsymbol{\S}$\textbf{50}.\quad  PROPERTY
$
\huge\bf{\binom{ A \ B}{C \ 0}}
$
}
\setlength\parindent{2em}
\setcounter{theoremn}{0}
\renewcommand{\thepage}{\S50-\arabic{page}}

\ \indent 
Let $\K$ be a field, $\bk \subset \K$ a subfield.
\vspace{0.5cm}

\begin{x}{\small\bf DEFINITION} \ 
Two $m \times n$ matrices \mM and \mN with entries in $\K$ are 
\un{\bk-equivalent}
\index{\bk-equivalent} 
if there exist nonsingular matrices \mP and \mQ with entries in \bk such that $N = P M Q$.
\vspace{0.2cm}

[Note: \ 
The dimension of the $\Q$-subspace of $\K^n$ generated by the rows of \mM (or \mN) is the same as the dimension of the $\Q$-subspace of $\K^m$ generated by the columns of \mM (or \mN).]
\end{x}
\vspace{0.3cm}

\begin{x}{\small\bf \un{N.B.}} \ 
The rank of \mM equals the rank of \mN, this being the largest integer $r$  for which there exists a nonsingular $r \times r$ submatrix of \mM (or \mN) (cf. \S35, \#8).
\end{x}
\vspace{0.3cm}

\begin{x}{\small\bf THEOREM} \ 
Let $\sE$ be a $\bk$-vector subspace of $\K$ which is spanned by a family (finite or infinite) of elements of $\K$ which are algebraically independent over $\bk$ $-$then every matrix \mM with entries in $\sE$ is $\bk$-equivalent to a matrix of the form
\[
\begin{pmatrix}
A &&& B\\[12pt]
C &&& 0\\
\end{pmatrix}
,
\]
where \mA is either zero-size or nonsingular.
\end{x}
\vspace{0.3cm}

To orient ourselves, here are two examples of the overall structural setup (ignoring for the time being the validity of the assumption on $\sE$).
\vspace{0.3cm}

\begin{x}{\small\bf EXAMPLE} \ 
Take $\K = \C$, $\bk = \Q$, let $\sE_0$ be the $\Q$-vector space $\fL$ of logarithms
of algebraic numbers, and put $\sE = \Q + \fL$.
\vspace{0.2cm}

[Note: \
The sum is direct.  In fact, 
\[
\Qbar \hsx  \cap   \hsx \fL= \ \{0\} \qquad \text{(cf. \S31, \#3)} \quad 
\implies \Q \hsx \cap \hsx \fL \ = \ \{0\}.]
\]
\end{x}
\vspace{0.3cm}

\begin{x}{\small\bf EXAMPLE} \ 
Take $\K = \C$, $\bk = \Qbar$, let $\sE_0$ be the $\Qbar$-vector space of homogeneous linear combinations of elements of $\fL$ with coefficients in $\Qbar$, and  put $\sE = \Qbar + \sE_0$ (hence $\sE = \fL^*$).
\vspace{0.2cm}

[Note: \ 
The sum $\Qbar + \sE_0$ is direct (cf. \S39, \#14).]
\end{x}
\vspace{0.3cm}

\begin{x}{\small\bf LEMMA} \ 
Suppose that $\sE$ is a \bk-vector subspace of $\K$ $-$then the following conditions are equivalent.
\vspace{0.2cm}

\qquad (i) \quad
$\sE$ is spanned by a family (finite or infinite) of elements of $\K$ which are algebraically independent over \bk.
\vspace{0.2cm}

\qquad (ii) \quad
Subsets of $\sE$ which are linearly independent over \bk are algebraically independent over \bk.
\vspace{0.2cm}

\qquad (iii) \quad
If $\sE^\prime$ is a vector subspace of $\sE$ and $x$ is an element of $\sE$ which does not belong to $\sE^\prime$, then $x$ is transcendental over $\bk (\sE^\prime).$
\vspace{0.2cm}

\vspace{0.2cm}

PROOF 

\qquad (i) $\implies$ (ii) \quad
Per the assumption, fix a basis \mB for $\sE$ over \bk consisting of elements of $\K$ which are algebraically independent over \bk.  
Let $x_1, \ldots, x_m$ be a set of \bk-linearly independent elements of $\sE$ and write each $x_i$ $(1 \leq i \leq m)$ as a linear combination with coefficients in \bk of elements $y_j \in B$ $(1 \leq j \leq n)$, say
\[
x_i 
\ = \ 
\sum\limits_{j=1}^n \hsx a_{i j} y_j.
\]
Since the matrix $[a_{i j}]$ has rank $m$, it follows that there is a subset $\{z_1, \ldots, z_{n-m}\}$ of $\{y_1, \ldots, y_n\}$ such that
\[
\bk(y_1, \ldots, y_n) 
\ = \ 
\bk(x_1, \ldots, x_m, z_1, \ldots, z_{n-m}).
\]
And this relation implies that $x_1, \ldots, x_m$ are algebraically independent over \bk.
\vspace{0.2cm}

\qquad (ii) $\implies$ (iii) \quad
Assume instead that $x \in \sE$, $x \notin \sE^\prime$ is algebraic over $\bk(\sE^\prime)$.  
Choose $y_1, \ldots, y_n$ in $\sE^\prime$, linearly independent over \bk, such that $x$ is algebraic over $\bk(y_1, \ldots, y_n)$ $-$then
$y_1, \ldots, y_n, x$ are algebraically dependent over \bk, hence by (ii), are linearly dependent over \bk, say
\[
a_1 y_1 + \cdots + a_n y_n - ax \ = \ 0.
\]
But $a$ cannot be zero (since otherwise $a = 0$ would force  $y_1, \ldots, y_n$ to be linearly dependent over \bk), hence
\[
x 
\ = \ 
\frac{a_1}{a} y_1 + \cdots + \frac{a_n}{a} y_n \in \sE^\prime, 
\]
contradicting $x \notin \sE^\prime$.
\vspace{0.2cm}

\qquad (iii) $\implies$ (i) \quad 
Let \mB be a basis for $\sE$ over \bk.  
Claim: \ Any subset $\{y_1, \ldots, y_n\} \subset B$ of \bk-linearly independent elements of \mB consists of \bk-algebraically independent elements.  
To establish this, proceed by induction on $n$.
\vspace{0.2cm}

\qquad \textbullet \quad
$n = 1$: \ Use (iii) with $\sE^\prime = \{0\}$: 
\[
y_1 \neq 0 \implies y_1 \notin \sE^\prime.
\]
Therefore $y_1$ is transcendental over \bk.
\vspace{0.2cm}

\qquad \textbullet \quad
$n \geq 2$: \ Assume the result holds at level $n - 1$ and let $y_1, \ldots, y_n$ be \bk-linearly independent elements of \mB.  
Denote by $\sE^\prime$ the vector subspace of $\sE$ over \bk
spanned by $y_1, \ldots, y_{n-1}$.  
Owing to the induction hypothesis, $y_1, \ldots, y_{n-1}$ are algebraically independent over \bk.  
But $y_n \notin \sE^\prime$, so by (iii), $y_n$ is transcendental over the field $\bk(y_1, \ldots, y_{n-1})$ from which  $y_1, \ldots, y_n$ are algebraically independent over \bk.  
\vspace{0.2cm}

[Note: \
There is yet another equivalent condition that can be added to this list, viz:
\vspace{0.2cm}

\qquad (iv) \quad 
For any nonzero polynomial $P \in \bk[X_1, \ldots, X_n]$, 
\[
Z(P) \hsx \cap \hsx \sE^n \ = \ \bigcup\limits_{\sV} \hsx \sV \hsx \cap \hsx \sE^n,
\]
where $\sV$ ranges over the $\K$-vector subspaces of $\K^n$ rational over \bk and contained in 
\[
Z(P) \ = \ \{\bx \in \K^n : P(\bx) = 0\}.]
\]
\end{x}
\vspace{0.3cm}

\begin{x}{\small\bf NOTATION} \ 
Let $\sE_0$ be the \bk-vector subspace of $\sE$ spanned by the entries of \mM.
\end{x}
\vspace{0.3cm}

The proof of \#3 goes via induction in the dimension $n$ of $\sE_0$.
\vspace{0.3cm}

\qquad \textbullet \quad 
$n = 1$: \quad 
Write $M = N x$, where \mN has entries in \bk and $x \in \sE$, $x \neq 0$.  
Let $r$ be the rank of \mN and let \mP and \mQ be nonsingular matrices with entries in \bk such that 
\[
P N Q \ = \ 
\begin{pmatrix}
I_r &&&0\\[12pt]
0 &&&0
\end{pmatrix}
.
\]
Then
\[
P M Q \ = \ 
\begin{pmatrix}
I_r x &&&0\\[12pt]
0 &&&0
\end{pmatrix}
,
\]
so matters are satisfied with the choices
\[
A \ = \ I_r x, \quad B \ = \ 0, \quad C \ = \ 0.
\]
\\[-1cm]

\qquad \textbullet \quad 
$n = 2$: \quad 
Write
\[
M \ = \ M_1 x_1 + M_2 x_2,
\]
where $M_1$ and $M_2$ are matrices with entries in \bk and where $x_1$, $x_2 \in \sE$ are linearly independent over \bk 
(hence algebraically independent over \bk (cf. \#6 (ii)).  
Denote by $r_1$ the rank of $M_1$.  
Choose nonsingular matrices $P_1$ and $Q_1$ with entries in \bk such that
\[
P_1 M_1 Q_1 \ = \ 
\begin{pmatrix}
I_{r_1} &&&0\\[12pt]
0 &&&0\\
\end{pmatrix}
.
\]
Denote by $A_2$, $B_2$, $C_2$, $D_2$ the matrices with entries in \bk such that 

\[
P_1 M_2 Q_1 \ = \ 
\begin{pmatrix}
A_2 &&&B_2\\[12pt]
C_2 &&&D_2\\
\end{pmatrix}
,
\]\\
where $A_2$ is an $r_1 \times r_1$ matrix.  
Then\\
\[
P_1 M Q_1 \ = \ 
\begin{pmatrix}
I_{r_1} x_1 + A_2 x_2 &&&B_2 x_2\\[12pt]  
C_2 x_2 &&&D_2 x_2\\
\end{pmatrix}
.
\]
Choose nonsingular matrices $P_2$ and $Q_2$ with entries in \bk such that
\[
P_2 D_2 Q_2 \ = \ 
\begin{pmatrix}
I_{r_2} &&&0\\[12pt]
0 &&&0\\
\end{pmatrix}
,
\]\\
where $r_2$ is the rank of $D_2$.  
Then\\
\[
\begin{pmatrix}
I_{r_1} &&&0\\[12pt]
0 &&&P_2\\
\end{pmatrix}
\  
P_1 M Q_1 
\ 
\begin{pmatrix}
I_{r_1} &&&0\\[12pt]
0 &&&Q_2\\
\end{pmatrix}
\]\\
equals\\
\[
\begin{pmatrix}
I_{r_1} x_1 + A_2 x_2 &&B_2^{\prime} x_2 &&B_2^{\prime\prime} x_2\\[12pt]
C_2^\prime x_2 &&I_{r_2} x_2 && 0\\[12pt]
C_2^{\prime\prime} x_2 &&0 &&0\\
\end{pmatrix}
,
\]\\
where $B_2^\prime$, $B_2^{\prime\prime}$, $C_2^\prime$, $C_2^{\prime\prime}$ have entries in \bk.  
Put now

\[
A \ = \ 
\begin{pmatrix}
I_{r_1} x_1 + A_2 x_2 &&&B_2^\prime x_2\\[12pt]
C_2^\prime x_2 &&&I_{r_2} x_2\\
\end{pmatrix}
\]
and take for \mB, \mC what remains.  
To check that \mA is nonsingular, note that the
determinant of \mA is a polynomial in $x_1$ and $x_2$ and the coefficient of $x_1^{r_1} x_2^{r_2}$ is 1.  
Therefore
\[
\det A \neq 0.
\]

\begin{spacing}{1.55}
\qquad \textbullet \quad 
$n > 2$: \quad 
Fix a nonzero element $x \in \sE_0$.  
Let $\sE_1$ be a subspace of $\sE_0$ such that $\sE_0 = \sE_1 \oplus \bk x$.  
Write $M = x N + M_1$, where \mN has entries in \bk and $M_1$ has entries in $\sE_1$.  
Denote by $r$ the rank of \mN.  
Choose nonsingular matrices \mP and \mQ with entries in \bk such that
\end{spacing}
\[
P N Q
\ = \ 
\begin{pmatrix}
I_r &&&0\\[12pt]
0 &&& 0\\
\end{pmatrix}
.
\]\\
Then\\
\[
P M Q
\ = \ 
\begin{pmatrix}
x I_r + A_1 &&&B_1\\[12pt]
C_1 &&&D_1
\end{pmatrix}
,
\]
where $A_1$, $B_1$, $C_1$, $D_1$ have their entries in $\sE_1$.  
Apply now the induction hypothesis to $D_1$:
\[
P^\prime D_1 Q^\prime
\ = \ 
\begin{pmatrix}
A^\prime &&&B^\prime\\[12pt]
C^\prime &&&0\\
\end{pmatrix}
.
\]
Here $A^\prime$ is nonsingular with enteries in $\sE_1$.  Next\\[8pt]
\[
\begin{pmatrix}
I_r && 0\\[12pt]
0 && P^\prime\\
\end{pmatrix}
\quad
\begin{pmatrix}
x I_r + A_1 && B_1\\[12pt]
C_1 && D_1\\
\end{pmatrix}
\quad
\begin{pmatrix}
I_r && 0\\[12pt]
0 && Q^\prime\\
\end{pmatrix}
\]\\
equals\\
\[
\begin{pmatrix}
A &&& B\\[12pt]
C &&& 0\\
\end{pmatrix}
.
\]\\
Here\\
\[
A 
\ = \ 
\begin{pmatrix}
x I_r + A_1 &&&B^{\prime\prime}\\[12pt]
C^{\prime\prime} &&&A^\prime
\end{pmatrix}
\]\\[-8pt]
\begin{spacing}{1.55}
\noindent and the entries of $B^{\prime\prime}$, $C^{\prime\prime}$ are in $\sE_1$.  
To assertain that \mA is nonsingular, note that the determinant of \mA is a polynomial in $x$ with coefficients in $\bk(\sE_1)$ whose 
term of highest degree is $x^r \det A^\prime$.  
Since $x \notin \sE_1$, it follows from \#6 (iii) that $x$ is transcendental over $\bk(\sE_1)$ and since $A^\prime$ is nonsingular, the bottom line is that
\end{spacing}
\[
\det A \ \neq \ 0.
\]
\vspace{0.1cm}

\begin{x}{\small\bf DEFINITION} \ 
Let $\sE$ be a \bk-vector subspace of $\K$ $-$then by property 
$\binom{ A \ B}{C \ 0}$
\index{property $\binom{ A \ B}{C \ 0}$} 
we shall understand the following: \ 
Any nonzero matrix \mM with enteries in $\sE$ is \bk-equivalent to a matrix of the form
\[
\begin{pmatrix}
A && B\\[12pt]
C && 0
\end{pmatrix}
,
\]
where \mA is nonsingular.
\vspace{0.2cm}

[Note: \
Stricly speaking this is a property of the triple
\[
(\bk, \K, \sE)
\]
but usually one abuses the language and simply says that $\sE$ has property $\binom{ A \ B}{C \ 0}$.]
\end{x}
\vspace{0.3cm}

\begin{x}{\small\bf \un{N.B.}}  \ 
The upshot of \#3 is that if $\sE$ is a \bk-vector subspace of $\K$ spanned by \bk-algebraically independent elements, then $\sE$ satisfies property $\binom{ A \ B}{C \ 0}$.
\end{x}
\vspace{0.3cm}

\begin{x}{\small\bf LEMMA} \ 
If $\sE_0$ is a \bk-vector subspace of $\K$ spanned by \bk-algebraically independent elements and if 
$\sE_0 \hsx \cap \hsx \bk = \{0\}$, then 
$\sE = \bk + \sE_0$ satisfies property $\binom{ A \ B}{C \ 0}$.
\vspace{0.2cm}

PROOF \ 
As a \bk-vector space, $\sE$ is isomorphic to the subspace $\sE^\prime = \bk X + \sE_0$ of $\K(X)$ and property 
$\binom{ A \ B}{C \ 0}$ holds for the triple 
\[
(\bk, \K(X), \sE^\prime).
\]
\end{x}
\vspace{0.3cm}

\begin{x}{\small\bf EXAMPLE} \ 
As in \#4, take $\K = \C$, $\bk = \Q$, $\sE_0 = \fL$, and admit FDC (cf. \S49, \#1) $-$then \#6 (ii) is in force which implies that \#6 (i) is in force.  
Accordingly, since $\sE_0 \cap \bk = \{0\}$, it follows that $\sE = \bk + \sE_0$ satisfies property $\binom{ A \ B}{C \ 0}$.
\vspace{0.2cm}

[Note: \ 
Of course, $\sE_0$ also satisfies property $\binom{ A \ B}{C \ 0}$.]
\end{x}
\vspace{0.3cm}

\begin{x}{\small\bf REMARK} \ 
The satisfaction of property $\binom{ A \ B}{C \ 0}$ is not automatic.
\vspace{0.2cm}

\begin{spacing}{1.55}
[To illustrate, choose elements $x$ and $u$ in $\K$ such that $u$, $ux$, $ux^2$ are
\bk-linearly independent ($\implies x \notin \bk$).  
Denote by $\sE$ the \bk-vector space $k u + k u x + k u x^2$ 
($\implies \dim_{\bk} (\sE) = 3$) $-$then the triple $(\bk, \K, \sE)$ does not satisfy property $\binom{ A \ B}{C \ 0}$.
Thus consider the line $\sV = \K(1,x)$ in $\K^2$ (the hyperplane defined by the equation $z_2 = x z_1$) and note that 
$\sV \hsx \cap \hsx \bk^2 = \{0\}$.  
Furthermore $\sV \hsx \cap \hsx \sE^2$ contains the \bk-linearly independent points
\end{spacing}
\[
(u, ux) , \quad (u x, u x^2)
\]
implying thereby that $\dim_{\bk} (\sV \hsx \cap \hsx \sE^2) \geq 2$.  
On the other hand, taking into account \S51, \#3 infra (with $d = 2$, $n = 1$),
\[
\dim_{\bk} (\sV \hsx \cap \hsx \sE^2) 
\ \leq \ 
1 (1 + 1) / 2
\ = \ 1.
\]
So, on the basis of this contradiction, the triple $(\bk, \K, \sE)$ does not satisfy property $\binom{ A \ B}{C \ 0}$.]
\end{x}
\vspace{0.75cm}

\[
\textbf{APPENDIX}
\]
\vspace{0.1cm}

Let $\K$ be a field, $\bk \subset \K$ a subfield.
\vspace{0.5cm}

{\small\bf LEMMA}
Suppose that $\sE$ is a \bk-vector subspace of $\K$ containing \bk $-$then the following conditions are equivalent (cf. \#6).
\vspace{0.2cm}

(i) \quad 
There exists a basis $\{x_i : i \in I\}$ for $\sE$ over \bk with $0 \in I$, $x_0 = 1$, and $\{x_i : i \in I, i \neq 0\}$ algebraically independent over \bk.
\vspace{0.2cm}

(ii) \quad 
If $x_1, \ldots, x_n$ are elements in $\sE$ such that $1, x_1, \ldots, x_n$ are linearly independent over \bk, then $x_1, \ldots, x_n$ are algebraically independent over \bk.
\vspace{0.2cm}

(iii) \quad 
For any tuple $(x_0, \ldots, x_n)$ consisting of \bk-linearly independent elements of $\sE$ and for any nonzero homogeneous polynomial 
$P \in \bk[X_0, \ldots, X_n]$, the number $P(x_0, \ldots, x_n)$ is not zero.
\vspace{0.2cm}

(iv) \quad 
If $P \in \bk[X_0, \ldots, X_n]$ is a nonzero homogeneous polynomial, then
\[
Z(P) \hsx \cap \hsx \sE^{n+1}
\ = \ 
\bigcup\limits_{\sV} \hsx \sV \hsx \cap \hsx \sE^{n+1},
\]
where $\sV$ ranges over the $\K$-vector subspaces of $\K^{n+1}$ rational over \bk and contained in
\[
Z(P) 
\ = \ 
\big\{\bx \in \K^{n+1} : P(\bx ) = 0\big\}.
\]

%% file: _51_vector_spaces_L_bis.tex
\chapter{
$\boldsymbol{\S}$\textbf{51}.\quad  VECTOR SPACES: $\fL$ (bis)}
\setlength\parindent{2em}
\setcounter{theoremn}{0}
\renewcommand{\thepage}{\S51-\arabic{page}}


\begin{x}{\small\bf RAPPEL} \ 
Let $\sV \subset \C^d$ be a $\C$-vector subspace such that $\sV \hsx \cap \hsx \Q^d = \{0\}$ $-$then
\[
\dim_{\Q} \hsx (\sV \hsx \cap \hsx \fL^d) 
\ \leq \ 
n (n+1)
\qquad \text{(cf. \S38, \#5)},
\]
where
\[
n \ = \ \dim _{\C} \hsx (\sV).
\]
\end{x}
\vspace{0.3cm}

\begin{x}{\small\bf \un{N.B.}} \ 
This result is unconditional.
\end{x}
\vspace{0.3cm}

Return now to the setup of \S50.
\vspace{0.3cm}

\begin{x}{\small\bf THEOREM} \ 
Let $\sE$ be a $\bk$-vector subspace of $\K$ satisfying property $\binom{ A \ B}{C \ 0}$.  
Let $\sV \subset \K^d$ be a $\K$-vector subspace $-$then 
\[
\dim_{\bk} \hsx (\sV \hsx \cap \hsx \sE^d) 
\ \leq \ 
n (n+1) / 2,
\]
where
\[
n \ = \ \dim _{\K} \hsx (\sV).
\]

PROOF \ 
When $d = 1$, $\sV = \{0\}$ and $\sV \hsx\hsx \cap \hsx\hsx \sE = \{0\}$.  
Assume now that $d \geq 2$ ($\implies n < d$).\\

\qquad \textbullet \quad 
By induction on $d$, if $r < d$ and if $\sW$ is a $\K$-vector subspace of $\K^r$ such that 
$\sW \hsx \cap \bk^r = \{0\}$, then the $\bk$-vector space $\sW \hsx \cap \hsx \sE^r$ is finite dimensional, in fact 
\[
\dim_{\K} \hsx (\sW \hsx \cap \hsx \sE^r) 
\ \leq \ 
r (r - 1) / 2 
\qquad \text{(see below)}.
\]

Take now $\ell$ elements $\bx_1, \ldots, \bx_\ell$ in $\sV \hsx \cap \hsx \sE^d$ which are linearly independent over $\bk$, the claim being that 
\[
\ell 
\ \leq \ 
n(n+1) / 2.
\]
Denote by \mM the $d \times \ell$ matrix whose columns are given by the coordinates of the $\bx_i$ $(i = 1, \ldots, \ell)$ 
$-$then the entries of \mM are in $\sE$, so \mM is $\bk$-equivalent to a matrix
\[
\begin{pmatrix}
A & &B \\[12pt]
C & &0
\end{pmatrix}
,
\]
where \mA is a nonsingular $r \times r$ matrix.  In addition
\[
d 
\ > \ 
n
\ \geq \ 
\rank M 
\ \geq \ 
r 
\implies 
r 
\ \leq \ 
n 
\ < \ 
d.
\]
Put $t = \ell - r$, thus \mB is an $r \times t$ matrix.  
Let $\sW$ be the $\K$-vector space spanned by the columns of \mB in $\K^r$.  
Since $\sV$ contains $\sW \times \{0\}^{d - r}$, we have $\sW \hsx \cap \hsx \bk^r = \{0\}$.  
On the other hand, the columns of \mM are $\bk$-linearly independent, hence the same is true of 
\[
\begin{pmatrix}
A & &B \\[12pt]
C & &0
\end{pmatrix}
,
\]
hence too for \mB.  
Therefore
\[
t \ = \ \dim_{\bk} \hsx (\sW \hsx \cap \hsx \sE^r)
\]
and by the induction hypothesis, 
\[
t 
\ \leq \ 
r (r - 1) / 2
\]
\qquad $\implies$
\begin{align*}
\ell \ 
&=\ 
r + t
\\[12pt]
&\leq \ 
r + r(r-1)/2
\\[12pt]
&\leq \ 
n + n(n-1)/2
\\[12pt]
&=\ 
n (n+1) /2.
\end{align*}
Finally 
\[
n \ \leq \ d - 1 
\implies
\ell \ \leq \ (d-1)(d-1+1) / 2 \ = \ d (d-1) /2
\]
which completes the induction.
\end{x}
\vspace{0.3cm}

\begin{x}{\small\bf APPLICATION} \ 
Take $\K = \C$, $\bk = \Q$, and $\sE_0 = \fL$.  
Admit FDC (cf. \S49, \#1) $-$then $\sE_0$ is a $\Q$-vector subspace of $\C$ satisfying property $\binom{ A \ B}{C \ 0}$ 
(cf. \S50, \#11), so for any $\C$-vector subspace $\sV \subset \C^d$ such that $\sV \hsx \cap \hsx \Q^d = \{0\}$ there follows
\[
\dim_{\Q} \hsx (\sV \hsx \cap \hsx \fL^d) 
\ \leq \ 
n(n+1) /2.
\]

[Note: \ 
It is not known if 
\[
\trdegQ \Q(\fL) \ \geq \ 2.
\]
However the mere presence of  property $\binom{ A \ B}{C \ 0}$ is not enough to imply that there exist two algebraically independent logarithms of algebraic numbers.]
\end{x}
\vspace{0.3cm}

\begin{x}{\small\bf \un{N.B.}}  \ 
The estimate

\[
\dim_{\Q} \hsx (\sV \hsx \cap \hsx \fL^d) 
\ \leq \ 
n(n+1) /2
\]
is sharp (cf. \S38, \#7).
\end{x}
\vspace{0.3cm}

\begin{x}{\small\bf IMPLICATION} \ 
\[
\text{FDC \ $\implies$ \ 4EC}.
\]

PROOF \ 
Refer back to the proof of \#1 in \S41.  
Follow it line by line, working with $\{x_1, x_2\}$ and $\{y_1, y_2\}$ 
(drop the ``$y_3$'') $-$then $\sV = \C \bx$ contains two $\Q$-linearly independent points 
(viz. $y_1 \bx, \hsx y_2 \bx$), hence
\[
2 \ \ \leq \ \dim_{\Q} \hsx (\sV \hsx \cap \hsx \fL^2).
\]
On the other hand (cf. \#4), 
\[
\dim_{\Q} \hsx (\sV \hsx \cap \hsx \fL^2)
\ \leq \ 
1 (1 + 1) /2 
\ = \ 
1.
\]
Contradiction.
\vspace{0.3cm}

[Note: \ Recall that
\[
\text{SCHC} \ \implies \ \text{4EC}  \qquad \text{(cf. \S47, \#17)}
\]
and
\[
\text{SCHC} \ \implies \ \text{FDC} \qquad \text{(cf. \S47, \#7 and \S49, \#1)}.]
\]
\end{x}
\vspace{0.3cm}

\begin{x}{\small\bf REMARK} \ 
Under SCHC, it can be shown that a finite subset of $\fL^*$ 
consisting of $\Qbar$-linearly independent elements along with 1 is $\Qbar$-algebraicallly independent.  
Agreeing to denote this property by the symbol SFDC, we therefore have the implication
\[
\text{SCHC} \ \implies \ \text{SFDC}.
\]
One can then work with the triple $(\Qbar, \C, \fL^*)$, which thus satisfies property $\binom{ A \ B}{C \ 0}$.  
So, for any $\C$-vector subspace $\sV \subset \C^d$ of dimension $n$ such that 
$\sV \cap \Qbar^{\hsx d} = \{0\}$, the $\Qbar$-vector subspace $\sV \hsx \cap \hsx \fL^{* \hsx d}$ 
has dimension $\leq n(n+1)/2$. 
\end{x}
\vspace{0.3cm}

\begin{x}{\small\bf \un{N.B.}} \ 
\[
\text{SCHC} \ \implies \ \text{S4EC}.
\]
\end{x}
\vspace{0.3cm}


%% file: _52_on_the_equation_z_plus_exp_z.tex
\chapter{
$\boldsymbol{\S}$\textbf{52}.\quad  ON THE EQUATION $\boldsymbol{z + e^z = 0}$}
\setlength\parindent{2em}
\setcounter{theoremn}{0}
\renewcommand{\thepage}{\S52-\arabic{page}}

\ \indent 
This equation has exactly one real root.  
Can it be expressed in ``elementary'' terms?
\vspace{0.3cm}

\begin{x}{\small\bf DEFINITION} \ 
A subfield $\F$ of $\C$ is closed under $\exp$ and $\Log$ if 
\vspace{0.2cm}

\qquad \textbullet \quad 
$z \in \F \implies \exp z \in \F$
\vspace{0.2cm}

\qquad \textbullet \quad 
$z \in \F - \{0\} \implies \Log z \in \F$.
\end{x}
\vspace{0.3cm}

\begin{x}{\small\bf NOTATION} \ 
Write $\E$ for the intersection of all subfields of $\C$ that are closed under $\exp$ and $\Log$, the elements of $\E$ being the 
\un{elementary numbers}.
\index{elementary numbers}
\end{x}
\vspace{0.3cm}

\begin{x}{\small\bf CONSTRUCTION} \ 
Set $E_0 = \{0\}$ and for each $n > 0$, let $E_n$ be the set of all complex numbers obtained by applying a field operation to a pair of elements of $E_{n-1}$ or by applying $\exp$ or $\Log$ to an element of $E_{n-1}$.
\vspace{0.2cm}

[Note: \ 
Division by zero or taking the logarithm of zero are not, of course, permitted.]
\end{x}
\vspace{0.3cm}

\begin{x}{\small\bf \un{N.B.}} \ 
Therefore
\[
\Q \subset \E.
\]
\end{x}
\vspace{0.3cm}

\begin{x}{\small\bf LEMMA} \ 
\[
\E 
\ = \ 
\bigcup\limits_{n=0}^\infty \hsx E_n.
\]
\vspace{0.2cm}

[Note: \ 
Consequently, $\E$  is countable.]
\end{x}
\vspace{0.3cm}

\begin{x}{\small\bf EXAMPLE} \ 
\[
e 
\ = \ 
\exp(\exp 0) \in \E.
\]
\end{x}
\vspace{0.3cm}
\begin{x}{\small\bf EXAMPLE} \ 
\[
\sqrt{-1} 
\ = \ 
\exp\bigg(\frac{\Log (-1)}{2}\bigg) \in \E.
\]
\end{x}
\vspace{0.3cm}

\begin{x}{\small\bf EXAMPLE} \ 
\[
\pi 
\ = \ 
- \sqrt{-1} \hsx \Log (-1) \hsx \in \hsx \E.
\]
\end{x}
\vspace{0.3cm}

\begin{x}{\small\bf EXAMPLE} \ 
\[
\sqrt{2} 
\ = \ 
\exp\bigg(\frac{\elln (2)}{2}\bigg) \hsx \in \hsx \E .
\]
\end{x}
\vspace{0.3cm}

\begin{x}{\small\bf THEOREM} \ 
(Admit SCHC) \ 
The real root $\rho$ of the equation $z + e^z = 0$ is not in $\E$.
\end{x}
\vspace{0.3cm}

This is definitely not obvious and it will first be necessary to  step through some preliminaries.
\vspace{0.3cm}

\begin{x}{\small\bf NOTATION} \ 
Given a finite set
\[
A \ = \ \{\alpha_1, \ldots, \alpha_n\}
\]
of nonzero complex numbers, if $A = \emptyset$ put $A_0 =  \Q$ and  if $A \neq \emptyset$, put
\[
A_i \ = \ 
\Q(\alpha_1, e^{\alpha_1} , \ldots, \alpha_i, e^{\alpha_i}) \qquad (i \in \{1, \ldots, n\}).
\]
\end{x}
\vspace{0.3cm}

\begin{x}{\small\bf \un{N.B.}}  \ 
Each element of $A_i$ is a rational function (with rational coefficients) of the numbers
\[
\alpha_1, e^{\alpha_1} , \ldots, \alpha_i, e^{\alpha_i}.
\]
\end{x}
\vspace{0.3cm}


\begin{x}{\small\bf DEFINITION} \ 
A 
\un{tower}
\index{tower} 
is a finite set
\[
A \ = \ \{\alpha_1, \ldots, \alpha_n\}
\]
of nozero complex numbers with the property that for each $i \in \{1, \ldots, n\}$ there exists an integer $m_i > 0$ such that 
$\alpha_i^{m_i} \in A_{i-1}$ or $e^{\alpha_i m_i} \in A_{i-1}$ (or both).
\end{x}
\vspace{0.3cm}

\begin{x}{\small\bf EXAMPLE} \ 
\[
A 
\ = \ 
(\alpha_1, \alpha_2, \alpha_3)
\ = \ 
\big(\elln(2), \hsy \elln(2) / 3, \hsy \elln\big(1 + e^{(\elln(2))/ 3}\big)\big)
\]
is a tower.
\vspace{0.2cm}

[One can choose
\[
m_1 = 1, \quad m_2 = 1, \quad m_3 = 1
\]
because
\[
e^{\alpha_1} = 2 \in A_0, \quad 
e^{\alpha_2} \in A_1, \quad 
e^{\alpha_3} \in A_2.]
\]

\end{x}
\vspace{0.3cm}

\begin{x}{\small\bf DEFINITION} \ 
A 
\un{reduced tower}
\index{reduced tower}
is a tower
\[
A \ = \ \{\alpha_1, \ldots, \alpha_n\}
\]
such that $\{\alpha_1, \ldots, \alpha_n\}$ is linearly independent over $\Q$.
\end{x}
\vspace{0.3cm}

\begin{x}{\small\bf \un{N.B.}}  \ 
The tower figuring in \#14 is not reduced (in fact $\alpha_1 - 3 \alpha_2 = 0$).  
\end{x}
\vspace{0.3cm}

\begin{x}{\small\bf LEMMA} \ 
Let
\[
A \ = \ \{\alpha_1, \ldots, \alpha_n\}
\]
be a tower and suppose that $q_1, \ldots, q_n$ are nonzero integers.  Set
\[
B \ = \ \{\beta_1, \ldots, \beta_n\},
\]
where
\[
\beta_i \ = \ \frac{\alpha_i}{q_i} \qquad (i = 1, \ldots, n).
\]
Then $\forall \ i$, 
\[
A_i \subset B_i
\]
and \mB is a tower.
\vspace{0.2cm}

PROOF \ 
Since 
\[
\alpha_i \ = \ \beta_i q_i 
\quad \text{and} \quad
e^{\alpha_i} \ = \ \big(e^{\beta_i}\big)^{q_i}, 
\]
it follows that every element of $A_i$ is a rational function (with rational coefficients) of the numbers
\[
\beta_1, e^{\beta_1}, \ldots, \beta_i, e^{\beta_i}, 
\]
hence $\forall \ i$, 
\[
A_i \subset B_i \qquad \text{(cf. \#12)}.
\]
This said, let $i \in \{1, \ldots, n\}$, thus $\alpha_i^{m_i} \hsx \in \hsx A_{i-1}$ or $e^{\alpha_i m_i} \hsx \in \hsx A_{i-1}$ (or both) and put 
$n_i = m_i q_i$.
\vspace{0.2cm}

\qquad \textbullet \ 
Suppose that $\alpha_i^{m_i} \in A_{i-1}$ $-$then
\[
\beta_i^{n_i} \ = \ 
\bigg(
\frac{\alpha_i^{m_i}}{q_i^{m_i}}
\bigg)^{q_i}
\hsx \in \hsx A_{i-1} \hsx \subset \hsx B_{i-1}.
\]
\\[-.75cm]

\qquad \textbullet \ 
Suppose that $e^{\alpha_i m_i} \hsx \in \hsx A_{i-1}$ $-$then
\[
e^{\beta_i n_i} \ = \ 
e^{\alpha_i m_i} 
\hsx \in \hsx A_{i-1} \hsx \subset \hsx B_{i-1}.
\]
Therefore \mB is a tower.
\end{x}
\vspace{0.3cm}

\begin{x}{\small\bf REDUCTION PRINCIPLE} \ 
Given $\gamma \in \E$, there is a reduced tower
\[
A \ = \ \{\alpha_1, \ldots, \alpha_n\}
\]
such that $\gamma \in A_n$.
\vspace{0.2cm}

PROOF \ 
If $\gamma \in \Q$, take for \mA the empty sequence.  
If $\gamma \notin \Q$, let $\Tee(\gamma)$ be the set of all towers
\[
A \ = \ \{\alpha_1, \ldots, \alpha_n\}
\]
with the property that $\gamma \in A_n$ $-$then $\Tee(\gamma)$ is not empty and, as will now be shown, the assumption that every element of $\Tee(\gamma)$ is not reduced is a non sequitur.  
So choose a tower
\[
A \ = \ \{\alpha_1, \ldots, \alpha_n\} \in \Tee(\gamma)
\]
and take $n$ minimal $(n \geq 1)$.  
Let $i$ be the smallest integer such that $\{\alpha_1, \ldots, \alpha_i\}$ is linearly dependent over $\Q$, hence
\[
\alpha_i \ = \ \sum\limits_{j=1}^{i-1} \hsx\hsx \frac{p_j}{q_j} \alpha_j
\]
for certain integers $p_1, q_1, \ldots, p_n, q_n$.  
Consider the sequence
\[
A^\prime \ = \ 
\bigg\{
\frac{\alpha_1}{q_1}, \ldots, \frac{\alpha_{i-1}}{q_{i-1}}, \alpha_{i+1} , \ldots, \alpha_n
\bigg\}.
\]
Then the claim is  that $A^\prime \in \Tee(\gamma)$, which contradicts the minimality of $n$.  
To establish this, note that the sequence
\[
\bigg\{
\frac{\alpha_1}{q_1}, \ldots, \frac{\alpha_{i-1}}{q_{i-1}}\bigg\}
\]
is a tower (cf. \#17).  In addition, 
\[
\alpha_i \in A_{i-1}^\prime \qquad \text{(by the formula above for $\alpha_i$)}
\]
and
\[
e^{\alpha_i} \in A_{i-1}^\prime \qquad 
\text{(it is a polynomial in the numbers $e^{\alpha_1 / q_1}, \ldots,  e^{\alpha_{i-1} / q_{i-1}}$)}.
\]
But
\[
A_{i-1} \ 
\subset \ 
A_{i-1}^\prime \qquad \text{(cf. \#17)}
\]
\qquad $\implies$
\[
A_i = A_{i-1}(\alpha_i, e^{\alpha_i}) \subset A_{i-1}^\prime.
\]
Therefore the tower condition for $A^\prime$ is satisfied at the boundary between $\ds\frac{\alpha_{i-1}}{q_{i-1}}$ and $\alpha_{i+1}$ and 
\[
\gamma \in A_n \subset A_{n-1}^\prime \implies A^\prime \in \Tee(\gamma),
\]
as claimed.
\end{x}
\vspace{0.3cm}

\begin{x}{\small\bf SUBLEMMA} \ 
Suppose that
\[
A \ = \ \{\alpha_1, \ldots, \alpha_n\}
\]
is a tower $-$then $\forall \ i$,
\[
\trdegQ A_i \ \leq \ i.
\]
\vspace{0.2cm}

PROOF \ 
Start with the situation when $n = 1$, say $\{\alpha, e^\alpha\}$, and for sake of argument, assume that 
$\alpha^m \in \Q$ $-$then $\alpha$ is algebraic (consider $X^m - \alpha^m$), hence
\allowdisplaybreaks
\begin{align*}
\trdegQ \Q(\alpha, e^\alpha) \ 
&=\ 
\trdegQ \Q(e^\alpha)
\\[12pt]
&\leq \ 1.
\end{align*}
Proceed from this point by induction, the underlying hypothesis being that
\[
\trdegQ A_{i-1} \ \leq \  i - 1.
\]
Let $r_i$ stand for $\alpha_i$ or $\ds e^{\alpha_i}$ $-$then
\allowdisplaybreaks
\begin{align*}
A_i 
&=\ 
A_{i-1} (\alpha_i, e^{\alpha_i})
\\[12pt]
&=\ 
A_{i-1} (r_i).
\end{align*}
However, on general grounds (cf. \S46, \#20),
\[
\trdegQ A_{i-1}(r_i) 
\ = \ 
\trdeg_{A_{i-1}} A_{i-1}(r_i)  + \trdegQ A_{i-1},
\]
or still, 
\[
\trdegQ A_{i-1}(r_i) 
\ \leq \ 
1 + i - 1 
\ = \ 
i.
\]
I.e. : 
\[
\trdegQ A_i \ \leq \ i.
\]
\end{x}
\vspace{0.3cm}

\begin{x}{\small\bf LEMMA} \ 
(Admit SCHC) \ 
Suppose that
\[
A \ = \ \{\alpha_1, \ldots, \alpha_n\}
\]
is a reduced tower $-$then not both $\alpha_i$ and $\ds e^{\alpha_i}$ are algebraic over $A_{i-1}$.
\vspace{0.2cm}

PROOF \ 
In the notation of \S46, \#20,
\[
\trdegQ (A_i / \Q) 
\ = \ 
\trdeg_{A_{i-1}} (A_i / A_{i-1}) + \trdegQ (A_{i-1} / \Q).
\]
To get a contradiction, suppose that both $\alpha_i$ and $\ds e^{\alpha_i}$ are algebraic over $A_{i-1}$, thus 
\[
A_{i-1} (\alpha_i, e^{\alpha_i})
\]
is an algebraic extension of $A_{i-1}$, so $A_i$ is an algebraic extension of $A_{i-1}$, hence 
\[
\trdeg_{A_{i-1}} (A_i / A_{i-1}) \ = \ 0 \qquad \text{(cf. \S46, \#18)}.
\]
Therefore
\[
\trdegQ (A_i / \Q) \ = \ \trdegQ (A_{i-1} / \Q).
\]
Owing now to Schanuel, 
\[
\trdegQ (A_i / \Q) \ \geq \ i.
\]
On the other hand (cf. \#19), 
\[
\trdegQ (A_{i-1} / \Q) \ \leq \  i - 1.
\]
Contradiction.
\end{x}
\vspace{0.3cm}

\begin{x}{\small\bf \un{N.B.}} \ 
$\forall \ i$, 
\[
\trdegQ A_i \ = \ i.
\]
\end{x}
\vspace{0.3cm}

Turning finally to the proof of \#10, suppose that $\rho \in \E$ $-$then in view of \#18, there is a reduced tower
\[
A \ = \ \{\alpha_1, \ldots, \alpha_n\}
\]
such that $\rho \in A_n$.  
Obviously $\rho \notin \Q$ and it can be assumed without loss of generality that $\rho \notin A_i$ if $i < n$.
\vspace{0.1cm}

Put
\[
A^\prime \ = \ \{\alpha_1, \ldots, \alpha_n, \rho\}.
\]
Then
\[
\rho \in A_n^\prime \ = \ \Q \big(\alpha_1, e^{\alpha_1}, \ldots, \alpha_n, e^{\alpha_n}\big) \ = \ A_n
\]
and
\[
\rho + e^\rho \ = \ 0 
\implies 
e^\rho \in A_n^\prime.
\]
Accordingly $A^\prime$ (which is clearly a tower) cannot be reduced (cf. \#20).  
On the other hand, \mA is reduced, thus
\[
\rho 
\ = \ 
\sum\limits_{i=1}^n \hsx 
\frac{p_i}{q_i} \alpha_i
\]
for certain integers $p_1, q_1, \ldots, p_n q_n$.  
Here $p_n \neq 0$ since $\rho \notin A_i$ for $i < n$.  
In terms of this data
\[
\rho + e^\rho \ =  \ 0 
\implies 
\sum\limits_{i=1}^n \hsx 
\frac{p_i}{q_i} \alpha_i \ + \ 
\prod\limits_{i=1}^n \hsx
\big(e^{\alpha_i / q_i}\big)^{p_i} 
\ = \ 0.
\]
Let 
\[
B \ = \ \{\alpha_1 / q_1, \ldots, \alpha_n / q_n\}.
\]
Then \mB is a tower (cf. \#17) and since \mA is reduced, the same is true of \mB.  
But $p_n \neq 0$, hence
\[
\alpha_n / q_n \  \text{algebraic over} \ B_{n-1} \implies e^{\alpha_n / q_n}  \  \text{algebraic over} \ B_{n-1}
\]
and vice versa.  
It therefore follows that \mB cannot be reduced (cf. \#20).  
Consequently the supposition that $\rho \in \E$ has led to a contradiction.
\vspace{0.3cm}

\begin{x}{\small\bf NOTATION} \ 
Write $\ov{\E}$ for the smallest algebraically closed subfield of $\C$ that is closed under $\exp$ and $\Log$.
\end{x}
\vspace{0.2cm}

\begin{x}{\small\bf \un{N.B.}}\ 
Evidently
\[
\E \subset \ov{\E}.
\]
\end{x}
\vspace{0.2cm}

\begin{x}{\small\bf THEOREM} \ 
(Admit SCHC) \ 
Suppose that $P(X,Y) \in \Qbar[X,Y]$ is an irreducible polynomial such that
\[
\begin{cases}
\ \deg_X P \geq 1 \quad \text{per} \quad \C[Y][X]\\
\ \deg_Y P \geq 1 \quad \text{per} \quad \C[X][Y]
\end{cases}
.
\]
Assume: \ For some nonzero $\alpha \in \C$, 
\[
P(\rho,e^\alpha) \ = \ 0.
\]
\\[-21PT]
Then $\alpha \notin \ov{\E}$.
\vspace{0.1cm}

[Note: \ 
$\alpha$ is necessarily transcendental.  
For if $\alpha$ was algebraic, then the relation 
\[
P(\rho,e^\alpha) \ = \ 0
\]
\\[-21PT]
implies that $e^\alpha$ would also be algebraic, which contradicts Hermite-Lindemann (cf. \S21, \#4).]
\end{x}
\vspace{0.1cm}

\begin{x}{\small\bf APPLICATION} \ 
Take $P(X,Y) = X + Y$ and take $\alpha = \rho$ $-$then 
\allowdisplaybreaks
\begin{align*}
P(\rho,e^\rho) \ 
&=\ 
\rho + e^\rho 
\\[5pt]
&=\ 
0
\\[5pt]
&\implies
\rho \notin \ov{\E}
\\[5pt]
&\implies
\rho \notin \E,
\end{align*}
thereby recovering \#10.
\end{x}

%% file: _53_on_the_equation_P_z_exp_z.tex
\chapter{
$\boldsymbol{\S}$\textbf{53}.\quad  ON THE EQUATION $\boldsymbol{P(z,e^z) = 0}$}
\setlength\parindent{2em}
\setcounter{theoremn}{0}
\renewcommand{\thepage}{\S53-\arabic{page}}


\begin{x}{\small\bf RAPPEL} \ 
Let $f$ be an entire function.  
Assume: \ $f$ has no zeros $-$then there is an entire function $g$ such that $f = e^g$.
\vspace{0.2cm}

[Note: \ 
If $f$ is of finite order, then $g$ is a polynomial (and the order of $f$ is equal to the degree of $g$).]
\end{x}
\vspace{0.3cm}

\begin{x}{\small\bf RAPPEL} \ 
Let $f$ be an entire function.  
Assume: \ $f$ has finitely many zeros $z_1 \neq 0, \ldots, z_n \neq 0$ (each counted with multiplicity), as well as a zero of order 
$m \geq 0$ at the origin $-$then
\[
f(z) 
\ = \ 
z^m \hsx e^{g(z)} \hsx 
\prod\limits_{k=1}^n \hsx \bigg(1 - \frac{z}{z_k}\bigg),
\]
where $g(z)$ is entire.
\vspace{0.2cm}

[Note: \ 
If $f$ is of finite order, then $g$ is a polynomial (and the order of $f$ is equal to the degree of $g$).]
\end{x}
\vspace{0.3cm}

\begin{x}{\small\bf DEFINITION} \ 
A polynomial $P \in \C[X,Y]$ satisfies the \un{standard conditions} if $P$ is irreducible and 

\[
\begin{cases}
\ \deg_X P \hsx \geq \hsx 1 \quad \text{per} \quad \C[Y][X]\\[8pt]
\ \deg_Y P \hsx \geq \hsx 1 \quad \text{per} \quad \C[X][Y]
\end{cases}
.
\]

Given such a $P$, let
\[
f(z) \ = \ P(z,e^z).
\]
Then $f(z)$ has order 1.
\end{x}
\vspace{0.3cm}

\begin{x}{\small\bf LEMMA} \ 
$f(z)$ has infinitely many zeros.
\vspace{0.2cm}

PROOF \ 
Suppose that $f(z)$ has finitely many zeros $-$then there exist complex constants $A$, $B$ and a polynomial $p(X) \in \C[X]$ 
such that 
\begin{align*}
f(z) \ 
&=\ 
e^{A z + B} \hsx p(z)
\\[12pt]
&=\ 
e^{A z} e^B\hsx p(z)
\\[12pt]
&=\ 
e^{A z} q(z),
\end{align*}
where
\[
q(z) \ = \ e^B p(z) \hsx \in \hsx \C[X].
\]
But the relation 
\[
P(z,e^z) \hsx - \hsx e^{A z} q(z) \ = \ 0
\]
is possible only if $A \in \N$ (expand the data and compare coefficients), hence
\[
P(X,Y) \ = \ Y^A q(X).
\]
Since $P$ depends on both \mX and \mY, neither $Y^A$ nor $q(X)$ are equal to 1, thus $P(X,Y)$ is reducible, which contradicts the fact that $P(X,Y)$ is irreducible.
\vspace{0.2cm}

[Note: \ 
To rule out from first principles the possibility that $A = 0$, observe that the relation
\[
P(z,e^z)  \ = \ q(z)
\]
would imply that $e^z$ is algebraic (cf. \S20, \#13), whereas $e^z$ is transcendental (cf. \S20, \#15).]
\end{x}
\vspace{0.3cm}

We come now to the main result which is an illustration of the old adage 
``assume more, get more'', there being, however, a price to pay, viz. the imposition of SCHC.
\vspace{0.5cm}

\begin{x}{\small\bf THEOREM} \ 
(Admit SCHC) \ 
Suppose that $P$ satisfies the standard conditions.   
Suppose in addition that $P \in \Q[X,Y]$ $-$then 
\[
f(z) \ = \ P(z,e^z)
\]
has infinitely many $\Q$-algebraically independent zeros.
\end{x}
\vspace{0.3cm}

The proof is lengthy and will be developed in the lines that follow.
\vspace{0.3cm}

\begin{x}{\small\bf DEFINITION} \ 
A zero $\alpha \neq 0$ of $f(z)$ is said to be 
\un{generic}
\index{zero \\ generic} 
if
\[
\trdegQ \Q(\alpha, e^\alpha) \ = \ 1.
\]

[Note: \ 
Therefore the point $(\alpha, e^\alpha)$ is a generic point of the curve $C \subset \C \times \C^\times$ given by 
$P(X,Y) = 0$.]
\end{x}
\vspace{0.3cm}

\begin{x}{\small\bf LEMMA} \ 
Every zero $\alpha \neq 0$ of $f(z)$ is generic.
\vspace{0.2cm}

PROOF \ 
According to \S52, \#24, $\alpha$ is necessarily transcendental, hence
\[
\trdegQ \Q(\alpha) \ = \ 1.
\]
But
\[
P(\alpha, Y) \hsx \in \hsx \Q(\alpha)[Y],
\]
so $e^\alpha$ is algebraic over $\Q(\alpha)$, which implies that
\[
\trdegQ \Q(\alpha, e^\alpha) \ = \ 1.
\]
\end{x}
\vspace{0.3cm}

\begin{x}{\small\bf \un{N.B.}} \ 
Distinct nonzero $\alpha, \beta$ with $f(\alpha) = 0$, $f(\beta) = 0$ are not automatically algebraically independent over $\Q$,
\vspace{0.2cm}

[Take
\[
P(X,Y) \ = \ 1 + X^2 Y + Y^2.
\]
Then
\[
P(\alpha, e^\alpha) \ = \ 0 
\implies
P(-\alpha, e^{-\alpha}) \ = \ 0.]
\]
\end{x}
\vspace{0.3cm}

However: 
\vspace{0.3cm}

\begin{x}{\small\bf SUBLEMMA} \ 
(Admit SCHC) \ 
Suppose that 
\[
\begin{cases}
\ f(\alpha) \ = \ 0 \qquad (\alpha \neq 0)\\[8pt]
\ f(\beta) \ = \ 0 \qquad (\beta \neq 0)
\end{cases}
\qquad \text{and} \quad \alpha \neq \pm \beta.
\]
Then $\alpha$ and $\beta$ are algebraically independent over $\Q$.
\vspace{0.2cm}

PROOF \ 
Bear in mind that $\alpha \neq 0$, $\beta \neq 0$ are transcendental and generic (cf. \#7).  
This said, assume that $\alpha$ and $\beta$ are algebraically dependent over $\Q$ $-$then
\[
\trdegQ \Q(\alpha, \beta, e^\alpha, e^\beta) 
\ = \ 
\trdegQ \Q(\alpha, \beta) 
\ = \ 
\trdegQ \Q(\alpha) 
\ = \ 
1.
\]
Owing now to Schanuel's conjecture, $\alpha$ and $\beta$ are linearly dependent over $\Q$: 
Linear independence over $\Q$ would imply that
\[
\trdegQ \Q(\alpha, \beta, e^\alpha, e^\beta) 
\ \geq \ 
2.
\]
Accordingly choose relatively prime integers $m$ and $n$ such that $m \alpha = n \beta$ (take $n > 0$ and suppose 
momentarily that $m > 0)$.  Put $\gamma = \ds \frac{\alpha}{n}$, hence
\[
e^\alpha
\ = \ 
\big(e^\gamma\big)^n 
\quad \text{and} \quad 
e^\beta
\ = \ 
\big(e^\gamma\big)^m .
\]
For every positive integer $j$, let
\[
C_j \subset \C \times \C^\times
\]
be the curve given by
\[
 P(j X, Y^j) \ = \ 0.
\]
Then
\[
\begin{cases}
\ 0 \ = \ f(\alpha) \ = \ P(\alpha,e^\alpha) \ = \ P(n \gamma, (e^\gamma)^n)\\[8pt]
\ 0 \ = \ f(\beta) \ = \ P(\beta,e^\beta) \ = \ P(m \gamma, (e^\gamma)^m)
\end{cases}
\]
\qquad $\implies$
\[
\big(\gamma, e^\gamma) \hsx \in \hsx C_n \hsx \cap \hsx C_m.
\]
Since $C_n$ and $C_m$ have a nonempty intersection, it follows that they have a common irreducible component and this means that
\[
P(n X, Y^n) \quad \text{and} \quad  P(m X, Y^m)
\]
have a common irreducible factor.
\vspace{0.3cm}

{\small\bf FACT} \ 
The $n^\nth$ roots of unity operate transitively on the irreducible components of $C_n$ 
and the  $m^\nth$ roots of unity operate transitively on the irreducible components of $C_m$.
\vspace{0.3cm}

\qquad\textbullet \quad
Factor $P(n X, Y^n)$ into relatively prime irreducibles:
\[
P(n X, Y^n) 
\ = \ 
\prod\limits_{j = 1}^k \hsx U_j (X,Y)^{s_j}.
\]
Then it can be shown that each $U_j(X,Y)$ is of the form $U_1(X, \omega Y)$ for some $n^\nth$ root of unity $\omega$ and 
$s_1 = \cdots = s_k$, call their common value $s$, hence
\[
\deg_X P 
\ = \ 
k  s \hsx \deg_X U_1
\]
and
\[
n \hsx \deg_Y P 
\ = \ 
k  s \hsx \deg_Y U_1.
\]
\vspace{0.1cm}

\qquad\textbullet \quad
Factor $P(m X, Y^m)$ into relatively prime irreducibles:
\[
P(m X, Y^m) 
\ = \ 
\prod\limits_{i = 1}^\ell \hsx V_i (X,Y)^{t_i}.
\]
Then it can be shown that each $V_i(X,Y)$ is of the form $V_1(X, \omega Y)$ for some $m^\nth$ root of unity $\omega$ and 
$t_1 = \cdots = t_\ell$, call their common value $t$, hence
\[
\deg_X P 
\ = \ 
\ell  t \hsx \deg_X V_1
\]
and
\[
m \hsx \deg_Y  P
\ = \ 
\ell  t \hsx \deg_Y V_1.
\]

It can be assumed that 
\[
U_1(X,Y) \ = \ V_1(X,Y),
\]
the common irreducible factor of $P(n X, Y^n)$ and $P(m X, Y^m)$ $-$then
\begin{align*}
k  s \hsx \deg_X U_1 \ 
&=\ 
\deg_X P
\\[12pt]
&=\ 
\ell  t \hsx \deg_X V_1
\\[12pt]
&=\ 
\ell  t \hsx \deg_X U_1.
\end{align*}
But
\[
\deg_X P \ \neq \  0 
\implies 
k s \ = \ \ell  t 
\ \neq \ 
0.
\]
Next
\begin{align*}
n \hsx \deg_Y P \ 
&=\ 
k  s \deg_Y U_1
\\[12pt]
&=\ 
\ell  t \hsx \deg_Y U_1
\\[12pt]
&=\ 
\ell  t \hsx \deg_Y V_1
\\[12pt]
&=\ 
m \hsx \deg_Y P.
\end{align*}
But
\[
\deg_Y P \ \neq \ 0
\implies 
n \ = \ m,
\]
contradicting the assumption that $m$, $n$ are relatively prime.
\vspace{0.2cm}

[Note: \ 
To treat the case when $m < 0$, consider the polynomial 
\[
T(X,Y) 
\ = \ 
Y^{-m \hsx \deg_Y P} \hsx P(m X, Y^m).
\]
Then
\[
\deg_X T \ = \ \deg_X P
\]
and 
\[
\deg_Y T \ = \ -m \deg_Y P.
\]
So as above, 
\[
m \alpha \ = \ n \beta 
\implies
-n \alpha \ = \ n \beta 
\implies
- \alpha \ = \ \beta 
\implies
\alpha \ = \ -\beta ,
\]
which is forbidden by hypothesis.]
\end{x}
\vspace{0.3cm}

\begin{x}{\small\bf DEFINITION} \ 
Under the assumptions of \#5, $P$ is said to be 
\un{primitive}
\index{primitive polynomial} 
if $\forall \ n \in \N$, the curve $C_n$ given by 
\[
P(n X, Y^n)
\ = \ 0
\]
is irreducible.
\end{x}
\vspace{0.3cm}

\begin{x}{\small\bf LEMMA} \ 
(Admit SCHC) \ 
Suppose that \mP is primitive and let $\alpha_1, \ldots, \alpha_n$ be nonzero zeros of 
$f(z) = P(z,e^z)$ subject to $\alpha_i \neq \pm \alpha_j$ for all $i \neq j$ 
$-$then $\alpha_1, \ldots, \alpha_n$ are algebraically  independent over $\Q$.
\vspace{0.2cm}

\begin{spacing}{1.3}
PROOF \ 
Searching for a contradiction, the first step is to tabulate the data.  
So assume that over $\Q$ there exists an algebraically dependent collection 
$\alpha_1, \ldots, \alpha_n, \alpha_{n+1}$ of $n + 1$ nonzero zeros of $f$ such that 
$\alpha_i \neq \pm \alpha_j$ for all $i \neq j$ and take $n$ minimal.  
In view of \#9, two such zeros are algebraically  independent over $\Q$, hence $n \geq 2$, and, by the minimality of $n$, 
the collection  $\alpha_1, \ldots, \alpha_n$ is algebraically  independent over $\Q$, hence
\end{spacing}
\[
\trdegQ \Q\big(\alpha_1, \ldots, \alpha_{n+1}, e^{\alpha_1}, \ldots, e^{\alpha_{n+1}}\big) 
\ = \ 
n 
\ < \ 
n + 1.
\]
Meanwhile, by Schanuel, if $\alpha_1, \ldots, \alpha_n, \alpha_{n+1}$ were linearly independent over $\Q$, then 
\[
\trdegQ \Q\big(\alpha_1, \ldots, \alpha_{n+1}, e^{\alpha_1}, \ldots, e^{\alpha_{n+1}}\big) 
\ \geq \ 
n + 1.
\]
Since this cannot be, it follows that there exist nonzero integers $m_1, \ldots, m_n, m$ with no common factor such that
\[
\sum\limits_{k=1}^n \hsx m_k \alpha_k 
\ = \ 
m \alpha_{n + 1} 
\qquad (m > 0).
\]
\begin{spacing}{1.3}
Put $\gamma_k = \ds \frac{\alpha_k}{m}$.  
Let $C \subset \C \times \C^\times$ be the curve defined by $P(X,Y) = 0$ and let 
$C_m \subset \C \times \C^\times$ be the curve defined by $P(m X, Y^m) = 0$.  
Since \mP is primitive, $C_m$ is irreducible and since $\alpha_1, \ldots, \alpha_n$ are algebraically  independent over $\Q$, 
the same is true of $\gamma_1, \ldots, \gamma_n$.  
Therefore $\big(\gamma_1, e^{\gamma_1}\big), \ldots, \big(\gamma_n, e^{\gamma_n}\big)$ are generic points in $C_m$.  
Moreover
\end{spacing}
\[
\trdegQ \Q(\gamma_1, e^{\gamma_1}, \ldots, \gamma_n, e^{\gamma_n}) 
\ = \ 
n.
\]

\end{x}
\vspace{0.3cm}

\begin{x}{\small\bf CONSTRUCTION} \ 
Define a map
\[
\phi : \big(\C \times \C^\times\big)^n \ra \C \times \C^\times
\]
by the prescription
\[
(x_1, y_1, \ldots, x_n, y_n) 
\ra 
\bigg(
\sum\limits_{k=1}^n \hsx m_k x_k, \hsx \prod\limits_{k=1}^n \hsx y_k^{m_k}
\bigg).
\]
Then
\begin{align*}
\phi(\gamma_1, e^{\gamma_1}, \ldots, \gamma_n, e^{\gamma_n}) \ 
&=\ 
\bigg(
\sum\limits_{k=1}^n \hsx m_k \gamma_k, \hsx \prod\limits_{k=1}^n \hsx e^{\gamma_k m_k}
\bigg)
\\[12pt]
&=\ 
\bigg(
\sum\limits_{k=1}^n \hsx \frac{m_k \alpha_k}{m}, \hsx \prod\limits_{k=1}^n \hsx e^{\textstyle\frac{m_k \alpha_k}{m}}
\bigg)
\\[12pt]
&=\ 
\big(\alpha_{n+1}, e^{\alpha_{n+1}}\big),
\end{align*}
a generic point in \mC, hence $\phi$ maps $\big(C_m\big)^n$ to \mC.  
So if $z_1, \ldots, z_n$ are zeros of $f$, then the pairs
\[
\bigg(\frac{z_1}{m}, e^{\frac{z_1}{m}}\bigg), \ldots, \bigg(\frac{z_n}{m}, e^{\frac{z_n}{m}}\bigg)
\]
lie in $C_m$, from which it follows that the sum
\[
\sum\limits_{k=1}^m \hsx \frac{m_k}{m} z_k
\]
is a zero of $f$.  
In particular:
\[
\alpha \ \equiv \ 
\frac{m_1 + m_2}{m}\hsx  \alpha_1 + \frac{m_3}{m}\hsx  \alpha_3 + \cdots + \frac{m_n}{m} \hsx \alpha_n
\]
is a zero of $f$ (take $z_1 = z_2 = \alpha_1$ and $z_k = \alpha_k$ $(k > 2)$).   
\vspace{0.3cm}

\qquad \un{$n > 2$:} \quad \ 
In this situation, the collection $\alpha_1, \alpha_3, \ldots, \alpha_n, \alpha$ is algebraically dependent over $\Q$ and consists of $n$ nonzero zeros of $f$, contradicting the minimality of $n$.
\vspace{0.2cm}

[Note: \ 
The condition $n > 2$ implies that $\alpha$ is nonzero and $\alpha \neq \pm \alpha_i \ \forall \ i$.]
\vspace{0.3cm}

\qquad \un{$n = 2$:} \quad \ 
It is a question of dealing with the collection $\alpha_1$, $\alpha_2$, $\alpha_3$ of $\Q$-algebraically dependent nonzero zeros of $f$ such that $\alpha_i \neq \pm \alpha_j$ for all $i \neq j$ satisfying
\[
m_1 \alpha_1 + m_2 \alpha_2 \ = \ m \alpha_3,
\]
where, as above,
\[
\alpha 
\ = \ 
\frac{m_1 + m_2}{m} \hsx \alpha_1
\]
is a zero of $f$.  
The claim then is that such a scenario is impossible.  
To this end, it will be shown below that each of the following conditions leads to a contradiction.
\[
(1) \quad m_1 + m_2 = 0;
\quad 
(2) \quad m_1 + m_2 = m;
\quad
(3) \quad m_1 + m_2 = -m.
\]
Therefore
\[
\alpha \neq 0 
\  \text{(cf. (1))}; \quad
\alpha \neq \alpha_1 
\  \text{(cf. (2))}; \quad
\alpha \neq -\alpha_1 
\  \text{(cf. (3))}. \quad
\]
Consequently $\alpha$ and $\alpha_1$ are algebraically independent over $\Q$ (cf. \#9).  
But this is nonsense since $\alpha$ and $\alpha_1$ are linearly dependent over $\Q$:
\[
1 \cdot \alpha - q \cdot \alpha_1 \ = \ 0 
\qquad \bigg(q = \frac{m_1 + m_2}{m} \in \Q\bigg).
\]

\qquad \un{Ad(1)}  \ 
$(m_1 + m_2 = 0)$: \quad 
To begin with, note that $\ds\frac{m_1}{m}\alpha_1$ and $\alpha_1$ are nonzero $\Q$-algebraicallly dependent zeros  of $f$, hence by \#9, 
\[
\frac{m_1}{m} \alpha_1 
\ = \ 
\pm \alpha_1 
\implies
m_1 \ = \ \pm m.
\]
To pin things down, take $m_1 = +m$ $-$then
\allowdisplaybreaks
\begin{align*}
m_1 \alpha_1 + m_2 \alpha_2 = m \alpha_3 \ 
&\implies 
m \alpha_1 - m \alpha_2 = m \alpha_3 
\\[12pt]
&\implies 
\alpha_1 = \alpha_2 + \alpha_3.
\end{align*}
Now interchange the roles of $\alpha_1$ and $\alpha_2$ to get
\[
\alpha_2 \ = \ \alpha_1 + \alpha_3
\]
or still, 
\[
\alpha_2 
\ = \ 
\alpha_2 + \alpha_3 + \alpha_3 
\implies
0 
\ = \ 
2 \alpha_3.
\]
Contradiction.
\vspace{0.5cm}

\qquad \un{Ad(2)}  \ 
$(m_1 + m_2 = m)$: \quad 
By switching the roles of the variables and multiplying by $-1$ if necessary, it can be assumed that 
$\abs{m} \geq \abs{m_1}$, $\abs{m_2}$ and $m > 0$, $m_1 > 0$.  
Construct a sequence $\{z_k\}$ of zeros of $f$ by the following procedure: 
Take $z_1 = \alpha_1$ and via recursion, take
\[
z_{k+1} 
\ = \ 
\frac{m_1}{m} \hsx z_k \hsx + \hsx \frac{m_2}{m} \alpha_2.
\]
Then the fact that
\[
\frac{m_2}{m} 
\ = \ 
1 - \frac{m_1}{m}
\]
leads to the relation
\[
z_{k+1} 
\ = \ 
\bigg(\frac{m_1}{m}\bigg)^k \alpha_1 \ + \ \bigg(1 - \bigg(\frac{m_1}{m}\bigg)^k\bigg) \alpha_2.
\]
Since
\[
0 
\ < \ 
\frac{m_1}{m}
\ < \ 1,
\]
the coefficient $\bigg(\ds\frac{m_1}{m}\bigg)^k$ of $\alpha_1$ takes a different value for each $k$, thus thanks to the $\Q$-algebraic independence of $\alpha_1$ and $\alpha_2$, the sequence $\{z_k\}$ assumes infinitely many distinct values.  
Put
\[
M 
\ = \ 
\max \{\abs{\alpha_1}, \abs{\alpha_2}\}.
\]
Then
\[
\abs{z_{k+1}}
\ \leq \ 
\abs{\bigg(\frac{m_1}{m}\bigg)^k} \hsx M 
\hsx + \hsx 
\abs{1 - \bigg(\frac{m_1}{m}\bigg)^k} \hsx M 
\ \leq \
2 M.
\]
But this means that the entire function $f$ has infinitely many zeros in the disc of radius $2M$ centered at the origin, so $f = 0$, a contrdiction.
\vspace{0.5cm}

\qquad \un{Ad(3)}  \ 
$(m_1 + m_2 = -m)$: \quad 
Let $s = \ds \frac{m_1}{m}$ $-$then 
\[
\frac{m_2}{m}
\ = \ 
-(1 + s)
\]
and
\[
m \alpha_3 
\ = \ 
m_1 \alpha_1 + m_2 \alpha_2
\]
\qquad $\implies$
\begin{align*}
\alpha_3 \ 
&=\ 
\frac{m_1}{m} \hsx \alpha_1 + \frac{m_2}{m} \hsx \alpha_2
\\[12pt]
&=\ 
s \hsx \alpha_1 - (1 + s) \alpha_2.
\end{align*}
On the other hand,
\[
s \hsx \alpha_3 - (1 + s) \alpha_2
\]
is a zero of $f$.  And
\begin{align*}
w \
&\equiv 
s \hsx \alpha_3 - (1 + s) \alpha_2
\\[12pt]
&=\ 
s(s \alpha_1 - (1 + s)\alpha_2) - (1 + s)\alpha_2
\\[12pt]
&=\ 
s^2 \hsx \alpha_1 - s(1 + s)\alpha_2  - (1 + s)\alpha_2 
\\[12pt]
&=\ 
s^2 \hsx \alpha_1 - (1 + s) (s \alpha_2 + \alpha_2)
\\[12pt]
&=\ 
s^2 \hsx \alpha_1 - (1 + s)^2 \alpha_2.
\end{align*}
Now treat $\alpha_1$, $\alpha_2$, $w$ as a collection of $\Q$-algebraically dependent nonzero zeros of $f$.  
Invoking the earlier analysis, we thus have
\[
s^2 - (1 + s)^2 
\ = \ 
-2s - 1
\ = \ 
0 \ \text{or} \ \pm1.
\]

\qquad \textbullet \quad 
If $- 2s - 1 = 1$, then
\begin{align*}
s = -1 
&\implies
-1 = \frac{m_1}{m}
\\[12pt]
&\implies
- m = m_1
\\[12pt]
&\implies
m_1 + m_2 = m_1
\\[12pt]
&\implies
m_2 = 0.
\end{align*}
So
\begin{align*}
m \alpha_3 \ 
&=\ 
m_1 \alpha_1 + m_2 \alpha_2 
\\[12pt]
&=\ 
m_1 \alpha_1 
\\[12pt]
&=\ 
-m \alpha_1
\\[12pt]
\implies 
\alpha_3 &= -\alpha_1.
\end{align*}
Contradiction.
\vspace{0.3cm}

\qquad \textbullet \quad 
If $- 2s - 1 = -1$, then
\[
s \ = \ 0 
\implies m_1 \ = \ 0.
\]
So
\begin{align*}
m \alpha_3 \ 
&=\ 
m_1 \alpha_1 + m_2 \alpha_2 
\\[12pt]
&=\ 
m_2 \alpha_2
\\[12pt]
&=\ 
-m \alpha_2
\\[12pt]
\implies
\alpha_3 &= -\alpha_2.
\end{align*}
Contradiction.
\vspace{0.3cm}

\qquad \textbullet \quad
If $- 2s - 1 = 0$, then
\begin{align*}
s^2 - (1 + s)^2 \ 
&=\ 
s^2 - (1 + 2s + s^2)
\\[12pt]
&=\ 
-1 -2s 
\\[12pt]
&=\ 
0.
\end{align*}
So matters reduce to ``$m_1 + m_2 = 0$'' $\ldots$ .
\end{x}
\vspace{0.3cm}

\begin{x}{\small\bf \un{N.B.}} \ 
It won't hurt to repeat: \ $P \in \Q[X,Y]$ satisfies the standard conditions and 
\[
f(z) \ = \ 
P(z,e^z)
\]
has infinitely many zeros (cf. \#4).
\end{x}
\vspace{0.3cm}

\begin{spacing}{1.3}
PROOF \ 
OF \ \#5 \quad In view of \#11, it can be assumed that \mP is not primitive.  
Choose, accordingly, an $n \in \N$ such that $C_n$ is reducible (cf. \#10) $-$then $C_n$ has an irreducible component defined by some polynomial $P_n(X,Y) \in \Q[X,Y]$ depending on both \mX and \mY and
\end{spacing} 
\[
0 \ < \ 
\deg_X P_n 
\ < \ 
\deg_X P.
\]
Noting that $\deg_X P > 1$, proceed by induction on $\deg_X P$, supposing that for all irreducible polynomials $T(X,Y) \in \Q[X,Y]$ satisfying the standard conditions such that 
\[
\deg_X T 
\ < \ 
\deg_X P
\]
the entire function 
\[
T(z, e^z)
\]
has infinitely many $\Q$-algebraically independent zeros $-$then by hypothesis, the entire function 
\[
f_n(z) 
\ = \ 
P_n(z,e^z)
\]
has infinitely many $\Q$-algebraically independent zeros, say $z_1, z_2, \ldots$, \ .  But
$P_n(X,Y)$ is a factor of $P(n X,Y^n)$, hence 
\[
f(n z_k) 
\ = \ 
P(n z_k, e^{n z_k}) 
\ = \ 0 
\qquad 
(k = 1, 2, \ldots).
\]
Therefore 
\[
n z_1, n z_2, \ldots 
\]
is an infinite collection of $\Q$-algebraically independent zeros of $f$.

\begin{x}{\small\bf REMARK} \ 
The result remains valid if $\Q$ is replaced by $\Qbar$, i.e., granted SCHC, if $P \in \Qbar[X,Y]$ satisfies the standard conditions, then
\[
f(z) \ = \ 
P(z,e^z)
\]
has infinitely many $\Q$-algebraically independent zeros.
\end{x}
\vspace{0.3cm}

\begin{x}{\small\bf EXAMPLE} \ 
(Admit SCHC) \ 
Consider $P(X,Y) = X - Y$ $-$then the entire function
\[
f(z) \ = \ 
P(z,e^z) \ = \ z - e^z
\]
has infinitely many $\Q$-algebraically independent zeros, thus the exponential function $e^z$ has infinitely many 
$\Q$-algebraically independent fixed points (cf. \S52, \#10).
\end{x}
\vspace{0.3cm}

\begin{x}{\small\bf THEOREM} \ 
(Admit SCHC) \ 
Suppose that $\K \subset \C$ is a finitely generated field $-$then for any $P \in \K[X,Y]$ satisfying the standard conditions, the equation 
\[
P(z,e^z) \ = \ 0
\]
has a solution $\alpha$ generic over $\K$:
\[
\trdeg_{\K} \hsx \K(\alpha, e^\alpha) \ = \ 1.
\]

[This was proved in 2014 by V. Mantova.]
\end{x}
\vspace{0.3cm}

\begin{x}{\small\bf APPLICATION} \ 
(Admit SCHC) \ 
\[
\#16 \implies \#14.
\]

[Start with the field $\K$ obtained by adjoining the coefficients of \mP to \mQ.  Choose $\alpha$ per supra.]
\end{x}
\vspace{0.3cm}

Here is a word or two on the proof of \#16.   
The key is to show that $P(z,e^z)$ has only finitely many zeros in $\ov{\K}$, the algebraic closure of $\K$ 
(this forces the other zeros to be generic over $\K$).  
The point of departure for this is the following result.
\vspace{0.3cm}

\begin{x}{\small\bf LEMMA} \ 
(Admit SCHC) \ 
There exists a finite dimensional $\Q$-vector space $F \subset \ov{\K}$ containing all the zeros of $P(z,e^z)$ in $\ov{\K}$.
\vspace{0.2cm}

[Without loss of generality, add to $\K$ the coefficients of \mP so that \mP is defined over $\K$.  
Recall that for any 
$\bz = (z_1, \ldots, z_n)$, 
\[
\trdegQ \Q(\bz,e^{\bz}) 
\ \geq \ 
\lindim_{\Q} \hsx \bz 
\qquad \text{(cf. \S47, \#24)}.
\]
If now each $z_i \in \bz$ is in $\ov{\K}$ and $P(z_i, e^{z_i}) = 0$, then $e^{z_i} \in \ov{\K}$.]
\end{x}
\vspace{0.3cm}

%% file: _54_zilber_fields.tex
\chapter{
$\boldsymbol{\S}$\textbf{54}.\quad  ZILBER FIELDS}
\setlength\parindent{2em}
\setcounter{theoremn}{0}
\renewcommand{\thepage}{\S54-\arabic{page}}

\ \indent 
These are fields subject to the following conditions.
\vspace{0.2cm}

\qquad \textbullet \quad (EAC)
\vspace{0.2cm}

\qquad \textbullet \quad (STD)
\vspace{0.2cm}

\qquad \textbullet \quad (SCHP)
\vspace{0.2cm}

\qquad \textbullet \quad (SEACP) $\subset$ (EACP)
\vspace{0.2cm}

\qquad \textbullet \quad (CCP)
\vspace{0.2cm}

The meaning of these abbreviations will be explained below.
\vspace{0.3cm}

\begin{x}{\small\bf DEFINITION} \ 
An 
\un{\mE-field}
\index{\mE-field} 
is a field $(\K, \hsy +, \hsy\cdot,\hsy  0, \hsy 1)$ of characteristic 0 equipped with a surjective homomorphism \mE from its additive group $(\K, \hsy +)$ to its multiplicative group $(\K^\times, \hsy \cdot)$, thus
\[
\forall \ x, y \in \K, \quad E(x+y) \ = \ E(x) \cdot E(y)
\]
and $E(0) = 1$.
\end{x}
\vspace{0.3cm}

\begin{x}{\small\bf EXAMPLE} \ 
\vspace{0.2cm}
To exhibit an \mE-field, take $\K = \R$, take $a > 0$, and equip it with the exponential function to base $a$, i.e., 
\[
E(x) \ = \exp_a(x) \ = \ a^x \qquad (x \in \R) \qquad \text{(cf. \S5, \#1)}.
\]
\vspace{0.2cm}

[Note: \ 
Denote this setup by the symbol $\R_{\exp}$ when $a = e$.] 
\end{x}
\vspace{0.3cm}

\begin{x}{\small\bf DEFINITION} \ 
An $E$-field $\K$ is an 
\un{EAC-field}
\index{EAC-field} 
if $\K$ is algebraically closed.
\end{x}
\vspace{0.3cm}

\begin{x}{\small\bf EXAMPLE} \ 
To exhibit an EAC-field, take $\K = \C$ and equip it with the usual exponential function $z \ra e^z$.
\vspace{0.2cm}

[Note: \ 
Denote this setup by the symbol $\C_{\exp}$.] 
\end{x}
\vspace{0.3cm}

\begin{x}{\small\bf \un{N.B.}}  \ 
If $\K$ is an \mE-field, then $\Q$ can be considered as a subfield of $\K$, since $\K$ has characteristic 0.
\end{x}
\vspace{0.3cm}

\begin{x}{\small\bf DEFINITION} \ 
Suppose that $\K$ is an \mE-field $-$then the kernel of the exponential map, i.e., 
\[
\{x \in \K : E(x) = 1\}, 
\]
is said to be 
\un{standard}
\index{standard (kernel of the exponential map of an \mE field)} 
(STD)
\index{STD} 
if it is an infinite cyclic group generated by a transcendental element $\alpha$, thus
\[
\Ker(E) \ = \ \alpha \Z.
\]

[Note: \ 
$\alpha$ is transcendental provided that it is not the root of a nonzero polynomial with coefficients in the copy of $\Q$ in $\K$.]
\end{x}
\vspace{0.3cm}

\begin{x}{\small\bf EXAMPLE} \ 
Take $\K = \C_{\exp}$ $-$then the kernel of the exponential map is $2\pi \sqrt{-1} \hsx\hsy \Z$, hence is standard (take $\alpha = 2 \pi \sqrt{-1}$).
\end{x}
\vspace{0.3cm}

\begin{x}{\small\bf DEFINITION} \ 
Suppose  that $\K$ is an \mE-field  $-$then to say that $\K$ has 
\un{Schanuel's property}
\index{Schanuel's property} 
(SCHP)
\index{SCHP} 
means that if $x_1, \ldots, x_n$ are elements of $\K$ which are linearly independent over $\Q$, then the field 
\[
\Q(x_1, \ldots, x_n, E(x_1), \ldots, E(X_n))
\]
has transcendence degree $\geq n$ over $\Q$.
\vspace{0.2cm}

[Note: \ 
When $\K = \C$, SCHP is, of course, conjectural (SCHC).]
\end{x}
\vspace{0.3cm}

\begin{x}{\small\bf NOTATION} \ 
Given an \mE-field $\K$, transcribe \S47, \#20 from $\C$ to $\K$ and given $\bx$, put
\[
\delta_A (\bx) \ = \ \trdegQ \Q(\bx, E(\bx)) - \lindim_{\Q} \bx,
\]
\end{x}
the 
\un{predimension}
\index{predimension} 
of $\bx$ (cf. \S47, \#26).
\vspace{0.5cm}

Therefore SCHP per $\K$ is the claim that $\forall \ \bx$, 
\[
\delta_A (\bx) \ \geq \ 0.
\]

\begin{x}{\small\bf NOTATION} \ 
(Admit SCHP) \ 
Given an \mE-field $\K$ and a finite set $X \subset \K$, view \mX as a tuple $-$then $\delta_{\K}(X) \geq 0$ and the 
\un{dimension of \mX in $\K$}
\index{dimension of \mX in $\K$}
is 
\[
\dim_{\K} (X) 
\ = \ 
\inf\limits_{\#Y < \infty} \ \{\delta_{\K}(Y) : X \subset Y \subset \K\}.
\]
\end{x}
\vspace{0.3cm}

\begin{x}{\small\bf DEFINITION} \ 
(Admit SCHP) \ 
Let $\K$ and $\LL$ be \mE-fields $-$then $\LL$ is a 
\un{strong extension}
\index{strong extension (of \mE-fields}
of $\K$ if $\K \subset \LL$ and 
\[
\dim_{\K}(X) \ = \ \dim_{\LL}(X) 
\]
for all $X \subset \K$, where \mX is finite.
\end{x}
\vspace{0.3cm}

\begin{x}{\small\bf THEOREM} \ 
(Admit SCHP) \ 
$\C_{\exp}$ is not a strong extension of $\R_{\exp}$.
\vspace{0.3cm}

PROOF \ 
It will be shown that
\[
\dim_{\R}(\pi) \ \neq \ \dim_{\C}(\pi).
\]
Owing to Nesterenko (cf. \S20, \#10):
\allowdisplaybreaks
\begin{align*}
\text{\textbullet} \quad \delta_{\R}(\pi) \ 
&=\ 
\trdegQ \Q(\pi, e^\pi) - \lindim_{\Q} (\pi) \hspace{5cm}
\\[12pt]
&=\ 
2 - 1
\\[12pt]
&=\ 
1.
\end{align*}
\allowdisplaybreaks
\begin{align*}
\text{\textbullet} \quad \delta_{\C}(\pi,\pi \sqrt{-1}) \ 
&=\ 
\trdegQ \Q(\pi,\pi \sqrt{-1}, e^\pi, e^{\pi \sqrt{-1}}) - \lindim_{\Q} (\pi, \pi \sqrt{-1}) 
\\[12pt]
&=\ 
\trdegQ \Q(\pi,\pi \sqrt{-1}, e^\pi, -1) - 2
\\[12pt]
&=\ 
\trdegQ \Q(\pi,\pi \sqrt{-1}, e^\pi) - 2
\\[12pt]
&=\ 
\trdegQ \Q(\pi,e^\pi) - 2
\\[12pt]
&=\ 
2 - 2
\\[12pt]
&=\ 
0,
\end{align*}
$\pi \sqrt{-1}$ being algebraic over $\Q(\pi)$.  Therefore
\[
\dim_{\C} (\pi) \ = \ 0.
\]

If now $\C_{\exp}$ was a strong extension of $\R_{\exp}$, then we'd have
\[
\dim_{\R} (\pi) \ = \ 0,
\]
so there would be a finite subset $X \subset \R$ with $\pi \in X$ such that $\delta_{\R}(X) = 0$.  
Explicate: 
\[
X \ = \ \{\pi, x_1, \ldots, x_n\}
\]
and suppose that 
\[
\lindim_{\Q} X \ = \ k + 1.
\]
Write
\allowdisplaybreaks
\begin{align*}
0 \ 
=\ 
&\delta_{\R} (\pi, x_1, \ldots, x_n)
\\[12pt]
=\ 
&\trdegQ \Q(\pi, x_1, \ldots, x_n, e^{\textstyle\pi}, e^{\textstyle x_1}, \ldots, e^{\textstyle x_n})
- 
\lindim_{\Q} (\pi, x_1, \ldots, x_n)
\\[12pt]
\implies \hspace{1cm}&
\\[12pt]
&\trdegQ \Q(\pi, x_1, \ldots, x_n, e^{\textstyle\pi}, e^{\textstyle x_1}, \ldots, e^{\textstyle x_n}) \ = \ k + 1
\\[12pt]
\implies \hspace{1cm}&
\\[12pt]
&\trdegQ \Q(\pi, \pi \sqrt{-1}\hsx, x_1, \ldots, x_n, e^{\textstyle\pi}, e^{\textstyle \pi \sqrt{-1}}, e^{\textstyle x_1}, \ldots, e^{\textstyle x_n}) 
\ = \ k + 1.
\end{align*}
On the other hand, thanks to Schanuel, 
\allowdisplaybreaks
\begin{align*}
\trdegQ &\Q(\pi, \pi \sqrt{-1} \hsx, x_1, \ldots, x_n, e^{\textstyle\pi}, e^{\pi \sqrt{-1}}, e^{\textstyle x_1}, \ldots, e^{\textstyle x_n}) \  
\\[12pt]
&\hspace{2cm} \geq \ \lindim_{\Q} (\pi, \pi \sqrt{-1}\hsx, x_1, \ldots, x_n)
\\[12pt]
&\hspace{2cm} =\ 
k + 2.
\end{align*}
Contradiction.
\end{x}
\vspace{0.3cm}

The next definition, viz. that of strong exponential closure, is on the technical side.
\vspace{0.5cm}

Let $\K$ be an EAC-field.  Put $G = \K \times \K^\times$ $-$then \mG is a $\Z$-module: 
\vspace{0.2cm}

\qquad(\textbullet): \quad $\Z \times G \ra G$
\vspace{0.2cm}

\hspace{1.8cm}  $m \cdot (x,y) = (m x, y^m)$.
\vspace{0.2cm}

This action can be generalized to matrices with integer coefficients:
\vspace{0.2cm}

\qquad(\textbullet): \quad
$M_{n \times n} (\Z) \times G^n \ra G^n$,
\vspace{0.2cm}

\noindent where a matrix $M = [m_{i j}]$ sends
\[
(x_1, \ldots, x_n, y_1, \ldots, y_n)
\]
to
\[
\bigg( 
\sum\limits_{j=1}^n \hsx m_{1 j} x_j, \ldots, \sum\limits_{j=1}^n \hsx m_{n j} x_j; \ 
\prod\limits_{j=1}^n \hsx y_j^{m_{1 j}}, \ldots, \prod\limits_{j=1}^n \hsx y_j^{m_{n j}}
\bigg) .
\]
\vspace{0.3cm}

\begin{x}{\small\bf NOTATION} \ 
If $V \subset G^n$, write $M \cdot V$ for its image and observe that if \mV is a subvariety of $G^n$, then so is $M \cdot V$.
\end{x}
\vspace{0.3cm}

\begin{x}{\small\bf DEFINITION} \ 
A subvariety $V \subset G^n$ satisfies the 
\un{dimension condition}
\index{dimension condition (subvariety)}
if for all $M \in M_{n \times n} (\Z)$, 
\[
\dim M \cdot V \ \geq \ \rank M.
\]

[Note: \ 
In particular, $\dim V \geq n$.]
\end{x}
\vspace{0.3cm}

\begin{x}{\small\bf DEFINITION} \ 
A subset \mV of $G^n$ is 
\un{additively free}
\index{additively free} 
if \mV is not contained in a set given by equations of the form
\[
\bigg\{(\bx,\by) : \hsx \sum\limits_{i=1}^n \hsx m_i x_i \hsx = \hsx a\bigg\},
\]
where the $m_i \in \Z$ are not all zero and $a \in \K$.
\end{x}
\vspace{0.3cm}

\begin{x}{\small\bf DEFINITION} \ 
A subset \mV of $G^n$ is 
\un{multiplicatively free}
\index{multiplicatively free} 
if \mV is not contained in a set given by equations of the form
\[
\bigg\{(\bx,\by) : \hsx \prod\limits_{i=1}^n \hsx y_i^{m_i} \hsx = \hsx b\bigg\},
\]
where the $m_i \in \Z$ are not all zero and $b \in \K^\times$.
\end{x}
\vspace{0.3cm}

\begin{x}{\small\bf \un{N.B.}} \ 
Call \mV 
\un{free}
\index{free}
if \mV is both additively and multiplicatively free.
\end{x}
\vspace{0.3cm}

\begin{x}{\small\bf DEFINITION} \ 
A subvariety $V \subset G^n$ is 
\un{admissible}
\index{admissible (subvariety)} 
if \mV is irreducible, satisfies the dimension condition, and is free.
\end{x}
\vspace{0.3cm}

\begin{spacing}{1.4}
\begin{x}{\small\bf DEFINITION} \ 
Suppose that $\K$ is an EAC-field $-$then $\K$ has the 
\un{exponential} \un{algebraic closure property}
\index{exponential algebraic closure property}
(EACP) 
\index{EACP} 
if for all admissible subvarieties \mV of $G^n$ that are defined over $\K$ and of dimension $n$, there is an $\bx$ in $\K^n$ such that 
$(\bx,E(\bx)) \in V$.  
\vspace{0.2cm}

[Note: \ 
Therefore $\K$ is exponentially algebraically closed iff each such variety \mV intersects the graph of exponentiation.]
\end{x}
\end{spacing}
\vspace{0.3cm}

\begin{x}{\small\bf REMARK} \ 
(Admit EACP) \ 
It can be shown that there are infinitely many $\Q$-algebraically independent $\bx$ such that $(\bx,E(\bx)) \in V$.  
\end{x}
\vspace{0.3cm}

\begin{spacing}{1.4}
\begin{x}{\small\bf EXAMPLE} \ 
(Admit SCHC) \ 
Take $\K = \C_{\exp}$ $-$then it is unknown whether EACP obtains in general but the simplest case, namely when $n = 1$, can be dealt with.  
To see how this goes, recall that a variety \mV in $\C^2$ is the set of common zeros of a collection of polynomials in $\C[X,Y]$ and, in fact, is the zero set of a single 
polynomial, i.e., given \mV, there is a polynomial $P(X,Y) \in \C[X,Y]$ such that
\[
V \ = \ Z(P) \ = \ \{(X,Y) \in \C \times \C : P(X,Y) = 0\}.
\]
And \mV is irreducible iff this is so of \mP.  
Working with $V \subset \C \times \C^\times$ (being interested only in solutions to $P(z,e^z) = 0$), transfer matters from \mV to \mP by imposing the standard conditions on \mP (cf. \S53, \#3) $-$then \mV is admissible.  
E.g.: \ To check freeness, $\forall$ nonzero $m \in \Z$,
\[
\begin{cases}
\ V \not\subset \{(X,Y) \in \C \times \C^\times : m X = a\}\\
\ V \not\subset \{(X,Y) \in \C \times \C^\times : Y^m = b \neq 0\}
\end{cases}
.
\]
Proceeding, to produce a point $(z, e^z) \in V$, what has been established in \S53, \#5 serves to settle things if $P \in \Q[X,Y]$ or if instead 
$P \in \Qbar[X,Y]$ (cf. \S53, \#13) and the general situation can be handled by an appeal to cf. \S53, \#15.
\end{x}
\end{spacing}
\vspace{0.3cm}

\begin{x}{\small\bf REMARK} \ 
There is a reinforcement of EACP to SEACP, where the ``S'' stands for ``strong''.  
This is done by demanding that the outcomes $(\bx,E(\bx)) \in V$ be generic in a suitable sense.
\vspace{0.2cm}

[Note: \ 
The discussion in \#21 is actually strong.]
\end{x}
\vspace{0.3cm}

Agreeing to admit SCHP, recall the notation of \#10.
\vspace{0.3cm}

\begin{x}{\small\bf NOTATION} \ 
Let $\K$ be an \mE-field with Schanuel's property.  
Given a finite set $X \subset \K$, put
\[
\ecl_{\K}(X) \ = \ \{x \in \K : \dim_{\K} (X \cup \{x\}) = \dim_{\K}(X)\}. 
\]
\end{x}
\vspace{0.3cm}

\begin{x}{\small\bf \un{N.B.}} \ 
$\ecl_{\K}(X) $ is called the 
\un{exponential closure}
\index{exponential closure} 
of \mX.
\end{x}
\vspace{0.3cm}

\begin{x}{\small\bf DEFINITION} \ 
(Admit SCHP) \ 
An \mE-field $\K$ has the
\un{countable closure} 
\un{property}
\index{countable closure property} 
 (CCP) \index{CCP} 
if for any finite set $X \subset \K$, $\ecl_\K(X)$ is countable.
\end{x}
\vspace{0.3cm}

There is another approach to exponential closure which forgoes SCHP and has the merit that it can be used to establish that $\C_{\exp}$ has the CCP.
\vspace{0.3cm}

\begin{x}{\small\bf DEFINITION} \ 
An 
\un{exponential polynomial}
\index{exponential polynomial} 
is a function of the form
\[
f(\bx) \ = \ P(\bx, E(\bx)),
\]
where
\[
P \in \K[X_1, \ldots, X_n, Y_1, \ldots, Y_n].
\]
\end{x}
\vspace{0.1cm}

\begin{spacing}{1.75}
\begin{x}{\small\bf \un{N.B.}} \ 
Formal differentiation of polynomials can be extended to exponential polynomials in a unique way such that 
$\ds\frac{\partial e^X}{\partial X} \ = \ e^X$.
\end{x}
\end{spacing}
\vspace{0.3cm}

\begin{x}{\small\bf DEFINITION} \ 
A 
\un{Khovanskii system of width $n$}
\index{Khovanskii system of width $n$} 
consists of exponential polynomials $f_1, \ldots, f_n$ with equations
\[
f_i(x_1, \ldots, x_n) \ = \ 0 \qquad (i = 1, \ldots, n)
\]
and the inequation\\
\[
\abs{
\begin{matrix}
\ \ds\frac{\partial f_1}{\partial x_1}
&\cdots
&\ds\frac{\partial f_1}{\partial x_{n}}  \\[8pt]
\vdots
&&\vdots\\[8pt]
\ds\frac{\partial f_n}{\partial x_1}
&\cdots
&\ds\frac{\partial f_n}{\partial x_n}\\[8pt]
\end{matrix}
}
\quad (x_1, \ldots, x_n)  \ \neq \ 0,
\]\\
the differentiation being the formal differentiation of exponential polynomials.
\end{x}
\vspace{0.3cm}

\begin{spacing}{1.4}
\begin{x}{\small\bf LEMMA} \ 
\ (Admit SCHP) \ 
Let $\K$ be an $E$-field, $X \subset \K$ a finite subset $-$then
$\ecl_\K(X)$ consists of those points $x \in \K$ with the property that there are $n \in \N$, $x_1, \ldots, x_n \in \K$, and exponential polynomials 
$f_1, \ldots, f_n$ with coefficients from $\Q(X)$ such that $x = x_1$ and $(x_1, \ldots, x_n)$ is a solution to the Khovanskii system given by the $f_i$.  
\end{x}
\end{spacing}
\vspace{0.3cm}

Now drop SCHP and for any $E$-field $\K$ take for the definition of $\ecl_\K(X)$ the property figuring in \#29, thereby extending the definition of CCP to all $E$-fields $\K$.
\vspace{0.3cm}

\begin{spacing}{1.4}
\begin{x}{\small\bf THEOREM} \ 
$\C_{\exp}$ has the countable closure property.
\vspace{0.2cm}

PROOF \ 
Given a finite subset $X \subset \C_{\exp}$, there are only countably many Khovanskii systems with coefficients from $\Q(X)$.  
The inequation 
in Khovanskii system amounts to saying that the Jacobian of the functions $f_1, \ldots, f_n$ does not vanish, so by the implicit function theorem, solutions to a Khovanskii system are isolated, hence there are but countably many solutions to each system, thus implying that
\[
\ecl_{\C_{\exp}}(X)
\]
is countable.
\end{x}
\end{spacing}
\vspace{0.3cm}

We come now to the fundamental definition: \ A 
\un{Zilber field}
is a field $\K$ subject to the conditions listed at the beginning.
\vspace{0.2cm}

[Note: \ 
Denote this setup by the symbol $\K_E$.]
\vspace{0.3cm}

\begin{x}{\small\bf THEOREM} \ 
For $\kappa$ uncountable, up to isomorphism there is a unique Zilber field of size $\kappa$.
\end{x}
\vspace{0.3cm}

\begin{x}{\small\bf CONJECTURE} \ 
The Zilber field of size continuum is isomorphic to $\C_{\exp}$.
\end{x}
\vspace{0.3cm}


%% file: _55_e_rings.tex
\chapter{
$\boldsymbol{\S}$\textbf{55}.\quad  $\boldsymbol{E}$-RINGS}
\setlength\parindent{2em}
\setcounter{theoremn}{0}
\renewcommand{\thepage}{\S55-\arabic{page}}


\begin{x}{\small\bf DEFINITION} \ 
An 
\un{\mE-ring}
\index{\mE-ring}
is a pair $(R,E)$, where \mR is a ring (commutative with 1) and 
\[
E: (R,+) \ra (U R, \cdot)
\]
is a map from the additive group of \mR to the multiplicative group of units of \mR such that
\[
\forall \ x, y \in R, \quad E(x + y) \  = \ E(x) \cdot E(y) 
\]
and $E(0) = 1$.

\vspace{0.2cm}

[Note: \ 
Every ring \mR becomes an \mE-ring via the stipulation
\[
E(x) \ = \ 1 \qquad (x \in R).]
\]
\end{x}
\vspace{0.3cm}

\begin{x}{\small\bf EXAMPLE} \ 
Every \mE-field is an \mE-ring (cf. \S54, \#1). 
\vspace{0.2cm}

[Note: \ 
By definition, an \mE-field has characteristic 0, matters being trivial in positive characteristic.  Thus suppose that $\K$ is a field of charcteristic $p > 0$ $-$then $\forall \ x \in \K$,
\allowdisplaybreaks
\begin{align*}
 1\
&=\ 
E(0)
\\[12pt]
&=\ 
E(x + x + \cdots + x) \qquad (p \ \text{terms})
\\[12pt]
&=\ 
E(x)^p 
\end{align*}
\qquad\qquad $\implies$
\allowdisplaybreaks
\begin{align*}
(E(x) - 1)^p
&=\ 
E(x)^p - 1^p
\\[12pt]
 &=\ 
E(x)^p - 1 
\\[12pt]
&=\ 0 
\\[12pt]
\implies 
\\[12pt]
&\ E(x) = 1.]
\end{align*}
\end{x}
\vspace{0.3cm}

\begin{x}{\small\bf EXAMPLE} \ 
Take $R = \Z$ and define \mE by the prescription
\[
E(x) \ = \ 1 \qquad (x \in \Z).
\]
Another possibility is the prescription
\[
E(x) \ = \ 
\begin{cases}
\ \hspace{.25cm} 1 \quad \text{if $x$ is even}\\[7pt]
\ -1 \quad \text{if $x$ is odd}\\
\end{cases}
(x \in \Z).
\]
\vspace{0.2cm}

[Note: \ 
These are the only two possibilities.]
\end{x}
\vspace{0.3cm}

\begin{spacing}{1.55}
\begin{x}{\small\bf RAPPEL} \ 
If \mG is a multiplicative group (finite or infinite) and \mR is a ring (commutative with 1), then the 
\un{group ring $R[G]$ of \mG over \mR}
\index{group ring $R[G]$ of \mG over \mR} 
is the set of all finite linear combinations of elements of \mG with coefficients in \mR, thus
\[
\sum\limits_{g \in G} \hsx r_g g,
\]
where $r_g = 0$ for all but finitely many elements of \mG and the ring operations are defined in the obvious way.
\vspace{0.2cm}

[Note: \ 
If 1 is the identity of \mR and $e$ is the identity of \mG, then $1 e$ is the identity of  $R[G]$.]
\end{x}
\vspace{0.3cm}

Let $X_1, \ldots, X_n$ be distinct indeterminants.
\vspace{0.5cm}

\begin{x}{\small\bf DEFINITION} \ 
The 
\un{free \mE-ring},
\index{free \mE-ring}
denoted
\[
[X_1, \ldots, X_n]^E,
\]
is an \mE-ring containing $X_1, \ldots, X_n$ as elements and having the property that for each \mE-ring \mR and elements 
$r_1, \ldots, r_n \in R$ there is one and only one \mE-ring morphism 
\[
f : [X_1, \ldots, X_n]^E \ra R
\]
such that
\[f(X_i) \ = \ r_i \qquad (i = 1, \ldots, n).
\]
\end{x}
\vspace{0.3cm}

\begin{x}{\small\bf \un{N.B.}} \ 
The free \mE-ring on no generators, denoted $[\emptyset]^E$ (``$n=0$''), is admitted.  
It has the property that for each \mE-ring \mR there is an \mE-morphism from $[\emptyset]^E$ to \mR.
\end{x}
\vspace{0.3cm}

The existence of 
\[
[X_1, \ldots, X_n]^E
\]
is established via an argument of recursion, itself a special case of the following considerations.  
Given an \mE-ring \mR, one can form the free \mE-ring extension of \mR on generators $X_1, \ldots, X_n$, denoted
\[
R[X_1, \ldots, X_n]^E,
\]
its elements being by definition the \mE-ring of exponential polynomials.
\vspace{0.2cm}

[Note: \ 
Take $R = \Z$ $(E \equiv 1)$ to recover 
\[
[X_1, \ldots, X_n]^E.]
\]
\end{spacing}
\vspace{0.3cm}

\begin{x}{\small\bf CONSTRUCTION} \ 
We shall construct three sequences:
\vspace{0.2cm}

\qquad \textbullet \quad $(R_k, +, \cdot)_{\hsx k \hsx \geq\hsx -1}$ are rings;
\vspace{0.2cm}

\qquad \textbullet \quad $(A_k, +)_{\hsx k \hsx\geq\hsx 0}$ are abelian groups;
\vspace{0.2cm}

\qquad \textbullet \quad $(E_k)_{\hsx k \hsx\geq\hsx -1}$ are \mE-morphisms from $R_k$ to $UR_{k+1}$.
\vspace{0.5cm}

\qquad \un{Initial Step:} \quad
Put $R_{\hsx -1} = R$, 
\[
R_0 \ = \ R[X_1, \ldots, X_n],
\]
and let $A_0$ be the ideal generated by $X_1, \ldots, X_n$.  So, as an additive group,
\[
R_0 \ = \ R \oplus A_0 \qquad (= R_{\hsx -1} \oplus A_0).
\]
Define the morphism
\[
E_{\hsx -1}:R_{\hsx -1} \ra R_0
\]
by the composition
\[
R_{\hsx -1} \ = \ 
R \overset{E}{\lra} R \overset{i}{\lra} R[X_1, \ldots, X_n] = R_0.
\]
\vspace{0.02cm}

\qquad \un{Inductive Step:} \quad
Suppose that $k \geq 0$ and $R_{k-1}$, $R_k$, $A_k$, and 
$E_{k-1}$ have been defined in such a way that
\[
R_k \ = \ R_{k-1} \oplus A_k, \ 
E_{k-1} : (R_{k-1},+) \ra (U R_k, \cdot).
\]
Let
\[
t: (A_k,+) \ra \big(t^{A_k},\cdot\big)
\]
be a formal isomorphism (additive $\ra$ multiplicative).  Define
\[
R_{k+1} \ = \ R_k[t^{A_k}].
\]
Therefore $R_k$ is a subring of $R_{k+1}$ and as an additive group
\[
R_{k+1} \ = \ R_k \oplus A_{k+1},
\]
where $A_{k+1}$ is the $R_k$-submodule of $R_{k+1}$ freely generated by the $t^a$ \ ($a \ \in\  A_k$, $a \ \neq \ 0$).  
Next extend
\[
E_k : (R_k,+) \ra (U R_{k+1}, \cdot)
\]
by
\[
E_k(x) 
\ = \ 
E_{k-1} (y) \cdot t^a \qquad (x = y + a, \ \text{with} \ y \in R_{k-1}, \ a \in A_k). 
\]
In this way there is assembled a chain of partial \mE-rings (the domain of exponentiation of $R_{k+1}$ is $R_k$):
\[
R_0 \subset R_1 \subset \cdots .
\]
Definition: 
\[
R [X_1, \ldots, X_n]^E 
\ = \ 
\bigcup\limits_{k=0}^\infty \hsx R_k,
\]
its \mE-ring morphism being the prescription
\[
E(x) \ = \ E_k(x) \qquad (x \in R_k).
\]
\end{x}
\vspace{0.3cm}

\begin{x}{\small\bf \un{N.B.}} \ 
$R_{k+1}$ as an additive group is the direct sum
\[
R \oplus A_0 \oplus A_1 \oplus \cdots \oplus A_{k+1}.
\]

\vspace{0.2cm}

[Note: \ 
The group ring $R_{k+1}$ is isomorphic to
\[
R_0[t^{ A_0 \hsx \oplus\hsx \cdots\hsx \oplus \hsx A_k}] 
\]
or still, is isomorphic to
\[
R_1[t^{ A_1  \hsx \oplus\hsx \cdots\hsx \oplus \hsx A_k}]
\]
\[
\cdots
\]
or still, is isomorphic to
\[
R_k[t^{A_k}] \hsx .]
\]
\end{x}
\vspace{0.3cm}

\begin{x}{\small\bf \un{N.B.}} \ 
\[
R [X_1, \ldots, X_n]^E 
\]
as an additive group is
\[
R \hsx \oplus\hsx A_0 \hsx \oplus\hsx A_1 \hsx\oplus\hsx \cdots \hsx\oplus\hsx A_k \hsx \oplus\hsx \cdots
\]
and as a group ring is
\[
R [X_1, \ldots, X_n] \hsy [t^{ A_0 \hsx \oplus A_1 \hsx \oplus \hsx \cdots \hsx \oplus A_k \hsx \oplus\hsx \cdots}].
\]
\end{x}
\vspace{0.3cm}

\begin{x}{\small\bf EXPONENTIATIONS} \ 
\vspace{0.3cm}

\qquad \textbullet \quad 
Let $P \in R_k \ (k \geq 0)$ $-$then \mP can be written uniquely as 
\[
P \ = \ P_0 + P_1 + \cdots + P_k,
\]
where $P_0 \in R_0$ and $P_\ell \in A_\ell$ $(\ell > 0)$.
\vspace{0.3cm}

\qquad \textbullet \quad 
Let $P \in A_k \ (k \geq 1)$ $-$then \mP can be written uniquely as 
\[
P \ = \ 
\sum\limits_{i=1}^N \hsx 
r_i E(a_i),
\]
where $a_i \in A_{k-1} - \{0\}$ and $a_i \neq a_j$ for $i \neq j$ and $r_1, \ldots r_N$ are nonzero elements of $R_{k-1}$.
\vspace{0.2cm}

[Note: \ 
The isomorphism $t:A_k \ra t^{A_k}$ is the restriction of the exponential map \mE to $A_k$:
\[
E(A_k) \ = \ t^{A_k}.]
\]
\end{x}
\vspace{0.3cm}

\begin{x}{\small\bf EXAMPLE} \ 
Take $n = 2$ and work with
\[
[X_1, X_2]^E \ \equiv \ [X,Y]^E \ \equiv \ \Z[X,Y]^E.
\]
Then $(k=2)$

\allowdisplaybreaks
\begin{align*}
P(X,Y) \ 
&=\ 
-3X^2 Y - X^5 Y^7 
\\[12pt]
&\hspace{.6cm}
+ (2XY + 5Y^2) \hsx E(-7X^3 + 11 X^5 Y^4)
\\[12pt]
&\hspace{.6cm}
+ (6 - 2XY^5) \hsx E((5X + 2X^7 Y^2) \hsx 
E(5X - 10Y^2))
\end{align*}
is an element of $R_2$ (per $\Z$):
\[
P \ = \ P_0 + P_1 + P_2.
\]
\end{x}
\vspace{0.3cm}

\begin{spacing}{1.6}
\begin{x}{\small\bf EXAMPLE} \ 
Consider the free \mE-ring $[\emptyset]^E$ on no generators $-$then the elements of $[\emptyset]^E$ are 
``exponential constants'', e.g., in suggestive notation,
\[
\text{\large{$e^{e^2 + 3} + 4 - 5 e^{3 + e^{-3}}$}}.
\]
\end{x}
\end{spacing}
\vspace{0.3cm}

\begin{x}{\small\bf LEMMA} \ 
Given an \mE-ring \mT and elements $t_1, \ldots, t_n \in T$, every \mE-ring morphism $\phi : R \ra T$ has a unique extension to an \mE-ring morphism
\[
\Phi: R [X_1, \ldots, X_n]^E  \ra T
\]
such that
\[
\Phi(X_i) \ = \ t_i \qquad (i = 1, \ldots, n).
\]

[Use the corresponding property of
\[
R [X_1, \ldots, X_n] \ = \ R_0
\]
and extend stepwise to each $R_k$ ($k > 0$).]
\vspace{0.5cm}

Suppose that $(R,E)$ is an \mE-ring.  
Given a set $I \neq \emptyset$, let $R^I$ be the set of functions $I \ra R$ $-$then 
$R^I$ is an \mE-ring.  
Let $f \in R^I$ and define $E f$ by the rule
\[
(E \hsy f) (i) \ = \ E(f(i)),
\]
i.e., operations are pointwise.
\vspace{0.3cm}

Take $I = R^n$ and consider $R^{R^n}$, the functions from $R^n$ to \mR.  
Define the coordinate functions $x_1, \ldots, x_n \in R^{R^n}$ by
\[
x_i(r_1, \ldots, r_n) \ = \ r_i \qquad (i = 1, \ldots, n).
\]
In \#13, take $T = R^{R^n}$.  
Embed \mR in $R^{R^n}$ by assigning to each $r \in R$ the constant function 
$C_r \hspace{.3cm} (C_r(r_1, \ldots, r_n) = r)$ $-$then the assignment
\[
C : \ 
\begin{cases}
\ R \ra R^{R^n}\\[7pt]
\ r \ra C_r
\end{cases}
\]
is an \mE-ring morphism, hence \mC admits a unique extension to an \mE-ring morphism
\[
R[X_1, \ldots, X_n]^{E} \ra R^{R^n}
\]
that sends each $X_i$ to $x_i$, the 
\un{canonical arrow},
\index{canonical arrow} 
call it $\Gamma$.

\end{x}
\vspace{0.3cm}

\begin{x}{\small\bf NOTATION} \ 
Write
\[
R[X_1, \ldots, X_n]^E
\]
in place of
\[
R^{R^n},
\]
its elements being by definition the \mE-ring of exponential polynomial functions.
\end{x}
\vspace{0.3cm}

\begin{spacing}{1.5}
\begin{x}{\small\bf LEMMA} \ 
If $(R,E)$ is an \mE-ring and if \mR is an integral domain of characteristic 0, then $R[X_1, \ldots, X_n]^{E}$ is an integral domain (and its units are of the form $u E (P)$, where $u$ is a unit of \mR and $P \in R[X_1, \ldots, X_n]^E$).
\vspace{0.2cm}

[Without going into detail, let us recall only that if \mR is an integral domain of characteristic 0 and \mG is  a multiplicative group, then the group ring $R[G]$ is an integral domain of characteristic 0 iff \mG is torsion free.]
\end{x}
\end{spacing}
\vspace{0.3cm}

\begin{x}{\small\bf \un{N.B.}} \ 
By induction on $k \geq 0$, assume that $R_k$ is an integral domain of characteristic 0 $-$then $A_k$ is torsion free.  
Therefore $t^{A_k}$ is torsion free, which implies that
\[
R_{k+1} \ = \ R_k[t^{A_k}]
\]
is an integral domain of characteristic 0.
\end{x}
\vspace{0.3cm}

In general, the canonical arrow
\[
\Gamma : R[X_1, \ldots, X_n]^E \ra R[x_1, \ldots, x_n]^E
\]
may have a nontrivial kernel.
\vspace{0.3cm}

\begin{x}{\small\bf EXAMPLE} \ 
Consider a ring \mR equipped with the trivial exponentiation, i.e., $E(x) = 1$ for all $x \in R$ $-$then $E(X_1) - 1$ is in the kernel of $\Gamma$.
\vspace{0.3cm}

[In fact, 
\allowdisplaybreaks
\begin{align*}
\Gamma(E(X_1) - 1) \ 
&=\ 
\Gamma E(X_1) - \Gamma 1
\\[12pt]
&=\ 
E(\Gamma X_1) - C_1
\\[12pt]
&=\ 
E(x_1) - C_1.
\end{align*}
And
\allowdisplaybreaks
\begin{align*}
E(x_1) (r_1, \ldots, r_n)\ 
&=\ 
E(x_1(r_1, \ldots, r_n))
\\[12pt]
&=\ 
E(r_1)
\\[12pt]
&=\ 
1
\\[12pt]
&=\ 
C_1(r_1, \ldots, r_n)
\end{align*}
\qquad\qquad $\implies$
\[
E(x_1) \ = \ C_1.
\]
Therefore
\allowdisplaybreaks
\begin{align*}
\Gamma(E(X_1) - 1) \ 
&=\ 
E(x_1) - C_1
\\[12pt]
&=\ 
C_1 - C_1
\\[12pt]
&=\ 
0.]
\end{align*}
\end{x}
\vspace{0.3cm}

\begin{x}{\small\bf THEOREM} \ 
Suppose that $(R,E)$ is an \mE-ring and \mR is an integral domain of characteristic 0.  
Make the following assumptions.
\vspace{0.2cm}

\qquad \textbullet \quad
There are derivations $\td_1, \ldots, \td_n$ of $R[x_1, \ldots, x_n]^{E}$ which are trivial on \mR and satisfy the condition 
$\td_i(x_j) = \delta_{i j}$ $(1 \leq i, j \leq n)$.
\vspace{0.2cm}

\qquad \textbullet \quad
There is a nonzero element $r \in R$ such that
\[
\td_i(E(f)) \ = \ r \td_i (f) E(f)
\]
for all $f$ in $R[x_1, \ldots, x_n]^{E}$ $(i = 1, \ldots, n)$.
\vspace{0.2cm}

Then $\Gamma$ is one-to-one.
\vspace{0.5cm}

Specialize now the theory outlined above and take $R = \C$, shifting matters to 
\[
\C[X_1, \ldots, X_n]^{\exp} \qquad (E = \exp), 
\]
which, as will be recalled, is a group ring (cf. \#9).  
Moreover, since $\C$ is an integral domain of charcteristic 0, it follows from \#15 that
\[
\C[X_1, \ldots, X_n]^{\exp} 
\]
is an integral domain.
\vspace{0.2cm}

[Note: \ 
While $\C[X_1, \ldots, X_n]$ is noetherian, this is definitely not the case of 
\[
\C[X_1, \ldots, X_n]^{\exp} .]
\]
\end{x}
\vspace{0.3cm}

\begin{x}{\small\bf THEOREM} \ 
The canonical arrow
\[
\Gamma : \C[X_1, \ldots, X_n]^{\exp} \ra \C[x_1, \ldots, x_n]^{\exp}
\]
is one-to-one.
\vspace{0.2cm}

[Apply \#18 (take $\td_1, \ldots, \td_n$ as the partial derivatives $\partial/\partial x_1, \ldots, \partial/\partial x_n$ and 
choose $r = 1$).]
\end{x}
\vspace{0.3cm}

\begin{x}{\small\bf NOTATION} \ 
Put
\[
\EXP(\C^n) 
\ = \ 
\Gamma \C[X_1, \ldots, X_n]^{\exp}.
\]
\end{x}
\vspace{0.3cm}

\begin{x}{\small\bf LEMMA} \ 
(cf. \#10) \ Each function $f$ in $\EXP(\C^n)$ can be written as a finite sum
\[
f 
\ = \ 
\sum\limits_i \hsx
P_i \cdot \exp(g_i),
\]
where
\[
P_i \in \C[X_1, \ldots, X_n]
\quad \text{and} \quad
g_i \in \EXP(\C^n).
\]
\end{x}
\vspace{0.3cm}

\begin{x}{\small\bf EXAMPLE} \ 
Take $n = 1$ and let $X_1 = X$ $-$then the function $z \ra e^z$ belongs to $\EXP(\C)$.
\vspace{0.2cm}

[For
\[
X \in A_0 \implies \EX \in A_1. 
\]
And $\Gamma X = x$, where $x:\C \ra \C$ is the function $z \ra z$ (i.e., $x(z) = z$), hence
\[
\Gamma \hsx \EX 
\ = \ 
\exp \Gamma X 
\ = \ 
\exp x,
\]
the function $\C \ra \C$ that sends $z$ to $\exp x(z) = \exp z$.
\end{x}
\vspace{0.3cm}

\begin{x}{\small\bf EXAMPLE} \ 
The function
\[
(z_1, z_2) \ra z_1 z_2 \cdot \exp(\exp(z_1 + z_2))
\]
belongs to $\EXP(\C^2)$.
\end{x}
\vspace{0.3cm}


%% file: _56_shanuel_implies_shapiro.tex
\chapter{
$\boldsymbol{\S}$\textbf{56}.\quad  SCHANUEL $\boldsymbol{\implies}$ SHAPIRO}
\setlength\parindent{2em}
\setcounter{theoremn}{0}
\renewcommand{\thepage}{\S56-\arabic{page}}


\begin{x}{\small\bf DEFINITION} \ 
Working over $\C$, an 
\un{exponential polynomial}
\index{exponential polynomial} 
is an entire function $f$ of the form
\[
f(z) 
\ = \ 
\lambda_1 e^{\textstyle\mu_1 z} + \cdots + \lambda_n e^{\textstyle\mu_n z},
\]
where $\lambda_1, \ldots, \lambda_n$ and $\mu_1, \ldots, \mu_n$ are complex numbers.
\end{x}
\vspace{0.3cm}

Under addition and multiplication, the set of all such functions form a commutative ring $\sE$ with 1.
\vspace{0.2cm}

[Note: \ 
The units are the elements of the form $\lambda e^{\mu z}$ ($\lambda \neq 0)$.]
\vspace{0.3cm}

\begin{x}{\small\bf REMARK}  \ 
This is the simplest situation since one could, e.g., allow $\lambda_1, \ldots, \lambda_n$ to be complex polynomials.
\end{x}
\vspace{0.3cm}

\begin{x}{\small\bf SHAPIRO'S CONJECTURE} \ 
If $f$, $g$ are two exponential polynomials with infinitely many zeros in common, then there exists an exponential polynomial $h$ such that $h$ is a common divisor of $f$, $g$ in the ring $\sE$ and $h$ has infinitely many zeros in $\C$.
\end{x}
\vspace{0.3cm}

As will be seen below, the proof of Shapiro's conjecture breaks up into two cases (terminology per infra).
\vspace{0.2cm}

\quad \un{Case 1:} \quad
Either $f$ or $g$ is simple.
\vspace{0.2cm}

\quad \un{Case 2:} \quad
Both $f$ and $g$ are irreducible.
\vspace{0.3cm}

\begin{x}{\small\bf \un{N.B.}} \ 
It turns out that the proof of Case 1 does not require Schanuel but the proof of Case 2 does require Schanuel, hence the rubric 
\[
\text{Shanuel $\implies$ Shapiro.}
\]
\end{x}
\vspace{0.3cm}


To prepare for the case distinction, we shall need some definitions and a few classical facts.
\vspace{0.3cm}

\begin{x}{\small\bf DEFINITION} \ 
Let
\[
f(z) 
\ = \ 
\lambda_1 e^{\textstyle\mu_1 z} + \cdots + \lambda_n e^{\textstyle\mu_n z}
\]
be an exponential polynomial $-$then its 
\un{support}, 
\index{support (exponential polynomial)}
denoted $\spt(f)$, is the vector space over $\Q$ generated by $\mu_1, \ldots, \mu_n$.
\end{x}
\vspace{0.3cm}

\begin{x}{\small\bf DEFINITION} \ 
An exponential polynomial $f$ is said to be 
\un{simple}
\index{silmple (exponential polynomial)} 
if
\[
\dim_{\Q} \spt (f) \ = \ 1.
\]
\end{x}
\vspace{0.3cm}

\begin{x}{\small\bf EXAMPLE}\ 
\[
f(z) \ = \ \sin z \ = \ \frac{e^{\sqrt{-1}\hsx \textstyle z}  - e^{-\sqrt{-1}\hsx \textstyle z}}{2 \sqrt{-1}}
\]
is simple.
\end{x}
\vspace{0.3cm}

\begin{x}{\small\bf DEFINITION} \ 
An exponential polynomial $f$ is said to be 
\un{irreducible}
\index{irreducible (exponential polynomial)}
if it is not a unit and has no divisors in the ring $\sE$ other than associates.
\end{x}
\vspace{0.3cm}

Here is Ritt's factorization theorem.
\vspace{0.3cm}

\begin{x}{\small\bf THEOREM} \ 
Every exponential polynomial $f$ can be written uniquely up to order and multiplication by a unit as a product in $\sE$ of the form
\[
S_1 \cdots S_c \ I_1 \cdots I_d,
\]
where all the $S_j$ are simple with
\[
\spt(S_j) \ \cap \ \spt(S_{j^\prime}) \ = \ \{0\}
\]
for $j \neq j^\prime$ and all the $I_k$ are irreducible.
\end{x}
\vspace{0.3cm}

\begin{spacing}{1.5}
Suppose that $f$, $g$ are two exponential polynomials with infinitely many zeros in common and neither one is simple.  
Write per Ritt:
\[
\begin{cases}
\ f = S_1 \cdots S_c \ I_1 \cdots I_d\\[8pt]
\ g = T_1 \cdots T_u \ J_1 \cdots J_v
\end{cases}
.
\]
Then a common zero of $f$, $g$ must be a zero of a factor of each function, thus two factors 
$\widetilde{f}$, 
$\wtilde[2.5pt]{g}$ 
of $f$, $g$ respectively have infinitely many zeros in common, thus if 
$\widetilde{f}$, 
$\wtilde[2.5pt]{g}$ 
have a common divisor $h$ in $\sE$ with infinitely many zeros, then $h$ is the common divisor of $f$, $g$ postulated in Shapiro's conjecture.
\end{spacing}
\vspace{0.1cm}

Matters have accordingly been reduced to Case 1 and Case 2 formulated at the beginning.
\vspace{0.5cm}

\[
\textbf{APPENDIX}
\]
\vspace{0.5cm}

Let \mR be a commutative ring with 1.
\vspace{0.5cm}

\qquad{\small\bf DEFINITION 1} \ 
Let $x$, $y \in R$ $-$then $y$ 
\un{divides} 
\index{divides (commutative ring)} 
$x$
(or $y$ is a 
\un{divisor}
\index{divisor (commutative ring)} 
of $x$) and $x$ is 
\un{divisible} 
\index{divisible (commutative ring)} 
by $y$
(or $x$ is a 
\un{multiple}
\index{multiple (commutative ring)}
of $y$) 
if there exists $z \in R$ such that $x = y z$.
\vspace{0.2cm}

[Note: \ 
The only elements of \mR which are divisors of 1 are the units of \mR, i.e., the elements of $U R$.]
\vspace{0.5cm}

\qquad{\small\bf DEFINITION 2} \ 
If $x$, $y \in R$ and if $x = y u$, where $u \in U R$, then $x$ and $y$ are said to be 
\un{associates}.
\index{associates}
\vspace{0.2cm}

[Note: \ 
Therefore $y$ divides $x$.  But also $y = x u^{-1}$, thus $x$ divides $y$.]
\vspace{0.5cm}

\qquad{\small\bf DEFINITION 3} \ 
The associates of an element $x \in R$ are the 
\un{improper divisors}
\index{improper divisors (commutative ring)}
of $x$.
\vspace{0.2cm}

[Note: \ 
A unit $u \in U R$ divides every element $x$ of \mR: \ $x = u(u^{-1}x)$.  
Still, the convention is not to include $U R$ in the set of divisors of $x$.]
\vspace{0.5cm}

\qquad{\small\bf DEFINITION 4} \ 
An element $x \in R$ is 
\un{irreducible}
\index{irreducible} 
if it is not a unit and its only divisors are associates, i.e., are improper.
\vspace{0.5cm}

\qquad{\small\bf DEFINITION 5} \ 
Irreducible elements $x$, $y \in R$ are 
\un{distinct}
\index{distinct (irreducible elements of a commutative ring)}
if they are not unit multiples of one another.
\vspace{0.5cm}

\qquad{\small\bf LEMMA} \ 
Distinct irreducibles $x$, $y \in R$ do not have a common divisor.
\vspace{0.2cm}

PROOF \ 
Suppose that $a$ is a common divisor: 
\[
\begin{cases}
\ x = au\\[8pt]
y = av
\end{cases}
\qquad (u, v \in U R).
\]
Then
\[
a \ = \ y v^{-1} 
\implies 
x \ = \ y v^{-1} u,
\]
i.e., $x$ is a unit multiple of $y$.  Contradiction.


%% file: _57_shapiro_conjecture_case_1.tex
\chapter{
$\boldsymbol{\S}$\textbf{57}.\quad  SHAPIRO'S CONJECTURE: \ CASE 1}
\setlength\parindent{2em}
\setcounter{theoremn}{0}
\renewcommand{\thepage}{\S57-\arabic{page}}

\ \indent 
Recall the setup: $f$, $g$ are two exponential polynomials with infinitely many zeros in common and either $f$ or $g$ is simple (cf. \#3).
\vspace{0.3cm}

\begin{spacing}{1.5}
\begin{x}{\small\bf THEOREM} \ 
(Skolem-Mahler-Lech) \ 
Let $f \in \sE$ and let $A \subset \Z$ be the set of integers on which $f$ vanishes $-$then \mA is the finite union of arithmetic progressions, i.e., sets of the form $\{m + k d: k \in \Z\}$ for some $m, \hsx d \in \Z$.  
Moreover, if \mA is infinite, then at least one of these arithmetic progressions has a nonzero difference $d$.
\end{x}
\end{spacing}
\vspace{0.2cm}

This is a wellknown result on the distribution of zeros of exponential polynomials and will be taken without proof.
\vspace{0.3cm}

\begin{x}{\small\bf LEMMA} \ 
Let $f \in \sE$.  
Suppose that $f(k) = 0$ $\forall \ k \in \Z$ $-$then $\sin (\pi z)$ divides $f$ in the ring $\sE$.
\vspace{0.2cm}

\begin{spacing}{1.5}
PROOF \ 
Let
\[
f(z) 
\ = \ 
\lambda_1 e^{\textstyle\mu_1 z} + \cdots + \lambda_n e^{\textstyle\mu_n z},
\]
with $\lambda_1, \ldots, \lambda_n \neq 0$.  
It can be assumed that $f$ is not identically zero and that $n \geq 2$ 
(since $\ds\lambda_1 e^{\textstyle\mu_1 z} = 0$ only if 
$\lambda_1 = 0$).  
Proceed by induction on the length $n$ of $f$.
\end{spacing}
\vspace{0.5cm}

\qquad \textbullet \quad \un{$n = 2$:} 
\[
f(z) 
\ = \ 
\lambda_1 e^{\textstyle\mu_1 z}  + \lambda_2 e^{\textstyle\mu_2 z}
\]
with $\lambda_1$, $\lambda_2 \neq 0$.  
Put $z = 0$ to get
\[
\lambda_1 + \lambda_2 \ = \ 0 
\implies
f(z) = 
\lambda_1(e^{\textstyle\mu_1 z} - e^{\textstyle\mu_2 z}).
\]
Put $z = 1$ to get
\[
e^{\textstyle\mu_1} - e^{\textstyle\mu_2} \ = \ 0
\]
\qquad\qquad $\implies$
\[
\mu_2 = \mu_1 + 2 k \pi \sqrt{-1} \qquad (\exists \ k \in \Z - \{0\})
\]
\qquad\qquad $\implies$
\[
f(z)
\ = \ 
\lambda_1 e^{\textstyle\mu_1 z} \hsx \big(1 - e^{2 k \pi \sqrt{-1} \hsx z} \big).
\]
Without loss of generality, take \ $k\  > \ 0$ (otherwise switch the roles of \ $\mu_1$ \ and \  $\mu_2$).  \\[5pt]
\noindent Next
\[
\sin z 
\ = \ 
\frac{e^{\sqrt{-1} \hsx z} - e^{- \sqrt{-1} \hsx z}}{2 \sqrt{-1}}
\]
\qquad\qquad $\implies$
\[
-2 \sqrt{-1} \hsx e^{\pi \sqrt{-1} \hsx z} \sin (\pi z) 
\ = \ 
1 - e^{2 \pi \sqrt{-1} \hsx z}
\]
\qquad\qquad $\implies$
\allowdisplaybreaks
\begin{align*}
\big(
1 
+ &e^{2 \pi \sqrt{-1} \hsx z}
+ e^{4 \pi \sqrt{-1} \hsx z}
+ \cdots + 
e^{2(k-1) \pi \sqrt{-1} \hsx z}
\big)
\big(
-2 \sqrt{-1} \hsx e^{\pi \sqrt{-1} \hsx z} \sin (\pi z) 
\big)
\\[12pt]
&=\ 
\big(
1 
+ e^{2 \pi \sqrt{-1} \hsx z}
+ e^{4 \pi \sqrt{-1} \hsx z}
+ \cdots + 
e^{2(k-1) \pi \sqrt{-1} \hsx z}
\big)
\big(
1 - e^{2 \pi \sqrt{-1} \hsx z}
\big)
\\[12pt]
&=\ 
1 
+ e^{2 \pi \sqrt{-1} \hsx z}
+ e^{4 \pi \sqrt{-1} \hsx z}
+ \cdots + 
e^{2(k-1) \pi \sqrt{-1} \hsx z}
\\[12pt]
&\hspace{1cm}
- e^{2 \pi \sqrt{-1} \hsx z}
- e^{4 \pi \sqrt{-1} \hsx z}
- \cdots - 
e^{2(k-1) \pi \sqrt{-1} \hsx z}
\ - \ e^{2 k \pi \sqrt{-1} \hsx z}
\\[12pt]
&=\ 
1 - e^{2 k \pi \sqrt{-1} \hsx z}
\end{align*}
\qquad $\implies$
\allowdisplaybreaks
\begin{align*}
f(z)\ 
&=\ 
\lambda_1 e^{\textstyle\mu_1 z} \hsx \big(1 - e^{2 k \pi \sqrt{-1} \hsx z} \big)
\\[12pt]
&=\ 
\lambda_1 e^{\textstyle\mu_1 z}F(z) (-2 \sqrt{-1} e^{\pi \sqrt{-1} \hsx z} \sin(\pi z))
\end{align*}

if 
\[
F(z) 
\ = \ 
1 + 
e^{2 \pi \sqrt{-1} \hsx z} 
+ 
e^{4 \pi \sqrt{-1} \hsx z} 
+ \cdots + 
e^{2 (k-1) \pi \sqrt{-1} \hsx z}.
\]
Therefore $\sin(\pi z)$ divides $f(z)$.
\vspace{0.3cm}

\qquad \textbullet \quad \un{$n > 2$:} \quad 
Suppose now that for all exponential polynomials $h(z)$ of length $\leq n - 1$ which vanish at the integers, $\sin(\pi z)$ divides $h(z)$.  
Setting $z = 1, 2, \ldots, n$ in $f(z)$ leads to the relations
\allowdisplaybreaks
\begin{align*}
&
\lambda_1 e^{\textstyle\mu_1} + \cdots + \lambda_n e^{\textstyle\mu_n} \ = \ 0
\\[12pt]
&
\lambda_1 \big(e^{\textstyle\mu_1}\big)^2 + \cdots + \lambda_n  \big(e^{\textstyle\mu_n}\big)^2 \ = \ 0
\\[12pt]
&
\hspace{2cm} \vdots
\\[12pt]
&
\lambda_1 \big(e^{\textstyle\mu_1}\big)^n + \cdots + \lambda_n \big(e^{\textstyle\mu_n}\big)^n \ = \ 0.
\end{align*}
Let $\delta_j = e^{\textstyle\mu_j}$ $(j = 1, \ldots, n)$, hence in matrix notation
\[
\begin{pmatrix}
\delta_1 &\delta_2 &\cdots &\delta_n\\[8pt]
\delta_1^2 &\delta_2^2 &\cdots &\delta_n^2\\[8pt]
&&\vdots\\[8pt]
\delta_1^n &\delta_2^n &\cdots &\delta_n^n\\[8pt]
\end{pmatrix}
\ 
\begin{pmatrix}
\lambda_1\\[8pt]
\lambda_2\\[8pt]
\vdots\\[8pt]
\lambda_n\\[8pt]
\end{pmatrix}
\ = \ 
\begin{pmatrix}
\text{ 0 }\\[8pt]
0\\[8pt]
\vdots\\[8pt]
0\\[8pt]
\end{pmatrix}
.
\]
Since $\lambda_1, \ldots, \lambda_n \neq 0$, they constitute a nontrivial solution of the corresponding system of linear equations, thus the determinant of the matrix vanishes:\\
\[
\abs{\ 
\begin{matrix}
\delta_1 &\delta_2 &\cdots &\delta_n\\[8pt]
\delta_1^2 &\delta_2^2 &\cdots &\delta_n^2\\[8pt]
&&\vdots\\[8pt]
\delta_1^n &\delta_2^n &\cdots &\delta_n^n\\[8pt]
\end{matrix}
\ 
}
\quad = \quad 0
\]
or still, 
\[
\delta_1 \hsx \delta_2 \cdots \delta_n \ 
\abs{\ 
\begin{matrix}
1 &1 &\cdots &1\\[8pt]
\delta_1 &\delta_2 &\cdots &\delta_n\\[8pt]
&&\vdots\\[8pt]
\delta_1^{n-1} &\delta_2^{n-1} &\cdots &\delta_n^{n-1}\\[8pt]
\end{matrix}
\ 
}
\quad = \quad 0.
\]
This is a Vandermonde determinant, so we are led to 
\[
\delta_1 \cdots \delta_n \hsx
\prod\limits_{1 \hsy \leq \hsy i  \hsy < \hsy j  \hsy \leq \hsy n} \hsx 
(\delta_i - \delta_j) 
\ = \
0.
\]
Since all the $\delta_i$ are nonzero, it must be the case that $\delta_i = \delta_j$ for some $i < j$.  
Without loss of generality, assume $\delta_1 = \delta_2$, i.e., $e^{\textstyle\mu_1} = e^{\textstyle\mu_2}$.  
Put
\[
h(z) 
\ = \ 
(\lambda_1 + \lambda_2) e^{\textstyle\mu_1 z} 
+ 
\sum\limits_{j=3}^n \hsx
\lambda_j e^{\textstyle\mu_j z}.
\]
Then
\allowdisplaybreaks
\begin{align*}
h(z) - \lambda_2\big(e^{\textstyle \mu_1 z} - e^{\textstyle\mu_2 z}\big)\ 
&=\ 
\lambda_1 e^{\textstyle\mu_1 z}
\hsx +\hsx
\lambda_2 e^{\textstyle\mu_1 z}
\hsx - \hsx
\lambda_2 e^{\textstyle\mu_1 z}
\hsx + \hsx
\lambda_2 e^{\textstyle\mu_2 z}
\hsx + \ 
\sum\limits_{j=3}^n \hsx
\lambda_j e^{\textstyle\mu_j z}
\\[12pt]
&=\ 
\lambda_1 e^{\textstyle\mu_1 z}
\hsx + \hsx
\lambda_2 e^{\textstyle\mu_2 z}
\hsx + \ 
\sum\limits_{j=3}^n \hsx
\lambda_j e^{\textstyle\mu_j z}
\\[12pt]
&=\ 
f(z).
\end{align*}
And $\forall \ k \in \Z$, 
\allowdisplaybreaks
\begin{align*}
h(k) \ 
&=\ 
f(k) + \lambda_2\big(e^{\hsx\textstyle\mu_1 k} - e^{\hsx\textstyle\mu_2 k}\big)
\\[12pt]
&=\ 
0.
\end{align*}
\begin{spacing}{1.5}
Consequently $h(z)$ vanishes at the integers.  
But its length is $< n$, hence by the induction hypothesis, $\sin(\pi z)$ divides $h(z)$.  
On the other hand, arguing as in the case $n=2$, $\sin (\pi z)$ divides 
$\lambda_2 \big(e^{\textstyle\mu_1 z} - e^{\textstyle\mu_2 z}\big)$.  
So finally $\sin(\pi z)$ divides $f(z)$. 
\end{spacing}
\vspace{0.2cm}

[Note: \ 
\[
e^{\mu z} 
\ = \ 
\sum\limits_{n=0}^\infty \hsx
\frac{\big(\mu z\big)^n}{n!}
\]
is, in general, not the same as
\[
\big(e^\mu\big)^z 
\ = \ 
e^{z \Log e^\mu} 
\ = \ 
e^{z (\mu + 2 \pi \sqrt{-1} \hsx m)}.
\]
But they are the same if $z = k \in \Z$:
\[
\big(e^\mu\big)^k 
\ = \ 
e^{k (\mu + 2 \pi \sqrt{-1} \hsx m)} 
\ = \ 
e^{k \mu} 
\ = \ 
e^{\mu k}.]
\]
\end{x}
\vspace{0.3cm}

\begin{spacing}{1.5}
\begin{x}{\small\bf THEOREM} \ 
If $f$, $g$ are two exponential polynomials with infinitely many zeros in common such that at least one of $f$, $g$ are simple, then there exists an exponential polynomial $h$ such that $h$ is a common divisor of $f$, $g$ in the ring $\sE$ and $h$ has infinitely many zeros in $\C$.
\vspace{0.2cm}

PROOF \ 
Take $f$ simple and write
\[
f(z) 
\ = \ 
u(z) \hsx
\prod\limits_{\ell = 1}^L \hsx
\big(1 - \alpha_\ell e^{\rho z}\big),
\]
where $\alpha_1, \ldots, \alpha_L, \rho$ are nonzero complex numbers and $u(z) \in \sE$ is a unit (the simplicity of $f$ implies that there is a nonzero $\kappa \in \C$ and $s_1, \ldots, s_n \in \Z$ such that 
$\mu_1 = s_1 \kappa, \ldots, \mu_n = s_n \kappa$).  
Since this is a finite product, $g$ must have infinitely many zeros in common with one of the factors, say 
$1 - \alpha_1 e^{\rho z}$.  
So suppose that
\[
1 - \alpha_1 e^{\rho z}
\ = \ 
0.
\]
Then
\[
\frac{1}{\alpha_1}
\ = \ 
e^{\rho z}
\]
\qquad\qquad $\implies$
\[
e^{\Log \frac{1}{\alpha_1}} \ = \ 
e^{\rho z}
\]
\qquad\qquad $\implies$
\[
\rho z - \Log \frac{1}{\alpha_1}
\ = \ 
2 k \pi \sqrt{-1} \qquad (\exists \ k \in \Z)
\]
\qquad\qquad $\implies$
\[
z 
\ = \ 
\frac{\raisebox{.15cm}{$\hsx\Log \fraOneAlphaone + 2 k \pi \sqrt{-1}$}}{\rho}.
\]
Therefore the exponential polynomial
\[
G(z) 
\ = \ 
g
\Bigg(\frac{\raisebox{.15cm}{$\hsx\Log \fraOneAlphaone + 2 z \pi \sqrt{-1}$}}{\rho}\Bigg)
\]
vanishes at infinitely many integers.  
Now apply \#1 $-$then for some $m_0$, $d_0 \in \Z$ $(d_0 \neq 0)$, \mG vanishes on $\{m_0 + kd_0 : k \in \Z\}$, 
thus $G(m_0 + z d_0)$ is an exponential polynomial which vanishes at all the integers,  
so $\sin(\pi z)$ divides $G(m_0 + z d_0)$ (cf. \#2).  
Moving on, any integer is a zero of the exponential polynomial
\[
F(z) 
\ = \ 
f
\Bigg(\frac{\raisebox{.15cm}{$\hsx\Log \fraOneAlphaone + 2 z \pi \sqrt{-1}$}}{\rho}\Bigg).
\]
Therefore $F(m_0 + z d_0)$ is an exponential polynomial which vanishes at all the integers, 
so $\sin(\pi z)$ divides $F(m_0 + z d_0)$ (cf. \#2).  
To conclude, consider
\[
h(z) 
\ = \ 
\sin\bigg( \frac{\pi}{d_0}
\bigg(\frac{\raisebox{.15cm}{$\rho z - \Log \fraOneAlphaone$}}{2 \pi \sqrt{-1}} - m_0\bigg)\bigg).
\]
To analyze \mG (ditto for \mF), start from 
\[
G(m_0 + z d_0) 
\ = \ 
\sin(\pi z) G_0(z).
\]
Then
\allowdisplaybreaks
\begin{align*}
G\bigg(m_0 + \frac{1}{d_0}& 
\bigg(\frac{\raisebox{.15cm}{$\rho z - \Log\fraOneAlphaone$}}{2 \pi \sqrt{-1}} - m_0\bigg)d_0\bigg) \ 
\\[15pt]
&=\
G\Bigg(m_0 + \frac{\raisebox{.15cm}{$\rho z - \Log\fraOneAlphaone$}}{2 \pi \sqrt{-1}} - m_0\Bigg) 
\\[15pt]
&=\ 
G\Bigg(\frac{\raisebox{.15cm}{$\rho z - \Log\fraOneAlphaone$}}{2 \pi \sqrt{-1}}\Bigg) 
\\[15pt]
&=\ 
g
\begin{pmatrix}
\frac{\raisebox{.001cm}{$\Log\fraOneAlphaone 
+ 
2$} \hsx
\bigg(
\frac
{\raisebox{-.2cm}{\large$\rho z - \Log\fraOneAlphaone$}}
{\raisebox{-.2cm}{\large$2 \pi \sqrt{-1}$}}
\bigg)
\raisebox{.001cm}{\large$\pi \sqrt{-1}$}
}
{\raisebox{-.2cm}{\large$\rho$}}
\end{pmatrix}
\\[15pt]
&=\ 
g
\begin{pmatrix}
\frac{\raisebox{-.2cm}{$\Log\fraOneAlphaone + \rho z - \LogfraOneAlphaone$}}{\raisebox{-.2cm}{\large$\rho$}}
\end{pmatrix}
\\[15pt]
&=\ 
g(z)
\\[15pt]
&=\ 
\sin\bigg( \frac{\pi}{d_0}
\bigg(\frac{\raisebox{.15cm}{$\rho z - \Log \fraOneAlphaone$}}{2 \pi \sqrt{-1}} - m_0\bigg)\bigg)
\hsx G_0 (\ldots)
\\[12pt]
&=\ 
h(z) \hsx G_0 (\ldots) \ .
\end{align*}
\end{x}
\end{spacing}


%% file: _58_shapiro_conjecture_case_2.tex
\chapter{
$\boldsymbol{\S}$\textbf{58}.\quad  SHAPIRO'S CONJECTURE: \ CASE 2}
\setlength\parindent{2em}
\setcounter{theoremn}{0}
\renewcommand{\thepage}{\S58-\arabic{page}}

\begin{spacing}{1.5}
\ \indent 
In this situation, both $f$, $g$ are irreducible.  
If $f = g u$ for some unit $u \in \sE$, (technically, $f$, $g$ are associates), then $g$ can serve as the ``$h$'' in \S56, \#3.  
On the other hand, if $f$, $g$ are distinct irreducibles (meaning that they are not unit multiples of one another), 
then they cannot have a common divisor (see the Lemma in the Appendix to \S56). 
Matters thus reduce to the following statement.
\end{spacing}
\vspace{0.3cm}

\begin{x}{\small\bf THEOREM} \ 
(Admit SCHC) \ 
Let $f$, $g$ be distinct irreducible exponential polynomials in $\sE$ $-$then $f$, $g$ have at most a finite number of zeros in common.
\end{x}
\vspace{0.3cm}

The proof is difficult and lengthy, thus an outline of the argument will have to do.
\vspace{0.3cm}

\begin{x}{\small\bf REMARK} \ 
Let $f$, $g$ be exponential polynomials and assume that $f$ is irreducible.  
Suppose further that $f$, $g$ have infinitely many zeros in common $-$then $f$ divides $g$ in the ring $\sE$ (i.e., $g / f$ is entire).

\vspace{0.2cm}

[Note: \ 
This assertion is equivalent to \#1.]
\end{x}
\vspace{0.3cm}

Proceeding to \#1, assume that $f$, $g$ are distinct irreducibles with infinitely many zeros in common, the objective being to show that this forces a contradiction (namely that $g$ divides $f$).
\vspace{0.2cm}

[Note: \ 
If $g$ divides $f$, then $g$ must be an associate of $f$, say $f = g u$ $(u \in U \sE)$, thereby forcing $f$ to be a unit multiple of $g$, contradicting the supposition of ``distinct''.]
\vspace{0.3cm}

\begin{x}{\small\bf NOTATION} \ 
Let \mS be the infinite set of nonzero common zeros of $f$, $g$.
\end{x}
\vspace{0.3cm}

\begin{x}{\small\bf MAIN LEMMA} \ 
(Admit SCHC) \ 
There exists an infinite subset $S^\prime$ of \mS such that the $\Q$-vector space spanned by $S^\prime$ is finite dimensional.
\end{x}
\vspace{0.3cm}

Without changing the notation, assume henceforth that \mS spans a finite dimensional vector space over $\Q$.
\\[0pt]

Write
\[
f(z) 
\ = \ 
\lambda_1 e^{\mu_1 z} + \cdots + \lambda_n e^{\mu_n z}
\]
and let $\Gamma$ be the 
divisible hull
\index{divisible hull of a multiplicative group} 
of the multiplicative group generated by
\[
\big\{e^{\mu_j s} : 1 \leq j \leq n, s \in S\big\},
\]
that is, $\forall \ \gamma \in \Gamma$ and any nonzero integer $\ell$, $\exists \ \zeta \in \Gamma$ such that $\zeta^{\ell} = \gamma$ and $\Gamma$ is the smallest such group containing
\[
\big\{e^{\mu_j s} : 1 \leq j \leq n, s \in S\big\}.
\]
Since $\spanxs{\Q} S$ is finite dimensional, $\Gamma$ has finite rank.
\vspace{0.3cm}

\begin{x}{\small\bf DEFINITION} \ 
A solution $\alpha_1, \ldots, \alpha_N$ of the linear equation
\[
a_1 x_1 + \cdots + a_N x_N
\ = \ 
1
\]
over $\C$ is 
\un{nondegenerate}
\index{nondegenerate} 
if for every proper nonempty subset \mJ of $\{1, \ldots, N\}$, 
\[
\sum\limits_{j \in J} \hsx a_j \alpha_j \ \neq \ 0.
\]
\end{x}
\vspace{0.3cm}

\begin{x}{\small\bf THEOREM} \ 
(Evertse-Schlickewei-Schmidt) \ 
Let \mN be a positive integer and let $\Lambda$ be a subgroup of $\big(\C^\times \big)^N$ of finite rank $r$ $-$then for any linear equation
\[
a_1 x_1 + \cdots + a_N x_N
\ = \ 
1
\]
over $\C$ with $a_1, \ldots, a_N \neq 0$ has at most
\[
\exp \big(\big( 6N\big)^{3N} (r + 1)\big)
\]
many nondegenerate solutions in $\Lambda$.

\vspace{0.2cm}

[Note: \ Only the fact that there exists a finite upper bound on the number of nondegenerate solutions in $\Lambda$ will actually be used.]
\end{x}
\vspace{0.3cm}

\begin{x}{\small\bf DISCUSSION} \ 
Let 
$q = \lindim_\Q S$
and fix a $\Q$-basis $\{s_1, \ldots, s_q\}$ of $\spanxs{\Q} S$.  
Let $s \in S$ $-$then there exist $c_1, \ldots, c_q \in \Q$ such that
\[
s \ = \ \sum\limits_{i=1}^q \hsx c_i s_i
\]
\qquad\qquad $\implies$
\[
0 
\ = \ 
f(s) 
\ = \ 
\lambda_1 \hsx
\prod\limits_{i=1}^q \hsx
e^{\textstyle\mu_1 c_i s_i} 
+
\cdots 
+
\lambda_n \hsx
\prod\limits_{i=1}^q \hsx
e^{\textstyle\mu_n c_i s_i}
\]
\qquad\qquad $\implies$
\[
\bigg(
\prod\limits_{i=1}^q \hsx
e^{\textstyle\mu_1 c_i s_i},
\ldots, 
\prod\limits_{i=1}^q \hsx
e^{\textstyle\mu_n c_i s_i}
\bigg)
\in \Gamma
\]
is a solution of the equation
\[
\lambda_1 x_1 + \cdots + \lambda_n x_n \ = \ 0.
\]
Put
\[
\lambda_j^\prime
\ = \ 
\bigg(
-\lambda_n \hsx
\prod\limits_{i=1}^q \hsx
e^{\textstyle\mu_n c_i s_i} 
\bigg)^{-1} \hsx \lambda_j
\qquad (1 \leq j \leq n - 1).
\]
Then
\allowdisplaybreaks
\begin{align*}
\lambda_1^\prime \hsx
\prod\limits_{i=1}^q \hsx
e^{\textstyle\mu_1 c_i s_i} 
&+
\cdots 
+
\lambda_{n-1}^\prime \hsx
\prod\limits_{i=1}^q \hsx
e^{\textstyle\mu_{n-1} c_i s_i}
\\[12pt]
&=\ 
\bigg(
-
\lambda_n \hsx
\prod\limits_{i=1}^q \hsx
e^{\textstyle\mu_n c_i s_i} 
\bigg)^{-1}
\hsx
\lambda_1 \hsx
\prod\limits_{i=1}^q \hsx
e^{\textstyle\mu_1 c_i s_i} 
\\[12pt]
&\hspace{1cm}
+ \cdots +  
\bigg(
-\lambda_n \hsx
\prod\limits_{i=1}^q \hsx
e^{\textstyle\mu_n c_i s_i} 
\bigg)^{-1}
\hsx
\lambda_{n-1} \hsx
\prod\limits_{i=1}^q \hsx
e^{\textstyle\mu_{n-1} c_i s_i}
\\[12pt]
&=\ 
\frac{\lambda_1}{-\lambda_n} \ 
\frac{\prod\limits_{i=1}^q \hsx e^{\textstyle\mu_1 c_i s_i}}{\prod\limits_{i=1}^q \hsx e^{\textstyle\mu_n c_i s_i}} 
\ 
+ \hsx \cdots \hsx +  
\ 
\frac{\lambda_{n-1}}{-\lambda_n} \ 
\frac{\prod\limits_{i=1}^q \hsx e^{\textstyle\mu_{n-1} c_i s_i}}{\prod\limits_{i=1}^q \hsx e^{\textstyle\mu_n c_i s_i}} 
\\[12pt]
&=\ 
-
\frac
{
\lambda_1 \hsx \prod\limits_{i=1}^q \hsx e^{\textstyle\mu_1 c_i s_i}
\ 
+ \hsx \cdots \hsx +  
\  
\lambda_{n-1} \hsx \prod\limits_{i=1}^q \hsx e^{\textstyle\mu_{n-1} c_i s_i}
}
{\lambda_n \hsx \prod\limits_{i=1}^q \hsx e^{\textstyle\mu_n c_i s_i}}
\\[12pt]
&=\ 
-
\frac{-\lambda_n \hsx \prod\limits_{i=1}^q \hsx e^{\textstyle\mu_n c_i s_i}} 
{\lambda_n \hsx \prod\limits_{i=1}^q \hsx e^{\textstyle\mu_n c_i s_i}} 
\\[12pt]
&= \ 
1
\end{align*}
\qquad\qquad$\implies$
\[
\bigg(
\prod\limits_{i=1}^q \hsx
e^{\textstyle\mu_1 c_i s_i}, 
\hsx \ldots, \hsx
\prod\limits_{i=1}^q \hsx
e^{\textstyle\mu_{n-1} c_i s_i}
\bigg)
\]
is a solution of the equation
\[
\lambda_1^\prime y_1 + \cdots + \lambda_{n-1}^\prime y_{n-1} \ = \ 1,
\]
all solutions which lie in some group $\Gamma_0$, a subgroup of $\Gamma$ of finite rank.  
Now 
apply \#6 to conclude that there are only finitely many nondegenerate solutions of 
\[
\lambda_1^\prime y_1 + \cdots + \lambda_{n-1}^\prime y_{n-1} \ = \ 1
\]
in $\Gamma_0$.
\end{x}
\vspace{0.3cm}

\begin{x}{\small\bf LEMMA} \ 
Let $\alpha$, $\beta \in S$ $(\alpha \neq \beta)$.  Suppose that
\[
\ba \ = \ (a_1, \ldots, a_n)
\]
is the solution of 
\[
\lambda_1 x_1 + \cdots + \lambda_n x_n \ = \ 0
\]
corresponding to $\alpha$ and
\[
\bb \ = \ (b_1, \ldots, b_n)
\]
is the solution of 
\[
\lambda_1 x_1 + \cdots + \lambda_n x_n \ = \ 0
\]
corresponding to $\beta$.  Then
\[
\ba \ \neq \ \bb.
\]
\\[-21pt]

PROOF \ 
If $\ba = \bb$, then for $j = 1, \ldots, n$, 
\[
\prod\limits_{i=1}^q \hsx
\big(e^{\textstyle\mu_j s_i}\big)^{c_{\ba,i}}
\ = \ 
\prod\limits_{i=1}^q \hsx
\big(e^{\textstyle\mu_j s_i}\big)^{c_{\bb,i}}
\]
iff
\[
\prod\limits_{i=1}^q \hsx
\big(e^{\textstyle\mu_j s_i}\big)^{\textstyle c_{\ba,i} - c_{\bb,i}} 
\ = \ 
1
\]
iff
\[
\mu_j
\sum\limits_{i=1}^q \hsx
s_i \big(\textstyle c_{\ba,i} - \textstyle c_{\bb,i}\big) \ 
\in \  2 \pi \sqrt{-1} \ \hsx \Z.
\]
So, for any $j = 1, \ldots, n$, 
\[
\sum\limits_{i=1}^q \hsx
s_i \big(\textstyle c_{\ba,i} - \textstyle c_{\bb,i}\big) 
\ = \ 
\frac{\textstyle 2 \pi \sqrt{-1}}{\textstyle \mu_j} \hsx N_j,
\]
where $N_j \in \Z$.  Therefore
\[
\frac{2 \pi \sqrt{-1}}{\mu_1} \hsx N_1
\ = \ 
\frac{2 \pi \sqrt{-1}}{\mu_2} \hsx N_2
\ = \ 
\cdots
\ = \ 
\frac{2 \pi \sqrt{-1}}{\mu_n} \hsx N_n
\]
\qquad\qquad $\implies$
\allowdisplaybreaks
\begin{align*}
&\mu_2 \ =\ \frac{\mu_1}{N_1} \hsx N_2
\\[12pt]
&\mu_3 \ =\ \frac{\mu_1}{N_1} \hsx N_3
\\[12pt]
&\hspace{.75cm} \vdots
\\[12pt]
&\mu_n \ =\ \frac{\mu_1}{N_1} \hsx N_n.
\end{align*}
\begin{spacing}{1.5}
\noindent Now put $\gamma = \ds\frac{\mu_1}{N_1}$ $-$then $f(z)$ is a polynomial in $e^{\gamma z}$, i.e., $f$ is simple, a contradiction since $f$ is not simple.
\vspace{0.3cm}

With this preparation, we are ready to tackle the proof of \#1 (as reformulated at the beginning: $f$, $g$ are distinct irreducibles with infinitely many zeros in common).  
It will be shown by induction on the length $n$ of $f$ that $g$ divides $f$.  
Since $f$, $g$ are distinct irreducibles, this is a contradiction.
\end{spacing}
\vspace{0.2cm}

\qquad \un{$n = 2$:} \quad
Suppose that
\[
f(z) 
\ = \ 
\lambda_1 e^{\mu_1 z} + \lambda_2 e^{\mu_2 z} 
\]
or still, 
\[
f(z) 
\ = \ 
\lambda_1 e^{\mu_1 z} \bigg(1 \hsx + \hsx \lambda_1^{-1} \lambda_2 \hsx e^{(\mu_2 - \mu_1)z}\bigg).
\]
Then $g(z)$ has infinitely many zeros in common with
\[
\bigg(1 \hsx + \hsx \lambda_1^{-1} \lambda_2 \hsx e^{(\mu_2 - \mu_1)z}\bigg)
\]
and as in \S57 there is an exponential polynomial of the form $\sin (T(z))$ dividing both $f(z)$ and $g(z)$.  
Since $g$ is irreducible, this implies that $g$ divides $f$.  
\\[12pt]
Proof:
\[
\begin{cases}
\ f = \sin(T) u\\[8pt]
\ g = \sin(T) v
\end{cases}
\qquad (u, v \in U \sE)
\]
\qquad\qquad $\implies$
\[
g v^{-1} \ = \ \sin(T)
\]
\qquad\qquad $\implies$
\[
f \ = \ g v^{-1} u.
\]
\vspace{0.2cm}

\begin{spacing}{1.5}
\qquad \un{Induction Hypothesis} \quad
Assume that for every exponential polynomial $h \neq g$ and of length $< n$, if $h$ and $g$ have infinitely many zeros in common, then $g$ divides $h$.
\end{spacing}
\vspace{0.5cm}

\qquad \un{$n > 2$:} \quad
Let as above
\[
\lambda_1^\prime y_1 + \cdots + \lambda_{n-1}^\prime y_{n-1} \ = \ 1
\]
be the linear equation associated with
\[
f(z) 
\ = \ 
\lambda_1 e^{\mu_1 z} + \cdots +  \lambda_n e^{\mu_n z}.
\]
Then $\Gamma_0$ contains just a finite number of nondegenerate solutions of this equation (cf. \#7).  
Consider the equation
\[
\lambda_1 x_1 + \cdots + \lambda_n x_n \ = \ 0.
\]
Then each $s \in S$ gives rise to a solution and since \mS is infinite, it follows from \#8 that this equation has infinitely many distinct solutions
\[
\bomega_s 
\ \equiv \ 
\big(\omega_1^{(s)}, \ldots, \omega_n^{(s)}\big) \in \Gamma,
\]
where 
\[
\omega_1^{(s)} 
\ = \ 
\prod\limits_{i=1}^q \hsx
e^{\mu_1 c_i s_i}, \hsx
\ldots, \hsx
\omega_n^{(s)} 
\ = \ 
\prod\limits_{i=1}^q \hsx
e^{\mu_n c_i s_i} \hsx .
\]
Each $\bomega_s$ can be turned into a solution of
\[
\lambda_1^\prime y_1 + \cdots + \lambda_{n-1}^\prime y_{n-1} \ = \ 1
\]
by simply removing its last component.  
Bottom line: \ There are an infinity of distinct solutions to
\[
\lambda_1^\prime y_1 + \cdots + \lambda_{n-1}^\prime y_{n-1} \ = \ 1,
\]
any such being determined by an $s \in S$.  
Moreover all but finitely many are degenerate (cf. \#6) and for a degenerate $\bomega_s$  there exists a proper nonempty $J_s \subset \{1, \ldots, n\}$ such that
\[
\sum\limits_{j \in J_s} \ 
\lambda_j \omega_j^{(s)} 
\ = \ 0.
\]
In fact, if
\[
\sum\limits_{j \in J_s} \hsx
\lambda_j^\prime \ 
\prod\limits_{i=1}^q \ 
e^{\mu_j c_i s_i}
\ = \ 
0,
\]
then
\[
\sum\limits_{j \in J_s} \hsx
\bigg(
- \lambda_n \hsx
\prod\limits_{i=1}^q \hsx
e^{\mu_n c_i s_i}
\bigg)^{-1}
\lambda_j \hsx
\prod\limits_{i=1}^q \hsx
e^{\mu_j c_i s_i}
\ = \ 0
\]
\qquad\qquad $\implies$
\[
\sum\limits_{j \in J_s} \hsx
\lambda_j \hsx
\prod\limits_{i=1}^q \hsx
e^{\mu_j c_i s_i}
\ = \ 0
\]
\qquad\qquad $\implies$
\[
\sum\limits_{j \in J_s} \hsx
\lambda_j \hsx
\omega_j^{(s)}
\ = \ 0.
\]
\\[-.25cm]

\noindent Owing now to the Box Principle (cf. \S7, \#15), we can find a proper nonempty subset
\[
T 
\ = \ 
\{j_1, \ldots, j_t\} \  \subset \  \{1, \ldots, n\}
\]
such that for infinitely many $s \in S$, 
\[
\sum\limits_{j \in T} \ 
\lambda_j \hsx
\omega_j^{(s)}
\ = \ 0.
\]
Therefore the equation
\[
\lambda_{j_1} x_{j_1} + \cdots + \lambda_{j_t} x_{j_t} 
\ = \ 
0
\]
has infinitely many solutions corresponding to common zeros of $f$, $g$.
\end{x}
\vspace{0.3cm}

\begin{x}{\small\bf LEMMA} \ 
$g$ divides $f$.
\vspace{0.5cm}

PROOF \ 
Put
\[
f_T(z) 
\ = \ 
\lambda_{j_1} e^{\mu_{j_1}z} + \cdots + \lambda_{j_t} e^{\mu_{j_t}z}.
\]
\begin{spacing}{1.5}
\noindent Then $g$ has infinitely many zeros in common with $f_T$ which are also zeros of $f$, thus also zeros of $f - f_T$.  
Both $f_T$ and $f - f_T$ are elements of $\sE$ of length strictly less
than $n$ (the length of $f$).  
Thanks to \S56, \#9, $g$ has infinitely many zeros in common with either an irreducible or a simple factor of $f_T$ in $\sE$, call this factor $h_T$.  
If $h_T$ is simple, then we are in Case 1 and $g$, $h_T$ must have a common divisor.  
Since $g$ is irreducible, it then divides $h_T$ ($g = a u$, $h_T = a b$, $g u^{-1} = a$, $h_T = g u^{-1} b$).  
If $h_T$ is irreducible, then it is either a unit multiple of $g$, in which case $g$ divides $h_T$, or $g$ and $h_T$ are distinct irreducibles, in which case $g$ divides $h_T$ (induction hypothesis).  
So, in all cases $g$ divides $h_T$, thus it also divides $f_T$.  
Analogously, $g$ divides $f - f_T$.  Therefore $g$ divides $f$.
\end{spacing}
\end{x}
\vspace{0.1cm}

\begin{x}{\small\bf \un{N.B.}} \ 
\#9 is the sought for contradiction.
\end{x}
\vspace{0.3cm}


%% file: _59_differential_algebra.tex
\chapter{
$\boldsymbol{\S}$\textbf{59}.\quad  DIFFERENTIAL ALGEBRA}
\setlength\parindent{2em}
\setcounter{theoremn}{0}
\renewcommand{\thepage}{\S59-\arabic{page}}

\ \indent 
Let $\K/\bk$ be fields of characteristic 0, where \bk is algebraically closed in $\K$.
\vspace{0.5cm}

\begin{x}{\small\bf DEFINITION} \ 
Suppose that \ \mV is a \ $\K$-vector space $-$then a linear map \ $\td:\K \ra V$ is a 
\un{\bk-derivation}
if $\forall \ x, y \in \K$, 
\[
\td(x y) \ = \ x \td(y) + y \td(x)
\]
and if $\forall \ a \in \bk$, 
\[
\td(a) \ = \ 0.
\]

[Note: \ 
In particular, $\td(1) = 0$.]
\end{x}
\vspace{0.3cm}

\begin{spacing}{1.5}
\begin{x}{\small\bf RAPPEL}  \ 
There is a $\K$-vector space $\Omega_{\K / \bk}$ and a \bk-derivation $\td_{\K / \bk} : \K \ra \Omega_{\K / \bk}$ with the property that for any $\K$-vector space \mV and any \bk-derivation $\td:\K \ra V$ there is a unique $\K$-linear map 
$\xi: \Omega_{\K / \bk} \ra V$ such that $\td = \xi \circ \td_{\K / \bk}$:
\[
\begin{tikzcd}[sep=huge]
\K \ar{d}[swap]{\td}\ar{rr}{\td_{\K/\bk}} &&\Omega_{\K/\bk} \arrow[d,shift right=0.5,dash] \arrow[d,shift right=-0.5,dash]\\
V &&\Omega_{\K/\bk} \ar{ll}{\xi}
\end{tikzcd}
.
\]
\end{x}
\end{spacing}
\vspace{0.3cm}

\begin{x}{\small\bf SCHOLIUM} \ 
Associated with every \bk-derivation $\td: \K \ra \K$ there is a unique derivation 
$\tD : \Omega_{\K / \bk} \ra \Omega_{\K / \bk}$ such that $\forall \ x_1, x_2 \in \K$, 
\[
\tD(x_1 \td_{\K / \bk} (x_2)) 
\ = \ 
\td(x_1) \td_{\K / \bk} (x_2) + x_1 \td_{\K / \bk} (\td(x_2)).
\]
\end{x}
\vspace{0.3cm}

\begin{x}{\small\bf SUBLEMMA} \ 
Suppose given a \bk-derivation $\td:\K \ra V$ $-$then for $y \in \K$,
$z \in \K^\times$,
\[
\tD \bigg( \td_{\K / \bk}  (y) \hsx - \hsx \frac{\td_{\K / \bk} (z)}{z}\bigg) 
\ = \ 
0
\]
if 
\[
\td(y) \ = \ \frac{\td(z)}{z}.
\]
\vspace{0.2cm}

PROOF \ 
The LHS equals
\[
\td_{\K / \bk}  (\td(y)) - \frac{1}{z} \td_{\K / \bk} (\td(z)) + \frac{\td(z)}{z^2} \td_{\K / \bk} (z)
\]
or still,
\allowdisplaybreaks
\begin{align*}
\td_{\K / \bk}  \bigg(\frac{\td(z)}{z}\bigg) - 
&\frac{1}{z} \td_{\K / \bk}  (\td(z)) + \frac{\td(z)}{z^2} \td_{\K / \bk} (z) \ 
\\[12pt]
&=\ 
\frac{z \td_{\K / \bk} (\td(z)) - (\td(z)) \td_{\K / \bk}(z)}{z^2} 
- \frac{1}{z} \td_{\K / \bk} (\td(z)) + \frac{\td(z)}{z^2} \td_{\K / \bk} (z)
\\[12pt]
&=\ 
0.
\end{align*}
\end{x}
\vspace{0.3cm}

\begin{x}{\small\bf SUBLEMMA} \ 
Suppose given a \bk-derivation $\td:\K \ra V$ $-$then for $y \in \K$, 
\[
\tD ( \td_{\K / \bk} (y)) \ = \ 0
\]
if $\td(y) = 1$.
\vspace{0.3cm}

PROOF \ 
The LHS equals
\allowdisplaybreaks
\begin{align*}
\tD(1 \td_{\K / \bk}(y)) \ 
&=\ 
\td(1) \td_{\K / \bk}(y) + 1 \td_{\K / \bk}(\td(y)) 
\\[12pt]
&=\ 
0 + \td_{\K / \bk}(1)
\\[12pt]
&=\ 
0.
\end{align*}
\end{x}
\vspace{0.3cm}

\begin{x}{\small\bf NOTATION} \ 
Given $y_i \in \K$, $z_i \in \K^\times$ $(i = 1, \ldots, n)$, put
\[
\omega_i 
\ = \ 
\td_{\K / \bk}(y_i) - \frac{d_{\K / \bk}(z_i)}{z_i} \hsx \in \hsx  \Omega_{\K / \bk}.
\]
\end{x}
\vspace{0.3cm}

\begin{x}{\small\bf LEMMA} \ 
Suppose that $\td:\K \ra V$ is a $\bk$-derivation.  
Assume that $\td(y_1) = 1$ and that $y_i \in\K$, $z_i \in \K^\times$ are connected by the relation
\[
\td(y_i) \ = \ \frac{\td(z_i)}{z_i} \qquad (i = 1, \ldots, n).
\]
Then for $f_1, \ldots, f_n, \hsx g \in \K$, 
\allowdisplaybreaks
\begin{align*}
\tD\big(\sum\limits_i \hsx f_i \hsx \omega_i + g \hsx \td_{\K/\bk}(y_1)\big)\ 
&=\ 
\sum\limits_i \hsx (\td(f_i) \hsx \omega_i  + f_i \hsx \tD  \omega_i ) + \td(g) \hsx \td_{\K/\bk}(y_1)+ g \hsx \tD (\td_{\K/\bk}(y_1))
\\[12pt]
&=\ 
\sum\limits_i \hsx (\td(f_i) \hsx \omega_i  + f_i \hsx 0) + \td(g) \hsx \td_{\K/\bk}(y_1)+ g \hsx 0
\\[12pt]
&=\ 
\sum\limits_i \hsx \td(f_i) \omega_i  + \td(g) \td_{\K/\bk}(y_1).
\end{align*}
\end{x}
\vspace{0.75cm}

In what follows, $\td:\K \ra \K$ is a derivation such that
\[
\Ker \hsx \td \ = \ \bk \quad (\supset \Q).
\]
\vspace{0.3cm}

\begin{x}{\small\bf CRITERION} \ 
Let $\K \supset \F \supset \bk$, where $\F$ is a field and
\[
\trdeg_{\bk} \hsx \F \ < \ \infty.
\]
Denote by \mE the $\K$-vector subspace of $\Omega_{\K / \bk}$ generated by $\td_{\K / \bk} \F$ $-$then
\[
\dim_{\K} E \ = \ \trdeg_{\bk} \hsx \F.
\]
\end{x}
\vspace{0.3cm}


\begin{spacing}{1.5}
\begin{x}{\small\bf EXAMPLE} \ 
Take $\F = \K$ $-$then
\[
\dim_{\K} \Omega_{\K / \bk} \ = \ \trdeg_{\bk} \hsx \K.
\]

[\quad \textbullet \quad
If $x_1, \ldots, x_n \in \K$ are algebraically dependent over \bk, then 
$\td_{\K / \bk}(x_1), \ldots,$ $\td_{\K / \bk}(x_n)$ $\in \Omega_{\K / \bk}$ 
are linearly dependent over $\K$.
\vspace{0.2cm}

\quad \hsx  \textbullet \quad
If $x_1, \ldots, x_n \in \K$ are algebraically independent over \bk, then $\td_{\K / \bk}(x_1), \ldots,$ $\td_{\K / \bk}(x_n)$ 
$\hsx \in \hsx \Omega_{\K / \bk}$ 
are linearly independent over $\K$.]
\vspace{0.5cm}

[Note: \ 
Therefore $\td_{\K / \bk} = 0$ iff $x$ is algebraic over \bk.]
\\[.0cm]

Keep to the setup of \#7 and in \#8, let
\[
\F = \bk(y_1, \ldots, y_n, z_1, \ldots, z_n)
\]
and suppose that $\trdeg_{\bk}\hsx \F < n+1$ $-$then there are elements $f_1, \ldots, f_n, g \in \K$ not all zero such that
\[
\sum\limits_i \hsx f_i \omega_i + g \td_{\K/\bk}(y_1) \ = \ 0.
\]
It can be assumed that  $f_1, \ldots, f_n, g \in \K$ have been chosen so that a minimal number of them are nonzero and at least one of them is 1.
\vspace{0.1cm}

Write
\allowdisplaybreaks
\begin{align*}
0\ 
&=\ 
\tD \hsx 0
\\[12pt]
&=\ 
\tD\hsx\big(\sum\limits_i \hsx f_i \omega_i + g \td_{\K/\bk}(y_1) \big)
\\[12pt]
&=\ 
\sum\limits_i \hsx \td(f_i) \omega_i + \td(g) \td_{\K/\bk}(y_1) 
\end{align*}
to conclude by minimality that
\[
\td(f_1) \ = \ 0, \ldots, \ \td(f_n) \ = \ 0, \ \td(g) \ = \ 0,
\]
thus
\[
f_1 \in \bk, \ldots, \ f_n \in \bk, \ g \in \bk,
\]
the field of constants of $\td$ being \bk (by hypothesis).  
Bearing in mind that
\[
\sum\limits_i \hsx f_i \omega_i + g \td_{\K/\bk}(y_1) \ = \ 0,
\]
let $c_i = f_i$, $c_0 = g$, hence
\[
\sum\limits_i \hsx c_i \omega_i +c_0 \td_{\K/\bk}(y_1) \ = \ 0.
\]
\end{x}
\end{spacing}
\vspace{0.3cm}

\begin{x}{\small\bf NOTATION} \ 
Put
\[
C \ = \ c_0 + c_1 y_1 + \cdots + c_n y_n.
\]
\end{x}
\vspace{0.3cm}

\begin{x}{\small\bf LEMMA} \ 
\[
\td_{\K/\bk} \hsx (C) \ = \ 
\sum\limits_i \hsx 
c_i \hsx \frac{\td_{\K/\bk}(z_i)}{z_i}.
\]
\vspace{0.2cm}

PROOF \ 
In fact,
\[
\sum\limits_i \hsx 
c_i \hsx \omega_i + c_0 \hsx \td_{\K/\bk} \hsx (y_1) \ = \ 0
\]
or still, 
\[
\sum\limits_i \hsx 
c_i \hsx \bigg(\td_{\K/\bk} \hsx (y_i) - \frac{\td_{\K/\bk}(z_i)}{z_i}\bigg) 
+ c_0 \td_{\K/\bk} \hsx (y_1) \ = \ 0
\]
\qquad\qquad $\implies$
\[
\sum\limits_i \hsx 
c_i \hsx \td_{\K/\bk} \hsx (y_i) 
+ c_0 \hsx \td_{\K/\bk} \hsx (y_1) 
\ = \
\sum\limits_i \hsx 
c_i \hsx \frac{\td_{\K/\bk}(z_i)}{z_i}
\]
\qquad\qquad $\implies$
\[
 c_0 \hsx \td_{\K/\bk} \hsx (y_1) 
 + 
 \sum\limits_i \hsx 
c_i \hsx \td_{\K/\bk} \hsx (y_i) 
\ = \
\sum\limits_i \hsx 
c_i \hsx \frac{\td_{\K/\bk}(z_i)}{z_i}
\]
\qquad\qquad $\implies$
\[
\td_{\K/\bk} \hsx (C) 
\ = \ 
\sum\limits_i \hsx 
c_i \hsx \frac{\td_{\K/\bk}(z_i)}{z_i}.
\]
\vspace{0.1cm}

Suppose that $c_1, \ldots, c_L$ is a $\Q$-basis for $c_1, \ldots, c_n$, hence
\[
c_i 
\ = \ 
\sum\limits_{\ell = 1}^L \hsx
q_{\ell, i} c_\ell \qquad (i = 1, \ldots, n).
\]
Here, at least a priori, the $q_{\ell, i} \in \Q$ but there is no loss of generality in taking $q_{\ell, i} \in \Z$.
\\

Accordingly
\allowdisplaybreaks
\begin{align*}
\td_{\K/\bk} \hsx (C) \ 
&=\  
\sum\limits_{i=1}^n \  
c_i \hsx 
\frac{\td_{\K/\bk}(z_i)}{z_i}
\\[12pt]
&=\  
\sum\limits_{i=1}^n \hsx 
\sum\limits_{\ell = 1}^L \ 
q_{\ell, i} \hsx c_\ell \hsx
\frac{\td_{\K/\bk}(z_i)}{z_i}
\\[12pt]
&=\  
\sum\limits_{\ell = 1}^L \ 
c_\ell \hsx \bigg(
\sum\limits_{i=1}^n \  
q_{\ell, i} \hsx
\frac{\td_{\K/\bk}(z_i)}{z_i}
\bigg)
\\[12pt]
&=\  
\sum\limits_{\ell = 1}^L \ 
c_\ell \hsx
\frac{\td_{\K/\bk}(w_\ell)}{w_\ell},
\end{align*}
where
\[
w_\ell
\ = \ 
\prod\limits_{i=1}^n \ 
z_i^{q_{\ell, i}}.
\]
\end{x}
\vspace{0.3cm}

\begin{x}{\small\bf LEMMA} \ 
Let $a_1, \ldots, a_L \in \bk$ be linearly independent over $\Q$, let $u_1, \ldots, u_L \in \K^\times$, let $v \in \K$, and assume that
\[
\td_{\K/\bk} \hsx (v) 
\ = \ 
\sum\limits_{\ell = 1}^L \ 
a_\ell \hsx
\frac{\td_{\K/\bk}(u_\ell)}{u_\ell}.
\]
Then
\[
\td_{\K/\bk} \hsx (u_1) 
\ = \ 0, 
\ldots, 
\td_{\K/\bk} \hsx (u_L) 
\ = \ 0.
\]
\end{x}
\vspace{0.3cm}

\begin{x}{\small\bf APPLICATION} \ 
Take $a_1 = c_1, \ldots, a_L = c_L$, take $v = C$, and take
\[
u_1 = w_1, \ldots, u_L = w_L.
\]
Then
\[
\td_{\K/\bk} \hsx (w_1) = 0, 
\ldots, 
\td_{\K/\bk} \hsx (w_L) = 0.
\]
\end{x}
\vspace{0.3cm}

\begin{x}{\small\bf \un{N.B.}} \ 
Since the standing assumption is that \bk is algebraically closed in $\K$, each $w_\ell \in \bk$ (cf. \#9).
\end{x}
\vspace{0.3cm}

\begin{x}{\small\bf APPLICATION} \ 
For $\ell = 1, \ldots, L$, 
\[
\prod\limits_{i=1}^n \ 
z_i^{q_{\ell, i}} \in \bk.
\]
Finally
\[
w_\ell \in \bk \implies \td(w_\ell) \ = \ 0
\]
\qquad\qquad $\implies$
\allowdisplaybreaks
\begin{align*}
0\ 
&=\ 
\frac{\td(w_\ell)}{w_\ell}
\\[12pt]
&=\
\sum\limits_{j=1}^L \ 
q_{\ell, j} \hsx
\frac{\td(z_j)}{z_j} 
\\[12pt]
&=\ 
\sum\limits_{j=1}^L \ 
q_{\ell, j} \hsx \td(y_j)
\\[12pt]
&=\ 
\td\hsx \bigg(
\sum\limits_{j=1}^L \ 
q_{\ell, j}  y_j \bigg)
\end{align*}
\qquad\qquad $\implies$
\[
\sum\limits_{j=1}^L \ 
q_{\ell, j} y_j \in \bk.
\]
\end{x}
\vspace{0.3cm}

\begin{x}{\small\bf SCHOLIUM} \ 
There exist integers $m_1, \ldots, m_n$ not all zero such that
\[
\sum\limits_{i=1}^n \  m_i y_i \in \bk.
\]

Recall: 
\vspace{0.2cm}

\qquad \textbullet \quad
$y_i \in \K$, $z_i \in \K^\times$, and 
\[
\td(y_i) \ = \ \frac{\td(z_i)}{z_i} \qquad (i = 1, \ldots, n).
\]
\vspace{0.2cm}

\qquad \textbullet \quad
$\F = \bk(y_1, \ldots, y_n, z_1, \ldots, z_n)$ and
\[
\trdeg_{\bk} \hsx \F < n + 1.
\]

Then under these assumptions: 
\vspace{0.2cm}

\qquad\qquad (1) \quad
There are $m_1, \ldots, m_n \in \Z$ not all zero such that
\[
\prod\limits_{i=1}^n \ 
z_i^{m_i} \in \bk.
\]
\vspace{0.2cm}

\qquad\qquad (2) \quad
There are $m_1, \ldots, m_n \in \Z$ not all zero such that
\[
\sum\limits_{i=1}^n \ 
m_i y_i \in \bk.
\]
\end{x}
\vspace{0.3cm}

\begin{x}{\small\bf STATEMENT} \ 
Maintain the supposition that
\[
\td(y_i) \ = \ 
\frac{\td(z_i)}{z_i} 
\qquad (i = 1, \ldots, n)
\]
but assume that the $y_i$ are $\Q$-linearly independent modulo $\bk$, i.e., 
\[
\sum\limits_{i = 1}^n \  q_i y_i \hsx \in \hsx \bk
\ \implies \ 
q_i \hsx = \hsx 0 
\qquad (i = 1, \ldots, n).
\]
Then
\[
\trdeg_{\bk} \hsx \F \ \geq \  n + 1.
\]
\end{x}
\vspace{0.3cm}


%% file: _60_formal_schanuel.tex
\chapter{
$\boldsymbol{\S}$\textbf{60}.\quad  FORMAL SCHANUEL}
\setlength\parindent{2em}
\setcounter{theoremn}{0}
\renewcommand{\thepage}{\S60-\arabic{page}}

\ \indent 
This is a version of Schanuel that can be established rigorously.  
However, before proceeding to the particulars, let us review the situation.

As it is usually formulated, Schanuel's conjecture is the following statement (cf. \S47, \#1).
\vspace{0.5cm}

\begin{x}{\small\bf CONJECTURE} \ 
Suppose that $x_1, \ldots, x_n$ are $\Q$-linearly independent complex numbers $-$then among the $2n$ numbers
\[
x_1, \ldots, x_n, \hsx e\strutv^{x_1}, \ldots, e\strutv^{x_n}, 
\]
at least $n$ are algebraically independent over $\Q$, i.e., 
\[
\trdegQ \Q(x_1, \ldots, x_n, \hsx e\strutv^{x_1}, \ldots, e\strutv^{x_n}) \ \geq \ n.
\]
\end{x}
\vspace{0.3cm}

There are other equivalent formulations.  E.g.: \ $\forall \ \bx$, 
\[
\delta(\bx) \ \geq \ 0 \qquad \text{(cf. \S47, \#24 and \#27)}.
\]

Here are two more.
\vspace{0.5cm}

\begin{x}{\small\bf CONJECTURE} \ 
Suppose that $x_1, \ldots, x_n$ are complex numbers such that
\[
\trdegQ \Q(x_1, \ldots, x_n, e\strutv^{x_1}, \ldots, e\strutv^{x_n})
\]
is $< n$ $-$then there are integers $m_1, \ldots, m_n$ not all zero such that 
\[
\sum\limits_{i=1}^n \hsx m_i x_i \ = \ 0.
\]
\end{x}
\vspace{0.3cm}

\begin{x}{\small\bf CONJECTURE} \ 
Suppose that $x_1, \ldots, x_n$ are complex numbers such that
\[
(x_1, \ldots, x_n, e\strutv^{x_1}, \ldots, e\strutv^{x_n})
\]
lie 
in an algebraic subvariety \mV of $\C^{2n}$ defined over $\Q$ and of dimension strictly less than $n$ $-$then there are integers 
$m_1, \ldots, m_n$ not all zero such that 
\[
\sum\limits_{i=1}^n \hsx m_i x_i \ = \ 0.
\]

\vspace{0.2cm}

[The assumption that
\[
(x_1, \ldots, x_n, e\strutv^{x_1}, \ldots, e\strutv^{x_n}) \in V \qquad (\dim V < n)
\]
forces
\[
\trdegQ \Q(x_1, \ldots, x_n, e\strutv^{x_1}, \ldots, e\strutv^{x_n}) \ < \ n.]
\]
\end{x}
\vspace{0.3cm}

We shall now turn to a setting in which an analog of Schanuel's conjecture is true.
\vspace{0.3cm}

\begin{x}{\small\bf NOTATION} \ 
Let \mR be a commutative ring with 1 $-$then
\[
R[[X]]
\]
is the ring of formal power series over \mR, a typical element of which is denoted by
\[
f(X) 
\ = \ 
\sum\limits_{n=0}^\infty \hsx a_n X^n 
\qquad (\forall \ n, \ a_n \in R).
\]
\end{x}
\vspace{0.3cm}

\begin{x}{\small\bf \un{N.B.}}  \ 
If \mR is an \mE-ring, then $R[[X]]$ is also an \mE-ring.
\vspace{0.2cm}

[Given $f \in R[[X]]$, write
\[
f \ = \ a_0 + g 
\qquad \bigg(g(X) = \sum\limits_{n=1}^\infty \hsx a_n X^n\bigg)
\]
and put
\[
\exp (f) \ = \ E(a_0) \exp (g),
\]
where $E(a_0)$ is that derived from \mR and
\[
\exp (g) 
\ = \ 
\sum\limits_{n=0}^\infty \hsx \frac{(g)^n}{n!}.]
\]
\end{x}
\vspace{0.3cm}

\begin{x}{\small\bf CONSTRUCTION} \ 
Let
\[
\begin{cases}
\ f(X) \ = \ \sum\limits_{n=1}^\infty \hsx  a_n X^n \ = \ a_1 X + a_2 X^2 + \cdots\\[15pt]
\ g(X) \ = \ \sum\limits_{n=0}^\infty \hsx  b_n X^n \ = \ b_0 + b_1 X^1 + b_2 X^2 + \cdots
\end{cases}
.
\]
Then their
\un{composite}
\index{composite} 
$g \circ f$ is the formal power series
\[
g(f(X)) 
\ = \ 
\sum\limits_{n=0}^\infty \hsx b_n (f(X))^n 
\ = \ 
\sum\limits_{n=0}^\infty \hsx c_n X^n.
\]
\end{x}
\vspace{0.3cm}

\begin{x}{\small\bf REMARK}\ 
The foregoing operation is valid only when $f(X)$ has no constant term (for then each $c_n$ depends on but a finite number of coefficients of $f(X)$ and $g(X)$).
\vspace{0.2cm}

[To illustrate, let
\[
\exp (X) 
\ = \ 
1 + X + \frac{X^2}{2!} + \frac{X^3}{3!} + \frac{X^4}{4!} + \cdots\ .
\]
Then it makes sense to form
\[
\exp (\exp (X) - 1) \ = \ 1 + X + X^2 + \frac{5}{6} X^3 + \frac{5}{8} X^4 + \cdots
\]
but
\[
\exp (\exp (X) ) \cdots \c ?]
\]

\vspace{0.2cm}

[Note: \ 
If $f(X)$ has no constant term, then $E(a_0) = E(0) = 1$ and one can form 
\[
\exp \circ f,
\]
which agrees with \#5.]
\end{x}
\vspace{0.3cm}

\begin{x}{\small\bf LEMMA} \ 
If \mR is an integral domain, then so is $R[[X]]$.
\end{x}
\vspace{0.3cm}

\begin{x}{\small\bf DEFINITION} \ 
A 
\un{formal Laurent series}
\index{formal Laurent series} 
over \mR is a series of the form
\[
f(X) 
\ = \ 
\sum\limits_{n \in \Z} \hsx a_n X^n,
\]
where $a_n = 0$ for all but finitely many negative indices $n$.
\end{x}
\vspace{0.3cm}

\begin{x}{\small\bf \un{N.B.}}  \ 
The formal Laurent series form a ring, denoted by 
$R((X))$.
\end{x}
\vspace{0.3cm}

\begin{x}{\small\bf LEMMA} \ 
If $R = \K$ is a field, then $\K((X))$ is a field.
\vspace{0.2cm}

[Note: \ 
$\K((X))$ can be identified with the field of fractions of the integral domain $\K[[X]]$.]
\end{x}
\vspace{0.3cm}

\begin{x}{\small\bf DEFINITION} \ 
Take $R = \K$ of characteristic 0 $-$then the formal derivative of the formal Laurent series
\[
f(X) 
\ = \ 
\sum\limits_{n \in \Z} \hsx a_n X^n
\]
is
\[
f^\prime 
\ = \ 
\partial f 
\ = \ 
\sum\limits_{n \in \Z} \hsx na_n X^{n-1}.
\]
\end{x}
\vspace{0.3cm}

\begin{x}{\small\bf \un{N.B.}} \ 
\[
\partial : \K((X)) \ra \K((X))
\]
is a $\K$-derivation $(\Ker \partial = \K)$.
\end{x}
\vspace{0.3cm}

Having dispensed with the formalities, specialize and take per \S59,
\[
\K \ = \ \C((X)), \quad \bk \ = \ \C, \quad \td \ = \ \partial.
\]
Let
\[
y_1 \in  X \C[[X]], \ldots, y_n \in X \C[[X]]
\]
be $\Q$-linearly independent and put
\[
z_1 = \exp(y_1), \ldots, z_n = \exp(y_n).
\]
\vspace{0.1cm}

\begin{x}{\small\bf THEOREM} \ 
\[
\trdeg_{\C} \hsx \C(y_1, \ldots, y_n, z_1, \ldots, z_n) \ \geq \ n + 1.
\]

[Quote \S59, \#17 (obviously, if the $y_i$ are $\Q$-linearly independent, then they are $\Q$-linearly independent modulo $\C$).]
\end{x}
\vspace{0.3cm}

This result can be rephrased.
\vspace{0.5cm}

\begin{x}{\small\bf RAPPEL} \ 
(cf. \S46, \#20) \ Given fields $\bk \hsx \subset \hsx  \K \subset \LL$,
\[
\trdeg_{\bk} \hsx (\LL / \bk) \ = \ \trdeg_{\K} \hsx (\LL / \K)  + \trdeg_{\bk} \hsx (\K / \bk).
\]
Abbreviate
\[
(y_1, \ldots, y_n, z_1, \ldots, z_n)
\]
to
\[
(\by, \bz).  
\]

Take in \#15
\[
\bk = \C, \quad \K = \C(X), \quad \LL = \C(X) (\by, \bz).  
\]
Then
\[
\trdeg_{\C} \hsx \C(X) (\by, \bz) 
\ = \ 
\trdeg_{\C(X)} \hsx \C(X) (\by, \bz)  + \trdeg_{\C} \hsx \C(X).  
\]
From \#14
\[
\trdeg_{\C} \hsx \C(X) (\by, \bz) 
\ > \ 
\trdeg_{\C} \hsx \C (\by, \bz) 
\ \geq \ 
n + 1.
\]
And
\[
\trdeg_{\C} \hsx \C(X)  \ = \ 1.
\]
Therefore
\allowdisplaybreaks
\begin{align*}
n + 1 \ 
&\leq\ 
\trdeg_{\C} \hsx \C(X) (\by, \bz) 
\\[12pt]
&=\ 
\trdeg_{\C(X)} \hsx \C(X) (\by, \bz) + 1
\end{align*}
\qquad\qquad $\implies$
\[
n \ \leq \ \trdeg_{\C(X)} \hsx \C(X) (\by, \bz).
\]
\end{x}
\vspace{0.3cm}

\begin{x}{\small\bf SUMMARY} \ 
The fact that
\[
\trdeg_{\C(X)} \hsx \C(X) (y_1, \ldots, y_n, z_1, \ldots, z_n) \ \geq \ n
\]
is 
\un{formal Schanuel},
\index{formal Schanuel} 
a result due to J. Ax.  
It is the power series analog of \#1 (which remains conjectural).
\end{x}
\vspace{0.3cm}

\begin{x}{\small\bf \un{N.B.}} \ 
\allowdisplaybreaks
\begin{align*}
\C \ \subset \ \ &\C[X] \ \subset \ \ \C[[X]]
\\[12pt]
&\ \ \cap \hspace{1.5cm} \cap \qquad .
\\[12pt]
 &\C[X] \ \subset \ \ \C((X))
\end{align*}
\end{x}
\vspace{0.3cm}


%% file: _61_an_arithmetic_criterion.tex
\chapter{
$\boldsymbol{\S}$\textbf{61}.\quad  AN ARITHMETIC CRITERION}
\setlength\parindent{2em}
\setcounter{theoremn}{0}
\renewcommand{\thepage}{\S61-\arabic{page}}

\ \indent 
Recall:
\vspace{0.3cm}

\begin{x}{\small\bf SCHANUEL'S CONJECTURE} \ 
Suppose that $x_1, \ldots, x_n$ are $\Q$-linearly independent complex numbers $-$then 
\[
\trdegQ \Q(x_1, \ldots, x_n, e\strutv^{x_1}, \ldots, e\strutv^{x_n}) 
\ \geq \ 
n.
\]
\end{x}
\vspace{0.3cm}

\begin{x}{\small\bf NOTATION}\ 
The symbol $\sD$ stands for the derivation
\[
\sD 
\ = \ 
\frac{\partial}{\partial X_0} + X_1 \frac{\partial}{\partial X_1}
\]
in the ring $\C[X_0,X_1]$.
\end{x}
\vspace{0.3cm}

\begin{x}{\small\bf DEFINITION} \ 
The 
\un{height}
\index{height (of a polynomial)} 
$H(P)$ of a polynomial $P \in \C[X_0,X_1]$ is the maximum of the absolute values of its coefficients.
\end{x}
\vspace{0.3cm}

\begin{x}{\small\bf DATA} \ 
Let $n$ be a positive integer, let $x_1, \ldots, x_n$ be $\Q$-linearly independent complex numbers, and let 
$\alpha_1 \in \C^\times, \ldots, \alpha_n \in \C^\times$.
\end{x}
\vspace{0.3cm}

\begin{x}{\small\bf PARAMETERS} \ 
Let $s_0$, $s_1$, $t_0$, $t_1$, $u$ be positive real numbers subject to
\[
\max\{1, t_0, 2 t_1\} 
\ < \ 
\min\{s_0, 2 s_1\}
\]
and
\[
\max\{s_0, s_1 + t_1\} 
\ < \ 
u
\ < \ 
\frac{1}{2} \hsx (1 + t_0 + t_1).
\]
\end{x}
\vspace{0.3cm}

\begin{spacing}{1.5} 
\begin{x}{\small\bf ROY'S CONJECTURE} \ 
In the presence of \#4 and \#5, assume that for any sufficiently large positive integer \mN, there exists a nonzero polynomial 
$P_N \in \Z[X_0,X_1]$ with partial degree $\leq \ N\strutw^{t_0}$ in $X_0$, with partial degree 
$\leq \ N\strutw^{t_1}$ in $X_1$, and with height $\leq \ e^N$ which satisfies
\[
\abs{\big(\sD^k P_N \big) \hsx 
\bigg(\sum\limits_{j=1}^n \hsx m_j x_j, \  
\prod\limits_{j=1}^n \hsx \alpha_j^{m_j} \bigg)} 
\ \leq \ 
\exp \big(- N^u\big)
\]\\
for all nonnegative integers $k, m_1, \ldots, m_n$, where
\[
k 
\ \leq \ 
N\strutw^{s_0} 
\quad \text{and} \quad 
\max\{m_1, \ldots, m_n\} 
\ \leq \  
N\strutw^{s_1}.
\]
Then
\[
\trdegQ \Q(x_1, \ldots, x_n, \alpha_1, \ldots, \alpha_n) 
\ \geq \ 
n.
\]
\end{x}
\vspace{0.3cm}

\begin{x}{\small\bf THEOREM} \ 
Roy's conjecture is equivalent to Schanuel's conjecture.
\end{x}
\vspace{0.3cm}

This result is due to Damien Roy.  
While we shall omit the proof, some hints will be given below.
\vspace{0.2cm}

[Note: \ 
Spelled out: \ If Roy's conjecture is true for some $n$ and some choice of $s_0$, $s_1$, $t_0$, $t_1$, $u$ (per \#5), then Schanuel's conjecture is true for this value of $n$.  
Conversely, if Schanuel's conjecture is true for some $n$, then Roy's conjecture is true for the same value of $n$ and any choice of $s_0$, $s_1$, $t_0$, $t_1$, $u$ (per \#5).]
\vspace{0.3cm}

In one direction, assume that the conditions in Roy's conjecture are in force $-$then it can be shown that there exists an integer $K \geq 1$ with the property that
\[
\alpha_j^K
\ = \ 
e^{Kx_j}
\qquad (j = 1, \ldots, n).
\]
Since $x_1, \ldots, x_n$ are $\Q$-linearly independent, the same is true of $K x_1, \ldots, K x_n$, hence
by Schanuel
\[
\trdegQ \Q(Kx_1, \ldots, Kx_n, e\strutv^{Kx_1}, \ldots, e\strutv^{Kx_n}) 
\ \geq \ 
n
\]
or still, 
\[
\trdegQ \Q(Kx_1, \ldots, Kx_n, \alpha_1^K, \ldots, \alpha_n^K) 
\ \geq \ 
n
\]
or still, 
\[
\trdegQ \Q(x_1, \ldots, x_n, \alpha_1, \ldots, \alpha_n) 
\ \geq \ 
n.
\]
Therefore
\[
\text{SCHANUEL $\implies$ ROY}.
\]
\vspace{0.2cm}

In the other direction, take the data as in \#4 and put $\alpha_j = e\strutv^{x_j}$ $(j = 1, \ldots, n)$.  
Take the parameters $s_0$, $s_1$, $t_0$, $t_1$, $u$ as in \#5 and impose the inequalities to be found there.
\end{spacing} 
\vspace{0.3cm}

\begin{x}{\small\bf NOTATION} \ 
Given $R > 0$, let
\[
B(0,R) 
\ = \ 
\{(z_1, z_2) \in \C^2 \ : \ \abs{z_1} \leq R, \ \abs{z_2} \leq R\}
\]
and for any continuous function $F:B(0,R) \ra \C$, put
\[
\abs{F}_R 
\ = \ 
\sup\{\abs{F(z_1,z_2)} \ : \ \abs{z_1} = R, \ \abs{z_2} = R\}.
\]

[Note: \ 
By the maximum modulus principle, when \mF is holomorphic in the interior of $B(0,R)$, $\abs{F}_R$ is the supremum of $\abs{F}$ on $B(0,R)$.]
\end{x}
\vspace{0.3cm}

\begin{x}{\small\bf EXAMPLE} \ 
Let \mL be a positive integer, let $r_0$, $r$, \mR be positive real numbers with $r \geq r_0$, $R \geq 2r$ $-$then
\[
\abs{F}_r 
\ \leq \ 
\sum\limits_{j, k \hsx \geq \hsx 0} \  \frac{1}{j! \hsy k!} \hsx \
\abs{\frac{\partial^{j+k}}{\partial z^j \hsy \partial w^k} \hsx (0,0)} \hsx r^{j + k}
\]
or still, 
\[
\abs{F}_r 
\ \leq \ 
\sum\limits_{j + k \hsx < \hsx L} \hsx \bigg(\frac{r}{r_0}\bigg)^{j+k} \abs{F}_{r_0} + 
\sum\limits_{j + k \hsx \geq \hsx L} \hsx \bigg(\frac{r}{R}\bigg)^{j+k} \hsx \abs{F}_R
\]
or still, 
\[
\abs{F}_r 
\ \leq \ 
\binom{L+1}{2} \hsx \bigg(\frac{r}{r_0}\bigg)^L \hsx \abs{F}_{r_0} + (2L + 4) 
\hsx \bigg(\frac{r}{R}\bigg)^L \hsx \abs{F}_R, 
\]
where
\[
\sum\limits_{j \hsx + \hsx k \hsx \geq \hsx L} \hsx 2^{L - j - k} 
\ = \ 
2L + 4.
\]
\vspace{0.2cm}

[Note: \ 
The conditions on \mF are, of course, the obvious ones \ldots \  .]
\end{x}
\vspace{0.3cm}

\begin{x}{\small\bf LEMMA} \ 
For any sufficiently large postive integer \mN, there exists a nonzero polynomial $P_N \in Z[X_0,X_1]$ with partial degree 
$\leq \ N\strutw^{t_0}$ in $X_0$,  with partial degree $\leq \ N\strutw^{t_1}$ in $X_1$, 
and with height $\leq \ e^N$ such that the function 
\[
f_N(z) 
\ = \ 
P_N(z, e^z)
\]
satisfies
\[
\abs{f_N}_r 
\ \leq \ 
\exp(-2 N\strutw^u).
\]

\vspace{0.2cm}

[Note: \ 
Here
\[
r \ = \ 1 + A N\strutw^{s_1}, 
\]
where
\[
A 
\ = \ 
\abs{x_1} + \cdots + \abs{x_n}.]
\]
\vspace{0.2cm}

The claim now is that
\[
\text{ROY $\implies$ SCHANUEL}.
\]
To verify that this is so, let $k, m_1, \ldots, m_n$ be nonnegative integers, where
\[
k 
\ \leq \ 
N\strutw^{s_0} 
\quad \text{and} \quad
\max\{m_1, \ldots, m_n\} 
\ \leq \ 
N\strutw^{s_1}.
\]
Then
\allowdisplaybreaks
\begin{align*}
\abs{
\big(
\sD^k P_N\big) \hsx 
\bigg( \sum\limits_{j=1}^n \hsx m_j x_j, \ 
\prod\limits_{j=1}^n \hsx \alpha_j^{m_j} 
\bigg)
}\ 
&=\ 
\abs{
\frac{\td^k f_N}{\td z^k} \hsx \bigg( \sum\limits_{j=1}^n \hsx m_j x_j \bigg)
}
\\[12pt]
&\leq\ 
k! \hsx \abs{f_N}_r
\\[12pt]
&\leq\
\exp (-N\strutw^u)
\end{align*}
if \mN is sufficiently large.  Consequently
\[
\trdegQ \Q(x_1, \ldots, x_n, e\strutv^{x_1}, \ldots, e\strutv^{x_n}) 
\ \geq \ 
n,
\]
thus
\[
\text{ROY $\implies$ SCHANUEL}.
\]
as claimed.
\end{x}
\vspace{0.3cm}

\begin{x}{\small\bf \un{N.B.}} \ 
Consider the situation when $n = 1$ $-$then
\allowdisplaybreaks
\begin{align*}
\abs{\big(\sD^k P_N\big) \hsx (m x, \alpha^m)}\ 
&=\ 
\abs{\frac{\td^k f_N}{\td z^k} (mx)}
\\[12pt]
&\leq\
k! \abs{f_N}_{\abs{mx} + 1}.
\end{align*}
Next
\[
\abs{mx} + 1 
\ \leq \ 
\abs{x} N\strutw^{s_1} + 1 
\ = \ 
r
\]
\qquad\qquad $\implies$
\[
\abs{\big(\sD^k P_N\big) \hsx (m x, \alpha^m)}
\ \leq \ 
k! \abs{f_N}_r.
\]
Since $s_0 < u$, it can be assumed that 
\[
N\strutw^{s_0} \hsx \elln(N\strutw^{s_0}) \ \leq \ N\strutw^u,
\]
from which
\allowdisplaybreaks
\begin{align*}
\exp (N\strutw^u) \ 
&\geq\ 
\exp \big( N\strutw^{s_0} \hsx \elln \big(N\strutw^{s_0}\big)\big)
\\[12pt]
&=\ 
\exp \big(\elln\big( \big(N\strutw^{s_0}\big)^{N\strutw^{s_0}}\big)\big)
\\[12pt]
&=\ 
\big(N\strutw^{s_0}\big)^{N\strutw^{s_0}}
\end{align*}
\qquad\qquad $\implies$
\[
k! \ \leq \ k^k \ \leq \ \big(N\strutw^{s_0}\big)^{N\strutw^{s_0}} \ \leq \ \exp (N\strutw^u)
\]
\qquad\qquad $\implies$
\allowdisplaybreaks
\begin{align*}
\abs{\big(\sD^k P_N\big) \hsx (m x, \alpha^m)} \ 
&\leq\ 
\exp (N\strutw^u) \abs{f_N}_r
\\[12pt]
&\leq\ 
\exp (N\strutw^u) \exp (- 2N\strutw^u)
\\[12pt]
&=\ 
\exp ( -N\strutw^u).
\end{align*}
\end{x}
\vspace{0.3cm}

\begin{x}{\small\bf REMARK} \ 
When $n = 1$, Schanuel is an acquired fact: \ If $x \in \C^\times$, then at least one of the two numbers $x$, $e^x$ is transcendental (Hermite-Lindemann), hence
\[
\trdegQ \Q(x, e^x) \ \geq \ 1,
\]
so Roy is automatic in this case.
\end{x}

\newpage
\[
\textbf{APPENDIX}
\]
\vspace{0.3cm}

{\small\bf PRETHEOREM} \ 
Let $(x, \alpha) \in \C \times \C^\times$ and let $s_0$, $s_1$, $t_0$, $t_1$, $u$ be positive real numbers satisfying the inequalities of \#5 $-$then the following conditions are equivalent: \\

\qquad\qquad (i) \quad 
There exists an integer $K \geq 1$ such that $\alpha^K = e^{K x}$.
\\[-.2cm]

\qquad\qquad (ii) \quad 
For any sufficiently large positive integer \mN, there exists a non-zero polynomial $P_N \in \Z[X_0,X_1]$ with partial degree $\leq \ N\strutw^{t_0}$ in $X_0$, with partial degree $\leq \ N\strutw^{t_1}$ in $X_1$, and with height $\leq \ e^N$ which satisfies
\[
\abs{\big(\sD^k P_N \big) \hsx (mx, \alpha^m)} \ \leq \  \exp (-N\strutw^u)
\]
for all nonnegative integers $k$, $m$ with 
\[
k \leq N\strutw^{s_0} \quad \text{and} \quad m \leq N\strutw^{s_1}.
\]
\vspace{0.1cm}

In what follows, we shall sketch the proof that
\[
\text{(ii)} 
\implies 
\text{(i) \qquad or }
\ \neg \  
\text{(i)} 
\implies
\neg \ \text{(ii)}.
\]
Now $\neg$ (i) means that $\forall \ K \in \N$, $\alpha^K \neq e\strutv^{K x}$, hence $\alpha e^{-x}$ is not a root of unity:
\[
\alpha e^{-x} 
\ = \ 
\zeta \qquad (\zeta^K = 1) 
\]
\qquad\qquad $\implies$
\[ 
\alpha^K  = \zeta^K e^{K x} 
\ = \ 
e^{K x} .
\]
\vspace{0.3cm}

{\small\bf OBJECTIVE} \ 
Let $(x, \alpha) \in \C \times \C^\times$ and let $s_0$, $s_1$, $t_0$, $t_1$, $u$ be positive real numbers such that
\[
\max\{1, t_0, 2t_1\} \ < \ \min\{s_0, 2s_1\} \ < \ u.
\]
Suppose that $\alpha e^{-x}$ is not a root of unity $-$then condition (ii) does not hold
for the pair $(x, \alpha)$.
\vspace{0.2cm}

[Note: \ 
The stated assumption on the parameters $s_0$, $s_1$, $t_0$, $t_1$, $u$ is weaker than that of \#5.  
Observe too that there is no restriction from above on $u$.]
\vspace{0.3cm}

{\small\bf NOTATION} \ 
Given $\gamma \in \C - \Q$ and a positive integer $N$, put
\[
\Gamma_\gamma (N) 
\ = \ 
\min\{ \abs{m + n \gamma} : m, n \in \Z, \ 0 < \max\{\abs{m}, \abs{n}\} \ < \ N. 
\]
\vspace{0.3cm}

{\small\bf LEMMA} \ 
For infinitely many \mN,
\[
\Gamma_\gamma(N) 
\ \geq \ 
\frac{1}{2N},
\]
i.e., for infinitely many \mN, 
\[
\abs{m + n \gamma} 
\ \geq \ 
\frac{1}{2N}
\]
for any pair $(m, n) \in \Z^2$ with 
\[
0 \ < \ \max\{\abs{m}, \abs{n}\} \ < \ N.
\]

PROOF \ 
Assume to the contrary that for any integer \mN larger than some $N_0$, there are integers $m(N)$ and $n(N)$ such that
\[
0 \ < \ \max\{\abs{m(N)}, \abs{n(N)}\} \ < \ N
\]
and 
\[
\abs{m(N) + n(N) \hsx \gamma}
\ < \ 
\frac{1}{2N}.
\]
Then $n(N) \neq 0$ and 
\allowdisplaybreaks
\begin{align*}
\abs{m(N) n(N + 1) - m(N + 1) n(N) }\ 
&\leq \
\abs{m(N) + n(N) \gamma} \cdot \abs{n(N+1)} 
\\[12pt]
&\hspace{1.5cm}
+ \abs{m(N+1) + n(N+1)\gamma} \cdot \abs{n(N)}
\\[12pt]
&<\ 1,
\end{align*}
\vspace{0.3cm}
\qquad\qquad $\implies$
\[
m(N) n(N + 1) - m(N + 1) n(N) 
\ = \ 0.
\]
Therefore the ratio
\[
\frac{m(N)}{n(N)}
\]
is a constant $q \in \Q$.  But
\allowdisplaybreaks
\begin{align*}
\abs{q + \gamma} \ 
&=\ 
\abs{m(N) + n(N) \gamma} \big/ \abs{n(N)}
\\[12pt]
&<\ 
\frac{1}{2N}
\end{align*}
for any $N > N_0$, hence $\gamma = -q$, a contradiction.
\vspace{0.5cm}

One can thus attach to each $\gamma \in \C - \Q$ an infinite subset $S_\gamma$ of $\N$, where the elements of $S_\gamma$ 
are the $N$ figuring in the definition of $\Gamma_\gamma (N)$.
\vspace{0.5cm}

{\small\bf \un{N.B.}} \ 
Choose $\lambda$ such that $e^\lambda = \alpha$ $-$then the ratio
\[
\gamma 
\ = \ 
\frac{\lambda - x}{2 \pi \sqrt{-1}} \ \in \C - \Q.
\]
\vspace{0.2cm}

[
Suppose instead that
\[
\frac{\lambda - x}{2 \pi \sqrt{-1}} \ = \ q \quad (\in \Q), 
\]
say $q = \ds\frac{m}{n}$ $(n > 0)$, so
\[
\lambda - x 
\ = \ 
q (2 \pi \sqrt{-1})
\ = \ 
\frac{m}{n} (2 \pi \sqrt{-1})
\]
\qquad\qquad $\implies$
\[
e^{\lambda - x}
\ = \ 
\exp \bigg(\frac{m}{n} \hsx 2 \pi \sqrt{-1}\bigg)
\]
\qquad\qquad $\implies$
\[
\alpha e^{-x} 
\ = \ 
\exp \bigg(\frac{m}{n} \hsx 2 \pi \sqrt{-1}\bigg)
\]
\qquad\qquad $\implies$
\[
\big(\alpha e^{-x}\big)^n
\ = \ 
\exp (m \hsx 2 \pi \sqrt{-1})
\ = \ 
1.]
\]
\vspace{0.3cm}

{\small\bf NOTATION} \ 
Let
\[
\bu 
\ = \ 
(0, 2 \pi \sqrt{-1}), 
\quad
\bv 
\ = \ 
(x, \lambda), 
\quad
\bw 
\ = \ 
(1, 1).
\]
\\[-1cm]

[Note: \ 
\allowdisplaybreaks
\begin{align*}
\bv  - \gamma \bu \ 
&=\ 
(x, \lambda) - \gamma \hsx (0,  2 \pi \sqrt{-1})
\\[12pt]
&=\ 
(x, \lambda) - \frac{\lambda - x}{ 2 \pi \sqrt{-1}} \hsx (0,  2 \pi \sqrt{-1})
\\[12pt]
&=\ 
(x, \lambda) - (\lambda - x) (0,1)
\\[12pt]
&=\ 
(x, \lambda) + (0, x - \lambda)
\\[12pt]
&=\ 
(x, \lambda + x - \lambda)
\\[12pt]
&=\ 
(x,x)
\\[12pt]
&=\ 
x \bw.]
\end{align*}
\vspace{0.3cm}

{\small\bf FACT} \ 
There exists a constant $C \geq 1$ (with $\bu$, $\bv \in B(0,C)$) such that for any $N \in S_\gamma$ and for any pair of real numbers $r, \ R$ with $R \geq 2r$ and $r \geq C N$ and for any continuous function $F: B(0,R) \ra \C$ which is holomorphic in the interior of $B(0,R)$, the estimate
\allowdisplaybreaks
\begin{align*}
\abs{F}_r \ 
&\leq \ 
\bigg(\frac{Cr}{N}\bigg)^{N^2}
\\[12pt]
&\qquad
\times \hsx \max\bigg\{\frac{1}{k!} 
\abs{D_{\bw}^k F(m\bu + n \bv)} N^k : 0 \leq k < N^2, \ 0 \leq m, n < N\bigg\} 
\\[12pt]
&\hspace{3cm}
+ \bigg(\frac{Cr}{R}\bigg)^{N^2} \hsx \abs{F}_R
\end{align*}
obtains.
\vspace{0.2cm}

[Note: \ 
Here
\[
D_{\bw} \ = \ \frac{\partial }{\partial z_1} + \frac{\partial}{\partial z_2} .]
\]
\vspace{0.3cm}

To establish our objective, proceed in steps.
\vspace{0.5cm}

\qquad \un{Step 1:} \quad
Take
\[
\gamma 
\ = \ 
\frac{\lambda - x}{2 \pi \sqrt{-1}} \hsx \in \hsx  \C - \Q.
\]
Then $S_\gamma$ is an infinite subset of $\N$, a generic element $N \in S_\gamma$ being allowed to ``float''.
\vspace{0.5cm}

\qquad \un{Step 2:} \quad
Put
\[
s 
\ = \ 
\min\{s_0 / 2, s_1\}
\]
and let \mM denote the smallest positive integer such that $N \leq M\strutw^s$ (tacitly, $N \in S_\gamma)$.  
Note that \mM depends on \mN (but \mM need not belong to $S_\gamma$) and we shall actually work with \mM rather than \mN in the statement of the objective.
\vspace{0.5cm}

\qquad \un{Step 3:} \quad
Choose a nonzero polynomial $Q_M \in \Z[X_0,X_1]$ with partial degree $\leq \ M\strutw^{t_0}$ in $X_0$, with partial degree 
$\leq \ M\strutw^{t_1}$ in $X_1$, and with height $\leq \ e^M$.
\vspace{0.5cm}

\qquad \un{Step 4:} \quad
Let
\[
\begin{cases}
\ 0 \leq k \leq M\strutw^{s_0} \\[8pt]
\ 0 \leq m \leq M\strutw^{s_1}
\end{cases}
\]
and put
\[
A 
\ = \ 
\max\limits_{k, m} \abs{\big(\sD^k Q_M\big) \hsx (mx, \alpha^m)},
\]
the claim being that if \mN is sufficiently large, then
\[
A 
\ > \ 
\exp (- M^u) \qquad (\exists u \gg 0),
\]
hence for some $k$, for some $m$,
\[
\abs{\big(\sD^k Q_M\big) \hsx (mx, \alpha^m)} 
\ > \ 
\exp (- M\strutw^u),
\]
thereby completing the proof.
\vspace{0.5cm}

\qquad \un{Step 5:} \quad
Define an entire function $G_M : \C^2 \ra \C$ by the prescription
\[
G_M(z, w)
\ = \ 
Q_M (z, e^w).
\]
Let
\[
\partial \ = \ \partial / \partial z + \partial / \partial w.
\]
Then
\[
\partial^k G_M(z,w) 
\ = \ 
\big(\sD^k Q_M \big) \hsx (z, e^w)
\]
for any integer $k \geq 0$ and any $(z,w) \in \C^2$.
\vspace{0.5cm}

\qquad \un{Step 6:} \quad
For any $(n, m) \in \Z^2$, 
\[
\big(\partial^k G_M\big) (n \bu + m \bv) 
\ = \ 
\big(\sD^k Q_M \big) (m x, \alpha^m).
\]
\vspace{0.5cm}

\qquad \un{Step 7:} \quad 
$N^2 \leq M\strutw^{s_0}$, \hsx $N \leq M\strutw^{s_1}$
\vspace{0.2cm}

\qquad $\implies$
\allowdisplaybreaks
\begin{align*}
&\max\bigg\{\frac{1}{k!} \abs{\big(\partial^k G_M\big) (m \bu + n \bv) }  N^k : 0 \leq k < N^2, 0 \leq m, \hsy n < N \bigg\}
\\[12pt]
&\hspace{2cm} \leq \ A \  \sum\limits_{k=0}^\infty \hsx \frac{N\strutw^k}{k!} 
\\[12pt]
&\hspace{2cm} =\ 
A \hsx e^N.
\end{align*}


\qquad \un{Step 8:} \quad
Introduce the constant $C \geq 1$ as above and specialize $r$, $R$ by taking $r = C N$ and $R = e C R$ $-$then in review
\allowdisplaybreaks
\begin{align*}
\abs{F}_r \ 
&\leq \ 
\bigg(\frac{Cr}{N}\bigg)^{N^2}
\\[12pt]
&\qquad
\times \hsx \max\bigg\{\frac{1}{k!} 
\abs{\big(\partial^kF\big) \hsx (m\bu + n \bv)} N^k : 0 \leq k < N^2, \ 0 \leq m, n < N\bigg\} 
\\[12pt]
&\hspace{3cm}
+ \bigg(\frac{Cr}{R}\bigg)^{N^2} \hsx \abs{F}_R
\end{align*}
and in the situation at hand $(F = G_M)$
\\[.2cm]
\[
\begin{cases}
\ \bigg(\ds\frac{Cr}{N}\bigg)^{N^2} = (C \cdot C)^{N^2} = (C^2)^{N^2} = C^{2N^2}\\[15pt]
\ \bigg(\ds\frac{Cr}{R}\bigg)^{N^2} =\bigg(\ds\frac{Cr}{eCr}\bigg)^{N^2} = e\strutv^{-N^2} 
\end{cases}
,
\]
so
\[
\abs{G_M}_r 
\ \leq \ 
C^{2N^2} \hsx A e^N + e^{-N^2} \abs{G_M}_R.
\]
\vspace{0.5cm}

\qquad \un{Step 9:} \quad
Since $\max\{1, t_0, s + t_1\} < 2s$, the definitions imply that
\allowdisplaybreaks
\begin{align*}
\abs{G_M}_R 
&\leq \ 
(M\strutw^{t_0} + 1) \hsx (M\strutw^{t_1} + 1) \times \exp (M + M\strutw^{t_0} \elln(R) + R M\strutw^{t_1} )
\\[12pt]
&\leq \
 e^{N^2} / 2
\end{align*}
provided $N$ is sufficiently large.
\vspace{0.5cm}

\qquad \un{Step 10:} \quad
$Q_M$ is a nonzero polynomial with integral coefficients, hence
\[
1 
\ \leq \ 
H(Q_M)
\ \leq \ 
\abs{Q_M}_1
\ \leq \ 
\abs{G_M}_\pi
\ \leq \ 
\abs{G_M}_r
\]
if $r \leq \pi$.
\vspace{0.5cm}

\qquad \un{Step 11:} \quad
Explicate the relation
\[
\abs{G_M}_r 
\ \leq \ 
C^{2N^2} A e^{N} + e^{-N^2} \hsx \abs{G_M}_R
\]
to arrive at
\[
1 
\ \leq \ 
\abs{G_M}_{C  N}
\ \leq \ 
C^{2N^2} A e^{N} + e^{-N^2} \hsx \big(e^{N^2}/2\big)
\]
for \mN large enough.\\[12pt]

I.e.: 
\[
1 \leq \ C^{2N^2} e^{N} A + \frac{1}{2}.
\]

I.e.: 
\[
\frac{1}{2} \ \leq \ C^{2N^2} e^{N} A.
\]

I.e.: 
\[
A \ \geq \ \frac{1}{2} C^{-2N^2} e^{-N}.
\]
\vspace{0.5cm}

\qquad \un{Step 12:} \quad
Apart from the restriction that
\[
\min\{s_0, 2s_1\} \ < \ u,
\]
the parameter $u \gg 0$ is at our disposal and can be chosen as large as we please.  
Bearing in mind that $2s$ is $<  u$, or now, as will be notationally convenient, $2s$ is $< v$, write
\allowdisplaybreaks
\begin{align*}
N \ \leq \  M\strutw^s 
&\implies 
N^2 \leq M\strutw^{2s} < \ M\strutw^v
\\[12pt]
&\implies
e^{N^2} < \exp (M\strutw^v).
\end{align*}
Consequently for some $u > v \gg 0$, 
\[
A \ \geq \ \frac{1}{2} C^{-2N^2} e^{-N} \ > \ \exp(-M\strutw^u).
\]
\vspace{0.2cm}

[To see this, ignore the $\ds\frac{1}{2}$ and for simplicity take $C = e$ $-$then
\allowdisplaybreaks
\begin{align*}
N^2 \ < \ M\strutw^v \ \implies 2 N^2
&< \ 
2M\strutw^v
\\[12pt]
&<\ 
2^v M\strutw^v 
\\[12pt]
&=\ 
(2M)\strutw^v
\\[12pt]
&=\ 
M^w.
\end{align*}
Here
\[
w \ = \ v \hsx \frac{\elln (2M)}{\elln (M)} \ > \ v.
\]
In fact,
\allowdisplaybreaks
\begin{align*}
(2M)^v \ = \ M\strutw^w 
&\implies
\elln((2M)^v) \ = \ \elln (M\strutw^w)
\\[12pt]
&\implies
v \elln(2M) \ = \ w \elln(M).
\end{align*}
Therefore
\allowdisplaybreaks
\begin{align*}
e^{2N^2} \hsx e^N \ 
&<\ 
\exp(M\strutw^w ) \hsx \exp(M\strutw^v)
\\[12pt]
&=\
\exp(M\strutw^w + M\strutw^v)
\\[12pt]
&<\ 
\exp (2M\strutw^w)
\\[12pt]
&<\ 
\exp (2^w M\strutw^w)
\\[12pt]
&=\
\exp ((2M)\strutw^w)
\\[12pt]
&=\
\exp(M\strutw^u)
\end{align*}
if 
\[
u 
\ = \ 
w \hsx \frac{\elln (2M)}{\elln (M)} 
\ > \ 
w \ (> v).
\]
Accordingly
\[
e^{-2N^2} \hsx e^{-N} 
\ > \ 
\exp (-M\strutw^u).]
\]

%% file: _62_real_numbers_bis.tex
\chapter{
$\boldsymbol{\S}$\textbf{62}.\quad  REAL NUMBERS (bis)}
\setlength\parindent{2em}
\setcounter{theoremn}{0}
\renewcommand{\thepage}{\S62-\arabic{page}}

\ \indent 
``Few mathematical structures have undergone as many revisions or have been presented in as many guises as the real numbers.  
Every generation re-examines the reals in the light of its \ldots mathematical objectives.''
\vspace{0.2cm}

[F. Faltin et al., Advances in Mathematics 16 (1975), p. 278.]
\vspace{0.3cm}

\[
*\hsx *\hsx*\hsx*\hsx*\hsx*\hsx*\hsx*\hsx*\hsx*\hsx*
\]
\\[-1cm]

``How do we get future generations to take the validity of real numbers for granted?  
We indoctrinate them early in their careers when they are eager but impressionable undergraduates.  
Here's how we do it.  
First we soften them up with a ``Constructing the Real Numbers'' blurb in their first calculus course.  
Needless to say we don't really construct the real numbers as they are by definition unconstructible.  
But the phrase sticks in their minds long after the details are forgotten.''
\vspace{0.2cm}

[N. J. Wildberger, The Mathematical Intelligencer 21 (1999), pp. 4-7.]
\vspace{0.2cm}

\[
*\hsx *\hsx*\hsx*\hsx*\hsx*\hsx*\hsx*\hsx*\hsx*\hsx*
\]

``How real are the real numbers? \ldots 
The frightening features are the unsolvability of the halting problem (Turing, 1936), the fact that most reals are uncomputable, and last but not least, the halting probability $\Omega$, which is irreducibly complex (algorithmically random), maximally unknowable, and dramatically illustrates the limits of reason.''
\vspace{0.2cm}

[Gregory Chaitin, arXiv:math/0411418 v 3 [math.HO] 29 Nov 2004.] 

\[
*\hsx *\hsx*\hsx*\hsx*\hsx*\hsx*\hsx*\hsx*\hsx*\hsx*
\]
\vspace{0.5cm}

For a systematic survey of the various constructions which lead to the real numbers, consult 
\vspace{0.2cm}

[Ittay Weiss, arXiv:math/1506.03467 v1 [math.HO] 18 May 2015.]
\[
*\hsx *\hsx*\hsx*\hsx*\hsx*\hsx*\hsx*\hsx*\hsx*\hsx*
\]

\vspace{0.5cm}
\[
\textbf{APPENDIX}
\]
\vspace{0.5cm}

In algorithmic information theory, a 
\un{halting probability}
\index{halting constant} 
(or 
\un{Chaitin constant})
\index{Chaitin constant} 
is a real number $\Omega$ which represents the probability that a randomly constructed program will halt.

To be precise, let $P_F$ be the domain of a prefix-free universal computable function \mF $-$then the halting probability 
$\Omega_F$ of $P_F$ is by definition
\[
\Omega_F 
\ = \ 
\sum\limits_{p \in P_F} \hsx 2^{-\abs{p}},
\]
where $\abs{p}$ denotes the length of a binary string $p$.  
The sum defining $\Omega_F$ is infinite and converges to a real number lying between 0 and 1.
\vspace{0.75cm}

{\small\bf FACT} \ 
$\Omega_F$ is transcendental.
\vspace{0.75cm}

There is a probabilistic interpretation of $\Omega_F$, from which the terminology.  
Thus let $(X, \mu)$ be the Cantor space and suppose that \mF is a prefix-free universl computable function $-$then the domain $P_F$ of \mF consists of an infinite set of binary strings: 
\[
P_F 
\ = \ 
\{p_1, p_2, \ldots\}.
\]
Each of these strings $p_i$ determines a subset $S_i$ of the Cantor space (viz. all sequences in Cantor space that begin with $p_i$).
Moreover the $S_i$ are pairwise disjoint and 
\[
\Omega_F 
\ = \ 
\mu \hsx \bigg( \hsx \bigcup\limits_{i \in \N} \hsx S_i \bigg).
\]
\vspace{0.3cm}

{\small\bf REMARK} \ 
$\Omega_F$ is not computable, i.e., there is no algorithm which, given 
$n$, returns the first $n$ digits of $\Omega_F$.
\vspace{0.3cm}

For more information on this material, consult George Barmpalias 
(arXiv:1707.08109 v 3 [math.LO]).

%% file: _XX_Supp_trans_of_series.tex
\centerline{\textbf{\large SUPPLEMENT I}}
\vspace{0.75cm}
\renewcommand{\thepage}{Transcendence of Series - 0}

\centerline{\textbf{TRANSCENDENCE OF SERIES}}
\vspace{0.5cm}

The overall theme is to discuss the transcendence of numbers of the form
\\[2pt]
\[
\sum\limits_{n\hsx=\hsx 1}^\infty \  \frac{A(n)}{B(n)} 
\quad \bigg(\text{or} \quad 
\sum\limits_{n\hsx =\hsx 0}^\infty \  \frac{A(n)}{B(n)}  \bigg)
\]
or
\[
\sum\limits_{n\hsx =\hsx -\infty}^\infty \  \frac{A(n)}{B(n)} 
\ \equiv \ 
\lim\limits_{N \ra \infty}
\sum\limits_{\abs{n} \hsx < \hsx N} \  \frac{A(n)}{B(n)}.
\]
\vspace{.2cm}

The literature on this subject is extensive and no attempt will be made at a systematic exposition.  
Foregoing this, we shall first examine a number of instructive special cases and then take a look at the general picture.
\vspace{.2cm}

[Note: \ 
Omitted details are to be regarded as excercises ad libitum.]
\vspace{.2cm}

\allowdisplaybreaks
\begin{align*}
\S1.  \qquad &\text{CANONICAL ILLUSTRATIONS} 
\\[9pt]
\S2.  \qquad &\text{THE ROLE OF THE COTANGENT}  
\\[9pt]
\S3.  \qquad &\text{APPLICATION OF NESTERENKO}
\\[9pt]
\S4.  \qquad &\text{INTRODUCTION OF SCHC}
\\[9pt]
\S5.  \qquad &\text{INTRODUCTION OF SCHC (bis)}
\\[9pt]
\S6.  \qquad &\text{CONSOLIDATION}
\\[9pt]
\S7.  \qquad &\text{CONSIDERATION OF $\frac{A}{B}$}
\\[9pt]
\S8.  \qquad &\text{AN ALGEBRAIC SERIES}
\end{align*}

%% file: _01_canonical_illustrations.tex
\chapter{
$\boldsymbol{\S}$\textbf{1}.\quad  CANONICAL ILLUSTRATIONS}
\setlength\parindent{2em}
\setcounter{theoremn}{0}
\renewcommand{\thepage}{Transcendence of Series \S1-\arabic{page}}


\begin{x}{\small\bf EXAMPLE} \ 
\[
\sum\limits_{n=1}^\infty \hsx 
\frac{1}{n (n+1)}
\ = \ 
\sum\limits_{n=1}^\infty \hsx 
\bigg(\frac{1}{n} - \frac{1}{n+1}\bigg) 
\ = \ 
1.
\]
\end{x}
\vspace{0.3cm}

\begin{x}{\small\bf EXAMPLE} \ 
\[
\sum\limits_{n=0}^\infty \hsx 
\frac{1}{n!}
\ = \ e 
\quad \text{and} \quad
\sum\limits_{n=1}^\infty \hsx 
\frac{(-1)^{n-1}}{2n - 1}
\ = \ 
\frac{\pi}{4},
\]
both of which are transcendental.
\end{x}
\vspace{0.3cm}

\begin{x}{\small\bf EXAMPLE} \ 
\[
\sum\limits_{n=1}^\infty \hsx 
\frac{(-1)^{n-1}}{n}
\ = \ 
\elln(2),
\]
a transcendental number (cf. \S21, \#9).
\end{x}
\vspace{0.3cm}

\begin{x}{\small\bf EXAMPLE} \ 
\[
\sum\limits_{n=1}^\infty \hsx 
\frac{1}{n^3}
\ = \ 
\zeta(3),
\]
an irrational number, the transcendence of which has yet to be shown.
\end{x}
\vspace{0.3cm}

\begin{x}{\small\bf EXAMPLE} \ 
\[
\sum\limits_{n=1}^\infty \hsx 
\bigg(\frac{1}{n} - \elln \bigg(1 + \frac{1}{n}\bigg) \bigg) 
\ = \ 
\gamma,
\]
\end{x}
$\gamma$ being Euler's constant, which is not known to be irrational, let alone transcendental.
\vspace{0.3cm}

\begin{x}{\small\bf EXAMPLE} \ 
\[
\sum\limits_{n=0}^\infty \  
\frac{(-1)^n}{(2n + 1)^2}
\ = \ 
G,
\]
\mG being Catalan's constant, whose irrationality status is unknown.
\vspace{0.2cm}

[Note: \ 
By comparison, 
\[
\sum\limits_{n=0}^\infty \  
\frac{(-1)^n}{(2n + 1)^3}
\ = \ 
\frac{\pi^3}{32}.]
\]
\end{x}
\vspace{0.3cm}

\begin{x}{\small\bf LEMMA} \ 
The zeros of the polynomial $X^2 - X - 1$ are 
$\phi = \ds\frac{1 + \sqrt{5}}{2}$ 
(the 
\un{golden ratio})
\index{golden ratio} 
and 
$\psi = \ds\frac{1 - \sqrt{5}}{2}$ 
$(= 1 - \phi = -\ds\frac{1}{\phi})$.
\vspace{0.3cm}

[Note: \ 
$\phi$ and $\psi$ are quadratic irrationals (cf. \S8, \#4).]
\end{x}
\vspace{0.3cm}

\begin{x}{\small\bf EXAMPLE} \ 
\[
\sum\limits_{n=-\infty}^\infty \hsx 
\frac{2n - 1}{n^2- n - 1}
\ = \ 
\sum\limits_{n=-\infty}^\infty \hsx 
\bigg(
\frac{1}{n - \phi} + \frac{1}{n - \psi}
\bigg)
\ = \ 
0.
\]
\end{x}
\vspace{0.3cm}

\begin{x}{\small\bf DEFINITION} \ 
The integers 0, 1, 1, 2, 3, 5, 8, 13, 21, 34, 55, 89, 144, \ldots are the 
\un{Fibonacci numbers}:
\index{Fibonacci numbers}
\[
F_0 = 0, \quad 
F_1 = 1, \quad 
F_n = F_{n-1} + F_{n-2} 
\qquad (n \geq 2).
\]
\end{x}
\vspace{0.3cm}

\begin{x}{\small\bf LEMMA} \ 
\[
F_n 
\ = \ 
\frac{\phi^n - \psi^n}{\phi - \psi}
\ = \ 
\frac{\phi^n - \psi^n}{\sqrt{5}}.
\]
\end{x}
\vspace{0.3cm}

\begin{x}{\small\bf \un{N.B.}} \ 
$\phi$ and $\psi$ are both solutions to the equations
\[
X^n \ = \ X^{n-1} + X^{n-2},
\]
hence
\[
\begin{cases}
\ \phi^n = \phi^{n-1} + \phi^{n-2}\\[8pt]
\ \psi^n = \psi^{n-1} + \psi^{n-2}
\end{cases}
.
\]
\end{x}
\vspace{0.3cm}

\begin{x}{\small\bf EXAMPLE} \ 
\allowdisplaybreaks
\begin{align*}
\sum\limits_{n=1}^\infty \  \frac{F_n}{n 2^n} \
&=\ 
\frac{1}{\sqrt{5}} \hsx\hsx 
\sum\limits_{n=1}^\infty \  
\frac{1}{n} \hsx 
\bigg(
\bigg(
\frac{\phi}{2}\bigg)^n - \bigg( - \frac{1}{2 \phi}\bigg)^n\bigg)
\\[12pt]
&=\ 
\frac{1}{\sqrt{5}} \hsx\hsx \elln(1 + \phi) - \frac{1}{\sqrt{5}} \hsx\hsx \elln(2 - \phi),
\end{align*}
a transcendental number (cf. \S31, \#11).
\end{x}
\vspace{0.3cm}

\begin{x}{\small\bf EXAMPLE} \ 
\[
\sum\limits_{n=1}^\infty \hsx 
\frac{1}{F_n F_{n+2}}
\ = \ 
1.
\]
\end{x}
\vspace{0.3cm}

\begin{x}{\small\bf EXAMPLE} \ 
\[
\sum\limits_{n=1}^\infty \  
\frac{(-1)^n}{F_n F_{n+1}}
\ = \ 
\frac{1 - \sqrt{5}}{2} 
\qquad (= \psi).
\]
\end{x}
\vspace{0.3cm}

\begin{x}{\small\bf LEMMA} \ 
If $\alpha_1, \ldots, \alpha_n$ are positive algebraic numbers and if $\beta_0, \beta_1, \ldots, \beta_n$ are algebraic numbers with $\beta_0 \neq 0$, then
\[
\beta_0 \pi \hsx + \hsx
\sum\limits_{j=1}^n \hsx \beta_j \elln(\alpha_j)
\]
is a transcendental number.
\vspace{0.2cm}

PROOF \ 
Replace $\pi$ by $-\sqrt{-1} \hsx \Log (-1)$ and quote \S31, \#11.
\vspace{0.2cm}

[The underlying supposition is that
\[
\beta_0 \hsx \pi \hsx + \hsx 
\sum\limits_{j=1}^n \hsx \beta_j \hsx \elln(\alpha_j)
\]
is nonzero.  
To see this, let $\{\elln(\alpha_j) : j \in S\}$ be a maximal $\Q$-linearly independent subset of
\[
\elln(\alpha_1), \ldots, \elln(\alpha_n),
\]
hence
\[
\beta_0 \hsx \pi \hsx + \hsx 
\sum\limits_{j=1}^n \hsx \beta_j \elln(\alpha_j)
\ = \ 
-\sqrt{-1} \hsx \beta_0 \hsx \Log (-1) \hsx + \hsx 
\sum\limits_{j \in S} \hsx C_j \elln(\alpha_j)
\]
for algebraic numbers $C_j$.  The claim now is that
\[
\Log (-1), \quad \elln(\alpha_j) \qquad (j \in S)
\]
are linearly independent over $\Q$, thus are linearly independent over $\Qbar$ (homogeneous Baker), thereby implying that
\[
-\sqrt{-1} \hsx \beta_0 \Log (-1) \hsx + \hsx 
\sum\limits_{j \in S} \hsx C_j \elln(\alpha_j)
\]
is nonzero.  So consider a rational dependence relation
\[
q_0 \hsx \Log (-1) \hsx + \hsx \sum\limits_{j \in S} \hsx q_j \elln(\alpha_j) \ = \ 0.
\]
The sum over $j \in S$ is a real number, while $\Log(-1)$ is pure imaginary, which forces $q_0 = 0$.  
But then $q_j = 0$ $\forall \ j \in S$.]
\end{x}
\vspace{0.3cm}

\begin{x}{\small\bf EXAMPLE} \ 
(Lehmer) \ 
\[
\sum\limits_{n=0}^\infty \ \hsx
\prod\limits_{j=1}^6 \  
\frac{1}{6n + j} 
\ = \ 
\frac{1}{4320} \hsx \big(192 \hsx \elln(2) - 81 \hsx \elln(3) + 7 \sqrt{3} (-\pi)\big),
\]
a transcendental number.
\end{x}
\vspace{0.3cm}


%% file: _02_the_role_of_the_cotangent.tex
\chapter{
$\boldsymbol{\S}$\textbf{2}.\quad  THE ROLE OF THE COTANGENT}
\setlength\parindent{2em}
\setcounter{theoremn}{0}
\renewcommand{\thepage}{Transcendence of Series \S2-\arabic{page}}


\begin{x}{\small\bf RAPPEL} \ 
$\forall \ z \in \C - \Z$, 
\[
\pi \hsx \cot (\pi z) \ 
\ = \ 
\sum\limits_{n\hsx = \hsx -\infty}^\infty \hsx \frac{1}{n + z}.
\]
\end{x}
\vspace{0.1cm}

\begin{x}{\small\bf THEOREM} \ 
Let $C \in \Q - \Z$ $-$then the series
\[
\sum\limits_{n=-\infty}^\infty \hsx \frac{1}{n + C}
\]
is transcendental if $C \not\equiv \ds\frac{1}{2} \hsx \modx \Z$.
\vspace{0.3cm}

PROOF \ 
Write
\allowdisplaybreaks
\begin{align*}
\pi \hsx \cot (\pi C) \ 
&=\ 
\pi \hsx \sqrt{-1} \  
\frac{e^{\pi \sqrt{-1} \hsx C}  + e^{-\pi \sqrt{-1} \hsx C}}{e^{\pi \sqrt{-1} \hsx C}  - e^{-\pi \sqrt{-1} \hsx C}}
\\[12pt]
&=\ 
\pi \hsx \sqrt{-1} \  
\frac{e^{2 \pi \sqrt{-1} \hsx C}  +1}{{e^{2 \pi \sqrt{-1} \hsx C}  - 1} }
\\[12pt]
&\neq 0.
\end{align*}
Let $C = \ds\frac{p}{q}$:
\[
\implies \ 
e^ {2 \pi \sqrt{-1} \hsx C} 
\ = \ 
\big(e^ {2 \pi \sqrt{-1} \hsx / q} \big)^p \in \Qbar.
\]
Therefore
\[
\sum\limits_{n = -\infty}^\infty \hsx \frac{1}{n + C}
\]
is transcendental (being $\pi$ times a nonzero algebraic number).
\vspace{0.2cm}

[Note: \ 
If $C \equiv \ds\frac{1}{2} \hsx \modx \Z$, then the series vanishes.  
In fact, $\forall \ m \in \Z$, 
\[
e^{2 \pi \sqrt{-1} \hsx \bigl(\frac{1}{2} + m\bigr)} 
\ = \ 
e^{\pi \sqrt{-1}}
\ = \ -1.
\]
One can also argue directly without an appeal to the formula: \ $\forall \ m \in \Z$, 
\allowdisplaybreaks
\begin{align*}
\sum\limits_{n\hsx =\hsx -\infty}^\infty \hsx \frac{1}{n + \frac{1}{2} + m} \
&=\ 
\sum\limits_{n\hsx =\hsx -\infty}^\infty \hsx \frac{1}{n - 1 - 2m + \frac{1}{2} + m}
\\[12pt]
&=\ 
\sum\limits_{n\hsx =\hsx -\infty}^\infty \hsx \frac{1}{n - \frac{1}{2} - m}
\\[12pt]
&=\ 
\sum\limits_{n\hsx =\hsx -\infty}^\infty \hsx \frac{1}{-n - \frac{1}{2} - m}
\\[12pt]
&=\ 
-\sum\limits_{n\hsx =\hsx -\infty}^\infty \hsx \frac{1}{n + \frac{1}{2} + m}. \hsx]
\end{align*}
\end{x}
\vspace{0.3cm}

\begin{x}{\small\bf LEMMA} \ 
$\forall \ k \geq 2$, $\forall \ z \in \C - \Z$,
\[
\frac{\td^{k-1}}{\td z^{k-1}} \hsx 
\bigg(\sum\limits_{n=-\infty}^\infty \hsx \frac{1}{n + z}\bigg) 
\ = \ 
(-1)^{k-1} \hsx (k - 1)! \hsx 
\sum\limits_{n=-\infty}^\infty \hsx \frac{1}{(n +z)^k}.
\]
\vspace{0.01cm}

\noindent Therefore
\[
\sum\limits_{n=-\infty}^\infty \hsx \frac{1}{(n +z)^k}
\ = \ 
\frac{(-1)^{k-1} \hsx (\pi \cot (\pi z))^{(k-1)}}{(k - 1)!}.  
\]
\end{x}
\vspace{0.3cm}

\begin{x}{\small\bf LEMMA} \ 
$\forall \ k \geq 2$, $\forall \ z \in \C - \Z$,
\[
\frac{\td^{k - 1}}{\td z^{k - 1}} (\pi \cot(\pi z))
\ = \ 
(2 \pi \hsx \sqrt{-1})^{k} \hsx 
\bigg(
\frac{A_{k,1}}{e^{2 \pi \sqrt{-1} \hsx z}  - 1} 
+ \cdots +
\frac{A_{k,k}}{(e^{2 \pi \sqrt{-1} \hsx z}  - 1)^k}
\bigg),
\]
where $A_{i,j} \in \Z$ and $A_{k,1} \neq 0$, $A_{k,k} \neq 0$.
\vspace{0.2cm}

PROOF \ 
Write
\allowdisplaybreaks
\begin{align*}
\pi \hsx \cot (\pi z) \ 
&=\ 
\pi \hsx \sqrt{-1} \hsx \hsx
\frac{e^{2 \pi \sqrt{-1} \hsx z}  + 1}{e^{2 \pi \sqrt{-1} \hsx z}  - 1} 
\\[12pt]
&=\ 
\pi \hsx \sqrt{-1} \hsx \hsx
\frac{e^{2 \pi \sqrt{-1} \hsx z}  - 1 + 1 + 1}{e^{2 \pi \sqrt{-1} \hsx z}  - 1} 
\\[12pt]
&=\ 
\pi \hsx \sqrt{-1} \hsx \hsx
\bigg(
\frac{e^{2 \pi \sqrt{-1} \hsx z}  - 1}{e^{2 \pi \sqrt{-1} \hsx z}  - 1} + 
\frac{2}{e^{2 \pi \sqrt{-1} \hsx z}  - 1}
\bigg)
\\[12pt]
&=\ 
\pi \hsx \sqrt{-1} \hsx \hsx
\bigg(1 + \frac{2}{e^{2 \pi \sqrt{-1} \hsx z}  - 1}\bigg).
\end{align*}
Differentiating this gives the result for $k= 2$.  
Proceeding by induction, assume matters have been established at level $\ell - 1$, hence
\[
A_{\ell-1,1}, \ldots, A_{\ell-1,\ell-1} \in \Z
\]
with $A_{\ell-1, 1} \neq 0$, $A_{\ell-1, \ell-1} \neq 0$ and 
\[
\frac{\td}{\td z} \bigg(\frac{\td^{\ell - 2}}{\td z^{\ell - 2}} (\pi \cot(\pi z) \bigg)
\ = \ 
(2 \pi \hsx \sqrt{-1})^{\ell - 1} \hsx 
\frac{\td}{\td z}
\Bigg(
\frac{A_{\ell-1,1}}{e^{2 \pi \sqrt{-1} \hsx z}  - 1} + \cdots +
\frac{A_{\ell-1,\ell-1}}{(e^{2 \pi \sqrt{-1} \hsx z}  - 1)^{\ell-1}}
\Bigg)
\]
\\[-5pt]
or still, 
\[
(2 \pi \hsx \sqrt{-1})^{\ell} \hsx 
\Bigg(
-A_{\ell-1,1} \hsx \frac{e^{2 \pi \sqrt{-1} \hsx z}}{(e^{2 \pi \sqrt{-1} \hsx z}  - 1)^2} 
- \cdots - 
(\ell - 1) A_{\ell-1,\ell-1} \hsx \frac{e^{2 \pi \sqrt{-1} \hsx z}}{\big(e^{2 \pi \sqrt{-1} \hsx z}  - 1\big)^{\ell}}
\Bigg)
\]
or still, 
\[
(2 \pi \hsx \sqrt{-1})^{\ell} \hsx 
\Bigg(
-A_{\ell-1,1} \hsx \frac{e^{2 \pi \sqrt{-1} \hsx z} - 1 + 1}{(e^{2 \pi \sqrt{-1} \hsx z}  - 1)^2} 
- \cdots - 
(\ell - 1) A_{\ell-1,\ell-1} \hsx \frac{e^{2 \pi \sqrt{-1} \hsx z} - 1 + 1}{\big(e^{2 \pi \sqrt{-1} \hsx z}  - 1\big)^{\ell}}
\Bigg),
\]
which equals $(2 \pi \hsx \sqrt{-1})^{\ell}$ times
\[
- 
\frac{A_{\ell-1,1}}{e^{2 \pi \sqrt{-1} \hsx z}  - 1} 
- 
\frac{A_{\ell-1,1}}{\big(e^{2 \pi \sqrt{-1} \hsx z}  - 1\big)^2} 
- \cdots - 
\frac{(\ell - 1) \hsx A_{\ell-1,\ell-1}}{\big(e^{2 \pi \sqrt{-1} \hsx z}  - 1\big)^{\ell-1}} 
-
\frac{(\ell - 1) \hsx A_{\ell-1,\ell-1}}{\big(e^{2 \pi \sqrt{-1} \hsx z}  - 1\big)^\ell},
\]
thereby leading to the result at level $\ell$.

\vspace{0.2cm}

[Note: \ 
To see the pattern, take $\ell = 3$ and put $w = e^{2 \pi \sqrt{-1} \hsx z} - 1$ $-$then
\[
\frac{A_{2,1}}{w} + \frac{A_{2,1}}{w^2} + \frac{2 A_{2,2}}{w^2} + \frac{2 A_{2,2}}{w^3} 
\ = \ 
\frac{A_{2,1}}{w} + \frac{A_{2,1} + 2 A_{2,2}}{w^2} + \frac{2 A_{2,2}}{w^3} \hsx .]
\]
Therefore
\[
\sum\limits_{n=-\infty}^\infty \hsx \frac{1}{(n+z)^k} 
\ = \ 
\frac{(-1)^{k-1}}{(k-1)!} \hsx \bigl(2\pi  \sqrt{-1}\bigr)^k  \hsx
\times \hsx 
\bigg(
\frac{A_{k,1}}{e^{2 \pi \sqrt{-1} \hsx z} - 1} + \cdots + 
\frac{A_{k,k}}{(e^{2 \pi \sqrt{-1} \hsx z} - 1)^k}
\bigg).
\]
\end{x}
\vspace{0.3cm}

\begin{x}{\small\bf NOTATION} \ 
Put
\[
A_k(z) 
\ = \ 
\frac{(-1)^{k-1}}{(k-1)!} \hsx \bigl(2 \pi \sqrt{-1}\bigr)^k  \hsx 
\bigg(
\frac{A_{k,1}}{e^{2 \pi \sqrt{-1} \hsx z} - 1} + \cdots + 
\frac{A_{k,k}}{(e^{2 \pi \sqrt{-1} \hsx z} - 1)^k}
\bigg).
\]
\\[-11pt]
Therefore
\[
\sum\limits_{n=-\infty}^\infty \hsx \frac{1}{(n+z)^k} \ = \ \pi^k \hsx A_k(z).
\]
\end{x}
\vspace{0.3cm}

\begin{x}{\small\bf \un{N.B.}} \ 
$\forall \ C \in \Q - \Z$, $A_k(C)$ is an algebraic number.
\end{x}
\vspace{0.3cm}

\begin{x}{\small\bf THEOREM} \ 
$\forall \ k \geq 2$, $\forall \ C \in \Q - \Z$, the series
\[
\sum\limits_{n=-\infty}^\infty \hsx \frac{1}{(n+C)^k}
\]
is either transcendental or zero.
\end{x}
\vspace{0.3cm}

\begin{x}{\small\bf REMARK} \ 
It can happen that 
\[
\frac{\td^{k-1}}{\td z^{k-1}} (\pi \cot(\pi z))  \Big|_{z = C} \quad (k \geq 2) \ = \ 0.
\]
To see this, take $k$ odd and observe that $\forall \ m \in \Z$, 
\allowdisplaybreaks
\begin{align*}
\sum\limits_{n=-\infty}^\infty \hsx \frac{1}{(n + \frac{1}{2} + m)^k}\ 
&=\ 
\sum\limits_{n=-\infty}^\infty \hsx \frac{1}{(- n - \frac{1}{2} - m)^k}
\\[12pt]
&=\ 
(-1)^k \hsx 
\sum\limits_{n=-\infty}^\infty \hsx \frac{1}{(n+\frac{1}{2} + m)^k}.
\end{align*}
\vspace{0.2cm}

[Note: \ 
The series does not vanish if $k$ is even and in that case we have transcendence.]
\end{x}
\vspace{0.3cm}


%% file: _03_application_of_nesterenko.tex
\chapter{
$\boldsymbol{\S}$\textbf{3}.\quad  APPLICATION OF NESTERENKO}
\setlength\parindent{2em}
\setcounter{theoremn}{0}
\renewcommand{\thepage}{Transcendence of Series \S3-\arabic{page}}


\begin{x}{\small\bf CRITERION} \ 
For any positive integer \mD, $\pi$ and $e^{\pi \sqrt{D}}$ are algebraically independent over $\Q$ (cf. \S20, \#10) (proof omitted).
\vspace{0.2cm}

[Note: \ 
In particular, $\pi$ and $e^{\pi}$ are algebraically independent over $\Q$.]
\end{x}
\vspace{0.3cm}

\begin{x}{\small\bf \un{N.B.}} \ 
If $r$ and $s$ are nonzero rational numbers, then $\pi^r$ and $\big(e^{\pi \sqrt{D}}\big)^s$ are algebraically independent over $\Q$ 
(cf. \S46, \#26).
\end{x}
\vspace{0.3cm}

\begin{x}{\small\bf THEOREM} \ 
Let $C \in \Q - \{0\}$ $-$then the series
\[
\sum\limits_{n=-\infty}^\infty \hsx 
\frac{1}{n^2+ C^2}
\] 
is transcendental.
\vspace{0.2cm}

PROOF \ 
Take $C > 0$ and let
\[
f(x) 
\ = \ 
\frac{\pi}{C} \hsx 
e^{- 2 \pi C \abs{x}}.
\] 
Then, using Poisson summation,
\[
\sum\limits_{n\hsx =\hsx -\infty}^\infty \hsx 
f(n + t)
\ = \ 
\sum\limits_{n\hsx =\hsx -\infty}^\infty \hsx 
\widehat{f} (n) \hsx e^{2 \pi \sqrt{-1} \hsx t \hsx n}.
\] 
Now put $t = 0$ to get
\[
\frac{\pi}{C} \hsx 
\sum\limits_{n\hsx =\hsx -\infty}^\infty \hsx 
e^{- 2 \pi C \abs{n}}
\ = \ 
\sum\limits_{n\hsx =\hsx -\infty}^\infty \hsx 
\frac{1}{n^2+ C^2}
\] 
or still, 
\[
\sum\limits_{n \hsx = \hsx -\infty}^\infty \hsx 
\frac{1}{n^2+ C^2}
\ = \ 
\frac{\pi}{C} \hsx 
\bigg(\frac{e^{2 \pi C} + 1}{e^{2 \pi C} - 1}\bigg),
\] 
a transcendental number (cf. infra).
\vspace{0.2cm}

\begin{spacing}{1.55}
[Note: \ 
Let
\[
C 
\ = \ 
\frac{p}{q} \qquad (p, q \in \N)
\] 
and write
\[
2 C 
\ = \ 
2 \hsx \frac{p}{q}
\ = \ 
\frac{\sqrt{4 p^2}}{q}
\ \equiv \ 
\frac{\sqrt{D}}{q}.
\] 
If
\[
\frac{\pi}{C} \hsx 
\bigg(\frac{e^{2 \pi C} + 1}{e^{2 \pi C} - 1}\bigg)
\ = \ 
\alpha \in \Qbar - \{0\},
\] 
then
\[
\frac{\pi}{C} \hsx  (e^{2 \pi C} + 1) - \alpha (e^{2 \pi C} - 1) 
\ = \ 
0.
\] 
Define a polynomial $P \in \Qbar[X,Y]$ by the prescription
\[
P(X,Y) 
\ = \ 
\frac{X}{C} (Y + 1) - \alpha (Y - 1).
\] 
Then
\[
P(\pi, e^{\pi \sqrt{D} \hsx /  q}) 
\ = \ 0.
\] 
But $\pi$ and $e^{\pi \sqrt{D} \hsx /  q}$ are algebraically independent over $\Q$ (cf. \#2), hence are algebraically indpendent over $\Qbar$ (cf. \S20, \#7).]
\end{spacing}
\end{x}
\vspace{0.3cm}

\begin{x}{\small\bf \un{N.B.}} \ 
For any positive real number \mC (not necessarily rational),
\allowdisplaybreaks
\begin{align*}
\sum\limits_{n=-\infty}^\infty \hsx \frac{1}{n^2+ C^2} \ 
&=\ 
\frac{\pi}{C}  \hsx
\bigg(
\frac{e^{2 \pi C} + 1}{e^{2 \pi C} - 1}
\bigg)
\\[12pt]
&=\ 
\frac{\pi}{C}  \hsx
\bigg(
\frac{e^{\pi C} + e^{-\pi C}}{e^{\pi C} - e^{-\pi C}}
\bigg).
\end{align*}
\end{x}
\vspace{0.3cm}

\begin{x}{\small\bf RAPPEL} \ 
\[
\coth z 
\ = \ 
\frac{\cosh z}{\sinh z}
\ = \ 
\frac{e^z + e^{-z}}{e^z - e^{-z}}.
\] 
\end{x}
\vspace{0.3cm}

\begin{x}{\small\bf \un{N.B.}} \ 
So, for any positive real number \mC (not necessarily rational),
\[
\sum\limits_{n=-\infty}^\infty \hsx 
\frac{1}{n^2+ C^2}
\ = \ 
\frac{\pi}{C} \coth (\pi C).
\] 
\vspace{0.2cm}

[There is another approach to this result using complex variables.  
Thus let
\[
f(z) 
\ = \ 
\frac{1}{z^2+ C^2}
\qquad (C > 0).
\] 
Then $f(z)$ has simple poles at $z = \pm C \sqrt{-1}$.
\vspace{0.3cm}

\qquad \textbullet \quad The residue of 
\[
\frac{\pi \hsx \cot (\pi z)}{z^2+ C^2}
\]
at $z = C \sqrt{-1}$ is
\allowdisplaybreaks
\begin{align*}
\lim\limits_{z \hsx \ra\hsx  C \sqrt{-1}} (z - C \sqrt{-1}) \ \frac{\pi \cot (\pi z)}{(z -C \sqrt{-1}) (z + C \sqrt{-1})} \ 
&=\ 
\frac{\pi \cot (\pi C \sqrt{-1})}{2 C \sqrt{-1}}
\\[12pt]
&=\ 
- \frac{\pi}{2 C} \hsx \coth (\pi C).
\end{align*}
\vspace{0.3cm}

\qquad \textbullet \quad The residue of 
\[
\frac{\pi \hsx \cot (\pi z)}{z^2+ C^2}
\]
at $z = - C \sqrt{-1}$ is
\[
- \frac{\pi}{2 C} \hsx \coth (\pi C).
\]
Since the sum of the residues is
\[
- \frac{\pi}{C} \hsx \coth (\pi C),
\]
it follows that
\allowdisplaybreaks
\begin{align*}
\sum\limits_{n=-\infty}^\infty \hsx \frac{1}{n^2+ C^2}\ 
&=\ 
-\text{(sum of residues)}
\\[12pt]
&=\ 
\frac{\pi}{C} \hsx \coth (\pi C).]
\end{align*}
\vspace{0.2cm}

[Note: \ 
The formalism here is that
\[
\sum\limits_{n=-\infty}^\infty \hsx f(n) \ = \ -S,
\]
where \mS is the sum of the residues of $\pi \cot (\pi z)$ $f(z)$ at the poles of $f(z)$.]
\end{x}
\vspace{0.3cm}

\begin{x}{\small\bf LEMMA} \ 
For any positive real number \mC (not necessarily rational),
\[
\sum\limits_{n=1}^\infty \hsx \frac{1}{n^2+ C^2}
\ = \ 
\frac{\pi}{2C} \hsx \coth (\pi C) - \frac{1}{2 C^2}.
\]
\vspace{0.2cm}

PROOF \
Write
\[
\sum\limits_{n\hsx = \hsx -\infty}^{-1} \hsx \frac{1}{n^2+ C^2} 
\hsx + \hsx \frac{1}{C^2} 
\hsx + \hsx \sum\limits_{n\hsx = \hsx 1}^\infty \hsx \frac{1}{n^2+ C^2}
\ = \ 
\frac{\pi}{C} \hsx \coth (\pi C). 
\]
\end{x}
\vspace{0.3cm}

\begin{x}{\small\bf EXAMPLE} \ 
Take $C = 1$ $-$then
\[
\sum\limits_{n=0}^\infty \  \frac{1}{n^2+ 1}
\ = \ 
\frac{1}{2} + \frac{\pi}{2} \hsx \frac{e^\pi + e^{-\pi}}{e^\pi - e^{-\pi}}.
\]
By comparison,
\[
\sum\limits_{n=0}^\infty \  \frac{1}{n^2 - 1}
\ = \ 
\frac{3}{4}.
\]
\vspace{0.2cm}

[Note: \ 
For the record,
\[
\sum\limits_{n=1}^\infty \  \frac{1}{n^2} \ = \ \frac{\pi^2}{6}.]
\]
\end{x}
\vspace{0.3cm}


\begin{x}{\small\bf REMARK} \ 
It is also possible to sum the series
\[
\sum\limits_{n\hsx =\hsx 0}^\infty \  \frac{(-1)^n}{n^2+ 1},
\]
the result being
\[
\frac{2 \pi}{e^\pi - e^{-\pi}}.
\]
\end{x}
\vspace{0.3cm}

\begin{x}{\small\bf THEOREM} \ 
Let $C \in \Q_{> 0}$ $-$then the series
\[
\sum\limits_{n\hsx = \hsx -\infty}^\infty \hsx \frac{1}{n^2+ C} 
\]
is transcendental.
\vspace{0.2cm}

PROOF \ 
Write
\[
\sum\limits_{n \hsx = \hsx -\infty}^\infty \hsx \frac{1}{n^2+ C} 
\ = \ 
\frac{\pi}{\sqrt{C}} \hsx 
\bigg(
\frac{e^{2 \pi \sqrt{C}} + 1}{e^{2 \pi \sqrt{C}} - 1}
\bigg)
\]
and let
\allowdisplaybreaks
\begin{align*}
C = \frac{p}{q} \quad (p, q \in \N) 
&\implies
\sqrt{C} = \bigg(\frac{p}{q}\bigg)^{1/2} = \frac{\sqrt{p q}}{q}
\\[12pt]
&\implies
2 \pi \sqrt{C} = \pi \sqrt{4} \hsx \frac{\sqrt{p q}}{q}
\\[12pt]
&\hspace{2.1cm}
= \ \pi \frac{\sqrt{4 p q}}{q}.
\end{align*}
Now apply \#2.
\end{x}
\vspace{0.3cm}

\begin{x}{\small\bf EXAMPLE} \ 
Take $C = 3$ $-$then
\[
\sum\limits_{n \hsx = \hsx 0}^\infty \hsx \frac{1}{n^2+ 3} 
\ = \ 
\frac{\pi}{2 \hsx \sqrt{3}} \  
\frac{e^{2 \pi \sqrt{3}} + 1}{e^{2 \pi \sqrt{3}} - 1} + \frac{1}{6}.
\]
\end{x}
\vspace{0.3cm}


\begin{x}{\small\bf THEOREM} \ 
Let $C \in \Q - \{0\}$ $-$then for every positive integer $k$, the series
\[
\sum\limits_{n=-\infty}^\infty \hsx 
\frac{1}{(n^2+ C^2)\strutv^k \upx}
\]
is transcendental.
\vspace{0.2cm}

PROOF \ 
Write
\[
\frac{1}{(n^2+ C^2)\strutv^k \upx}
\ = \ 
\frac{1}{(n + \sqrt{-1} \hsx C)\strutv^k\hsx (n - \sqrt{-1} \hsx C)\strutv^k \upx}
\]
and decompose the term on the right into partial fractions:
\[
\sum\limits_{j=1}^k \hsx 
\frac{\alpha_j}{(n + \sqrt{-1} \hsx C)\strutv^j} 
+ 
\frac{\beta_j}{(n - \sqrt{-1} \hsx C)\strutv^j}
\qquad (\alpha_j, \beta_j \in \Qbar).
\]
Proceed \ldots \ .
\end{x}
\vspace{0.3cm}

\begin{x}{\small\bf EXAMPLE} \ 
Take $C = 1$ $-$then
\[
\sum\limits_{n \hsx = \hsx 0}^\infty \hsx 
\frac{1}{(n^2+ 1)^2}
\ = \ 
\frac{\pi}{4} \hsx 
\frac{e^{2 \pi} + 1}{e^{2 \pi} - 1} 
+ 
\frac{\pi^2}{4} \hsx 
\frac{e^{2 \pi}}{(e^{2 \pi} - 1)^2} 
+ \frac{1}{2}.
\]
\vspace{0.2cm}

[Consider
\[
R(X,Y) 
\ = \ 
\frac{X}{4} \hsx 
\frac{Y + 1}{Y - 1} 
+ \frac{X^2}{4} \hsx 
\frac{Y}{(Y-1)^2} 
+ 
\frac{1}{2}
\]
and write
\[
e^{2 \pi} 
\ = \ 
e^{\pi \sqrt{4}}
\qquad \text{(so $D = 4$)}.]
\]
\end{x}
\vspace{0.3cm}

\begin{x}{\small\bf THEOREM} \ 
Let $C \in \Q_{> 0}$ $-$then for every positive integer $k$, the series
\[
\sum\limits_{n=-\infty}^\infty \hsx 
\frac{1}{(n^2+ C)^k}
\]
is transcendental.

\end{x}
\vspace{0.3cm}


%% file: _04_introduction_of_schc.tex
\chapter{
$\boldsymbol{\S}$\textbf{4}.\quad  INTRODUCTION OF SCHC}
\setlength\parindent{2em}
\setcounter{theoremn}{0}
\renewcommand{\thepage}{Transcendence of Series \S4-\arabic{page}}


\begin{x}{\small\bf THEOREM} \ 
Let $C \in \Q - \Z$ $-$then the series
\[
\sum\limits_{n \hsx = \hsx -\infty}^\infty \hsx 
\frac{1}{n^3+ C^3}
\]
is transcendental.
\vspace{0.2cm}

PROOF \ 
Since
\[
\sum\limits_{n \hsx = \hsx -\infty}^\infty \hsx 
\frac{1}{n^3 - C^3}
\ = \ \sum\limits_{n \hsx = \hsx -\infty}^\infty \hsx 
\frac{1}{(-n)^3 - C^3}
\ = \ 
-\sum\limits_{n \hsx = \hsx -\infty}^\infty \hsx 
\frac{1}{n^3+ C^3},
\]
it can be assumed that \mC is positive.  This said, write
\[
n^3 + C^3 
\ = \ 
(n + C) (n + C \rho) (n + C \rho^2), 
\]
where
\[
\rho \ = \ (-1 - \sqrt{-1} \hsx \sqrt{3}) / 2
\]
is a primitive cube root of unity.  
Decompose $\ds\frac{1}{n^3+ C^3}$ into partial fractions: 
\[
\frac{1}{3 C^2} \hsx \frac{1}{n + C} 
+ 
\frac{\rho}{3 C^2} \hsx \frac{1}{n + C \rho} 
+ 
\frac{\rho^2}{3 C^2} \hsx \frac{1}{n + C \rho^2}.
\]\\[0pt]
Then
\[
\sum\limits_{n \hsx = \hsx -\infty}^\infty \hsx 
\frac{1}{n^3+ C^3}
\]
equals\\
\[
\frac{\pi \sqrt{-1}}{3 \hsx C^2} \ 
\bigg[
\frac{e^{2 \pi \sqrt{-1} \hsx C} + 1}{e^{2 \pi \sqrt{-1} \hsx  C} - 1}
+ 
\rho \hsx
\frac{e^{2 \pi \sqrt{-1} \hsx  C\rho} + 1}{e^{2 \pi \sqrt{-1} \hsx  C \rho} - 1}
+ 
\rho^2 \hsx
\frac{e^{2 \pi \sqrt{-1} \hsx  C\rho^2} + 1}{e^{2 \pi \sqrt{-1} \hsx  C \rho^2} - 1}
\bigg]
.
\]
Here we have used the formula for the cotangent in terms of exponentials (see \S2, 
\#2) (proof thereof).  
Expand the data to arrive at a fraction
\[
2 \pi \sqrt{-1} \hsx \frac{A}{B},
\]
where \mA equals\\
\allowdisplaybreaks
\begin{align*}
\Big(e^{-2 \pi \sqrt{-1} \hsx C} + e^{2 \pi \sqrt{-1} \hsx C}\Big)
&+ 
\rho
\Big(e^{\pi \sqrt{-1} \hsx C} e^{\pi  C \sqrt{3}}  
+ e^{-\pi \sqrt{-1} \hsx C} e^{-\pi  C \sqrt{3}}\Big)
\\[12pt]
&+
\rho^2
\Big(e^{\pi \sqrt{-1} \hsx C} e^{-\pi  C \sqrt{3}}  
+ e^{-\pi \sqrt{-1} \hsx C} e^{\pi  C \sqrt{3}}\Big)
\end{align*}
and \mB equals \\
\[
3 C^2 
\big(e^{2 \pi \sqrt{-1} \hsx C} - 1\big)
\big(e^{2 \pi \sqrt{-1} \hsx C \rho} - 1\big)
\big(e^{2 \pi \sqrt{-1} \hsx C \rho^2} - 1\big).
\]
Owing now to \S3, \hsx \#2, \hsx $\pi$ and $\Big(e^{\pi \sqrt{3}}\Big)^C = e^{\pi C \sqrt{3}}$ are algebraically independent over $\Q$, hence the numerator is either transcendental or zero.  
If the numerator is zero, then the algebraic coefficients of $e^{\pi C \sqrt{3}}$ and $e^{-\pi C \sqrt{3}}$ must both be zero, which implies that
\[
\begin{cases}
\ \rho e^{\pi \sqrt{-1} \hsx C} + \rho^2 e^{-\pi \sqrt{-1} \hsx C} = 0\\[8pt]
\ \rho^2 e^{\pi \sqrt{-1} \hsx C} + \rho e^{-\pi \sqrt{-1} \hsx C} = 0
\end{cases}
.
\]
The first equation implies that
\[
C 
\ = \ 
\frac{1}{6}  + K_1 \qquad (\exists \ K_1 \in \Z)
\]
and the second equation implies that
\[
C 
\ = \ 
-\frac{1}{6}  + K_2 \qquad (\exists \ K_2 \in \Z)
\]
\qquad\qquad $\implies$
\[
\frac{1}{6} + K_1 
\ = \ 
-\frac{1}{6} + K_2
\implies
\frac{1}{3} \ = \  K_2 - K_1,
\]
a contradiction.  Therefore the series is transcendental.
\end{x}
\vspace{0.3cm}

\begin{x}{\small\bf REMARK} \ 
At least one of 
\[
\sum\limits_{n \hsx = \hsx 1}^\infty \hsx 
\frac{1}{n^3+ C^3}
\quad \text{and} \quad 
\sum\limits_{n \hsx = \hsx 1}^\infty \hsx 
\frac{1}{n^3- C^3}
\]
is transcendental.
\end{x}
\vspace{0.3cm}

\begin{x}{\small\bf THEOREM} \ 
Let $C \in \Q - \Z$ $-$then for every positive integer $k$, the series
\[
\sum\limits_{n \hsx = \hsx -\infty}^\infty \hsx 
\frac{1}{(n^3+ C^3)\strutv^k \upx}
\]
is transcendental or zero (transcendental if $k$ is even).
\vspace{0.2cm}

[Start by decomposing
\[
\frac{1}{(n + C)\strutv^k (n + \rho C)\strutv^k (n + \rho^2 C)\strutv^k \upx}
\]
into partial fractions.]
\end{x}
\vspace{0.3cm}

\begin{x}{\small\bf CRITERION} \ 
(Admit SCHC) \ 
If $\alpha_1, \ldots, \alpha_n$ are algebraic numbers such that $\sqrt{-1}, \alpha_1, \ldots, \alpha_n$ are linearly independent over $\Q$, then
\[
\pi, e\strutv^{\textstyle\pi \alpha_1}, \ldots, e\strutv^{\textstyle\pi \alpha_n}
\]
are algebraically independent over $\Q$.
\end{x}
\vspace{0.3cm}

\begin{x}{\small\bf \un{N.B.}} \ 
Take $n = 1$, $\alpha_1 = 1$ $-$then the conclusion is that $\pi$ and $e^\pi$ are algebraically independent over $\Q$ 
(cf. \S3, \#1) (no need for SCHC in this situation).
\end{x}
\vspace{0.3cm}

\begin{x}{\small\bf EXAMPLE} \ 
(Admit SCHC) \ 
Take $n = 2$, $\alpha_1 = \sqrt[\leftroot{0}\uproot{3}3]{C} \hsx\sqrt{3}$, $\alpha_2 = \sqrt{-1} \hsx \sqrt[\leftroot{0}\uproot{3}3]{C}$, where
\[
C \in \Q - \Z, \ C \neq D^3 \quad (D \in \Q).
\]
Then
\[
\pi, \quad e^{\pi \sqrt[\leftroot{0}\uproot{3}3]{C} \hsx\sqrt{3}}, \quad e^{\pi \sqrt{-1} \sqrt[\leftroot{0}\uproot{3}3]{C}}
\]
are algebraically independent over $\Q$.
\vspace{0.2cm}

[To check that $\sqrt{-1}$, $\alpha_1$, $\alpha_2$ are linearly independent over $\Q$, consider a rational dependence relation
\allowdisplaybreaks
\begin{align*}
r \sqrt{-1} + s \alpha_1 + t \alpha_2 \ 
&=\ 
r \sqrt{-1} + s \sqrt[\leftroot{0}\uproot{3}3]{C} \hsx\sqrt{3} + t \sqrt{-1} \sqrt[\leftroot{0}\uproot{3}3]{C} 
\\[12pt]
&=\ 0.
\end{align*}
Then $s = 0$, leaving
\[
r \sqrt{-1} + t \sqrt{-1}  \hsx \sqrt[\leftroot{0}\uproot{3}3]{C} \ = \ 0
\]
or still, 
\allowdisplaybreaks
\begin{align*}
r + t  \hsx \sqrt[\leftroot{0}\uproot{3}3]{C} = 0
&\implies
\sqrt[\leftroot{0}\uproot{3}3]{C} = -\frac{r}{t}
\\[12pt]
&\implies C = \bigg(-\frac{r}{t}\bigg)^3.]
\end{align*}

\end{x}
\vspace{0.3cm}

\begin{x}{\small\bf THEOREM} \ 
(Admit SCHC) \ 
Suppose that $C \in \Q - \Z$ is not a cube in $\Q$ $-$then the series
\[
\sum\limits_{n \hsx = \hsx -\infty}^\infty \hsx \frac{1}{n^3 + C}
\]
is transcendental.
\vspace{0.2cm}

PROOF \ 
The verification is an elaboration of that used in \#1 (which considers the situation when ``\mC'' is a cube).  
So, to begin with, recast matters into the form\\
\[
\frac{\pi \sqrt{-1}}{3 \hsx \sqrt[\leftroot{0}\uproot{3}3]{C^2}} \ 
\bigg[
\frac{e^{2 \pi \sqrt{-1} \hsx \sqrt[\leftroot{0}\uproot{3}3]{C}} + 1}{e^{2 \pi \sqrt{-1} \hsx \sqrt[\leftroot{0}\uproot{3}3]{C}} - 1}
+ 
\rho \hsx
\frac{e^{2 \pi \sqrt{-1} \hsx \sqrt[\leftroot{0}\uproot{3}3]{C}\hsx \rho} + 1}{e^{2 \pi \sqrt{-1} \hsx \sqrt[\leftroot{0}\uproot{3}3]{C}\hsx\rho} - 1}
+ 
\rho^2 \hsx
\frac{e^{2 \pi \sqrt{-1} \hsx \sqrt[\leftroot{0}\uproot{3}3]{C}\hsx\rho^2} + 1}{e^{2 \pi \sqrt{-1} \hsx \sqrt[\leftroot{0}\uproot{3}3]{C}\hsx\rho^2} - 1}
\bigg]
.
\]\\
This done, combine terms in the sum to form a fraction and, using \#6, check that its numerator is not zero.
\end{x}
\vspace{0.3cm}


%% file: _05_introduction_of_schc_bis.tex
\chapter{
$\boldsymbol{\S}$\textbf{5}.\quad  INTRODUCTION OF SCHC (bis)}
\setlength\parindent{2em}
\setcounter{theoremn}{0}
\renewcommand{\thepage}{Transcendence of Series \S5-\arabic{page}}


\begin{x}{\small\bf EXAMPLE} \ 
\[
\sum\limits_{n=0}^\infty \hsx 
\frac{1}{n^4 + 4}
\ = \ 
\frac{\pi}{8} \hsx 
\frac{e^{4 \pi} - 1}{e^{4 \pi} - e^{2 \pi} + 1} + \frac{1}{8}.
\]
\vspace{0.2cm}

[To  ascertain that the right hand side is transcendental, suppose that
\[
\pi \hsx 
\frac{e^{4 \pi} - 1}{e^{4 \pi} - e^{2 \pi} + 1} 
\ = \ 
\alpha \in \Qbar - \{0\}.
\]
Then
\[
\pi (e^{4 \pi} - 1) - \alpha (e^{4 \pi} - e^{2 \pi}  + 1) \ = \ 0.
\]
Define a polynomial $P \hsx \in \hsx \Qbar[X,Y]$ by the prescription
\[
P(X,Y) 
\ = \ 
X (Y^4 - 1) - \alpha(Y^4 - Y^2 + 1) 
\ = \ 
0.
\]
Then
\[
P(\pi,e^\pi) 
\ = \ 
\pi (e^{4 \pi} - 1) - \alpha (e^{4 \pi} - e^{2 \pi}  + 1) \ = \ 0,
\]
which contradicts the fact that $\pi$ and $e^\pi$ are algebraically independent over $\Qbar$.]
\end{x}
\vspace{0.3cm}

\begin{x}{\small\bf LEMMA} \ 
(Admit SCHC) \ 
\[
\pi, \quad e^{\pi \sqrt{2}}, \quad e^{\pi \sqrt{-1} \hsx \sqrt{2}}
\]
are algebraically independent over $\Q$.
\vspace{0.2cm}

PROOF \ 
In \S4, \#4, take $n = 2$, \  $\alpha_1 = \sqrt{2}$,  \ $\alpha_2 = \sqrt{-1} \hsx \sqrt{2}$.
\end{x}
\vspace{0.3cm}

\begin{x}{\small\bf THEOREM} \ 
(Admit SCHC) \ 
Let $C \in \Qbar - \{0\}$ $-$then the series
\[
\sum\limits_{n=-\infty}^\infty \hsx 
\frac{1}{n^4 + C^4}
\]
is transcendental.
\vspace{0.2cm}

PROOF \ 
Write
\[
\frac{1}{n^4 + C^4}
\ = \ 
\frac{1}{n^4 - (\xi C)^4},
\]
where
\[
\xi 
\ = \ 
e^{\pi \sqrt{-1} / 4}
\ = \ 
\sqrt{2} / 2 + \sqrt{-1} \hsx \sqrt{2} / 2.
\]
Then
\[
\sum\limits_{n=-\infty}^\infty \hsx 
\frac{1}{n^4 + C^4}
\]
equals
\[
\frac{\pi}{2 \xi^3 C^3} \ 
\bigg[
\frac{
\big(e^{2\pi \sqrt{-1} \hsx \xi C} + 1\big)
\big(e^{2\pi \xi C} - 1\big) 
 - 
\sqrt{-1}
\big(e^{2\pi \xi C} +1)
\big(e^{2\pi \sqrt{-1} \hsx \xi C} - 1\big)
}
{\sqrt{-1} \hsx
\big(e^{2\pi \sqrt{-1} \hsx \xi C} - 1\big)
\big(e^{2\pi \xi C} - 1\big)
\upx}
\bigg].
\]
Note that
\[
e^{2\pi \sqrt{-1} \hsx \xi C}
\ = \ 
e^{\pi \sqrt{-1} \hsx C \sqrt{2}} e^{-\pi C / \sqrt{2}}
\]
and use the fact that
\[
\pi, \quad e^{\pi \sqrt{2}},\quad e^{\pi \sqrt{-1} \sqrt{2}}
\]
are algebraically independent over $\Q$ (cf. \#2).
\end{x}
\vspace{0.3cm}


%% file: _06_consolidation.tex
\chapter{
$\boldsymbol{\S}$\textbf{6}.\quad  CONSOLIDATION}
\setlength\parindent{2em}
\setcounter{theoremn}{0}
\renewcommand{\thepage}{Transcendence of Series \S6-\arabic{page}}

\ \indent 

Our objective here is to analyze the series
\[
\sum\limits_{n \hsx = \hsx -\infty}^\infty \hsx 
\frac{1}{n^p + C^p},
\]
where $p = 1$ or $p$ is a prime $\geq 2$ and $C \in \Q - \Z$.
\vspace{0.5cm}

\qquad \textbullet \quad \un{$p = 1:$}
\[
\sum\limits_{n \hsx = \hsx -\infty}^\infty \hsx 
\frac{1}{n + C}
\]
is transcendental or zero (cf. \S2, \#2).
\vspace{0.5cm}

\qquad \textbullet \quad \un{$p = 2:$}
\[
\sum\limits_{n \hsx = \hsx -\infty}^\infty \hsx 
\frac{1}{n^2 + C^2}
\]
is transcendental or zero (cf. \S3, \#3).
\vspace{0.5cm}

\qquad \textbullet \quad \un{$p = 3:$}
\[
\sum\limits_{n \hsx = \hsx -\infty}^\infty \hsx 
\frac{1}{n^3 + C^3}
\]
is transcendental or zero (cf. \S4, \#1).
\vspace{0.5cm}

\begin{x}{\small\bf THEOREM} \ 
(Admit SCHC) \ 
Let $p$ be a prime $\geq 5$ and let $C \in \Q - \Z$ $-$then the series
\[
\sum\limits_{n \hsx = \hsx -\infty}^\infty \hsx 
\frac{1}{n^p + C^p}
\]
is transcendental or zero.
\vspace{0.2cm}

PROOF \ 
Let 
\[
\zeta
\ = \ 
e^{2 \pi \sqrt{-1} \hsx / p}
\]
be a primitive $p^\nth$ root of unity $-$then
\[
1, \ \zeta, \ldots, \zeta^{p-2}
\]
are linearly independent over $\Q$, thus
\[
\sqrt{-1}, \sqrt{-1} \hsx \zeta, \ldots, \sqrt{-1}\hsx \zeta^{p-2}
\]
are also linearly independent over $\Q$.  Therefore
\[
\pi, \ e\strutv^{\textstyle\pi \sqrt{-1}\hsx \zeta}, \ldots, e\strutv^{\textstyle\pi \sqrt{-1}\hsx \zeta^{p-2}}
\]
are algebraically independent over $\Q$ (cf. \S4, \#4).  Write
\[
n^p + C^p 
\ = \ 
(n + C) \cdots (n + \zeta^{p-1} C)
\]
to arrive at\\
\[
\pi \sqrt{-1} \  
\bigg(
\alpha_0 \hsx 
\frac{e^{2 \pi \sqrt{-1} \hsx C} + 1}{e^{2 \pi \sqrt{-1} \hsx C} - 1}
+ \cdots + 
\alpha_{p-1} \hsx 
\frac{e^{2 \pi \sqrt{-1} \hsx C \zeta^{p-1}} + 1}{e^{2 \pi \sqrt{-1} \hsx C \zeta^{p-1}} - 1}
\bigg),
\]\\
where the $\alpha_i \in \Qbar$.  Using the fact that
\[
\zeta^{p-1} 
\ = \ 
- 1 - \zeta - \cdots - \zeta^{p-2},
\]
the sum inside the parenthesis can be reduced to a rational function in algebraically independent terms which can be transcendental, zero, or algebraic nonzero but the $\pi$ out in front rules out the last possibility.
\end{x}
\vspace{0.3cm}


%% file: _07_consideration_of_a_over_b.tex
\chapter{
$\boldsymbol{\S}$\textbf{7}.\quad  CONSIDERATION OF $\boldsymbol{\frac{A}{B}}$ }
\setlength\parindent{2em}
\setcounter{theoremn}{0}
\renewcommand{\thepage}{Transcendence of Series \S7-\arabic{page}}

\ \indent 
Let $A(X)$, $B(X)$ be elements of $\Qbar[X]$ with 
\[
\deg A \ < \ \deg B.
\]

Assume: 
\[
B(X) 
\ = \ 
(X + \alpha_1)^{m_1} \cdots (X + \alpha_k)^{m_k},
\]
where $\alpha_1, \ldots, \alpha_k$ are algebraic, nonintegral, and such that 
\[
1, \alpha_1, \ldots, \alpha_k
\]
are linearly independent over $\Q$.
\vspace{0.3cm}

\begin{x}{\small\bf THEOREM} \ 
(Admit SCHC) \ 
The series
\[
\sum\limits_{n = -\infty}^\infty \hsx \frac{A(n)}{B(n)}
\]
is transcendental or zero.
\end{x}
\vspace{0.3cm}

\begin{x}{\small\bf RAPPEL}  \ 
(cf. \S2, \#3) $\forall \ j \geq 2$, $\forall \ z \in \C - \Z$, 
\[
\sum\limits_{n \hsx = \hsx  -\infty}^\infty \hsx \frac{1}{(n + z)^j \upx}
\ = \ 
\frac{(-1)^{j-1} (\pi \cot(\pi z))^{(j-1)}}{(j-1)! \upx}.
\]
\end{x}
\vspace{0.3cm}

\begin{x}{\small\bf \un{N.B.}} \ 
When $j = 1$, 
\[
\sum\limits_{n \hsx = \hsx  -\infty}^\infty \hsx \frac{1}{n + z}
\ = \ 
\pi \cot (\pi z).
\]
Using partial fractions, write
\[
\frac{A(n)}{B(n)}
\ = \ 
\sum\limits_{i \hsx = \hsx 1}^k \  
\sum\limits_{j \hsx = \hsx 1}^{m_i} \  
C_{i j} \hsx \frac{1}{(n + \alpha_i)^j \upx}.
\]
Then 
\allowdisplaybreaks
\begin{align*}
\sum\limits_{n \hsx = \hsx -\infty}^\infty \hsx \frac{A(n)}{B(n)} \ 
&=\ 
\sum\limits_{n \hsx = \hsx -\infty}^\infty \hsx
\bigg(
\sum\limits_{i\hsx= \hsx 1}^k \  
\sum\limits_{j\hsx = \hsx 1}^{m_i} \  
C_{i j} \hsx 
\frac{1}{(n + \alpha_i)^j \upx}
\bigg)
\\[12pt]
&=\ 
\sum\limits_{i=1}^k \  
\sum\limits_{j=1}^{m_i} \  
C_{i j} \hsx 
\sum\limits_{n \hsx = \hsx -\infty}^\infty \hsx
\frac{1}{(n + \alpha_i)^j \upx }
\\[12pt]
&=\ 
\sum\limits_{i\hsx = \hsx1}^k \  
\sum\limits_{j\hsx =\hsx 1}^{m_i} \  
C_{i j} \hsx 
\frac{(-1)^{j-1} (\pi \hsy\cot (\pi \alpha_i))^{(j-1)}}{(j-1)! \upx}
\\[12pt]
&=\ 
\pi \ 
\sum\limits_{i\hsx = \hsx 1}^k \hsx 
\sum\limits_{j\hsx = \hsx 1}^{m_i} \hsx 
D_{i j} \hsx  (\cot (\pi \alpha_i))^{(j-1)},
\end{align*}
where
\[
D_{i j}
\ = \ 
C_{i j} \ \frac{(-1)^{j-1}}{(j-1)!}.
\]
\vspace{0.2cm}

\qquad \textbf{FACT} \ 
For any integer $m > 1$, 
\[
\bigg(\frac{\td}{\td z}\bigg)^m \cot z
\]
is a polynomial in $\cot z$.
\vspace{0.2cm}

[The formula is
\[
\bigg(\frac{\td}{\td z}\bigg)^m \cot z
\]
equals 
\[
(2 \sqrt{-1})\strutv^m \hsx 
(\cot z - \sqrt{-1}) \ 
\sum\limits_{\ell \hsx = \hsx 1}^m \  
\frac{\ell !}{2^\ell} \hsx S(m,\ell) \hsx (\sqrt{-1} \hsx \cot z - 1)^\ell.
\]
Here the $S(m,\ell) \in \Z$ are the Stirling subset numbers (a.k.a. the Stirling numbers of the second kind).]
\vspace{0.2cm}

[Note: \ 
$\forall \ k \geq 2$, $\forall \ z \in \C - \Z$, 
\[
\sum\limits_{n \hsx = \hsx -\infty}^\infty \hsx
\frac{1}{(n+z)^k \upx}
\ = \ 
\frac{(-2 \pi \sqrt{-1})^k}{(k-1)! \upx} 
\sum\limits_{\ell \hsx = \hsx  1}^k \ 
\frac{(\ell - 1)! \hsx S(k,\ell)}{(e^{-2 \pi \sqrt{-1} \hsx z} - 1)\strutv^{\ell} \upx} 
\qquad \text{(cf. \S2, \#3)}.]
\]
\end{x}
\vspace{0.3cm}

\begin{x}{\small\bf RAPPEL} \ 
\[
\cot (\pi z) 
\ = \ 
\sqrt{-1} \ 
\frac{e^{2 \pi \sqrt{-1} \hsx z} + 1}{e^{2 \pi \sqrt{-1} \hsx z} - 1}.
\]
\end{x}
\vspace{0.3cm}

\begin{x}{\small\bf APPLICATION} \ 
\[
(\cot (\pi \alpha_i))^{(j-1)}
\]
is an algebraic linear combination of rational functions evaluated at $e^{2 \pi \sqrt{-1} \hsx \alpha_i}$.
\vspace{0.5cm}

The assumption on the $\alpha_i$ is that 
\[
1, \ \alpha_1, \ldots, \alpha_k
\]
are linearly independent over $\Q$ or still, that
\[
\sqrt{-1}, \ \sqrt{-1} \hsx \alpha_1, \ldots, \sqrt{-1} \hsx \alpha_k
\]
are linearly independent over $\Q$ or still, that
\[
\sqrt{-1}, \ 2\sqrt{-1} \hsx \alpha_1, \ldots, 2\sqrt{-1} \hsx \alpha_k
\]
are linearly independent over $\Q$.  Therefore
\[
\pi, \ e^{\textstyle 2 \pi \sqrt{-1} \hsx \alpha_1}, \ldots, e^{\textstyle 2 \pi \sqrt{-1} \hsx \alpha_k}
\]
are algebraically independent over $\Q$ (cf. \S4, \#4).  \\

To finish the proof, rearrange the sum so as to form a  polynomial in $\pi$, 
the coefficients of a given power of $\pi$ being a rational expression in 
\[
e^{\textstyle 2 \pi \sqrt{-1} \hsx \alpha_1}, \ldots, e^{\textstyle 2 \pi \sqrt{-1} \hsx \alpha_k}.
\]
Complete the argument by citing algebraic independence over $\Q$ (which eliminates the algebraic nonzero possibility).
\end{x}
\vspace{0.3cm}

There is one set of circumstances under which the series
\[
\sum\limits_{n \hsx = \hsx  -\infty}^\infty \hsx \frac{A(n)}{B(n)}
\]
is transcendental (thereby ruling out the zero contingency).\\

Assume: \ The roots of $B(X)$ are simple, hence
\[
m_1 = 1, \ldots, m_k = 1.
\]

To proceed, write
\[
\sum\limits_{n \hsx = \hsx -\infty}^\infty \hsx \frac{A(n)}{B(n)}
\ = \ 
\pi \hsx 
\sum\limits_{i \hsx = \hsx 1}^k \hsx
C_i \hsx \cot(\pi \alpha_i)
\]
or still, 
\[
\sum\limits_{n \hsx = \hsx -\infty}^\infty \hsx \frac{A(n)}{B(n)}
\ = \ 
\pi \sqrt{-1} \  
\sum\limits_{i \hsx = \hsx 1}^k \ 
C_i \hsx 
\frac{e^{2 \pi \sqrt{-1} \hsx \alpha_i} + 1}{e^{2 \pi \sqrt{-1} \hsx \alpha_i} - 1},
\]
the claim being that the expression on the right is nonzero, thus that the series
\[
\sum\limits_{n \hsx =  \hsx -\infty}^\infty \hsx \frac{A(n)}{B(n)}
\]
is transcendental.

Rewrite the expression as
\[
\frac{\pi \sqrt{-1}}{\prod\limits_{i \hsx = \hsx 1}^k \hsx \big(e^{2 \pi \sqrt{-1} \hsx \alpha_i} - 1\big) } \ 
\sum\limits_{i \hsx = \hsx  1}^k \hsx 
C_i 
\big(e^{2 \pi \sqrt{-1} \hsx \alpha_i} + 1 \big) \ 
\prod\limits_{a \hsx \neq  \hsx i} \hsx
\big(e^{2 \pi \sqrt{-1} \hsx \alpha_a} - 1 \big).
\]
Matters then reduce to showing that the polynomial 
\[
\sum\limits_{i \hsx = \hsx  1}^k \ 
C_i (X_i + 1) \ 
\prod\limits_{a \hsx \neq  \hsx i} \  
(X_a - 1)
\]
is not identically zero.  
Suppose it were identically zero.  
Given $i$, take
\[
X_i = 0, \ 
X_j = -1, \  (j \neq i), \ 
X_a = 2 \ (a \neq i)
\]
to see that $C_i = 0$.  
But $i$ is arbitrary, so $C_i = 0$ $\forall \ i$, contradicting the tacit assumption that $A \neq 0$.


%% file: _08_an_algebraic_series.tex
\chapter{
$\boldsymbol{\S}$\textbf{8}.\quad  AN ALGEBRAIC SERIES}
\setlength\parindent{2em}
\setcounter{theoremn}{0}
\renewcommand{\thepage}{Transcendence of Series \S8-\arabic{page}}

\ \indent 
Instead of looking for a transcendental series, this time we shall exhibit an algebraic series.
\vspace{0.3cm}

\begin{x}{\small\bf THEOREM} \ 
Suppose that $P(X) \in \Qbar[X]$ and $z \in \Qbar$ $(0 < \abs{z} < 1)$ $-$then the series
\[
\sum\limits_{n \hsx = \hsx 0}^\infty \  z^n P(n)
\]
is algebraic.
\vspace{0.2cm}

PROOF \ 
First of all, the manipulations infra are justified by the absolute convergence of our series, so if
\[
P(X) 
\ = \
\sum\limits_{i \hsx = \hsx 0}^k \hsx a_i X^i, 
\]
then
\[
\sum\limits_{n \hsx = \hsx 0}^\infty \  z^n P(n)
\ = \
\sum\limits_{i \hsx = \hsx 0}^k \  a_i 
\sum\limits_{n \hsx= \hsx 0}^\infty \  z^n n^i.
\]
Write
\[
X^i 
\ = \ 
\sum\limits_{j \hsx = \hsx 0}^i \  S(i,j) \hsx (X)_j,
\]
where $(X)_0 = 1$ and for $j \geq 1$, 
\[
(X)_j
\ = \ 
X (X - 1) \cdots (X - j + 1).
\]
Inserting this data leads to
\[
\sum\limits_{i \hsx = \hsx 0}^k \  a_i \ 
\sum\limits_{j \hsx = \hsx 0}^i \ S(i,j) \
\sum\limits_{n \hsx = \hsx 0}^\infty \ (n)_j z^n 
\]
or still, 
\[
\sum\limits_{i \hsx = \hsx 0}^k \ a_i \
\sum\limits_{j \hsx = \hsx 0}^i \ S(i,j) \
\sum\limits_{n \hsx = \hsx 0}^\infty \ n(n-1) \cdots(n - j + 1) z^n
\]
or still, 
\[
\sum\limits_{i \hsx = \hsx 0}^k \ a_i \
\sum\limits_{j \hsx = \hsx 0}^i \ S(i,j) \
\sum\limits_{n \hsx = \hsx 1}^\infty \hsx n(n-1) \cdots(n - j + 1) z^n
\]

\[
\vdots
\]
or still, 
\[
\sum\limits_{i \hsx = \hsx 0}^k \hsx a_i \hsx
\sum\limits_{j \hsx = \hsx 0}^i \hsx S(i,j) \hsx
\sum\limits_{n \hsx = \hsx j-1}^\infty \hsx n(n-1) \cdots(n - j + 1) z^n
\]
or still, 
\[
\sum\limits_{i \hsx = \hsx 0}^k \hsx a_i \hsx
\sum\limits_{j \hsx = \hsx 0}^i \hsx S(i,j) \hsx
\sum\limits_{n \hsx = \hsx j}^\infty \hsx n(n-1) \cdots(n - j + 1) z^n
\]
or still, 
\[
\sum\limits_{i \hsx = \hsx 0}^k \hsx a_i \hsx
\sum\limits_{j \hsx = \hsx 0}^i \hsx S(i,j) \hsx z^j \hsx
\sum\limits_{n \hsx = \hsx 0}^\infty \hsx (n+1) \cdots (n+j) z^n
\]
or still, 
\[
\sum\limits_{i \hsx = \hsx 0}^k \hsx a_i \hsx
\sum\limits_{j \hsx = \hsx 0}^i \hsx S(i,j) \hsx z^j \hsx \bigg(\frac{z^j}{1 - z} \bigg)^{(j)}
\]
or still, 
\[
\sum\limits_{i \hsx = \hsx 0}^k \hsx a_i \hsx
\sum\limits_{j \hsx = \hsx 0}^i \hsx \frac{S(i,j) \hsx j! \hsx z^j}{(1 - z)^{j+1}\upy},
\]
an algebraic number.
\end{x}
\vspace{0.3cm}


%% file: _XX_Supp_zeta_function_values.tex
\centerline{\textbf{\large SUPPLEMENT II}}
\vspace{1cm}
\renewcommand{\thepage}{Zeta Function Values - 0}
\centerline{\textbf{ZETA FUNCTION VALUES}}


\allowdisplaybreaks
\begin{align*}
\S1.  \qquad &\text{BERNOULLI NUMBERS} 
\\[10pt]
\S2.  \qquad &\text{$\zeta(2n)$}  
\\[10pt]
\S3.  \qquad &\text{$\zeta(2)$}
\\[10pt]
\S4.  \qquad &\text{$\zeta(2)$ (bis)}
\\[10pt]
\S5.  \qquad &\text{$\zeta(2n)$ (bis)}
\\[10pt]
\S6.  \qquad &\text{$\zeta(3)$}
\\[10pt]
\S7.  \qquad &\text{CONJUGATE BERNOULLI NUMBERS}
\\[10pt]
\S8.  \qquad &\text{$\zeta(2n+1)$}
\end{align*}

%% file: _01_Bernoulli_Numbers.tex
\chapter{
$\boldsymbol{\S}$\textbf{1}.\quad  BERNOULLI NUMBERS}
\setlength\parindent{2em}
\setcounter{theoremn}{0}
\renewcommand{\thepage}{Zeta Function Values I \S1-\arabic{page}}

\indent Define the Bernoulli polynomials $B_n(x)$ $(n = 0, 1, 2, \ldots)$ via the generating function\\[5pt]
\[
\frac{t e^{xt}}{e^t - 1} 
\ = \ 
\sum\limits_{n \hsx  = \hsx  0}^\infty \hsx B_n(x) \hsx \frac{t^n}{n!}.
\]
\vspace{0.2cm}

[Note: \ 
\[
B_0(x) = 1, \ 
B_1(x) = x - \frac{1}{2}, \ 
B_2(x) = x^2 - x + \frac{1}{6}.]
\]
\\[-20pt]

There are two sign conventions at play here.
\vspace{0.5cm}

\qquad $(+)$ \ 
Define the Bernoulli numbers $B_n^+$ $(n = 0, 1, 2, \ldots)$ by taking $x = 1$, hence the generating function \\
\[
\frac{t e^{t}}{e^t - 1} 
\ = \ 
\sum\limits_{n \hsx  = \hsx  0}^\infty \hsx B_n^+ \hsx \frac{t^n}{n!}.
\]
\vspace{0.2cm}

[Note: \ 
$B_0^+= 1, \ 
B_1^+ = \ds\frac{1}{2}, \ 
B_2^+ = \ds\frac{1}{6}$.]
\vspace{0.5cm}

\qquad $(-)$ \ 
Define the Bernoulli numbers $B_n^-$ $(n = 0, 1, 2, \ldots)$ by taking $x = 0$, hence the generating function\\[3pt]
\[
\frac{t}{e^t - 1} 
\ = \ 
\sum\limits_{n \hsx  = \hsx  0}^\infty \hsx B_n^- \hsx \frac{t^n}{n!}.
\]
\vspace{0.2cm}

[Note: \ 
$B_0^- = 1, \ 
B_1^- = -\ds\frac{1}{2}, \ 
B_2^- = \ds\frac{1}{6}$.]
\vspace{0.5cm}

\begin{x}{\small\bf REMARK} \ 
A Bernoulli number is real and rational.
\end{x}
\vspace{0.3cm}


\begin{x}{\small\bf LEMMA} \ 
\[
B_n^+
\ = \
(-1)^n B_n^-.
\]
\end{x}
\vspace{0.3cm}

\begin{x}{\small\bf LEMMA} \ 
If $n$ is an odd integer $\geq 3$, then
\[
B_n^+
 = 
0, \quad
B_n^-
 = 
0.
\]
\end{x}
\vspace{0.3cm}

\begin{x}{\small\bf \un{N.B.}} \ 
In formulas involving even index Bernoulli numbers, it is permissible to drop the $\pm$ and simply use the symbol $B_n$.
\end{x}
\vspace{0.3cm}

\begin{x}{\small\bf EXAMPLE} \ 
\[
x \cot x 
\hsx = \hsx
\sum\limits_{n \hsx  = \hsx  0}^\infty \hsx (-1)^n \hsx \frac{2^{2 n}}{(2 n)!} B_{2n} \hsx x^{2n} 
\qquad (0 < \abs{x} < \pi).
\]
\end{x}
\vspace{0.3cm}

\begin{x}{\small\bf LEMMA} \ 
$\forall \ n \geq 1$,
\[
\int\limits_0^1 \hsx B_n(x) \hsx dx
\hsx = \hsx
0.
\]
\end{x}
\vspace{0.3cm}

\begin{x}{\small\bf LEMMA} \ 
$\forall \ n \geq 1$,
\[
\int\limits_0^1 \hsx B_n(x) \hsx B_m(x) \hsx dx
\hsx = \hsx
(-1)^{n - 1} \hsx \frac{m! \hsx n!}{(m + n)!} \hsx B_{m + n}^-.
\]
\end{x}
\vspace{0.3cm}

\begin{x}{\small\bf LEMMA} \ 
$\forall \ n \geq 1$,
\[
\frac{d}{dx} \hsx B_n(x) 
\hsx = \hsx
n B_{n - 1}(x).
\]
\end{x}
\newpage

\[
\textbf{APPENDIX}
\]
\vspace{0.5cm}

{\small\bf LEMMA} (MULTIPLICATION FORMULA) \ 
\[
B_n(m x)
\ = \
m^{n - 1} \hsx \sum\limits_{k \hsx  = \hsx  0}^{m - 1} \hsx B_n \bigg(x + \frac{k}{m}\bigg).
\]
\vspace{0.3cm}

{\small\bf APPLICATION} \ 
Take $x = 0$, $m = 2$ $-$then
\[
B_{2n} \bigg(\frac{0}{2} \bigg) + B_{2n} \bigg(\frac{1}{2} \bigg)
\ = \
2^{1 - 2n} B_{2n}(0),
\]
i.e.,
\allowdisplaybreaks
\begin{align*}
 B_{2n} \bigg(\frac{1}{2} \bigg) \ 
 &=\ 2^{1 - 2n} B_{2n} - B_{2n}\\[12pt]
  &=\ (2^{1 - 2n} - 1) B_{2n}.
\end{align*}
\vspace{0.3cm}

{\small\bf LEMMA} (ADDITION FORMULA) \ 
\[
B_n(x + y)
\hsx = \hsx
\sum\limits_{k \hsx  = \hsx  0}^n \hsx {n\choose k} \hsx B_k(x) y^{n - k}.
\]
\vspace{0.3cm}


%% file: _02_Zeta_2n.tex
\chapter{
$\boldsymbol{\S}$\textbf{2}.\quad  $\zeta(2 n)$}
\setlength\parindent{2em}
\setcounter{theoremn}{0}
\renewcommand{\thepage}{Zeta Function Values \S2-\arabic{page}}


\begin{x}{\small\bf THEOREM} \ 
$\forall \ n \geq 1$, 
\allowdisplaybreaks
\begin{align*}
\zeta (2 n ) \ 
&\equiv\ \sum\limits_{k \hsx  = \hsx  1}^\infty \hsx  \frac{1}{k^{2n}}\\[12pt]
&=\ (-1)^{n-1} \frac{(2 \pi)^{2n}}{2 (2n)!} \hsx B_{2n}
\end{align*}
or still, 
\[
\zeta (2 n )
\hsx = \hsx 
(-1)^{n - 1} \hsx \frac{2^{2n - 1}}{(2n)!} \hsx B_{2n} \pi^{2n}.
\]
\end{x}
\vspace{0.3cm}

\begin{x}{\small\bf APPLICATION} \ 
$\zeta (2 n )$ is transcendental.
\vspace{0.2cm}

[Recall that $\pi$ is transcendental, hence $\pi^{2n}$ is transcendental.]
\end{x}
\vspace{0.3cm}

The stated formular for $\zeta (2 n )$ can now be proved in many different ways.  
What follows is one of them.
\vspace{0.3cm}

\begin{x}{\small\bf NOTATION} \ 
Given an $f \in L^1[0,1]$, put
\[
\widehat{f} (k)
\hsx = \hsx 
\int\limits_0^1 \hsx f(x) e^{-2 \pi \sqrt{-1} \hsx k x} \hsx dx \qquad (k \in \Z).
\]
\end{x}
\vspace{0.3cm}

\begin{x}{\small\bf PLANCHEREL} \ 
Given an $f \in L^2[0,1]$, 
\[
\int\limits_0^1 \hsx \abs{f(x)}^2 \hsx dx
\ = \ 
\sum\limits_{-\infty}^\infty \  \abs{\widehat{f}(k)}^2.
\]
\vspace{0.1cm}

[Note: \ 
Recall that 
\[
L^2[0,1] \subset L^1[0,1] \hsx .]
\]
\end{x}
\vspace{0.3cm}

\begin{x}{\small\bf LEMMA} \ 
Take $f(x) = B_n(x)$ $-$then
\[
\widehat{B}_n(k) 
\hsx = \hsx 
- \frac{n!}{(2 \pi \sqrt{-1} \hsx k)^n}
\]
if $k \neq 0$ while $\widehat{B}_n(0) = 0.$
\vspace{0.2cm}

PROOF \ 
The second point is covered by \S1, \#6.  
As for the first point, take $n \geq 1$ and write 
\allowdisplaybreaks
\begin{align*}
\widehat{B}_n(k) \
&=\ \int\limits_0^1 \hsx B_n(x) e^{-2 \pi \sqrt{-1} \hsx k x} dx \\[12pt]
&=\ -\frac{1}{2 \pi \sqrt{-1} \hsx k} \hsx \int\limits_0^1  \hsx B_n(x) \frac{d}{dx}e^{-2 \pi \sqrt{-1} \hsx k x} dx \\[12pt]
&=\ -\frac{1}{2 \pi \sqrt{-1} \hsx k} \hsx B_n(x) e^{-2 \pi \sqrt{-1} \hsx k x} \bigg\arrowvert_0^1 
+
\frac{1}{2 \pi \sqrt{-1} \hsx k} \hsx \int\limits_0^1  \hsx\frac{d}{dx}  B_n(x) e^{-2 \pi \sqrt{-1} \hsx k x} dx.
\end{align*}
\vspace{0.2cm}

\un{$n = 1$}:
\allowdisplaybreaks
\begin{align*}
\widehat{B}_1(k) \ 
&=\ -\frac{1}{2 \pi \sqrt{-1} \hsx k} \hsx \bigg(x - \frac{1}{2}\bigg) e^{-2 \pi \sqrt{-1} \hsx k x} \bigg\arrowvert_0^1
-
\frac{1}{2 \pi \sqrt{-1} \hsx k} \hsx \int\limits_0^1  \hsx 1 \cdot e^{-2 \pi \sqrt{-1} \hsx k x} dx\\[12pt]
&=\ -\frac{1}{2 \pi \sqrt{-1} k} \hsx \bigg(\frac{1}{2} + \frac{1}{2}\bigg) 
- 
\frac{1}{2 \pi \sqrt{-1} \hsx k} \  0 \qquad (k \neq 0)\\[12pt]
&=\ -\frac{1}{2 \pi \sqrt{-1} \hsx k}.
\end{align*}
\vspace{0.2cm}

\un{$n > 1$}: \ To begin with 
\[
-\frac{1}{2 \pi \sqrt{-1} \hsx k} \hsx B_n(x) \hsx e^{-2 \pi \sqrt{-1} \hsx k x} \bigg\arrowvert_0^1 
\hsx = \hsx
-\frac{1}{2 \pi \sqrt{-1} \hsx k} \hsx (B_n(1) - B_n(0)).
\]
And 
\allowdisplaybreaks
\begin{align*}
B_n(1) - B_n(0) \ 
&=\ B_n^+ - B_n^- \\[12pt]
&=\ (-1)^n B_n^- - B_n^- \qquad (\text{cf. \S1, \#2})
\\[12pt]
&=\ B_n^- ((-1)^n - 1).
\end{align*}
But 
\[
\begin{cases}
\ n \ \text{even}, \ \geq 2 \implies (-1)^n = 1 \implies B_n^- ((-1)^n - 1) = 0\\[8pt]
\ n \ \text{odd}, \ \geq 3 \implies B_n^- = 0 \ \text{(cf. \S1, \#3}) \implies  B_n^- ((-1)^n - 1) = 0\\
\end{cases}
.
\]
Therefore
\[
B_n(1) - B_n(0) \hsx = \hsx 0,
\]
leaving 
\[
\frac{1}{2 \pi \sqrt{-1} \hsx k} \hsx  \int\limits_0^1 \hsx \frac{d}{dx} B_n(x) e^{-2\pi \sqrt{-1} \hsx kx} dx.
\]
Using \S1, \#8, replace $\ds\frac{d}{dx} B_n(x)$ by $n B_{n-1}(x)$ to arrive at
\allowdisplaybreaks
\begin{align*}
\widehat{B}_n(k) \ 
&=\ \frac{n}{2 \pi \sqrt{-1} \hsx k} \hsx \int\limits_0^1 \hsx B_{n-1} (x) e^{-2 \pi \sqrt{-1} \hsx k x} dx\\[12pt]
&=\ \frac{n}{2 \pi \sqrt{-1} \hsx k} \hsx \widehat{B}_{n-1}(k),
\end{align*}
so, inductively, 
\allowdisplaybreaks
\begin{align*}
\widehat{B}_n(k) \ 
&=\ \frac{n}{2 \pi \sqrt{-1} \hsx k} \cdot \frac{n-1}{2 \pi \sqrt{-1} \hsx k} \hsx \widehat{B}_{n-2}(k)\\[12pt]
&\hspace{.2cm}  \vdots\\[12pt]
&=\ \frac{n (n-1) \cdots 2}{2 \pi \sqrt{-1} \hsx k)^{n-1}} \hsx \widehat{B}_1(k)\\[12pt]
&=\ \frac{n!}{(2 \pi \sqrt{-1} \hsx k)^{n-1}} \hsx \bigg( -\frac{1}{2 \pi \sqrt{-1} \hsx k}\bigg)\\[12pt]
&=\ -\frac{n!}{(2 \pi \sqrt{-1} \hsx k)^{n}}.
\end{align*}
Hence the lemma.
\vspace{0.5cm}

To prove the theorem, take $f = B_n$ $(n \geq 1)$ in Plancherel:
\[
\int\limits_0^1 \hsx \abs{B_n(x)}^2 dx 
\ = \ 
\sum\limits_{-\infty}^\infty \hsx \abs{\widehat{B}_n(k)}^2.
\]
Here
\allowdisplaybreaks
\begin{align*}
\int\limits_0^1 \hsx \abs{B_n(x)}^2 dx  \ 
&=\ \int\limits_0^1 \hsx {B_n(x)}\hsx {B_n(x)} dx \\[12pt]
&=\ (-1)^{n-1} \hsx \frac{(n!)^2}{(2n)!} \hsx B_{2n}^- \qquad (\text{cf. \S1, \#7})\\[12pt]
&=\ (-1)^{n-1} \hsx \frac{(n!)^2}{(2n)!} \hsx B_{2n} \qquad (\text{cf. \S1, \#4}).
\end{align*}
On the other hand, 
\allowdisplaybreaks
\begin{align*}
\sum\limits_{-\infty}^\infty \hsx \abs{\widehat{B}_n(k)}^2 \ 
&=\ \sum\limits_{k \hsx \neq \hsx 0} \hsx \abs{-\frac{n!}{(2\pi \sqrt{-1} \hsx k)^n}}^2\\[12pt]
&=\ 2 \hsx \sum\limits_{k \hsx  = \hsx  1}^\infty \hsx \frac{(n!)^2}{(2 \pi k)^{2n}}\\[12pt]
&=\ 2 \hsx \frac{(n!)^2}{(2 \pi)^{2n}} \hsx \sum\limits_{k \hsx  = \hsx  1}^\infty  \hsx \frac{1}{k^{2n}}\\[12pt]
&=\ 2 \hsx \frac{(n!)^2}{(2 \pi)^{2n}} \hsx \zeta(2n).
\end{align*}

Now cancel the $(n!)^2$ to get 
\[
\zeta(2n) 
\hsx = \hsx 
(-1)^{n-1} \hsx \frac{(2\pi)^{2n}}{2 (2n)!} \hsx B_{2n}.
\]
\end{x}
\vspace{0.3cm}

\begin{x}{\small\bf SCHOLIUM} \ 
\[
\Q[\zeta(2), \ \zeta(4),\ \zeta(6), \ldots] 
\hsx = \hsx 
\Q[\pi^2].
\]
\end{x}

\vspace{0.1cm}


%% file: _03_Zeta_2.tex
\chapter{
$\boldsymbol{\S}$\textbf{3}.\quad  $\zeta(2)$}
\setlength\parindent{2em}
\setcounter{theoremn}{0}
\renewcommand{\thepage}{Zeta Function Values \S3-\arabic{page}}

\ \indent 
In \S2, \#1, take $n = 1$ to get 
\[
\zeta(2) 
\hsx = \hsx
\frac{\pi^2}{6}.
\]
Of course there are a ``million'' proofs of this result but for motivational purposes we shall single out one of these.
\vspace{0.3cm}

\begin{x}{\small\bf NOTATION} \ 
The symbol
\[
\ds\int\limits_0^1 \hsx
\ds\int\limits_0^1 \hsx
f(x,y) \hsx dx dy
\]
stands for a double integral over the unit square $[0,1] \hsx \times \hsx [0,1]$, possibly improper.
\end{x}
\vspace{0.3cm}

\begin{x}{\small\bf SUBLEMMA} \ 
\[
\frac{3}{4} \hsx \zeta(2)
\ = \ 
\sum\limits_{n \hsx = \hsx 0}^\infty \hsx \frac{1}{(2n+1)^2}.
\]
\vspace{0.2cm}

PROOF \ 

\allowdisplaybreaks
\begin{align*}
\zeta(2)\
&=\ 
\sum\limits_{n \hsx = \hsx 1}^\infty \hsx \frac{1}{n^2}\\[15pt]
&=\ \sum\limits_{n \hsx = \hsx 0}^\infty \hsx \frac{1}{(2n+1)^2} \hsx + \hsx \sum\limits_{n=1}^\infty \hsx \frac{1}{(2n)^2}\\[15pt]
&=\ \sum\limits_{n \hsx = \hsx 0}^\infty \hsx \frac{1}{(2n+1)^2} \hsx + \hsx 
\frac{1}{4} \hsx \sum\limits_{n=1}^\infty \hsx \frac{1}{n^2}\\[15pt]
&=\ \sum\limits_{n \hsx = \hsx 0}^\infty \hsx \frac{1}{(2n+1)^2} \hsx + \hsx  \frac{1}{4} \hsx \zeta(2)
\end{align*}

\qquad\qquad $\implies$ 
\[
\frac{3}{4} \hsx \zeta(2) 
\ = \ 
\sum\limits_{n \hsx = \hsx 0}^\infty \hsx \frac{1}{(2n+1)^2}.
\]
\end{x}
\vspace{0.3cm}

\begin{x}{\small\bf LEMMA} \ 
\allowdisplaybreaks
\begin{align*}
\int\limits_0^1 \hsx
\int\limits_0^1 \hsx
\frac{1}{1 - x^2 y^2} dx dy
&=\ 
\int\limits_0^1 \hsx
\int\limits_0^1 \ 
\sum\limits_{n \hsx = \hsx 0}^\infty \hsx (x y)^{2n} \hsx dx dy \\[15pt]
&=\ \sum\limits_{n \hsx = \hsx 0}^\infty \hsx \frac{1}{(2n+1)^2}\\[15pt]
&=\ \frac{3}{4} \hsx \zeta(2).
\end{align*}
\vspace{0.2cm}

[Note: \ The singularity at the corner $(x,y) = (1,1)$ can be safely ignored \ldots \hsx .]
\vspace{0.5cm}

Define a bijective map from 
\[
\Pi_2 
\hsx \equiv \hsx 
\bigg\{(u,v):u > 0, v > 0, u + v < \frac{\pi}{2}\bigg\}
\]
to $]0,1[ \hsx \times \hsx ]0,1[$ by the prescription
\[
(u, v) 
\ra 
\bigg(\frac{\sin u}{\cos v}, \frac{\sin v}{\cos u} \bigg)
\]
with Jacobian
\allowdisplaybreaks
\begin{align*}
\frac{\partial(x,y)}{\partial(u,v)} \ 
&=\ 
\begin{pmatrix}
\cos u / \cos v & \sin u \sin v / \cos^2 v\\[12pt]
\sin u \sin v / \cos^2 u &\cos v / \cos u\\
\end{pmatrix}\\[15pt]
&=\ 1 - \frac{\sin^2 u \sin^2 v}{\cos^2 u \cos^2 v} \\[15pt]
&=\ 1 - x^2 y^2.
\end{align*}

[Note: \  The details are in the Appendix to this \S.]
\end{x}
\vspace{0.3cm}

Therefore
\allowdisplaybreaks
\begin{align*}
\frac{3}{4} \zeta(2) \ 
&=\ 
\int\limits_0^1 \hsx
\int\limits_0^1 \hsx
\frac{1}{1 - x^2 y^2} \hsx dx dy \\[15pt]
&=\ \text{Area}(\Pi_2) \\[15pt]
&=\ \frac{\pi^2}{8}
\end{align*}
\qquad\qquad $\implies$
\[
\zeta(2) 
\hsx = \hsx 
\frac{\pi^2}{6}.
\]
\vspace{0.2cm}

\begin{x}{\small\bf LEMMA} \ 
\[
\zeta(2) 
\hsx = \hsx 
\ds\int\limits_0^1 
\ds\int\limits_0^1
\frac{1}{1 - x y} dx dy.
\]
\vspace{0.2cm}

PROOF \ 
The RHS equals
\[
\int\limits_0^1 \hsx
\int\limits_0^1 \  
\sum\limits_{n \hsx = \hsx 0}^\infty \hsx x^n y^n dx dy
\]
or still, 
\[
\sum\limits_{n \hsx = \hsx 0}^\infty \hsx 
\bigg( \int\limits_0^1 \hsx x^n dx \bigg) 
\hsx \cdot \hsx 
\bigg( \int\limits_0^1 \hsx y^n dy \bigg) 
\]
or still, 
\[
\sum\limits_{n \hsx = \hsx 0}^\infty \hsx \frac{x^{n+1}}{n+1} \bigg\arrowvert_0^1 
\ \cdot \ 
 \frac{y^{n+1}}{n+1} \bigg\arrowvert_0^1
\]
or still, 
\[
\sum\limits_{n \hsx = \hsx 0}^\infty \hsx \frac{1}{(n+1)^2} 
\ = \ 
\sum\limits_{n \hsx = \hsx 1}^\infty \hsx \frac{1}{n^2} 
\ = \ 
\zeta(2).
\]
\\[-.5cm]

To establish the connection between \#3 and \#4, write

\allowdisplaybreaks
\begin{align*}
\text{\textbullet} \quad 
\ds\int\limits_0^1 
\ds\int\limits_0^1
\bigg( \frac{1}{1 - xy} - \frac{1}{1 + xy}\bigg) dx dy\
&=\ 
\ds\int\limits_0^1 
\ds\int\limits_0^1
\bigg( \frac{2 xy}{1 - x^2y^2} \bigg) dx dy 
\hspace{1.75cm}
\\[15pt]
&=\ 
\frac{1}{2} \hsx 
\ds\int\limits_0^1 
\ds\int\limits_0^1
\frac{1}{1 - xy} dx dy .
\end{align*}


\qquad \textbullet \quad
$
\ds\int\limits_0^1 
\ds\int\limits_0^1
\bigg( \frac{1}{1 - xy} + \frac{1}{1 + xy}\bigg) dx dy 
\hsx = \hsx 
2
\ds\int\limits_0^1 
\ds\int\limits_0^1
\frac{1}{1 - x^2 y^2} dx dy.\\
$
\\
\\
\vspace{0.2cm}
Then\\
\vspace{0.2cm}
\[
2
\ds\int\limits_0^1 
\ds\int\limits_0^1
\frac{1}{1 - x y} dx dy
\hsx = \hsx
\frac{1}{2}
\ds\int\limits_0^1 
\ds\int\limits_0^1
\frac{1}{1 - x y} dx dy
\hsx + \hsx
2
\ds\int\limits_0^1 
\ds\int\limits_0^1
\frac{1}{1 - x^2 y^2} dx dy
\]
\\

\qquad\qquad $\implies$
\[
2 \zeta(2) 
\hsx = \hsx
\frac{1}{2} \zeta(2) \hsx + \hsx 
2
\ds\int\limits_0^1 
\ds\int\limits_0^1
\frac{1}{1 - x^2 y^2} dx dy
\]
\\

\qquad\qquad $\implies$
\[
\frac{3}{4} \hsx \zeta(2)
\hsx = \hsx 
\ds\int\limits_0^1 \hsx
\ds\int\limits_0^1 \hsx
\frac{1}{1 - x^2 y^2} dx dy.
\]
\end{x}
\vspace{0.75cm}

\[
\textbf{APPENDIX}
\]
\vspace{0.25cm}

{\small\bf NOTATION}
\\[-.2cm]

\[
\Pi_n 
\ = \
\bigg\{(u_1, u_2, \ldots, u_n) \in \R^n: u_i > 0, \ u_i + u_{i+1} < \frac{\pi}{2} \ (1 \leq i \leq n)\bigg\}.
\]
\\[-.7cm]

[Note: \ 
In what follows the indices $i$ of the $n$ coordinates of a point in $\R^n$ are to be regarded as integers modulo $n$, thus
\[
x_i \hsx = \hsx \frac{\sin u_i}{\cos u_{i+1}} \qquad (i \in \N \ \modx \ n).]
\]
\\[-20pt]


Introduce
\[
x_1 \hsx = \hsx \frac{\sin u_1}{\cos u_2}, \ 
x_2 \hsx = \hsx \frac{\sin u_2}{\cos u_3}, \ 
\ldots, 
x_{n-1} \hsx = \hsx \frac{\sin u_{n-1}}{\cos u_n}, \ 
x_n \hsx = \hsx \frac{\sin u_n}{\cos u_1}
\]
to get an arrow $\Pi_n \ra \R^n$.
\vspace{0.5cm}

{\small\bf LEMMA 1}
The arrow $\Pi_n \ra \R^n$ is one-to-one and its range is the open unit cube $(]0,1[)^n$.
\vspace{0.5cm}

{\small\bf LEMMA 2}
The Jacobian

\[
\frac{\partial(x_1, \ldots, x_n)}{\partial(u_1, \ldots, u_n)}
\]
equals
\[
1 \ \pm \ (x_1 \cdots x_n)^2, 
\]
the sign $-$ or $+$ according to whether $n$ is even or odd.
\vspace{0.3cm}

The volume of $\Pi_n$ is 
\[
\int\limits_{\Pi_n} \hsx 1 d u_1 \cdots d u_n
\]
or still, 
\[
\int\limits_0^1 \cdots \int\limits_0^1  \ 
\frac{1}{1 \ \pm \ (x_1 \cdots x_n)^2} \hsx d x_1 \cdots d x_n
\]
or still, 
\[
\int\limits_0^1 \cdots \int\limits_0^1  \ 
\sum\limits_{k \hsx = \hsx 0}^\infty 
\hsx 
(-1)^{n k} (x_1 \cdots x_n)^{2 k} \hsx d x_1 \cdots d x_n.
\]
\vspace{0.2cm}

[Note: \ 
When $n$ is even, the integrand in the second integral is singular at 
\[
(x_1, \ldots, x_n) 
\hsx = \hsx 
(1, \ldots, 1)
\]
but the change of variables remains valid since the integrand is elsewhere positive.]
\vspace{0.2cm}

Take now $n \geq 2$ $-$then in view of absolute convergence, the third integral equals
\[
\sum\limits_{k \hsx = \hsx 0}^\infty \hsx (-1)^{nk} 
\int\limits_0^1 \hsx
\ldots 
\int\limits_0^1 \hsx
(x_1, \ldots, x_n)^{2k} \hsx dx_1 \ldots dx_n.
\]
But 
\allowdisplaybreaks
\begin{align*}
\int\limits_0^1 \hsx
\ldots 
\int\limits_0^1 \hsx
(x_1, \ldots, x_n)^{2k} \hsx dx_1 \ldots dx_n \ 
&=\ 
\bigg( \int\limits_0^1 \hsx x_1^{2k}\hsx dx_1 \bigg) 
\ldots 
\bigg( \int\limits_0^1 \hsx x_n^{2k}\hsx dx_n \bigg)  \\[15pt]
&=\ \frac{1}{(2k+1)^n}.
\end{align*}
Therefore the volume of $\Pi_n$ is 

\[
\sum\limits_{k \hsx = \hsx 0}^\infty \hsx \frac{(-1)^{nk}}{(2k + 1)^n},
\]
a rational multiple of $\pi^n$.
\vspace{0.5cm}

{\small\bf \un{N.B.}} \ 
When $n = 1$, $\Pi_n$ reduces to the line segment $0 < u_1 < \pi/4$ and the bottom line is the wellknown formula
\[
\frac{\pi}{4} 
\hsx = \hsx 
1 - \frac{1}{3} + \frac{1}{5} - \frac{1}{7} + \cdots, 
\]
the value of 
\[
\int\limits_0^1 \hsx \frac{1}{1 + x^2} \hsx dx.
\]
\vspace{0.3cm}

{\small\bf REMARK} Take $n$ even $-$then 
\[
\sum\limits_{k \hsx = \hsx 0}^\infty \hsx \frac{(-1)^{nk}}{(2k + 1)^n} 
\ = \ 
(1 - 2^{-n}) \hsx \zeta(n).
\]


%% file: _04_Zeta_2_bis.tex
\chapter{
$\boldsymbol{\S}$\textbf{4}.\quad  $\boldsymbol{\zeta(2)}$\  (bis)}
\setlength\parindent{2em}
\setcounter{theoremn}{0}
\renewcommand{\thepage}{Zeta Function Values \S4-\arabic{page}}

\begin{spacing}{1.5}
\ \indent 
Since $\zeta(2) = \ds\frac{\pi^2}{6}$, it follows that $\zeta(2)$ is transcendental, hence irrational.
But let's ignore this, the objective being to prove from first principles that $\zeta(2)$ is irrational, the point being that the methods utilized can be extended in the next \S \  to establish that $\zeta(3)$ is irrational.
\end{spacing}
\vspace{0.5cm}

\begin{x}{\small\bf NOTATION} \ 
Let $d_n$ be the least common multiple of $1, 2, \ldots, n$ and set $d_0 = 1$.
\end{x}
\vspace{0.3cm}

\begin{x}{\small\bf LEMMA} \ 
$\forall \ K > e$, 
\[
d_n < K^n \qquad \text{if} \ n \gg 0.
\]

PROOF \ 
\allowdisplaybreaks
\begin{align*}
d_n \ 
&=\ \prod\limits_{p \hsx  \leq \hsx n} \hsx p^{[\elln(n) / \elln(p)]}\\[15pt]
&\leq\ \prod\limits_{p \hsx  \leq \hsx n} \hsx p^{\elln(n) / \elln(p)}\\[15pt]
&=\ \prod\limits_{p \hsx  \leq \hsx n} \hsx n\\[15pt]
&=\ n^{\pi(n)},
\end{align*}
$\pi(n)$ the prime counting function.  
Owing to the prime number theorem, 
\[
\lim\limits_{n \ra \infty} \hsx \frac{\pi(n) \hsy \elln(n)}{n} \hsx = 1,
\]
so if $A > 1$, then 
\[
n \gg 0 \ \implies  \frac{\pi(n) \hsy \elln(n)}{n} < A
\]
\end{x}
\vspace{0.3cm}
or still, 
\allowdisplaybreaks
\begin{align*}
n \gg 0 \ 
&\implies \pi(n) \elln(n) < nA\\[12pt]
&\implies n^{\pi(n)} < (e^A)^n = K^n,
\end{align*}
where $K = e^A > e$, i.e., 
\[
n \gg 0 \ \implies \ d_n = n^{\pi(n)}  < K^n.
\]

\begin{x}{\small\bf \un{N.B.}} \ 
In particular, 
\[
n \gg 0 \ \implies \ d_n < 3^n.
\]
\end{x}
\vspace{0.3cm}

\begin{x}{\small\bf NOTATION} \ 
Let
\[
P_n(x) 
\ = \ 
 \frac{1}{n!} \hsx \frac{d^n}{dx^n} \bigl(x^n (1 - x)^n \bigr). 
\]
Then 
\[
P_n(x) 
\hsx = \hsx 
\sum\limits_{k=0}^n \hsx (-1)^k {n\choose k} {n+k \choose k} 
\hsx x^k,
\]
a polynomial of degree $n$ with integral coefficients.
\end{x}
\vspace{0.3cm}

\begin{x}{\small\bf SUBLEMMA} \ 
For $i \leq n - 1$, 
\[
\begin{cases}
\ \ds\frac{d^i}{dx^i} \hsx \bigl(x^n (1 - x)^n \bigr)(0) \hsx = \hsx 0\\[15pt]
\ \ds\frac{d^i}{dx^i} \hsx \bigl(x^n (1 - x)^n \bigr)(1) \hsx = \hsx 0
\end{cases}
.
\]
\end{x}
\vspace{0.3cm}

\begin{x}{\small\bf LEMMA} \ 
Suppose that $f(x)$ is sufficiently differentiable $-$then
\[
\abs{\hsx \int\limits_0^1 \hsx P_n(x) f(x) dx \hsx} 
\  = \  
\abs{\hsx \int\limits_0^1 \hsx \frac{1}{n!} \hsx x^n (1 - x)^n \frac{d^n}{dx^n} f(x) dx \hsx}.
\]
\vspace{0.2cm}

PROOF \ 
Write
\allowdisplaybreaks
\begin{align*}
\int\limits_0^1 \hsx P_n(x) f(x) dx \ 
&=\ \int\limits_0^1 \hsx \frac{1}{n!} \frac{d^n}{dx^n} (x^n (1 - x)^n ) f(x) dx\\[15pt]
&=\ \frac{1}{n!} 
\frac{d^{n-1}}{dx^{n-1}} (x^n (1 - x)^n ) f(x) \bigg\arrowvert_0^1
\hsx - \hsx 
\int\limits_0^1 \hsx \frac{1}{n!} \frac{d^{n-1}}{dx^{n-1}} (x^n (1 - x)^n ) 
\hsx 
\frac{d}{dx}f(x) dx\\[15pt]
&=\ - \hsx \int\limits_0^1 \hsx \frac{1}{n!} \frac{d^{n-1}}{dx^{n-1}} (x^n (1 - x)^n ) 
\hsx 
\frac{d}{dx}f(x) dx.
\end{align*}
Proceed from here by iteration.
\end{x}
\vspace{0.3cm}

\begin{x}{\small\bf INTEGRAL FORMULAS} \ 
\vspace{0.5cm}

\qquad \textbullet \quad Let $r$ be a nonnegative integer $-$then 

\[
\int\limits_0^1 \hsx 
\int\limits_0^1 \hsx 
\frac{x^r y^r}{1 - xy} \hsx dx dy
\ = \ 
\sum\limits_{n=1}^\infty \frac{1}{(n+r)^2}.
\]
So
\allowdisplaybreaks
\begin{align*}
r = 0 \ 
&\implies 
\int\limits_0^1 \hsx 
\int\limits_0^1 \hsx 
\frac{1}{1 - xy} \hsx dx dy \ = \ \zeta(2) \qquad \text{(cf. \S3, \#4)}.\\[15pt]
r > 0 \ 
&\implies 
\int\limits_0^1 \hsx 
\int\limits_0^1 \hsx 
\frac{x^r y^r}{1 - xy} \hsx dx dy 
\  =\  \zeta(2) - \bigg( \frac{1}{1^2} + \frac{1}{2^2} + \cdots + \frac{1}{r^2} \bigg).
\end{align*}
\vspace{0.2cm}

\qquad \textbullet \quad Let $r, s$ be a nonnegative integers with $r > s$$-$then 

\[
\int\limits_0^1 \hsx 
\int\limits_0^1 \hsx 
\frac{x^r y^s}{1 - xy} \hsx dx dy
\hsx = \hsx 
\frac{1}{r - s} \hsx \bigg\{\frac{1}{s+1} + \frac{1}{s+2} + \cdots + \frac{1}{r}\bigg\}.
\]
\end{x}
\vspace{0.3cm}

\begin{x}{\small\bf APPLICATION} \ 
\[
\int\limits_0^1 \hsx 
\int\limits_0^1 \hsx 
\frac{x^r y^r}{1 - xy} \hsx dx dy
\hsx = \hsx 
\zeta(2) - \frac{a}{d_r^{\raisebox{1.5pt}{$\scriptstyle 2$}}}
\]
and
\[
\int\limits_0^1 \hsx 
\int\limits_0^1 \hsx 
\frac{x^r y^s}{1 - xy} \hsx dx dy
\hsx = \hsx 
\frac{b}{d_r^{\raisebox{1.5pt}{$\scriptstyle 2$}}},
\]
where $a, b$ are integers.
\end{x}
\vspace{0.3cm}

Therefore:
\vspace{0.3cm}

\begin{x}{\small\bf LEMMA} \ 
If $P(x)$, $Q(y)$ are polynomials of degree $n$ with integer coefficients, then 
\[
\int\limits_0^1 \hsx 
\int\limits_0^1 \hsx 
\frac{P(x) Q(y)}{1 - xy} \hsx dx dy
\ = \
\frac{A \zeta(2) + B}{d_n^{\raisebox{1.5pt}{$\scriptstyle 2$}}},
\]
where $A$, $B$ are integers.
\end{x}
\vspace{0.3cm}

\begin{x}{\small\bf NOTATION} \ 
Put 
\[
I_n 
\ = \ 
\int\limits_0^1 \hsx 
\int\limits_0^1 \hsx 
\frac{P_n(x) (1 - y)^n}{1 - xy} \hsx dx dy.
\]
Take $Q(y) = (1 - y)^n$ to get 
\[
I_n
\hsx = \hsx 
\frac{A_n \zeta(2) + B_n}{d_n^2}, 
\]
where $A_n$, $B_n$ are integers depending on $n$.
\end{x}
\vspace{0.3cm}


\begin{x}{\small\bf LEMMA} \ 
\[
\abs{I_n} 
\ = \ 
\int\limits_0^1 \hsx 
\int\limits_0^1 \hsx 
\frac{x^n (1 - x)^n y^n (1 - y)^n}{(1 - xy)^{n+1}} \hsx dx dy.
\]
\vspace{0.2cm}

PROOF \ Taking into account \#6, 

\allowdisplaybreaks
\begin{align*}
\abs{I_n}  \ 
&=\ 
\abs{\hsx
\int\limits_0^1 \hsx \frac{x^n (1 - x)^n}{n!} 
\frac{d^n}{dx^n}
\bigg(
\int\limits_0^1 \hsx \frac{(1 - y)^n}{1 - xy}  
\hsx dy \bigg) dx \hsx}\\[15pt]
&=\ 
\abs{\hsx
\int\limits_0^1 \hsx \frac{x^n (1 - x)^n}{n!} 
\bigg(
\int\limits_0^1 \hsx \frac{d^n}{dx^n} \hsx \bigg( \frac{(1 - y)^n}{1 - xy}  
\bigg)
\hsx dy \bigg) dx \hsx}\\[15pt]
&=\ 
\abs{\hsx
\int\limits_0^1 \hsx \frac{x^n (1 - x)^n}{n!} 
\bigg(
\int\limits_0^1 \hsx  \frac{n! \hsy y^n(1 - y)^n}{(1 - xy)^{n+1}\hsx}  
\hsx dy \bigg) dx \hsx}\\[15pt]
&=\ 
\int\limits_0^1 \hsx
\int\limits_0^1 \hsx
\frac{x^n (1 - x)^n y^n (1 - y)^n}{(1 - xy)^{n+1}} 
\hsx dx dy.
\end{align*}
\end{x}
\vspace{0.3cm}

\begin{x}{\small\bf \un{N.B.}} \ 
$I_n$ is nonzero (the integrand is positive for all $x, y \in ]0,1[\hsx)$.
\end{x}
\vspace{0.3cm}

The function
\[
f(x,y) 
\hsx = \hsx 
\frac{x (1 - x) \hsy y \hsy (1 - y)}{1 - xy} 
\qquad (0 \leq x < 1, \ 0 \leq y < 1)
\]
vanishes on the boundary of $[0,1] \times [0,1]$ and, although not defined at $(1,1)$, it does however tend to 0 as 
$x, y \uparrow 1$.
\vspace{0.3cm}

\begin{x}{\small\bf LEMMA} \ 
The maximum of $f(x,y)$ in $0 < x < 1$, $0 < y < 1$ is 
\[
\bigg(\frac{\sqrt{5} \hsx -1}{2}\bigg)^5.
\]
\vspace{0.2cm}

PROOF \ 
Consider the relations
\[
\frac{\partial}{\partial x} f(x,y) \hsx = \hsx 0,
\quad 
\frac{\partial}{\partial y} f(x,y) \hsx = \hsx 0,
\]
i.e., 
\[
1 - 2x + yx^2 \hsx = \hsx 0, \quad 1 - 2y + xy^2 \hsx = \hsx 0.
\]
Then
\allowdisplaybreaks
\begin{align*}
y \hsx = \hsx \frac{2x - 1}{x^2} 
&\implies 
1 - 2\bigg( \frac{2x - 1}{x^2} \bigg) + x \bigg( \frac{2x - 1}{x^2} \bigg)^2 
\hsx = \hsx
0\\[15pt]
&\implies x^3 - 2x + 1 
\hsx = \hsx
0,
\end{align*}
the roots of which are 
\[
1, \ \frac{-1 \pm \sqrt{5}}{2}, 
\quad \text{so} \quad 
x 
\hsx = \hsx
\frac{\sqrt{5} \hsx -1}{2}.
\]
Analogously
\[
y 
\hsx = \hsx
\frac{\sqrt{5} \hsx -1}{2}.
\]
Therefore $f(x,y)$ achieves its maximum at 
\[
\bigg(\frac{\sqrt{5} \hsx -1}{2}, \frac{\sqrt{5} \hsx -1}{2}\bigg),
\]
the value being
\[
\bigg(\frac{\sqrt{5} \hsx -1}{2}\bigg)^5.
\]
\end{x}
\vspace{0.3cm}

\begin{x}{\small\bf APPLICATION} \ 
\allowdisplaybreaks
\begin{align*}
\abs{I_n} \ 
&=\
\int\limits_0^1 \hsx 
\int\limits_0^1 \hsx 
f(x,y) \hsx \frac{1}{1 - xy} \hsx dx dy\\[15pt]
&\leq\ \bigg(\frac{\sqrt{5} \hsx - 1}{2} \bigg)^{5n} \hsx 
\int\limits_0^1 \hsx 
\int\limits_0^1 \hsx 
\frac{1}{1 - xy} \hsx dx dy\\[15pt]
&=\ \bigg(\frac{\sqrt{5} \hsx - 1}{2} \bigg)^{5n} \hsx \zeta(2) \qquad \text{(cf. \S3, \#4)}.
\end{align*}
\end{x}
\vspace{0.3cm}


\begin{x}{\small\bf \un{N.B.}} \ 
\[
\frac{\sqrt{5} \hsx - 1}{2}
\hsx < \hsx
0.62,
\]
And
\[
(0.62)^5 \ < \  \frac{1}{10} 
\]

\qquad$\implies $
\allowdisplaybreaks
\begin{align*}9^n \hsx \bigg(\frac{\sqrt{5} \hsx - 1}{2} \bigg)^{5n}
&=\ \bigg(9 \hsx \cdot \hsx \bigg(\frac{\sqrt{5} \hsx - 1}{2} \bigg)^5 \bigg)^n\\[15pt]
&<\ \bigg(9 \hsx \cdot \hsx \frac{1}{10} \bigg)^n\\[15pt] 
&=\  \bigg(\frac{9}{10}\bigg)^n  
\\[15pt] 
&\ra \  0 
\qquad\qquad (n \ra \infty).
\end{align*}
\end{x}
\vspace{0.3cm}

\begin{x}{\small\bf THEOREM} \ 
$\zeta(2)$ is irrational.
\vspace{0.2cm}

PROOF \ 
Suppose instead that $\zeta(2)$ is rational, say $\zeta(2) = \ds\frac{a}{b}$ $(a,\ b\in \N)$.  
Write
\allowdisplaybreaks
\begin{align*}
I_n \ 
&=\ \frac{A_n \zeta(2) + B_n}{d_n^2} \qquad \text{(cf. \#10)}\\[15pt]
&=\ \frac{A_n \bigg(\ds\frac{a}{b}\bigg) + B_n}{d_n^2}
\end{align*}
\qquad\qquad $\implies$
\[
\abs{A_n \bigg(\frac{a}{b}\bigg) + B_n} \hsx \leq \hsx d_n^2 \hsx \abs{I_n}
\]
\qquad\qquad $\implies$ $(n \gg 0)$
\allowdisplaybreaks
\begin{align*}
\abs{A_n \bigg(\frac{a}{b}\bigg) + B_n}\ 
&\leq \ 9^n \hsx \abs{I_n} \qquad \text{(cf. \#3)}\\[15pt]
&\leq \ 9^n \hsx \bigg(\frac{\sqrt{5} \hsx - 1}{2} \bigg)^{5n} \hsx \zeta(2)
\end{align*}
\qquad\qquad $\implies$ $(n \gg 0)$
\allowdisplaybreaks
\begin{align*}
\abs{A_n a + B_n b}\ 
&\leq \ 9^n \hsx \bigg(\frac{\sqrt{5} \hsx - 1}{2} \bigg)^{5n} \hsx b\\[15pt]
&\approx\ b\bigg(\frac{9}{10}\bigg)^n \ 
\\[15pt]
&\ra \ 0.
\end{align*}
But $I_n$ is nonzero (cf. \#12), hence
\[
0 
\hsx < \hsx 
\abs{A_n a + B_n b} \ra 0 \qquad (n \ra \infty),
\]
a contradiction (a sequence of positive integers cannot tend to 0).
\end{x}
\vspace{0.3cm}


%% file: _05_Zeta_2n_bis.tex
\chapter{
$\boldsymbol{\S}$\textbf{5}.\quad   $\boldsymbol{\zeta(2 n)}$\  (bis)}
\setlength\parindent{2em}
\setcounter{theoremn}{0}
\renewcommand{\thepage}{Zeta Function Values \S5-\arabic{page}}


\begin{x}{\small\bf RAPPEL} \ 
\[
\pi x \cot(\pi x) 
\hsx = \hsx 
1 + 2 x^2 \hsx \sum\limits_{k \hsx = \hsx 1}^\infty \hsx \frac{1}{x^2 - k^2}.
\]
\end{x}
\vspace{0.3cm}

\begin{x}{\small\bf RAPPEL} \ 
\[
\pi x \cot(\pi x) 
\hsx = \hsx 
1 + \sum\limits_{n \hsx = \hsx 1}^\infty \hsx 
(-1)^n \hsx \frac{2^{2n} B_{2n}}{(2n)!}  \hsx\pi^{2n} \hsx x^{2n}.
\]
\end{x}
\vspace{0.3cm}

\begin{x}{\small\bf \un{N.B.}} \ 
These expansions are valid for $\abs{x}$ sufficiently small.
\vspace{0.3cm}

Given $k$, expand
\[
\frac{2 x ^2}{x^2 - k^2}
\]
in powers of $x$:
\[
\frac{2 x ^2}{x^2 - k^2}
\hsx = \hsx 
-2 \hsx \sum\limits_{n \hsx = \hsx 1}^\infty \hsx \bigg(\frac{x^2}{k^2}\bigg)^n.
\]
Therefore the coefficient of $x^{2n}$ is
\[
-2 \hsx \sum\limits_{n \hsx = \hsx 1}^\infty \hsx \frac{1}{k^{2n}}.
\]
And then
\allowdisplaybreaks
\begin{align*}
\sum\limits_{k \hsx = \hsx 1}^\infty \hsx \frac{2x^2}{x^2 - k^2} \ 
&=\ 
-2 \hsx 
\sum\limits_{k \hsx = \hsx 1}^\infty \hsx 
\sum\limits_{n \hsx = \hsx 1}^\infty \hsx 
\bigg(\frac{x^2}{k^2}\bigg)^n
\\[15pt]
&=\ 
-2 \hsx 
\sum\limits_{n \hsx = \hsx 1}^\infty \hsx 
\sum\limits_{k \hsx = \hsx 1}^\infty \hsx 
\frac{x^{2n}}{k^{2n}}
\\[15pt]
&=\ 
-2 \hsx 
\sum\limits_{n \hsx = \hsx 1}^\infty \hsx 
\bigg( \sum\limits_{k \hsx = \hsx 1}^\infty \hsx 
\frac{1}{k^{2n}}\bigg) x^{2n}
\\[15pt]
&=\ -2 \hsx 
\sum\limits_{n \hsx = \hsx 1}^\infty \hsx  
\zeta(2 n) x^{2n},
\end{align*}
i.e.,  $-2 \zeta(2n)$ is the coefficient of $x^{2n}$.  
But the coefficient of $x^{2n}$ is also 
\[
(-1)^{n} \hsx \frac{2^{2n} \hsx B_{2n}}{(2n)!} \hsx \pi^{2n}.
\]
Consequently 

\[
\zeta(2n)
\hsx = \hsx 
(-1)^{n-1} \hsx \frac{2^{2n-1}}{(2n)!} \hsx B_{2n} \hsy \pi^{2n}
\]
\\
as predicted by the considerations of \#2.
\end{x}
\vspace{0.3cm}


%% file: _06_Zeta_3.tex
\chapter{
$\boldsymbol{\S}$\textbf{6}.\quad  $\boldsymbol{\zeta(3)}$}
\setlength\parindent{2em}
\setcounter{theoremn}{0}
\renewcommand{\thepage}{Zeta Function Values \S6-\arabic{page}}


\begin{x}{\small\bf THEOREM} \ 
$\zeta(3)$ is irrational.
\vspace{0.5cm}

The proof is similar to that for $\zeta(2)$ (cf. \S4, \#16), albeit technically more complicated.  
In outline form, here is how it goes.
\vspace{0.5cm}

\qquad \un{Step \ 1}:
\vspace{0.3cm}

\qquad\textbullet \quad Let $r$ be a nonnegative integer $-$then
\[
-
\int\limits_0^1
\int\limits_0^1 \hsx
\frac{\ell n(xy)}{1 - xy} x^r y^r dx dy 
\  = \ 
2\bigg(\zeta(3) - \sum\limits_{k = 1}^r \frac{1}{k^3}\bigg) \ \in \ 2 \zeta(3) + \frac{1}{d_r^{\raisebox{1.5pt}{$\scriptstyle 3$}}} \Z.
\]
In particular: 
\[
-
\int\limits_0^1
\int\limits_0^1 \hsx
\frac{\ell n(xy)}{1 - xy} dx dy 
\  = \  
2 \hsy \zeta(3).
\]
\\[-15pt]

\qquad\textbullet \quad Let $r$, $s$ be nonnegative integers with $r > s$ $-$then
\[
-
\int\limits_0^1
\int\limits_0^1 \hsx
\frac{\ell n(xy)}{1 - xy} x^r y^s dx dy 
\hsx = \hsx 
\frac{1}{r - s} \hsx \bigg( \frac{1}{(s + 1)^2} + \cdots + \frac{1}{r^2} \bigg) \ \in \ \frac{1}{d_r^{\raisebox{1.5pt}{$\scriptstyle 3$}}} \Z.
\]
\\[-15pt]

\qquad \un{Step \ 2}:
\allowdisplaybreaks
\begin{align*}
I_n \ 
&\equiv\  
-
\int\limits_0^1
\int\limits_0^1 \hsx
\frac{P_n(x) P_n(y)}{1 - xy} \hsx \ell n (xy) \hsx dx dy\\[15pt]
&=\ \frac{A_n \zeta(3) + B_n}{d_n^{\raisebox{1.5pt}{$\scriptstyle 3$}}},
\end{align*}
where $A_n$, $B_n \in \Z$.
\vspace{0.3cm}

\qquad \un{Step \ 3}:
\[
-\frac{\ell n (xy)}{1 - xy} 
\  = \  
\int\limits_0^1 \hsx \frac{1}{1 - (1 - xy)z} \hsx dz.
\]
\vspace{0.3cm}

\qquad \un{Step \ 4}:
\allowdisplaybreaks
\begin{align*}
\abs{I_n} \ 
&=\ \abs{\hsx
\int\limits_0^1 \hsx
\int\limits_0^1 \hsx
\int\limits_0^1 \hsx
\frac{P_n(x)P_n(y)}{1 - (1 - xy)z} 
\hsx dz dx dy \hsx}\\[15pt]
&=\ \abs{\hsx
\int\limits_0^1 \hsx \frac{x^n(1 - x)^n}{n!} \hsx \frac{d^n}{dx^n}
\bigg(
\int\limits_0^1 \hsx
\int\limits_0^1 \hsx
\frac{P_n(y)}{1 - (1 - xy) z} \hsx dy dz \bigg) dx \hsx} \\[15pt]
&\hspace{0.15cm} \vdots\\[15pt]
&=\ \abs{\hsx
\int\limits_0^1 \hsx P_n(y) 
\bigg(
\int\limits_0^1 \hsx
\int\limits_0^1 \hsx
\frac{x^n (1 - x)^n y^n z^n}{(1 - (1 - xy) z)^{n+1}} 
\hsx dx dz \bigg) dy\hsx}.
\end{align*}
\vspace{0.3cm}

\qquad \un{Step \ 5}:
Let $D = \{(u, v, w): u, v, w \in ]0,1[\}$ $-$then the map 
\[
(u, v, w) \ra (x, y, z)
\]
defined by $x = u$, $y = v$ and 
\[
z 
\hsx = \hsx 
\frac{1 - w}{1 - (1 - uv) w}
\]
from \mD to \mD is one-to-one and onto.  
In addition, 
\[
\frac{\partial (x, y, z)}{\partial(u, v, w)}
\hsx = \hsx 
-\frac{uv}{(1 - (1 - uv) w)^2}.
\]
\vspace{0.3cm}

\qquad \un{Step \ 6}:
The function 
\[
\frac{u (1 - u) v (1 - v) w (1 - w)}{1 - (1 - uv) w}
\]
is bounded above by $\ds\frac{1}{27}$ in the region \mD.
\vspace{0.3cm}

\qquad \un{Step \ 7}:
In $I_n$ make a change of variable and use the relations
\[
z^n 
\hsx = \hsx 
\frac{(1 - w)^n}{(1 - (1 - uv) w)^n}
\]
\allowdisplaybreaks
\begin{align*}
(1 - (1 - xy) z)^{n+1} \ 
&=\ \bigg(1 - (1 - u v) \hsx \frac{1 - w}{1 - (1 - uv)w}\bigg)^{n+1}
\\[15pt]
&=\ \frac{(uv)^{n+1}}{(1 - (1 - uv) w)^{n+1}}
\end{align*}
to get 
\[
\abs{I_n} 
\ = \ 
\int\limits_0^1
\int\limits_0^1
\int\limits_0^1
\frac{u^n (1 - u)^n v^n (1 - v)^n w^n (1 - w)^n}{(1 - (1 - uv) w)^{n+1}}
\hsx 
du dv dw.
\]
\vspace{0.3cm}

\qquad \un{Step \ 8}:\ 
Therefore
\allowdisplaybreaks
\begin{align*}
0 \
&<\  \abs{I_n} \\
&\leq\ \bigg(\frac{1}{27}\bigg)^n \hsx
\int\limits_0^1
\int\limits_0^1
\int\limits_0^1
\hsx
\frac{1}{1 - (1 - uv)w} 
\hsx du dv dw
\\[15pt]
&=\ \bigg(\frac{1}{27}\bigg)^n \hsx 
\int\limits_0^1
\int\limits_0^1
\hsx 
-\frac{\ell n (uv)}{1 - uv} 
\hsx du dv
\\[15pt]
&=\ 2 \bigg(\frac{1}{27}\bigg)^n \hsx \zeta(3).
\end{align*}
\vspace{0.3cm}

\qquad \un{Step \ 9}:\ 
\allowdisplaybreaks
\begin{align*}
0 \
&<\  \abs{I_n} 
\\[15pt]
&=\ \frac{\abs{A_n \zeta(3) + B_n}}{d_n^{\raisebox{1.5pt}{$\scriptstyle 3$}}} 
\\[15pt]
&\leq\ 2 \bigg(\frac{1}{27}\bigg)^n \hsx \zeta(3).
\end{align*}
\vspace{0.3cm}

\qquad \un{Step \ 10}: \ 
To derive a contradiction, suppose that $\zeta(3)$ is rational, say $\zeta(3) = \ds\frac{a}{b}$ $(a, b \in \N)$ $-$then 
\[
0 
\hsx < \hsx 
\abs{A_n \bigg(\frac{a}{b}\bigg) + B_n}
\hsx \leq \hsx 
2 \bigg(\frac{1}{27}\bigg)^n \hsx \zeta(3) d_n^3
\]

$\implies$
\allowdisplaybreaks
\begin{align*}
0 \ 
&<\ \abs{A_n a + B_n b} 
\\[15pt]
&\leq\ 2 b \bigg(\frac{1}{27}\bigg)^n \hsx  d_n^3
\\[15pt]
&<\ 2 b \bigg(\frac{1}{27}\bigg)^n (2.8)^{3n} \qquad (\text{cf. \S4, \#2 (take $K = 2.8$)})
\\[15pt]
&=\ 2 b \bigg(\frac{(2.8)^3}{27}\bigg)^n
\\[15pt]
&<\ 2 b (0.9)^n 
\\[15pt]
&\ra 0 \qquad (n \ra \infty).
\end{align*}
\end{x}
\vspace{0.3cm}

\begin{x}{\small\bf \un{N.B.}} \ 
The irrationality of $\zeta(3)$ is thereby established but the issue of its transcendence remains open.
\end{x}
\vspace{0.3cm}

\begin{x}{\small\bf REMARK} \ 
It was shown by T. Rivoal that the $\Q$-vector space generated by 
\[
1,\  \zeta(3),\  \zeta(5),\  \zeta(7), \ldots
\]
is infinite dimensional, hence there exist infinitely many $n$ such that $\zeta(2n+1)$ is irrational (but it is unknown whether 
$\zeta(5)$ is irrational).
\vspace{0.2cm}

[Note: \ 
For an account, consult S. Fischler (arXiv:math.0303066).]
\end{x}
\vspace{0.3cm}

In the book ``Zeta and $q$-Zeta Functions and Associated Series and Integrals'' by H.M. Srivastava and Junesang Choi, the reader will find a large collection of formulas for $\zeta(2n +1)$.


%% file: _07_Conjugate_Bernoulli_Numbers.tex
\chapter{
$\boldsymbol{\S}$\textbf{7}.\quad  CONJUGATE BERNOULLI NUMBERS}
\setlength\parindent{2em}
\setcounter{theoremn}{0}
\renewcommand{\thepage}{Zeta Function Values \S7-\arabic{page}}


\begin{x}{\small\bf DEFINITION} \ 
If $f$ is a 1-periodic function, then its 
\un{periodic Hilbert transform}
\index{periodic Hilbert transform}
$\sH[f]$ is given by
\[
\sH[f] (x)
\hsx = \hsx 
\PV \hsx \int\limits_{-1/2}^{1/2} \hsx f(x - y)  \cot(\pi y) dy.
\]
\end{x}
\vspace{0.3cm}

\begin{x}{\small\bf CONSTRUCTION} \ 
Start with the Bernoulli polynomial $B_n(x)$ and put
\[
\sB_n(x)
\hsx = \hsx 
B_n(x - [x]),
\]
a so-called 
\un{Bernoulli function}.  
\index{Bernoulli function}
It is 1-periodic and 
\[
\frac{\sB_n(x)}{n!} 
\hsx = \hsx 
- \hsx \sum\limits_{\substack{k \in \Z\\ k \neq 0}} \hsx \frac{e^{2 \pi \sqrt{-1} \hsx k x}}{(2 \pi \sqrt{-1} \hsx k)^n},
\]
a formula which holds for all real $x$ if $n \geq 2$ and for all $x \notin \Z$ if $n = 1$.
\end{x}
\vspace{0.3cm}

\begin{x}{\small\bf DEFINITION} \ 
The 
\un{conjugate Bernoulli functions}
\index{conjugate  Bernoulli functions}  \ 
$\widetilde{B}_n(x)$ are defined for $x \in [0,1[$ $(x \neq 0$ if $n = 1$) by the restriction of $\sH[\sB_n]$ to $[0,1[$.
\end{x}
\vspace{0.3cm}

\begin{x}{\small\bf EXAMPLE} \ 
For $0 < x < 1$, 
\[
\widetilde{B}_1(x)
\  = \  
-\frac{1}{\pi} \hsx \elln (2 \sin(\pi x)).
\]
\end{x}
\vspace{0.3cm}

\begin{x}{\small\bf EXAMPLE} \ 
\allowdisplaybreaks
\begin{align*}
\widetilde{B}_{2n+1} \bigg(\frac{1}{2}\bigg) \ 
&=\ \sH[\sB_{2n+1}]  \bigg(\frac{1}{2}\bigg)\\
&=\ \PV \hsx \int\limits_{-1/2}^{1/2} \hsx \sB_{2n+1}\bigg(\frac{1}{2} - y\bigg)  \cot(\pi y) dy \\
&=\ \PV \hsx \int\limits_{-1/2}^{1/2} \hsx B_{2n+1}\bigg(\frac{1}{2} - y\bigg) \cot(\pi y) dy.
\end{align*}
\vspace{0.2cm}

[Note: \ By definition
\[
\sB_{2n+1}\bigg(\frac{1}{2} - y\bigg)
\hsx = \hsx 
B_{2n+1}\bigg(\frac{1}{2} - y - \bigg[ \frac{1}{2} - y \bigg]\bigg) .
\]
But
\allowdisplaybreaks
\begin{align*}
-\frac{1}{2} < y < \frac{1}{2} 
&\implies \frac{1}{2} > -y > -\frac{1}{2}
\\[12pt]
&\implies \frac{1}{2} + \frac{1}{2} > \frac{1}{2} - y > \frac{1}{2} - \frac{1}{2} 
\\[12pt]
&\implies 1 > \frac{1}{2} - y > 0
\\[12pt]
&\implies \bigg[\frac{1}{2} - y\bigg] = 0 \hsx.]
\end{align*}
\end{x}
\vspace{0.3cm}

\begin{x}{\small\bf \un{N.B.}} \ 
\[
\widetilde{B}_n(x) 
\ = \ 
-2(n!) \hsx \sum\limits_{k \hsx = \hsx 1}^\infty \hsx 
\frac{\sin(2 \pi k x - n \pi / 2)}{(2 \pi k)\strutv^n \upx} 
\qquad (x \neq 0 \ \text{if} \ n = 1).
\]
\end{x}
\vspace{0.3cm}

\begin{x}{\small\bf LEMMA} \ 
$\forall \ n \in \N$,
\[
\widetilde{B}_n(1-x) 
\hsx = \hsx 
(-1)^{n+1} \widetilde{B}_n(x) \qquad (0 < x < 1).
\]
\vspace{0.2cm}

PROOF \ From \#6,
\[
\widetilde{B}_n(1-x) 
\hsx = \hsx 
-2(n!) \hsx 
\sum\limits_{k \hsx = \hsx 1}^\infty \hsx 
\frac{\sin(2 \pi k (1-x)  - n \pi / 2)}{(2 \pi k)^n}.
\]
Write
\allowdisplaybreaks
\begin{align*}
\sin(2 \pi k (1 - x) - n \pi/2) \ 
&=\ \sin(2 \pi k - 2 \pi k x - n \pi / 2 + n \pi / 2 - n \pi / 2)
\\[12pt]
&=\ \sin((-2 \pi k x + n \pi / 2) + (2 \pi k - n \pi))
\\[12pt]
&=\ \sin(-2 \pi k x + n \pi / 2) \cos(2 \pi k - n \pi ) 
\\[12pt]
& \hspace{2.5cm} + \hsx
\sin(2 \pi k - n \pi) \cos(-2 \pi k x + n \pi /2)
\\[12pt]
&=\ -\sin(2 \pi k x - n \pi/2) \cos(-n \pi) 
\hsx + \hsx 
\sin(-n\pi) \cos(-2\pi k x + n \pi /2)
\\[12pt]
&=\ \sin(2 \pi k x - n \pi/2) (-1) \cos(n \pi)
\hsx + \hsx 
(0) \cos(-2 \pi k x + n \pi /2)
\\[12pt]
&=\ \sin(2 \pi k x - n \pi /2) (-1) (-1)^n
\\[12pt]
&=\ (-1)^{n+1} \hsx \sin(2 \pi k x - n \pi /2),
\end{align*}
matters then being manifest.
\end{x}
\vspace{0.3cm}

\begin{x}{\small\bf APPLICATION} \ 
Take $x = \ds\frac{1}{2}$ $-$then
\[
\widetilde{B}_{2n}\bigg(\frac{1}{2}\bigg) \ 
=\ (-1)^{2n+1} \hsx \widetilde{B}_{2n}\bigg(\frac{1}{2}\bigg) \ = \ -\widetilde{B}_{2n}\bigg(\frac{1}{2}\bigg)
\]
\qquad\qquad$\implies$
\[
\widetilde{B}_{2n}\bigg(\frac{1}{2}\bigg) = 0.
\]
\end{x}
\vspace{0.3cm}

\begin{x}{\small\bf DEFINITION} \ 
The 
\un{conjugate Bernoulli numbers}
\index{conjugate Bernoulli numbers}  
$\widetilde{B}_n$ are defined by
\[
\widetilde{B}_n 
\hsx = \hsx 
\widetilde{B}_n(0) \qquad (n > 1).
\]
\end{x}
\vspace{0.3cm}

\begin{x}{\small\bf RAPPEL} \ 
$\forall \ n > 1$, 
\[
\sum\limits_{k \hsx = \hsx 1}^\infty \hsx \frac{(-1)^{k+1}}{(2 \pi k)^n} 
\hsx = \hsx 
(2 \pi)^{-n} (1 - 2^{1 - n}) \hsx \zeta(n).
\]
\end{x}
\vspace{0.3cm}

\begin{x}{\small\bf LEMMA} \ 
$\forall \ n > 1$, 
\[
\widetilde{B}_n\bigg(\frac{1}{2}\bigg)
\ = \ 
\bigl( 2^{1 - n} - 1 \bigr) \widetilde{B}_n.
\]
\vspace{0.2cm}

PROOF \ From \#6,
\[
\widetilde{B}_n\bigg(\frac{1}{2}\bigg)
\ = \ 
-2(n!) \hsx \sum\limits_{k \hsx = \hsx 1}^\infty \hsx 
\frac{\sin(\pi k - n \pi / 2)}{(2 \pi k)^n}.
\]
But
\allowdisplaybreaks
\begin{align*}
\sin(\pi k - n \pi /2) \ 
&=\ \sin(\pi k) \cos \bigg(\frac{n \pi}{2}\bigg) - \sin\bigg(\frac{n \pi}{2}\bigg) \cos(\pi k)
\\[12pt]
&=\ - \sin\bigg(\frac{n \pi}{2}\bigg) \cos(\pi k)
\\[12pt]
&=\ - \sin\bigg(\frac{n \pi}{2}\bigg) (-1)^{k}
\\[12pt]
&=\ \sin\bigg(\frac{n \pi}{2}\bigg) (-1)^{k+1}.
\end{align*}
Therefore
\allowdisplaybreaks
\begin{align*}
\widetilde{B}_n \bigg(\frac{1}{2}\bigg)  \ 
&=\ -2(n!) \hsx\sin\bigg(\frac{n \pi}{2}\bigg) \hsx 
\sum\limits_{k \hsx = \hsx 1}^\infty \hsx 
\frac{(-1)^{k+1}}{(2 \pi k)^n}
\\[12pt] 
&=\ -2(n!) \hsx\sin\bigg(\frac{n \pi}{2}\bigg) \hsx (2 \pi)^{-n} \hsx \bigl(1 - 2^{1 - n}\bigr) \hsx \zeta(n)
\\[12pt]
&=\ \bigl(2^{1 - n} - 1\bigr) \hsx 2 (n!) \hsx  \sin\bigg(\frac{n \pi}{2}\bigg) (2 \pi)^{-n} \zeta(n).
\end{align*}
However
\allowdisplaybreaks
\begin{align*}
\widetilde{B}_n \hsx \ 
&=\  \hsx \widetilde{B}_n(0) 
\\[12pt]
&=\ -2(n!) \hsx 
\sum\limits_{k \hsx = \hsx 1}^\infty \hsx 
\frac{\sin\big(-\frac{n \pi}{2}\big)}{(2 \pi k)^n} 
\\[12pt]
&=\ 2(n!) \hsx \sin\bigg(\frac{n \pi}{2}\bigg) \hsx 
\sum\limits_{k \hsx = \hsx 1}^\infty \hsx \frac{1}{(2 \pi k)^n}
\\[12pt]
&=\ 2(n!) \hsx \sin\bigg(\frac{n \pi}{2}\bigg) \hsx (2 \pi)^{-n} \zeta(n).
\end{align*}
Therefore
\[
\widetilde{B}_n\bigg(\frac{1}{2}\bigg) 
\hsx = \hsx 
\bigl(2^{1 - n} - 1\bigr) \widetilde{B}_n.
\]
\end{x}
\vspace{0.3cm}

\begin{x}{\small\bf DEFINITION} \ 
Given $x \in \R$, put
\[
\Omega(x)
\hsx = \hsx 
\PV \hsx \int\limits_{-1/2}^{1/2} \hsx e^{x y}  \cot(\pi y) dy,
\]
the 
\un{omega function}.
\index{omega function}
\end{x}
\vspace{0.3cm}

\begin{x}{\small\bf \un{N.B.}} \ 
Therefore the omega function is the periodic Hilbert transform at 0 of the 1-periodic function $f$ defined by periodic extension of 
$f(y) = e^{-xy}$ $(y \in \bigl[-\frac{1}{2},\frac{1}{2}\bigr[ \hsx)$:
\allowdisplaybreaks
\begin{align*}
\Omega(x) \ 
&=\ \PV \hsx \int\limits_{-1/2}^{1/2} \hsx e^{-(0-y)x}  \cot(\pi y) dy\\[15pt]
&=\ \sH[e^{- \cdot x}] \hsx (0).
\end{align*}
\end{x}
\vspace{0.3cm}

\begin{x}{\small\bf LEMMA} \ 
There is an expansion
\[
\Omega(x)
\ = \ 
\sum\limits_{j \hsx = \hsx 0}^\infty \hsx \frac{\Omega_j}{j!} \hsx x^j,
\]
where
\[
\Omega_j \ = \ 
D_x^j \Omega(x) \bigg\arrowvert_{x=0} 
\ = \ 
\PV \hsx \int\limits_{-1/2}^{1/2} \hsx y^j  \cot(\pi y) dy.
\]
\vspace{0.2cm}

The omega function figures in the generating function for the $\widetilde{B}_n \bigg(\ds\frac{1}{2}\bigg)$.
\end{x}
\vspace{0.3cm}

\begin{x}{\small\bf THEOREM} \ 
For $\abs{x} < 2\pi$, 
\[
-\frac{x e^{x/2}}{e^x - 1} \hsx \Omega(x) 
\ = \ 
\sum\limits_{k \hsx = \hsx 0}^\infty \hsx \widetilde{B}_k \bigg(\frac{1}{2}\bigg)  \frac{x^k}{k!}.
\]
\vspace{0.2cm}

PROOF \ 
Ignoring the minus sign, on the LHS, it is a question of the Cauchy 
product of two infinite series: 

\[
\bigg( 
\sum\limits_{k \hsx = \hsx 0}^\infty \hsx B_k \bigg(\frac{1}{2}\bigg)  \frac{x^k}{k!}\bigg)
\hsx \times \hsx 
\bigg( 
\sum\limits_{k \hsx = \hsx 0}^\infty \hsx \frac{\Omega_k}{k!} x^k \bigg),
\]
a generic term being 
\[
\sum\limits_{j \hsx = \hsx 0}^k \hsx
B_{k-j} \bigg(\frac{1}{2}\bigg) \frac{x^{k-j}}{(k-j)!} \hsx  \Omega_j \frac{x^j}{j!}
\]
or still, 
\[
\bigg( 
\sum\limits_{j \hsx = \hsx 0}^k \hsx \binom{k}{j} B_{k-j} \bigg(\frac{1}{2}\bigg) \Omega_j\bigg) \frac{x^k}{k!}.
\]
Owing to the addition formula (see the Appendix to \S1), 

\[
B_k \bigg( \frac{1}{2} - y\bigg)
\ = \ 
\sum\limits_{j \hsx = \hsx 0}^k \hsx \binom{k}{j} B_{k-j} \bigg(\frac{1}{2}\bigg) (-y)^j.
\]
On the other hand, 

\[
\Omega_j
\hsx = \hsx 
\PV \hsx \int\limits_{-1/2}^{1/2} \hsx y^j \cot(\pi y) dy.
\]
And $\Omega_{2j} = 0$.  So in the sum

\[
\sum\limits_{j \hsx = \hsx 0}^k \hsx \binom{k}{j} B_{k-j} \bigg(\frac{1}{2}\bigg)  \Omega_j,
\]
only the odd $j$ contribute.  This said, consider

\[
\PV \hsx \int\limits_{-1/2}^{1/2} \  
\sum\limits_{j \hsx = \hsx 0}^k \hsx \binom{k}{j} B_{k-j} \bigg(\frac{1}{2}\bigg) y^j \cot(\pi y) dy
\]
or still, 

\[
-\PV \hsx \int\limits_{-1/2}^{1/2} \hsx 
\sum\limits_{j \hsx = \hsx 0}^k \hsx \binom{k}{j} B_{k-j} \bigg(\frac{1}{2}\bigg) (-1) y^j \cot(\pi y) dy.
\]
Assume that $j$ is odd, say $j = 2 \ell + 1$ $-$then 
\allowdisplaybreaks
\begin{align*}
(-y)^j \ 
&=\ 
(-y)^{2 \ell + 1} 
\\[12pt]
&=\ 
(-1)^{2 \ell + 1} \hsx (y)^{2 \ell + 1} 
\\[12pt]
&=\ 
(-1)^1 (y)^{2 \ell + 1}
\\[12pt]
&=\ 
(-1)y^j.
\end{align*}
The data thus reduces to 
\[
-\PV \hsx \int\limits_{-1/2}^{1/2} \hsx 
B_k \bigg(\frac{1}{2} - y\bigg) \cot (\pi y) \hsx dy 
\ \equiv \ 
-\widetilde{B}_k \bigg(\frac{1}{2}\bigg),
\]
from which the result.
\end{x}
\vspace{0.3cm}

\begin{x}{\small\bf THEOREM} \ 
\[\Omega(2 \pi x)
\hsx = \hsx 
\frac{1}{\pi} \bigg(e^{-\pi x} - e^{\pi x}\bigg) \hsx
\sum\limits_{k \hsx = \hsx 1}^\infty \hsx (-1)^k \frac{k}{x^2 + k^2}.
\]
\vspace{0.2cm}

[It can be shown that 
\[
2 \hsx \sum\limits_{k \hsx = \hsx 1}^\infty \hsx (-1)^{k+1} \hsx 
\int\limits_0^1 \hsx e^{2 \pi x y} \sin(2 \pi k y) dy
\hsx = \hsx 
\frac{1}{\pi} \bigg(e^{2\pi x} - 1\bigg) \hsx \sum\limits_{k \hsx = \hsx 1}^\infty \hsx (-1)^k \frac{k}{x^2 + k^2}
\]
or still, 
\[
e^{\pi x} \Omega(-2 \pi x) 
\hsx = \hsx 
\frac{1}{\pi} \bigg(e^{2\pi x} - 1\bigg) \hsx \sum\limits_{k \hsx = \hsx 1}^\infty \hsx (-1)^k \frac{k}{x^2 + k^2}
\]
or still, 
\[
\Omega(-2 \pi x)
\hsx = \hsx 
\frac{1}{\pi} \bigg(e^{\pi x} - e^{-\pi x}\bigg) \hsx
\sum\limits_{k \hsx = \hsx 1}^\infty \hsx (-1)^k \frac{k}{x^2 + k^2} \hsx.]
\]
\end{x}
\vspace{0.3cm}

\begin{x}{\small\bf REMARK} \ 
By way of comparison, recall that

\[
\frac{\pi}{\sin(\pi x)} 
\hsx = \hsx 
\frac{1}{x} + 2 \hsx \sum\limits_{k \hsx = \hsx 1}^\infty \hsx (-1)^k \frac{x}{x^2 - k^2}.
\]
\end{x}
\vspace{0.3cm}


%% file: _08_Zeta_2n_1.tex
\chapter{
$\boldsymbol{\S}$\textbf{8}.\quad  $\boldsymbol{\zeta(2n + 1)}$}
\setlength\parindent{2em}
\setcounter{theoremn}{0}
\renewcommand{\thepage}{Zeta Function Values \S8-\arabic{page}}

\ \indent 
The formula for $\zeta(2n)$ in terms of Bernoulli numbers (cf. \S2, \#1) admits an analog for $\zeta(2n+1)$ in terms of conjugate Bernoulli numbers.
\vspace{0.5cm}

\begin{x}{\small\bf THEOREM} \ 
\[
\zeta(2n+1) 
\ = \
(-1)^n \hsx 2^{2n} \hsx \pi^{2n+1} \hsx \frac{\widetilde{B}_{2n+1}}{(2n+1)!}.
\]

PROOF \ 
\vspace{0.3cm}

\qquad \un{Step 1:} \ $\abs{x} \ < \ 1$\\

\qquad\qquad $\implies$
\allowdisplaybreaks
\begin{align*}
\hspace{1.5cm} \sum\limits_{k \hsx = \hsx 1}^\infty \hsx (-1)^k \hsx \frac{k}{x^2 + k^2} \ 
&=\ 
\sum\limits_{k \hsx = \hsx 1}^\infty \hsx \frac{(-1)^k}{k}
\sum\limits_{n \hsx = \hsx 0}^\infty \hsx (-1)^n \bigg(\frac{x}{k}\bigg)^{2n}
\\[15pt] 
&=\ 
\sum\limits_{k \hsx = \hsx 1}^\infty \hsx \frac{(-1)^k}{k}
\sum\limits_{n \hsx = \hsx 1}^\infty \hsx (-1)^n \bigg(\frac{x}{k}\bigg)^{2n}
+ \ 
\sum\limits_{k \hsx = \hsx 1}^\infty \hsx \frac{(-1)^k}{k}
\\[15pt] 
&=\ 
\sum\limits_{n \hsx = \hsx 0}^\infty \hsx \bigg(
\sum\limits_{k \hsx = \hsx 1}^\infty \hsx (-1)^k \frac{2}{k^{2n+1}} \bigg) (-1)^n x^{2n}.
\end{align*}
\vspace{0.2cm}

\qquad \un{Step 2:} \ Write (cf. \S7, \#15)
\\[-20pt]

\allowdisplaybreaks
\begin{align*}
 \sum\limits_{k \hsx = \hsx 0}^\infty \hsx 
\frac{\widetilde{B}_k \bigl(\frac{1}{2}\bigr)}{k!} \hsx (2\pi x)^k \  
&=\ 
-2 \pi x \
\frac{e^{\pi x}}{e^{2\pi x} -1} \hsx \Omega(2 \pi x)
\\[15pt]
&=\ 
-2 \pi x \
\frac{e^{\pi x}}{e^{2\pi x} -1} \hsx \frac{e^{-\pi x}}{e^{-\pi x}} \hsx \Omega(2 \pi x)
\\[15pt]
&=\ 
-2 \pi x \
\frac{1}{e^{\pi x} - e^{-\pi x}} \hsx \Omega(2 \pi x)
\\[15pt]
&=\ 
2 \pi x \
\frac{1}{e^{-\pi x} - e^{\pi x}} \hsx \Omega(2 \pi x)
\\[15pt]
&=\ 
2 x \
\frac{\pi}{e^{-\pi x} - e^{\pi x}} \hsx \Omega(2 \pi x)
\\[15pt]
&=\ 
2 x \ 
\sum\limits_{k \hsx = \hsx 1}^\infty \hsx (-1)^k \frac{k}{x^2 + k^2} \qquad (\text{cf. \S7, \#26})
\\[15pt]
&=\ 
2 x \hsx 
\sum\limits_{n \hsx = \hsx 0}^\infty \hsx
\bigg( \sum\limits_{k \hsx = \hsx 1}^\infty \hsx (-1)^k \hsx \frac{1}{k^{2n+1}}\bigg) (-1)^n x^{2n}.
\end{align*}
Accordingly

\[
\frac{1}{2x} \hsx \sum\limits_{k \hsx = \hsx 0}^\infty \hsx 
\frac{\widetilde{B}_k \bigl(\frac{1}{2}\bigr)}{k!} \hsx (2\pi)^k x^k 
\hsx = \hsx 
\sum\limits_{n \hsx = \hsx 0}^\infty \hsx 
\bigg(\sum\limits_{k \hsx = \hsx 1}^\infty \hsx (-1)^k \frac{1}{k^{2n+1}}\bigg) (-1)^n x^{2n}.
\]

\noindent So, comparing coefficients, 

\[
\widetilde{B}_{2 n} \bigg(\frac{1}{2}\bigg) 
\hsx = \hsx 
0 \qquad (\text{cf. \S7, \#8}),
\]

\noindent and 

\[
\frac{\widetilde{B}_{2n+1} \bigl(\frac{1}{2}\bigr)}{(2n+1)!} \hsx 2^{2n} \hsx \pi^{2n+1} 
\hsx = \hsx 
(-1)^n \hsx \sum\limits_{k \hsx = \hsx 1}^\infty \hsx (-1)^k \hsx \frac{1}{k^{2n+1}}.
\]
\vspace{0.1cm}

\qquad \un{Step 3:} \ First (cf. \S7, \#10)
\[
\sum\limits_{k \hsx = \hsx 1}^\infty \hsx (-1)^k \frac{1}{k^{2n+1}} 
\hsx = \hsx 
(2^{-2n} - 1) \zeta(2n+1).
\]
Therefore
\[
\zeta(2n+1)
\hsx = \hsx 
\frac{1}{2^{-2n} - 1} \hsx (-1)^n \hsx 2^{2n} \hsx \pi^{2n+1} \hsx 
\frac{\widetilde{B}_{2n+1} (\frac{1}{2})}{(2n+1)!}.
\]
But (cf. \S7, \#11)
\[
\widetilde{B}_{2n+1} \bigg(\frac{1}{2}\bigg)
\hsx = \hsx
(2^{-2n} - 1) \widetilde{B}_{2n+1},
\]
thus 
\allowdisplaybreaks
\begin{align*}
\zeta(2n + 1) \ 
&=\ 
\frac{1}{2^{-2n} - 1} \hsx (-1)^n \hsx 2^{2n} \hsx \pi^{2n+1} \hsx 
\frac{(2^{-2n} - 1) \widetilde{B}_{2n+1}}{(2n+1)!}
\\[15pt]
&=\ 
\hsx (-1)^n \hsx 2^{2n} \hsx \pi^{2n+1} \hsx \frac{\widetilde{B}_{2n+1}}{(2n+1)!},
\end{align*}
the statement of \#1.
\end{x}
\vspace{0.3cm}

Question: \ Is
\[
\frac{\zeta(2n+1)}{\pi^{2n+1}}
\]
rational or irrational? Answer: Nobody knows.  Of course, part of the problem is the structure of $\widetilde{B}_{2n+1}$ which appears to be complicated.  
E.g.: 
\allowdisplaybreaks
\begin{align*}
\widetilde{B}_3 \bigg(\frac{1}{2}\bigg) \ 
&=\ 
\frac{\ell n(2)}{4\pi} - 2 \int\limits_{0^+}^{1/2} \hsx y^3 \cot(\pi y) dy
\\[15pt]
&=\ 
(2^{-2} - 1) \hsx \widetilde{B}_3.
\end{align*}
\vspace{0.1cm}

\begin{x}{\small\bf THEOREM} \ 
\[
\zeta(2n+1) 
\ = \
(-1)^{n+1} \hsx \frac{2^{2n}\pi^{2n+1}}{(2n+1)!} \hsx 
\int\limits_0^1 \hsx {B_{2n+1}} (y) \cot(\pi y) dy. 
\]
\vspace{0.2cm}

PROOF \ 
In fact
\allowdisplaybreaks
\begin{align*}
\widetilde{B}_{2n+1} \ 
&\equiv\ 
\widetilde{B}_{2n+1} (0) \qquad (\text{cf. \S7, \#9})
\\[15pt]
&=\ -\PV \int\limits_0^1 \hsx B_{2n+1}(y) \cot(\pi y) dy
\\[15pt]
&=\ 
-\int\limits_0^1 \hsx B_{2n+1}(y) \cot(\pi y) dy
\end{align*}
after replacing $y$ by $-y$ and taking into account the 1-periodicity.
\vspace{0.2cm}

[Note: \ The PV is not necessary since 
\[
\lim\limits_{x \hsx \ra \hsx  0} \hsx x \cot x \hsx = \hsx 1.]
\]
\end{x}
\vspace{0.3cm}

\begin{x}{\small\bf REMARK} \ 
In a similar vein, 
\[
\zeta(2n) 
\ = \
(-1)^{n+1} \hsx \frac{2^{2n-1}\pi^{2n}}{(2n)!} \hsx 
\int\limits_0^1 \hsx \widetilde{B}_{2n} (y) \cot(\pi y) dy. 
\]
\end{x}
\vspace{0.3cm}


%% file: __refs.tex
\centerline{\textbf{\large REFERENCES}}
\setcounter{page}{1}
\setcounter{theoremn}{0}
\renewcommand{\thepage}{References-\arabic{page}}
\vspace{0.75cm}

\begin{rf}
P. D'Aquino, A. Macintyre, G. Terzo, From Schanuel's Conjecture to Shapiro's Conjecture, arXiv:1206.6747 [math.NT].
\end{rf}

\begin{rf}
J. Ax, On Schanuel's Conjectures, Ann. of Math. {\bf 93}, (1971), p. 252-268.
\end{rf}

\begin{rf}
E. Burger and R. Tubbs, Making Transcendence Transparent, Springer-Verlag, 2004.
\end{rf}

\begin{rf}
Gregory Chaitin, How real are the real numbers?, arXiv:math/0411418 v 3 [math.HO] 29 Nov 2004.
\end{rf}

\begin{rf}
Cristian S. Calude, Michael J. Dinneen, Chi-Kou Shu, Computing A Glimpse of Randomness, 
arXiv:nlin/0112022 [nlin.CD].
\end{rf}

\begin{rf}
Guy Diaz, La conjecture des quatre exponentielles et les conjectures de D. Bertrand su la fonction modulaire, 
Journal de Th\'eorie des Nombres de Bordeaux, {\bf 9} (1997) p. 229-245.
\end{rf}

\begin{rf}
Guy Diaz, 
Utilisation de la conjugaison complexe dans l'\'etude de la transcendance de valeurs de la fonction exponentielle usuelle, 
Journal de Th\'eorie des Nombres de Bordeaux, {\bf 16} (2004), p. 535-553.
\end{rf}


\begin{rf}
P. Erd\' os, Representations of real numbers as sums and products of Liouville numbers, 
Michigan Math. J. {\bf 9}, p. 59-60 (1962).
\end{rf}

\begin{rf}
G. R. Everest and J. Van Der Poorten, Factorisation in the Ring of Exponential Polynomials, Proc. Amer. Math. Soc. {\bf 125}, (1997), p. 1293-1298.
\end{rf}

\begin{rf}
F. Faltin, N.  Metropolis, B. Ross, G. C. Rota, The real numbers as a wreath product, Advances in Mathematics {\bf 16} (1975) p. 278-304.
\end{rf}

\begin{rf}
N. I. Fel'dman and Yu. V. Nesterenko, Transcendental Numbers, Springer-Verlag, 1998.
\end{rf}

\begin{rf}
Lothar Sebastian Krapp, Schanuel’s Conjecture and Exponential Fields, Dissertation Univ. Oxford, (2015), 
\url{http://www.math.uni-konstanz.de/~krapp/research/Schanuels_Conjecture_and_Exponential_Fields_Errata.pdf} .
\end{rf}

\begin{rf}
F. M. S. Lima and Diego Marques, Some transcendental functions with an empty exceptional set, 
arXiv:1010.1668 [math.NT].
\end{rf}

\begin{rf}
Vincenzo Mantova, Umberto Zannier, Polynomial exponential equations and Zilber's conjeture, arXiv:1402.0685 [math.NT].
\end{rf}

\begin{rf}
Deigo Marques and Jonathon Sondow, The Schanuel Conjecture Implies Gelfond's Power Tower Conjecture, 
arXiv:1212.6931 [math.NT] (2012).
\end{rf}

\begin{rf}
Deigo Marques and Jonathon Sondow, Schanuel's Conjecture and Algebraic Powers $z^w$ and $w^z$ with $z$ and $w$ Transcendental, 
arXiv:1010.6216 [math.NT].
\end{rf}

\begin{rf}
L. A. MacColl, A factorization theory for polynomials in $x$ and in functions $e^{\alpha x}$, 
Bull. Amer. Math. Soc. {\bf 41} (2): 104-109 (February 1935).
\end{rf}

\begin{rf} 
E. Maillet, Sur quelques propri\'et\'es des nombres transcendants de Liouville, 
Bulletin de la S. M. F., \textbf{50} (1922), p. 74-99.
\end{rf}

\begin{rf} 
David Masser, Auxiliary Polynomials in Number Theory, Cambridge University Press, 2016.
\end{rf}

\begin{rf}
M. Ram Murty and Purusottam Rath, Transcendental Numbers, Springer-Verlag, 2014.
\end{rf}

\begin{rf}
Yu. V. Nesterenko, Algebraic Independence, Narosa Publishing House, 2009.
\end{rf}

\begin{rf}
J. R. Ritt, On the zeros of exponential polynomials, 
Trans, Amer. Math. Soc. \textbf{31}, (1929), p. 680-686.
\end{rf}

\begin{rf}
Damien Roy, An arithmetic criterion for the values of the exponential function, Acta Arithmetica \textbf{97} (2001), p. 183-194.
\end{rf}

\begin{rf}
R. Tijdeman, On the number of zeros of general exponential polynomials, 
Nederl. Akad. Wetensch. Proc. Set. A \textbf{74}, Indag. Math. \textbf{33} (1971), p. 1 - 7.
\end{rf}

\begin{rf}
Lou van Den Dries, Exponential Rings, Exponential Polynomials and Exponential Functions, Pac. J. Math. \textbf{113} (1984).
\end{rf}

\begin{rf}
Michel Waldschmidt, Linear Independence of logarithms and algebraic numbers, 
\url{https://webusers.imj-prg.fr/~michel.waldschmidt/articles/pdf/LIL.pdf}.
\end{rf}

\begin{rf}
Michel Waldschmidt, Schanuel's conjecture: algebraic independence of transcendental numbers, Colloquium de Giorgi, (2014).
\end{rf}

\begin{rf}
Michel Waldschmidt, Nombres Transcendants, Springer-Verlag, 1974.
\end{rf}

\begin{rf}
Michel Waldschmidt, Diophantine Approximation on Linear Algebraic Groups,  Springer-Verlag, 2000.
\end{rf}

\begin{rf}
Michel Waldschmidt, Elliptic Functions and Transcendence, 
Surveys in Number Theory, Springer-Verlag, (2008), p. 143-188.
\end{rf}

\begin{rf}
N. J. Wildberger, Real fish, real numbers, real jobs, The Mathematical Intelligencer \textbf{21} (1999), pp. 4-7.
\end{rf}